\documentclass{amsart}

\usepackage[utf8]{inputenc}
\usepackage[english]{babel}
\usepackage[T1]{fontenc}
\usepackage{appendix}

\usepackage{hyperref}

\usepackage{dirtytalk}

\usepackage{amsmath, amssymb, amsthm}   
\usepackage{amstext, amsfonts, mathrsfs}
\usepackage{epsfig}
\usepackage{epic, eepic}

\usepackage[all]{xy} 
\usepackage{xcolor}
\usepackage{centernot}
\usepackage{array}
\usepackage{enumerate}
\usepackage{babelbib}
\usepackage{graphics,graphicx}
\usepackage{tabularx}
\usepackage{cases}
\usepackage[shortlabels]{enumitem}

\usepackage{float}
\usepackage{eso-pic}
\usepackage{pst-all}
\usepackage{caption}
\usepackage{subcaption}

\usepackage[top=3cm, bottom=3cm, left=3cm, right=3cm]{geometry}

\newtheorem{thm}{Theorem}[section]
\newtheorem{prop}[thm]{Proposition}
\newtheorem{claim}[thm]{Claim}
\newtheorem{propdefi}[thm]{Proposition and Definition}
\newtheorem{lem}[thm]{Lemma}
\newtheorem{lemdef}[thm]{Lemma and Definition}
\newtheorem{coro}[thm]{Corollary}

\newtheorem{mainthm}{Theorem}

\newtheorem{thmintro}{Theorem}
\newtheorem{propintro}[thmintro]{Proposition}

\theoremstyle{remark}
\newtheorem{rmk}[thm]{Remark}

\theoremstyle{definition}
\newtheorem{defi}[thm]{Definition}
\newtheorem{nota}[thm]{Notation}
\newtheorem{example}[thm]{Example}

\numberwithin{equation}{section}

\def\l{\ensuremath{\left}}
\def\r{\ensuremath{\right}}

\def\nn{\noindent}

\def\qq{\nn\quad}

\def\emptyset{{\varnothing}}

\def\ssm{{\smallsetminus}}
\def\inv{\ensuremath{{-1}}}

\newcommand{\mb}[1]{\mbox{#1}}

\newcommand{\res}[2]
{\ensuremath{#1}|_{\ensuremath{#2}}}

\def\N{{\mathbb N}}    
\def\Z{{\mathbb Z}}    
\def\R{{\mathbb R}}

\def\H{{\mathbb H}}

\def\S{{\mathbb S}}   
\def\T{{\mathbb T}}    
\def\D{{\mathbb D}}

\DeclareMathAlphabet{\mathcal}{OMS}{cmsy}{m}{n}

\def\cA{{\mathcal A}}  \def\cG{{\mathcal G}} \def\cM{{\mathcal M}}    \def\cH{{\mathcal H}} \def\cN{{\mathcal N}} \def\cT{{\mathcal T}} \def\cC{{\mathcal C}}   \def\cO{{\mathcal O}} \def\cU{{\mathcal U}} \def\cD{{\mathcal D}}    \def\cV{{\mathcal V}}   \def\cK{{\mathcal K}} \def\cQ{{\mathcal Q}} \def\cW{{\mathcal W}} \def\cF{{\mathcal F}}  \def\cL{{\fontfamily{cmr}\selectfont{\mathcal L}}} \def\cR{{\mathcal R}}

\def\scA{{\mathscr A}}         \def\scT{{\mathscr T}}

        \def\bN{{\mathbf N}}

\def\gA{{\mathfrak A}}    \def\gS{{\mathfrak S}}                        \def\gR{{\mathfrak R}}

\newcommand{\Id}{\operatorname{Id}}

\newcommand{\Vect}{\operatorname{Vect}}
\newcommand{\adh}{\operatorname{adh}}

\def\intr{\operatorname{int}} 
\def\dist{\operatorname{dist}}
\def\mod{\operatorname{mod}}
\def\cst{\operatorname{\textit{cst}}}
\def\max{{\operatorname{max}}}

\def\qms{quasi-Morse--Smale}

\def\qt{quasi-transverse}
\def\qts{quasi-transverse}
\def\bb{building block}
\def\bbs{building blocks}
\def\vf{vector field}
\def\vfs{vector fields}

\def\minc{minimal continuation}
\def\ncs{normalized coordinate system}
\def\ncss{normalized coordinate systems}
\def\lcs{linearizing coordinate system}

\def\acs{affine coordinate system}

\def\msis{maximal saddle invariant set}

\def\nbh{neighborhood}
\def\nbhs{neighborhoods}

\newcommand{\iin}{\ensuremath{\mathrm{in}}}
\newcommand{\out}{\ensuremath{\mathrm{out}}}
\def\tan{\ensuremath{\mathrm{tan}}}

\def\uu{\ensuremath{{uu}}}
\def\ss{\ensuremath{{ss}}}
\def\cu{\ensuremath{{cu}}}
\def\cs{\ensuremath{{cs}}}

\def\bby{B\'eguin--Bonatti--Yu}
\def\sqt{strongly quasi-transverse}
\def\st{strongly transverse}

\def\diff{diffeomorphism}
\def\diffs{diffeomorphisms}
\def\homeo{homeomorphism}
\def\homeos{homeomorphism}

\def\paif{pair of affine invariant foliations}

\def\cc{connected component}
\def\ccs{connected components}

\def\foutin{\ensuremath{f_{\out, \iin}}}

\def\pP{{\ensuremath{\partial P}}}

\def\Pout{\ensuremath{P^\out}}
\def\Pin{\ensuremath{P^\iin}}

\def\tPout{\ensuremath{\hat P^\out}}
\def\tPin{\ensuremath{\hat P^\iin}}

\def\ppP{{\ensuremath{\partial_+ P}}}

\def\loc{\ensuremath{\mathrm{loc}}}
\def\Mod{{\ensuremath{\mathrm{Mod}}}}

\newcommand{\led}{{\ensuremath{\lambda, \epsilon, \delta}}}
\newcommand{\ed}{{\ensuremath{\epsilon, \delta}}}

\def\ra{\ensuremath{\rightarrow}}
\def\la{\ensuremath{\leftarrow}}
\def\ua{\ensuremath{\uparrow}}
\def\da{\ensuremath{\downarrow}}

\def\imb{intermediate block}
\def\imbs{intermediate blocks}

\begin{document}


\title[Anosov flows in dimension 3 from gluing building blocks]{Anosov flows in dimension 3 from gluing building blocks with quasi-transverse boundary}
\author{Neige Paulet}
\address{Department of Math \& Stats, Queen's University, 48 University Ave., Jeffery Hall, Kingston, ON K7L 3N6, Canada}
\email{neige.paulet@queensu.ca}

\begin{abstract}
We prove a new result allowing to construct Anosov flows in dimension 3 by gluing \textit{building blocks}.
By a building block, we mean a compact 3-manifold with boundary $P$, equipped with a $\cC^1$ vector field $X$, such that the maximal invariant set $\bigcap_{t \in \R} X^t (P)$ is a saddle hyperbolic set, and the boundary $\partial P$ is \textit{quasi-transverse} to $X$, i.e., transverse except for a finite number of periodic orbits contained in $\partial P$.
Our gluing theorem is a generalization of a recent result of F.~B\'eguin, C.~Bonatti, and B.~Yu who only considered the case where the block does not contain attractors nor repellers, and the boundary $\pP$ is transverse to $X$.
The quasi-transverse setting is much more natural.
Indeed, our result can be seen as a counterpart of a theorem by Barbot and Fenley which roughly states that every 3-dimensional Anosov flow admits a canonical decomposition into building blocks (with quasi-transverse boundary).
We will also show a number of applications of our theorem. 
\end{abstract}


\vspace*{1em}
\maketitle

\setcounter{tocdepth}{1}

\tableofcontents

\section*{Introduction}

The definition of Anosov flows appeared in the 1960's, when D. Anosov studied the qualitative dynamical properties of the geodesic flow on Riemannian manifolds with negative curvature in his famous article \cite{anosovGeodesicFlowsClosed1967}, which became the prototype example of an Anosov flow.
The understanding of Anosov flows, together with their generalizations and their relationships to other fields and structures, has led to many discoveries in dynamical systems, geometry and geometric topology.
As Anosov flows are structurally stable, one can hope to obtain a complete classification of their orbital equivalence classes on a given manifold by a finite number of combinatorial invariants.
Such a classification in dimension~3 (the minimal dimension for an Anosov flow) is still today the motivation of many works.
They revealed remarkable interactions between the dynamics of the flow and the topology of the underlying manifold.

For certain classes of 3-manifolds (torus bundles over the circle, or Seifert manifolds for example), the topology of the manifold almost completely determines the dynamics of the Anosov flow it can carry, up to orbital equivalence.
On the other hand, there are 3-manifolds which carry several non-orbit equivalent Anosov flows.

The main object of this paper is the construction of new examples of Anosov flows in 3-manifolds.

\subsection*{Previous constructions of Anosov flows}
\label{sec: construction decomposition anosov}
The two standard examples of Anosov flows are geodesic flows of a negatively curved manifolds and suspensions of hyperbolic toral automorphisms. 
These were the only examples of Anosov flows for decades. 
These two flows can also be described as the action of a one-parameter subgroup of a Lie group $G$, acting on a quotient of $G$ by a co-compact lattice.
Such a flow is called an \emph{algebraic flow}.
In the 1980's the first examples of non-algebraic Anosov flows appeared.
These examples show that the problem of classification of Anosov flows in dimension~3 cannot be easily reduced, because of \emph{flexibility} phenomena: there are many Anosov flows which do not have the qualities common to algebraic flows, and there are manifolds carrying many Anosov flows whose dynamics are not equivalent, or even very distant.

\begin{itemize}[leftmargin=0.8cm]
    \item J. Franks and B. Williams construct in \cite{franksAnomalousAnosovFlows1980} the first non-transitive Anosov flow in dimension~3 by the \emph{Blow up -- Excise -- Glue} surgery.
    
    \item  M. Handel and W.P. Thurston construct in \cite{handelAnosovFlowsNew1980} the first non-algebraic transitive Anosov flow by performing a \emph{cut-and-paste} surgery on the geodesic flow of a hyperbolic surface.

    \item C. Bonatti and R. Langevin construct in \cite{bonattiExempleFlotAnosov1994} the first example of transitive Anosov flow, transverse to an embedded torus, but not equivalent to a suspension.
    The two previous examples are obtained by gluing the boundary components of \bbs{} that have been \emph{cut out of a (blown-up) Anosov flow}. 
    The Bonatti--Langevin example is the first Anosov flow constructed from a \bb{} which is not obtained by surgery on an initial standard Anosov flow.
    T.~Barbot generalizes this construction in~\cite{barbotGeneralizationsBonattiLangevinExample1998}.

    \item F. B\'eguin, C. Bonatti, and B. Yu have developed a general procedure for constructing Anosov flows by gluing ``abstract building blocks'' in \cite{beguinBuildingAnosovFlows2017}, in the spirit of the Bonatti--Langevin construction but where the manifold and the dynamics are not explicit.
    A \emph{\bby{} block} is a pair $(P,X)$ where $P$ is a compact manifold with boundary equipped with a $\cC^1$ vector field $X$, transverse to $\partial P$, and such that the maximal invariant set $\bigcap_{t \in \R} X^t (P)$ forms a hyperbolic set with one-dimensional strong stable and unstable bundles.
    The authors show that under very general conditions, there is a way to glue the boundary components of $P$ via a \diff{} $\varphi \colon \partial P \to \partial P$ which pairs the components of the \emph{exit boundary} (along which the flow exits $P$) with the components of the \emph{entrance boundary} (along which the flow enters $P$) to obtain a closed manifold $P_\varphi := P/\varphi$ equipped with a vector field $X_\varphi$ induced by $X$ which is Anosov.
\end{itemize}

The gluing procedure of \bby{} is a powerful technique to show the flexibility of Anosov flows, and allows to build Anosov dynamics on manifolds with a rich and complicated topology.
They construct for example the first family of manifolds $\cM_N$ carrying $N$ pairwise non-orbit equivalent Anosov flows for every integer $N$.
However for Anosov flows constructed with this procedure, there always exist embedded tori which are transverse to the flow. 
The existence of such tori is not a common property.
For example, the unit tangent bundle of a closed hyperbolic surface contain plenty of incompressible tori, but one can easily prove that none of these is isotopic to a torus transverse to the geodesic flow.

\vspace*{-0.8em}
\subsection*{Decomposition of an Anosov flow: the modified JSJ decomposition}
The counterpart of the \bby{} gluing procedure is to find a good way to decompose a given Anosov vector field $Z$ on a closed orientable 3-manifold $\cM$ along tori into building blocks $(P_i, X_i)$.

The Jaco--Shalen--Johannson (JSJ) theorem states that any irreducible connected orientable closed manifold $\cM$ of dimension~3 can be cut along a minimal finite collection, unique up to isotopy, of incompressible tori $\cT = \{T_1, \dots, T_n \}$ into pieces $P_1, \dots, P_m$ such that each $P_i$ is either atoroidal or admits a Seifert fibration. 
Any 3-manifold $\cM$ which carries an Anosov vector field $X$ is irreducible (since its universal covering is $\R^3$), and thus admits a JSJ-decomposition along incompressible tori.
In \cite{barbotPseudoAnosovFlowsToroidal2013}, the authors study in detail the ``optimal'' position of an incompressible torus $T$ embedded in an orientable manifold $\cM$ with respect to the Anosov vector field $X$ on $\cM$.
We will say that the torus $T$ is \emph{quasi-transverse} to $X$ if $T$ contains a finite (possibly zero) number of periodic orbits $\cO_* = \{\cO_1, \dots, \cO_n\}$ of the flow,
it is transverse to $X$ on the complementary of the orbits $\cO_*$, and the transverse orientation given by the vector field $X$ on two adjacent components of $T \ssm \cO_*$ never coincide.

T. Barbot and S. Fenley show (\cite[Theorem 6.10]{barbotPseudoAnosovFlowsToroidal2013}) that any incompressible torus embedded in a 3-manifold $\cM$ carrying an Anosov vector field $X$ is homotopic to a torus quasi-transverse to $X$ and \emph{weakly embedded} (embedded outside the periodic orbits of $X$ contained in the torus).
This result gives a \emph{modified JSJ decomposition} of an Anosov flow (see \cite[Section 2.2]{barthelmeCountingPeriodicOrbits2017} for a precise statement) which is unique up to homotopy along the flow, and only depends on the topology of the underlying manifold.
With the help of this decomposition, S. Fenley and T. Barbot have started the intensive study of Anosov flows on toroidal 3-manifolds (\cite{barbotPseudoAnosovFlowsToroidal2013}) and their classification in restriction to JSJ Seifert pieces and to some graph manifolds (\cite{barbotClassificationRigidityTotally2015, barbotFreeSeifertPieces2021, barbotOrbitalEquivalenceClasses2022}).
Let us stress that Anosov flows on atoroidal manifolds or in restriction to atoroidal pieces of the JSJ decomposition are still very poorly understood.

In conclusion, up to some technical details, any Anosov flow in a toroidal manifold can be canonically decomposed into blocks $(P_i,X_i)$ where $P_i$ is a manifold with boundary, and $X_i$ is a $\cC^1$ vector field on $P_i$ which is quasi-transverse to the boundary $\partial P_i$.
It is therefore natural to try to generalize the \bby{} construction for such ``quasi-transverse building blocks.''
This is the goal of this paper.

\subsection*{The Gluing Theorem}

In our setting, a \emph{building block} $(P,X)$ is a compact 3-manifold $P$ with boundary equipped with a $\cC^1$ vector field $X$, such that the maximal invariant set $\Lambda$ of the generated flow is (saddle) hyperbolic, and whose boundary $\pP$ is \emph{quasi-transverse} to the vector field.
We provide very general sufficient conditions to glue the boundary components of building blocks, so that the resulting manifold is closed and equipped with an Anosov flow induced by the initial vector field.
A key element is the one-dimensional lamination $\cL$ induced by the trace of the stable and unstable manifold of $\Lambda$ on the boundary $\pP$, and its image by a \emph{gluing map} $\varphi: \pP \to \pP$.
We say that the block is \emph{filled} if $\cL$ is a \emph{filling lamination} (i.e., it extends to a foliation on $\pP$ with a technical additive assumption) and that the gluing map is \emph{\sqt{}} if the pair of laminations $(\cL, \varphi_* \cL)$ extend to a pair of foliation transverse except along the boundary orbits $\cO_*$ where they coincide. 
These assumption (except the technical one) are necessary to obtain an Anosov flow by such gluing procedure.
To do so, we want to perturb the initial \vf{} and the boundary of the block by isotopy, and to perturb the initial gluing map by isotopy so as to preserve the \say{drawing} of the transverse intersection of the pair of laminations $(\cL, \varphi_* \cL)$.
We will say that such perturbation is a \emph{strong isotopy of the triple} $(P,X,\varphi)$.

\begin{thmintro}[Gluing Theorem, Theorem~\ref{thm: gluing theorem}] \label{thmintro: gluing theorem}
Let $(P,X)$ be a filled building block, and $\varphi$ be a \sqt{} gluing map of $(P,X)$.
There exists a triple $(P_1, X_1, \varphi_1)$ strongly isotopic to $(P,X,\varphi)$ such that
$X_1$ induces an Anosov vector field on the closed 3-manifold $P_1/\varphi_1$.
\end{thmintro}

This theorem is the analog of the \bby{} gluing theorem (\cite[Theorem~1.5]{beguinBuildingAnosovFlows2017}) in the case where the set $\cO_*$ of periodic orbits of $X$ contained in $\pP$ is empty, and extends it to blocks containing attractors or repellers.
Moreover it includes the Franks--Williams surgery, the Handel--Thurston surgeries and their generalizations. 
Finally, it allows to consider the most natural building blocks for Anosov flows in dimension~3, in the sense that the quasi-transverse position is the ``optimal'' position of an incompressible closed surface embedded in an Anosov flow according to the work of T. Barbot and S. Fenley mentioned above.

\subsubsection*{Transitivity criterion}

We show a transitivity criterion of an Anosov vector field $X_\varphi$ on a closed manifold $P_\varphi$ obtained by gluing a \bb{} $(P,X)$ by a \sqt{} gluing map $\varphi$.
We associate an oriented graph $G = G(P,X,\varphi)$, analogous to the Smale's graph in the following way.
The vertices are the basic pieces $\Lambda_i$ of the maximal hyperbolic invariant set $\Lambda$ of $(P,X)$, and there exists a directed edge from $\Lambda_i$ to $\Lambda_j$ if and only if $\cW^u(\Lambda_i)$ intersects $\cW^s(\Lambda_j)$, or $\varphi(\cW^u(\Lambda_i))$ intersects $\cW^s(\Lambda_j)$.
We say that an oriented graph is \textit{strongly connected} if each pair of vertex can be connected by a oriented path of edges.

\begin{propintro}[Proposition~\ref{prop: transitivity criterion}] \label{propintro: transitivity criterion}
If the graph $G(P,X,\varphi)$ is strongly connected, then the Anosov vector field $X_\varphi$ on $P_\varphi$ is transitive.
\end{propintro}

\subsection*{Applications of the Gluing Theorem}
We show a number of application of Theorem~\ref{thmintro: gluing theorem}.

\subsubsection*{Realizing bi-foliation in Anosov flows}

As first application, we prove that any type of \emph{quasi-transverse bi-foliation} on a torus can be realized as the trace of the stable and unstable foliation of a transitive Anosov flow on an embedded \qt{} torus.
If $(P,X)$ is a \bb{}, then the boundary lamination $\cL$ of $\pP$ is a generalization of a \emph{Morse--Smale type} lamination on a surface, called \emph{quasi-Morse--Smale}.
It contains a finite number of compact leaves, and each half non-compact leaf accumulates on a compact leaf.
Among the compact leaves, we distinguish the compact leaves induced by the periodic orbits of $X$ in $\pP$, which we call \emph{marked compact leaves}.
If $(\cF_1, \cF_2)$ is a pair of quasi-Morse--Smale foliations on the oriented torus $\T^2$, such that $\cF_1$ and $\cF_2$ are tangent along their marked leaves and transverse to each other on the complementary of the marked leaves, they are said to be \emph{\qt{}}.
Then we can associate to the pair a finite combinatorial data $\sigma = \sigma(\cF_1, \cF_2)$, which encodes the (cyclic) order of the compact leaves of $\cF_1$ and $\cF_2$ on $\T^2$, and the \say{type} of holonomy of each one.
We say that $\sigma$ is a \emph{combinatorial type of the bifoliation}.
We prove

\begin{thmintro}[Theorem~\ref{thm: realize qms bifoliation}] \label{thmintro: realize qms bifoliation}
Let $\sigma$ be a combinatorial type of a \qt{} bifoliation.
There exists a transitive Anosov vector field $Z$ on an oriented $3$-manifold $\cM$ and an incompressible torus $T$ embedded in $\cM$, \qt{} to $Z$, such that the trace of the stable and unstable foliation $\cF^s$ and $\cF^u$ on $T$ induces a bifoliation $(\cF_1, \cF_2)$ on $T$ of combinatorial type $\sigma$.
\end{thmintro}

\subsubsection*{Embedding block in Anosov flows}
We prove that any filled building block can be embedded in an Anosov flow.

\begin{thmintro}[Theorem~\ref{thm: embed block in anosov}] \label{thmintro: embed block in anosov}
For any (transitive) orientable filled block $(P,X)$, there exists a (transitive) Anosov vector field $Z$ on a closed orientable 3-manifold $\cM$, such that $(P,X)$ is embedded in $(\cM, Z)$.
More precisely, there exists a finite collection of incompressible tori $\cT$ embedded in $\cM$, \qt{} to $Z$, such that the closure of one \cc{} of $\cM \ssm \cT$ is a compact submanifold diffeomorphic to $P$ and such that the restriction of $Z$ on $P$ is orbit equivalent to~$X$.
\end{thmintro}

\subsubsection*{Realizing geometric type in building blocks}

We study the dynamical properties of building blocks, giving necessary and sufficient criteria for an abstract geometric type to be realized on a Markov partition in a nice building block. 
The \emph{geometric type of a Markov partition}, introduced by C. Bonatti and R. Langevin, is a finite combinatorial data $\scT$ which encodes the geometry of the intersection of the rectangles of a Markov partition $\cR$ under the action of the dynamic.
F. B\'eguin and C. Bonatti show (\cite{beguinFlotsSmaleDimension2002}) that this data characterizes the \emph{germ} of a vector field $X$ on a \emph{saddle saturated} hyperbolic compact invariant set $K$.
They show that the germ of $X$ along $K$ determines a unique (up to orbit equivalence) orientable manifold with boundary $U$, provided with a vector field $Y$ transverse to the boundary $\partial U$, which realizes the germ of $X$ along $K$, and this manifold is in a sense the \say{simplest} topologically speaking.
The pair $(U,Y)$ is a building block, and one will then speak of \emph{model block}.
B\'eguin, Bonatti and Veitiez construct in \cite{beguinConstructionFlotsSmale1999}, for a given abstract geometric type $\scT$, a model block $(U,Y)$ whose maximal invariant admits a Markov partition $\cR$ of geometric type $\scT$, which is unique up to orbit equivalence.
We then speak of \emph{the model of the geometric type $\scT$}.
We are interested in the question of the realizability of an abstract geometrical type $\scT$ by a filled \bb{}.
We show the following necessary and sufficient criterion.

\begin{thmintro}[Theorem~\ref{thm: geometric type in block}] \label{thmintro: geometric type in block}
Let $\scT$ be an abstract geometric type.
There exists an orientable filled block $(P,X)$, admitting a Markov partition of geometric type $\scT$ if and only if the model $(U, Y)$ of $\scT$ satisfies the following conditions;
\begin{enumerate}
    \item 
    $\partial U$ is a union of tori and spheres, each sphere contains two open disks $D_i$ and $D_j$, disjoint from the boundary lamination $\cL_Y$, bounded by two distinct compact leaves of $\cL_Y$;
    \item
    $\cL_Y$ is a pre-foliation on the complementary $\partial U \ssm \bigcup_i D_i$.
\end{enumerate}
Such a block $(P,X)$ is then unique up to isotopy.
\end{thmintro}


The two conditions can be checked using a simple algorithm for attaching the rectangles of the associated Markov partition $\cR$, which allows to compute the boundary of the model $U$ and the boundary lamination $\cL_Y$.

\subsubsection*{Orbit complement as JSJ pieces of Anosov flows}
We show that one can realize periodic orbits complement of transitive (pseudo)-Anosov flow as JSJ piece of transitive Anosov flow.
A \emph{pseudo-Anosov flow} is a generalization of an Anosov flow where we allow a finite number of singularities of stable and unstable foliations of $p$-prong type, $p \geq 3$.

\begin{thmintro}[Theorem~\ref{thm: orbit complement and JSJ piece}]
\label{thmintro: orbit complement and JSJ piece}
Let $\Gamma = \{ \gamma_1, \dots, \gamma_n \}$ be a finite collection of periodic orbits of a transitive pseudo-Anosov vector field $X$ on an orientable 3-manifold $\cM$.
Assume that all the singular orbits of $X$ are contained in $\Gamma$ and that the complementary $\cM \ssm \Gamma$ is atoroidal.
Then there exists an orientable 3-manifold $\cN$ carrying a transitive Anosov vector field $Y$ such that the JSJ decomposition of $\cN$ is made of two atoroidal pieces $P$ and $P'$, both homeomorphic to $\cM \ssm \Gamma$, and a periodic Seifert piece.
The restriction of $Y$ to $P$ and $P'$ is obtained from $X$ by a DA\footnote{Derived from Anosov.} bifurcation on the orbits of $\Gamma$.
\end{thmintro}

A Seifert piece in a 3-manifold $\cM$ carrying an Anosov vector field $X$ is said to be \emph{periodic} if there exists a Seifert fibration for which the regular fiber is homotopic to a power of a periodic orbit of the flow of $X$.
Those pieces are well understood from works of Barbot and Fenley (\cite{barbotClassificationRigidityTotally2015}).

As a corollary of this theorem, we obtain sufficient conditions for the complementary of a \emph{hyperbolic knot} in $\S^3$ to be realized as an atoroidal JSJ piece of a transitive Anosov flow.
Here is an example of such knots.
A \textit{plumbing} (or Murasugi sum, see \cite{gabaiPseudoAnosovMapsSurgery1990}) $S = S_1 \# S_2$ is a Seifert surface obtained by gluing two Seifert surfaces $S_i$ of unlinked knots $\partial S_i$ in $\S^3$ along a disk.

\begin{thmintro}[Corollary~\ref{coro: plumbing of figure eight knot}]
\label{thmintro: plumbing of figure eight knot}
Let $K = \partial (S_1 \# S_2 \# \dots \# S_n)$ be a plumbing of $n$ copies of the Seifert surface of the figure eight knot.
Then $K$ is a hyperbolic fibered knot and the complementary of $K$ is an atoroidal JSJ piece of a manifold carrying a transitive Anosov flow.
\end{thmintro}

\subsubsection*{Gluing pieces of skewed Anosov flows}
We show that, under minimal conditions, one can glue pieces cut along a collection of incompressible tori embedded in a skewed $\R$-covered Anosov flow.
An Anosov flow on a 3-manifold $\cM$ is said to be \emph{$\R$-covered} if the leaf space of the lifted stable foliations $\widetilde \cF^s$ on the universal cover $\widetilde \cM$ is separated (hence homeomorphic to $\R$).
It is said to be \emph{skewed $\R$-covered} if it is moreover not orbit equivalent to a suspension.
A result of Barbot (\cite{barbotMisePositionOptimale1995}) allows us to cut building blocks out of a a skewed $\R$-covered Anosov flow along a collection of incompressible tori. 
We call such piece a \emph{skewed $\R$-covered Anosov block}.

\begin{thmintro}[Theorem~\ref{thm: gluing skewed blocks}] \label{thmintro: gluing skewed blocks}
Let $(P,X)$ and $(P', X')$ be two skewed $\R$-covered Anosov blocks and $\varphi \colon \pP \to \pP'$ a gluing map.
There is a gluing map $\psi$ isotopic to $\varphi$ among gluing maps such that the \vf{} $Z$ induced by $X$ and $X'$ on $P \cup P' / \psi$ is Anosov.
\end{thmintro}

Note that this statement does not require any assumption on the action of the gluing map on the boundary laminations.
Any \emph{piece of geodesic flow} or finite covering of geodesic flow is a skewed $\R$-covered Anosov block.
Hence Theorem~\ref{thmintro: gluing skewed blocks} generalizes the Handel--Thurston construction and their generalization and gets rid of the \emph{positivity constraint} on the isotopy class of the gluing map.
Note, however, that it allows to use blocks that are much more general than geodesic flow pieces, and that are not \emph{a priori} cut out from the same Anosov flow, inside the very rich family of skewed $\R$-covered Anosov flows.

\subsection*{Structure of the paper}

The paper is organized as follows.
The first part is devoted to the proof of Theorem~\ref{thmintro: gluing theorem}.
In the preliminary Section~\ref{sec: preliminaries}, we formally define building blocks, gluing maps, and fundamental properties.
In Section~\ref{sec: normalization}, we show that we can put a triple $(P,X,\varphi)$ candidate for Theorem~\ref{thmintro: gluing theorem} in a \say{normalized} form.
In Section~\ref{sec: crossing map}, we study the hyperbolic properties of the crossing map of the flow from the entrance boundary to the exit boundary of the block.
In Section~\ref{sec: spreading}, we show how to \say{spread hyperbolicity} by a process of coordinate change on the boundary of $P$.
In Section~\ref{sec: parameters and cones}, we show that we can use this change of coordinates to modify the gluing map $\varphi$ in order to create hyperbolicity along the new recurrent orbits of the flow obtained after gluing, and this in a way compatible with the natural hyperbolicity of the initial flow.
In Section~\ref{sec: proof gluing thm}, we show that for such a choice of gluing map, the flow induced by the initial flow on the glued manifold is Anosov, which completes the proof of Theorem~\ref{thmintro: gluing theorem}. We also show the transitivity criterion stated in Proposition~\ref{propintro: transitivity criterion}.

The second part is devoted to the applications of the Gluing Theorem.
In Section~\ref{sec: prescribed boundary lamination}, we show Theorem~\ref{thmintro: realize qms bifoliation}.
In Section~\ref{sec: embed block in anosov flow}, we show Theorem~\ref{thmintro: embed block in anosov}.
In Section~\ref{sec: geometric type}, we show Theorem~\ref{thmintro: geometric type in block}.
In Section~\ref{sec: orbit complement}, we show Theorem~\ref{thmintro: orbit complement and JSJ piece}.
Finally, we show Theorem~\ref{thmintro: gluing skewed blocks} in Section~\ref{sec: gluing skewed blocks}.

\section{Preliminaries: building blocks and gluing maps}
\label{sec: preliminaries}
\subsection{Building block}
\label{sec: preli; subsec: building block}

A (local) \textit{quadrant} of a hyperbolic periodic orbit $\cO$ is a connected component of the complement of the local stable and unstable manifolds of $\cO$ in a \nbh{} of $\cO$.
Two quadrants are \textit{opposite} if the intersection of their closure is $\cO$.

\begin{defi} [Surface quasi-transverse to a vector field]\label{def: qt surface}
Let $S$ be a transversely orientable surface $\cC^1\!$-immersed in a $3$-dimensional manifold $M$ \linebreak[4] equipped with a $\cC^1$ vector field $X$.
We say that $S$ is \emph{\qt{} to~$X$}~if
\begin{enumerate}
    \item \label{def: qt surface; it: tangent orbit}
    $S$ contains a finite collection \mb{$\cO_* = \{ \cO_1, \dots, \cO_n \}$} of hyperbolic periodic orbits of $X$;
    \item \label{def: qt surface; it: transverse}
    $X$ is transverse to $S \ssm \cO_*$;
    \item \label{def: qt surface; it: in and out}
    each orbit $\cO_i \in \cO_*$ is two-sided in $S$, and the local sides lie in opposite quadrants.
\end{enumerate}
\end{defi}

Item \ref{def: qt surface; it: in and out} implies that the orbits of $X$ cross the two (local) sides in two opposite directions, or more formally the transverse orientation induced by the \vf{} $X$ on the two sides of $\cO_i$ never simultaneously coincide with a global transverse orientation of the surface $S$.

\begin{figure}[htb]
    \centering
    \vspace*{-2em}
    \includegraphics[height=0.25\textheight]{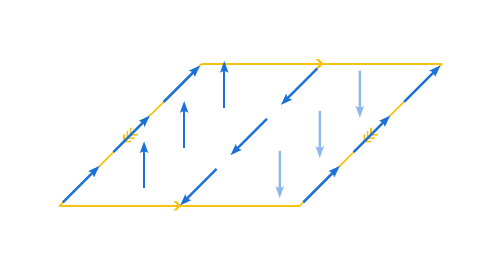}
    \vspace*{-2em}
    \caption{A torus quasi-transverse to a vector field $X$ containing two periodic orbits}
\end{figure}

In particular, a surface $S$ transverse to a vector field $X$ is a \qt{} surface with $\cO_* = \emptyset$.
We denote $X^t$ the flow generated by $X$.
Recall that a compact set $\Lambda$ invariant by the flow $X^t$ on a manifold $M$ is said to be \emph{hyperbolic of index $(1,1)$} for $X$ if there exists an $X^t$-invariant splitting of the tangent space of $M$ over $\Lambda$ into a sum
\mb{$\res{TM}{\Lambda} = E^\ss \oplus \R. X \oplus E^\uu$} of 1-dimensional subbundles,
and constants $\lambda >1$, and $C>0$ such that for a Riemannian metric on $M$, we have
\begin{itemize}[--]
    \item $ \forall v \in E^\uu$, $\forall t \geq 0$, $\Vert (X^t)_* v \Vert \geq C \lambda^t v \Vert$;
    \item $\forall v \in E^\ss$, $\forall t \leq 0$, $\Vert (X^t)_* v \Vert \geq C \lambda^{-t} \Vert v \Vert$.
\end{itemize}

\begin{defi}[Building block]
\label{def: building block}
Let $P$ be a compact 3-dimensional manifold with boundary, equipped with a $\cC^1$ vector field $X$.
We say that the pair $(P,X)$ is a \emph{building block} (or more simply a \emph{block}) if
\begin{enumerate}
    \item the boundary $\partial P$ is \qt{} to the vector field $X$;
    \item the maximal invariant set of the flow of
    $X$ in $P$, denoted $\Lambda := \bigcap_{t \in \R} X^t (P)$, is an index $(1,1)$ hyperbolic set for the flow of $X$.
\end{enumerate}
\end{defi}

\begin{rmk}
The manifold $P$ is not necessarily connected.
Building blocks should be thought as the basic pieces of a building game, our goal being to build some Anosov flows by gluing a collection of such blocks together. 
From a formal point of view, a finite collection of \bbs{} can always be viewed as a single non-connected \bb{}.
\end{rmk}

\begin{rmk} \label{rmk: block and bby block}
Let $(P,X)$ be a building block.
In the case where the collection $\cO_*$ of periodic orbit contained in $\pP$ is empty, the boundary is transverse to the vector field $X$ and one recover the definition of a \emph{hyperbolic plug} in the sense of \cite[Definition~3.1~and~Definition~3.2]{beguinBuildingAnosovFlows2017}.
    We will then say that $(P,X)$ is a \emph{\bby{} block}.
    In this case the maximal invariant set $\Lambda$ is contained in the interior of $P$, which is no longer true if the collection $\cO_*$ is non-empty.
\end{rmk}

\begin{defi}[Entrance boundary and exit boundary]
Let $(P,X)$ be a \bb{} and $\cO_*$ be the set of periodic orbits of $X$ contained in $\pP$.
We call the \emph{entrance boundary} of $(P,X)$, denoted $\Pin$, the subset of $\pP \ssm \cO_*$ where the flow enters $P$, and the \emph{exit boundary} of $(P,X)$, denoted $\Pout$, the subset of $\pP \ssm \cO_*$ where the flow exits $P$.
\end{defi}

\begin{figure}[htb]
    \centering
    \vspace*{-1em}
    \includegraphics[height=0.3\textheight]{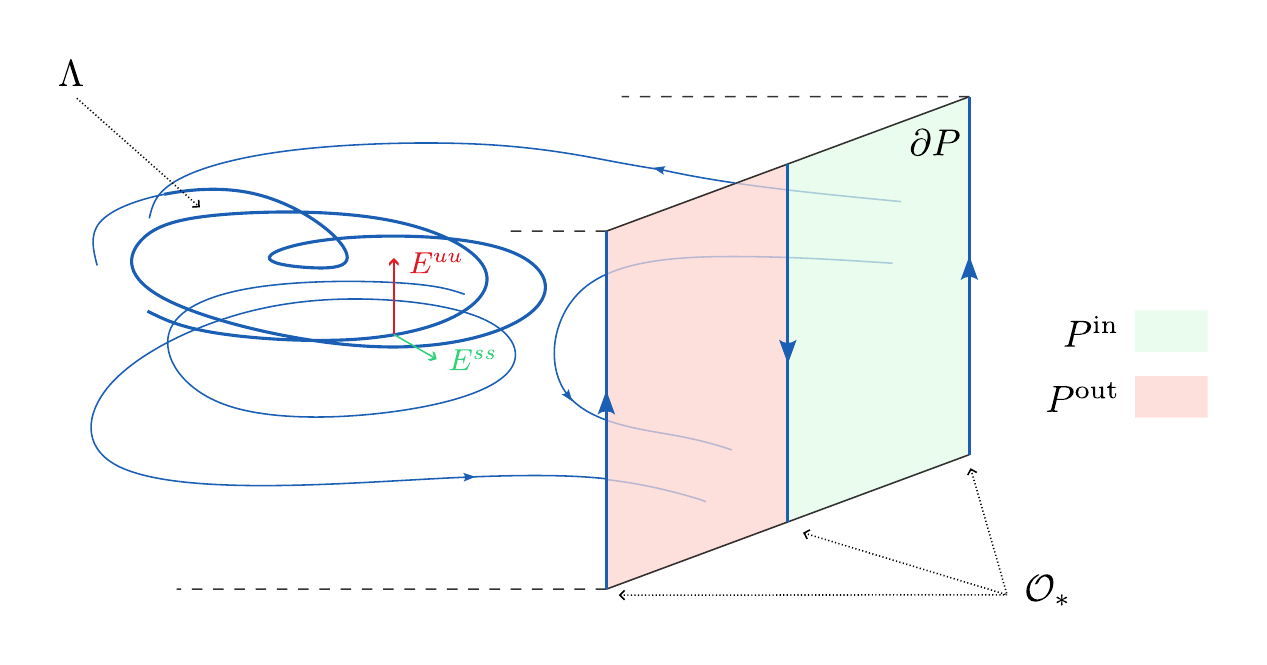}
    \vspace*{-1em}
    \caption{A building block $(P,X)$}
    \label{fig: building block}
\end{figure}

\begin{figure}[htb]
    \centering
    \vspace*{-2em}
    \includegraphics[height=0.32\textheight]{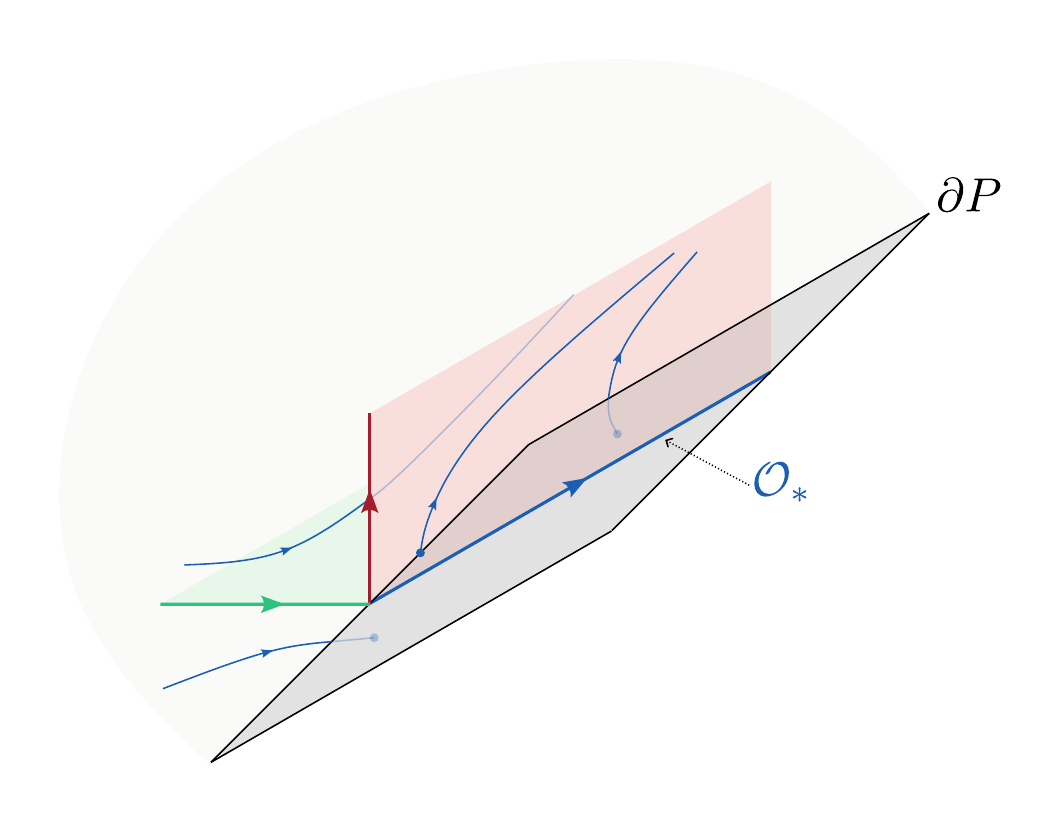}
    \vspace*{-1em}
    \caption{Boundary of a building block in a \nbh{} of a periodic orbit contained in the boundary}
    \label{fig: boundary building block local}
\end{figure}

We have the splitting $\pP \ssm \cO_* = \Pin \sqcup \Pout$, and it follows from Definition~\ref{def: qt surface}, Item~\ref{def: qt surface; it: in and out} that the \ccs{} of $\Pin$ and $\Pout$ alternate along $\cO_*$.
We refer to Figures~\ref{fig: building block} and~\ref{fig: boundary building block local}.

\subsection{Boundary lamination}
\label{sec: preli; subsec: boundary lamination}

Let $(P,X)$ be a \bb{}, and $\Lambda := \bigcap_{t\in \R} X^t (P)$ the maximal invariant set.
Let $\cW^s$ and $\cW^u$ be respectively the stable and the unstable manifolds of $\Lambda$.
The hyperbolicity of index $(1,1)$ of $\Lambda$ implies that $\cW^s$ and $\cW^u$ are 2-dimensional laminations, transverse to each other.
We suppose that the reader is familiar with basic hyperbolic dynamic theory and refer to ~\cite{hirschInvariantManifolds1977} for more details.
Let $\cO$ be a periodic orbit of $X$. 
The \emph{multipliers} of $\cO$ are the eigenvalues of the first return map on a Poincar\'e section for $\cO$.
The orbit $\cO$ has positive multipliers if and only if its (local) invariant manifolds are cylinders. 
In that case, the union $\cW^s_\loc(\cO) \cup \cW^u_\loc(\cO)$ of the local invariant manifolds of $\cO$ disconnects a \nbh{} of $\cO$ into four open regions, which are called \emph{quadrants}.
We will say that $(\tilde P, \tilde X)$ is a \emph{continuation} of $(P,X)$ if $\tilde P$ is a 3-manifold, $\tilde X$ is a $\cC^1$ vector field on $\tilde P$, and there exists an embedding $h \colon P \to \tilde P$ which maps the vector field $X$ to the vector field $\tilde X$ restricted to $h(P)$, and such that $h(\Lambda) \subset \intr \tilde P$.
If $\cO \in \cO_*$ is a periodic orbit contained in the boundary $\pP$ of a \bb{} $(P,X)$, the previous definitions hold in a continuation of $(P,X)$.

\begin{claim}\label{claim: boundary quadrant and multipliers}
Let $(P,X)$ be a \bb{} and $\cO$ a periodic orbit of $X$ contained in $\pP$.
The multipliers of $\cO$ are positive and $\partial P$ crosses two opposite quadrants of $\cO$.
\end{claim}

\begin{figure}[htb]
    \centering
    \vspace*{-1em}
    \captionsetup{width=.83\linewidth}
    \includegraphics[height=0.28\textheight]{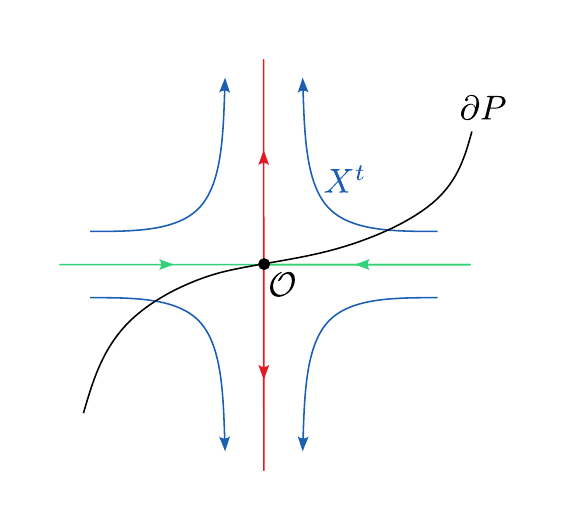}
    \vspace*{-1.5em}
    \caption{The boundary $\partial P$ crosses two opposite quadrants of~$\cO$}
    \label{fig: boundary crosses quadrant}
\end{figure}

\begin{proof}[Proof of Claim~\ref{claim: boundary quadrant and multipliers}]
Let $A$ be a tubular neighborhood of $\cO$ in $\partial P$.
It follows from Item~\ref{def: qt surface; it: in and out} of Definition~\ref{def: qt surface} that $A$ is an annulus.
Let $p$ be a point in $\cO$ and consider $\cC^1$-coordinates $(x, y, \theta)$ in the neighborhood of $p$ (in a continuation $(\tilde P, \tilde X)$ of the block) in which the orbit $\cO$ coincides with $x=y=0$ and the stable and unstable manifolds of $\cO$ are respectively straightened on $y=0$ and $x=0$.
The tangent plane to $A$ at $p$ is transverse to at least one of the two planes $x=0$ and $y=0$.
Let us assume for concreteness that it is transverse to $x=0$.
Since $A$ contains the periodic orbit $\cO: x=y=0$, we deduce that in the neighborhood of $p$, $A$ is a smooth graph $y= \phi(x, \theta)$, and the two connected components of $A \ssm \cO$ correspond to $x<0$ and $x>0$.
Orbits in a neighborhood of $\cW^s(\cO) \cup \cW^u(\cO)$ enter a neighborhood of $p$ near the unstable manifold $y=0$ and exit in a neighborhood of $p$ near the stable manifold $x=0$.
Thus, where $\phi>0$, these orbits go from the hypograph of $\phi$ to the epigraph of $\phi$, and where $\phi<0$, they go from the epigraph to the hypograph of $\phi$.
Now we know that $A$ is transverse to the orbits of $\tilde X$ and that the orbits cross $A$ in two different directions on one side and on the other side of $\cO$ by definition of a \qt{} surface.
Therefore, up to changing the coordinate $y$ by $-y$, we have $\phi(x, \theta)<0$ for $x<0$ and $\phi(x, \theta)>0$ for $x>0$.
In other words, $A$ is topologically transverse to the stable manifold of $\cO$ in a \nbh{} of $p$, and this is true for any point $p$.
We deduce that $\cW^s_\loc(\cO) \ssm \cO$ is composed of two connected components located on both sides of the annulus $A$.
Therefore $\cW^s_\loc(\cO)$ is an annulus.
The argument is the same for $\cW^u_\loc(\cO)$.
We deduce that the multipliers of $\cO$ are positive and $A$ crosses two opposite quadrants (Figure~\ref{fig: boundary crosses quadrant}).
\end{proof}

\begin{propdefi}[Entrance, Exit, Boundary lamination] \label{prop: existence boundary lamination}
\mb{}
\begin{enumerate}[i)]
    \item $\cL := (\cW^u \cup \cW^s) \cap \pP$ is a 1-dimensional lamination on $\pP$, called \emph{boundary lamination of $(P,X)$};

    \item $\cL^\iin := \cL \cap \Pin = \cW^s \cap \Pin$ is a 1-dimensional lamination on \Pin{}, called \emph{entrance lamination of $(P,X)$};
    
    \item $\cL^\out := \cL \cap \Pout = \cW^u \cap \Pout$ is a 1-dimensional lamination on \Pout{}, called \emph{exit lamination of $(P,X)$}.

\end{enumerate}
\end{propdefi}

\begin{rmk} \label{rmk: entrance and exit lamination are quasi-disjoint}
The lamination $\cW^s$ does not intersect the exit boundary $\Pout{}$ because any orbit in $\cW^s$ converges to $\Lambda$ in the future so never exits the block $(P,X)$, and similarly the lamination $\cW^u$ does not intersect the entrance boundary $\Pin{}$ because any orbit in $\cW^u$ converges to $\Lambda$ in the past.
The trace $\cW^s \cap \pP$ and the trace $\cW^u \cap \pP$ are therefore disjoint on $\pP$ except along the periodic orbits $\cO_*$ of $X$ contained in $\pP$ where they coincide.
The previous lemma states that the traces $\cW^s \cap \pP$ and $\cW^u \cap \pP$ come together along $\cO_*$ in a (single) lamination $\cL = (\cW^s \cup \cW^u) \cap \pP$ of dimension 1 on $\pP$.
\end{rmk}

\begin{proof}[Proof of Proposition~\ref{prop: existence boundary lamination}]
The lamination $\cW^s$ is tangent to the \vf{} $X$, and the \vf{} $X$ is transverse to $\Pin$, so $\cW^s$ intersects transversely $\Pin$.
Therefore $\cL^\iin :=\cW^s \cap \Pin$ is a one-dimensional lamination on $\Pin$.
Similarly for the lamination $\cL^\out = \cW^u \cap \Pout$.

According to Remark~\ref{rmk: entrance and exit lamination are quasi-disjoint},
$(\cW^s \cup \cW^u) \cap \pP = \cL^\iin \cup \cO_* \cup \cL^\out$.
It suffices to show that the union of $\cL^\iin$ and $\cL^\out$ reconnects along $\cO_* = \partial \Pin = \partial \Pout$ in a lamination $\cL$ on $\pP$.

Let $(\tilde P, \tilde X)$ be an extension of $(P, X)$.
Let $p$ be a point of $\cO \in \cO_*$, and $\cU_p$ a neighborhood of $p$ in $\tilde P$ with coordinates $(x, y, \theta) \in \R^3$ of class $\cC^1$, such that $\cO = \{x=y= 0\}$.
Let $(\partial_x, \partial_y, \partial_\theta)$ be the associated vector fields and consider the usual scalar product in these coordinates.
Since $\partial P$ is a smooth submanifold, it is locally transversely orientable by a normal vector field which we denote $\vec n$ in $\cU_p$. 
Since $\partial P$ contains $\cO$, we deduce that $\vec n (r)$ converges to a vector in the normal plane to $\partial_\theta$ when $r$ tends to $p$.
Similarly, the lamination $\cW^s$ is locally transversely orientable by a normal vector field which we denote $\vec m$ in $\cU_p$.
Since the leaf passing through $p$ contains $\cO$, we deduce that the vector $\vec m (r)$ converges to a vector in the plane normal to the vector $\partial_\theta$ when $r$ tends to~$p$.

A vector field tangent to the lamination $\cL^\iin$ is defined by the vector product $\vec v = \vec n \times \vec m$ at each intersection point $\cW^s \cap \Pin$.
It is a vector field normal to the plane generated by $\vec n$ and $\vec m$.
We deduce that $\vec v (r)$ converges to a vector tangent to $\R . \partial_\theta$ for $r \in \cL^\iin$ which tends to $p$, in other words $\vec v (r)$ converges to a vector tangent at $T_p \cO$. Doing this at any point $p$ of the set $\cO_*$, we show that the lamination $\cL^\iin$ on $\Pin$ reconnects with $\cO_* = \partial P^\iin$ in a lamination on \linebreak[4]$\overline \Pin = \Pin \cup \cO_*$.

A similar proof shows that the lamination $\cL^\out$ on $\Pout$ reconnects with $\cO_* = \partial P^\out$ in a lamination $\overline \Pout = \Pout \cup \cO_*$.
Therefore, $$\cL :=  (\cW^s \cup \cW^u) \cap \pP = \cL^\iin \cup \cO_* \cup \cL^\out$$ is a lamination of dimension~1 in $\partial P = \Pin \cup \cO_* \cup \Pout$.
\end{proof}

\vspace*{-1.5em}

\subsubsection*{Quasi-Morse--Smale laminations}

\begin{defi}[Quasi-Morse--Smale lamination] \label{def: qms lam}
Let $\cL$ be a 1-dimensional lamination on a closed surface $S$.
We say that $\cL$ is a \emph{quasi-Morse--Smale} lamination if it satisfies the following conditions:
\begin{enumerate}
    \item \label{def: qms lam; it: finite compact}
    There exists a finite number of compact leaves $\Gamma = \{\gamma_1, \dots, \gamma_N\}$;
    \item \label{def: qms lam; it: half non-compact}
    Each half non-compact leaf accumulates on a single compact leaf;
    \item \label{def: qms lam; it: holonomy}
    Each compact leaf $\gamma$, given an orientation, has a contracting or expanding holonomy on each side.\footnote{If $\gamma$ admits two sides in $S$, the holonomy of $\gamma$ may be contracting on both sides, expanding on both sides, or contracting on one side and expanding on the other. If it has a single side in $S$, its holonomy is either contracting or expanding.}
\end{enumerate}
The (possibly empty) set of elements of $\Gamma$ such that the holonomy is contracting on one side and expanding on the other is denoted $\Gamma_*$, and called the set of \emph{marked compact leaves}.
We further require:
\begin{enumerate}[resume]
    \item \label{def: qms lam; it: in out}
    There exists a splitting $S \ssm \Gamma_* = S^\iin \, \sqcup \, S^\out$ into two disjoint open sets such that each leaf $\gamma_* \in \Gamma_*$ is adjacent to a connected component of $S^\iin$ and to a connected component of $S^\out$.
    We call it the \emph{(in,out)-splitting of $S$ for $\cL$.}
\end{enumerate}
\end{defi}

\begin{figure}[htb]
    \centering
    \vspace*{-1.3em}
    \includegraphics[height=0.27\textheight]{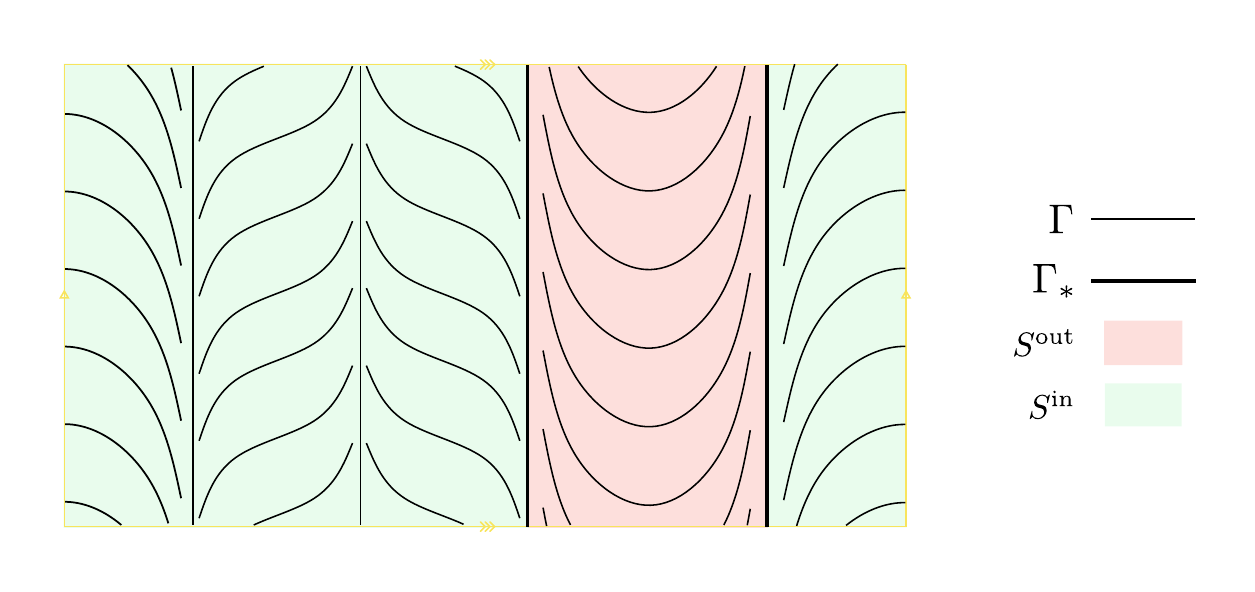}
    \vspace*{-1em}
    \caption{An example of \qms{} lamination on the torus}
    \label{fig: qms lamination}
\end{figure}

\begin{rmk} \label{rmk: QMS lam and MS lam}
A \qms{} lamination $\cL$ such that the set of marked leaves $\Gamma_*$ is empty is a \emph{Morse--Smale lamination} in the sense of \cite[Definition~3.9]{beguinBuildingAnosovFlows2017}, which justifies this terminology.
\end{rmk}

\begin{prop} \label{prop: block boundary lam are qms}
Let $(P,X)$ be a \bb{}.
Then the boundary lamination $\cL$ is a \qms{} lamination on $\pP$, and the set $\Gamma_*$ of marked leaves of $\cL$ coincides with the set $\cO_*$ of periodic orbits of $X$ contained in $\pP$.
\end{prop}

In the case where the boundary $\partial P$ is transverse to the vector field $X$, i.e. if $(P,X)$ is a \bby{} block, this property is known.
More precisely, we have the following lemma.

\begin{lem}[{\cite[Proposition~3.8]{beguinBuildingAnosovFlows2017}}]
\label{lem: boundary lamination of BBY block is MS}
Let $(U,Y)$ be a \bb{} whose boundary $\partial U$ is transverse to the vector field $Y$.
Let $\cL_Y$ be the boundary lamination of $(U,Y)$.
Then it satisfies the following properties:
\begin{enumerate}
    \item there exists a finite number of compact leaves;
    \item each non-compact half-leaf accumulates on a single compact leaf;
    \item each compact leaf has an orientation for which its holonomy is contracting (on each side).
\end{enumerate}
\end{lem}

In other words, the boundary lamination $\cL_Y$ of a \bby{} block $(U,Y)$ is a \qms{} lamination for which the set of marked leaves is empty, namely a \emph{Morse--Smale lamination}.
The following lemma reduces the study of the boundary lamination of a building block $(P,X)$ to the case where the boundary of the block is transverse to the vector field.
We denote by $\overline A$ the closure of a set $A$.
Let $(P,X)$ be a \bb{}.
Let $\cO_*$ be the set of periodic orbits of $X$ contained in $\pP$, and $\Lambda_X$ the maximal invariant set of $(P,X)$, $\Pin$ the entrance boundary and $\Pout$ the exit boundary.

\begin{lem}[From building block to \bby{} block] \label{lem: from bb to bby}
There exists a \bby{} block $(U,Y)$ such that, if $\Lambda_Y$ denotes the maximal invariant set of $(U,Y)$, $(\cW^s_Y, \cW^u_Y)$ the pair of stable and unstable lamination of $\Lambda_Y$, $U^\iin$ the entrance boundary, $U^\out$ the exit boundary, and $\cL_Y$ the boundary lamination, then
\begin{enumerate}
    \item \emph{(embedding)} \label{lem: block to BBY; it: embedding}
    There exists an embedding $h \colon P \to U$ which maps the vector field $X$ on the vector field $Y$ in $h(P)$ and the set $\Lambda_X$ on $\Lambda_Y$.
    \item \emph{(lamination)} \label{lem: block to BBY; it: lamination}
    There exists a finite collection $D_* = \{D_1, \dots, D_n\}$ of open disks contained in $\partial U$, disjoint from $\cL_Y$, bounded by pairwise distinct compact leaves of $\cL_Y$, such that $ U^\iin \ssm D_*$ is isotopic to $\overline \Pin$ along the lamination $\cW^s_Y$ and $U^\out \ssm D_*$ is isotopic to $ \overline \Pout$ along the lamination $\cW^u_Y$.
\end{enumerate}
\end{lem}

Recall that a \emph{linearizing neighborhood} of a periodic orbit $\cO$ of $X$ is a tubular neighborhood $\cV$ of $\cO$ on which the flow of $X$ is (locally) orbit equivalent to the suspension flow of the diagonal linear map $(x, y) \mapsto (\lambda x, \mu y)$ on $\R^2$ with $0 < \vert \lambda \vert < 1$ and $\vert \mu \vert>1$.
Every (saddle) hyperbolic periodic orbit has a linearizing neighborhood.
A $\cC^1$-coordinate system $\xi = (x,y,\theta) \in \R^2 \times \R/\Z$ on $\cV$ in which the flow is linear is a \emph{\lcs{} of $\cO$}. 

\begin{proof}[Proof of Lemma~\ref{lem: from bb to bby}]
Let $(\tilde P, \tilde X)$ be a continuation of $(P,X)$.
Let $\cO_* = \{ \cO_1, \dots,  \cO_n \}$ be the collection of periodic orbits of $X$ contained in $\pP$.
For each $\cO_i \in \cO_*$, we can consider (up to local orbit equivalence) a linearizing neighborhood $\cV_i$ of $\cO_i$ for the flow of $\tilde X$ in $\tilde P$, provided with a \lcs{} $(x,y,\theta) \in \R^2 \times \R/\Z$.
The boundary $\partial P$ crosses the opposite quadrants $\{x > 0, y > 0\}$ and $\{x < 0, y < 0\}$ (Claim~\ref{claim: boundary quadrant and multipliers}).
Let $S$ be a topological surface, smooth outside a finite number of simple closed curves $c_1^\iin, \dots, c_n^\iin$ and $c_1^\out, \dots, c_n^\out$, which coincides with $\pP$ outside the $\cV_i$ neighborhoods, and decomposes into the union $S = S^\iin \cup (\bigcup_i A_i) \cup S^\out$, where
\begin{itemize}
    \item $S^\iin$ is a surface with boundary transverse to the \vf{} $\tilde X$ and bounded by $c_1^\iin, \dots, c_n^\iin$, which coincides with $\Pin$ outside the union of the neighborhoods $\cV_i$;
    \item $S^\out$ is a surface with boundary transverse to the \vf{} $\tilde X$ and bound\-ed by $c_1^\out, \dots, c_n^\out$, which coincides with $\Pout$ outside the union of the neighborhoods $\cV_i$;
    \item $A_i$ is an annulus tangent to the \vf{} $\tilde X$, contained in the linearizing neighborhood $\cV_i$ of $\cO_i$ in the quadrant $\{ x > 0, y < 0\}$, and bounded by the curves $c^\iin_i$ and $c^\out_i$, such that all orbits of $\tilde X$ on $A_i$ go from $c^\iin_i$ to $c^\out_i$.
\end{itemize}
See Figure~\ref{fig: intermediate block}.
\begin{figure}[htb]
    \centering
    \vspace*{-2em}
    \includegraphics[height=0.27\textheight]{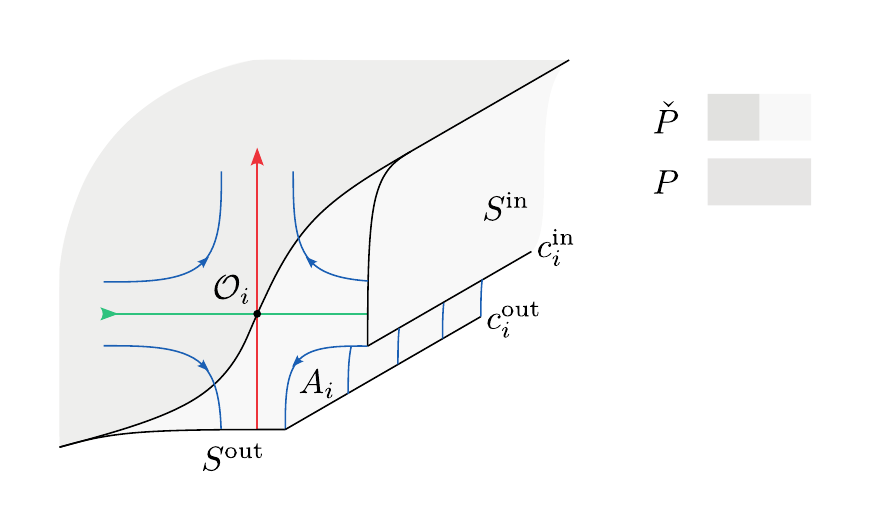}
    \vspace*{-1em}
    \caption{Manifold $\check P$ in a linearizing \nbh{} $\cV_i$ of $\cO_i$}
    \label{fig: intermediate block}
\end{figure}

The surface $S$ cuts in $\tilde P$ a compact submanifold $\check P$ with boundary, which contains $P$ and coincides with $P$ outside the union of the neighborhoods $\cV_i$.
Let $\check X$ be the restriction of $\tilde X$ to $\check P$.
The \vf{} $\check X$ restricted to $A_i$ is orbit equivalent to the vertical field $\partial t$ on $\S^1 \times I$, where $t$ is the coordinate on $I=[0,1]$.
We glue on each annuli $A_i \simeq \S^1 \times I$ a cylinder $C_i = \D^2 \times I$ provided with the vertical vector field $\partial_t$, along $B_i = \partial \D_2 \times I \subset \partial C_i$.
Let $f_i \colon B_i = \partial \D^2 \times I \to \S^1 \times I \simeq A_i$ a \diff{} which preserve the coordinate $t \in I$ (Figure~\ref{fig: glue cylinder for bby block}).

\begin{figure}[htb]
    \centering
    \vspace*{-1.5em}
    \includegraphics[height=0.32\textheight]{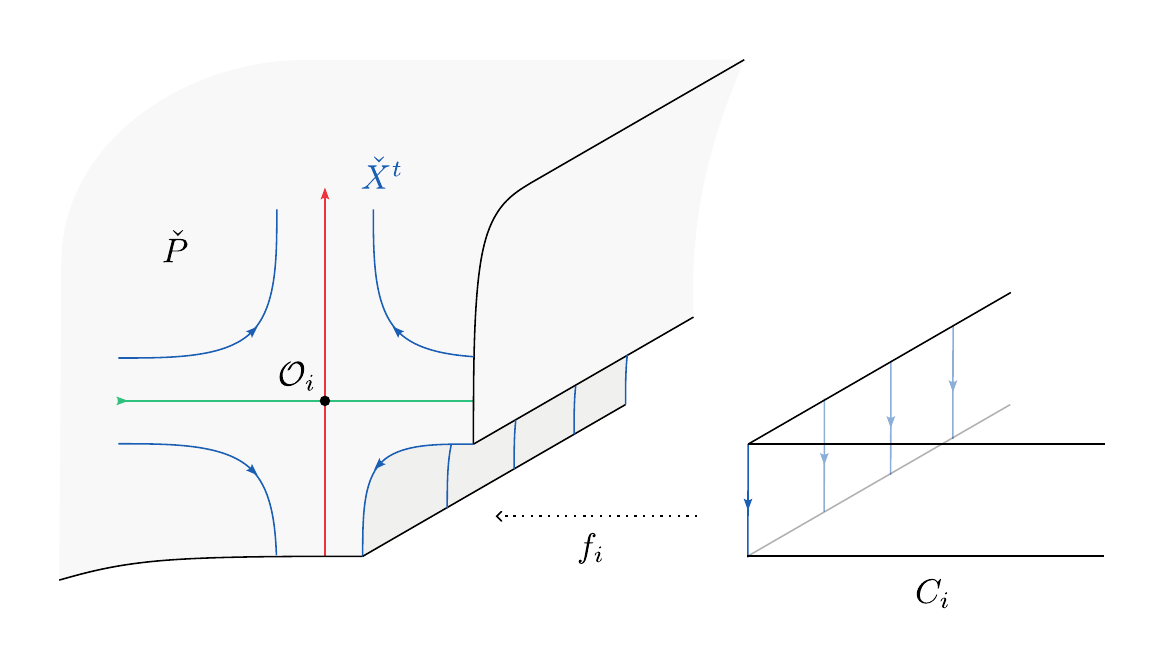}
    \vspace*{-1em}
    \caption{Gluing of a cylinder $C_i$ along the annulus $A_i$ tangent to the \vf{} $\check X$}
    \label{fig: glue cylinder for bby block}
\end{figure}

Define the quotient manifold $ U := \check P \cup (\bigcup_i C_i)/ f $ where $f$ is the product of $f_i$.
Each $f_i$ maps the vertical vector field of $C_i$ on the vector field $\check X$, it follows that $U$ is equipped with a $\cC^1$ \vf{} $Y$ for a differentiable structure on $U$, induced by the union of $\check X$ and the vertical \vf{} on each $C_i$.

Let us show that the pair $(U,Y)$ satisfies Lemma~\ref{lem: from bb to bby}.
The manifold $U$ is a compact smooth manifold of dimension~3 with boundary.
The boundary $\partial U$ is the quotient $\partial U = S \cup (\bigcup_i \partial C_i)/f $.
A point of $\partial U$ is either in $S \ssm \bigcup_i A_i$, where $Y$ coincides with the \vf{} $\check X$, or in the boundary $\D^2 \times \{ 0 \}$ or $\D^2 \times \{ 1 \}$ of a cylinder $C_i = \D^2 \times [0,1]$, where $Y$ is the vertical vector field.
It follows that $Y$ is transverse to $\partial U$.
It is clear by construction that the manifold $P$ is embedded in the manifold $U$ via an embedding which maps $X$ to $Y$.
Moreover, the maximal invariant set of $Y$ in $P$ coincides with the maximal invariant set of $X$ in $P$ via this embedding.
Indeed, an orbit of $Y$ in $C_i$ enters and exits $C_i$, hence $U$, in uniformly bounded time, and an orbit of $Y$ which intersects $\check P$ is either an orbit of the flow of $X$ contained in $P$, or intersects the boundary of $\check P$.
It follows that $(U,Y)$ is a \bb{} which satisfies Item~\ref{lem: block to BBY; it: embedding} of Lemma~\ref{lem: from bb to bby}.

\smallskip

It remains to check Item~\ref{lem: block to BBY; it: lamination}.
The boundary lamination $\cL_Y$ does not intersect the boundary of $U$ contained in $\partial C_i$ because all orbits in $C_i$ enter and exit $U$.
In $\check P$, the orbits of $Y$ coincide with the orbits of $\check X$.
We denote $\check \cW^s$ and $\check \cW^u$ the stable and unstable manifolds of the set $\Lambda_X$ for the flow of $\check X$ in $\check P$.
Then $\check \cW^s$ is the union of the saturation of the laminations $\cW^s_X$ by the flow of $\check X$ and the local stable manifolds of the periodic orbits $\cO_i$ in the linearizing neighborhoods $\cV_i$.
Similarly $\check \cW^u$ is the union of the saturation of the laminations $\cW^u_X$ by the flow of $\check X$ and of the local unstable manifolds of the periodic orbits $\cO_i$ in the linearizing neighborhoods $\cV_i$.
Define
$ \check \cL^\iin := \cW^s \cap \partial \check P = \cW^s \cap S$ and  $\cL^\out := \cW^u \cap \partial \check P = \cW^s \cap S $.
By definition, the orbits of the points of $\check \cL^\iin$ are the orbits of the flow of $\check X = \res{Y}{\check P}$ that accumulate on the set $\Lambda_X = \Lambda_Y$ in the future, in other words $\check \cL^\iin = \cL^\iin_Y$ is the entrance lamination of $(U,Y)$.
The local stable manifold of the orbit $\cO_i$ in the neighborhood $\cV_i$ intersects $\partial \check P$ along a compact leaf $\gamma^s_i$, which bounds a connected component $D_i^\iin$ of $U^\iin \ssm \cL^\iin_Y$, where $U^\iin$ is the union of the \ccs{} of $\partial U$ along which the \vf{} $Y$ points into $U$.
This component is a disk which contains the entrance boundary $\D^2 \times \{ 0 \}$ of the cylinder $C_i$ (Figure~\ref{fig: lamination on bby block}).

\begin{figure}[htb]
    \centering
        \vspace*{-1em}
    \includegraphics[height=0.27\textheight]{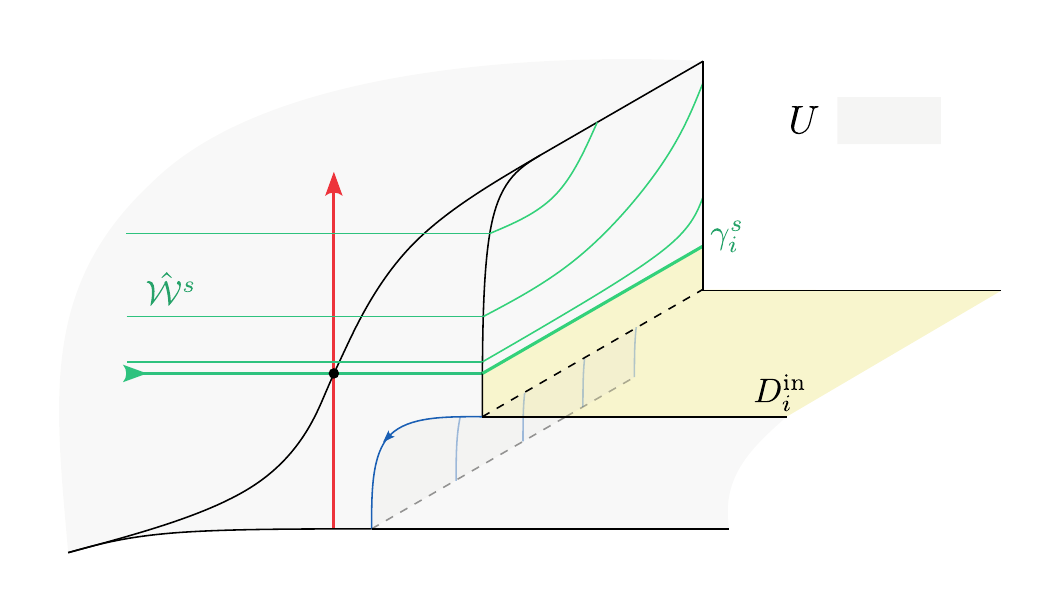}
        \vspace*{-1em}
    \caption{Disk $D^\iin_i$ in $U^\iin$, disjoint from $\cL^\iin_Y$ and bounded by a compact leaf $\gamma^s_i$}
    \label{fig: lamination on bby block}
\end{figure}

Similarly, $\check \cL^\out =\cL^\out_Y$ is the exit lamination of $(U,Y)$, and the local unstable manifold of the orbit $\cO_i$ in the neighborhood $\cV_i$ intersects $\partial \check P$ along a compact leaf $\gamma^u_i$.
This leaf borders a connected component $D_i^\out$ of $U^\out \ssm \cL^\out_Y$, where $U^\out$ is the union of the \ccs{} of $\partial U$ along which the vector field $Y$ points inwards.
This component is a disk which contains the exit boundary $\D^2 \times \{ 1\}$ of the cylinder $C_i$.
Let $D_* = (\bigcup_i D_i^\iin) \cup (\bigcup_i D_i^\out)$.
The surfaces $U^\iin \ssm D$ and $\overline \Pin$ coincide outside the union of neighborhoods $\cV_i$, and in each $\cV_i$, they are annuli whose interior is transverse to the lamination $\check \cW^s$. {Thus there is an isotopy between $U^\iin \ssm D_*$ and $\Pin$, supported in the linearizing  neighborhoods $\cV_i$, which preserves each leaf of the lamination $\check \cW^s$.
The argument is symmetric for the isotopy between $U^\out \ssm D_*$ and $\overline \Pout$.}
\end{proof}

\begin{figure}[htb]
    \centering
        \vspace*{-1.5em}
    \includegraphics[height=0.265\textheight]{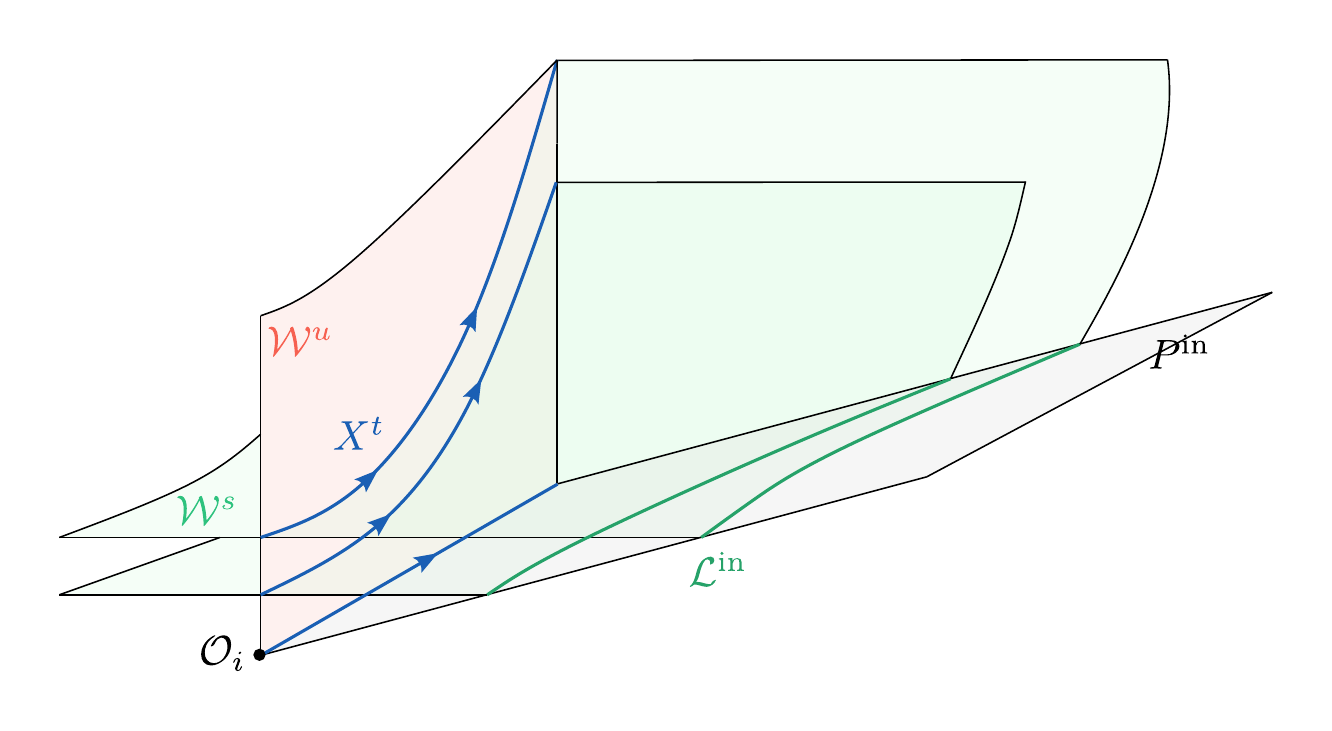}
        \vspace*{-1em} 
    \caption{Holonomy of an orbit $\cO_i \in \cO_*$ along the boundary lamination on the entrance boundary $\Pin$}
    \label{fig: holonomy periodic orbit boundary}
\end{figure}

\vspace*{-0.5em}
\begin{proof}[Proof of Proposition~\ref{prop: block boundary lam are qms}]
Let $(P,X)$ be a \bb{}  and $(U,Y)$ be the \bb{} associated by Lemma~\ref{lem: from bb to bby} whose boundary $\partial U$ is transverse to the vector field $Y$.
Let $\cL_X = \cL^\iin_X \cup \cO_* \cup \cL^\out_X$ denote the boundary lamination of $(P,X)$ and $\cL_Y = \cL^\iin_Y \cup \cL^\out_Y$ the boundary lamination of $(U,Y)$.
Let $D_* = \{D_1, \dots, D_n \}$ be the collection of (open) disks in $\partial U$, disjoint from $\cL_Y$ and bordered by compact leaves of $\cL_Y$ given by Lemma~\ref{lem: from bb to bby}, Item~\ref{lem: block to BBY; it: lamination}.
It follows from this item the existence of a \diff{} $H^s \colon U^\iin \ssm D_* \to \overline{\Pin}$ which maps the restriction of the lamination $\cL^\iin_Y$ to the lamination $\cL_X^s \cup \cO_*$.
We conclude that the lamination $\cL^\iin_X \cup \cO_*$ satisfies the first three items of Definition~\ref{def: qms lam} of a \qms{} lamination on $\overline{\Pin}$.
The compact leaves bordering the disks $D_i \in D_* \cap U^\iin$ correspond via this diffeomorphism to the compact leaves of $\cL^\iin_X \cup \cO_*$ on the boundary of $\overline{\Pin}$, in other words to the periodic orbits $\cO_i \in \cO_*$ of $X$ on $\pP$.
The similar argument for $\cL^\out_X \cup \cO_*$ implies that the lamination $\cL_X = \cL^\iin_X \cup \cO_* \cL^\out_X$ on $\pP = \overline{\Pin} \cup \overline{\Pout}$ satisfies the first three items of Definition~\ref{def: qms lam}, and each compact leaf of $\cL^\iin_X$ and $\cL^\out_X$ admits an orientation for which its holonomy is contracting (on both sides).

It remains to show that the holonomy of the oriented elements of $\cO_*$ is contracting on one side and expanding on the other, i.e. $\cO_*$ is the set of marked leaves of the lamination $\cL_X$.
Item~\ref{def: qms lam; it: in out} of Definition~\ref{def: qms lam} is then satisfied by definition of a \qt{} surface (Definition~\ref{def: qt surface}, Item~\ref{def: qt surface; it: in and out}) for the decomposition $\pP \ssm \cO_* = \Pin \cup \Pout$.
Let $\cO_i \in \cO_*$ be a periodic orbit of $X$ on $\pP$.
It is a compact leaf of $\cL_X$, adjacent on one side to $\Pin$ and on the other to $\Pout$.
From the previous proof, there exists an orientation of $\cO_i$ such that its holonomy on the side of $\Pin$ is contracting, and there exists an orientation such that its holonomy on the side of $\Pout$ is contracting.
Let us show that these two orientations cannot coincide.
Let $l^s$ be a leaf of $\cL^\iin_X$ which accumulates on $\cO_i$.
Then the leaf $l^s$ intersects a local product \nbh{} of a point $p \in \cO_i$,
and it follows that $l^s$ is the transverse intersection of the stable manifold $\cW^s_X(\gamma)$ of an orbit $\gamma \in \cW^u_X(\cO_i)$ with the surface $\Pin$.
As the surfaces $\Pin$ and $\cW^u_X(\cO_i)$ are both transverse to $\cW^s_X$, we deduce that the holonomy of $\cL^\iin_X \cup \cO_i$ along $\cO_i$ is conjugate to the holonomy of the lamination induced by $\cW^s_X$ on $\cW^u_X(\cO_i)$ along $\cO_i$.
This holonomy is conjugate to the first return map of the flow of $X$ on a transversal of $\cO_i$ in $\cW^u_X(\cO_i)$, and is therefore expanding for the orientation of $\cO_i$ given by the flow.

Similarly, the holonomy of $\cL^\out_X \cup \cO_i$ along $\cO_i$ is conjugate to the holonomy of the lamination induced by $\cW^u_X$ on $\cW^s_X(\cO_i)$ along $\cO_i$, and this holonomy is contracting for the orientation of $\cO_i$ induced by the orientation of the flow of $X$ (Figure~\ref{fig: holonomy periodic orbit boundary}).
Finally, each compact leaf $\cO_i \in \cO_*$ of the lamination $\cL_X$ has an expanding holonomy on one side and a contracting one on the other.
We conclude that the lamination $\cL_X = \cL^\iin_X \cup \cO_* \cL^\out_X$ is a \qms{} lamination, which proves Proposition~\ref{prop: block boundary lam are qms}.
\end{proof}

\vspace*{-0.8em}

\subsubsection*{Prefoliation and filling lamination}

\begin{defi}[Prefoliation] \label{def: prefoliation}
A 1-dimensional lamination $\cL$ on a surface $S$ is said to be a \emph{prefoliation} if there exists a 1-dimensional foliation $\cF$ on $S$ which contains $\cL$ as a sublamination, i.e. such that the leaves of $\cL$ are leaves of $\cF$.
\end{defi}

\begin{rmk} \label{rmk: pre lamination necessary condition}
This is a necessary condition on the boundary lamination $\cL$ of a block $(P,X)$ for the existence of an Anosov vector field $Z$ on a closed manifold $\cM$ obtained by gluing the boundary of a block $(P,X)$.
Indeed, the surface $\pP$ then projects in $\cM$ onto a surface $S$ \qt{} to the \vf{} $Z$ and provided with two foliation $(\zeta^u, \zeta^s)$ of dimension~1, whose union contains the boundary lamination $\cL$ as sublamination.
These are the traces on $S$ of the unstable and stable foliations $\cF^u$ and $\cF^s$ of the Anosov flow generated by $Z$.
It follows that each \cc{} of $\pP$ is a closed connected surface of zero Euler characteristic, hence a torus or a Klein bottle.
\end{rmk}

\begin{prop} \label{prop: complementary of a prefoliation}
A \qms{} lamination $\cL$ on a closed surface $S$ is a prefoliation on $S$ if and only if each connected component of $S \ssm \cL$ is either
\begin{enumerate}[i)]
    \item homeomorphic to an annulus or a Moebius strip bounded by compact leaves, or
    \label{prop: complement prefol; it: annulus and mobius}
    \item homeomorphic to $\R^2$, and the accessible boundary\footnote{Those are the points of the boundary which are endpoints of a segment contained in the interior of the component.} is the union of two non-compact leaves $l$ and $l'$ which are asymptotic to each other at both ends.
    \label{prop: complement prefol; it: strip}
\end{enumerate}
\end{prop}

\begin{defi}[Strip, Filling lamination] \label{def: strip and filling lam}
Let $B$ be a \cc{} of $S \ssm \cL$.
We say that $B$ is a \emph{strip} if $B$ is homeomorphic to $\R^2$ and the accessible boundary of $B$ is the union of two non-compact leaves $l$ and $l'$ of $\cL$ which are asymptotic to each other at both ends.
If all \cc{} of $S \ssm \cL$ are strips, we say that the lamination $\cL$ is \emph{filling on $S$}.
\end{defi}

\begin{proof}[Proof of Proposition~\ref{prop: complementary of a prefoliation}]
A lamination contains a finite number of compact leaves (Definition~\ref{def: qms lam}, Item~\ref{def: qms lam; it: finite compact}).
The proposition is then a direct consequence of the following fact, shown in \cite{beguinBuildingAnosovFlows2017}.
\end{proof}
\begin{lem}[{\cite[Lemme~3.12]{beguinBuildingAnosovFlows2017}}]
Let $S$ be a closed surface and $\cL$ be a 1-dimen\-sional lamination on $S$ which contains a finite number of compact leaves.
Then there exists only a finite number of \ccs{} of $S \ssm \cL$ which are not strips.
\end{lem}

\begin{figure}[htb]
    \centering
    \vspace*{-1.3em}
    \includegraphics[height=0.28\textheight]{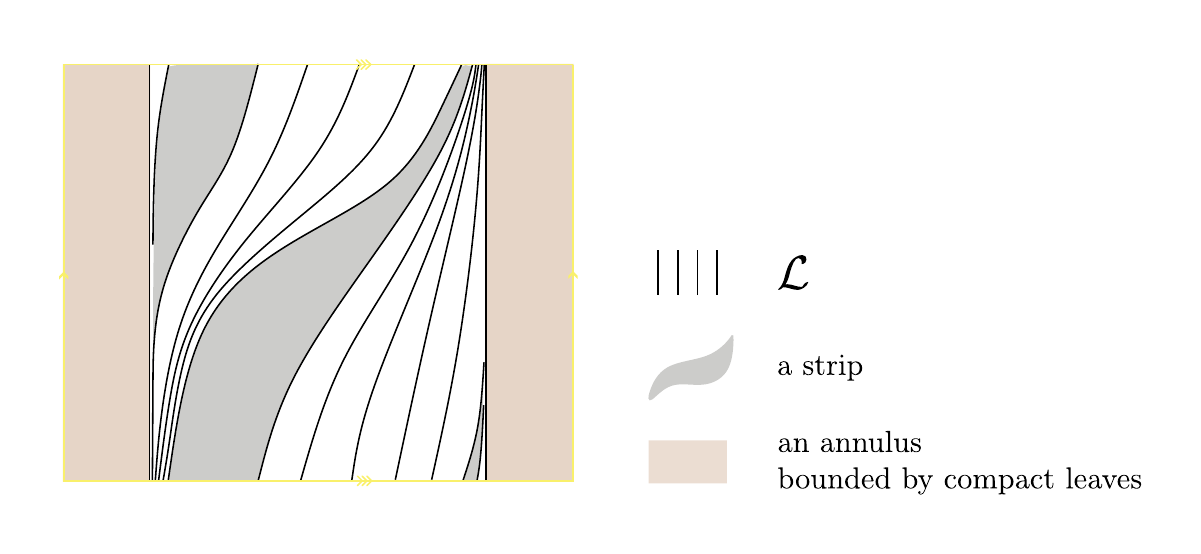}
    \vspace*{-1em}
    \caption{A prefoliation $\cL$ on the torus $S = \T^2$}
    \label{fig: prefoliation}
\end{figure}

Figure \ref{fig: prefoliation} shows a type~\ref{prop: complement prefol; it: annulus and mobius} (in brown) and a type~\ref{prop: complement prefol; it: strip} (in grey) for a prefoliation $\cL$ on $\T^2$. 

We will have to make the following additional assumption for technical reasons.

\begin{defi}[Filled block]
\label{def: filled block}
We say that a block $(P,X)$ is a \emph{filled block} if the boundary lamination $\cL$ is filling on $\pP$.
\end{defi}

\begin{figure}[htb]
    \centering
    \vspace*{-1em}
    \includegraphics[height=0.28\textheight]{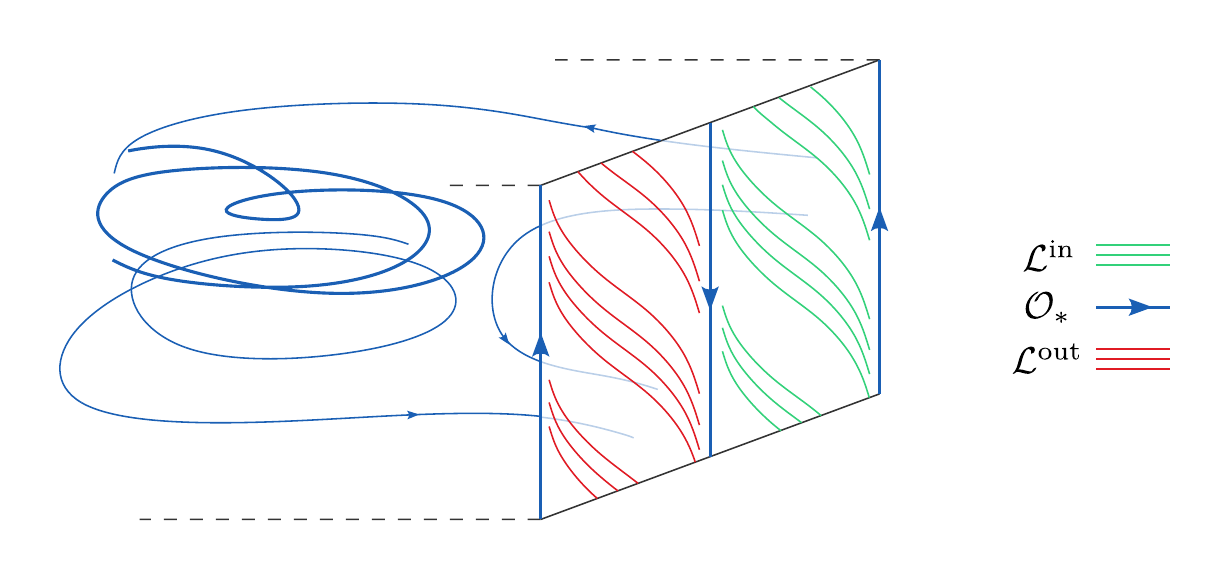}
    \vspace*{-1em}
    \caption{A \bb{} $(P,X)$ with a torus $\T^2$ as boundary and a filling boundary lamination}
    \label{fig: block with filling lam}
\end{figure}

In other words, $(P,X)$ is a filled block if the boundary lamination $\cL$ is a prefoliation and no \cc{} of $\pP \ssm \cL$ is an annulus or a Moebius strip.
Figure~\ref{fig: block with filling lam} represents an example of a \bb{} whose boundary lamination is filling.

\subsection{Gluing map}

We want to choose a \diff{} $\varphi \colon \pP \to \pP$ in such a way that the quotient $P/\varphi = P/(x \sim \varphi(x))$ is a closed smooth manifold, and the vector field $X$ induces a vector field of class $\cC^1$ on $P/ \varphi$ which is Anosov.
We then note $P_\varphi := P/\varphi$ and $X_\varphi$ the vector field induced by $X$ on $P_\varphi$.

\subsubsection*{Gluing maps}

\begin{defi}[Gluing map of a building block]
\label{def: gluing map}
Let $(P,X)$ be a \bb{}, $\cO_*$ the periodic orbits of $X$ contained in $\pP$, $\Pin$ the entrance boundary and $\Pout$ the exit boundary of $(P,X)$.
A \emph{gluing map of $(P,X)$} is a $\cC^1$-involution $\varphi \colon \pP \to \pP$ such that:
\begin{enumerate}
    \item \label{def: gluing map; it: pair cc}
    there exists a partition $\pP = \partial_1 P \,\sqcup \, \partial_2 P$, where $\partial_1 P$ and $\partial_2 P$ are unions of \ccs{} of $\pP$, and such that $\varphi (\partial_1 P) = \partial_2 P$;
    \item\label{def: gluing map; it: orbit on orbit}
    $\varphi(\cO_*) \!=\! \cO_*$ and $\varphi$ preserves the orientation of the flow on the orbits of $\cO_*$;
    \item \label{def: gluing map; it: out on in}
    $\varphi (\Pout{}) = \Pin{}$.
\end{enumerate}
\end{defi}

\begin{rmk}$ $
\begin{itemize}
    \item One could equivalently define a gluing map $\varphi$ of a block $(P,X)$ as a \diff{} from $\partial_1 P$ to $\partial_2 P$, for some partition $\pP = \partial_1P \sqcup \partial_2 P$ in two collections of \ccs{}, and which satisfies Items~\ref{def: gluing map; it: orbit on orbit} and~\ref{def: gluing map; it: out on in} of Definition~\ref{def: gluing map}.
    This definition has the drawback of arbitrarily breaking the symmetry between $\partial_1P$ and $\partial_2P$, so we prefer to see a gluing map of a block as an involution of $\pP$ which matches the \ccs{} pairwise.
    In the case where the boundary $\pP$ is transverse to the vector field $X$, in other words if $(P,X)$ is a \bby{} block, we can take $\partial_1 P = \Pout$, $\partial_2 P = \Pin$, and we already have a natural distinction between $\partial_1 P$ and $\partial_2 P$.
    It is then more natural to define a gluing map as a \diff{} $\varphi \colon \Pout \to \Pin$ as the authors of \cite{beguinBuildingAnosovFlows2017} do.
    
    \item If $\varphi$ is a gluing map of a block $(P,X)$, there exists a smooth structure on the quotient space $P / \varphi$ which makes it a closed smooth manifold denoted $P_\varphi$.
    If the boundary of $P$ is transverse to the vector field $X$, then $X$ always induces a vector field of class $\cC^1$ on $P_\varphi$ which is denoted $X_\varphi$.
    In the case where there are periodic orbits contained in the boundary, this is not always true.
    Indeed, it is necessary to be able to identify the dynamics in the neighborhood of two orbits which are paired, in particular the periodic orbits must have the same multipliers.
\end{itemize}
\end{rmk}

\begin{defi}[Dynamical gluing map]
\label{def: dynamical gluing map}
A \emph{dynamical gluing map} of $(P,X)$ is a gluing map which maps the vector field $X$ to itself.
Formally, we require that there exists a continuation $(\tilde P, \tilde X)$ of $(P,X)$ which contains $P$ in its interior, and a tubular neighborhood $U$ of $\pP$ in $\tilde P$, such that $\varphi$ is the restriction of a \diff{} $\tilde \varphi \colon U \to U$ which preserves the vector field $\tilde X$.
\end{defi}

\begin{claim}
\label{claim: dynamical gluing map induce vector field}
If $\varphi$ is a dynamical gluing map of $(P,X)$, then there exists a differentiable structure on the quotient manifold $P_\varphi = P / \varphi$ which makes it a closed smooth manifold equipped with a vector field $X_\varphi$ of class $\cC^1$ induced by the vector field $X$ on $P$.
\end{claim}

\begin{proof}
This is a consequence of the definition of a dynamical gluing map, and of the natural differentiable structure on $P_\varphi$.
\end{proof}

\subsubsection*{Strong (quasi)-transversality}
We are led to make additional assumptions about the gluing map concerning the way the boundary lamination $\cL$ and the image $\varphi_* \cL$ intersects each other.

\begin{defi}[Strongly quasi-transverse laminations] \label{def: sqt laminations}
Let $S$ be a surface (possibly with boundary) and $(\cL_1,\cL_2)$ a pair of laminations on $S$.
\begin{itemize}
\item $\cL_1$ and $\cL_2$ are said to be \emph{\st} if they are transverse to each other on $S$ and the closure of a connected component of the complementary $S \ssm (\cL_1 \cup \cL_2)$ is a compact rectangle
bounded by two disjoint arcs of leaves of $\cL_1$ and two disjoint arcs of leaves of $\cL_2$.
\end{itemize}
Suppose further that $\cL_1$ and $\cL_2$ are \qms{} laminations.
Denote $\Gamma_{\cL_i}$ the compact leaves of $\cL_i$, $\Gamma_{\cL_i, *}$ the marked compact leaves, and $S \ssm \Gamma_{\cL_i, *} = S^\iin_i \cup S^\out_i$ the decomposition of $S \ssm \Gamma_{\cL_i,*} $ corresponding to $\cL_i$, for $i=1,2$.
\begin{itemize}
\item They are said to be \emph{quasi-transverse} if they coincide on the set of marked leaves $\Gamma_{\cL_1, *} = \Gamma_{\cL_2, *} =: \Gamma_*$, are transverse to each other on the complementary $S \ssm \Gamma_*$, and if $S^\iin_1 = S^\out_2$.
\item They are said to be \emph{\sqt{}} if they are \qt{} and if they are \st{} on $S \ssm \Gamma_*$.
\end{itemize}
\end{defi}

\begin{figure}[htb]
    \centering
    \captionsetup{width=.83\linewidth}
    \vspace*{-1.5em}
    \hspace*{-1.5em}
    \includegraphics[height=0.23\textheight]{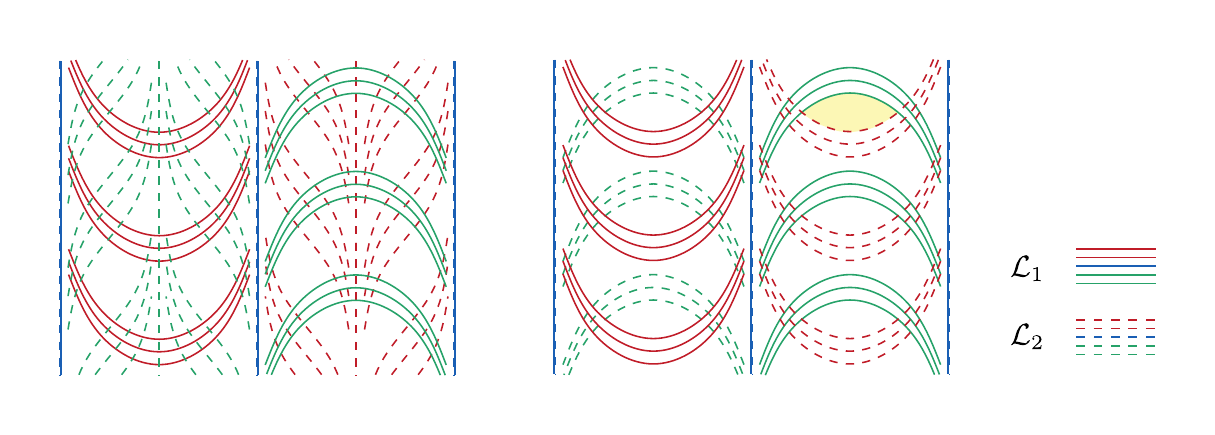}
    \hspace*{-2em}\vspace*{-1.7em}
    \caption{Pairs of quasi-transverse laminations on the torus $\T^2$}
    \label{fig: lam fqt}
\end{figure}

On the left of Figure~\ref{fig: lam fqt}, the laminations $\cL_1$ (in full line) and $\cL_2$ (in dashed line) are \sqt{}{}.
On the right, they are \qt{}, but not \sqt{}.
Indeed, some components of the complementary $\T^2 \ssm (\cL_1 \cup \cL_2)$, for example the one colored in yellow, are not bounded by two disjoint arcs of leaves of $\cL_1$ and two disjoint arcs of leaves of $\cL_2$.

\begin{defi}[Strongly quasi-transverse gluing map] \label{def: sqt gluing map}
Let $(P,X)$ be a \bb{} and $\varphi$ be a gluing map of $(P,X)$.
We say that $\varphi$ is a \emph{\sqt{} gluing map} if $(\cL, \varphi_* \cL)$ is a pair of \sqt{} laminations on $\partial P$.
\end{defi}

\begin{figure}[htb]
    \centering
    \vspace*{-1em}
    \includegraphics[height=0.48\textheight]{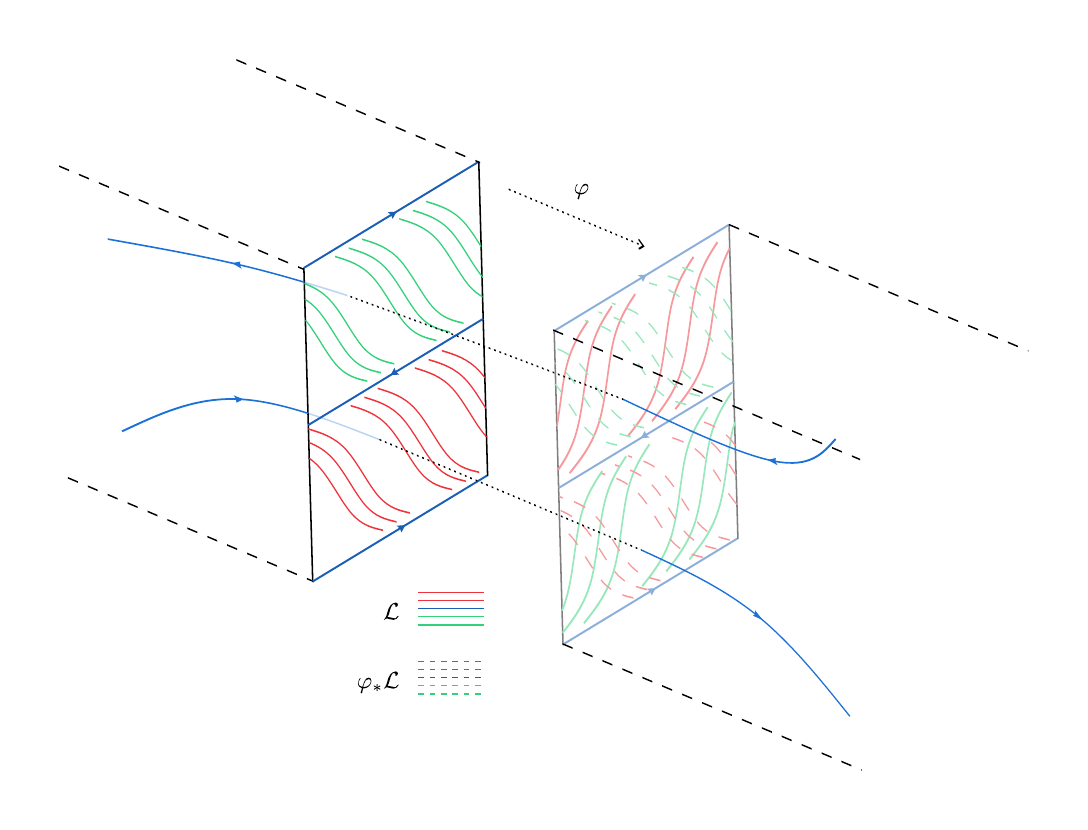}
    \vspace*{-1em}
    \caption{Building block $(P,X)$ with filling boundary lamination $\cL$, and a strongly quasi-transverse gluing map $\varphi$}
    \label{fig: filled block with sqt gluing map}
\end{figure}

\vspace*{-2.5em}
\begin{rmk}
The conditions from Definitions~\ref{def: prefoliation} and\ref{def: sqt gluing map} ensure that we can complete the pair $(\cL, \varphi_* \cL)$ into a pair of \qts{} foliations on $\pP$.
This is a necessary condition for our construction.
Indeed, if a block $(P,X)$ has a dynamical gluing map $\varphi$ such that the vector field $X_\varphi$ on $P_\varphi := P/\varphi$ is Anosov, then the manifold $P$ projects onto the closed manifold $P_\varphi$ and the image of the boundary is an embedded \qt{} surface $S$, equipped with a pair of \qts{} foliations $(\xi^u, \xi^s)$ of dimension~1 whose union contains the projection of the boundary lamination $\cL$ as sublamination.
These are the intersections with $S$ of the unstable and stable foliation $\cF^u$ and $\cF^s$ of the Anosov flow generated by $X_\varphi$.
It follows that each \cc{} of $\pP$ is a closed connected surface of Euler characteristic zero, hence a torus or a Klein bottle.
\end{rmk}

\subsection{Strong isotopy}
\label{sec: preli; subsec: strong isotopy}

\subsubsection*{Equivalence of blocks}
We define the following equivalence relation on \bbs{}.

\begin{defi}[Isotopic \bbs{}]
\label{def: isotopic blocks}
Let $(P_0,X_0)$ and $(P_1, X_1)$ be two \bbs{}.
They are said to be \emph{isotopic} if there exists a family $(P_t, X_t)$ of \bbs{} and a 3-manifold $\tilde P$ such that:
\begin{enumerate}
    \item There exists a family of embeddings $\{ H_t \colon P_0 \to \tilde P \}_{t \in[0,1]}$, continuous in the $\cC^1$ topology, such that $H_0$ is the inclusion and $H_t(P) = P_t$.
    \item There exists a family of vector fields $\{ \tilde X_t \}_{t \in [0,1]}$ on $\tilde P$, continuous in the $\cC^1$ topology, where $\tilde X_t$ is an extension of $X_t$.
\end{enumerate}
We say that $\{(P_t, X_t) \}_{t \in [0,1]}$ is an \emph{isotopy} of \bbs{}.
\end{defi}

\begin{figure}[htb]
    \centering
    \vspace*{-2em}
    \hspace*{-2em}
    \includegraphics[height=0.28\textheight]{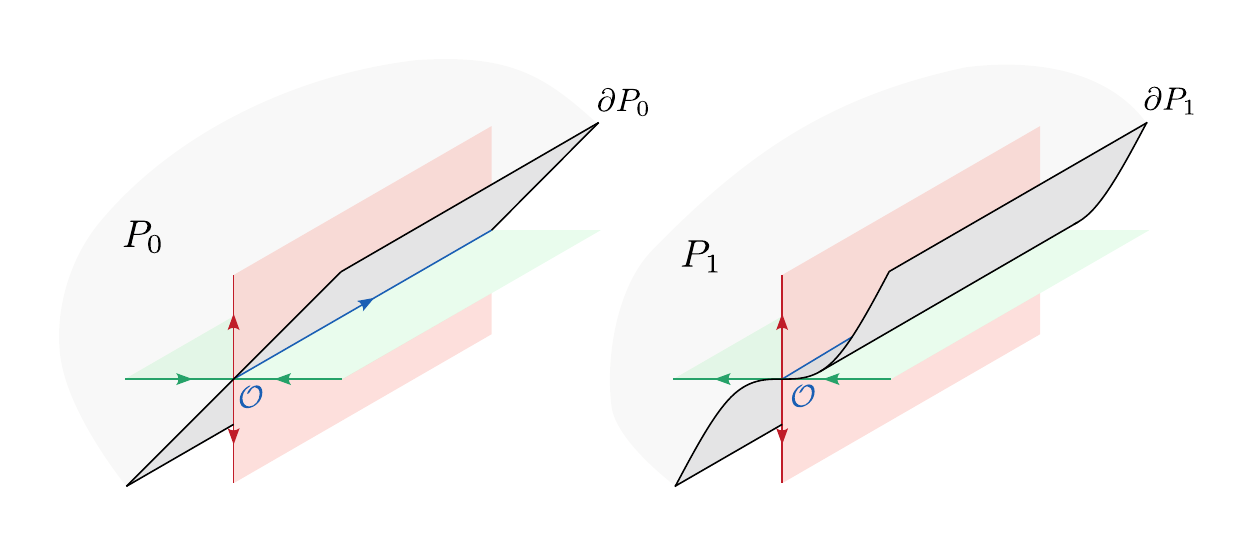}
    \vspace*{-1em}
    \caption{Two isotopic but non-orbit equivalent blocks in the \nbh{} of a common periodic orbit $\cO$ contained in the boundary}
    \label{fig: isotopic but not orbit equivalent block}
\end{figure}

\begin{rmk} \label{rmk: different isotopy definition bby}
In \cite[Definition~3.28]{beguinBuildingAnosovFlows2017}, the authors use a simpler notion of a building block isotopy, for which we have a continuous family of vector fields $\{ X_t \}_{t \in [0,1]}$ all carried by the same manifold $P$ with boundary.
In this article, we will really need to consider families $\{ P_t \}_{t \in [0,1]}$ of manifolds with boundary. More precisely, we will have to consider cases where the position of the boundary of $P_0$ and the position of the boundary of $P_1$ in the neighborhood of a periodic orbit is not the same, as on Figure~\ref{fig: isotopic but not orbit equivalent block}: $\partial P_1$ is transverse to the stable and unstable manifolds of a periodic orbit $\cO \in \cO_*$ contained in $\pP_0$, and $\pP_1$ is topologically transverse but tangent to the stable manifold $\cW^s(\cO)$ along $\cO$.
In the case where the collection $\cO_*$ is empty, two blocks which are isotopic according to Definition~\ref{def: isotopic blocks} will be orbit equivalent to blocks which are isotopic in the sense of~\cite{beguinBuildingAnosovFlows2017}.
\end{rmk}

We will see that this isotopy relation is equivalent to a \say{weakened} orbit equivalence relation.
Let us recall some elementary definitions. An \emph{orbit equivalence} between a vector field $X$ on a manifold $M$ and a vector field $Y$ on a manifold $N$ is a homeomorphism $h \colon M \to N$, which maps the oriented orbits of the flow of $X$ to the oriented orbits of the flow of $Y$. A vector field $X$ on a manifold $M$ is said to be \emph{$\cC^1$-structurally stable} if, for any perturbation $Y$ of $X$ small enough in the $\cC^1$-topology on $M$, $Y$ is orbit equivalent to $X$ on~$M$.

We define a simple continuation of a building block $(P,X)$, which consists in gluing linearizing neighborhoods on the periodic orbits contained in the boundary $\pP$.

\begin{defi}[Minimal continuation]
\label{def: minc}
Let $(P,X)$ be a \bb{} and $\cO_*$ be the collection of periodic orbits of $X$ contained in $\pP$.
We say that $(\tilde P, \tilde X)$ is a \emph{\minc{} of $(P,X)$} if
\begin{enumerate}
    \item $P \subset \tilde P$ and $X = \res{\tilde X}{P}$;
    \item $\cO_* \subset \intr \tilde P$;
    \item there exists a disjoint union $\cV_*$ of linearizing neighborhoods of periodic orbits $\cO_*$ in $\tilde P$ such that
    $\tilde P \ssm P \, \subset \, \cV_*$.
\end{enumerate}
\end{defi}

The following proposition shows that the isotopy relation is equivalent to having a common \minc{}, up to orbit equivalence.

\begin{prop} \label{prop: isotopy vs orbit eq}
Let $(P_0, X_0)$ and $(P_1, X_1)$ be two building blocks.
\begin{enumerate}
    \item \label{prop: isotopy vs orbit eq; it: isotopy imples same minc}
    If $(P_0, X_0)$ and $(P_1, X_1)$ are isotopic, then there exists an orbit equivalence $(P_0',X_0')$ of $(P_0, X_0)$, such that $(P_0',X_0')$ and $(P_1, X_1)$ have a common \minc{} $(\tilde P, \tilde X)$.
    \item \label{prop: isotopy vs orbit eq; it: same minc implies isotopy}
    Conversely, if $(P_0,X_0)$ and $(P_1, X_1)$ are contained in the same \minc{} $(\tilde P, \tilde X)$, then there exists a block isotopy $(P_t, X_t)$ such that $P_t \subset \tilde P$, $X_t = \res{\tilde X}{P_t}$, and $P_t \ssm \cO_*$ is obtained by pushing $P_0 \ssm \cO_*$ along the flow of $\tilde X$.
\end{enumerate}
\end{prop}

Two blocks contained in the same \minc{} have the same collection $\cO_*$ of periodic orbits contained in the boundary, which justifies the previous notation.
We use the following result.

\begin{lem} \label{lem: bby block structurally stable}
Let $(P,X)$ be a \bb{} such that the boundary $\partial P$ is transverse to the vector field $X$.
Then $(P,X)$ is structurally stable, i.e. for any perturbation $Y$ of $X$ small enough in the $\cC^1$-topology on $P$, $(P,Y)$ is a \bb{} orbit equivalent to $(P,X)$.
\end{lem}

\begin{proof}
We use the following result which allows us to embed $P$ into a closed manifold $\cM$ with a structurally stable extension of $X$.
A neighborhood $U \subset \cM$ is said to be \emph{filtrating} for a vector field $Z$ on $\cM$ if the boundary $\partial U$ is transverse to $Z$ and the intersection with any orbit of the flow of $Z$ is connected.

\begin{claim}[{\cite[Appendix~A]{beguinFlotsSmaleDimension2002}}]
\label{claim: smale extension}
There exists a closed 3-manifold $\cM$ equipped with a structurally stable vector field $\bar X$, such that $P$ is embedded in $\cM$, $\bar X$ is an extension of $X$, and $P$ is a filtrating neighborhood of $\Lambda$ for the flow of $\bar X$.
\end{claim}

A $\cC^1$-small perturbation $Y$ of the vector field $X$ extends to a $\cC^1$-small perturbation $\bar Y$ of the \vf{} $\bar X$ on $\cM$.
By structural stability, there exists an orbit equivalence $h: \cM \to \cM$ which maps the oriented orbits of $\bar X$ to the oriented orbits of $\bar Y$, and which is $\cC^0$-close to the identity.
The orbit equivalence $h$ maps $P$ to a filtrating neighborhood of a locally maximal compact hyperbolic invariant set of index $(1,1)$ for the flow of $\bar Y$.
Up to choose a small perturbation, and make an isotopy along the orbits of the flow, we can assume that $h(P)=P$ and the \vf{} $\bar Y$ is transverse to the boundary $\partial P$.
By setting $Y:= \res{\bar Y}{P}$, it follows that $(P,Y)$ is a building block orbit equivalent to $(P,X)$.
\end{proof}

\begin{rmk} \label{rmk: block not structurally stable}
This result is no longer true if the boundary of $P$ is quasi-trans\-verse to the vector field $X$ and contains a nonempty set of periodic orbits $\cO_*$.
Indeed, a small perturbation of such a \bb{} is not always a \bb{} because the condition of quasi-transversality of the boundary to the vector field is not stable by small perturbation.
Moreover, the maximal invariant set $\Lambda$ intersects the boundary $\pP$ along the periodic orbits $\cO_*$.
For this reason one cannot require \emph{a priori} that the orbit equivalence preserves the boundary of $P$.
\end{rmk}

\begin{proof}[Proof of Proposition~\ref{prop: isotopy vs orbit eq}]
Let us show the first item.
Let $\{ (P_t, X_t)\}_{t \in [0,1]}$ be an isotopy of \bbs{}.
Let $(U_t,Y_t)$ be the \bby block associated to the block $(P_t, X_t)$ by Lemma~\ref{lem: from bb to bby} whose boundary is transverse to the vector field, and such that $(P_t, X_t)$ is embedded in $(U_t, Y_t)$.
Each step in the construction of $(U_t, Y_t)$ in the proof of Lemma~\ref{lem: from bb to bby} can be done in a continuous way in $t$, i.e. so that $ \{(U_t, Y_t)\}_{t \in [0,1]}$ is an isotopy of \bbs{}.
The blocks $(U_t, Y_t)$ are structurally stable by Lemma~\ref{lem: bby block structurally stable}.
By Definition~\ref{def: isotopic blocks}, there exists a $\cC^1$-continuous family of embeddings $H_t \colon U_0 \to \tilde U$ such that $H_0$ is the inclusion and $H_t (U_0) = U_t$.
Let $Y'_t$ be the conjugate of $Y_t$ by $H_t$.
It is a vector field on $U_0$ arbitrarily $\cC^1$-close to $Y'_0 = Y_0$.
By structural stability of $(U_0, Y_0)$,
the \vfs{} $Y'_t$ and $Y_0$ are orbit equivalent via a homeomorphism $h'_t :U_0 \to U_0$ with $h'_0$ the identity on $U_0$. Conjugating again by $H_t$ we deduce that the vector fields $Y_t$ and $Y_0$ are orbit equivalent via a homeomorphism $h_t \colon U_0 \to U_t$.
Let $(P_t', X_t')$ be the image by $h_t$ of the block $(P_t, X_t)$ in $U_0$.
Then $X_t'$ coincides with the restriction of $Y_t$ to $P_t'$, and the family $\{ (P_t', X_t') \}_{t \in [0,1]}$ is an isotopy of orbit equivalent \bbs{}.

Let $\cO_*$ be the periodic orbits of $X_0$ contained in $\pP_0$.
Each block $(P_t', X_t')$ has naturally the same collection $\cO_*$ of periodic orbits contained in the boundary.
Let $\cV_*$ be a disjoint union of small tubular neighborhoods of the orbits $\cO_*$ in $U_0$ contained in a union of linearizing \nbhs{} of the orbits $\cO_*$.
Then the family $\{ \pP'_t \ssm \cV_* \}_t$ is an isotopy of compact surfaces uniformly transverse to $\tilde X$ embedded in $\tilde P$.
If we push along the orbits of $\tilde X$, we can suppose that $\pP'_t \ssm \cV_* = \pP'_0 \ssm \cV_* = \pP_0 \ssm \cV_*$.
This operation amounts to an orbit equivalence on each block $(P_t', X_t')$.
We deduce that (up an orbit equivalence) the blocks $P_t'$ coincide outside $\cV_*$, and consequently the blocks $(P_t', X_t')$ are all contained in the same \minc{} $(\tilde P_0, \tilde X_0)$, with $\tilde P_0 = P_0 \cup \cV_*$, and the vector field is the restriction of the vector field $\tilde X$.

 Let us show the second item.
Let $\cO \in \cO_*$ be a periodic orbit of $X$ on $\pP$ and $\cV$ a linearizing neighborhood of $\cO$ in $\tilde P$.
According to Claim~\ref{claim: boundary quadrant and multipliers}, the surfaces $(\partial P_0 \ssm \cO) \cap \cV$ and $(\partial P_1 \ssm \cO) \cap \cV$ are annuli contained in the same quadrant and which coincide in their boundary $\cO \in \cO_*$.
As the flow is trivializable in a quadrant, we deduce the result.
\end{proof}

\begin{rmk}
\label{rmk: not flow isotopic}
In Item~\ref{prop: isotopy vs orbit eq; it: same minc implies isotopy} of Proposition~\ref{prop: isotopy vs orbit eq}, we cannot have an isotopy along the orbits of the flow of $\tilde X$ from $P_t$ to $P_0$.
To see this, let us consider again the case of the blocks $(P_0, X_0)$ and $(P_1, X_1)$ of Remark~\ref{rmk: different isotopy definition bby}.
The isotopy along the flow between $\partial P_0 \ssm \cO$ and $\partial P_1 \ssm \cO$ in this neighborhood does not extend along the periodic orbit $\cO$: the transit time tends to infinity when approaching the orbit $\cO$.
\end{rmk}

\subsubsection*{Strongly isotopic triples}

\begin{defi}[Strongly isotopic triples] \label{def: strongly isotopic triples}
Let $(P_0, X_0)$ and $(P_1, X_1)$ be two \bbs{} with gluing maps $\varphi_0$ and $\varphi_1$.
The triples $(P_0, X_0, \varphi_0)$ and \linebreak[4]$(P_1, X_1, \varphi_1)$ are said to be \emph{isotopic} if
\begin{enumerate}
    \item \label{def: strongly isotopic triples; it: block isotopy}
    There exists an isotopy $\{(P_t, X_t) \}_{t \in [0,1]}$ of \bbs{}.
    \item \label{def: strongly isotopic triples; it: gluing maps isotopy}
    There exists a continuous family $\{ \varphi_t \colon \pP_t \to \pP_t \}_{t \in [0,1]}$ of gluing maps of $(P_t, X_t)$ (in the sense of Remark~\ref{rmk: continuous family of gluing maps}).
\end{enumerate}
They are said to be \emph{strongly isotopic} if in addition
\begin{enumerate}[resume]
    \item \label{def: strongly isotopic triples; it: lamination}
    There exists a continuous family of homeomorphisms
    $$\left\{ h_t \colon \pP_0 \ssm \cO_{0,*} \to \pP_t \ssm \cO_{t,*}
    \right\}_{t \in [0,1]}$$ such that:
    \begin{itemize}[--]
        \item $h_0= \Id$,
        \item $h_t$ maps the lamination $\cL_0 \ssm \cO_{0,*}$ on the lamination $\cL_t \ssm \cO_{t,*}$,
        \item $h_1$ maps the lamination $(\varphi_0)_* (\cL_0 \ssm \cO_{0,*})$ on the lamination \linebreak[4]$(\varphi_1)_* (\cL_1 \ssm \cO_{1,*})$.
    \end{itemize}
\end{enumerate}
\end{defi}

\begin{rmk} \label{rmk: continuous family of gluing maps}
The gluing maps $\varphi_t$ do not share the same domain and destination, so we have to specify what is formally meant by a continuous family.
By definition of a building blocks isotopy $(P_t, X_t)$ (Definition~\ref{def: isotopic blocks}), there exists a $\cC^1$-continuous family of embeddings $H_t \colon P_0 \to \tilde P$ in dimension 3, such that $H_t (P) = P_t$.
Consider $\psi_t \colon \pP_0 \to \pP_0$ the conjugate of $\varphi_t$ by $H_t$, and ask that the family be continuous in the $\cC^1$-topology.
The same remark holds for the family of homeomorphisms $h_t$, where we require a continuous family in the $\cC^0$-topology.
\end{rmk}

\begin{rmk}
In more intuitive terms, the strong isotopy of triples (Definition~\ref{def: strongly isotopic triples}) means that, in addition to a classical isotopy of building blocks and gluing maps, one requires a $\cC^0$-isotopy of the intersection patterns of the laminations.
This last isotopy is only continuous because we will lose regularity of the gluing maps in the course of the proof by conjugating by orbit equivalences which are only continuous and non-trivial on the boundary of the blocks.
\end{rmk}

The definition of strongly isotopic triples ensures that the \say{glued manifolds} are the same, and that the \say{intersection pattern of the laminations} is the same. 
More formally, we have the following proposition.

\begin{prop}\label{prop: properties of strongly isotopic triples}
$\!$Let $(P_0, X_0, \varphi_0)$ and $(P_1, X_1, \varphi_1)$ be two strongly isotopic triples. Then
\begin{enumerate}
    \item \label{prop: strongly isotopic trip, it: same minc}
    Up to orbit equivalence, the blocks $(P_0,X_0)$ and $(P_1, X_1)$ are contained in the same \minc{}.
    
    \item \label{prop: strongly isotopic trip, it: glued manifold homeo}
    The quotient manifolds $P_0 / \varphi_0$ and $P_1 / \varphi_1$ are homeomorphic.
    
    \item \label{prop: strongly isotopic trip, it: lamination pattern}
    There exists a \homeo{} $h \colon \pP_0 \ssm \cO_{*,0} \to \pP_1 \ssm \cO_{*, 1}$ such that
    \begin{itemize}[--]
        \item $h_* (\cL_0 \ssm \cO_{*,0}) = \cL_1 \ssm \cO_{*,1}$,
        \item $h_* (\varphi_0)_* (\cL_0 \ssm \cO_{*,0}) = (\varphi_1)_*( \cL_1 \ssm \cO_{*,1})$.
    \end{itemize}
\end{enumerate}
\end{prop}

\begin{proof}
The first item is a consequence of Proposition~\ref{prop: isotopy vs orbit eq}, Item~\ref{prop: isotopy vs orbit eq; it: isotopy imples same minc}.
The others follow directly from Definition~\ref{def: strongly isotopic triples}
\end{proof}

\begin{figure}[htb]
    \centering
    \vspace*{-1em}
    \includegraphics[height=0.21\textheight]{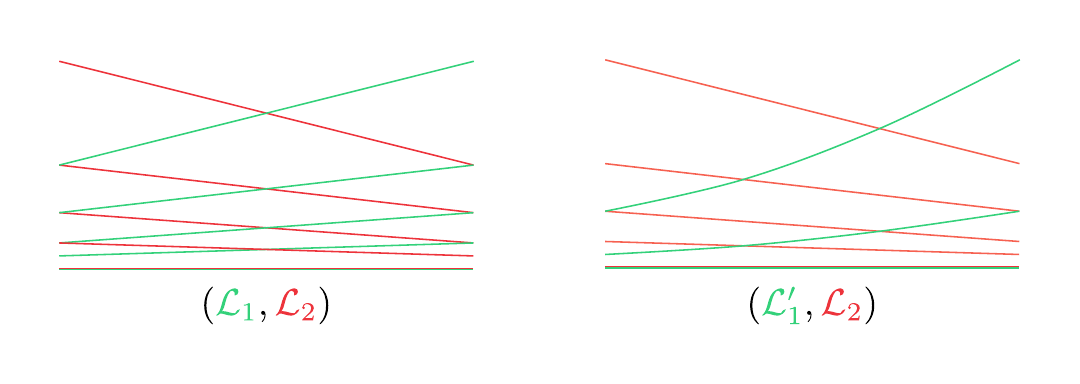}
    \vspace*{-2em}
    \caption{Two laminations $\cL_1$ and $\cL_1'$ which are not isotopic via a family of laminations \sqt{} to~$\cL_2$}
    \label{fig: sqt laminations not strongly isotopic}
\end{figure}

\begin{rmk} \label{rmk: no usual sqt isotopy for equivalent triples}
In Definition~\ref{def: strongly isotopic triples}, one cannot simply ask for the existence of a continuous family of \sqt{} gluing maps $\varphi_t \colon \partial P_t \to \partial P_t$.
Indeed, if $(P_0, X_0)$ and $(P_1,X_1)$ are two blocks contained in the same \minc{} $(\tilde P, \tilde X)$, we want to identify as equivalent the following gluing maps:
\begin{itemize}
    \item $\varphi_0 \colon \partial P_0 \to \partial P_0$, where $\pP_0$ is identified with $x=y$ in linearizing coordinates $(\cV_i, \xi_i = (x,y,\cO_i))$ in the neighborhood of $\cO_i \in \cO_*$, and $\varphi_0$ is defined by $\varphi_0(x,y,\theta) = (-x, -y, \theta)$ in coordinates $\xi_i$ and $\xi_j$;
    \item $\varphi_1: \partial P_1 \to \partial P_1$, where $\pP_1$ is identified with $y = \sqrt{x}$ in linearizing coordinates $(\cV_i, \xi_i = (x,y,\theta))$ in the neighborhood of $\cO_i \in \cO_*$, and $\varphi_1$ is defined $\varphi_1(x,y,\theta) = (-x, -y, \theta)$ in coordinates $\xi_i$ and $\xi_j$.
\end{itemize}

Now these two gluing \diffs{} are not isotopic via a family of \sqt{} gluing maps $\varphi_t \colon \pP_t \to \pP_t$.
Indeed, let us assume by contradiction that there is such an isotopy.
Let $W^s \in \cW^s$ and $W^u \in \cW^u$ be respectively a leaf of the stable and unstable manifold of $\Lambda$, such that if $l^s_0 = W^s \cap \pP_0$ and $l^u_0 = W^u \cap \pP_0$ are the leaves induced on the boundary lamination $\cL_0$ of $\pP_0$, $l^s_0$ and $\varphi_0(l^u_0)$ accumulate on a same periodic orbit $\cO_i \in \cO_*$.
We can without loss of generality suppose that $W^s$ is defined by the equation $\theta = -\ln(x) (\mod \Z)$ in the coordinates $\xi_i$ and $W^u$ is defined by the equation $\theta = \ln(x) (\mod \Z)$ in the corresponding coordinates $\xi_j$.
The set $\{ p_{n,0} \} \subset \pP_0$ of intersection points of $l^s_0 \cap \varphi_0 (l^u_0)$ in the neighborhood $\cV_i$ is ordered (the first intersection, then the second one, etc), and we easily check that there are two convergent subsequences to two distinct points of $\cO_i$.
Let $l^s_t = W^s \cap \pP_t$ and $l^u_t = W^u \cap \pP_t$.
Then, by strong transversality, we can consider the set $\{ p_{n,t} \} \subset \pP_t$ of the intersection points of $l^s_t \cap \varphi_t (l^u_t)$, and by continuity (in $t$) it also contains two convergent subsequences.
However, we can easily check that the gluing map $\varphi_1$ does not satisfy this property:
the sequence of points $\{ p_{n,1} \} \subset \pP_1$ admits three convergent subsequences at three distinct points of $\cO_i$ (Figure~\ref{fig: sqt laminations not strongly isotopic}).
\end{rmk}

\subsection{Statement of the gluing theorem}
We recall that an Anosov flow is a flow $Y^t$ on a closed manifold $\cM$, such that the whole manifold is a hyperbolic set for the flow. More precisely,

\begin{defi}[Anosov flow] \label{def: anosov flow}
Let $\cM$ be a closed manifold equipped with a vector field $Y$ of class $\cC^1$.
We say that $Y^t$ is an \emph{Anosov flow}, or equivalently $Y$ is an \emph{Anosov vector field}, if there exists a $Y^t$-invariant splitting of the tangent bundle into the sum $T\cM = E^\ss \oplus \R.Y \oplus E^\uu$ and constants $\lambda >1$, and $C>0$ such that
\begin{itemize}[--]
    \item $ \forall v \in E^\uu$, $\forall t \geq 0$, $\Vert (Y^t)_* v \Vert \geq C \lambda^t v \Vert$,
    \item $\forall v \in E^\ss$, $\forall t \leq 0$, $\Vert (Y^t)_* v \Vert \geq C \lambda^{-t} \Vert v \Vert$.
\end{itemize}
for a certain Riemannian metric on $\cM$.
\end{defi}

We state the main theorem of this paper.

\begin{mainthm}[Gluing theorem] \label{thm: gluing theorem}
Let $(P,X)$ be a filled building block, and $\varphi$ be a \sqt{} gluing map of $(P,X)$.
There exists a triple $(P_1, X_1, \varphi_1)$ strongly isotopic to $(P,X,\varphi)$ such that
$X_1$ induces an Anosov vector field on the closed 3-manifold $P_1/\varphi_1$.
\end{mainthm}

Sections~\ref{sec: normalization} to~\ref{sec: proof gluing thm} are devoted to the proof of this theorem, which is organized as follows.
In Section~\ref{sec: normalization}, we show that we can put the initial triple $(P,X, \varphi)$ of the theorem in a \say{normalized} form.
In Section~\ref{sec: crossing map}, we study the hyperbolic properties of the crossing map of the flow from the entrance boundary $\Pin$ to the exit boundary $\Pout$.
In Section~\ref{sec: spreading}, we show how to \say{spread hyperbolicity} by a process of coordinate change on the boundary of $P$.
In Section~\ref{sec: parameters and cones}, we show that we can use this change of coordinates to modify the gluing map $\varphi$ in order to create hyperbolicity along the new recurrent orbits of the flow obtained after the gluing process, and this in a way compatible with the natural hyperbolicity of the initial flow.
In Section~\ref{sec: proof gluing thm}, we show that for such a choice of gluing map, the flow induced by the initial flow on the glued manifold is Anosov, which completes the proof of Theorem~\ref{thm: gluing theorem}.

\section{Normalization}
\label{sec: normalization}

In this section, we consider a triple $(P,X,\varphi)$ consisting of a filled \bb{} $(P,X)$ and a \sqt{} gluing map $\varphi$ (see Section~\ref{sec: preliminaries} for the definitions).
Denote $\cO_* = \{ \cO_1, \dots, \cO_n\}$ the collection of periodic orbits of $X$ contained in $\pP$, $\Lambda$ the maximal invariant set of $(P,X)$, $\cW^s=\cW^s(\Lambda)$ and $\cW^u=\cW^u(\Lambda)$ the stable and unstable manifolds of $\Lambda$.

The goal of this section is to modify the triple $(P,X,\varphi)$ by strong isotopy in order to obtain a \say{normalized} form.
In short, we want to:
\begin{itemize}
    \item linearize by pieces the dynamics on a neighborhood of the \emph{saddle basic pieces} of $\Lambda$;
    \item adjust the multipliers of the orbits $\cO_i \in \cO_*$ to $\{ \frac{1}{2}, 2 \}$;
    \item construct a pair of foliations $(\cG^s, \cG^u)$ on the block which extends the laminations $(\cW^s, \cW^u)$;
    \item put the boundary of the block in canonical position in the neighborhood of the orbits $\cO_i \in \cO_*$;
    \item explicitly control the gluing map $\varphi$ in the neighborhood of the orbits $\cO_i \in \cO_*$.
\end{itemize}
To formalize this, we need to define a notion of normalized \bb{} (Definition~\ref{def: normalized block}) and of normalized gluing map (Definition~\ref{def: normalized gluing map}).
We can then state the normalization result (Proposition~\ref{prop: normalization of triple}).

\begin{rmk}
    In \cite{beguinBuildingAnosovFlows2017}, the building blocks which are considered are \emph{saddle}: they contains neither attractors or repellers.
    The reason why the authors only considered saddle (\bby{}) blocks is to be able to linearize the dynamic in the \nbh{} of the maximal invariant set, which is a technical tool for the proof.
    We will show here that we can deal separately with the attractors and the repellers while doing the same construction and use the same idea of proof.
\end{rmk}

\subsection{Definitions and statement}
\label{sec: normalization; subsec: statement}

\subsubsection*{Normalized block}$\!\!$Let $(\tilde P, \tilde X)$ be a \minc{} of $(P,X)$ (Definition~\ref{def: minc}).

\begin{nota}[Attracting, Repelling, Saddle basic pieces] \label{nota: repellers attractors saddle pieces}
    Denote $\gR$, $\gS$ and $\gA$ the union of repelling, saddle and attracting basic pieces respectively.
\end{nota}

\begin{defi}[Maximal saddle invariant set] \label{def: msis}
The \emph{maximal saddle invariant set} $\Lambda_s$ is the compact invariant set in $\Lambda$ containing neither attractors nor repellers and which is maximal for this property. Equivalently, $\Lambda_s$ is equal to the union of $\gS$ and the intersections of the invariant manifolds of $\gS$.
\end{defi}

\begin{defi}[Local transverse section] \label{def: local section}
We say that $\Sigma \subset \tilde P$ is a \emph{local transverse section} of $\Lambda_s$ for the flow of $\tilde X$, if
\begin{enumerate} 
    \item $\Sigma$ is a finite union of disjoint compact embedded disks transverse to $\tilde X$;
    \item $\Sigma$ intersects any orbit of $\Lambda_s$; 
    \item $\partial \Sigma \cap \Lambda_s = \emptyset$.
\end{enumerate}
\end{defi}

If $D$ is a compact subset of $\Sigma$, we say that the first return of $D$ is \emph{well defined} if the positive orbit of any point of $D$ intersects the interior of $\Sigma$ in continuous time. If $D$ has this property, we can \emph{renormalize} the vector field $\tilde X$ so that the return time of $D$ on $\Sigma$ is equal to $1$.

\begin{defi}[Normalized building block]
\label{def: normalized block}
We will say that a block $(P,X)$ is \emph{normalized} if there exists a \minc{} $(\tilde P, \tilde X)$ such that
\begin{enumerate}

    \item \label{def: normalized, it: affine section}
    \emph{(affine coordinates)}
    There exists a local transverse section $\Sigma$ of $\Lambda_s$ in $\tilde P$, and a cover of $\Sigma \cap \Lambda_s$ in $\Sigma$ into a finite collection of disjoint disks $D = \{ D_1, \dots, D_n \}$, provided with a coordinate system $\chi_i = (x, y) \colon D_i \to \R^2$ on each $D_i$ such that:
    \begin{enumerate}
        \item \label{def: normalized, it: affine section, it: return map}
        the first return of the disks of $D$ on $\Sigma$ is well defined,
        \item \label{def: normalized, it: affine section, it: unique disq}
        for any orbit $\cO_i \in \cO_*$, there exists a unique disk $D_i \in D$ which intersects $\cO_i$ at a single point,
        \item \label{def: normalized, it: affine section, it: affine coord}
        if $f$ is the first return map of the flow of $\tilde X$ on $\Sigma$, then the intersection $D_i \cap f_\Sigma^\inv(D_j)$ is connected and the restriction of $f$ to $D_i \cap f_\Sigma^\inv(D_j)$ is diagonal affine in the coordinates $\chi_i$ and $\chi_j$ (in the sense of Remark~\ref{rmk: diagonal affine return}).
    \end{enumerate}
    We say that $\Sigma$ is an \emph{affine section of $\Lambda_s$}, and $\{ (D_i, \chi_i)\}$ is a \emph{\acs{}} on $\Sigma$.

    \item \label{def: normalized, it: multipliers}
    \emph{(multipliers)}
    The multipliers of the orbits $\cO_i \in \cO_*$ are equal to $\{ \frac{1}{2}, 2\}$.

    \item \label{def: normalized, it: foliation}
    \emph{(foliation)}
    There exists a pair of smooth 2-dimensional foliation
    $(\cG^s, \cG^u)$ on $\tilde P$ satisfying :
\begin{enumerate}
    \item \label{def: normalized, it: foliation, it: complete lam}
    $\cG^s$ and $\cG^u$ are $\tilde X^t$-invariant, transverse to each other, and contain $\cW^s$ and $\cW^u$ as sublamination respectively,
    \item \label{def: normalized, it: foliation, it: affine}
    The foliations $\cG^s \cap \Sigma$ and $\cG^u \cap \Sigma$ coincide with $x= \cst$ and $y= \cst$ in the affine coordinates on $\Sigma$ given by Item~\ref{def: normalized, it: affine section}.
\end{enumerate}
    We say that $(\cG^s, \cG^u)$ is a \emph{pair of affine invariant foliations} of $(P,X)$.
    \item \label{def: normalized, it: straight boundary}
    \emph{(straightened boundary)}
    For any orbit $\cO_i \in \cO_*$, there exists a
    \lcs{} $(\cV_i, \xi_i = (x,y,\theta) \in \R^2 \times \R/\Z)$ of $\cO_i$ for the flow of $\tilde X$, compatible with the \acs{} of the section $\Sigma$, and in which the boundary coincides with the diagonal $\{x=y\}$. More precisely,
    \begin{enumerate}
        \item \label{def: normalized, it: straight boundary, it: compatible coord}
        if $D_i$ is the disk given by Item~\ref{def: normalized, it: affine section, it: unique disq} which intersects $\cO_i$, then it coincides with $\{ \theta = 0 \}$, and the coordinates $\xi_i$ and $\chi_i$ coincide on $D_i$;
        \item \label{def: normalized, it: straight boundary, it: diagonal}
        $P$ coincides with the region $\{ y \geq x \}$ in the chart $\xi_i$ on $\cV_i$.
    \end{enumerate}
    We say that $\cV_i$ is a \emph{normalized neighborhood} of $\cO_i$ and $(\cV_i, \xi_i)$ is a \emph{\ncs{}} of $\cO_i$.
\end{enumerate}
\end{defi}

\begin{figure}[htb]
     \centering
     \vspace*{-1em}
     \includegraphics[height=0.33\textheight]{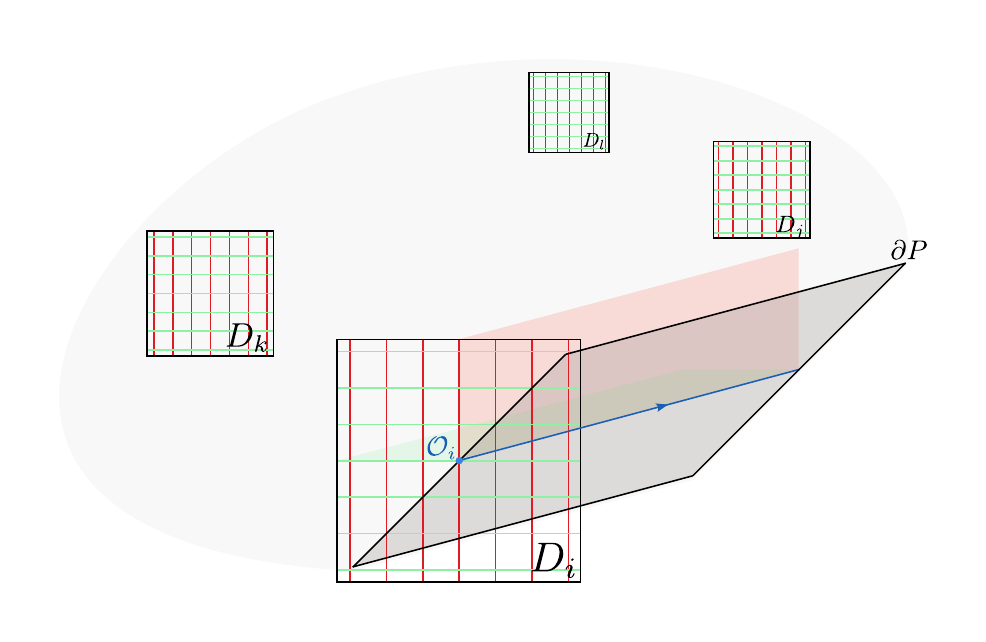} 
     \vspace*{-1em}
     \caption{Normalized block $(P,X)$ in the neighborhood of a periodic orbit $\cO_i \in \cO_*$ in a \ncs{} $(\cV_i, \xi_i)$}
     \label{fig: normalized block on snc}
 \end{figure}

We refer to Figure~\ref{fig: normalized block on snc}.

\begin{rmk} \label{rmk: diagonal affine return}
Item~\ref{def: normalized, it: affine section} is formally stated as follows.
For any pair $(i,j)$ such that $D_i \cap f^\inv(D_j) \neq \emptyset$, there exists  $(a_{i,j}, b_{i,j}) \in \R^2$ and real numbers $0< \vert \lambda_{i,j} \vert <1$ and $\vert \mu_{i,j} \vert >1$ such that for all $p \in D_i \cap f^\inv(D_j)$, if $(x,y) = \xi_i(p)$, we have
$$(\xi_j \circ f \circ \xi_i^\inv) (x, y) = (a_{i,j} + \lambda_{i,j}.x,\ b_{i,j} + \mu_{i,j}.y).$$
\end{rmk}

\pagebreak[4]
\subsubsection*{Properties of a normalized block}

\begin{coro} \label{coro: expression flow in normalized neighborhood}
Let $(P,X)$ be a normalized \bb{} and $(\tilde P, \tilde X)$ a \minc{}.
For any periodic orbit $\cO_i \in \cO_*$ contained in $\pP$, if $(\cV_i, \xi_i = (x,y,\theta))$ is a \ncs{} of $\cO_i$ in $\tilde P$ (Definition~\ref{def: normalized block}, Item~\ref{def: normalized, it: straight boundary}), then the expression of the flow of $\tilde X$ in these coordinates is
$$ \tilde X^t (x,y,\theta) = (2^{-t} x,\, 2^t y,\, \theta +t).$$
\end{coro}

\begin{proof}
By assumption (Definition~\ref{def: normalized block}, Item~\ref{def: normalized, it: straight boundary, it: compatible coord}), $\cV_i$ is contained in the suspension neighborhood of the disk $D_i$ which intersects $\cO_i$, and the flow of $\tilde X$ in coordinates $\xi_i$ coincides with the suspension flow of the first-return map $f$ on the disk $D_i$ in coordinates $\chi_i$.
According to Item~\ref{def: normalized, it: affine section, it: affine coord}, this map is affine in each coordinates of the system $\chi_i = (x,y)$ on $D_i \cap f^\inv(D_i)$,
and according to Item~\ref{def: normalized, it: multipliers}, it has multipliers $\frac{1}{2}$ and $2$.
The statement follows.
\end{proof}

\begin{rmk} \label{rmk: equation affine foliation}
Let $(P,X)$ be a normalized block, and $(\cG^s, \cG^u)$ a \paif{} on $(P,X)$ (Definition~\ref{def: normalized block}, Item~\ref{def: normalized, it: foliation}).
The \ncs{} $\xi_i = (x,y,\theta) \colon \cV_i \to \R^2 \times \R/\Z$ induce the coordinate system
$ \rho_i = (x, \theta) \colon \pP \cap \cV_i \to \R \times \R/\Z$ on the boundary of $P$ by \say{forgetting} the $y$ coordinate.
In these coordinates, if $\tilde \theta \in \R$ is a lift of the coordinate $\theta \in \R/\Z$, then the leaves of the foliation $\cG^u \cap \partial P$ are the curves
    $x = \cst 2^{-\tilde \theta}$ and the leaves of the foliation $\cG^s \cap \partial P$ are the curves
    $x = \cst 2^{\tilde \theta}$.
\end{rmk}

\begin{figure}[htb]
    \centering
    \vspace*{-1em}
    \includegraphics[width=0.55\textwidth]{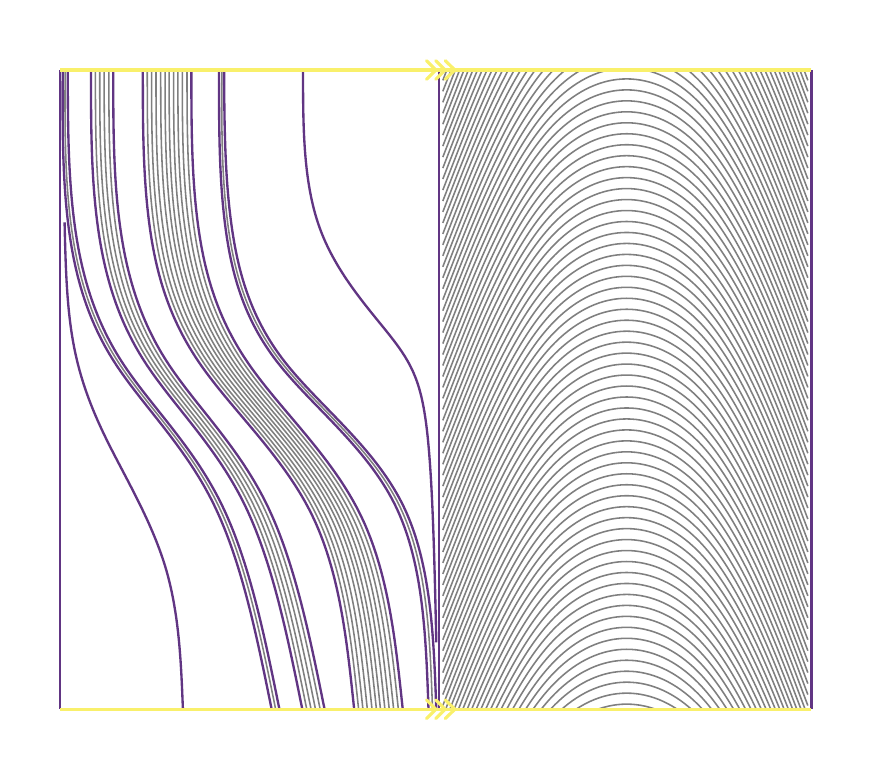}
    \vspace*{-1em}
    \caption{Boundary (in violet) of a filling \qms{} lamination on an annulus}
    \label{fig: boundary of lamination}
\end{figure}

\begin{defi}[Boundary of a lamination] \label{def: boundary of a lamination}
Let $\cL$ be a lamination on a surface $S$.
The \emph{boundary $\partial \cL$ of $\cL$} is the complementary $\cL \ssm \intr_S \cL$.
It is a closed sublamination of $\cL$ in $S$ with empty interior.
A leaf $l \in \partial \cL$ is called a \emph{boundary leaf} of~$\cL$.
\end{defi}
See Figure~\ref{fig: boundary of lamination}.

We state now the following result.

\begin{prop} \label{prop: boundary foliation induced by paif}
Let $(P,X)$ be a building block equipped with an affine section $\Sigma$ and a \paif{} $(\cG^s, \cG^u)$.
Then
\begin{enumerate}
    \item \label{prop: boundary foliation, it: out}
    $\cG^{u,\out}:= \cG^u \cap \Pout$ and $\cG^{s, \out} := \cG^s \cap \Pout$ are smooth one-dimensional foliations on $\Pout$ transverse to each other;
    
    \item \label{prop: boundary foliation, it: in}
    $\cG^{u,\iin} := \cG^u \cap \Pin$ and $\cG^{s, \iin} := \cG^s \cap \Pin$ are smooth one-dimensional foliations on $\Pin$ transverse to each other;
   
    \item \label{prop: boundary foliation, it: Gb}
    $\cF := \cG^{u, \out} \, \cup \, \cO_* \, \cup \, \cG^{s, \iin}$ is a smooth one-dimensional foliation of $\pP$, which contains the boundary lamination $\cL$ as sublamination.
    The holonomy of compact leaves of $\partial \cL$ for the foliation $\cF$ are conjugated to affine maps.\footnote{Of each of the two sides of the compact leaf when the compact leaf is two sided.}

\end{enumerate}
\end{prop}

\begin{rmk} \label{rmk: nbh of end of strip with affine holonomy}
    Every \cc{} of $\pP \ssm \cL$ which is a strip $B$ accumulates at both end on a compact boundary leaf (Definition~\ref{def: boundary of a lamination}). 
    It follows from Proposition~\ref{prop: boundary foliation induced by paif}, Item~\ref{prop: boundary foliation, it: Gb} that there is a \nbh{} of the ends of $B$ where the holonomy maps of the foliation $\cF$ are conjugated to affine maps.
    This is a key ingredient for the proof of \cite[Theorem~1.5]{beguinBuildingAnosovFlows2017}, which was obtained by considering only \emph{saddle} \bby{} blocks and linearizing globally the dynamic in the \nbh{} of the maximal invariant set.
    We bypass this problem by linearizing only in the \nbh{} of the maximal \emph{saddle} invariant set, and get the same property.
\end{rmk}

We will need general preliminary results, valid for any building block $(P,X)$.
Let $O$ be an orbit of $X$ and $\cW^s(O)$ the stable manifold of $O$.
A \emph{stable separatrix of $O$} is a \cc{} of $\cW^s(O) \ssm O$, denoted $\cW^s_+(O)$.
We say that a stable separatrix $\cW^s_+(O)$ of $O$ is \emph{free} if it is disjoint from $\Lambda$.
Similarly, we define an unstable separatrix and a free unstable separatrix of $O$.

\begin{lem}\label{lem: free separatrix}
Any compact leaf $\gamma$ of $\cL^\iin\,$\footnote{$\gamma$ is not a periodic orbit of $X$.} is the transverse intersection of a free stable separatrix $\cW^s_+(O)$ of a periodic orbit $O$ of $X$ with $\Pin$.
Similarly, any compact leaf $\gamma$ of $\cL^\out$ is the transverse intersection of a free unstable separatrix $\cW^u_+(O)$ of a periodic orbit $O$ with $\Pout$.
\end{lem}

\begin{proof}[Proof of Lemma~\ref{lem: free separatrix}]
Each orbit of $\cW^s \ssm \Lambda$ transversely intersects $\Pin$ at a single point.
Let $C$ be the \cc{} of $\cW^s \ssm \Lambda$ such that $\gamma = C \cap \Pin$.
The orbit space on $C$ is a circle, so $C$ is a cylinder trivially foliated by the flow.
It follows that $C$ is contained in the stable manifold of a periodic orbit $O$.
The orbit $O$ is the unique periodic orbit of $X$ in $\cW^s(O)$, and $C$ is flow invariant. It follows that $C$ is a \cc{} of $\cW^s(O) \ssm O$, in other words a stable separatrix of the periodic orbit $O$.
Since $C$ is disjoint from $\Lambda$, this separatrix is free.
\end{proof}

Recall that $\Lambda_s$ denote the \msis{} (Definition~\ref{def: msis}).

\vspace*{-0.2em}
\begin{lem} \label{lem: boundary leaf and periodic orbit are saddle}
    The boundary sublamination $\partial \cL$ of $\cL$ is contained in the lamination induced by the \msis{} $\Lambda_s$.
    Any periodic orbit $\cO_i \in \cO_*$ contained in $\pP$ belongs to $\Lambda_s$.
\end{lem}

\vspace*{-0.65em}
\begin{proof}
    The stable manifolds of the attractors $\gA$ in $P$ foliates an open set of $P$ (the union of basins of the attractors) with non-empty interior, the stable manifolds of $\Lambda_s$ forms a lamination with empty interior in $P$, and the stable manifolds of $\gR$ coincide with $\gR$ so never intersects the boundary. 
    It follows that the interior of $\cL^\iin$ in $\Pin$ contains the intersection $\cW^s(\gA) \cap \Pin$ and its complementary $\partial \cL^\iin$ contains $\cW^s(\Lambda_s) \cap \Pin$.
    Similarly, the interior of $\cL^\out$ in $\Pout$ contains the intersection $\cW^u(\gR) \cap \Pout$ and its complementary $\partial \cL^\out$ contains $\cW^u(\Lambda_s) \cap \Pout$.
    The boundary lamination $\cL = \cL^\iin \cap \cO_* \cap \cL^\out$ has boundary $\partial \cL$ contained in the union $\partial \cL^\iin \cup \cO_* \cup \partial \cL^\out$.
    Therefore it remains to show that $\cO_*$ is contained in $\Lambda_s$ in order to prove the lemma.
    Let $\cO \in \cO_*$ and consider $\cV$ a small tubular \nbh{} of $\cO$.
    Suppose that $\cO$ belongs to an attracting basic piece $\gA_i$. 
    Consider a point $p \in \gA_i$ in $\cV$ disjoint from $\cO$.
    The unstable manifold $\cW^u(p)$ crosses $\Pout$, so there are orbits of $\cW^u(\gA_i) =\gA_i \subset \Lambda$ which exit the block, this is a contradiction.
    We deduce that $\cO$ is an isolated periodic orbit, so it is a saddle basic piece and the first assumption is wrong.
    The argument are similar for proving that $\cO$ can not belong to a repelling basic piece.
    It follows that the second assertion, hence the lemma, is true.
\end{proof}

\vspace*{-0.65em}
\begin{proof}[Proof of Proposition~\ref{prop: boundary foliation induced by paif}]
The foliations $\cG^s$ and $\cG^u$ are smooth and transverse to each other, tangent to the \vf{} $X$, and the surfaces $\Pin$ and $\Pout$ are transverse to $X$.
We deduce the first two items.
Let $\cF := \cG^{u, \out} \cup \, \cO_* \, \cup \, \cG^{s, \iin}$.
The lamination $\cL^\out$ is a sublamination of $\cG^{u, \out}$ on $\Pout$ and $\cL^\iin$ is a sublamination of $\cG^{s, \iin}$ on $\Pin$.
The foliation $\cG^{u, \out}$ on $\Pout{}$ and $\cG^{s, \iin}$ on $\Pin{}$ reconnect with $\cO_*$ to form a $\cC^1$-foliation which contains the boundary lamination \mb{$\cL = \cL^\out \cup \cO_* \cup \cL^\iin$} as sublamination (see the equations of Remark~\ref{rmk: equation affine foliation}).

Let $\gamma$ be a compact boundary leaf of $\cL \subset \cF$.
If $\gamma$ is a leaf of $\cL^\iin$, then it is the intersection of a cylindrical free stable separatrix $\cW^s_+(O)$ of a periodic orbit $O \in \Lambda_s$ according to Lemma~\ref{lem: free separatrix} and Lemma~\ref{lem: boundary leaf and periodic orbit are saddle}.
The paths $\gamma$ and $O$ are freely homotopic in $\cW^s_+(O)$ (as non-oriented paths).
The holonomy of $O$ for the foliation $\cG^s$ is conjugate to the first return map of the flow on a transversal of $O$ in $\cW^u(O)$.
Since the flow of $X$ admits an affine section in the neighborhood of the \msis{} $\Lambda_s$, this map is conjugate to an affine map (expanding for the orientation induced by the flow).
We deduce that the holonomy maps of $\gamma$ for the foliation $\cG^{s, \iin}$ are conjugate to affine maps (contracting for the orientation induced by the flow).
The same is true for a compact leaf of $\cL^\out$.
In the case $\gamma$ is a periodic orbit $\cO_i \in \cO_*$, we consider the holonomy on one side and on the other of $\cO_i$.
One corresponds to the holonomy of the foliation $\cG^s$ along $\cO_i$ (therefore expanding for the flow orientation) and the other to the holonomy of the foliation $\cG^u$ along $\cO_i$ (therefore contracting), and they are conjugated to affine maps by the same proof. 
\end{proof}

\subsubsection*{Normalized gluing map}

\begin{defi}[Normalized gluing map and triple] \label{def: normalized gluing map}
Let $(P,X)$ be a normalized block with a \paif{} $(\cG^s, \cG^u)$. Denote $\cG^{s, \iin}$ and $\cG^{u, \out}$ the induced foliations on $\Pin$ and $\Pout$.
We say that a gluing map $\varphi \colon \pP \to \pP$ of $(P,X)$ is \emph{normalized} if
\begin{enumerate}
    \item \label{def: normalized gluing; it: trivial}
    For each periodic orbit $\cO_i$ in $\cO_*$, there exist \ncss{} $(\cV_i, \xi_i = (x,y,\theta))$ of $\cO_i$ and $(\cV_j, \xi_j = (x,y,\theta))$ of  $\cO_j = \varphi(\cO_i)$ such that the expression for $\varphi$ in these coordinates is
    \begin{equation*}
        \varphi(x,x,\theta) = (-x,-x,\theta).
    \end{equation*}
    \item \label{def: normalized gluing; it: foliation transverse}
    The foliation $\varphi_* \cG^{u, \out}$ is transverse to the foliation $\cG^{s, \iin}$ on $\Pin$.
\end{enumerate}
Then we say that $(P,X,\varphi)$ is a normalized triple.
\end{defi}

\begin{rmk}
The boundary $\pP$ has equation $x=y$ in a \ncs{} $(\cV_i, \xi_i)$, which justifies the expression of Item~\ref{def: normalized gluing; it: trivial}.
More formally (but we don't write it down for sake of simplicity), this expression means that for any $(x, x, \theta)$ in the image $\xi_i (\cV_i \cap \partial P)$, we have
$\xi_j \circ \varphi \circ \xi_i^\inv (x,x,\theta) = (-x,-x,\theta)$.
\end{rmk}

\begin{prop} \label{prop: normalized gluing map is dynamic}
A gluing map $\varphi$ of a normalized block $(P,X)$ satisfying Item~\ref{def: normalized gluing; it: trivial} of Definition~\ref{def: normalized gluing map} is a dynamical gluing map of $(P,X)$ (in the sense of Definition~\ref{def: dynamical gluing map}).
\end{prop}

\begin{proof}
Let $(\tilde P, \tilde X)$ be a \minc{} of $(P,X)$.
It must be shown that $\varphi$ is the restriction of a \diff{} $\tilde \varphi$ on a neighborhood of $\pP$ in $\tilde P$ which maps the vector field $\tilde X$ to itself.
Let $\cO_*$ be the collection of periodic orbits of $X$ contained in $\pP$.
For each periodic orbit $\cO_i \in \cO_*$, if $(\cV_i, \xi_i)$ and $(\cV_j, \xi_j)$ denote the \ncss{} in the neighborhood of $\cO_i$ and $\cO_j := \varphi(\cO_i)$, then the formula $(x,y,\theta) \mapsto (-x,-y,\theta)$ defines an extension $\tilde \varphi$ of $\varphi$ to a neighborhood of $\partial P$ in $\cV_i$ and $\cV_j$.
Let $\cV_*$ be the union of the normalized neighborhoods $\cV_i$ of $\cO_i$, and $\cV$ a neighborhood of $\cO_*$ strictly contained in $\cV_*$.
    Then $\pP \ssm \cV$ is a compact surface transverse to $X$.
    By the flow box theorem for $\cC^1$ vector fields, there exists a small tubular neighborhood $U$ of $\pP \ssm \cV$ with coordinates $(p,t) \in (\pP \ssm \cV) \times (-1, 1)$ in which the vector field $\tilde X$ is the trivial field $\partial/ \partial t$ and $P$ is identified with the set $\{ -1 < t \leq 0 \}$.
    We extend $\varphi$ by the map $(p,t) \mapsto (\varphi(p), -t)$.
    These two local definitions are compatible.
    Indeed the set $\pP \cap (\cV_* \ssm \cV)$ is uniformly transverse to $\tilde X$ and there is only one way to extend $\varphi$ on a tubular neighborhood by a \diff{} which maps the vector field $\tilde X$ to itself.
We thus obtain coordinate systems on a tubular neighborhood of $\pP$ in $\tilde P$ in which the \diff{} $\varphi$ extends to a \diff{} $\tilde \varphi$ which maps the \vf{} $\tilde X$ to itself.
\end{proof}

\subsubsection*{Statement}

The aim of this section is to show the following proposition.
Recall that a \bb{} $(P,X)$ is said to be saddle if its maximal invariant set $\Lambda$ contains neither attractors nor repellers, and filled if the \ccs{} of the complementary of the boundary lamination $\pP \ssm \cL$ are strips.

\begin{prop}[Normalization of a triple]
\label{prop: normalization of triple}
Let $(P_0,X_0, \varphi_0)$ be a filled building block, provided with a \sqt{} gluing map.
Then there exist a normalized triple $(P,X,\varphi)$ strongly isotopic to $(P_0, X_0, \varphi_0)$
\end{prop}

We refer to Definition~\ref{def: strongly isotopic triples} for the strong triple isotopy.

\begin{rmk} \label{rmk: normalized triple is filled saddle and sqt}
The \bb{} $(P,X)$ obtained by Proposition~\ref{prop: normalization of triple} is filled and the gluing map $\varphi$ is \sqt{} (indeed those properties are preserved under strong isotopy).
\end{rmk}

\begin{rmk}
We follow closely the techniques of \say{normalization} of blocks and gluing maps of \bby{} (\cite[Section~5]{beguinBuildingAnosovFlows2017}).
The first difference is that our blocks can contains attractors and repellers. 
We will isolate the attractors and repellers and we recover the same techniques of linearization by considering the dynamic near the \msis{}.
The other difference is that we want to \say{straighten} the boundary in the \nbh{} of the periodic orbits to put it in the canonical position (Item~\ref{def: normalized, it: straight boundary} of Definition~\ref{def: normalized block}).
This operation done at Subsection~\ref{sec: normalization; subsec: boundary straightening} is typically obtained by pushing the boundary along the flow for an infinite time in the neighborhood of the periodic orbits contained in the boundary.
This takes us out of the orbit equivalence class of the initial block, and the gluing map no longer extends to the periodic orbits.
We will say that we obtain an \emph{incomplete gluing map}.
Then we have to show that we can modify this incomplete gluing map in its strong isotopy class so that it extends into a true gluing map.
On the other hand, some technical difficulties are added by the fact that the blocks we work with are not structurally stable, unlike the \bby{} blocks.
Therefore, perturbing the vector field of a block, as it is done in the sections~\ref{sec: normalization; subsec: affine section} and~\ref{sec: normalization; subsec: multipliers} to linearize the flow and adjust the multipliers (Item~\ref{def: normalized, it: affine section} and~\ref{def: normalized, it: multipliers}, Definition~\ref{def: normalized block}), will force us to consider a new block carried by another manifold.
\end{rmk}

\pagebreak[4]
\subsubsection*{Proof summary}
The proof of this proposition is organized as follows.
\begin{itemize}[leftmargin=*]
    \item \emph{(affine section)} 
    In Subsection~\ref{sec: normalization; subsec: affine section}, we show that we can perturb a block $(P,X)$ in its orbit equivalence class, so as to linearize the dynamics in the neighborhood of the \msis{} $\Lambda_s$.
    More precisely, we show that there exists a block $(P_1, X_1)$ isotopic to $(P,X)$ among the orbit equivalent blocks which admits an affine section $\Sigma$ of $\Lambda_s$ (Proposition~\ref{prop: affine section}).
    This block then satisfies Item~\ref{def: normalized, it: affine section} of Definition~\ref{def: normalized block}.
    The key point is that the return map on a local transverse section $\Sigma$ of the \msis{} $\Lambda_s$ is arbitrarily close to an affine map if the section is small enough.
    One can then use a (weakened) structural stability result on building blocks (Proposition~\ref{prop: weak structural stability}).

    \item \emph{(multipliers)} 
    In Subsection~\ref{sec: normalization; subsec: multipliers}, we show that any block $(P,X)$ which admits an affine section is isotopic to an orbit equivalent block $(P_1,X_1)$ which admits an affine section and such that the multipliers of the periodic orbits contained in the boundary are equal to $\{\frac{1}{2}, 2\}$, in other words, which satisfies Items~\ref{def: normalized, it: affine section} and~\ref{def: normalized, it: multipliers} of Definition~\ref{def: normalized block}.
    This is Proposition~\ref{prop: multipliers}.
    The proof consists in continuously modifying the return map on a Markov partition in the neighborhood of the periodic orbits of the boundary.

    \item \emph{(pair of affine invariant foliations)} 
    In Subsection~\ref{sec: normalization; subsec: pair of affine invariant foliation}, we show that we can foliate a filled \bb{} $(P,X)$ which admits an affine section $\Sigma$ by a \paif{} $(\cG^s, \cG^u)$, in other words which satisfies Item~\ref{def: normalized, it: foliation} of Definition~\ref{def: normalized block}.
    This is Proposition~\ref{prop: existence paif}.
    The key point is that one can push through the flow of $X$ the horizontal and vertical foliations induced by the affine coordinates on $\Sigma$, in order to obtain a pair of transverse and $X^t$-invariant foliations on an invariant neighborhood of $\Lambda_s$ in $P$, which contain the stable and unstable manifolds of $\Lambda_s$.
    Since the block $(P,X)$ is a filled block, there is a unique way to complete these foliations in a \paif{} on $P$.
   
    \item \emph{(strongly isotopic incomplete gluing map)}
    In Subsection~\ref{sec: normalization; subsec: incomplete gluing maps}, we consider a block $(P_1,X_1)$ with a \paif{}, and isotopic to an orbit equivalent block $(P_0, X_0)$ with a \sqt{} gluing map $\varphi_0$.
    We show (Proposition~\ref{prop: induced gluing map on orbit equivalent block}) that we can obtain an \say{almost} gluing map of $(P,X)$: an involution $\hat \varphi \colon \pP \ssm \cO_* \to \pP \ssm \cO_*$ which maps the exit boundary $\Pout$ to the entrance boundary $\Pin$ but does not extend a priori to the periodic orbits $\cO_*$.
    We will say that it is a \emph{incomplete gluing map} (Definition~\ref{def: incomplete gluing map}).
    We require that $\hat \varphi$ and $\varphi_0$ respect a relation compatible with Definition~\ref{def: strongly isotopic triples}, Item~\ref{def: strongly isotopic triples; it: lamination}, in other words we want the pattern of the intersection of the laminations to be the same.
    We will say that they are \emph{strongly isotopic} (Definition~\ref{def: strongly isotopic incomplete gluing maps}).

    \item \emph{(boundary straightening)}
    In Subsection~\ref{sec: normalization; subsec: boundary straightening}, we show that if $(P,X)$ is \bb{} which admits an affine section $\Sigma$, then there exists a block $(P_1, X_1)$ in a \minc{} of $(P,X)$ which satisfies Item~\ref{def: normalized, it: straight boundary} of Definition~\ref{def: normalized block}, in other words the boundary $\pP_1$ is diagonal in the affine coordinates of the section $\Sigma$ in the neighborhood of the periodic orbits $\cO_*$ contained in the boundary.
    This is Proposition~\ref{prop: straight block associated}.
    If the initial block $(P,X)$ satisfies Items~\ref{def: normalized, it: affine section},~\ref{def: normalized, it: multipliers} and~\ref{def: normalized, it: foliation} of Definition~\ref{def: normalized block}, it follows that the block $(P_1, X_1)$ is a normalized block.
    In a final paragraph, we explain that we can recover an incomplete gluing map on the new block $(P_1, X_1)$ which is strongly isotopic to the gluing map of the initial block $(P,X)$ (Lemma~\ref{lem: induced incomplete gluing map on isotopic blocks}).
    Such an incomplete gluing map is typically obtained by pushing the initial gluing map by the flow on the boundary of $P_1$ for an infinite time, and will in general not extend to periodic orbits.

    \item \emph{(gluing map normalization)}
    In Subsection~\ref{sec: normalization; subsec: normalization of gluing map}, we consider a normalized block $(P,X)$ equipped with an incomplete gluing map $\hat \varphi$ which maps the boundary lamination $\cL$ to a strongly transverse lamination on $\pP \ssm \cO_*$.
    We show that there exists a normalized (complete) gluing map $\varphi_1$ of $(P,X)$ strongly isotopic to $\hat \varphi$.
    This is Proposition~\ref{prop: strongly isotopic normalized gluing map}.
    We first put the foliations transverse to each other in $\pP \ssm \cO_*$ (Proposition~\ref{prop: incomplete gluing map with transverse foliation}).
    We then modify the incomplete gluing map $\hat \varphi$ by strong isotopy on the complementary of the periodic orbits $\cO_*$ contained in $\pP$, so as to control the image by $\hat \varphi$ of the foliations induced on the boundary in the neighborhood of $\cO_*$, and such that the new incomplete gluing map extends trivially on $\cO_*$.

    \item \emph{(conclusion)}
    In a last Subsection~\ref{sec: normalization; subsec: proof of triple normalization}, we finally show that, starting from a filled \bb{} $(P,X)$ provided with a \sqt{} gluing map $\varphi$, we can put together the results of the previous subsections in order to construct a normalized triple $(P_1, X_1, \varphi_1)$ strongly isotopic to $(P,X, \varphi)$ which satisfies Proposition~\ref{prop: normalization of triple}. 
\end{itemize}

\subsection{Affine section}

\label{sec: normalization; subsec: affine section}

In this section we show that we can perturb a block $(P,X)$ by isotopy in its orbit equivalence class in a block which admits an affine section in the sense of Definition~\ref{def: normalized block}, Item~\ref{def: normalized, it: affine section}.

\begin{defi}[Affine block] \label{def: affine block}
    We say that a \bb{} $(P,X)$ is \emph{affine} if it admits an affine section of the \msis{}, in the sense of Definition~\ref{def: normalized block}, Item~\ref{def: normalized, it: affine section}.
\end{defi}

\begin{prop}
\label{prop: affine section}
Let $(P,X)$ be a \bb{}.
There exists an affine \bb{} $(P_1, X_1)$, isotopic to $(P,X)$ among orbit equivalent \bbs{}.
\end{prop}

It is crucial to work near the saddle basic pieces, as it ensures the existence of a local transverse section arbitrarily small, and allows us to linearize the return map and use structural stability.
We follow closely the proof of \cite[Lemma~5.3]{beguinBuildingAnosovFlows2017} for \emph{saddle hyperbolic plugs} which works the same when restricted to the \msis{}.
We perturb the return map on a local transverse section $\Sigma$ of the \msis{} $\Lambda_s$.
If the section is small enough, there exists a perturbation $X'$ arbitrarily $\cC^1$-close to $X$ such that the first-return map on $\Sigma$ is affine in some coordinates on a neighborhood of $\Sigma \cap \Lambda_s$.
A difference is that the blocks that we consider are not structurally stable.

\subsubsection*{Weak structural stability of building blocks}

The following proposition gives a \say{weakened} structural stability result for \bbs{}.

\begin{prop}
\label{prop: weak structural stability}
Let $(P,X)$ be a \bb{} and $(\tilde P, \tilde X)$ be a \minc{} (Definition~\ref{def: minc}).
If $\tilde Y$ is a vector field on $\tilde P$ close enough to $\tilde X$ in the $\cC^1$ topology, then there exists a block $(Q,Y) = (Q, \res{\tilde Y}{Q})$ contained in $\tilde P$ which is isotopic to $(P,X)$ among the orbit equivalent \bbs{}.
\end{prop}

\begin{proof}
According to Lemma~\ref{lem: from bb to bby}, up to shrink the \minc{} $(\tilde P, \tilde X)$, we can embed $(\tilde P, \tilde X)$ into a \bb{} $(\hat P,\hat X)$ whose boundary $\partial \hat P$ is transverse to the \vf{} $\hat X$.
A small enough perturbation $\tilde Y$ of $\tilde X$ on $\tilde P$ extends to $\hat P$ into a $\epsilon$-perturbation of $\hat X$, denoted $\hat Y$.
By structural stability (Lemma~\ref{lem: bby block structurally stable}), there exists a \homeo{} $h \colon \hat P \to \hat P$ isotopic to the identity which realizes an orbit equivalence between the flow of $\hat X$ and the flow of $\hat Y$.
Moreover, $h$ is $K. \epsilon$ close to the identity in the $\cC^0$ topology on $\hat P$, with $K$ a uniform constant.
Then $h$ induces an orbit equivalence between $(P,X)$ and $(Q, Y)$ where $Q := h(P)$ is contained in a $K. \epsilon$-neighborhood of $P$ in $\hat P$, and $Y := \res{\hat Y}{h(P)}$.
Let $\cV_0$ be a disjoint union of small tubular neighborhoods of periodic orbits of $X$ contained in $\pP$, contained in linearizing neighborhoods.
Then $\pP \ssm \cV_0$ is a compact surface uniformly transverse to $\hat X$.
For $\epsilon$ small enough, $h(\pP \ssm \cV_0)$ is a compact surface uniformly transverse to $\hat Y$ in a $\epsilon$-neighborhood of $\pP \ssm \cV_0$.
Up to push the surface along the orbits of $\hat X$ and modify $h$, we can assume that $h(\pP \ssm \cV_0) = \pP \ssm \cV_0$.
Then for $\epsilon$ small enough, $Q = h(P)$ is contained in the \minc{} $\tilde P$ of $(P,X)$.
Moreover (up to do an isotopy along the orbits of the flow), $Q$ is a surface \qt{} to $Y$, so $(Q, Y)$ is an orbit equivalent \bb{} to $(P,X)$.
Since $h$ is isotopic to the identity, the \vfs{} $\tilde X$ and $\tilde Y$ are isotopic.
Moreover, the boundary of $P$ and $Q$ are contained in the same \minc{}, so according to Proposition~\ref{prop: isotopy vs orbit eq}, Item~\ref{prop: isotopy vs orbit eq; it: same minc implies isotopy}, the blocks $(P,X)$ and $(h(P), \res{\tilde X}{h(P)})$ are isotopic among \bbs{}.
We deduce that $(P,X)$ and $(h(P), \res{\tilde Y}{h(P)}) = (h(P), Y)$ are isotopic among the orbit equivalent \bbs{}.
The block $(h(P),Y)$ satisfies Proposition~\ref{prop: weak structural stability}.
\end{proof}

In general one cannot have $h(P)=P$, i.e., $h$ does not preserve the boundary of $P$, and if $Y$ is a small perturbation of $X$ on $P$, the pair $(P,Y)$ is not always a \bb{} (Remark~\ref{rmk: block not structurally stable}).

\subsubsection*{Markov partition} 

We will need the notion of Markov partition.
We give the definition following \cite[Section~2.1]{beguinFlotsSmaleDimension2002}.
Let $\Sigma$ be a transverse local section of $\Lambda_s$ in a \bb{} $(\tilde P, \tilde X)$, $f$ the first return map of the flow of $\tilde X$ on $\Sigma$ (Definition~\ref{def: local section}).
\begin{itemize} 
    \item A \emph{rectangle} $R$ on $\Sigma$ is an embedding $h$ of $I \times J$ into $\Sigma$, where $I$ and $J$ are two non-trivial closed segments of $\R$.
    The \emph{horizontal boundary} of $R$, denoted $\partial^s R$ is the image of $I \times \partial J$ and the \emph{vertical boundary} of $R$, denoted $\partial^u R$ is the image of $\partial I \times J$.
    \item We will always suppose that a rectangle $R$ is equipped with a pair of trivial transversal $\cC^{0,1}$-foliation $(\zeta^s,\zeta^u)$, containing the trace of $(\cW^s, \cW^u)$, invariant by $f$, and such that the horizontal boundary of $R$ is formed by two leaves of $\zeta^s$ and the vertical boundary of $R$ is formed by two leaves of $\zeta^u$.
    \item A horizontal subrectangle $H$ of $R$ is a rectangle contained in $R$ such that $\partial^u H \subset \partial^u R$ and $\partial^s H$ is formed by two leaves of $\zeta^s$.
    A \emph{vertical subrectangle} $V$ of $R$ is a rectangle contained in $R$ such that $\partial^s V \subset \partial^s R$ and $\partial^u V$ is formed by two leaves of $\zeta^u$. \end{itemize}

\begin{defi}[{\cite[Definition~2.1]{beguinFlotsSmaleDimension2002}}] \label{def: markov partition}
A \emph{Markov partition} of $\Lambda_s$ is the data $(\Sigma, R = \bigcup_i R_i)$ of a transverse local section $\Sigma$ of $\Lambda_s$ and a finite collection $R_1, \dots, R_n$ of foliated rectangles on $\Sigma$ such that:
\begin{enumerate}
    \item The union of rectangles $R_i$ covers the set $\Lambda \cap \Sigma$
    \item For any $i$, the first return $f(R_i)$ of $R_i$ on $\Sigma$ is well defined.
    \item For any $i,j$ the intersection $f(R_i) \cap R_j$ has a finite number of components. 
    Each of these components is both a horizontal subrectangle of $f(R_i)$ and a vertical subrectangle of $R_j$. 
    Moreover, we require that this intersection be disjoint from the vertical boundary of $R_j$ and the horizontal boundary of $f(R_i)$.
    \item For any $x \in R_i \cap f^\inv (R_j)$, any vector $v$ tangent at $x$ to the foliation $\zeta^s$ is uniformly contracted by the differential $d_xf$ and any vector $w$ tangent at $x$ to the foliation $\zeta^u$ is uniformly expanded by the differential $d_xf$ for a metric on $\Sigma$.
\end{enumerate}
\end{defi}

If $\cR = (\Sigma, R)$ is a Markov partition, we will always assume that we have renormalized $\tilde X$ so that the first return time of $R$ on $\Sigma$ is equal to $1$.
We call \emph{suspension \nbh{} of $R$} the union $\cC = \bigcup_{t=0}^1 \tilde X^t(R)$ of the orbit segments between a point of $R$ and its first return on $\Sigma$.

\begin{lem} \label{lem: existence of markov partition for saddle block}
Let $(P,X)$ be a block, $(\tilde P, \tilde X)$ a \minc{}, and $\cO_*$ the collection of periodic orbits of $X$ contained in $\pP$.
Then there exists a Markov partition $\cR = (\Sigma, R)$ (in the sense of Definition~\ref{def: markov partition}) of $\Lambda_s$ such that $\Sigma$ intersects each orbit $\cO_i \in \cO_*$ at a single point $p_i$ and the intersection of an orbit of $\tilde X$ with the suspension \nbh{} $\cC$ of $R$ is connected.
\end{lem}

\begin{proof}[Proof of Lemma~\ref{lem: existence of markov partition for saddle block}]
According to Lemma~\ref{lem: from bb to bby}, up to shrink the \minc{}, we can embed $(P,X)$ into a \bb{} $(U,Y)$ whose boundary $\partial U$ is transverse to the vector field $Y$ and which contains the \minc{} $(\tilde P, \tilde X)$.
Then $\Lambda_s$ is a locally maximal saddle hyperbolic invariant compact set of index $(1,1)$ for the flow of $Y$ and $U$ contains a filtrating neighborhood of $\Lambda_s$ for $Y$ (Claim~\ref{claim: smale extension}).
The existence of a Markov partition $\cR = (\Sigma, R)$ of $\Lambda_s$ (in the sense of Definition~\ref{def: markov partition}) for the flow of $Y$ in $U$ and such that the intersection of the orbit segments with the suspension \nbh{} $\cC$ of $\cR$ is connected is shown in \cite[Proposition~2.8]{beguinFlotsSmaleDimension2002}.
Up to shrink the section  $\Sigma$, we can suppose that it is contained in the \minc{} $\tilde P$.
Every orbit $\cO_i \in \cO_*$ belongs to $\Lambda_s$ (Lemma~\ref{lem: boundary leaf and periodic orbit are saddle}), so up to cutting $\Sigma$ and removing some \ccs{} we can assume that $\cO_i$ intersects $\Sigma$ at a single point $p_i$.
\end{proof}

\begin{proof}[Proof of Proposition~\ref{prop: affine section}]
Let $(P,X)$ be a block and let $\cO_*$ be the collection of periodic orbits of $X$ contained in $\pP$.
Up to do a small perturbation in $\cC^1$ norm and use Proposition~\ref{prop: weak structural stability}, we can suppose that the \vf{} $X$ is of class $\cC^2$.
Let $\cR = (R, \Sigma)$ be a Markov partition of $\Lambda_s$ on a transverse local section $\Sigma$ given by Lemma~\ref{lem: existence of markov partition for saddle block}.
Denote $f \colon R \to \Sigma$ the first return map and $(\zeta^s, \zeta^u)$ the pair of $f$-invariant transverse foliations on $R$.
Any orbit $\cO_i \in \cO_*$ intersects $\Sigma$ at a single point $p_i$.
We will perturb the \vf{} $X$ into a \vf{} $Y$ arbitrarily $\cC^1$-close and such that $\Sigma$ is an affine section for $Y$.
Let $R^n$ be the \emph{order $n$ refinement} of the partition $R$, defined by
$$ R^n := \underset{k=-n}{\overset{n}{\bigcap}} f^k R \subset \Sigma. $$
Let $\epsilon_n$ be the diameter of the partition $R_n$, in other words the maximum of the diameter of the \ccs{} of $R^n$.
The diameter $\epsilon_n$ is bounded above (up to a multiplicative constant independent of $n$) by the maximum of the norm of $df^n$ in the stable direction (tangent to $\zeta^s$) and of $df^{-n}$ in the unstable direction (tangent to $\zeta^u$) on the rectangles, so $\epsilon_n \to 0$ when $n \to + \infty$.
Let $R^n_1, \dots, R^n_{k_n}$ be the \ccs{} of $R^n$.
For any $p \in R^n_i$, there exists a $\cC^1$-identification of the tangent space $T_p \Sigma$ with the rectangle $R^n_i$.
More precisely, there exists $\varphi_p \colon R^n_i \to T_p\Sigma$, such that
\begin{itemize}
    \item $\varphi_p$ is a \diff{} on its image $V := \varphi_p (R_i^n)$ which is a $\epsilon_n$-neighbor\-hood of the origin $O$ in $T_p \Sigma$, and maps $p$ to $O$;
    \item there exists a $\cC^0$-identification $\chi_{i,p} = (x, y) \colon T_p \Sigma \to \R^2$ such that the foliation $(\varphi_p)_* \zeta^u$ coincides with the vertical foliation $x = \cst$ and $(\varphi_p)_* \zeta^s$ coincides with the vertical foliation $y = \cst$.
\end{itemize}

For any $n \geq 1$, the intersection $R_{(i,j)}^n := R_i^n \cap f^\inv(R_j^n)$ is connected, and for any pair $(i,j)$ such that this intersection is nonempty we denote the corresponding restriction of $f$ by
$$f_{n, (i,j)} \colon R_i^n \cap f^\inv(R_j^n) \to R_j^n.$$
Let $\{ p_{n,(i,j)} \}_{(i,j)}$ be a set of points of $\Sigma \cap \Lambda$ such that $p_{n,(i,j)} \in R_{(i,j)}^n$, and such that if $R_{i,i}^n$ is the rectangle intersecting the orbit $\cO_i \in \cO_*$, then $p_{n,(i,i)} = p_i \in \cO_i$ is the intersection point of $\cO_i$ with $\Sigma$.
Let $\hat f_{n, (i,j)}$ be the conjugate of $d_{p_{n, (i,j)}} f$ by $\varphi_{p_{n,(i,j)}}$ and $\varphi_{f(p_{n,(i,j)})}$, in other words
\begin{equation}
    \label{eq: differential conjugate}
    \hat f_{n, (i,j)} := \varphi_{f(p_{n,(i,j)})}^\inv \circ d_{p_{n,(i,j)}} f \circ \varphi_{p_{n,(i,j)}} \colon R_{(i,j)}^n \to R_j^n.
\end{equation}
By construction, there exist coordinates $\chi_i = (x,y) \in \R^2$ on $R_i^n$ and $\chi_j = (x,y) \in \R^2$ on $R_j^n$ such that the map $\hat f_{n, (i,j)}$ is affine in coordinates $\chi_i$ and $\chi_j$, and the foliations $(\zeta^s, \zeta^u)$ coincide with the horizontal and vertical foliations in these coordinates.

\begin{claim} \label{claim: perturb f on sigma}
There exists $n \in \N^*$ and $\hat f = \hat f_n \colon R \to \Sigma$ which is a \diff{} onto its image, such that
\begin{enumerate}
    \item $\hat f$ coincides with $\hat f_{n, (i,j)}$ on $R_{(i,j)}^n$ for each pair $(i,j)$;
    \item $\hat f$ coincides with $f$ on a neighborhood of the boundary of $R$;
    \item $\hat f$ is arbitrarily $\cC^1$-close to $f$ for $n$ large enough.
\end{enumerate}
\end{claim}

\begin{proof}
Since $f$ is $\cC^2$, according to Equation~\eqref{eq: differential conjugate}, the map $\hat f_{n, (i,j)}$ differs from $f_{n, (i,j)}$ by a term uniformly proportional to $\epsilon_n^2$, i.e., there exists a constant $\cst >0$ such that for any $p \in R_{(i,j)}^n$,
$$ \vert \hat f_{n, (i,j)}(p) - f_{n, (i,j)}(p) \vert \leq \cst \epsilon_n^2. $$
For rectangles of diameter $\epsilon_n$, if the distance between the rectangles $R^n_i$ in $\Sigma$ is a $O(\epsilon_n)$, we can set $\hat f = \hat f_{n, (i,j)}$ on $R_{(i,j)}^n$ and extend on $R$ so as to reconnect this map with the map $f$ on the boundary of $R$, remaining $\epsilon_n$-close to $f$ in $\cC^1$ topology.
It is therefore sufficient to check that the distance in $\Sigma$ between the rectangles of a refinement $R^n$ does not decrease too fast.
Since the expansion and contraction of the differential of $f$ is uniform on a neighborhood of $\Lambda_s \cap \Sigma$, the distance in $\Sigma$ between the rectangles is also of order $\epsilon_n$.
\end{proof}

The suspension of $\hat f$ on $R$ defines a flow generated by a $\cC^1$-vector field $Y$ on a neighborhood of $\Lambda$ in $\tilde P$, which coincides with $\tilde X$ on a neighborhood of the boundary, and with first return $\hat f$ on $\Sigma$.
Let $\tilde Y$ denote the extension of $Y$ by $\tilde X$ on $\tilde P$.
Then $\tilde Y$ is arbitrarily $\cC^1$-close to $\tilde X$.
By Proposition~\ref{prop: weak structural stability}, there exists a block $(Q,Y)$ isotopic to $(P,X)$ among the orbit equivalent blocks in $\tilde P$, where $Y = \res{\tilde Y}{Q}$.
By construction, the flow of $Y$ admits an affine section $\Sigma$ of its \msis{} $\Lambda_s$.
Indeed, there exists a cover of $\Sigma \cap \Lambda_s$, given by the collection of rectangles $R^n = \{R^n_1, \dots, R^n_{k_n} \}$ of a large enough refinement of a Markov partition $\cR$ on $\Sigma$, and coordinates $\chi_i$ on $R^n_i$ such that the map $\hat f$ restricted to $R_{(i,j)}^n = R^n_i \cap f^\inv (R^n_j)$ is affine in coordinates $\chi_i$ and $\chi_j$.
Hence the block $(Q,Y)$ satisfies Proposition~\ref{prop: affine section}.
\end{proof}

\begin{rmk} \label{rmk: affine section and equal multipliers} Note that one can apply Proposition~\ref{prop: affine section} while preserving the multipliers of the periodic orbits contained in the boundary, which will be useful in the next subsection.
Indeed, in the neighbourhood of a periodic orbit $\cO_i \in \cO_*$ we have replaced the first return of $X$ on $\Sigma$ by its differential at the point $p_i = \cO_i \cap \Sigma$, which preserves the multipliers of the periodic orbits of the new flow in the boundary.
More precisely, using the notations from the proof, let $\cO_i$ be a periodic orbit of $X$ in $\pP$, and let $R_i^n$ be the rectangle of $R^n$ intersecting $\cO_i$. Let $R_{(i,i)}^n = R^n_i \cap f^\inv(R_i^n)$.
According to Claim~\ref{claim: perturb f on sigma} the return map $\hat f$ of the flow of $\tilde Y$ on $R_{(i,i)}^n$ coincides with the $\cC^1$-conjugate of the differential of $f$ in $p_i = \cO_i \cap \Sigma$, where $f$ is the return map of the flow of $\tilde X$ on $\Sigma$.
It follows that $p_i$ is a fixed point of $\hat f$ on $\Sigma$, contained in the boundary of $P$, and its multipliers are those of $d_p f$. The orbit equivalence between $X$ and $Y$ fixes the periodic orbits $\cO_i$, and the multipliers of these orbits for the flow of $X$ and the flow of $Y$ are the same. \end{rmk}

\subsection{Multipliers of periodic orbits}
\label{sec: normalization; subsec: multipliers}

In this subsection, we perturb the block by isotopy in its orbit equivalence class so that the multipliers of the periodic orbits contained in the boundary are $\{ \frac{1}{2}, 2\}$.
More precisely, we show

\begin{prop} \label{prop: multipliers} 
Let $(P,X)$ be an affine \bb{}.
There exists an affine block $(P_1, X_1)$ isotopic to $(P,X)$ among the orbit equivalent \bbs{}, such that the multipliers of the periodic orbits of $X_1$ contained in $\pP_1$ are equal to $\{ \frac{1}{2}, 2\}$. \end{prop}

\begin{proof}
Let $\cO_*$ be the collection of periodic orbits of $X$ contained in $\pP$.
Let $\Sigma$ be an affine section of $(P,X)$, and $\cR = (\Sigma, R)$ a Markov partition of $\Lambda_s$ in the sense of Definition~\ref{def: markov partition} which satisfies the assumptions of Lemma~\ref{lem: existence of markov partition for saddle block}.
Up to shrink $R$, we assume that the Markov partition is contained in a domain of affine coordinates on $\Sigma$ (Definition~\ref{def: normalized block}, Item~\ref{def: normalized, it: affine section}).
The map $f \colon R \to \Sigma$ is well defined and preserves a pair of transverse foliations $(\zeta^s, \zeta^u)$ on $R$ which are the horizontal and vertical foliations in affine coordinates. 
These foliations are now of class $\cC^1$ on $R$ and contains the trace of the stable and unstable manifold of $\Lambda$.
Each orbit $\cO_i$ is associated with a single fixed point $p_i$ of $f$ and positive multipliers according to Claim~\ref{claim: boundary quadrant and multipliers}. 
We want the eigenvalues of the differential $D_{p_i} f$ to be equal to $\{ \frac{1}{2}, 2\}$. 
We will therefore continuously modify the derivative at the points $p_i$, for $i=1, \dots, n$.

Let $R/\zeta^s$ be the quotient of the rectangles by the foliation $\zeta^s$.
The map $f \colon R \cap f^\inv (R) \to R \cap f(R)$ induces on $R/\zeta^s$ a map $f^u \colon \bigcup_{k,l} I_k^l \to \bigcup_k I_k$, where $I_1, \dots, I_N$ are disjoint closed ordered intervals of $\R$, and $I_k^1, \dots, I_k^{n_k}$ are closed disjoint ordered strict subintervals of $I_k$, and the restriction of $f^u$ to $I_k^l$ is a uniformly expanding affine map in some $I_i$.
Let $p^u_i \in I_{k_i}^{l_i}$ be the projection of the fixed point $p_i \in \cO_i$ of $f$.

\begin{claim} \label{claim: change multipliers dim 1} For any $\lambda >1$, there exists a diffeomorphism $g^u = g^u_\lambda \colon \bigcup_{k,l} I_k^l \to \bigcup_k I_k$ such that 
\begin{enumerate} 
    \item $g^u = f^u$ on $I_k^l$ for $(k,l) \neq (k_i, l_i)$, 
    \label{claim: change mult dim1, it: egal F} 
    \item the restriction $g^u \colon I_{k_i}^{l_i} \to I_{k_i}$ is a uniformly expanding \diff{} which coincides with $f^u$ near the boundary of $I_{k_i}^{l_i}$, 
    \label{claim: change mult dim1, it: on I0} 
    \item $p_i$ is a fixed point of $g^u$ and $g^u{}'(p) = \lambda$, 
    \label{claim: change mult dim1, it: derivative in p} 
    \item $g^u_\lambda$ depends continuously on $\lambda$ and is conjugate to $f^u$ by a \homeo{} \mb{$H \colon \bigcup_k I_k \to \bigcup_k I_k$}, equal to the identity outside $I_{k_i}^{l_i}$ and on the boundary of~$I_{k_i}^{l_i}$. 
    \label{claim: change mult dim1, it: continuity in lambda and conjugation} 
\end{enumerate} 
\end{claim}

\begin{proof}
To simplify the notation, let $(k_1, l_1) = (1,1)$, and let $p=p_1$. The set $K := \bigcap_{n \in \Z} (f^u)^n \big(\bigcup_k I_k\big)$ is a Cantor set disjoint from the boundary of $I_k^l$, and $\vert I_1 \vert > \vert I_1^1 \vert$.
Let $V_\lambda$ be a neighborhood of $p$ in $I_1^1$, whose boundary is disjoint from the boundary of $I_1^1$, small enough so that $ \vert I_1 \vert > \vert I_1^1 \vert - (\lambda-1) \vert V_\lambda \vert $.
Then $g^u(x) = \lambda (x - p) +p$ on $V_\lambda$.
Then $\vert I_1 \ssm g^u(V_\lambda) \vert > \vert I_1^1 \ssm V_\lambda \vert$, so we can extend $g^u$ on $I_1^1$ into a uniformly expanding \diff{} which we still denote $g^u$, and which coincides with $f^u$ on the boundary of $I_1^1$.
We can choose the neighborhood $V_\lambda$ which depends continuously on $\lambda$, and extend $g^u = g^u_\lambda$ on $I_1^1$ continuously into $\lambda$.
Then $g^u$ and $f^u$ are two expanding \diffs{} on $I_1^1 \to I_1$, which coincide on the boundary of $I_1^1$. Let $H \colon \bigcup_k I_k \ssm I_1^1 \to \bigcup_k I_k \ssm I_1^1$ be a map equal to the identity.
There is a unique way of extending $H$ on $I_1^1$ so that we have a conjugation $H$ between $f^u$ and $g^u$. 
We can repeat the same procedure on each fixed point $p^u_i$, by noticing that the map $g^u$ obtained still satisfies the hypotheses of the fact, and each $p^u_i$ is still a fixed point of~$g^u$.
\end{proof}

Let $g^u$ be the map given by Claim~\ref{claim: change multipliers dim 1} for $\lambda=2$.
In the same way, when quotienting $R$ by the foliation $\zeta^u$, the map $f^\inv$ induces an affine map $f^s \colon \bigcup_{k,l} J_k^l \to \bigcup_k J_k $ which satisfies the hypotheses of Claim~\ref{claim: change multipliers dim 1}.
Let $p^s_i$ be the projection in $R/\zeta^u$ of the fixed point $p_i \in \cO_i$ of $f :R \to \Sigma$ on the orbit $\cO_i \in \cO_*$.
There exists a map $g^s \colon \bigcup_{k,l} J_k^l \to \bigcup_k J_k$ expanding on each $J_k^l$, with $g^s(J_k^l) = f^s(J_k^l)$ for any pair $(k,l)$, and such that each $p^s_i \in J_{k_i}^{l_i}$ is a fixed point of $g^s$ with $(g^s)'(p^s_i) = 2$, and $g^s$ is conjugate to $f^s$.
Since the foliations $(\zeta^s, \zeta^u)$ are of class $\cC^1$, the product coordinates $R \to R/\zeta^s \times R/\zeta^u$ are of class $\cC^1$ and the product map $\l( (g^s)^\inv, g^u \r)$ induces a \diff{} $g$ of $R \cap f^\inv (R)$ with image in $R \cap f(R)$.

According to Claim~\ref{claim: change multipliers dim 1}, Item~\ref{claim: change mult dim1, it: continuity in lambda and conjugation}, the map $g$ is isotopic to $f$ via a continuous family of \diffs{} $f_t \colon R \cap f^\inv(R) \to R \cap f(R)$, with $f_0 = f$ and $f_1 = g$, satisfying:

\begin{itemize}[--] 
    \item each point $p_i \in \Sigma \cap \cO_i$ is a fixed point of $f_t$, and the multipliers of the differential of $d_{p_i} f_1 = d_{p_i} g$ are $\{\frac{1}{2}, 2\}$;\\[-0.5em]
    \item the map $f_t \colon R \cap f^\inv(R) \to R \cap f(R)$ is conjugate to $f$ by a \homeo{} $h_t: R \to R$, equal to the identity on the edge of $R$.
\end{itemize}

Let $\cC = \bigcup_{t=0}^1 \tilde X^t R$ be the suspension \nbh{} of the Markov partition associated with $\tilde X$ in $\tilde P$.
It is a neighbourhood of $\Lambda$, whose intersection with any orbit of $\tilde X$ is connected (Lemma~\ref{lem: existence of markov partition for saddle block}), homeomorphic to the suspension manifold (with edge and corner) $R_f := R \times [0,1] / (x,1) \sim (f(x), 0)$ by a \homeo{} which conjugates the flow of $\tilde X$ with the suspension flow induced by $\partial_s$, where $s$ is the coordinate on $[0,1]$.
Moreover, the suspension flow on $R_f$ is topologically conjugate to the suspension flow on $R_{f_t} := R \times [0,1] / (x,1) \sim (f_t(x), 0)$.
We deduce that $R_{f_t}$ is naturally embedded in $\tilde P$ via a homeomorphism which maps $R_{f_t}$ on $\cC$ and which maps the orbits of the flow of $\partial_t$ of $R_{f_t}$ on the orbits of the flow of a \vf{} $\tilde X_t$ on $\cC$ whose first return on $R \subset \cC$ is equal to $f_t$, and which coincides with $\tilde X$ on the boundary of $\cC$. 
The \vfs{} are orbit equivalent via a \homeo{} $h_t \colon \cC \to \cC$, equal to the identity on the boundary of $\cC$.

Extend $\tilde X_t$ by $\tilde X$ on the rest of the manifold, and we obtain a continuous family (for the $\cC^1$-norm) of vector fields on $\tilde P$, which we still denote $\tilde X_t$.
As any orbit which exits $\cC$ never enters it again, the orbit equivalence $h_t$ on $\cC$ extends by the identity into an orbit equivalence $h_t \colon \tilde P \to \tilde P$ between $\tilde X_t$ and $\tilde X$, with $h_0 = \Id$.
Let $P_t = h_t^\inv (P)$. 
Up to  an isotopy along the flow of $\tilde X_t$, the surface $\partial P_t$ is \qt{} to $\tilde X_t$.
The family $(P_t, X_t)$ is an isotopy of orbit equivalent \bbs{}, where $X_t = \res{\tilde X_t}{P_t}$.

Let $\cO_{*,t}$ be the collection of periodic orbits of $X_t$ contained in $\pP_t$. 
Each orbit $\cO_{i,t} \in \cO_{*,t}$ intersects the section $\Sigma$ at a single point $p_{i,t} = p_i$.
By the properties of $f_1 = g$, each $p_i$ is a fixed point of $p$ with multipliers $\{\frac{1}{2}, 2\}$. It follows that the periodic orbits of $X_1$ contained in $\pP_1$ have multipliers $\{\frac{1}{2}, 2\}$.

The section $\Sigma$ is a priori no longer an affine section for the flow of $X_1$, because we have modified the first return map on the neighbourhood of the maximal invariant set.
We can apply Proposition~\ref{prop: affine section} and Remark~\ref{rmk: affine section and equal multipliers} to the block $(P_1, X_1)$ to obtain a block $(P'_1, X'_1)$, isotopic to $(P_1, X_1)$ among the orbit equivalent blocks, which admits an affine section and whose periodic orbits contained in the boundary have multipliers $\{\frac{1}{2}, 2\}$.
This block satisfies Proposition~\ref{prop: multipliers}.
\end{proof}

\subsection{Pair of affine invariant foliation}
\label{sec: normalization; subsec: pair of affine invariant foliation}

In this subsection, we show that we can foliate a filled affine block $(P,X)$ with a \paif{} $(\cG^s, \cG^u)$ (Definition~\ref{def: normalized block}, Item~\ref{def: normalized, it: foliation}).

\begin{prop}
\label{prop: existence paif}
Let $(P,X)$ be a filled affine \bb{}, and $\Sigma$ an affine section of $(P,X)$.
There exists a \paif{}
$(\cG^s, \cG^u)$, in the sense of Definition~\ref{def: normalized block}, Item~\ref{def: normalized, it: foliation}.
\end{prop}

We repeat here the idea of the proof of \cite[Lemma~5.6,~Proposition~5.7]{beguinBuildingAnosovFlows2017}.
The key points of the proof are the following.
\begin{itemize}[leftmargin=*]
    \item We can push the horizontal and vertical foliations induced by the \acs{} on the section $\Sigma$, in order to obtain a pair of transverse $X^t$-invariant foliations on an invariant neighborhood of $\Lambda$ in $P$, which extends the stable and unstable manifolds of $\Lambda$ (Lemma~\ref{lem: paif local}).
    These foliations will naturally satisfy Item~\ref{def: normalized, it: foliation, it: affine} of Definition~\ref{def: normalized block}.
    \item We then show that there is a (unique) way to complete these foliations into a \paif{} on $P$. We use the fact that the boundary lamination is filling.
\end{itemize}

\subsubsection*{Affine invariant foliations on a \nbh{} of $\Lambda$}

\begin{lem} \label{lem: paif local}
    Let $(P,X)$ be an affine \bb{} equipped with an affine section $\Sigma$. There exist an invariant \nbh{} $\cU$ of $\Lambda$, endowed with a pair of smooth transverse invariant 2-dimensional foliations $(\hat \cG^s, \hat \cG^u)$, containing \linebreak[4]$(\cW^s, \cW^u)$ as sublaminations, and such that the trace on $\Sigma$ are the vertical and horizontal foliations in the affine coordinate. 
\end{lem}

\begin{proof}
   Let $(\zeta^s, \zeta^u)$ be the pair of transverse foliations, invariant by the first return map $f$, which coincide with the horizontal and vertical foliations on $\Sigma$ in the affine coordinate system on a neighborhood of $\Lambda_s \cap \Sigma$ in $\Sigma$.
According to Lemma~\ref{lem: existence of markov partition for saddle block}, there exists a Markov partition $\cR = (\Sigma, R)$ on $\Sigma$ in a \minc{} $(\tilde P, \tilde X)$ such that the intersection of the orbits of the flow of $\tilde X$ with the suspension \nbh{} $\cC = \bigcup_{t=0}^1 \tilde X^t (R)$ associated with the partition is connected, and such that (up to take a refinement of the partition) $R$ is contained in the domain of affine coordinates of class $\cC^1$ on $\Sigma$.
Recall that the map $f \colon R \to \Sigma$ is well defined and coincides with time 1 of the flow $\tilde X$ restricted to $R$.
Consider the union $\cU$ of the flow-saturated set of $\cC$ with the basins of $\gA$ and $\gR$.
It is an invariant \nbh{} of $\Lambda$, hence of $\cW^s\cup \cW^u$.
By pushing the horizontal and vertical foliations in affine coordinates on $R$ by the flow of $X$, containing the trace of the stable and unstable lamination, we obtain two foliations $\hat \cG^s$ and $\hat \cG^u$ of dimension $2$ and class $\cC^1$, transverse to each other, $X^t$-invariant, which extend $\cW^s$ and $\cW^u$ on $\cU$.
By the $\lambda$-lemma, they continuously extend $\cW^s$ and $\cW^u$ on $\cU$.
They satisfy Lemma~\ref{lem: paif local}.
\end{proof}

\subsubsection*{Adapted neighborhoods} 

Let $(P,X)$ a \bb{} satisfying Lemma~\ref{lem: paif local}.
The trace $U = \cU \cap \pP$ is a \nbh{} of the boundary lamination $\cL$ foliated by the trace of $\hat \cG^u$ and $\hat \cG^s$ on $\pP$.
Denote $(\hat \cG^{s, \iin}, \hat \cG^{u, \iin})$ the trace of $(\hat \cG^s, \hat \cG^u)$ on $\Pin$ and $(\hat \cG^{s, \out}, \hat \cG^{u, \out})$ the trace on $\Pout$.
We adapt the notions of \cite[Sections~5.3]{beguinBuildingAnosovFlows2017}.
The definitions of Item~\ref{def: adapted nbh; it: of L^in} and~\ref{def: adapted nbh; it: of L^out} are not used in the subsection but will be useful later.

\begin{defi}[Adapted \nbhs{}] \label{def: adapted nbh}
\mb{}
\begin{enumerate}
    \item \label{def: adapted nbh; it: of L}
    Let $U$ be a \nbh{} of $\cL$ on $\pP$.
    We say that a \nbh{} $U$ of $\cL$ in $\pP$ is an \emph{adapted \nbh{} of $\cL$} if the closure of \ccs{} $\pP \ssm \cU$ are squares $C$ diffeomorphic to $[0,1]^2$ with two sides on leaves of $\cG^{s, *}$ and two sides on leaves of $\hat \cG^{u, *}$ with $* = \iin$ or $\out$ depending on whether $C$ is contained in $\Pin$ or $\Pout$.
    We say that $C$ is an \emph{adapted square}.
    
    \item \label{def: adapted nbh; it: of L^in}
    Let $U^\iin$ be a \nbh{} of $\cL^\iin$ on $\Pin$.
    We say that $U^\iin$ is an \emph{adapted \nbh{} of $\cL^\iin$} if the closure of \ccs{} $\Pin \ssm U^\iin$ are either squares diffeomorphic to $[0,1]^2$, band $\R \times [0,1]$ or half-band $\R_+ \times [0,1]$, with horizontal sides on leaves of $\cG^{s, \iin}$ and vertical sides (when well defined) on leaves of $\hat \cG^{u, \iin}$, or ends (when well defined) accumulating on $\partial \Pin = \cO_*$.
    We say that $C$ is an \emph{adapted square, band or half-band}.
    
    \item \label{def: adapted nbh; it: of L^out}
    Let $U^\out$ be a \nbh{} of $\cL^\out$ on $\Pout$.
    We say that $U^\out$ is an \emph{adapted \nbh{} of $\cL^\out$} if the closure of \ccs{} $\Pout \ssm U^\out$ are either squares diffeomorphic to $[0,1]^2$, band diffeomorphic $\R \times [0,1]$ or half-band $\R_+ \times [0,1]$, with horizontal sides on leaves of $\hat \cG^{u, \out}$ and vertical sides (when well defined) on leaves of $\cG^{s, \out}$, or ends (when well defined) accumulating on $\partial \Pout = \cO_*$.
    We say that $C$ is an \emph{adapted square, band or half-band.}

\end{enumerate}
\end{defi}

\begin{figure}[htb]
    \centering
    \vspace*{-1em}
    \hspace*{-1.5em}
    \includegraphics[height=.283\textheight]{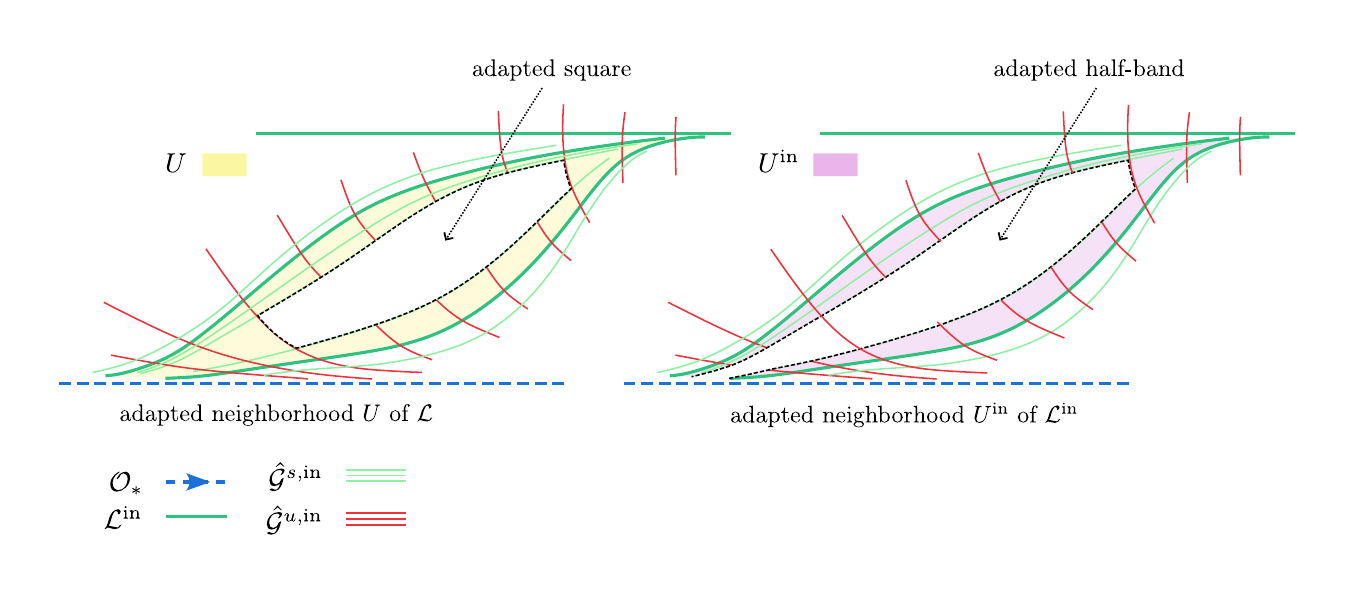}
    \hspace*{-1em}
      \vspace*{-1em}
    \caption{Adapted square of an adapted \nbh{} $U$ of $\cL$ (left) and adapted half-band of an adapted \nbh{} $U^\iin$ of $\cL^\iin$ (right)}
    \label{fig: adapted square half band}
\end{figure}

In order to point out the difference between adapted \nbhs{} of $\cL$ and of $\cL^\iin$ we pictured the \nbhs{} inside a \cc{} $B$ of $\Pin \ssm \cL^\iin$ in the case of Item~\ref{def: adapted nbh; it: of L} and in the case of Item~\ref{def: adapted nbh; it: of L^in} in Figure~\ref{fig: adapted square half band}.

One can easily prove the following lemma,  using the fact that $\cL$ is a filling \qms{} lamination.

\begin{lem} \label{lem: basis of adapted nbh}
    Adapted \nbhs{} of $\cL$ (resp. $\cL^\iin$, resp. $\cL^\out$) are a basis of \nbhs{} of $\cL$ (resp. $\cL^\iin$, resp. $\cL^\out$) in $\pP$ (resp. $\Pin$, resp. $\Pout$).
\end{lem}

\subsubsection*{Global affine invariant foliations}

\begin{proof}[Proof of Proposition~\ref{prop: existence paif}]
Let $(P,X)$ a filled affine \bb{}, and suppose that $(P,X)$ satisfied Lemma~\ref{lem: paif local}.
Up to shrink $\cU$, we can suppose according to Lemma~\ref{lem: basis of adapted nbh} that $\partial \cU \cap \pP$ is an adapted \nbh{} of $\cL$ (Definition~\ref{def: adapted nbh}).
The complementary $\Pin \ssm \cU$ is a finite number of adapted squares $C_i$ in $\Pin$.
We complete the pair of foliations trivially on each square $C_i$ of the complementary $\Pin \ssm \cU$ by the horizontal and vertical foliations in the horizontal and vertical coordinate.
We obtain a pair of transverse smooth foliations $\cG^{s, \iin}$ and $\cG^{u, \iin}$ on $\Pin$, respectively extending $\hat \cG^{s, \iin}$ and $\hat \cG^{u, \iin}$.
Any point of a component $C_i$ is uniformly far from the entrance lamination $\cL^\iin$.
Any future orbit of a point of $C_i$ intersects $\Pout$ at a single point, and the orbit segment is uniformly bounded in length.
There exists $X'$ a re-normalization of $X$ such that the $X'$-orbit of each point of $C_i$ intersects transversely $\Pout$ in time $t=1$.
It follows that the orbit $\cC_i = \bigcup_{t \in \R} X^t(C_i)$ of $C_i$ by the flow of $X$ in $P$ is homeomorphic to $C_i \times [0,1]$, on which the vertical field $\partial/\partial t$ is orbit equivalent to $X$.
We push the foliations $\cG^{s, \iin}$ and $\cG^{u, \iin}$ of $C_i$ by the flow of $X$ in each $\cC_i$ cylinder and we obtain two transverse foliations $\cG^s$ and $\cG^u$ of class $\cC^1$ on $P$.
These foliations are $X^t$-invariant, and extend the foliations $\hat \cG^s$ and $\hat \cG^u$ respectively.
They satisfy Proposition~\ref{prop: existence paif}.
\end{proof}

\subsection{Incomplete gluing maps and strong isotopy} 

\label{sec: normalization; subsec: incomplete gluing maps}
In this section, we show that we can obtain a weakened version of a gluing map on a \bb{} $(P,X)$, isotopic to a \bb{} $(P_0, X_0)$ equipped with a \sqt{} gluing map $\varphi_0$, and
which satisfies a relation compatible with the strong isotopy relation of triples.
Let us define it formally.

\begin{defi}[Incomplete gluing map] \label{def: incomplete gluing map}
Let $(P,X)$ be a \bb{}, $\cO_*$ the periodic orbits of $X$ contained in $\pP$, $\Pin$ the entrance boundary and $\Pout$ the exit boundary of $(P,X)$.
An \emph{incomplete gluing map of $(P,X)$} is a $\cC^1$-involution $\hat \varphi \colon \pP \ssm \cO_* \to \pP \ssm \cO_*$ such that
\begin{enumerate}
    \item there exists a partition $\pP = \partial_1 P \, \sqcup \, \partial_2 P$ where $\partial_1 P$ and $\partial_2 P$ are unions of \ccs{} of $\pP$, and such that $\hat \varphi (\partial_1P \ssm \cO_*) = \partial_2P \ssm \cO_*$;
    \item $\hat \varphi (\Pout) = \Pin$.
\end{enumerate}
\end{defi}

In particular, the restriction of a gluing map $\varphi$ of a block $(P,X)$ to the complementary of the periodic orbits $\cO_*$ in $\pP$ is an incomplete gluing map, but the interest of introducing this notion is to consider \diffs{} that do not extend to $\cO_*$.
We consider the following relation on incomplete gluing maps.

\begin{defi}[Strongly isotopic incomplete gluing maps] 
\label{def: strongly isotopic incomplete gluing maps}
We will say that two incomplete gluing maps $\hat \varphi_0$ of $(P_0, X_0)$ and $\hat \varphi_1$ of $(P_1, X_1)$ are strongly isotopic if there exists a block isotopy $\{(P_t, X_t)\}_{t \in [0,1]}$ and a continuous family of homeomorphisms $\{ h_t \colon \pP_0 \ssm \cO_{0,*} \to \pP_t \ssm \cO_{t,*} \}_{t \in [0,1]}$, such that
\begin{enumerate}
        \item $h_0= \Id$,
        \item $h_t$ maps the lamination $\cL_0 \ssm \cO_{0,*}$ on the lamination $\cL_t \ssm \cO_{t,*}$,
        \item $h_1$ maps the lamination $\!(\varphi_0)_{*}\!(\cL_0 \ssm \cO_{0,*}\!)\!$ on the lamination $\!(\varphi_1)_{*} \!(\cL_1 \ssm \cO_{1,*}\!)$.
\end{enumerate}
\end{defi}

\begin{rmk}
We recover Item~\ref{def: strongly isotopic triples; it: lamination} of Definition~\ref{def: strongly isotopic triples} of strongly isotopic triples, i.e., if $(P_0, X_0, \varphi_0)$ and $(P_1, X_1, \varphi_1)$ are two strongly isotopic triples in the sense of Definition~\ref{def: strongly isotopic triples}, then the restriction of $\varphi_0$ and the restriction of $\varphi_1$ are strongly isotopic incomplete gluing maps.
Intuitively, this relation ensures that the intersection patterns of the laminations are the same.
\end{rmk}

We will consider incomplete gluing maps which achieve a \say{good} lamination intersection pattern.

\begin{defi}[Strongly transverse incomplete gluing map] \label{def: strongly transverse incomplete gluing maps}
Let $(P,X)$ be a \bb{} and $\hat \varphi$ an incomplete gluing map of $(P,X)$.
We say that it is a \emph{strongly transverse incomplete gluing map} if $(\cL \ssm \cO_*, \hat \varphi_* (\cL \ssm \cO_*))$ is a pair of strongly transverse laminations on $(P,X)$ (Definition~\ref{def: sqt laminations}).
\end{defi}

In particular, the restriction of a \sqt{} gluing map of a block $(P,X)$ to the complementary of the periodic orbits $\cO_*$ in $\pP$ is a strongly transverse incomplete gluing map.
We show the following

\begin{prop} \label{prop: induced gluing map on orbit equivalent block}
Let $(P,X)$ a filled affine block
and $(P_0,X_0)$ a block isotopic to $(P,X)$ among orbit equivalent blocks.
Let $\hat \varphi_0$ be a \st{} incomplete gluing map of $(P_0,X_0)$.
Then there exists a strongly transverse incomplete gluing map $\hat \varphi$ of $(P, X)$ strongly isotopic to $\hat \varphi_0$.
\end{prop}

\begin{proof}
    Let $\{(P_t, X_t)\}_{t \in [0,1]}$ be the isotopy of orbit equivalent blocks, with $(P_1, X_1) = (P,X)$.
    Let $H_t : P_0 \to P_t$ be a continuous family of orbit equivalences between the flow of $X_0$ and the flow of $X_t$.
    Let $\hat \varphi$ be the conjugate of $\hat \varphi_0$ by the orbit equivalence $H$.
    It satisfies Definition~\ref{def: incomplete gluing map}, but is only continuous.
We will say that it is a \emph{topological incomplete gluing map} of $(P, X)$.
There exists a continuous family of homeomorphisms $$h_t \colon \pP_0 \ssm \cO_{*,0} \to \pP_t \ssm \cO_{*,t}$$ (this is the restriction of the $H_t$ orbit equivalence) which maps the lamination $\cL_0 \ssm \cO_{*,0}$ on the lamination $\cL_t \ssm \cO_{t, *}$ and such that $h_1$ maps the lamination $(\hat \varphi_0)_*(\cL_0 \ssm \cO_{*,0})$ on the lamination $(\hat \varphi_1)_*(\cL \ssm \cO_*)$.
The map $\hat \varphi_1$ satisfies the definition of strong isotopy (Definition~\ref{def: strongly isotopic incomplete gluing maps}) of incomplete gluing maps, except that it is not differentiable.
Let $(\cG^s, \cG^u)$ a \paif{} and $\cG^{s, \iin}$, $\cG^{u, \out}$ the trace on $\Pin$ and $\Pout$ respectively (Proposition~\ref{prop: boundary foliation induced by paif}).
We can smooth the image of the lamination $\cL^\out$ in $\Pin$ in small charts $\{D_i\}_i$ for the foliation $\cG^{s, \iin}$, so that in any chart $D_i$, every leaf of $(\hat \varphi_1)_* \cL^\out \cap D_i$ intersect every leaf of $\cL^\iin \cap D_i$ at a unique point in $D_i$, which is possible by (topological) transversality.
We can do the smoothing so that the new smooth lamination $\cK$ is smoothly transverse to $\cL^\iin$ and satisfies the same previous property of intersection in $D_i$.
It follows that the two laminations are strongly isotopic, in the sense that they are isotopic via \homeos{} of $\Pin$ preserving $\cL^\iin$.
Now, there exist a smooth map $\hat \varphi$ of $\pP \ssm \cO_*$, $\cC^0$-close to $\hat \varphi_1$, which map $\cL^\out$ to $\cK$.
By construction, this \diff{} is an incomplete gluing map strongly isotopic to $\hat \varphi$, hence to $\hat \varphi_0$.
\end{proof}

\subsection{Straightening the boundary of a \bb{}}
\label{sec: normalization; subsec: boundary straightening}

In this subsection, we will \say{straighten} the boundary of the block. More precisely, we show

\begin{prop} \label{prop: straight block associated}
Let $(P,X)$ an affine block and $(\tilde P, \tilde X)$ a \minc{}.
There exists a \bb{} $(P_1, X_1)$, such that $(\tilde P, \tilde X)$ is a \minc{} of $(P_1, X_1)$, and such that for any orbit $\cO_i$ of $X_1$ in $\pP_1$, there exists a
    \lcs{} $(\cV_i, \xi_i = (x,y,\theta) \in \R^2 \times \R/\Z)$ of $\cO_i$ for the flow of $\tilde X$, compatible with the \acs{} $(D_i, \chi_i)$ of the section $\Sigma$, and in which $\pP_1$ coincides with the diagonal $\{x=y\}$. More precisely,
    \begin{enumerate}
        \item if $D_i \subset \Sigma$ is the disk which intersects $\cO_i$, then it coincides with $\{ \theta = 0 \}$, and the coordinates $\xi_i$ and $\chi_i$ coincide on $D_i$;
        \label{prop: straight block; it: compatible coordinate}
        \item $P_1$ coincides with the region $\{ y \geq x \}$ in coordinates $\xi_i$ on $\cV_i$.
        \label{prop: straight block; it: diagonal boundary}
    \end{enumerate}
This block is unique up to isotopy along the flow.
\end{prop}

We refer to Figure~\ref{fig: normalized block on snc}.

\begin{coro} \label{coro: straight block is normalized and isotopic}
If $(P,X)$ satisfies Item~\ref{def: normalized, it: affine section},~\ref{def: normalized, it: multipliers},~\ref{def: normalized, it: foliation} of Definition~\ref{def: normalized block}, then the block $(P_1, X_1)$ given by Proposition~\ref{prop: straight block associated} is a normalized block isotopic to $(P,X)$.
\end{coro}

\begin{proof}[Proof of Corollary~\ref{coro: straight block is normalized and isotopic}]
The isotopy is a consequence of the fact that the blocks $(P_0, X_0)$ and $(P_1, X_1)$ have a common \minc{} and of Proposition \ref{prop: isotopy vs orbit eq}, Item~\ref{prop: isotopy vs orbit eq; it: same minc implies isotopy}.
The block $(P, X)$ satisfies Items~\ref{def: normalized, it: affine section},~\ref{def: normalized, it: multipliers} and~\ref{def: normalized, it: foliation} of Definition~\ref{def: normalized block}, hence so does $(P_1, X_1)$ by the property of having a common \minc{}.
According to Proposition~\ref{prop: straight block associated}, $(P_1, X_1)$ satisfies Item~\ref{def: normalized, it: straight boundary}.
The block $(P_1, X_1)$ thus satisfies Definition~\ref{def: normalized block} of a normalized block.
\end{proof}

\begin{figure}[htb]
    \centering
    \vspace*{-1em}
    \includegraphics[height=0.27\textheight]{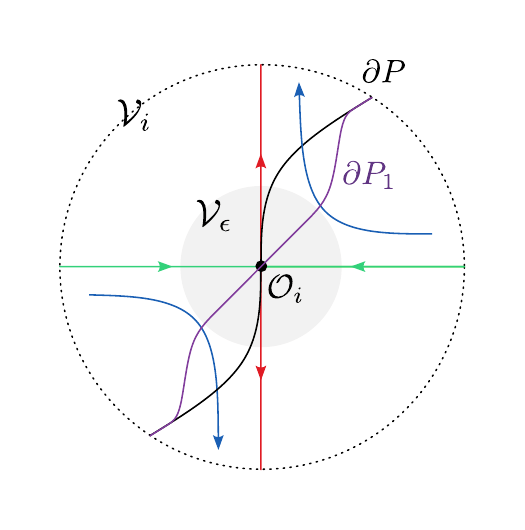}
    \vspace*{-1em}
    \centering \caption{Straightening the boundary of the block $(P,X)$ in a \minc{} $(\tilde P, \tilde X)$.}
    \label{fig: straightened boundary}
\end{figure}

\begin{proof}[Proof of Proposition~\ref{prop: straight block associated}]
Let $\cO_*$ be the collection of periodic orbits of $X$ contained in $\pP$.
For each orbit $\cO_i$, let $D_i$ be the unique disk which intersects the orbit $\cO_i$ with affine coordinates $\chi_i = (x,y) \in \R^2$ (Definition~\ref{def: normalized block}, Item~\ref{def: normalized, it: affine section})
Let $\cV_i$ be a tubular neighborhood of $\cO_i$ contained in the suspension neighborhood of $D_i$, in other words contained in $\bigcup_{t=0}^1 \tilde X^t (D_i)$ (we assume that we have renormalized the vector field $\tilde X$ such that every disk $D_i$ returns on $\Sigma$ in time $1$).
Let $\xi_i = (x,y,\theta) \in \R^2 \times \R/\Z$ be a \lcs{} of $\cO_i$ on $\cV_i$ such that $D_i$ coincides with $\{\theta=0\}$ in $\cV_i$ and $\xi_i$ and $\chi_i$ coincide on $D_i \cap \cV_i$.
Up to make a translation, we can suppose that $\cO_i = \{x=y=0\}$.
$\tilde X^t(x,y,\theta) = (2^{-t} x, 2^t y, \theta +t)$.
Two components of $\pP \ssm \cO_*$ adjacent in $\cO_i$ are contained in two opposite quadrants of $\cO_i$ (Claim~\ref{claim: boundary quadrant and multipliers}), namely (up to reverse the orientation) the quadrants $\{x \geq 0, y \geq 0 \}$ and $\{x \leq 0, y \leq 0\}$.
We can thus define a smooth surface $S$, which coincides with $\pP$ outside the \nbhs{} $\cV_i$, and such that in each $\cV_i$, we have
\begin{itemize}
    \item $S$ coincides with the diagonal annulus of equation $\{ x=y \}$ in an $\epsilon$-neighbor\-hood of $\cO_i$ strictly contained in $\cV_i$ which we denote $\cV_\epsilon$,
    \item $S$ is a annulus transverse to the vector field $\tilde X$ in the annulus $\cV_i \, \ssm \, \cV_\epsilon$, and coincides with $\pP$ in the neighborhood of the boundary of $\cV_i$.
\end{itemize}
Then $S$ is a closed surface, \qt{} to the \vf{} $\tilde X$, and cuts in $\tilde P$ a submanifold $P_1$ which coincides with $P$ outside the union of the \nbhs{} $\cV_i$.
If we denote $X_1$ the restriction of $\tilde X$ to $P_1$, then the pair $(P_1, X_1)$ is a \bb{} with \minc{} $(\tilde P, \tilde X)$.
By construction, it satisfies Proposition~\ref{prop: straight block associated}.
If $(P',X')$ is another block which satisfies the proposition, we can make a $\cC^1$ change of coordinate to get back to the same \lcs{} $(\cV_i, \xi_i)$.
The boundaries $\pP'$ and $\pP_1$ in this neighborhood are annuli transverse to the stable and unstable manifolds of $\cO_i$ contained in the same two opposite quadrants, and which coincide in their boundary.
We deduce that they are isotopic along the flow of $\tilde X$ in $\cV_i$.
\end{proof}

\subsubsection*{Strongly isotopic incomplete gluing maps}
In this paragraph, we explain how to obtain a strongly transverse incomplete gluing map (Definitions~\ref{def: incomplete gluing map} and~\ref{def: strongly transverse incomplete gluing maps}) on a block $(P,X)$ having a common \minc{} with a block $(P_0, X_0)$ equipped with a \sqt{} incomplete gluing map $\hat \varphi_0$, and so as to satisfy the strong isotopy relation (Definition~\ref{def: strongly isotopic incomplete gluing maps}).

\begin{lem} \label{lem: induced incomplete gluing map on isotopic blocks}
Let $(P, X)$ and $(P_0,X_0)$ be two blocks with a common \minc{} $(\tilde P, \tilde X)$.
Let $\hat \varphi_0$ be a strongly transverse incomplete gluing map of $(P_0,X_0)$.
Then there exists a strongly transverse incomplete gluing map $\hat \varphi$ of $(P, X)$ strongly isotopic to $\hat \varphi_0$.
\end{lem}

\begin{rmk} \label{rmk: gluing map does not extend 2}
Even if the initial incomplete gluing map $\hat \varphi_0$ extends into a gluing map, this will not be the case for the incomplete gluing map $\hat \varphi$ obtained by Lemma~\ref{lem: induced incomplete gluing map on isotopic blocks}.
Typically, one will obtain $\hat \varphi$ by pushing $\hat \varphi_0$ by an isotopy along the flow for a time which may tend to infinity as one approaches the periodic orbits contained in the boundary.
\end{rmk}

\begin{proof}
According to Proposition~\ref{prop: isotopy vs orbit eq}, there exists a continuous family $(P_t, X_t)$ of \bbs{} having a common \minc{} $(\tilde P, \tilde X)$, and which joins $(P_0,X_0)$ to $(P, X)$.
These blocks all coincide outside a union $\cV_*$ of linearizing neighborhoods of the periodic orbits $\cO_*$ of $X_0$ contained in $\pP_0$.
There exists $H_t \colon \pP_0 \ssm \cO_* \to \pP_t \ssm \cO_*$ a unique isotopy along the orbits of $X_0$ trivial on the complementary of $\cV_*$.
It does not extend in general along the orbits $\cO_*$.
The family
$\hat \varphi_t = (H_t)_* \hat \varphi \colon \pP_t \ssm \cO_* \to \pP_t \ssm \cO_*$ 
is a continuous family of incomplete gluing maps of $(P_t, X_t)$.
Since the isotopy $\{ H_t \}_t$ preserves the stable and unstable manifolds of the flow of $\tilde X$, we deduce that the incomplete gluing maps $\hat \varphi := \hat \varphi_1$ and $\hat \varphi_0$ are strongly isotopic in the sense of Definition~\ref{def: strongly isotopic incomplete gluing maps}.
Since the isotopy is of class $\cC^1$, the strong transversality of $\hat \varphi$ leads to the strong transversality of $\hat \varphi$.
\end{proof}

\subsection{Normalization of gluing map}
\label{sec: normalization; subsec: normalization of gluing map}
In this subsection,
$(P,X)$ denotes a normalized filled block (Definition~\ref{def: filled block} and~\ref{def: normalized block}).
We consider a strongly transverse incomplete gluing map $\hat \varphi \colon \pP \ssm \cO_* \to \pP \ssm \cO_*$.
Suppose that $\hat \varphi$ is strongly isotopic to the restriction of a \sqt{} gluing map $\varphi_0$ (in the sense of Definition~\ref{def: strongly isotopic incomplete gluing maps}).\footnote{$\hat \varphi$ does not extend a priori to $\cO_*$, but is strongly isotopic to an incomplete gluing map which extends to a gluing map $\varphi_0$. 
Note that $\varphi_0$ is not in general a gluing map of $(P,X)$, but only a gluing map of a block $(P_0, X_0)$ isotopic to $(P,X)$.
Typically, the block $P_0$ is obtained by pushing $\partial P$ along the orbits of the flow for an infinite time.}
The block $(P,X)$ is provided with a \paif{} $(\cG^s, \cG^u)$.
We denote $(\cG^{s, \out},\cG^{u, \out})$ the pair of induced foliations on $\Pout$ and $(\cG^{s, \iin},\cG^{u, \iin})$ the pair of induced foliations on $\Pin$.

The goal of this section is to construct a normalized gluing map $\varphi_1$ on $\pP$, so as to obtain a normalized triple $(P, X, \varphi_1)$ strongly isotopic to $(P_0,X_0,\varphi_0)$ in the sense of Definition~\ref{def: strongly isotopic triples}.

\begin{prop} \label{prop: strongly isotopic normalized gluing map}
There exists a normalized gluing map $\varphi_1$ of $(P, X)$ whose restriction to $\pP \ssm \cO_*$ is strongly isotopic to $\hat \varphi$.
\end{prop}

Recall that this means that $\varphi_1$ satisfies Definition~\ref{def: normalized gluing map}, in other words coincides with the reflection $(x,x,\theta) \mapsto (-x, -x, \theta)$ in \ncss{} of the orbits $\cO_i \in \cO_*$, and $\varphi_1$ maps the foliation $\cG^{u, \out}$ on a foliation transverse to $ \cG^{s, \iin}$ on $\pP$.

\subsubsection*{Summary of the proof}
The organization of the proof is as follows.
\begin{itemize}[leftmargin=*]
    \item In a first step, we show that we can make an isotopy among the strongly transverse incomplete gluing maps of $(P,X)$ so that $\hat \varphi$ maps the foliation $\cG^{u, \out}$ on a foliation transverse to $\cG^{s, \iin}$ in $\Pin$.
    This step will allow us to obtain Item~\ref{def: normalized gluing; it: foliation transverse} of Definition~\ref{def: normalized gluing map} of a normalized gluing map.
    This is Proposition~\ref{prop: incomplete gluing map with transverse foliation}.

    The key idea is to first locally straighten the image of the laminations in a \nbh{} of the boundary $\partial \cL^\iin$ and $\partial \cL^\out$ of $\cL^\iin$ and $\cL^\out$, which are sublaminations with empty interior (Lemma~\ref{lem: local straight foliation}).
    In a second step, we push the image of the laminations $\cL^\out$ transversally to the foliation $\cG^{s,  _iin}$ in the complementary of an adapted \nbh{} of $\cL^\iin$ (Lemma~\ref{lem: psi^in put L^out transverse to G^s,in}), and similarly put the image of $\cL^\iin$ transverse to $\cG^{u, \out}$ in the complementary of an adapted \nbh{} of $\cL^\out$ (Lemma~\ref{lem: psi^out put L^in transverse to G^u,out}).
    Once this is done, it remains to modify the mapped foliation on the complementary of $\hat \varphi_* \cL^\out \cup \cL^\iin$.
    By strong transversality of $\hat \varphi$, these complements are rectangles on which the foliations are trivial.

    \item In a second step, we show that we can modify the gluing map $\hat \varphi$ obtained in the previous step by strong isotopy so that $\hat \varphi$ is equal to the reflection $(x, x, \theta) \mapsto (-x, -x, \theta)$ in \ncss{} of the periodic orbits of $X$ contained in $\pP$.
    Such an incomplete gluing map then extends to a normalized gluing map $\varphi_1$ of $(P,X)$.    
    This result is summarized in Proposition~\ref{prop: extension trivial for incomplete gluing map}.
    
    The technique is to construct an isotopy which maps each leaf of $\hat \varphi_* \cG^{u, \out}$ to the \say{corresponding} leaf of $\cG^{u, \iin}$ for parameters given by the \ncs{} $\xi_i$ in the neighborhood of $\cO_i$, and an analogous isotopy for the foliation $\cG^s$.
    We then show that the composition extends trivially on $\cO_i$ in coordinates $\xi_i$.
    
\end{itemize}

\subsubsection*{Step 1: Transverse form for an incomplete gluing map}
A first step in the proof of Proposition~\ref{prop: strongly isotopic normalized gluing map} consists in making an isotopy of $\hat \varphi$ among the strongly transverse incomplete gluing maps of $(P,X)$, so that the image of $\cG^{u, \out}$ by the incomplete gluing map is transverse to $\cG^{s, \iin}$ on $\Pin$.

\begin{prop} \label{prop: incomplete gluing map with transverse foliation}
There exists an incomplete gluing map $\hat \varphi_1$ of $(P,X)$ strongly isotopic to $\hat \varphi$, such that the foliation $(\hat \varphi_1)_* \cG^{u, \out}$ is transverse to the foliation $\cG^{s, \iin}$ on $\Pin$.
\end{prop}

\begin{figure}[htb]
    \centering
    \vspace*{-2em}
    \includegraphics[height=0.27\textheight]{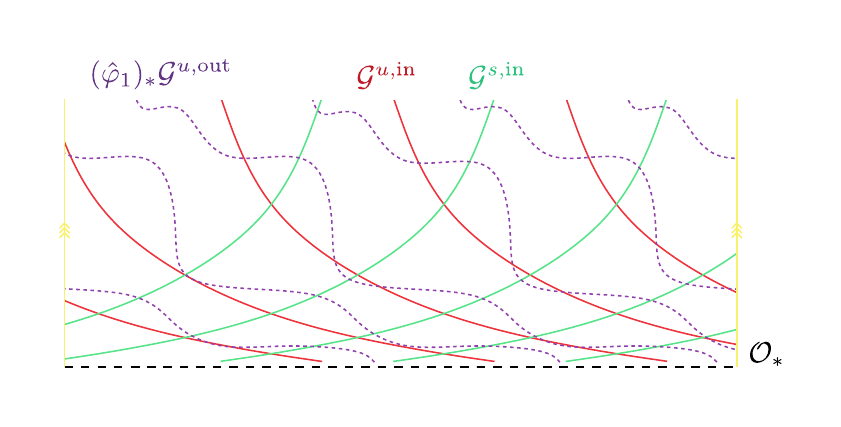}
    \vspace*{-1em}
    \centering \caption{Foliations on a neighborhood of $\cO_*$ in $\Pin$}
    \label{fig: transverse foliation gluing map}
\end{figure}

We repeat the proof of \cite[Proposition~5.2]{beguinBuildingAnosovFlows2017}, with the difference that we are working with laminations $\cL^\iin$ and $\cL^\out$ with non-empty interior due to the attractors and repellers in $\Lambda$, on a non-compact set $\pP \ssm \cO_*$, and the isotopies we make on $\pP \ssm \cO_*$ to straighten the image of the foliations will not extend on $\cO_*$.
The trick is to consider $\partial \cL^\iin$ and $\partial \cL^\out$ the boundary of $\cL^\iin$ and $\cL^\out$ respectively (Definition~\ref{def: boundary of a lamination}), which is a sublamination with empty interior.
We first show the following lemma.

\begin{lem} \label{lem: local straight foliation}
    There exist an incomplete gluing map $\hat \varphi'$ strongly isotopic to $\hat \varphi$ such that
    \begin{enumerate}
        \item \label{lem: local straight foliation; it: G^u,out}
        $\hat \varphi'_* \cG^{u, \out} = \cG^{u, \iin}$ in a \nbh{} of $\hat \varphi'_*( \cL^\out) \cap \partial \cL^\iin$, 
        \item \label{lem: local straight foliation; it: G^s,out}
        $\hat \varphi'_* \cG^{s, \out} = \cG^{s, \iin}$ in a \nbh{} of $\cL^\iin \cap \hat \varphi'_* (\partial \cL^\out)$.
    \end{enumerate}
\end{lem}

The key points of the proof are the following.
The transverse intersection of $\hat \varphi_* \cL^\out$ with the boundary $\partial \cL^\iin$ of $\cL^\iin$ has empty interior, and we can cover it by arbitrarily thin domains where the foliations are transverse.
We first straighten the image of the foliation $\cG^{u, \out}$ in the \nbh{} of the compact leaves of $\partial \cL^\iin$ along arbitrarily thin annuli by using the general Lemma~\ref{lem: general straightening on an annulus}.
Then we cover the other \ccs{} of the intersection by arbitrarily thin rectangles and (half)-bands and use the general Lemma~\ref{lem: general straightening on a disk} to straighten the image of the foliation in the \nbh{} of the intersection, which proves Item~\ref{lem: local straight foliation; it: G^u,out}.
We can do it so as to preserves the foliation $\cG^{s, \iin}$, which allows to proves Item~\ref{lem: local straight foliation; it: G^s,out} in the same way by reversing in $\Pout$ without destroying what have been done.

We first state two general lemmas.

\begin{figure}[h]
    \centering
    \vspace*{-1em}
    \includegraphics[width=0.95\textwidth]{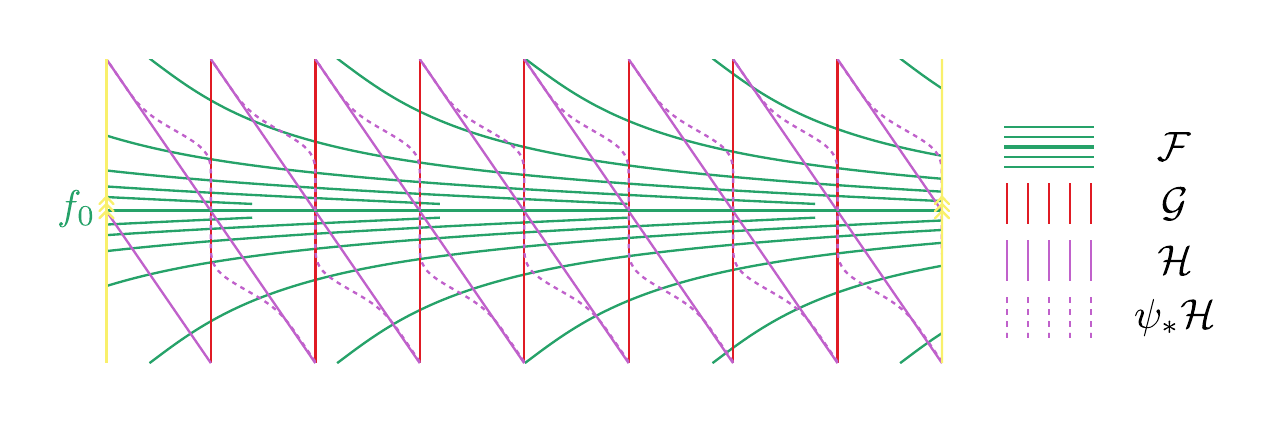}
    \vspace*{-1em}
    \caption{Action of the \diff{} $\psi: A \to A$ on the foliation $\cH$}
    \label{fig: general straightening on annulus}
\end{figure}

\begin{lem} \label{lem: general straightening on an annulus}
    Let $A =\R/\Z \times [-1, 1]$ an annulus, endowed with $3$ foliations $\cF, \cG, \cH$ such that $\cG$ and $\cH$ are trivial foliations transverse to $\cF$, and $\cF$ have a single compact leaf $f_0 := \R/\Z \times \{0\}$ with contracting holonomy.
    Then there exist a \diff{} $\psi : A \to A$, isotopic to the identity via \diffs{} preserving each leaf of $\cF$ and equal to the identity in a \nbh{} of $\partial A$, and such that leaves of $\psi_* \cH$ are leaves of $\cG$ in the \nbh{} of $f_0$.
\end{lem}
We refer to Figure~\ref{fig: general straightening on annulus}.

\vspace*{-1em}
\begin{proof}
Let $(\theta,x)$ the coordinate of $\R/\Z \times [-1, 1]$.
Up to change the coordinate system, we can suppose that $\cG$ is the vertical foliation $\theta = \cst$.
Let $\tilde A:= \R \times [-1,1]$ the universal cover of $A$ with coordinate $(\tilde \theta, x)$ and $\pi : \tilde A \to A$ the projection.
Denote $\tilde \cF$, $\tilde \cG$, $\tilde \cH$ the lifts of foliations in $\tilde D$.
Denote $\tilde f(\tilde \theta, x)$ the leaf of $\tilde \cF$ through $(\tilde \theta,x)$, and similarly $\tilde g(\tilde \theta, x)$ and $\tilde h(\tilde \theta,x)$ for $\cG$ and $\cH$.

\begin{claim}
There exist $\epsilon>0$ such that for all $\vert x \vert < \epsilon$ and $\tilde \theta \in \R$, 
the leaf $\tilde f(\tilde \theta, x)$ intersect the leaf $\tilde h(0, \tilde \theta)$ at a unique point $p = p(\tilde \theta, x) \in \tilde A$.
\end{claim}

\begin{proof}
The foliation $\tilde \cH$ is transverse to the leaf $\tilde f_0 = \R \times \{ 0 \}$ of $\tilde \cF$.
By continuity of the foliations, there exist $\epsilon>0$ such that every leaf $\tilde f(0, x)$ intersect $\tilde h(0, 0)$ for $\vert x \vert < \epsilon$.
As the holonomy of $\tilde \cF$ along $\tilde f_0$ is contracting, up to reduce $\epsilon$, we have that every leaf $\tilde f(0, x)$ crosses from $\tilde h(0, 0)$ to $\tilde h(1, 0)$ for $\vert x \vert < \epsilon$ (hence crosses every leaf $\tilde h(\tilde \theta, 0)$ for $\tilde \theta \in [0,1]$) and this holonomy map is contracting.
By invariance by integer translation $\tilde \theta \to \tilde \theta +1$, we deduce that the leaf $\tilde f(\tilde \theta, x)$ intersect the leaf $\tilde h(\tilde \theta, 0)$ for every $\vert x \vert < \epsilon$ and any $\tilde \theta \in\R$.
Now by transversality of the foliations, this intersection is unique.
\end{proof}

The map $p \colon \R \times (-\epsilon, \epsilon) \to \R \times [-1,1]$; $(\tilde \theta, x) \mapsto p(\tilde \theta, x)$ thus defined is a smooth injective map (because the foliations are smoothly transverse), commuting with the integer translation along $\tilde \theta$, which preserves each leaf of $\tilde \cF$, and map $\tilde \cH$ to $\tilde \cG$.
It remains to do a barycentric isotopy in the leaves of $\cF$ away from $x = 0$ to get the identity in the \nbh{} of $\vert x \vert = \epsilon$, and we extend $p$ by the identity on $\tilde A$.
More precisely let $\tau \colon [0,\epsilon] \to [0,1]$ be a smooth decreasing function which is equal to $1$ in a \nbh{} of $0$ and $0$ in a \nbh{} of $\epsilon$.
Define $\hat p (x, \tilde \theta)$ the weight barycenter by $\tau (\vert x \vert)$ between the point $p(x, \tilde \theta)$ and $(x, \tilde \theta)$ in the leaf $\tilde f(\tilde \theta, x)$.
Extend $\hat p$ by the identity on $\tilde A$.
This define a \diff{} $\hat p$, equal to the identity in the neighborhood of $\partial A$,
isotopic to the identity among \diffs{} which preserves the foliation $\tilde \cF$ leaf-to-leaf,
and which maps leaves of $\tilde \cH$ to leaves of $\tilde \cF$ on a neighborhood of $\tilde f_0$.
This \diff{} commutes with the integer translation along the coordinate $\tilde \theta$, so it is quotiented into a \diff{} $\psi : A \to A$ which satisfies the desired properties.
\end{proof}

\begin{lem} \label{lem: general straightening on a disk}
    Let $D$ a connected and simply connected domain of the plan endowed with 3 foliations $\cF, \cG, \cH$ such that $\cG$ and $\cH$ are transverse to $\cF$.
    Let $\cK$ be a sublamination of $\cF$ with empty interior and $\cL$ a sublamination of $\cH$.
    Suppose the intersection $\cK \cap \cL$ is disjoint from $\partial D$.
    Then there exist a \diff{} $\psi : D \to D$, isotopic to the identity via \diffs{} preserving each leaf of $\cF$ and equal to the identity in the \nbh{} of $\partial D$, and such that leaves of $\psi_* \cH$ are leaves $\cG$ in the \nbh{} of $\cK \cap \cL$.
\end{lem}
\begin{proof}
    The intersection $\Delta := \cK \cap \cL$ have empty interior and the \ccs{} $0$-dimensional and $1$-dimensional submanifold of $D$.
    As the foliations cannot have closed leaves by Poincar\'e--Hopf theorem, the \ccs{} of $\Delta$ are points, segments, lines and half-lines in $D$.
    We can cover $\Delta$ by a locally finite union $\{R_i\}$ such that:
    \begin{itemize}
        \item The sets $R_i$ are pairwise disjoint and diffeomorphic to either
        \begin{enumerate}
            \item a rectangle $[0,1]^2$, \label{it: rectangle}
            \item a half-band $\R_+ \times [0,1]$, or \label{it: half band}
            \item a band $\R \times [0,1]$. \label{it: band}
        \end{enumerate}
        \item The boundary of $R_i$ is disjoint from $\Delta$.
        \item The vertical segment of $R_i$ are on leaves of $\cG$.
        \item The horizontal segments of $R_i$ are on leaves of $\cF$ and the horizontal boundary of $R_i$ is disjoint from $\cK$.
    \end{itemize}
We refer to Figure~\ref{fig: general straightening on disk}.

\begin{figure}[htb]
    \centering
    \vspace*{-1.5em}
    \includegraphics[width=0.9\textwidth]{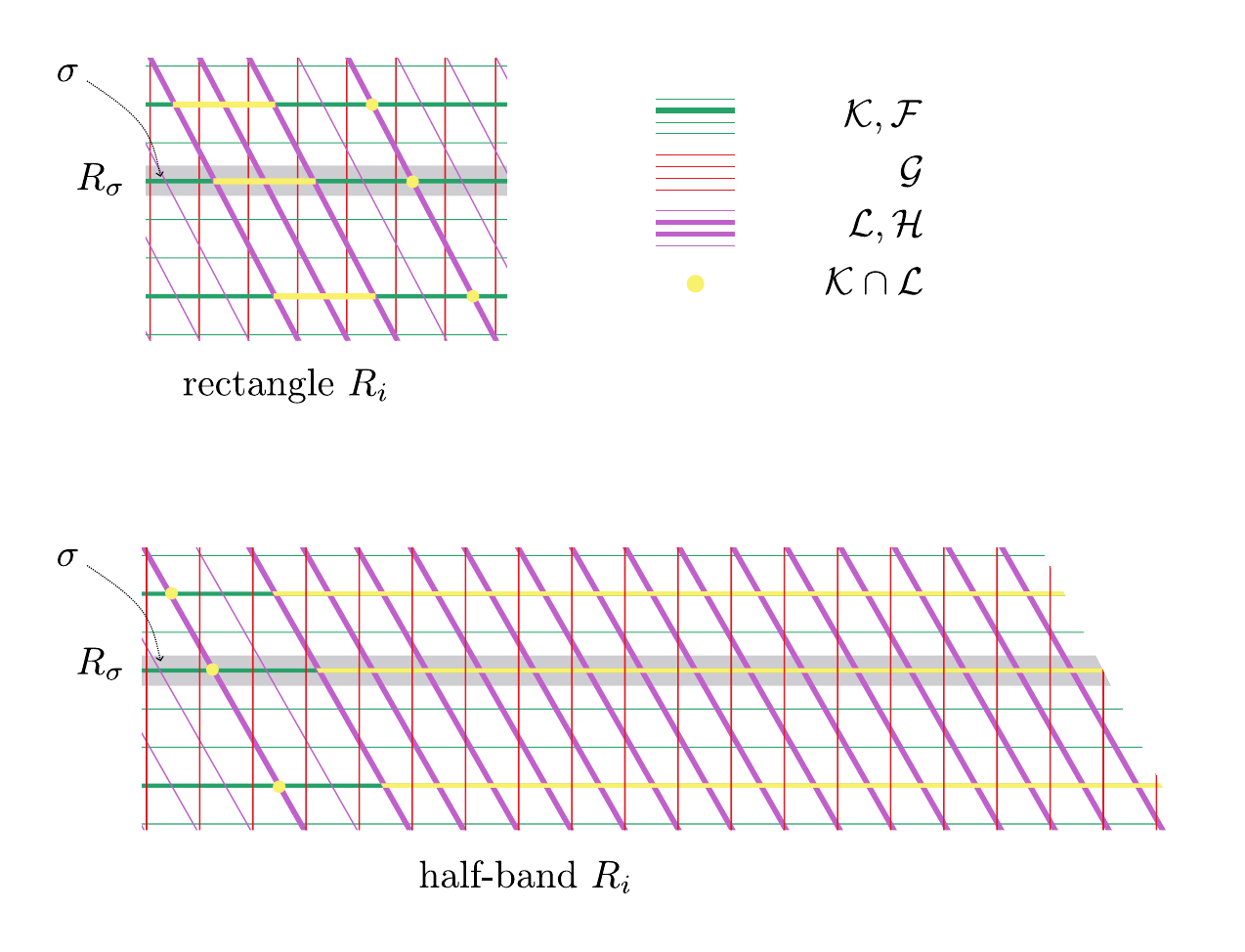}
    \vspace*{-1em}
    \caption{Foliated rectangle or half-band $R_i$}
    \label{fig: general straightening on disk}
\end{figure}

Consider Case~\ref{it: rectangle}.
Denote $R= R_i$.
Take $\sigma$ a \cc{} of $\cK$.
As $\cK$ is of empty interior, we can take an arbitrarily thin horizontal subrectangle $R_\sigma$ of $R$ containing $\sigma$ such that the horizontal boundary of $R_\sigma$ is disjoint from $\cK$.
Up to take $R_\sigma$ thin, every leaf of $\cH$ crosses $R_\sigma$ vertically in the \nbh{} of $\Delta \cap R_\sigma$, because $\Delta$ is disjoint from the horizontal boundary of $R_\sigma$.
One can locally straighten the foliation $\cH$ on the vertical foliation in a \nbh{} of $\Delta \cap R_\sigma$.
More precisely, there exist a \diff{} $\psi_\sigma: R \to R$ supported in the interior of $R_\sigma$ such that $\psi$ preserves the horizontal foliation $\cF$ and $\cH$ coincide with $\cG$ in the \nbh{} of $\Delta \cap R_\sigma$.
We cover the \ccs{} of $\cK$ by a finite disjoint union of such arbitrarily thin subrectangle $R_{\sigma_j}$ and we take $\psi_{R}$ as the product of the $\psi_{\sigma_j}$.
A barycentric isotopy in the leaves of $\cF$ shows that such a \diff{} is isotopic to the identity via \diff{} preserving each leaf of $\cF$, and equal to the identity in the \nbh{} of $D$.

We proceed the same in Case~\ref{it: half band} or Case~\ref{it: band}. Denote $R= R_i$ a band $\R \times [0,1]$ or a half-band $\R_+ \times [0,1]$.
Take $\sigma$ a \cc{} of $\cK$. It is a horizontal line or half-line.
As $\cK$ is of empty interior, we can take an arbitrarily thin horizontal sub-(half-)band $R_\sigma$ of $R$ containing $\sigma$ such that the horizontal boundary of $R_\sigma$ is disjoint from $\cK$.
Up to take $R_\sigma$ thin, every leaf of $\cH$ crosses $R_\sigma$ vertically in the \nbh{} of $\Delta \cap R_\sigma$, because $\Delta$ is disjoint from the horizontal boundary of $R_\sigma$.
The \diff{} $\psi_R$ is then defined as for Case~\ref{it: rectangle}.
\end{proof}

\begin{proof}[Proof of Lemma~\ref{lem: local straight foliation}]
We proceed in two steps.

\smallskip

\emph{Step 1: Push the image of $\cG^{u, \out}$ on $\cG^{u, \iin}$.}

\smallskip

The lamination $\partial \cL^\iin$ is a sublamination of $\cL^\iin$, hence of $\cG^{s, \iin}$, with empty interior, and $\hat \varphi_* \cL^\out$ is transverse to $\partial \cL^\iin$.
The \ccs{} of $\Delta := \partial \cL^\iin \cap \hat \varphi_* \cL^\out$ are points, segments, (half)-lines and circles.
There is a \nbh{} of $\Delta$ in $\Pin$  where $\cG^{s, \iin}$ and $\hat \varphi_* \cG^{u, \out}$ are transverse.
We first work in a \nbh{} of the circle component, which are compact leaves of $\partial \cL^\iin$ (hence of $\cG^{s, \iin}$).
They are in finite number with contracting holonomy (Proposition~\ref{prop: block boundary lam are qms}).
There exist a finite disjoint union of arbitrarily thin annuli $\{ A_i \}$, covering each circle component $\gamma_i$ of $\Delta$, such that $\hat \varphi_* \cG^{u, \out}$ is transverse to $\cG^{s, \iin}$ on $A_i$.
We use Lemma~\ref{lem: general straightening on an annulus} with $\cF = \cG^{s, \iin}$, $\cG = \cG^{u, \iin}$ and $\cH = \hat \varphi_* \cG^{u, \out}$ in each $A_i$.
As the collection $\{A_i\}_i$ is a finite collection of disjoint annuli, we can construct a \diff{} $\psi_A$ on $\pP \ssm \cO_*$ supported in the interior of the $A_i$'s, isotopic to the identity via \diff{} preserving each leaf of $\cG^{s, \iin}$, such that $(\psi_A \circ \hat \varphi)_* \cG^{u, \out}$ is transverse to $\cG^{s, \iin}$ on $A_i$ and coincide with $\cG^{u, \iin}$ in a \nbh{} $V$ of each $\gamma_i$.
Up to shrink $V$ we can suppose that the boundary of $V$ is transverse to both $\cG^{s, \iin}$ and $(\psi_A \circ \hat \varphi)_* \cG^{u, \out}$.

The composition $\psi_A \circ \hat \varphi$ is a strongly transverse incomplete gluing map, strong\-ly isotopic to $\hat \varphi$.
We work now in the complementary of the \nbh{} $V$ of the closed connected component $\gamma_i$ of $\Delta$.
Denote $\Delta'= ((\psi_A \circ \hat \varphi)_* \cL^\out \cap \partial \cL^\iin) \ssm V$.
The \ccs{} of $\Delta'$ are points, segments and (half)-lines. 
One can cover $\Delta'$ by a disjoint union of arbitrarily thin, connected and simply connected domains $\{ D_i \}$ in $\Pin$ such that the boundary of $D_i$ is disjoint from $\Delta$, and the foliation $(\psi_A \circ \hat \varphi)_* \cG^{u, \out}$ is transverse to $\cG^{s, \iin}$ on $D_i$.
We use Lemma~\ref{lem: general straightening on a disk} on each $D_i$, with $\cF = \cG^{s, \iin}$, $\cK = \partial \cL^\iin$, $\cG = \cG^{u, \iin}$, $\cH = (\psi_A \circ \hat \varphi)_* \cG^{u, \out}$, $\cL = \check \varphi_* \cL^\out$, but we replace the intersection $\cK \cap \cL$ by $(\cK \cap \cL) \ssm V$, which does not change the proof.
As the collection $\{D_i\}$ is a locally finite collection of disjoint domains, it gives a \diff{} $\psi_D$ on $\pP \ssm \cO_*$ supported in the interior of the $D_i$'s, isotopic to the identity via \diff{} preserving each leaf of $\cG^{s, \iin}$, equal to the identity in $V$, and such that leaves of $(\psi_D \circ \psi_A \circ \hat \varphi)_* \cG^{u, \out}$ are leaves of $\cG^{u, \iin}$ in a \nbh{} of $(\psi_D \circ \psi_A \circ \hat \varphi)_* \cL^\out \cap \partial \cL^\iin$.
Define $\psi^\iin : \pP \ssm \cO_* \to \pP \ssm \cO_*$, equal to the identity on $\Pout$ and to $\psi_D \circ \psi_A$ on $\Pin$.
Define $\hat \varphi_1 := \psi^\iin \circ \hat \varphi$.
It is a strongly transverse incomplete gluing map, strongly isotopic to $\hat \varphi$, such that $(\hat \varphi_1)_* \cG^{u, \out}$ coincide with $\cG^{u, \iin}$ in a \nbh{} of $(\hat \varphi_1)_* \cL^\out \cap \partial \cL^\iin$.

\smallskip

\emph{Step 2: Push the image of $\cG^{s, \out}$ on $\cG^{s, \iin}$.}

\smallskip
It remains now to push the image of $\cG^{s, \out}$ on $\cG^{s, \iin}$ in a neighborhood of the intersection $(\hat \varphi_1)_* \partial \cL^\out \cap \cL^\iin$, without destroying the previous work.
We consider $\Pout$ and essentially do the same work, considering the foliations $\cF = \cG^{u, \out}$, $\cG = \cG^{s, \out}$ and $\cH = (\hat \varphi_1)^\inv_*\cG^{s, \iin}$, and the sublaminations $\cK = \partial \cL^\out$ with empty interior and $\cL = (\hat \varphi_1)^\inv_* \cL^\iin$, and use Lemma~\ref{lem: general straightening on an annulus} then Lemma~\ref{lem: general straightening on a disk} as previously.
We construct a \diff{} $\psi^\out$ supported in $\Pout$, such that $\psi^\out$ is isotopic the identity, preserves each leaf of $\cG^{u, \out}$, and maps leaves of $(\hat \varphi_1)^\inv_* \cG^{s, \iin}$ to leaves of $\cG^{s, \out}$ in a \nbh{} of $(\psi^\out \circ \hat \varphi_1^\inv)_* \cL^\iin \cap \partial \cL^\out$.
Define $\hat \varphi_2 = \psi^\iin \circ \hat \varphi \circ (\psi^\out)^\inv$. It is a strongly transverse incomplete gluing map, strongly isotopic to $\hat \varphi$. 
Moreover, we check that $(\hat \varphi_2)_* \cG^{u, \out} = (\psi^\iin \circ \hat \varphi_1)_* \cG^{u, \out} = \cG^{u, \iin}$ in a \nbh{} of $(\hat \varphi_2)_* \cL^\out \cap \partial \cL^\iin = (\hat \varphi_1)_* \cL^\out \cap \partial \cL^\iin$, hence $\hat \varphi_2$ satisfies Lemma~\ref{lem: local straight foliation}.
\end{proof}

Before proving Proposition~\ref{prop: incomplete gluing map with transverse foliation} we need the general lemma.

\begin{lem} \label{lem: general transversalisation on a square}
    Let $C$ be either a square $[0,1]^2$, a half band $\R_+ \times [0,1]$ or a band $\R \times [0,1]$, endowed with a pair of foliation $\cF$ and $\cG$ such that the horizontal boundary of $C$ are leaves of $\cF$ and the vertical boundary of $C$ (when existing) are leaves of $\cG$.
    We suppose that $\cG$ is transverse to $\cF$ in the \nbh{} of $\partial C$ and every leaf of $\cG$ crosses $C$ from bottom to top.
    Then there is a \diff{} $\psi : C \to C$, isotopic to the identity, and such that $\psi_*\cG$ is transverse to $\cF$ on $C$.
\end{lem}

\begin{figure}[htb]
    \centering
    \vspace*{-1em}
    \includegraphics[width=0.66\textwidth]{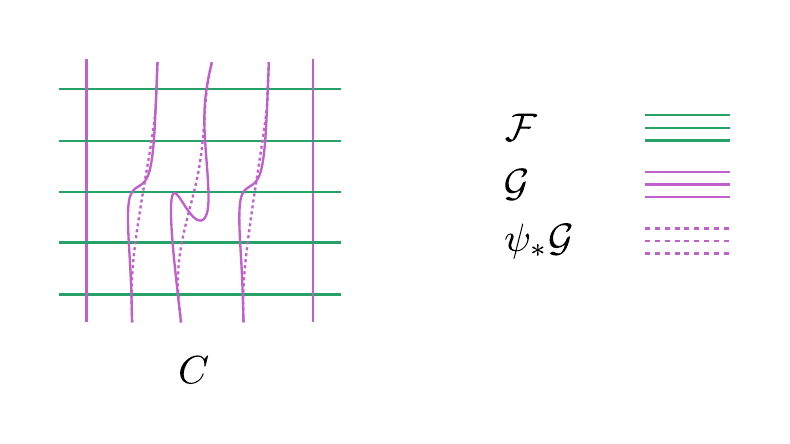}
    \vspace*{-1em}
    \caption{Isotopy of the foliation $\mathcal{G}$ in the square $C$}
    \label{fig: general transversalisation on square}
\end{figure}

We refer to Figure~\ref{fig: general transversalisation on square}.

\begin{proof}
The hypothesis imply that each of the foliations $\cF$ and $\cG$ are smoothly conjugated to trivial foliations, such that we may assume that $\cF$ is the horizontal foliation. 
Now, there is a foliation $\cH$ on $C$ transverse to the leaves of $\cF$ and coinciding with $\cG$ in a neighborhood of $\partial C$ and having the same holonomy from the horizontal bottom to the horizontal top of $C$. The foliations $\cG$ and $\cH$ are smoothly conjugated by a diffeomorphism $\psi$ which coincide with the identity close to $\partial C$, concluding.
\end{proof}

\begin{proof}[Proof of Proposition~\ref{prop: incomplete gluing map with transverse foliation}]
Up to make a strong isotopy among the incomplete gluing maps, we can assume that the \st{} incomplete gluing map $\hat \varphi$ of $(P,X)$ satisfies Lemma~\ref{lem: local straight foliation}. We proceed in two steps.

\begin{lem} \label{lem: psi^in put L^out transverse to G^s,in}
There exists a \diff{} $\psi^\iin \colon \pP \ssm \cO_* \to \pP \ssm \cO_*$, supported in $\Pin$, equal to the identity on a neighborhood of $\cL^\iin$, such that $(\psi^\iin \circ \hat \varphi)_* \cL^\out$ is transverse to $\cG^{s, \iin}$.
\end{lem}

\begin{proof}
First remark that there is a \nbh{} $U^\iin$ on $\cL^\iin$ such that $\hat \varphi_* \cL^\out$ is transverse to $\cG^{s, \iin}$.
Indeed, from Lemma~\ref{lem: local straight foliation}, there is a \nbh{} $\cO$ of $\partial \cL^\iin \cap \hat \varphi_* \cL^\out$ such that $\hat \varphi_* \cG^{u, \out}$ coincide with $\cG^{s, \iin}$.
We can choose a \nbh{} of $\partial \cL^\iin$ small enough such that the intersection with $\hat \varphi_* \cL^\out$ is contained in $\cO$.
Now the union of a \nbh{} of $\partial \cL^\iin$ with $\cL^\iin$ is a \nbh{} of $\cL^\iin$, satisfying the property.
Up to reduce, we can choose $U^\iin$ to be an adapted \nbh{} of $\cL^\iin$ (Definition~\ref{def: adapted nbh} and Lemma~\ref{lem: basis of adapted nbh}).
The complementary $\Pin \ssm U^\iin$ is a finite collection $\{C_i\}$ of pairwise disjoint adapted squares and (half)-bands (see Figure~\ref{fig: adapted square half band}).

The foliation $\hat \varphi_* \cL^\out$ coincide with $\cG^{u, \iin}$ in the \nbh{} of $\partial C_i$.
Moreover,

\begin{claim}
Each leaf of $\hat \varphi_* \cL^\out \cap C_i$ is a segment crossing $C_i$ from bottom to top.
\end{claim}

\begin{proof}
First remark that any leaf of $\hat \varphi_* \cL^\out \cap C_i$ 
crossing one horizontal boundary of $C_i$ cannot cross it twice.
Indeed, by extending this leaf on the strip $B_i \subset \Pin \ssm \cL^\iin$ containing $C_i$, it would give us a bi-gon $G$ with a boundary $g_1$ on a leaf of $\cL^\iin$ and a boundary $g_2$ on a leaf of $\hat \varphi_* \cL^\out$, with interior disjoint from $\cL^\iin$. 
Up to choose another smaller bi-gon inside $G$, either $G$ is foliated by $\hat \varphi_* \cL^\out$ either it is disjoint from $\cL^\iin \cup \hat \varphi_* \cL^\out$.
The first case is impossible by transversality of $\hat \varphi_* \cL^\out$ with $\cL^\iin$ and the second by strong transversality of $\hat \varphi$.
Suppose now the existence of a half-leaf $l$ of $(\hat \varphi_* \cL^\out) \cap C_i$ which remains in the interior of $C_i$, hence in the interior of the strip $B_i$ containing $C_i$.
As every laminations are pre-foliations, there is no compact leaf in $B_i$ therefore $l$ is asymptotic to an end of $B_i$.
Recall that $\hat \varphi$ is strongly isotopic to the restriction of a \sqt{} gluing map $\varphi$. 
In order to keep notations light, we suppose that $\varphi$ is a gluing map of the same building block $(P,X)$. The proof would be the same by considering an isotopic block.
By strong isotopy, there is a half-leave $l'$ of $\varphi_* \cL^\out$ which is contained in the interior of the strip $B_i$ and asymptotic to an end of $B_i$.
This end accumulates on a periodic orbit $\cO \in \cO_*$, and let $\cO' := \varphi(\cO)$ (where the orientation is preserved). 
Half-leaves of $\cL^\out$ in the \nbh{} of $\cO'$ in $\Pout$ accumulate on $\cO'$ with a contracting holonomy (for the orientation of $\cO$ by the flow), so the same occurs for half-leaves of $\varphi_* \cL^\out$ in the \nbh{} of $\cO = \varphi(\cO')$ in $\Pin$.
On the other hand, the holonomy of $l'$ on a small transversal $\tau$ to $\cO$ is bounded by the (expanding) holonomy of the leaves of $\cL^\iin$ bordering the strip $B_i$ containing $C_i$, hence cannot be contracting which is a contradiction.
\end{proof}

We can use Lemma~\ref{lem: general transversalisation on a square} on each pairwise disjoint $C_i$ with $\cF  = \cG^{s, \iin}$ and $\cG= \hat \varphi_* \cG^{u, \out}$ to construct a \diff{} $\psi_{C_i}$ of $\pP \ssm \cO_*$ supported in the interior of $C_i$.
The collection $\{C_i\}$ is pairwise disjoint and locally compact in $\Pin$, which allows to define the \diff{} $\psi^\iin$ by the product of $\psi_{C_i}$.
\end{proof}

Denote $\hat \varphi_1 = \psi^\iin \cap \hat \varphi$.

\begin{lem} \label{lem: psi^out put L^in transverse to G^u,out}
There exists a \diff{} $\psi^\out \colon \pP \ssm \cO_* \to \pP \ssm \cO_*$, supported in $\Pout$, equal to the identity on a \nbh{} of $\cL^\out$, such that $(\psi^\out \circ \hat \varphi_1^\inv)_* \cL^\iin$ is transverse to $\cG^{u, \out}$.
\end{lem}

\begin{proof}
    The proof is the same as for Lemma~\ref{lem: psi^in put L^out transverse to G^s,in}.
\end{proof}

Denote $\hat \varphi_2 = \psi^\iin \circ \hat \varphi \circ (\psi^\out)^\inv : \pP \ssm \cO_* \to \pP \ssm \cO_*$.
It is isotopic to $\hat \varphi$ through strongly transverse incomplete gluing maps, and $(\hat \varphi_2)_* \cL^\out = (\hat \varphi_1)_* \cL^\out$ is transverse to $\cG^{s, \iin}$ and $(\hat \varphi_2)_* \cG^{u, \out}$ is transverse to $\cL^\iin$.
So the foliations $\hat \varphi_* \cG^{u, \out}$ and $\cG^{s, \iin}$ may fail to be transverse only in the interior of a \cc{} of $\Pin \ssm (\cL^\iin \cup (\hat \varphi_2)_* \cL^\out)$.
By strong transversality of $\hat \varphi_2$, the closure of those \ccs{} are squares with two sides on leaves of $\cL^\iin$ and two sides on leaves of $(\hat \varphi_2)_* \cL^\out$.
We conclude by using Lemma~\ref{lem: general transversalisation on a square} with $\cF = \cG^{s, \iin}$ and $\cG = (\hat \varphi_2)_* \cG^{u, \out}$.
The result is a strongly transverse incomplete gluing maps, strongly isotopic to $\hat \varphi$ because it coincides with $\hat \varphi$ in a \nbh{} of the boundary lamination $\cL \ssm \cO_*$, and which satisfies Proposition~\ref{prop: incomplete gluing map with transverse foliation}.
\end{proof}

\subsubsection*{Step 2: Trivial extension of an incomplete gluing map.}
We keep the notations of the beginning of Subsection~\ref{sec: normalization; subsec: normalization of gluing map}.
In this step, up to make a strong isotopy of incomplete gluing maps, we assume that $\hat \varphi$ satisfies Proposition~\ref{prop: incomplete gluing map with transverse foliation}, i.e., $\hat \varphi_* \cG^{u, \out}$ is transverse to $\cG^{s, \iin}$ on $\Pin$.
For any $\cO_i \in \cO_*$, there exists a unique orbit $\cO_j \in \cO_*$ such that if $A$ and $A'$ are \ccs{} of $\pP \ssm \cO_*$ adjacent to $\cO_i$, then $\hat \varphi(A)$ and $\hat \varphi(A')$ are \ccs{} of $\pP \ssm \cO_*$ adjacent in $\cO_j$.
This pairing defines an involution without fixed point $\sigma \colon \{1,\dots, n\} \to \{ 1,\dots,n \}$, which determines the combinatorial gluing of~$\hat \varphi$.

We show the following proposition, which states that there exists a normalized gluing map $\varphi_1$ of $(P,X)$ whose restriction to the complementary of the orbits $\cO_*$ is strongly isotopic to $\hat \varphi$.
Recall that a normalized block $(P,X)$ is equipped with \ncss{} $(\cV_i ,\xi_i = (x,y,\theta))$ on tubular \nbhs{} of the periodic orbits $\cO_i \in \cO_*$ contained in $\pP$ (Definition~\ref{def: normalized block}, Item~\ref{def: normalized, it: straight boundary}).
In these coordinates, the boundary $\partial P$ has equation $x=y$, the exit boundary $\Pout$ is in the domain $\{ x<0 \}$, and the entrance boundary $\Pin$ in the domain $\{ x>0 \}$.

\begin{prop} \label{prop: extension trivial for incomplete gluing map}
There exists a gluing map $\varphi_1 \colon \pP \to \pP$
such that
\begin{enumerate}
    \item \label{prop: normalization gluing map; it: strong isotopy}
    The restriction of $\varphi_1$ to $\pP \ssm \cO_*$ is strongly isotopic to $\hat \varphi$
    (in the sense of Definition~\ref{def: strongly isotopic incomplete gluing maps});
    
    \item \label{prop: normalization gluing map; it: trivial snc}
    for any $\cO_i \in \cO_*$, there exist \ncss{} $(\cV_i, \xi_i = (x,y,\theta))$ and $(\cV_j, \xi_j = (x,y,\theta))$ of $\cO_i$ and $\cO_j= \varphi_1(\cO_i)$ such that the expression of $\varphi_1$ in these coordinates is $(x,x,\theta) \mapsto (-x,-x,\theta)$;
    
    \item \label{prop: normalization gluing map; it: foliation}
    the foliation $(\varphi_1)_* \cG^{s, \iin}$ is transverse to the foliation $\cG^{s, \iin}$ on $\Pin$.
\end{enumerate}
\end{prop}

In other words, $\varphi_1$ is a normalized gluing map of $(P,X)$ (Definition~\ref{def: normalized gluing map}) whose restriction is strongly isotopic to $\hat \varphi$.

Let us start by fixing some notations.
Let $(\cV_i, \xi_i)$ be a \ncs{} of the orbit $\cO_i \in \cO_*$.
The coordinates $\xi_i = (x,y,\theta) \in \R^2 \times \R/\Z$ on $\cV_i$ induce the coordinates $\rho_i = (x,\theta) \in \R \times \R/\Z$ on $\partial P \cap \cV_i$ by forgetting the coordinate $y$.
Up to multiply the coordinates $x$ and $y$ by a constant and then shrink $\cV_i$, we can assume that $\rho_i (\Pin \cap \cV_i) = ]0,1] \times \R/\Z$ and $\rho_i(\Pout \cap \cV_i) = [-1, 0[ \times \R/\Z$.
According to the foliations equation (Remark~\ref{rmk: equation affine foliation}) each leaf of $\cG^{s, \iin}$ and $\cG^{u, \iin}$ transversely intersects the section $\{ x=1 \}$ once, and each leaf of $\cG^{s, \out}$ and $\cG^{u, \out}$ transversely intersects the section $\{ x=-1 \}$ once.

\begin{figure}[htb]
    \centering
    \vspace*{-1em}
    \includegraphics[height=0.32\textheight]{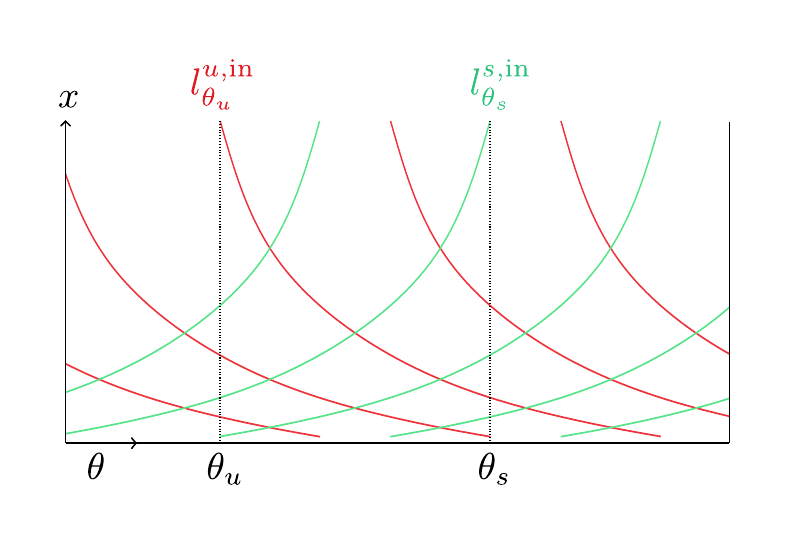}
    \vspace*{-1em}
    \centering \caption{Parameterization of the leaves of $\cG^{s, \iin}$ and $\cG^{u, \iin}$ by the coordinate $\theta$ on $\cV^\iin_i$}
    \label{fig: coord theta foliation}
\end{figure}

Denote $l^{s, \iin}_{i, \theta}$ and $l^{u, \iin}_{i, \theta}$ the leaves of $\cG^{s, \iin}$ and $\cG^{u, \iin}$ containing $(1,\theta)$ in coordinates $\rho_i$.
Similarly, denote $l^{s, \out}_{i, \theta}$ and $l^{u, \out}_{i, \theta}$ the leaves of $\cG^{s, \out}$ and $\cG^{u, \out}$ containing $(-1,\theta)$ in coordinates $\rho_i$ (see Figure~\ref{fig: coord theta foliation}).

We will show the following two lemmas, analogous to each other.

\begin{lem} \label{lem: identity cn psi in}
There exists a \diff{} $\psi^\iin \colon \pP \ssm \cO_* \to \pP \ssm \cO_*$ with support in $\Pin$ and a neighborhood $\cV^\iin$ of $\partial \Pin$ in $\Pin$ such that
\begin{enumerate}
    \item $\psi^\iin$ is isotopic to the identity among the \diffs{} of $\pP \ssm \cO_*$ which preserves the foliation $\cG^{s, \iin}$ leaf-to-leaf;
    \label{lem: identity psi in; it: preserve foliation}
    \item for any $i$ and $j=\sigma(i)$, for any $\theta \in \R/\Z$, the leaf $\psi^\iin ( \hat \varphi (l^{u,\out}_{j, \theta} ))$ coincides with the leaf $l^{u, \iin}_{i, \theta}$ on $\cV^\iin$.
    \label{lem: identity psi in; it: theta parameters}
\end{enumerate}
\end{lem}

\begin{lem} \label{lem: identity cn psi out}
There exists a \diff{} $\psi^\out \colon \pP \ssm \cO_* \to \pP \ssm \cO_*$ with support in $\Pout$ and a neighborhood $\cV^\out$ of $\partial \Pout$ in $\Pout$ such that
\begin{enumerate}
    \item $\psi^\out$ is isotopic to the identity among the \diffs{} of $\pP \ssm \cO_*$ which preserves the foliation $\cG^{u, \out}$ leaf-to-leaf;
    \label{lem: identity psi out; it: preserve foliation}
    \item for any $i$ and $j=\sigma(i)$, for any $\theta \in \R/\Z$, the leaf $\psi^\out ( \hat \varphi^\inv (l^{s,\iin}_{j, \theta} ))$ coincides with the leaf $l^{s, \out}_{i, \theta}$ on $\cV^\out$.
    \label{lem: identity psi out; it: theta parameters}
\end{enumerate}
\end{lem}

Let us show that Proposition~\ref{prop: extension trivial for incomplete gluing map} follows from Lemmas~\ref{lem: identity cn psi in} and~\ref{lem: identity cn psi out}.

\begin{proof}[Proof of Proposition~\ref{prop: extension trivial for incomplete gluing map}]
Denote $\hat \varphi_1 := \psi^\iin \circ \hat \varphi \circ (\psi^\out)^\inv: \pP \ssm \cO_* \to \pP \ssm \cO_*$.
According to Item~\ref{lem: identity psi in; it: preserve foliation} of Lemma~\ref{lem: identity cn psi in} and its analogue of Lemma~\ref{lem: identity cn psi out}, the \diff{} $\hat \varphi_1$ is an incomplete gluing map strongly isotopic to $\hat \varphi$.
Moreover $\hat \varphi$ maps the foliation $\cG^{u, \out}$ to a foliation transverse to $\cG^{s, \iin}$ on $\Pin$ by assumption, so it is the same for $\hat \varphi_1$, still according to Item~\ref{lem: identity psi in; it: preserve foliation} of Lemma~\ref{lem: identity cn psi in} for $\psi^\iin$, and its analogue for $\psi^\out$.
Let $\cW^\iin := \cV^\iin \cap \hat \varphi_1 (\cV^\out)$, where $\cV^\iin$ and $\cV^\out$ are the neighborhoods given by Lemmas~\ref{lem: identity cn psi in} and~\ref{lem: identity cn psi out}.
This is a neighborhood of $\partial \Pin$ in~$\Pin$.

Let $i$ and $j=\sigma(i)$ be $\theta_s \in \R/\Z$ and $\theta_u \in \R/\Z$.
Then the following equalities are true on $\cW^\iin$ (see Figure~\ref{fig: normalization image of foliation}):
\begin{align*}
    \hat \varphi_1 (l^{u, \out}_{j, \theta_u}) & = \psi^\iin \circ \hat \varphi \circ (\psi^\out)^\inv (l^{u, \out}_{j, \theta_u}) \\
    & = \psi^\iin \circ \hat \varphi (l^{u, \out}_{j,\theta_u})
    & \mb{(Lemma~\ref{lem: identity cn psi out}, Item~\ref{lem: identity psi out; it: preserve foliation})}  \\
    & = l^{u, \iin}_{i,\theta_u},& \mb{(Lemma~\ref{lem: identity cn psi in}, Item~\ref{lem: identity psi in; it: theta parameters})}
\end{align*}
\begin{align*}
    \hat \varphi_1 (l^{s, \out}_{j, \theta_s}) & = \psi^\iin \circ \hat \varphi \circ (\psi^\out)^\inv (l^{s, \out}_{j, \theta_s}) \\
    & = \psi^\iin \circ \hat \varphi (\hat \varphi^\inv (l^{s, \iin}_{i, \theta_s})) = \psi^\iin (l^{s, \iin}_{i,\theta_s})
    & \mb{(Lemma~\ref{lem: identity cn psi out}, Item~\ref{lem: identity psi out; it: theta parameters})}  \\
    & = l^{s, \iin}_{i,\theta_s}. & \mb{(Lemma~\ref{lem: identity cn psi in}, Item~\ref{lem: identity psi in; it: preserve foliation})}
\end{align*}

\begin{figure}[htb]
    \centering
    \vspace*{-1em}
    \includegraphics[height=0.27\textheight]{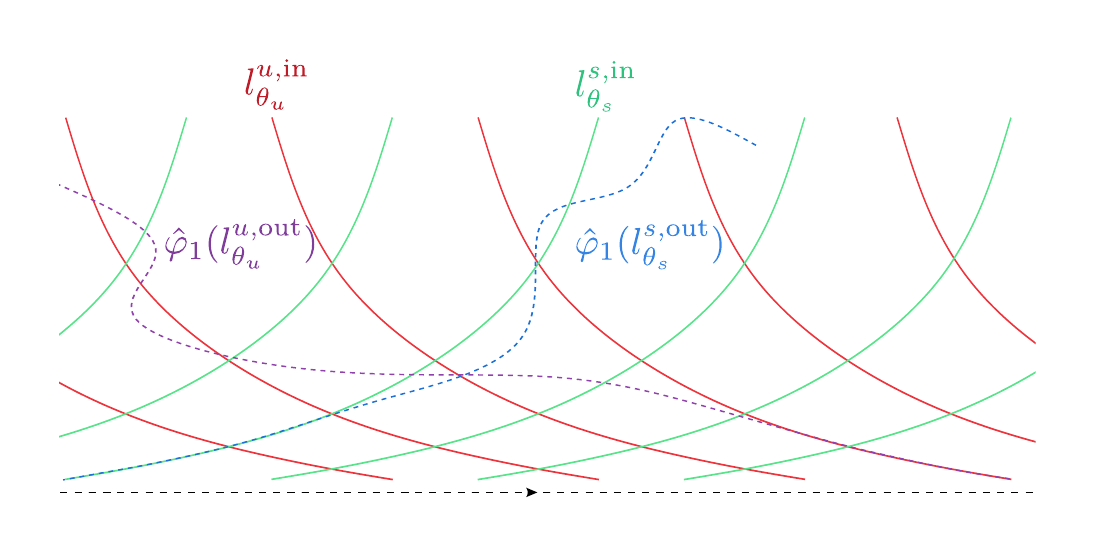}
    \vspace*{-1em}
    \centering \caption{Image of the foliations by $\hat \varphi_1$}
    \label{fig: normalization image of foliation}
\end{figure}

Recall that the flow of $X$ is linear in coordinates $\xi_i = (x, y, \theta)$ on $\cV_i$, given by the expression $X^t(x,y,\theta) = (2^{-t} x, 2^t y, \theta+t)$, and the boundary $\partial P$ coincides with the diagonal $x=y$.
The leaf $l^{s, \iin}_{i, \theta_s}$ thus has equation $x=2^{-\theta_s + \tilde \theta}$ and the leaf $l^{u, \iin}_{i, \theta_u}$ has equation $x=2^{-\theta_u - \tilde \theta}$ in coordinates $\rho_i = (x,\theta)$, where $\tilde \theta \in \R$ is a lift of $\theta \in \R/\Z$ (Remark~\ref{rmk: equation affine foliation}).
So the intersection points of these two leaves in coordinates $\rho_i = (x, \theta)$ are the following points, for $k \in \Z$, with $k$ large enough:
$$ p_k = \l(2^{\frac{\theta_u - \theta_s}{2} - \frac{k}{2}}, \frac{\theta_s + \theta_u}{2} + \frac{k}{2} \r).$$
Similarly, the intersection points of the leaves $l^{s, \out}_{j, \theta_s}$ and $l^{u, \out}_{j, \theta_u}$ in coordinates $\rho_j = (x, \theta)$ are the following points, for $k \in \Z$, with $k$ large enough:
$$ q_k = \l(- 2^{\frac{\theta_u - \theta_s}{2} - \frac{k}{2}}, \frac{\theta_s + \theta_u}{2} + \frac{k}{2} \r).$$
Therefore, for $x$ close enough to $0$ we have $\hat \varphi_1 (x, \theta) = \big(-x 2^{-\frac{k(x, \theta)}{2}}, \, \theta + \frac{k(x, \theta)}{2}\big)$ in coordinates $\rho_j$ and $\rho_i$, with $k(x, \theta) \in \Z$.
By continuity, $k(x, \theta)$ is constant, and there exists $k \in \Z$ such that for $x$ close enough to $0$, we have
$$ \hat \varphi_1 (x, \theta) = \l( -x 2^{-\frac{k}{2}}, \, \theta + \frac{k}{2} \r).$$

Let $\tau \colon [0,1] \to [0,1]$ be a smooth decreasing plateau function, such that $\tau(x)= 1$ on a neighborhood of $0$ and $\tau(x) = 0$ on a neighborhood of $1$.
Let $\psi \colon \pP \ssm \cO_* \to \pP \ssm \cO_*$ be the \diff{} with support in $\Pin$, equal to the identity outside $\bigcup_i \cV_i$, and such that the expression of $\psi$ in coordinates $\rho_i$ on each $\cV_i \cap \Pin$ is 
$$ (x, \theta) \mapsto \l( x 2^{\frac{k \tau(x)}{2}}, \theta- \frac{k \tau(x)}{2} \r). $$
Then $\psi$ is isotopic to the identity among the \diffs{} preserving $\cG^{u, \iin}$ leaf-to-leaf, supported in a neighborhood of $\partial \Pin$ where the foliations $\cG^{u, \iin}$ and $\varphi_* \cG^{u, \out}$ coincide (up to choose the support of $\tau$ close to $0$).
Therefore, the composition 
$\hat \varphi_1 := \psi \circ \hat \varphi_1' \colon \pP \ssm \cO_* \to \pP \ssm \cO_* $ 
is an incomplete gluing map strongly isotopic to $\hat \varphi_1$,
which maps the foliation $\cG^{u, \out}$ on a foliation transverse to $\cG^{s, \iin}$ on $\Pin$,
and for all $i$ and $j=\sigma(i)$, the expression of $\hat \varphi_1'$ in coordinates $\rho_j$ and $\rho_i$ is $\hat \varphi_1' (x, \theta) = (-x, \theta)$ for all $x$ close enough to $0$.
This \diff{} naturally extends on $\cO_*$ by $(0, \theta) \mapsto (0, \theta)$.
The \diff{} $\hat \varphi_1' \colon \pP \to \pP$ is thus a gluing map of $(P,X)$, whose restriction to $\pP \ssm \cO_*$ is strongly isotopic to $\hat \varphi$, hence satisfies Item~\ref{prop: normalization gluing map; it: strong isotopy} of Proposition~\ref{prop: extension trivial for incomplete gluing map}.
Recall the foliation $(\hat \varphi_1')_* \cG^{u, \iin}$ is transverse to $\cG^{s, \iin}$ on $\Pin$, hence $\hat \varphi_1'$ satisfies Item~\ref{prop: normalization gluing map; it: foliation} of Proposition~\ref{prop: extension trivial for incomplete gluing map}.
Finally, $\hat \varphi'_1$ satisfies Item~\ref{prop: normalization gluing map; it: trivial snc} of Proposition~\ref{prop: extension trivial for incomplete gluing map} by construction, in other words it is equal to the reflection $(x, x, \theta) \mapsto (-x,-x, \theta)$ in \ncss{} $(\cV_i', \xi_i = (x,y,\theta))$ and $(\cV_j', \xi_j = (x,y,\theta))$, for smaller normalized \nbhs{} $\cV_i' \subset \cV_i$ and $\cV_j' \subset \cV_j$.
The proposition is shown.
\end{proof}

It remains to prove Lemma~\ref{lem: identity cn psi in}, Lemma~\ref{lem: identity cn psi out} being completely analogous.

\begin{proof}[Proof of Lemma~\ref{lem: identity cn psi in}]
Let $\cV_i^\iin := \cV_i \cap \Pin$, and $\cV^\out_j := \cV_j \cap \Pout$, where $\cV_i$ and $\cV_j$ are tubular \nbhs{} of the periodic orbits $\cO_i$ and $\cO_j$ contained in $\pP$, with $j=\sigma(i)$.
It is sufficient to construct a \diff{} $\psi^\iin = \psi^\iin_i \colon \cV^\iin_i \to \cV^\iin_i$ which satisfies the properties stated on $\cV^\iin_i$, equal to the identity on the boundary of $\cV^\iin_i$.
We will work in the universal covering. 
Let $I = ]0,1]$ and $J=[-1, 0[$.
Denote $\pi \colon \R \times \R \to \R \times \R/\Z$ the projection.
Let $\pi^\iin_i \colon \widetilde \cV_i^\iin \to \cV^\iin_i$ be the universal covering of $\cV^\iin_i$ and $\tilde \rho_i =(x, \tilde \theta) \colon \widetilde \cV_i^\iin \to I \times \R$ the lift of coordinates $\rho_i$ such that $\pi \circ \tilde \rho_i = \rho_i \circ \pi^\iin_i$.
We also define $\widetilde \cV_j^\out$ the universal covering of $\cV_j^\out$ and $\tilde \rho_j = (x, \tilde \theta) \colon \widetilde \cV^\out_j \to J \times \R$ the lift of coordinates $\rho_j$.
Denote $(\tilde \cG^{s, \iin}$, $(\tilde \cG^{u, \iin})$ the lifts of the foliations in $\widetilde \cV^\iin_i$,
and $(\tilde \cG^{s, \out}$, $\tilde \cG^{u, \out})$ the lifts of the corresponding foliations in $\widetilde \cV^\out_i$.
For any $\tilde \theta \in \R$, denote $\tilde l^{s,\iin}_{\theta}$ and $\tilde l^{u, \iin}_{\theta}$ the leaves of $\tilde \cG^{s, \iin}$ and $\tilde \cG^{u, \iin}$ containing $(1, \tilde \theta)$ in the coordinate system $\tilde \rho_i$.
Similarly, for any $\tilde \theta \in \R$, denote $\tilde l^{s,\out}_{\theta}$ and $\tilde l^{u, \out}_{\theta}$ the leaves of $\tilde \cG^{s, \out}$ and $\tilde \cG^{u, \out}$ containing $(-1, \tilde \theta)$ in the coordinate system $\tilde \rho_j$.
The incomplete gluing map $\hat \varphi$ lift to a \diff{} $\tilde \varphi$, which is well defined from a neighborhood of $\{x=0\}$ in $\widetilde \cV^\out_j$ into a neighborhood of $\{x=0\}$ in $\widetilde \cV^\iin_i$.
For any $p \in \widetilde \cV^\iin_i \cap \tilde \varphi (\widetilde \cV^\out_j)$, let us denote $\tilde \theta_u(p)$ and $\tilde \theta_s(p)$ the real numbers such that $p  \in \tilde l^{s, \iin}_{\theta_s (p)} \cap \tilde \varphi (\tilde l^{s, \iin}_{\theta_s (p)})$.

\begin{lem} \label{lem: epsilon for intersection of leaves well defined}
There exists $\epsilon>0$ such that for any $p \in \widetilde \cV^\iin_i \cap \tilde \varphi (\widetilde \cV^\out_j)$ with coordinates $\tilde \rho_i(p) = (x(p), \tilde \theta(p))$ with $0 < x(p) < \epsilon$, the leaves $\tilde l^{s, \iin}_{\theta_s(p)}$ and $\tilde l^{u, \iin}_{\theta_u(p)}$ intersect at a single point $q \in \widetilde \cV_i$.
\end{lem}

\begin{proof}
Recall that, for every $\tilde \theta_s$ and $\tilde \theta_u$, the leaf $\tilde l^{s, \iin}_{\theta_s}$ has equation
$\tilde \theta = \log_2 (x) + \tilde \theta_s$, 
the leaf $\tilde l^{u, \iin}_{\theta_u}$ has equation
$\tilde \theta = -\log_2 (x) + \tilde \theta_u$, and 
the leaf $\tilde l^{u, \out}_{\theta_u}$ has equation
$\tilde \theta = -\log_2 (x) + \tilde \theta_u$ (in the corresponding coordinate systems $\rho_i$ and $\rho_j$).
Note that $\tilde l^{s, \iin}_{\theta_s}$ intersects $\tilde l^{u, \iin}_{\theta_u}$ if and only if $\tilde \theta_s < \tilde \theta_u$.
Orient each leaf of $\tilde l^{u, \out}_{\theta_u}$ in the increasing direction of the coordinate $\tilde \theta$ in $\rho_j$.
The leaf is then oriented in the direction of the lift of the orbit $\cO_j = \partial \cV^\out_j$.
Orient the leaf $\tilde \varphi_*(\tilde l^{u, \out}_{\theta_u})$ with the mapped orientation.
Each leaf $\tilde l^{s, \iin}_{\theta_s}$ of $\tilde \cG^{s, \iin}$ disconnects the plane $\tilde \cV^\iin_i$ into two components. 
We call \emph{left component of $\tilde l^{s, \iin}_{\theta_s}$} the component containing the leaves of $\tilde \cG^{s, \iin}$ parametrized by $\tilde \theta < \theta_s$, and \emph{right component of $\tilde l^{s, \iin}_{\theta_s}$} the other one.
With this orientation, the oriented leaves of $\tilde \cG^{u, \iin}$ intersect the leaves of $\tilde \cG^{s, \iin}$ from left to right.
Likewise,

\begin{claim} \label{claim: sense of intersection for foliation Gu}
The oriented leaves of $\tilde \varphi_* \tilde \cG^{u, \iin}$ transversely intersect the oriented leaves of $\tilde \cG^{s, \iin}$ from their left to their right on $\widetilde \cV^\iin_i$.
\end{claim}

\begin{proof}
By assumption $\hat \varphi$ is strongly isotopic to the restriction of a \sqt{} gluing map $\varphi_0$ of an isotopic \bb{} $(P_0, X_0)$.
We denote $\cL^\iin_0$ and $\cL^\out_0$ the boundary laminations on the entrance boundary $\Pin_0$ and the exit boundary $\Pout_0$ of $(P_0, X_0)$.
Let $\cO_{0,*} = \cO_{0,1}, \dots, \cO_{0,n}$ be the periodic orbits of $X_0$ contained in $\pP_0$, numbered compatibly with the numbering of the orbits $\cO_*$ (via block isotopy).
Then $\varphi_0$ maps the oriented orbit $\cO_{0,j}$ to the oriented orbit $\cO_{0,i}$, and maps $\Pout_0$ to $\Pin_0$.
It follows that the leaves of $(\varphi_0)_* \cL^\out_0$ accumulate on the orbit $\cO_{0,i}$ in a contracting way on $\Pin_0$, hence transversely intersect the leaves of $\cL^\iin_0$ from their left to their right, for the same orientation conventions.
By strong isotopy of $\hat \varphi$ with the restriction of $\varphi_0$ to $\pP_0 \ssm \cO_{0,*}$, it is the same for the intersection of oriented leaves of $\hat \varphi_* \cL^\out$ with oriented leaves of $\cL^\iin$ on $\Pin$.
By transversality of the foliations, it is the same for the leaves of $ \varphi_* \cG^{u, \iin}$ with $\cG^{s, \iin}$.
This property remains true for the universal covering.
\end{proof}

We don't know a priori if the leaf $\tilde l^{u, \iin}_{\theta_u(p)}$ intersects the leaf $\tilde l^{s, \iin}_{\theta_s(p)}$, but we have the following claim.

\begin{claim} \label{claim: convergence theta u - theta s}
If $x(p) \to 0$, then $(\tilde \theta_s(p) - \tilde \theta_u(p)) \to + \infty$ and this convergence is uniform in $\tilde \theta(p)$.
\end{claim}

\begin{proof}
According to Claim~\ref{claim: sense of intersection for foliation Gu}, if $p = (x(p), \tilde \theta(p))$ is a coordinate for the oriented leaf $\tilde \varphi (\tilde l^{u, \iin}_{\theta_u(p)})$, in other words $\tilde \theta_u(p)$ is fixed and $x(p) \to 0$, then $\tilde \theta_s (p) \to + \infty$.
By continuity (the foliations $\tilde \varphi_* \tilde \cG^{u, \iin}$ and $\tilde \cG^{s, \iin}$ are transverse), the convergence is locally uniform in $\tilde \theta(p)$.
Since everything commutes with the translation $\tilde \theta \mapsto \tilde \theta+1$, we deduce that the convergence is globally uniform.
\end{proof}

It follows from Claim~\ref{claim: convergence theta u - theta s} the existence of $\epsilon >0$ such that for all $p = (x, \theta) \in \widetilde \cV^\iin_i$ with $x < \epsilon$, we have $\tilde \theta_s(p) < \tilde \theta_u(p)$.
As noticed previously, the explicit equations of the leaves $\tilde l^{u, \iin}$ and $\tilde l^{s, \iin}$ imply that for all $p = (x, \tilde \theta) \in \widetilde\cV^\iin_i$ with $x < \epsilon$, the leaves
$\tilde l^{u, \iin}_{\theta_u(p)}$ and $\tilde l^{s, \iin}_{\theta_s(p)}$ intersect at a single point
\[ q(p) = \l( 2^{ \frac{\tilde \theta_u(p) - \tilde \theta_s(p)}{2} }, \frac{\tilde \theta_u(p) + \tilde \theta_s(p)}{2} \r).  \qedhere\]
\end{proof}

We can now prove Lemma~\ref{lem: identity cn psi in}.
Lemma~\ref{lem: epsilon for intersection of leaves well defined} gives a map
$p \mapsto q(p)$, such that, given a point $p \in \cV^\iin_i$ such that $x(p) < \epsilon$, it associates the point $q(p)$ of intersection of the leaf $l^{s, \iin}_{\theta_s(p)}$ with the leaf $\tilde l^{u, \iin}_{\theta_u(p)}$.
This map is well defined from $]0, \epsilon[ \times \R$ (in the coordinate system $\rho_i = (x, \theta)$) on its image $\cV$ where $\cV \subset \cV^\iin$ contains a neighborhood $]0, \eta] \times \R$ ($\eta < \epsilon$), and injective.
By smooth transversality of the foliations, it is differentiable.
Let $\tau \colon [0,1] \to [0,1]$ be a smooth decreasing function which is equal to $1$ on $[0, \frac{\eta}{2}]$ and $0$ on $[\eta, 1]$.
Define $\tilde \psi^\iin \colon \tilde \cV^\iin_i \to \cV^\iin_i$ in the following way (Figure~\ref{fig: straightening foliation image}).
\begin{figure}[htb]
    \centering
    \vspace*{-2em}
    \includegraphics[height=0.28\textheight]{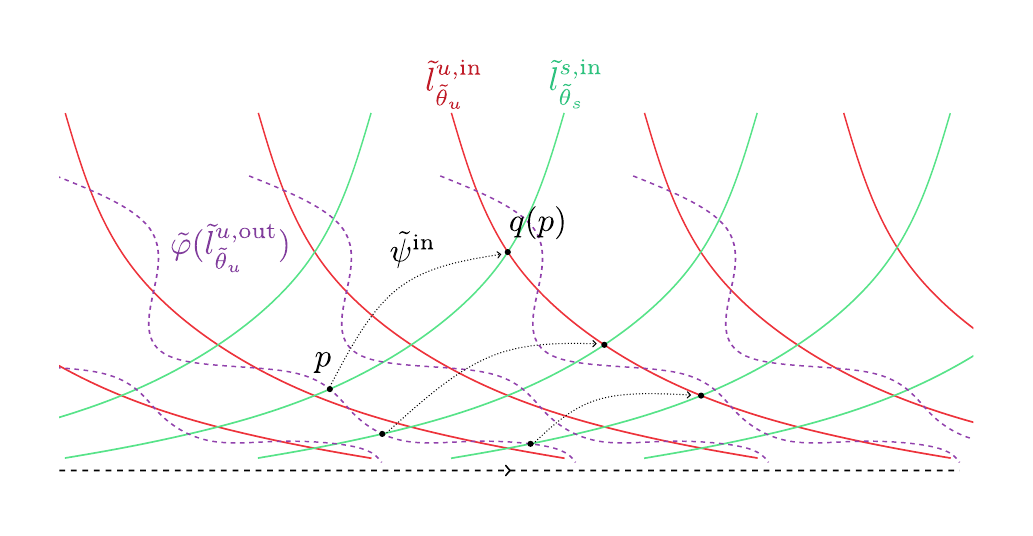}
    \vspace*{-1em}
    \centering \caption{Action of $\tilde \psi^\iin$ on $\tilde \cV^\iin_i$}
    \label{fig: straightening foliation image}
\end{figure}

Let $p= (x, \tilde \theta) \in \widetilde \cV^\iin_i$.
If $x \geq \eta$, define $\tilde \psi^\iin(x, \tilde \theta) = (x, \tilde \theta)$.
If $x \leq \eta$, let $q = q(p) \in \cV^\iin$ given by Lemma~\ref{lem: epsilon for intersection of leaves well defined}.
    Then $p$ and $q$ belong to the same leaf $l^{s, \iin}_{\theta_s}$ of $\tilde \cG^{s, \iin}$, and define $\tilde \psi^\iin(p)$ the weighted barycenter by $\tau(x)$ between the points $p$ and $q$ in the leaf $l^{s, \iin}_{\theta_s}$, which we denote by use of notation
    $$\tilde \psi^\iin (p = (x, \tilde \theta) ) = \tau(x) q + (1- \tau(x)) p .$$
This expression define a \diff{} $\tilde \psi^\iin_i$, equal to the identity in the neighborhood of the boundary $\{ x=1 \}$ of $\tilde \cV^\iin_i$,
isotopic to the identity among \diffs{} which preserves the foliation $\widetilde \cG^{s, \iin}$ leaf-to-leaf,
and which maps the leaf $\tilde \varphi(\tilde l^{u, \out}_{\theta_u})$ to the leaf $\tilde l^{u, \iin}_{\theta_u}$ on a neighborhood of the boundary $\{ x=0 \}$ of $\tilde \cV^\iin_i$.
This \diff{} commutes with the integer translation along the coordinate $\tilde \theta$, so it is quotiented into a \diff{} $\psi^\iin \colon \cV^\iin_i \to \cV^\iin_i$ which satisfies the desired properties.
\end{proof}

\subsubsection*{Conclusion}
\begin{proof}[Proof of Proposition~\ref{prop: strongly isotopic normalized gluing map}]
It is sufficient to apply successively Propositions \ref{prop: incomplete gluing map with transverse foliation} and~\ref{prop: extension trivial for incomplete gluing map}
\end{proof}

\subsection{Proof of Proposition~\ref{prop: normalization of triple}}
\label{sec: normalization; subsec: proof of triple normalization}

Let $(P_0,X_0)$ be a filled \bb{}, and $\varphi_0$ be a \sqt{} gluing map of $(P_0,X_0)$.
\begin{itemize}[leftmargin=*]
    \item Propositions~\ref{prop: affine section} and~\ref{prop: multipliers} applied successively to $(P_0,X_0)$ give the existence of a \bb{} $(P, X)$ which satisfies Items~\ref{def: normalized, it: affine section} and~\ref{def: normalized, it: multipliers} of the definition of a normalized block~\ref{def: normalized block}, in other words it admits an affine section and the multipliers of the orbits of $X$ in $\pP$ are equal to $\{\frac{1}{2}, 2\}$,
    and it is isotopic to $(P_0,X_0)$ among the orbit equivalent \bbs{}.
    The block $(P,X)$ satisfies the assumptions of Proposition~\ref{prop: existence paif}.
    Let $(\cG^s, \cG^u)$ be a \paif{} on $(P, X)$ given by Proposition~\ref{prop: existence paif}.
    The block $(P, X)$ equipped with the pair of foliations $(\cG^s, \cG^u)$ satisfies Items~\ref{def: normalized, it: affine section},~\ref{def: normalized, it: multipliers} and~\ref{def: normalized, it: foliation} of Definition~\ref{def: normalized block}.
    It follows from Proposition~\ref{prop: induced gluing map on orbit equivalent block} the existence of a strongly transverse incomplete gluing map $\hat \varphi$ of $(P, X)$, strongly isotopic to the restriction of the gluing map $\varphi_0$ to the complementary of the periodic orbits of $X_0$ in~$\pP_0$.
    
    \item Proposition~\ref{prop: straight block associated} applied to $(P, X)$ gives the existence of a filled normalized block $(P_1, X_1)$ isotopic to $(P, X)$.
    Lemma~\ref{lem: induced incomplete gluing map on isotopic blocks} applied to the incomplete gluing map $\hat \varphi$ of $(P, X)$ gives the existence of a strongly transverse incomplete gluing map $\hat \varphi_1$ of $(P_1, X_1)$ and strongly isotopic to $\hat \varphi$ and thus to the restriction of the \sqt{} gluing map $\varphi_0$ to the complementary of the periodic orbits of $X_0$ in $\pP_0$. 

    \item The block $(P_1, X_1)$ is a filled normalized block with a strongly transverse incomplete gluing map $\hat \varphi_1$, strongly isotopic to the restriction of the \sqt{} gluing map $\varphi_0$ of $(P_0,X_0)$.
    One can apply Proposition~\ref{prop: strongly isotopic normalized gluing map}, which gives the existence of a normalized gluing map $\varphi'_1$
    of $(P_1, X_1)$, whose restriction is strongly isotopic to $\hat \varphi_1$ and thus to the restriction of $\varphi_0$.
    To show that the triple $(P_1, X_1, \varphi'_1)$ is strongly isotopic to the triple $(P_0, X_0, \varphi_0)$, it remains to show that $\varphi'_1$ is isotopic to $\varphi_0$ (Item~\ref{def: strongly isotopic triples; it: gluing maps isotopy} of Definition~\ref{def: strongly isotopic triples}).
    The blocks $(P_0, X_0)$ and $(P_1, X_1)$ are isotopic.
    Up to make an orbit equivalence $(P_0', X_0')$ of $(P_0,X_0)$, they are contained in a common \minc{} (Proposition~\ref{prop: isotopy vs orbit eq}), and the boundaries $\pP'_0$ and $\pP_1$ are isotopic relatively to the periodic orbits contained in the boundary.
    Let $\varphi_0'$ be the conjugate of $\varphi_0$ by the orbit equivalence and $\varphi_1$ be the push of $\varphi_0'$ on $\pP_1$ by the isotopy between $\pP_0'$ and $\pP_1$.
    Then $\varphi_1$ is a gluing map, a priori only topological, that we will use only to adjust the isotopy class of $\varphi_1'$.
    Up to make a small perturbation, we can assume that it is differentiable.
    It is isotopic to $\varphi_0$ in the sense of Definition~\ref{def: strongly isotopic triples}.
    
    Let us show that we can adjust $\varphi'_1$ to make it isotopic 
    to $\varphi_1$ without destroying the strong isotopy.
    If $T$ is a \cc{} of $\pP$ that does not intersect the collection $\cO_*$, the restriction of $\hat \varphi_1$ to $T$ is a complete gluing map.
    The restrictions to $T$ of $\hat \varphi_1$ and $\varphi_0$, and of $\varphi'_1$ and $\hat \varphi_1$ are then isotopic.
    The result follows by transitivity.
Suppose $T$ contains a collection of periodic orbits of $X$, and let $T'$ be the \cc{} of $\pP$ matched with $T$ by $\varphi_1$.
Let $\cO_i$ be a periodic orbit contained in $T$, and $(\cV_i, \xi_i= (x,y,\theta))$ an \ncs{} of $\cO_i$.
As $\varphi'_1$ is a normalized gluing map, up to shrink $\cV_i$, we can assume that
$(\varphi'_1)_* \cG^{u, \out} = \cG^{u, \iin}$ and $(\varphi'_1)_* \cG^{s, \iin} = \cG^{s, \out}$ on $\pP \cap \cV_i$.
Let $A := \pP \cap \cV_i$ and $\rho_i = (x,\theta) \in I \times \R/\Z$ the coordinates induced on $A$.
There exists $\tau \colon T \to T$ a Dehn twist on the annulus $A$, equal to the identity on the boundary of $A$, and such that $\tau \circ \varphi'_1 \colon T' \to T$ is isotopic to $\varphi_1 \colon T' \to T$.
For the composition $\varphi''_1:= \tau \circ \varphi'_1$ to still be a normalized gluing map, it is necessary and sufficient that:
\begin{enumerate}
    \item $\tau$ is equal to the identity in the neighborhood of $\cO_i$,
    \label{it: twist id}
    \item $\tau$ maps the foliation $\cG^{u, \iin}$ to a foliation transverse to $\cG^{s, \iin}$,
    \label{it: twist in}
    \item $\tau$ maps the foliation $\cG^{s, \out}$ to a foliation transverse to $\cG^{u, \out}$.
    \label{it: twist out}
\end{enumerate}
Let $\sigma := \{\theta = 0\}$ and $n$ be the integer such that $\tau([\sigma]) = [\sigma] + n [\cO]$, where $[\cdot]$ denotes the homology class of the curve.
We have $\tau (x, \theta) = (x, \theta + f(x) \mod \Z)$, with $f \colon [-1, 1] \to \R$ a smooth monotone function, equal to $O$ on a neighborhood of $-1$, and equal to $n$ on a neighborhood of $1$.
According to Remark~\ref{rmk: equation affine foliation} there exists a constant $c>0$ such that, if $^*$ denotes $\iin$ or $\out$, the tangent directions to the foliations $(\cG^{s, *}, \cG^{u, *})$ at a point $p=(x, \theta) \in A$ are
\vspace*{-0.3em}
\[ T_{(x, \theta)} \cG^{s, *} \colon dx = cx . d\theta, \qq T _{(x, \theta)}\cG^{u, *}  \colon dx = - cx . d\theta.\]
\begin{itemize}
    \item If $n \geq 0$, we choose the support of $\tau$ in a compact set of $A \cap \Pout$.
    The direction tangent to $\cG^{s, \out}$ at the point $p=(x, \theta)$ is mapped on the direction $(cx.f'(x) + 1) dx = cx. d\theta$.
    Since $f'(x)\geq 0$ and $x \geq 0$, the direction of the image is transverse to the direction $dx = - cx . d\theta$ tangent to $\cG^{u, \out}$ at the point $\tau(p) = (x, \theta + f(x))$, which shows that Items~\ref{it: twist id} and~\ref{it: twist out} are true.
    Item~\ref{it: twist in} is trivially true because $\tau$ is the identity on $\Pin$.
    \item If $n \leq 0$, we choose the support of $\tau$ in a compact set of $A \cap \Pin$.
    The direction tangent to $\cG^{u, \iin}$ at the point $p=(x, \theta)$ is mapped on the direction $(-cx.f'(x)+1) . dx = -cx . d\theta$.
    Since $f'(x)\leq 0$ and $x \leq 0$, the direction of the image is transverse to the direction $dx = cx . d\theta$ tangent to $\cG^{s, \iin}$ at the point $\tau(p) = (x, \theta + f(x))$, which shows that Items~\ref{it: twist id} and~\ref{it: twist in} are true.
    Item~\ref{it: twist out} is trivially true because $\tau$ is the identity on $\Pout$.
\end{itemize}
    It follows that the triple $(P_1, X_1, \varphi''_1)$ is a normalized triple strongly isotopic to the triple $(P_0,X_0,\varphi_0)$.
\end{itemize}

\vspace*{-0.7em}

\section{Crossing map}
\label{sec: crossing map}

In this section $(P,X)$ denotes a normalized filled \bb{}.
We refer to Sections \ref{sec: preliminaries} and~\ref{sec: normalization} for notations, definitions and associated properties.
By definition of a normalized \bb{} (Definition~\ref{def: normalized block}), each orbit $\cO_i \in \cO_*$ \linebreak[4]contained in the boundary is provided with a \ncs{} which we denote $(\cV_i, \xi_i = (x,y,\theta))$.
For the rest of the proof, we fix $g= \langle \cdot, \cdot\rangle$ a Riemannian metric on a \minc{} $\tilde P$,
adapted to the hyperbolic splitting of the maximal invariant set $\Lambda$, and such that in a \ncs{} $(\cV_i, \xi_i = (x, y, \theta))$ of $\cO_i$, we have
$g = dx^2 + dy^2 + d \theta^2 $.

\subsection{Definition and statement of the main proposition}

\label{sec: crossing map; subsec: def and statement}

In this section, we will study some properties of the crossing map $\foutin$ of the flow of $X$ from $\Pin{}$ to $\Pout{}$.
For any $p \in P$, let $\tau(p)$ be the \emph{exit time} (when well defined) of the orbit of $p$ by the flow of $X$, in other words such that $X^{\tau(p)}(p) \in \Pout{}$.

\begin{defi}[Crossing map] \label{def: def foutin}
The map $\foutin{} \colon \Pin \to \Pout$ defined by $$\foutin(p) = X^{\tau(p)}(p)$$
is a \diff{} from $\Pin \ssm \cL^\iin$ to $\Pout \ssm \cL^\out$, called \emph{the crossing map of $(P,X)$}.
\end{defi}

Indeed, the time $\tau$ is well defined if and only if $p$ is in the complementary of the stable manifold $\cW^s$ of $\Lambda$.
The image of $\foutin$ corresponds to the set of points of $\Pout$ whose orbit intersects $\Pin$, and thus coincides with the complementary of the unstable manifold $\cW^u$ of $\Lambda$ in $\Pout$.
A crossing map of a $\cC^1$-flow between (smooth) surfaces quasi-transverse to the vector field is~$\cC^1$.




\smallskip

The purpose of this section is to prove the following main proposition, which states that \foutin{} expands the direction tangent to $\cG^{u, \iin}$ arbitrarily strongly in a small neighborhood of $\cL^\iin$, and the analogous result for the inverse map along the direction tangent to $\cG^{s, \out}$.

\begin{prop}[Expansion of the crossing map] 
\label{prop: crossing map} \mb{}
\begin{enumerate}
    \item For any $\lambda > 1$, there exists a neighborhood $\cW^\iin_\lambda$ of $\cL^\iin$ in \Pin{},
    such that the crossing map $\foutin{}$ expands by a factor $\lambda$ the direction tangent to $\cG^{u, \iin}$ on $\cW^\iin_\lambda \ssm \cL^\iin$:
    \vspace*{-0.3em}
    $$\forall p \in \cW^\iin_\lambda \ssm \cL^\iin, \quad \forall v \in T_p \cG^{u, \iin}, \quad \Vert (\foutin)_* v \Vert > \lambda \Vert v \Vert.$$
    \item For any $\lambda > 1$, there exists a neighborhood $\cW^\out_\lambda$ of $\cL^\out$ in \Pout{},
    such that the inverse $\foutin^\inv$ expands by a factor $\lambda$ the direction tangent to $\cG^{s, \out}$ on $\cW^\out_\lambda \ssm \cL^\out$:
    \vspace*{-0.3em}
    $$\forall p \in \cW^\out_\lambda \ssm \cL^\out, \quad \forall v \in T_p \cG^{s, \out}, \quad \Vert (\foutin^\inv)_* v \Vert > \lambda \Vert v \Vert.$$
\end{enumerate}
\end{prop}

\subsubsection*{Remarks}
\begin{figure}[htb]
    \centering
    \vspace*{-1em}
    \includegraphics[height=0.4\textheight]{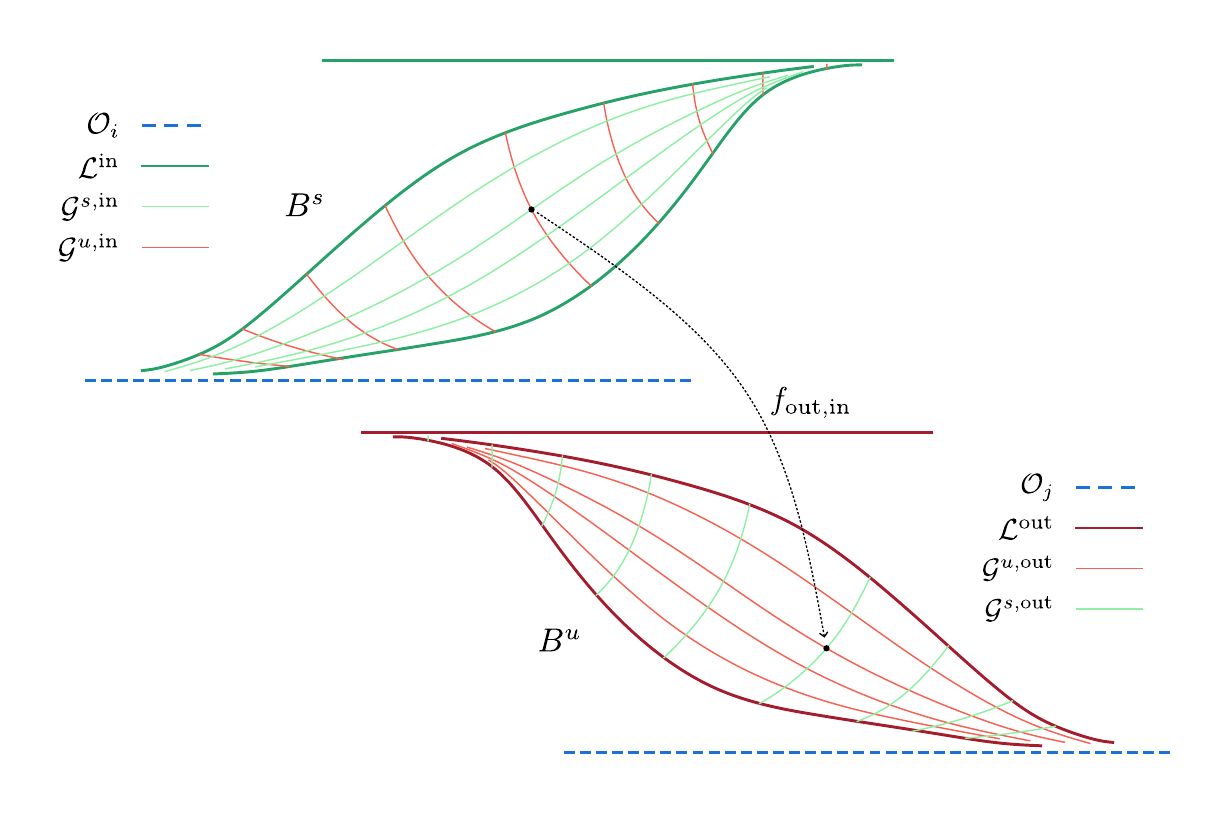}
    \vspace*{-1em}
    \caption{Action of the crossing map from a strip $B^s \subset \Pin$ to a strip $B^u \subset \Pout$}
    \label{fig: crossing map strip to strip}
\end{figure}

This result is valid for any normalized \bb{}, not necessary filled. 
However, in order to have pictures in mind and because it corresponds to our setting, let us describe the crossing map for a filled \bb{}.
The boundary lamination $\cL$ is filling hence the connected components of $\Pin \ssm \cL^\iin$ and $\Pout \cL^\out$ are \emph{strips}, i.e., topological disks whose accessible boundary is formed by two non-compact leaves asymptotic to each other at both end (Definition~\ref{def: strip and filling lam}).
The map $\foutin{}$ maps a strip $B^s \subset \Pin$ to a strip $B^u = \foutin (B^s) \subset \Pout{}$ (Figure~\ref{fig: crossing map strip to strip}).

If $B^s$ is a strip of $\Pin$, there are two possible cases for each end of $B^s$:
    \begin{enumerate}
        \item either the given end of $B^s$ accumulates on a compact leaf of $\cL^\iin$ in $\Pin$, or 
        \item the given end of $B$ accumulates on a component of $\partial \Pin$, in other words a periodic orbit $\cO_i \in \cO_*$ for the flow of $X$.
    \end{enumerate}

    Let $\cW^\iin_\lambda$ be the neighborhood given by Proposition~\ref{prop: crossing map}.

    If $B^s$ is a \cc{} of $\Pin \ssm \cL^\iin$, then $\cW^\iin_\lambda \cap B^s$ is a neighborhood of the accessible boundary of $B^s$, and
    \begin{enumerate}
        \item either we are in the previous case 1), and $\cW^\iin_\lambda \cap B^s$ contains a neighborhood of the corresponding end of $B^s$, or
        \item we are in the previous case 2), and $\cW^\iin_\lambda \cap B^s$ does not contain any neighborhood of the corresponding end of $B^s$.
    \end{enumerate}
    Up to shrinkink $\cW^\iin_\lambda$, we can suppose it is an adapted \nbh{} of $\cL^\iin$ (Definition~\ref{def: adapted nbh}, Item~\ref{def: adapted nbh; it: of L^in}) and we can refer to Figure~\ref{fig: adapted square half band}.

\subsubsection*{Proof summary}

Let us give the main ideas for the proof of the first item of Proposition~\ref{prop: crossing map}. The second item is symmetrical.

\begin{itemize}[leftmargin=*]
    \item We can extend the hyperbolic structure of $\Lambda^s$ on a small filtrating \nbh{} $\cU$ of $\Lambda^s$, in the sense that there are directions $E^\ss$ and $E^\uu$ in $\cU$, invariant by the differential of $X^t$ in $\cU$, and which are respectively contracted and expanded in the future (Lemma~\ref{lem: filtrating neighborhood hyp}).
    According to the $\lambda$-lemma (\cite{katokIntroductionModernTheory1995}), a vector transverse to the plane $E^\ss \oplus \R. X$ carried by the flow for large enough positive times is exponentially expanded and tends to the strong unstable direction $E^\uu$.
    
    \item A point in $\Pin{}$ close to the lamination $\cL^\iin$ in a \cc{} of $\Pin \ssm \cL^\iin$ is close to the stable manifold of the maximal saddle invariant set (Lemma~\ref{lem: boundary leaf and periodic orbit are saddle}) thus has an orbit which spends a long time in the neighborhood $\cU$ before intersecting $\Pout{}$.
    It follows that a vector $v$ on a point close to the lamination $\cL^\iin$ and far enough from the orbits $\cO_*$, tangent to $\cG^{u, \iin}$, is expanded by the differential of $X^t$ after a uniformly bounded time.
    
    \item The image of a vector $v$ by the differential of the crossing map $\foutin{}$ is the projection $\pi \colon TP \to TP^\out$ onto the tangent space $T \Pout{}$ parallel to $\R.X$ of the vector $X^\tau_*v$ pushed by the flow on \Pout{}, i.e., $
       (\foutin)_*v = \pi (X^\tau_*v)$
    where $\tau$ is the exit time of the orbit.
    We show that the projection of a vector $v = v^\uu + v^X + v^\ss \in E^\uu \oplus \R. X \oplus E^\ss$ on a surface $S$ on $\cU$ parallel to $\R.X$ does not contract the strong unstable component $v^\uu$ (Lemma~\ref{lem: projection bounded effect on strong unstable}).
    \item The tricky point occurs in the neighborhood of an orbit $\cO_i \in \cO_*$.
    Indeed a vector $v$ tangent to $\cG^{u, \iin}$ at a point $p \in \Pin$ very close to an orbit $\cO_i$ has a strong unstable component $v^\uu$ very small compared to its component $v^X$ tangent to $\R. X$, on which the action of the flow is isometric.
    It is therefore necessary to spend a very long time in a hyperbolic neighborhood for the vector to be exponentially expanded by the $X^t$ flow.
    However, the trajectory of an orbit of a point $p \in \Pin$ very close to an orbit $\cO_i \in \cO_*$ spends a very long time in a neighborhood of $\cO_i$, where the dynamics of the flow is hyperbolic.
    We will see that this competition balances.
    To do so, we decompose the crossing map \foutin{} from \Pin{} to an \say{intermediate} surface \tPin{} close to \Pin{} and uniformly transverse to the vector field, then from an \say{intermediate} surface \tPout{} close to \Pout{} and uniformly transverse to the vector field.
    We study the intermediate maps explicitly in a \ncs{} $(\cV_i, \xi_i)$ of $\cO_i$.
\end{itemize}

\subsection{Preliminary lemmas}
Before proving the proposition, we will need a number of lemmas which we now state and prove.
\label{sec: crossing map; subsec: preliminary lemmas}

\subsubsection*{Hyperbolic filtrating neighborhood}
A neighborhood of $\Lambda_s$ in $P$ is said to be \emph{filtrating} if the intersection of this neighborhood with any orbit of the flow of $X$ in $P$ is connected.
The hyperbolic structure of a locally maximal invariant set extend on small filtrating \nbh{}.
More precisely we have:

\begin{lem}
\label{lem: filtrating neighborhood hyp}
There exists a filtrating neighborhood $\cU$ of $\Lambda_s$ in $P$,
and an $X^t$-invariant hyperbolic splitting $\cU$, still denoted $\res{TP}{\cU} = E^\uu \oplus \R.X \oplus E^\ss$, which coincides with the hyperbolic splitting on $\Lambda_s$.
There exist constants $\lambda_0 > 1$ and $C_0 >0$ such that if the orbit of $p \in \cU$ is in $\cU$ at time $t>0$, then
\begin{itemize}[--]
    \item $ \forall v \in E^\uu(p)$, $\Vert (X^t)_* v \Vert \geq C_0 \lambda_0^t v \Vert$;
    \item $\forall v \in E^\ss(p)$, $\Vert (X^t)_* v \Vert \leq C_0^\inv \lambda_0^{-t} \Vert v \Vert$.
\end{itemize}
Moreover, the direction $E^\ss \oplus \R.X$ is tangent to $\cG^s$ and the direction $E^\uu \oplus \R.X$ is tangent to $\cG^u$.
\end{lem}

\begin{proof}
By Lemma~\ref{lem: existence of markov partition for saddle block}, there exists a Markov partition $\cR = (R, \Sigma)$ of $\Lambda_s$ on $\Sigma$, such that the intersection of any orbit of $X$ in $P$ with the suspension \nbh{} $ \cU := \bigcup_{t \in [0,1]} X^t (R) $ is connected (recall that the return time of $R$ on $\Sigma$ is equal to $1$).
Let $(\zeta^s, \zeta^u)$ be the pair of $f$-invariant foliations on $\Sigma$, such that $f$ uniformly expands the direction tangent to $\zeta^u$ and uniformly contracts the direction tangent to $\zeta^s$.
Let us show that the neighborhood $\cU$ satisfies Lemma~\ref{lem: filtrating neighborhood hyp}.
According to the properties of the Markov partition, it is a filtrating neighborhood of $\Lambda_s$.
We can push the decomposition $T\Sigma = T \zeta^s \oplus T \zeta^u$ on $\cU$ by the differential of the flow for $t \in [0,1]$.
By invariance and hyperbolicity of $f$, they are invariant and respectively contracted and expanded in the future by the flow, and extend the hyperbolic structure of $\Lambda^s$ on $\cU$. \end{proof}

We will need the following result.

\begin{claim} \label{claim: orbit length bounded outside nbh U}
    The segments of $X$-orbit from $\Pin$ to $\Pout$ which never enters $\cU$ have uniformly bounded length.
    The segment of $X$-orbit from $\Pin$ to $\cU$ and from $\cU$ to $\Pout$ have uniformly bounded length.
\end{claim}

\begin{proof}
Let us show the first assertion.
       If this were not the case, we would have a sequence of points $p_n \in \Pin$ and increasing times $t_n \to +\infty$ such that the orbit segments $[p_n, X^{t_n}p_n]$ are in the complementary of $\cU$ and $X^{t_n}p_n \in \Pout$.
    Up to take a subsequence the sequence $(p_n)_n$ converges to a point $q$ in the closure $\overline \Pin$. 
    Suppose that there exists a time $T>0$ such that $X^T(q)$ is in the interior of $\cU$.
By continuity of flow, for an arbitrarily large $n$, the distance $\operatorname{dist}( X^T(p_n), X^T(q)) $ is arbitrarily small, and $X^T(p_n)$ is in $\cU$ for $n$ large enough, which is a contradiction as soon as $t_n > T$.
Suppose now that there exists a time $T>0$ such that $X^T(q)$ crosses $\Pout$.
The distance $\operatorname{dist}( X^T(p_n), X^T(q))$ is arbitrarily small, but a point in $P$ arbitrarily close to $\Pout{}$ has a positive orbit which exits $P$ in arbitrarily small time, which is a contradiction for $t_n$ big in front of $T$.
Suppose now that the positive orbit of $q$ accumulate on an attractor in $\gA$. Hence $q$ belongs to the basin of the attractor, which is an open set, and so do $p_n$ for $n$ large enough, which is a contradiction.  
We deduce that the positive orbit of $q$ neither enters $\cU$, neither exits through $\Pout{}$, and neither accumulates on an attractor.
This is impossible.
The other assertions are shown similarly.
\end{proof}

\begin{rmk} \label{rmk: extend hyperbolic splitting}
    We can push the hyperbolic splitting of $\cU$ from Lemma~\ref{lem: filtrating neighborhood hyp} in the future along the orbits of $\cU$ exiting through $\Pout$ and in the past along the orbits of $\cU$ entering through $\Pin$ because the \nbh{} $\cU$ is filtrating.
    Moreover, the hyperbolic properties are still satisfied along these orbits segments according to Claim~\ref{claim: orbit length bounded outside nbh U}.
\end{rmk}

\subsubsection*{Projection of a center-unstable vector onto a transverse surface}
We will need the following lemma, which states that the projection parallel to $\R.X$ of a vector $v \in E^\uu \oplus \R.X$ onto a surface $S$ transverse to $X$, at a point $p \in \cU$ has a norm essentially less than the norm of the strong unstable component in $E^\uu$ of the vector $v$.

\begin{lem}
\label{lem: projection bounded effect on strong unstable}
Let $S$ be a surface embedded in $P$ transverse to the vector field $X$.
Let $p\in S \cap \cU$ and $\pi_p \colon T_p P \to T_p S$ be the projection on $T_pS$ parallel to $\R. X$ in the neighborhood $\cU$ provided with the hyperbolic splitting.
There exists a constant $\cst>0$ such that for any vector $v = v^\uu + v^X \in E^\uu \oplus \R.X$ tangent to the center unstable bundle on a point $p \in S$, we have
$\Vert \pi_p(v) \Vert \geq \cst \Vert v^\uu \Vert. $
\end{lem}

\begin{proof}
Recall the following result.
\begin{claim} \label{claim: projection bounded by angle}
Let $F,G$ be two vector planes in $\R^3$, and $H$ be a supplementary line of $F$ and $G$, and $\Vert \cdot \Vert_2$ the Euclidean norm.
Let $ \pi \colon F \to G$ be the projection of $F$ onto $G$ parallel to $H$, linear invertible. Then we have
$$
\dfrac{ \Vert \pi(v) \Vert_2 }{ \Vert v \Vert_2 }
=\l\vert \dfrac{ \sin \measuredangle(H,\R.v) }{ \sin \measuredangle(H,\R.\pi(v)) } \r\vert.
$$
\end{claim}
If $\Vert \cdot \Vert_2$ is the Euclidean norm associated to a basis of tangent vectors $(\partial_\uu$, $\partial_\ss, X)$ where $\partial_\uu \in E^\uu$ and $\partial_\ss \in E^\ss$,
then we have
$$ \frac{\Vert \pi (v) \Vert_2}{\Vert v \Vert_2} \geq \sin( \measuredangle(v,X) ) = \frac{\Vert v^\uu \Vert_2}{\Vert v \Vert_2}.$$
By equivalence of norms, we deduce the existence of the constant $\cst>0$ which satisfies the lemma.
\end{proof}

\subsubsection*{Angles between the foliations and vector field on $\partial P$}

The vector field $X$ becomes more and more tangent to the surface $\pP$ as we approach the set of periodic orbits $\cO_*$.
Similarly the $\cG^{s, \iin}$ and $\cG^{u, \iin}$ laminations on $\Pin{}$ are increasingly tangent to each other as one approaches the set $\cO_* \subset \pP$.
Using the expression for the vector field $X$ and the foliations in the normalized coordinates $(\cV_i, \xi_i)$, we show that the closing velocity of the angles is essentially the same, of the order of the distance to the orbit $\cO_i$.
This property will be useful in the rest of the proof.

\begin{lem}\label{lem: angle foliation and vector field on boundary}
There exists a constant $\cst>0$ such that for all $p \in \Pin{}$,
\begin{enumerate}
    \item \label{angle gu gs on Pin} $ \cst^\inv \dist(p, \cO_*) \leq \measuredangle(\cG^{s, \iin}, \cG^{u, \iin})_p \leq \cst \dist(p, \cO_*) $,
    \item \label{angle X Pin on Pin} $ \cst^\inv \dist(p, \cO_*) \leq \measuredangle(X, \Pin{})_p \leq \cst \dist(p, \cO_*) $,
    \item \label{angle gu X on Pin} $ \cst^\inv \dist(p, \cO_*) \leq \measuredangle(X, \cG^{u, \iin})_p \leq \cst \dist(p, \cO_*) $.
\end{enumerate}
\end{lem}

\begin{proof}
It is sufficient to check this in the neighborhood of a periodic orbit $\cO_i$ in the $(\cV_i, \xi_i = (x,y,\theta))$.
Let $(\partial_x, \partial_y, \partial_\theta)$ be the basis of vector fields associated with the coordinates $(x,y,\theta)$.
Recall that the norm $\Vert \cdot \Vert$ on $P$ coincides with the Euclidean norm in these coordinates.
By definition of a normalized block (Definition~\ref{def: normalized block} and Corollary~\ref{coro: expression flow in normalized neighborhood}), the surface $\Pin{}$ has local equation $x=y$ and, if $c := \ln(2)$, the vector field $X$ in these coordinates is
$X(x,y,\theta) = -c x \partial_x + c y \partial_y + \partial_\theta. $
The neighborhoods $\cV_i$ are contained in the neighborhood $\cU$ provided with a hyperbolic splitting of the tangent space $\res{TP}{\cU} = E^\uu \oplus \R.X \oplus E^\ss$ (Lemma~\ref{lem: filtrating neighborhood hyp}), up to restrict $\cV_*$, which we will always suppose.
The direction $E^\uu$ is generated by $\partial_y$ and the direction $E^\ss$ is generated by $\partial_x$.
Note that the distance from a point $p=(x,x,\theta) \in \Pin{}$ to $\cO_i$ is equal to $x$.
We show each of the three items.
\begin{enumerate}[leftmargin=*]
    \item The foliation $\cG^{s, \iin}$ is generated by the vector field $\partial_s = c x \, (\partial_x + \partial_y) + \partial_\theta$ and the foliation $\cG^{u, \iin}$ is generated by the vector field $\partial_u = - c x \, (\partial_x + \partial_y) + \partial_\theta$ \linebreak[4](Figure~\ref{fig: angle foliation Pin}). As the couple $(\partial_x + \partial_y, \partial_\theta)$ forms an orthogonal basis on $\Pin{} \cap \cV_i$, the angle between the two tangent directions at a point $p= (x,x,\theta) \in \Pin{}$ is proportional to $x$, hence to the distance from $p$ to $\cO_i$.
    
    \begin{figure}[t]
        \centering
        \vspace*{-1em}
        \includegraphics[height=0.3\textheight]{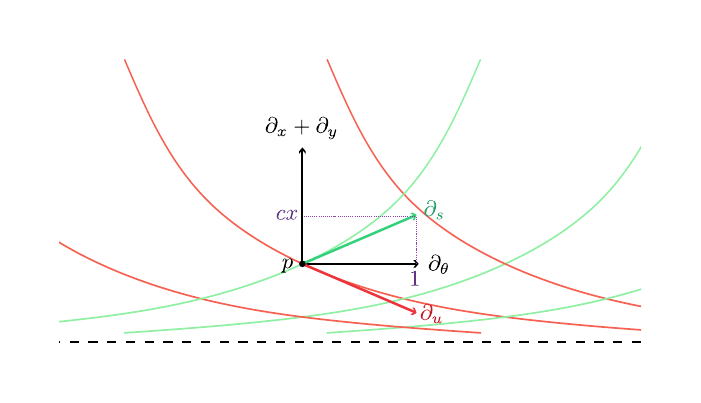}
        \vspace*{-1em}
        \caption{Coordinates of the vector fields $\partial_u$ and $\partial_s$ generating the foliation $\cG^{u, \iin}$ and $\cG^{s, \iin}$ at a point $p =(x, x, \theta) \in \Pin$}
        \label{fig: angle foliation Pin}
    \end{figure}
    
\pagebreak[4]
    
    \item The vector field $X$ on $\Pin{}$ is
    $X(x,x,\theta) =c x \l( \partial_y - \partial_x \r) + \partial_\theta $.
    It is orthogonal to $\partial_x + \partial_y$ at any point of the annulus and the angle between $X$ and $\partial_\theta$ at a point $p =(x,x,\theta) \in \Pin$ is proportional to $x$, thus to the distance from $p$ to $\cO_i$ (Figure~\ref{fig: coordinate vector X on Pin}).
    
    \begin{figure}[htb]
        \centering
        \captionsetup{width=.86\linewidth}
        \vspace*{-2em}
        \includegraphics[height=0.37\textheight]{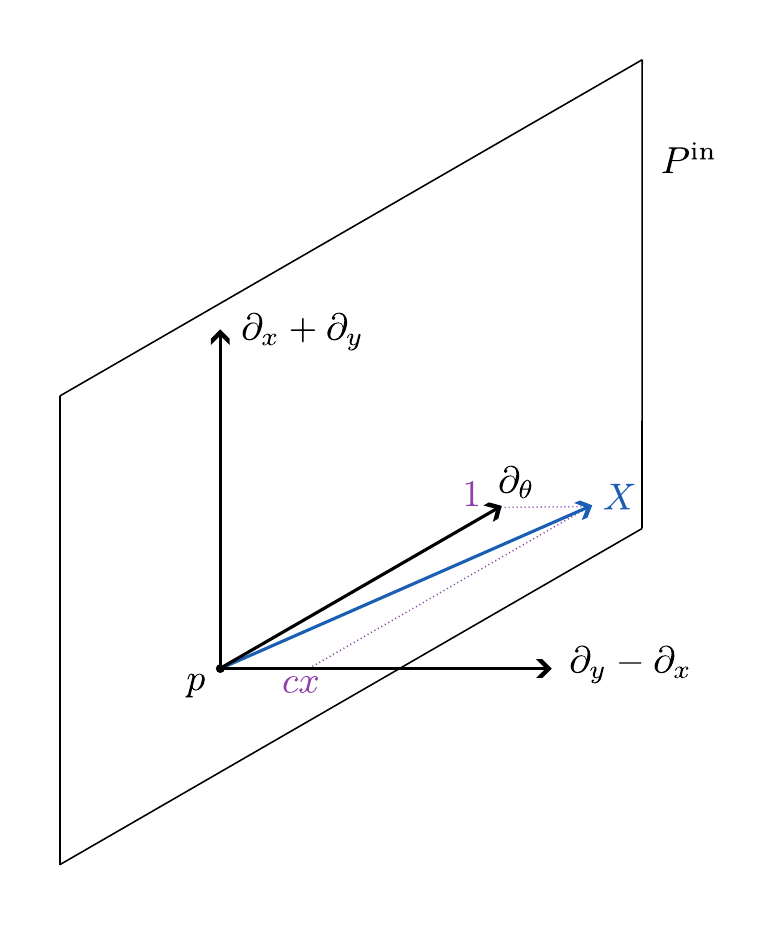}
        \vspace*{-1em}
        \centering \caption{Coordinates of the vector $X(p)$ at point $p=(x,x,\theta) \in \Pin$}
        \label{fig: coordinate vector X on Pin}
    \end{figure}
    \item The direction tangent to $\cG^{u, \iin}$ at the point $p = (x,x,\theta) \in \Pin$ is generated by the vector field 
    $\partial^u = -cx \, (\partial_x + \partial_y) + \partial_\theta = -2 c x \, \partial_y + X(x,x,\theta)$
    Since the vector $X$ is uniformly transverse to $\partial_y$, we deduce that the angle between $X$ and the vector $\partial^u$ at a point $p = (x,x,\theta) \in \Pin$ is proportional to $x$, hence to the distance from $p$ to $\cO_i$.\qedhere
\end{enumerate}
\end{proof}

\subsection{Decomposition of the crossing map}
\label{sec: crossing map; subsec: decomposition map}

In order to explain the crossing of an orbit close to the set $\cO_*$,
we will consider surfaces $\tPin{}$ and $\tPout{}$ which coincide with $\Pin{}$ and $\Pout{}$ except in the \nbh{} of the orbits $\cO_*$ and are uniformly transverse to the \vf{} $X$.
Let $\cV_* = \bigcup_i \cV_i$ be the union of the normalized \nbhs{} of the orbits $\cO_i \in \cO_*$ in a \minc{} $\tilde P$.
The expression of the flow in the coordinate system $(\cV_i, \xi_i = (x,y,\theta) \in \R^2 \times \R/\Z)$ is
$X^t(x,y,\theta) = (2^{-t} x, 2^{t} y, \theta + t)$,
and the boundary of $P$ is identified with $\{ x=y \}$ (Definition~\ref{def: normalized block}).
Let $r_0>0$ so that the point $(r_0, r_0, 0)$ is in the image of $\cV_i$ in coordinates $\xi_i$ for all $i$.

\begin{figure}[htb]
    \centering
    \vspace*{-2em}
    \includegraphics[height=0.3\textheight]{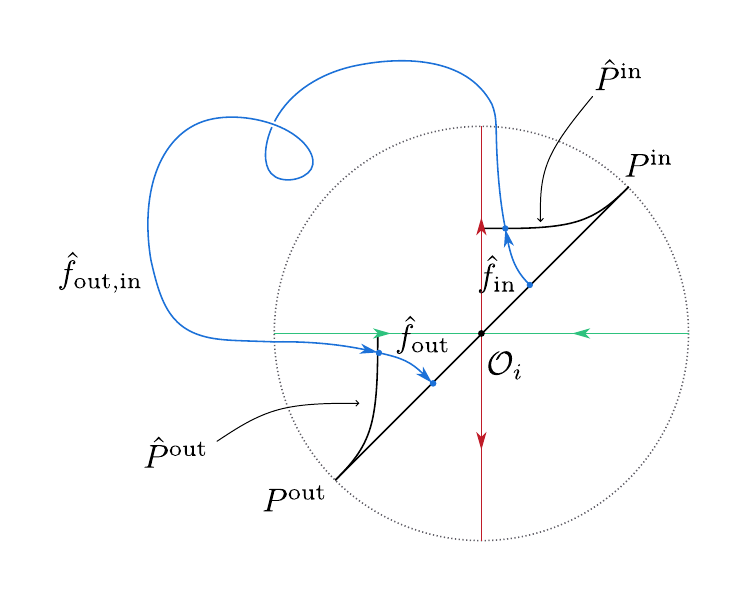}
    \vspace*{-1em}
    \caption{Decomposition of the crossing map $\foutin{}$ in a normalized neighborhood $(\cV_i, \xi_i)$.}
    \label{fig: decomposition crossing map}
\end{figure}

\begin{nota} \label{nota: transverse intermediate annuli}
Let \tPin{} be a smooth surface such that 
$\tPin{}$ coincides with $\Pin{}$ outside $\cV_*$,
$\tPin{} = \{ y=r_0, x >0\}$ close to $x=0$ in $\cV_i$, and $\tPin{}$ is uniformly transverse to the \vf{} $X$.
Similarly, let \tPout{} be a smooth surface such that $\tPout{}$ coincides with $\Pout{}$ outside $\cV_*$, $\tPout{}= \{ x=-r_0, y <0\}$ close to $y=0$ in $\cV_i$,
and $\tPout{}$ is uniformly transverse to the \vf{} $X$.
We refer to Figure~\ref{fig: decomposition crossing map}.
\end{nota}

\begin{rmk}$\,$
\begin{itemize}
    \item Each \cc{} of $\tPin$ (respectively $\tPout$) is either a closed surface or the interior of a surface whose boundary is formed by closed curves in the unstable (respectively stable) manifold of periodic orbits of $\cO_*$.
    \item Each point on \tPin{} has a past orbit which enters via \Pin{}, each point on \tPout{} has a future orbit which exits via \Pout{}.
    The time is not bounded because one can be arbitrarily close to the boundary of \Pin{} and \Pout{}, in which case the orbit spends an arbitrarily long time in the neighborhood of an orbit $\cO_i$ before intersecting \tPin{} or \tPout{}.
    \item By transversality of \tPin{} and \tPout{} with the vector field $X$, the lamination $\cW^s$ induces a lamination of dimension 1 on \tPin{} and the lamination $\cW^u$ induces a lamination of dimension 1 on \tPout{}.
    Denote $\hat \cL^\iin = \cW^s \cap \tPin$ and $\hat \cL^\out = \cW^u \cap \tPout$.
\end{itemize}
\end{rmk}

We write 
$\foutin = \hat f_\out \circ \hat f_{\out, \iin} \circ \hat f_\iin$
where
$\hat f_\iin \colon \Pin{} \to \tPin$, $\hat f_{\out, \iin} \colon \tPin \ssm \hat \cL^\iin \to \tPout \ssm \hat \cL^\out$ and  $\hat f_\out{} \colon \tPout \to \Pout$
are the crossing maps for the flow of $X$ between the corresponding transverse surfaces (see Figure~\ref{fig: decomposition crossing map}).
The composition map is well defined on $\Pin \ssm \cL^\iin$ with image in $\Pout \ssm \cL^\out$.
The local crossings $\hat f_\iin$ and $\hat f_\out$ are trivial outside the normalized neighborhoods $\cV_*$, on which we have explicit expressions.
The crossing map $\hat f_{\out, \iin}$ from \tPin{} to \tPout{} is a crossing map between two sections uniformly transverse to the \vf{}~$X$.

\subsubsection*{Intermediate transit time}
For any $p$ in $P$, recall that $\tau(p)$ is the exit time (if well defined) of the orbit of the point $p$ on $\Pout{}$, in other words such that $X^{\tau(p)}(p) \in \Pout{}$. We also define
\begin{itemize}
    \item $\hat \tau$ the time of intersection of the orbit (positive or negative) of a point with the surface $\tPout$ (when it is well defined), in other words such that $X^{\hat \tau(p)}(p) \in \tPout$.
    This time is unique because an orbit which intersects \tPout{} never passes through \tPout{} again.
    This time is well defined on the complementary of $\cW^s$ in $P$.
    Indeed, any point on $P$ has an orbit which intersects \Pout{}.
    Moreover, any orbit that intersects \Pout{} intersects \tPout{}.
    \item $\tau_\iin$ the time of intersection of the orbit (positive or negative) of a point with the surface $\tPin$ (when it is well defined), in other words such that $X^{\tau_\iin(p)}(p) \in \tPin$.
    This time is unique because an orbit which intersects \tPin{} never passes through \tPin{} again.
    This time is well defined on the complementary of $\cW^u$ in $P$.
    Indeed, any point on $P$ has an orbit which intersects \Pin{}.
    Moreover, any orbit that intersects \Pin{} intersects \tPin{}.
\end{itemize}

\subsection{Proof of Proposition~\ref{prop: crossing map}}
\label{sec: crossing map; subsec: proof}

We are now able to show Proposition~\ref{prop: crossing map}.
We are reduced to the study of the expanding effects on the vectors tangent to $\cG^u$ of the three maps $\hat f_\iin$, $\hat f_{\out, \iin}$, $\hat f_\out$.
First, we show that $\hat f_\iin$ does not contract too much the direction tangent to $\cG^u$ (Lemma~\ref{lem: non-contraction of enter crossing}) on $\Pin$.
In a second step, we show that the map $\hat f_\out \circ \hat f_{\out, \iin}$ arbitrarily strongly expands the direction tangent to $\cG^u$ in a small neighborhood of $\cW^s$ on $\tPin$ (Lemma~\ref{lem: expansion middle crossing} and~\ref{lem: expansion exit crossing}).
In a conclusion, we assemble the previous results to show that the crossing map $\foutin$ expands arbitrarily strongly the direction tangent to $\cG^u$ in a small neighborhood of $\cL^\iin$ on $\Pin$.

\subsubsection*{Step 1: Lower bound for the contraction of $\hat f_\iin$ in the unstable direction}
\begin{lem} \label{lem: non-contraction of enter crossing}
The derivative of $\hat f_\iin$ is uniformly bounded from below in the direction tangent to the foliation $\cG^{u, \iin}$.
In other words, there exists a constant $\cst>0$ such that
$$ \forall p \in \Pin{}, \qq \forall v \in T_p\cG^{u, \iin},\qq \Vert (\hat f_\iin)_* v \Vert \geq \cst \Vert v \Vert.$$
\end{lem}

This lemma tells us that the competition balances between the expansion of an unstable vector $v \in E^\uu \oplus \R.X$ under the effect of the flow, and the contraction under the effect of the projection on the surface $\tPin$.
More precisely, a point on $\Pin{}$ close to $\cO_i$ has an orbit whose local crossing from $\Pin{}$ to $\tPin{}$ is not trivial.
A vector tangent to $\cG^{u, \iin}$ at this point has a small strong unstable component compared to its tangent component at $X$.
The projection parallel to $X$ of such a vector pushed by the flow on the transverse surface $\tPin{}$ has a strong contracting effect on the tangent component to $X$ and a uniformly bounded effect on the strong unstable component.
But the time spent in the \nbh{} of the hyperbolic periodic orbit during the transition from \Pin{} to \tPin{} is long enough for the flow to expand the strong unstable component sufficiently so that the projection on $\tPin{}$ of the vector pushed by the flow is not (too) contracting.

\begin{proof}
We have $(\hat f_\iin)_* = \pi \circ X^{\tau_\iin}_*$,
where $\pi \colon TP \to T \tPin$ is the projection on the tangent space of \tPin{} parallel to $\R. X$.
The map $\hat f_\iin$ is equal to the identity everywhere outside $\cV_* = \bigcup_i \cV_i$, there is nothing to show.
Let us show the property on $\cV$, which admits an $X^t$-invariant hyperbolic splitting.
Let $(\cV_i, \xi_i = (x,y,\theta))$ be a \ncs{} of $\cO_i \in \cO_*$, provided with the hyperbolic splitting $E^\uu \oplus \R.X \oplus E^\ss$ of the tangent space (Lemma~\ref{lem: filtrating neighborhood hyp}).
The flow in these coordinates is
$X^t(x,y,\theta) = (2^{-t} \, x, 2^t \, y, \theta + t)$.
The lemma is based on the following two points.
\begin{enumerate}
    \item \label{it: relation uu X on Pin} Let $v = v^\uu + v^X \in E^\uu \oplus \R. X$ a tangent vector to $\cG^{u, \iin}$ above a point $p \in \Pin{} \cap \cV_i$.
    Then there exists a uniform constant $\cst$ such that
    $$\Vert v^\uu \Vert \geq \cst\, \dist(p, \cO_i) \, \Vert v \Vert.$$
    \item \label{it: crossing time Pin tPin} The transit time $\tau_\iin$ is well defined on $\Pin{}$ and satisfied :
    there exists a uniform constant $\cst{}>0$ such that for any $p \in \Pin{}$,
    $$ \tau_\iin(p) \geq \log_2 (\cst \dist(p, \cO_i)^\inv).$$
    A vector in $\cV_i$ pushed by the (linear) flow during a time $\tau_\iin(p)$ has thus its strong unstable component expanded by a factor bounded below by $\cst \dist(p, \cO_i)^\inv$.
\end{enumerate}

The first point is a direct consequence of Lemma~\ref{lem: angle foliation and vector field on boundary}, Item~\ref{angle gu X on Pin}, which states that the angle $\measuredangle(X, \cG^u_\iin)$ is of order of the distance to the orbit $\cO_i$.
For the second point, the orbit of a point $p=(x,x,\theta) \in \Pin{}$ intersects the section $\tPin{}$ at time $ \tau_\iin(p) = - \log_2 \l( 2x \r) $.
We conclude by noticing that $x$ is equal to the distance $\dist(p, \cO_i)$.

Let $v = v^\uu + v^X \in E^\uu \oplus \R. X$ a tangent vector to $\cG^{u, \iin}$ above a point $p \in \Pin{}$.
Then up to change the constant at each line, we have
{\renewcommand{\arraystretch}{1.5}
$$ \begin{array}{l@{\hspace{2cm}}r}
    \Vert (\hat f_\iin)_* v \Vert  =  \Vert \pi \circ (X^{\tau_\iin(p)})_* v \Vert \geq \cst\, \Vert (X^{\tau_\iin(p)})_* v^\uu \Vert 
    & \text{(Lemma~\ref{lem: projection bounded effect on strong unstable})} \\
    \Vert (\hat f_\iin)_* v \Vert  \geq  \cst\, 2^{\tau_\iin(p)} \Vert v^\uu \Vert \geq \cst\, \dfrac{1}{\dist(p, \cO_i)} \Vert v^\uu \Vert 
    & \text{(Item~\ref{it: crossing time Pin tPin})} \\
    \Vert (\hat f_\iin)_* v \Vert \geq \cst\,  \Vert v \Vert.
    & \text{(Item~\ref{it: relation uu X on Pin})}
\end{array} 
$$}
This completes the proof of Lemma \ref{lem: non-contraction of enter crossing}.
\end{proof}

\subsubsection*{Step 2: Expansion of the composition $\hat f_\out \circ \hat f_{\out, \iin}$ in the unstable direction}

We show that the composition $\hat f_\out \circ \hat f_{\out, \iin}$ expands arbitrarily strongly a vector tangent to $\cG^u \cap \tPin$ in a small neighborhood of $\cW^s \cap \tPin$.
For this we study the expansion property of $\hat f_{\out, \iin}$ on the one hand (Lemma~\ref{lem: expansion middle crossing}) and of $\hat f_\out$ on the other hand (Lemma~\ref{lem: expansion exit crossing}).

\begin{lem} \label{lem: expansion middle crossing}
There exists a constant $\cst>0$ such that for any $p \in \tPin \ssm \cW^s$ and for any vector $v$ tangent to $\cG^u \cap \tPin$ on $p$,
$\Vert (\hat f_{\out, \iin})_* v \Vert \geq \cst \,\lambda_0^{\hat \tau(p)} \Vert v \Vert$.
\end{lem}

\begin{proof}
We have
$(\hat f_{\out, \iin})_* = \hat{\pi} \circ  X^{\hat \tau}_*  $
where $\hat \pi \colon TP \to T \tPout$ is the projection on the tangent space of $\tPout{}$ parallel to $\R.X$.
This lemma is based on the following key points.
Let $\cU$ be the filtrating neighborhood, equipped with a hyperbolic splitting (Lemma~\ref{lem: filtrating neighborhood hyp}).
Let $p \in \tPin{}$, disjoint from $\cW^s$, and $v$ be a vector tangent to $\cG^u \cap \tPin$ at point $p$.
There are two situations.
\begin{enumerate}[leftmargin=*]
    \item The orbit of $p$ never enters $\cU$.
    The orbit of $p$ enters through $\Pin$ and escapes through $\Pout$ at time $\tau(p)$, and the total length, greater than $\tau(p)$, is uniformly bounded by a constant say $T_0$ from Claim~\ref{claim: orbit length bounded outside nbh U}.
    The set $\Pout \ssm \cU$ is contained in $P \ssm \cV_*$, where the surfaces \Pout{} and \tPout{} coincide.
    We deduce that the time $\hat \tau = \tau$ is uniformly bounded on the set of point of $\hat \Pin$ whose orbit never crosses $\cU$.
    Then
\begin{enumerate}
    \item \label{key point norm flow} The action of the differential $X^{\hat \tau}_*$ is uniformly bounded;
    \item \label{key point angle flow} For any vector $v$ tangent to $\tPin$, the angle between the vector $w= X^{\hat \tau}_*(v)$ and the direction tangent to $X$ is uniformly bounded from below;
    \item \label{key point projection} If the angle between a vector $w$ of norm 1 and the direction tangent to $X$ at a point $p\in \tPout$ is uniformly bounded from below, then the norm of the projected vector $\hat \pi (w) \in T \tPout$ is uniformly bounded.
\end{enumerate}
Item~\ref{key point angle flow} is a consequence of the uniform transversality of $\tPin$ to $X$, and Item~\ref{key point projection} is a consequence of the uniform transversality of $\tPout$ to $X$.
We deduce that there exists a constant $\cst'>0$ uniform in $p$ such that
    \begin{align*}
        \| d_p \hat f_{\out, \iin}(v) \| & = \| \hat \pi_p \circ d_p X^{\hat \tau(p)} (v) \| \\
         & \geq \cst' \| d_p X^{\hat \tau(p)} (v) \| &
        \mbox{(Item~\ref{key point angle flow} and~\ref{key point projection})} \\
         & \geq \cst' \| v \|. & \mbox{(Item~\ref{key point norm flow})}
    \end{align*}
    Lemma~\ref{lem: expansion middle crossing} is satisfied by setting $\cst = \cst' (\lambda_0^{T_0})^\inv $.
    
    \item The orbit of $p$ enters $\cU$.
    The orbit enters through $\Pin$ and exit through $\Pout$, hence according to~\ref{rmk: extend hyperbolic splitting}, we can extend the splitting $E^\uu \oplus \R. X \oplus E^\ss$ along the orbit, and the action of the flow is still uniformly hyperbolic because the length of segment of orbit disjoint from $\cU$ are uniformly bounded.  
    The tangent vector $v \in T_p \cap T_p \cG^u$ decomposes into the sum
    $v = v^\uu + v^X \in E^\uu \oplus \R.X$.
    By transversality of \tPin{} with $X$, there exists a uniform constant $\cst >0$ such that $ \Vert v^\uu \Vert \geq \cst \, \Vert v \Vert $.
    By hyperbolicity, the direction $E^\uu$ is preserved and the component $v^\uu$ is expanded by the differential of the flow along the orbit of $p$, i.e., we have $w= X^{\hat \tau}_* (v) = w^\uu + w^X \in E^\uu \oplus \R.X$, and $\| w^\uu \| \geq C_0 \lambda_0^{\hat \tau} \| v^\uu \| \geq \cst \lambda_0^{\hat \tau} \| v \|$.
    Moreover, we know from Lemma~\ref{lem: projection bounded effect on strong unstable} that the norm of the projection $\hat \pi(w)$ is essentially bounded from below by the strong unstable component $w^\uu$, hence
        \begin{align*}
        \| d_p \hat f_{\out, \iin}(v) \| & = \|  \hat \pi_p \circ d_p X^{\hat \tau(p)} (v) \| =
        \|  \hat \pi_p (w^\uu + w^X) \|  \\
          &\geq \cst \| w^\uu\| \geq \cst \lambda_0^{\hat \tau}  \| v \| .\qedhere
    \end{align*}
\end{enumerate}
\end{proof}

Finally, we show that the exit crossing map $\hat f_\out{}$ exponentially expands the unstable vectors as a function of the transit time from \tPout{} to \Pout{}.

\begin{lem} \label{lem: expansion exit crossing}
There exists a constant $\cst>0$ such that for any $p \in \tPout{}$, for any vector $v$ tangent to the foliation $\cG^u \cap \tPout$ on $p$, we have $\Vert (\hat f_\out{})_*v \Vert \geq \cst\, \lambda_0^{\tau(p)} \Vert v \Vert $.
\end{lem}

In this sense, there is no competition between expansion and contraction.
Indeed, $\hat f_\out$ is a crossing map of the flow starting from a surface \tPout{} uniformly transverse to $X$.
The projection on the destination $\Pout{}$ will not be too contracting.
Moreover, the transport by the flow in the neighborhood of $\cO_*$ contributes to exponentially expand an unstable vector.

\begin{proof}
We have $(\hat f_\out{})_* = \pi \circ X^{\tau}_* $
where $\pi$ is the projection on $\Pout{}$ parallel to $\R. X$, and the transit time $\tau$ is well defined everywhere on \Pout{}.
The argument is the same as the proof of the previous lemma.
We stand in a domain where the map $\hat f_\out$ is not trivial, in other words a neighborhood $\cV_i \subset \cV_*$ of a periodic orbit $\cO_i$.
There exists an invariant hyperbolic splitting $E^\uu \oplus \R. X \oplus E^\ss$ on $\cV_i$.
For $p \in \tPout \cap \cV_i$, a vector tangent to $v \in T_p \tPout \cap T_p \cG^u$ decomposes into the sum
    $v = v^\uu + v^X \in E^\uu \oplus \R.X$.
    By transversality of \tPout{} with $X$, there exists a uniform constant $\cst >0$ such that 
    $\Vert v^\uu \Vert \geq \cst \, \Vert v \Vert $.
    The orbit of $p$ is entirely contained in $\cV_i$, and the intersection time $\tau$ with $\Pout$ is strictly positive.
    By hyperbolicity, the direction $E^\uu$ is preserved and the component $v^\uu$ is expanded by the differential $X^\tau_*$ of the flow along the orbit of $p$, i.e., we have $w= X^{\tau}_* (v) = w^\uu + w^X \in E^\uu \oplus \R.X$, and $$\| w^\uu \| \geq C_0 \lambda_0^{\tau} \| v^\uu \| \geq \cst \lambda_0^{\tau} \| v \|.$$
    Moreover, we know from Lemma~\ref{lem: projection bounded effect on strong unstable} that the norm of the projection $\hat \pi(w)$ is essentially bounded below by the strong unstable component $w^\uu$, hence
        \begin{align*}
        \| d_p \hat f_\out(v) \| & = \|  \hat \pi_p \circ d_p X^{\tau(p)} (v) \| =
        \|  \hat \pi_p (w^\uu + w^X) \|  \\
         & \geq \cst \| w^\uu\|  \geq \cst \lambda_0^{\tau} \| v \|. \qedhere
    \end{align*} 
\end{proof}

\subsubsection*{Conclusion: Expansion of the crossing map $\foutin$ in the unstable direction}

\begin{proof}[Proof of Proposition~\ref{prop: crossing map}]
We show the first item.
Let $\lambda>1$, and $\cst>0$ be constants satisfying Lemmas~\ref{lem: non-contraction of enter crossing},~\ref{lem: expansion middle crossing} and~\ref{lem: expansion exit crossing}.

\begin{claim} \label{claim: T_lambda}
There exists $\hat \cW^\iin_\lambda$ a neighborhood of $\hat \cL^\iin = \cW^s \cap \tPin$ on $\tPin$ such that the orbit of any point $p \in \hat \cW_\lambda$ is defined for a time at least equal to
\begin{equation}\label{eq: T_lambda}
    T_\lambda = \log_{\lambda_0}(\lambda) - 3\log_{\lambda_0}(\cst).
\end{equation}
\end{claim}

\begin{proof}
This neighborhood exists because a point arbitrarily close to $\cW^s$ has an orbit which spends an arbitrarily long time in the neighborhood of $\Lambda$, thus in the \minc{} $\tilde P$, whose boundary is disjoint from $\Lambda$.
Moreover, an orbit of a point in \tPin{} defined at time $T$ in $\tilde P$ is defined at time $T' \geq \frac{T}{2}$ in $P$.
Indeed, the complementary $\tilde P \ssm P$ is contained in the neighborhood $\cV_*$ by definition of a \minc{} (Definition~\ref{def: minc}), and it is clear that an orbit that enters a neighborhood $\cV_i$ spends a time greater than half of $T$ in $P$, which concludes.
\end{proof}

Let $\cW^\iin_\lambda = (\hat f_\iin)^\inv(\hat{\cW}^\iin_\lambda)$.
This is a neighborhood of $\cL^\iin$ in \Pin{} (Figure~\ref{fig: expanding neighborhood}).

\begin{figure}[htb]
    \centering
    \vspace*{-1em}
      \includegraphics[height=0.34\textheight]{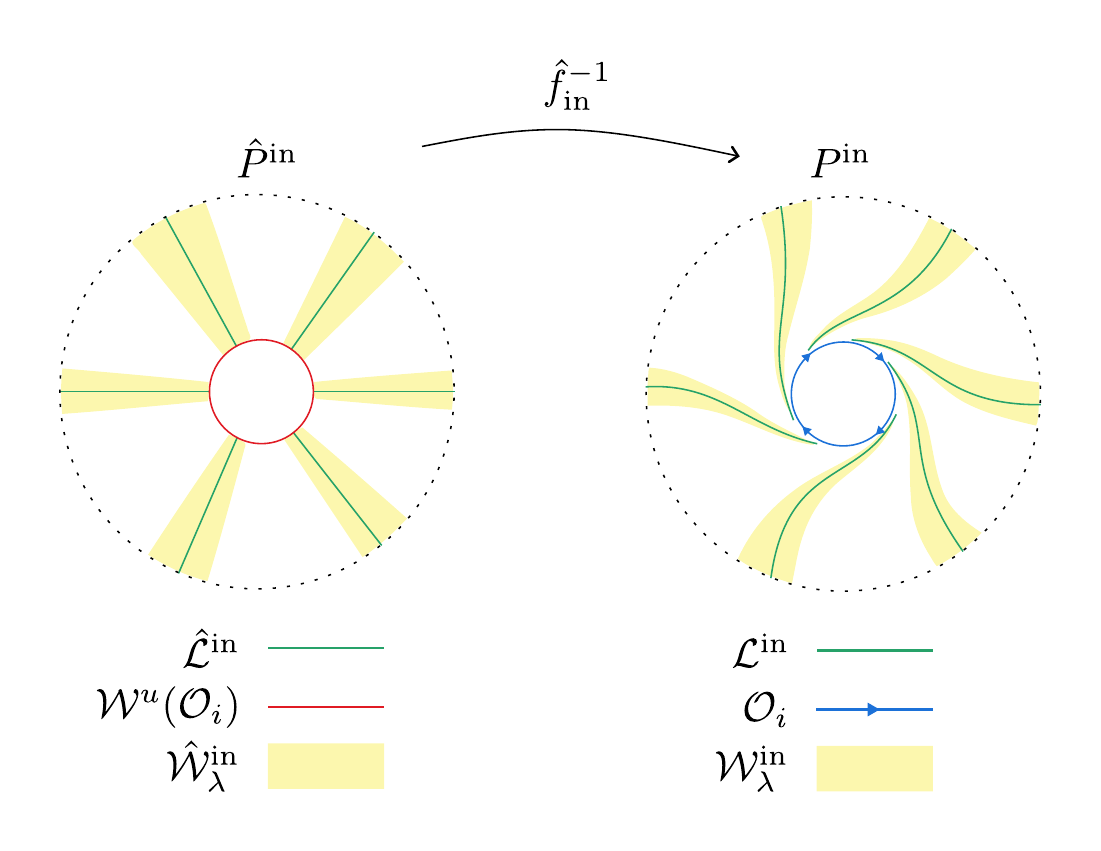}
       \vspace*{-1em}
    \caption{Construction of the neighborhood $\mathcal{W}^{\mathrm{in}}_\lambda$ of $\mathcal{L}^{\mathrm{in}}$ in $P^{\mathrm{in}}$}
    \label{fig: expanding neighborhood}
\end{figure}

Let
$p \in  \cW^\iin_\lambda \ssm \cL^\iin$,  $q = \hat f_\iin(p) \in \tPin{}$ and $ s = \hat f_{\out, \iin}(q) \in \tPout{}$.
Let $v$ be a vector tangent to $\cG^{u, \iin}$ at $p$.
Let us note that $\tau(s) = \tau(q) - \hat \tau(q)$ and $ \tau(q) \geq T_\lambda$ by Claim~\ref{claim: T_lambda}.
Then
\begin{align*}
\Vert (\foutin)_* v \Vert & =  \Vert (\hat f_\out{} \circ \hat f_{\out, \iin} \circ \hat f_\iin)_* v \Vert \\
& \geq \cst\, \lambda_0^{\tau(s)} \, \Vert (\hat f_{\out, \iin} \circ \hat f_\iin)_* v \Vert  
& \text{(Lemma~\ref{lem: expansion exit crossing})}\\
& \geq   \cst^2 \, \lambda_0^{\tau(q) -\hat \tau(q)} \, \lambda_0^{\hat \tau(q)} \, \Vert (\hat f_\iin)_* v \Vert  
& \text{(Lemma~\ref{lem: expansion middle crossing})}\\
& \geq  \cst^3 \, \lambda_0^{T_\lambda} \Vert v \Vert  
& \text{(Lemma~\ref{lem: non-contraction of enter crossing})} \\
& \geq  \lambda \, \Vert v \Vert.
& \text{(from Equation~\eqref{eq: T_lambda}}
\end{align*}

This concludes the proof of the first item of Proposition~\ref{prop: crossing map}. The second item is completely symmetrical.
\end{proof}

\section{Spread the hyperbolicity}
\label{sec: spreading}

In the previous section we have shown that the expansion of the crossing map $\foutin \colon \Pin{} \to \Pout{}$ in the direction tangent to $\cG^{u, \iin}$ is arbitrarily strong close to the lamination $\cL^\iin$ (Proposition~\ref{prop: crossing map}).
We refer to the beginning of Section~\ref{sec: crossing map} for a reminder of the notations.
We will now \say{spread} this expansion. 
The goal is that the map $\foutin$ expands the direction tangent to $\cG^{u, \iin}$, not only in the neighborhood of $\cL^\iin$, but almost all over \Pin{}.
This spreading is done with a \diff{} $\psi^\iin \colon \pP \to \pP$, supported on $\Pin$, which we can think of as a change of coordinates on \Pin{}.
The \diff{} $\psi^\iin$ will depend on three parameters \led{}.
Let us present them informally:
\begin{itemize}
    \item $\lambda > 1$ is a lower bound of the expansion factor of the composition \linebreak[4]\mb{$\foutin \circ \psi^\iin$} in the direction $\cG^{u, \iin}$ outside a neighborhood of the orbits~$\cO_*$;
    \item $\epsilon>0$ measures the perturbation of the foliation $\cG^{s, \iin}$ induced by the action of $\psi^\iin$;
    \item $\delta>0$ is the size of the neighborhood of the orbits $\cO_*$ on which we no longer control the expansion of the composition $\foutin \circ \psi^\iin$.
\end{itemize}

Similarly, the expansion of the inverse crossing map $\mb{$\foutin^\inv \colon \Pout{} \to \Pin{}$}$ in the direction $\cG^{s, \out}$ is arbitrarily strong close to the lamination $\cL^\out$ (Proposition~\ref{prop: crossing map}).
The symmetrical operation consists in \say{spreading} the expansion of $\foutin^\inv$ in the direction $\cG^{s, \out}$ far from the lamination $\cL^\out$, by means of a \diff{} $\psi^\out: \pP \to \pP$ supported in \Pout{}, which one can think of as a change of coordinates on \Pout{}.
These two symmetrical steps must be done in such a way that one does not destroy the other, in other words so that the composition $\psi^\out \circ \foutin \circ \psi^\iin$ continues to expand the direction $\cG^{u,\iin}$ and its inverse continues to expand the direction $\cG^{s,\out}$.
This will be possible thanks to the control of the perturbation of the foliation by the \diffs{} $\psi^\iin$ and $\psi^\out$.

Let us translate this result into the framework of the gluing procedure of the block $(P,X)$.
If $\varphi$ is a normalized gluing map of the block $(P,X)$, then:
\begin{itemize}[leftmargin=*]
    \item Replacing $\foutin$ by $\psi^\out \circ \foutin \circ \psi^\iin$ amounts to replacing the gluing map $\varphi$ by $\psi^\iin \circ \varphi \circ \psi^\out$.
    Indeed the first return map of the flow of $X_\varphi$ in the glued manifold $P/\varphi$ on the projection of the surface $\partial P$ is the composition of the crossing map $\foutin \colon \Pin \to \Pout$ and the gluing map $\varphi$.
    If $\psi = \psi^\iin \circ \varphi \circ \psi^\out$ is another gluing map, then the first-return map of $X_\psi$ on $\partial P$ is (up to conjugate) the map $\varphi \circ (\psi^\out \circ \foutin \circ \psi^\iin)$.
    \item Since $\psi^\iin$ and $\psi^\out$ have their support contained in the complementary of a neighborhood of the orbits $\cO_*$, the composition $\psi^\iin \circ \varphi \circ \psi^\out$ is a dynamical gluing map, because it satisfies Item~\ref{def: normalized gluing; it: trivial} of Definition~\ref{def: normalized gluing map}.
    Moreover, the properties of $\psi^\iin$ and $\psi^\out$ imply that the composition $\psi^\iin \circ \varphi \circ \psi^\out$ is isotopic to $\varphi$ among the \sqt{} gluing maps, which implies that the triples are strongly isotopic (Definition~\ref{def: strongly isotopic triples}).
\end{itemize}

\subsection{Statement of the main proposition}

\label{sec: spread; subsec: statement}

\begin{figure}[htb]
    \centering
    \vspace*{-1em}
    \includegraphics[height=0.28\textheight]{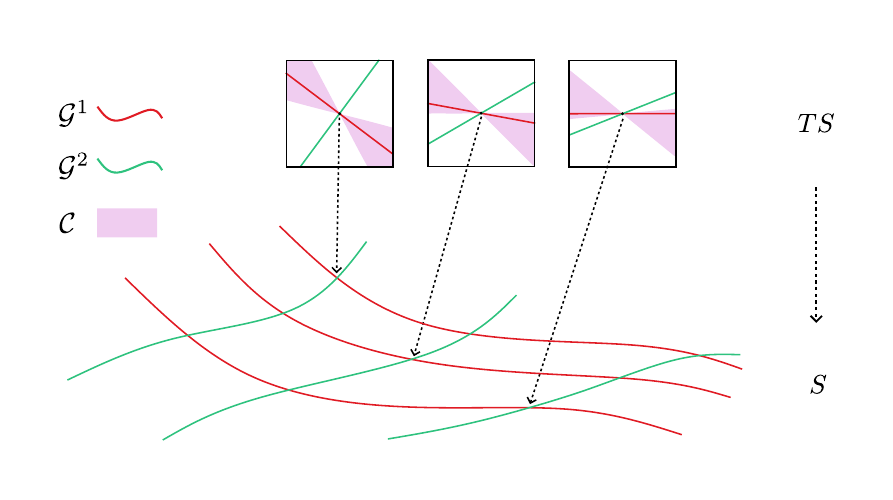}
    \vspace*{-1em}
    \caption{A $(K, \cG^1/\, \cG^2)$-cone field}
    \label{fig: cone field slope foliation}
\end{figure}

We will use cone fields on the surface $\partial P$.
We will talk about the slope of a cone field on a surface relatively to a pair of transverse foliations on the surface (Figure~\ref{fig: cone field slope foliation}).

\begin{defi}[Cones] \label{def: cone slope}
Let $S$ be a surface equipped with a Riemannian norm $\Vert \cdot \Vert$, and $\cG^1$ and $\cG^2$ be two transverse foliation on $S$, and $K \geq 0$.
\begin{itemize}
    \item We say that $C \in T_p S$ is a (closed) \emph{$(K,\cG^1/\,\cG^2)$-cone} if 
    $$C = \{ v_1 + v_2 \in T_p\cG^1 \oplus T_p\cG^2 \mid \Vert v_2 \Vert \leq K \Vert v_1 \Vert \},\footnote{The order is important: a $(K,\cG^1/\,\cG^2)$-cone contains the direction tangent to $\cG^1$ and does not contain the tangent direction to $\cG^2$.}$$
    where $T_p\cG^1 \oplus T_p\cG^2$ is the direct sum of the transverse directions tangent respectively to the foliation $\cG^1$ and $\cG^2$ on $p$.
    We say that $K$ is \emph{the slope of~$C$}.
    
    \item The \emph{interior of a cone} is the set
    $$\intr C := \{ v_1 + v_2 \in T_p\cG^1 \oplus T_p\cG^2 \mid \Vert v_2 \Vert < K \Vert v_1 \Vert \} \cup \{O\}.$$    
\end{itemize}
Let $C=\{C(p), p\in T_p S\}_{p\in S}$ be a cone field. 
\begin{itemize}
    \item We say that $C$ is a $(K(p),\cG^1/\,\cG^2)$-cone field if each $C(p)$ is a $(K(p),\cG^1/\,\cG^2)$-cone on $T_p S$.
    \item If $K(p)=K$ is constant, we will say that $C$ is a $(K,\cG^1/\,\cG^2)$-cone field.
    \item The $(1,\cG^1/\,\cG^2)$-cone field is called the \emph{bisector} cone field of $(\cG^1/\,\cG^2)$.
\end{itemize}
\end{defi}

The main result of the chapter is the following.

\begin{prop}[Spreading of the expansion] \label{prop: spread expansion}
For any $\lambda>1$, $\epsilon>0$ and $\delta>0$, there exists a diffeomorphism \mb{$\psi^\iin = \psi^\iin_\led \colon \pP \to \pP$} such that:
\begin{enumerate}
    \item \label{prop: spread; it: support} $\psi^\iin$ is the identity outside a compact contained in a finite number of \ccs{} of $\Pin \ssm \cL^\iin$.
    
    \item \label{prop: spread; it: foliation preserved} $\psi^\iin$ preserves the foliation $\cG^{u, \iin}$ leaf-to-leaf.
    
    \item \label{prop: spread; it: foliation perturbed} $\psi^\iin$ slightly perturbs the foliation $\cG^{s, \iin}$: the direction tangent to \linebreak[4]$(\psi^\iin)_{*}^\inv \cG^{s, \iin}$ lies inside a $(\epsilon,\cG^{s, \iin}/\cG^{u, \iin})$-cone field.

    \item \label{prop: spread; it: derivative bounded} The derivative of $\psi^\iin$ in the direction $\cG^{u, \iin}$ is uniformly bounded by a constant depending only on $\lambda$: for all $\lambda >1$, there exists \mb{$\cst = \cst(\lambda)>0$}, such that for all $\ed{} >0$ we have
    \[\forall p \in \Pin{}, \qq \forall v \in T_{p} \cG^{u, \iin}, \qq \cst^\inv \, < \frac{\Vert \psi^\iin_* v \Vert }{\Vert v \Vert} < \cst.\]
    \item \label{prop: spread; it: composition expansion} If $\cV_\delta \subset \cV_*$ denotes the $\delta$-neighborhood of the periodic orbits of $\cO_*$, then the composition \mb{$\foutin \circ \psi^\iin$} expands by a factor $\lambda$ the direction tangent to $\cG^{u, \iin}$ on $\Pin \ssm \cL^\iin$ in the complementary of $\cV_\delta$:
    $$\forall p \in (\Pin \ssm \cL^\iin) \cV_\delta, \qq \forall v \in T_p\cG^{u, \iin}, \qq \Vert (\foutin \circ \psi^\iin)_* v \Vert \geq \Vert \lambda \Vert.$$
    \item \label{prop: spread; it: independant delta}
    $\psi^\iin$ does not depend on $\delta$ on the complementary of $\cV_*$:
    for any $\lambda>1$, $\epsilon>0$, the family of \diffs{} $\big\{ \psi^\iin_{\lambda, \epsilon, \delta}\big\}_{\delta >0}$ coincide on $\pP \ssm \cV_*$.
    \end{enumerate}
\end{prop}

\begin{figure}[h]
    \centering
    \vspace*{-2.0em}
    \includegraphics[width=1.0\textwidth]{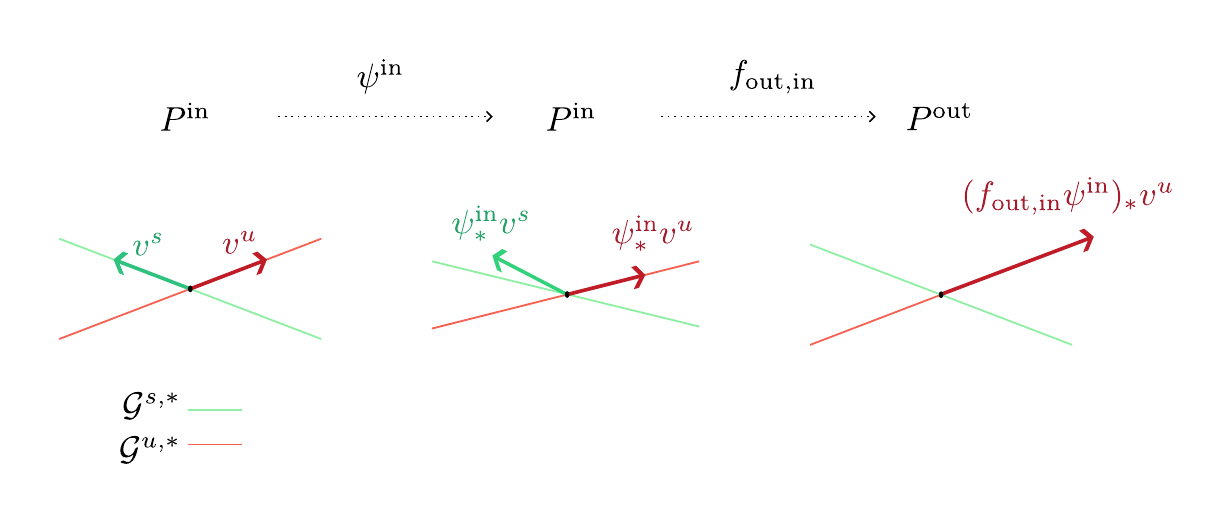}
    \vspace*{-2.0em}
    \caption{Action of the composition $\foutin \circ \psi^\iin$ on $\Pin \ssm \cV_\delta$}
    \label{fig: action psi in}
\end{figure}

We refer to Figure~\ref{fig: action psi in} to illustrate the properties of $\psi^\iin$.

We state the following symmetric result which is proved in an analogous way.

\begin{prop}[Spreading of the expansion] \label{prop sym: spread expansion}
For any $\lambda>1$, $\epsilon>0$ and $\delta>0$, there exists a diffeomorphism \mb{$\psi^\out = \psi^\out_\led \colon \pP \to \pP$} such that
\begin{enumerate}
    \item \label{prop sym: spread; it: support} $\psi^\out$ is the identity outside a compact contained in a finite number of \ccs{} of $\Pout \ssm \cL^\out$.
    
    \item \label{prop sym: spread; it: foliation preserved} $\psi^\out$ preserves the foliation $\cG^{s, \out}$ leaf-to-leaf.
    
    \item \label{prop sym: spread; it: foliation perturbed} $\psi^\out$ slightly perturbs the foliation $\cG^{u, \out}$: the direction tangent to \linebreak[4]$\psi^\out_* \cG^{u, \out}$ lies inside a $(\epsilon,\cG^{u, \out}/, \cG^{s, \out})$-cone field.

    \item \label{prop sym: spread; it: derivative bounded} The derivative of $\psi^\out$ in the direction $\cG^{u, \out}$ is uniformly bounded by a constant depending only on $\lambda$: for all $\lambda >1$, there exists \mb{$\cst = \cst(\lambda)>0$}, such that for all $\ed{} >0$ we have
    $$\forall p \in \Pout{}, \qq \forall v \in T_{p} \cG^{u, \out}, \qq \cst^\inv \, < \frac{\Vert \psi^\out_* v \Vert }{\Vert v \Vert} < \cst.$$
    \item \label{prop sym: spread; it: composition expansion} If $\cV_\delta \subset \cV_*$ denotes the $\delta$-neighborhood of the periodic orbits of $\cO_*$, then the composition \mb{$(\psi^\out \circ \foutin)^\inv$} expands by a factor $\lambda$ the direction tangent to $\cG^{s, \out}$ on $\Pout \ssm \cL^\out$ in the complementary of $\cV_\delta$:
    $$\forall p \in (\Pout{} \ssm \cL^\out) \ssm \cV_\delta, \qq \forall v \in T_{p} \cG^{s, \out}, \qq \Vert (\psi^\out \circ \foutin)^\inv_* v \Vert \geq \lambda \Vert v \Vert.$$
    
    \item \label{prop sym: spread; it: independant delta}
    $\psi^\out$ does not depend on $\delta$ on the complementary of $\cV_*$:
    for any $\lambda>1$, $\epsilon>0$, the family of \diffs{} $\{ \psi^\out_{\lambda, \epsilon, \delta} \}_{\delta >0}$ coincide on $\pP \ssm \cV_*$.
    \end{enumerate}
\end{prop}

We follow closely the proof of \cite[Proposition~6.2]{beguinBuildingAnosovFlows2017} to spread the expansion on the interior of the strips.
The difference is that the strips in our blocks have ends which can accumulate on periodic orbits.
The spreading in the \nbh{} of a strip end which accumulates on a periodic orbit will have to be controlled, which justifies the new parameter $\delta$.

\subsubsection*{Construction of the \diff{} on a strip}
In the previous section, we proved (Proposition~\ref{prop: crossing map}) that for any $\lambda>1$, there exists a neighborhood $\cW^\iin_\lambda$ of $\cL^\iin$ on $\Pin$ such that the crossing map $\foutin \colon \Pin{} \ssm \cL^\iin \to \Pout{} \ssm \cL^\out$ expands the direction $\cG^{u, \iin}$ by a factor $\lambda$ on $\cW^\iin_\lambda$.
Therefore the identity diffeomorphism $\psi^\iin = \Id$ satisfies all items of Proposition~\ref{prop: spread expansion} on the neighborhood $\cW^\iin_\lambda$, for the parameter $\lambda$ and for any parameter $\epsilon>0$ and $\delta>0$.
To prove Proposition~\ref{prop: spread expansion}, it is therefore sufficient to work on the complementary of $\cW^\iin_\lambda$ in $\Pin$.
The complementary $\Pin \ssm \cW^\iin_\lambda$ is contained in a finite number of \ccs{} of $\Pin \ssm \cL^\iin$ which we denote $B_1, \dots, B_n$.
Recall that the boundary lamination $\cL$ of $(P,X)$ is filling, so the \ccs{} of $\Pin \ssm \cL^\iin$ are strips bounded by two noncompact leaves of $\cL^\iin$ which are asymptotic to each other at both ends.
The proof of Proposition~\ref{prop: spread expansion} consists in building a diffeomorphism \mb{$ \psi_i: B_i \to B_i$} on each strip $B_i$ which satisfies the items of Proposition~\ref{prop: spread expansion} on this strip, and equal to the identity outside a compact of the strip.
We can extend $\psi_i$ by the identity on $\pP$ and their product $\psi^\iin = \prod_i \psi_i$ defines a diffeomorphism which satisfies the proposition.
In summary, it suffices to prove the following 

\begin{lem}\label{lem: spread expansion strip}
Let $B \subset \Pin \ssm \cO_*$ be a strip. 
For all $\lambda>1$, $\epsilon>0$ and $\delta>0$,
there exists a \diff{} $\psi = \psi_\led{} \colon B \to B$, such that:

\begin{enumerate}
    \item $\psi$ is the identity outside a compact subset of $B$.
    \label{lem: spread strip; it: support}
    
    \item $\psi$ preserves $\cG^{u, \iin}$ leaf-to-leaf.
    \label{lem: spread strip; it: foliation preserved}
    
    \item The direction tangent to $\psi^\inv_* (\cG^{s, \iin})$ lies inside a $(\epsilon,\cG^{s, \iin}/, \cG^{u, \iin})$-cone field.
    \label{lem: spread strip; it: foliation perturbed}
    
    \item The derivative of $\psi$ in the direction tangent to $\cG^{u, \iin}$ is uniformly bounded by a constant depending only on $\lambda$: for all $\lambda >1$, there exists $\cst = \cst(\lambda)>0$, such that for all $\ed{} >0$ we have
    \[ \forall p \in \Pin{}, \qq \forall v \in T_{p} \cG^{u, \iin}, \qq \cst^\inv \, < \frac{\Vert \psi_* v \Vert }{\Vert v \Vert} < \cst.\]
    \label{lem: spread strip; it: derivative bounded}\vspace*{-1em}
    \item If $\cV_\delta \subset \cV_*$ denotes the $\delta$-neighborhood of $\cO_*$ orbits, then the composition \mb{$\foutin \circ \psi$} expands by a factor $\lambda$ the direction tangent to $\cG^{u, \iin}$ on $B \ssm \cV_\delta$:
    \[\forall p \in B \ssm \cV_\delta, \qq \forall v \in T_p\cG^{u, \iin}, \qq \Vert (\foutin \circ \psi)_* v \Vert \geq \Vert \lambda \Vert.\]
    \label{lem: spread strip; it: composition expansion}\vspace*{-1em}
    \item 
    $\psi$ does not depend on $\delta$ on the complementary of $\cV_*$ :
    for any $\lambda>1$, $\epsilon>0$, the family of \diffs{} $\{ \psi_{\lambda, \epsilon, \delta} \}_{\delta >0}$ coincide on $B \ssm \cV_*$.
    \label{lem: spread strip; it: independant delta}
\end{enumerate}
\end{lem}

Let us draw attention to one point.
The support of $ \psi \colon B \to B$ is a compact contained in $B$, but the complementary $B \ssm \cW^\iin_\lambda$ is not necessarily a compact contained in $B$.
A strip $B$ accumulates at both ends on a compact leaf of the boundary lamination $\cL$.
If one of these compact leaves is a periodic orbit $\cO_j$, then the neighborhood $\cW^\iin_\lambda$ does not overlap the corresponding end of the strip $B$,
because there is no $\delta$-neighborhood of the $\cO_*$ periodic orbits contained in $\cW^\iin_\lambda$.
Therefore the expansion property of the composition $\foutin \circ \psi^\iin$ along $\cG^{u, \iin}$ is no longer true on a $\delta$-neighborhood of $\cO_*$.

Let us show that Proposition~\ref{prop: spread expansion} follows from Lemma~\ref{lem: spread expansion strip}.

\begin{proof}[Proof of Proposition~\ref{prop: spread expansion}]
Let $\lambda>1$, $\epsilon>0$, $\delta>0$.
We define the diffeomorphism $\psi^\iin = \psi^\iin_\led \colon \Pin{} \to \Pin{}$ as follows.
Let $\cW^\iin_\lambda$ be the neighborhood given by Proposition~\ref{prop: crossing map}.
The complementary of $\cW^\iin_\lambda$ is contained in a finite number of \ccs{} $B_1, \dots, B_n$ of $\Pin \ssm \cL^\iin$ (the integer $n$ depends on $\lambda$).
For $i= 1, \dots, n$, let be the \diff{} $\psi_i = (\psi_i)_\led \colon B_i \to B_i$ given by Lemma~\ref{lem: spread expansion strip}.
Item~\ref{prop: spread; it: support} of Lemma~\ref{lem: spread expansion strip} allows us to extend each $\psi_i$ by the identity on $\pP$, and we define
$$    \psi_\led^\iin = \underset{i=1}{\overset{m}{\prod}} { \psi_i } \colon \pP \to \pP.$$
The support of $\psi^\iin_\led$ is contained in a finite number of strips $B_i$ of the entrance boundary, and is compact in each of these strips, which gives Item~\ref{prop: spread; it: support}.
Items~\ref{prop: spread; it: foliation preserved}, ~\ref{prop: spread; it: foliation perturbed} and~\ref{prop: spread; it: independant delta} are direct consequences of Items~\ref{lem: spread strip; it: foliation preserved}, ~\ref{lem: spread strip; it: foliation perturbed} and~\ref{lem: spread strip; it: independant delta} of Lemma~\ref{lem: spread expansion strip}.
Item~\ref{prop: spread; it: derivative bounded} is true because the $\psi_i$ \diffs{} are in finite number.

Let us show Item~\ref{prop: spread; it: composition expansion}.
On a strip $B_i$, we have $\psi^\iin = \psi_i$. We use Item~\ref{lem: spread strip; it: composition expansion} of Lemma~\ref{lem: spread expansion strip}.
On the complementary of the strips ${}^c{(B_1 \cup \dots \cup B_m)} $ and by Proposition~\ref{prop: crossing map}, $\foutin \circ \psi^\iin = \foutin$ expands by a factor $\lambda$ the direction $\cG^{u, \iin}$.
\end{proof}

\begin{figure}[htb]
    \centering
    \vspace*{-1em}
    \includegraphics[height=0.28\textheight]{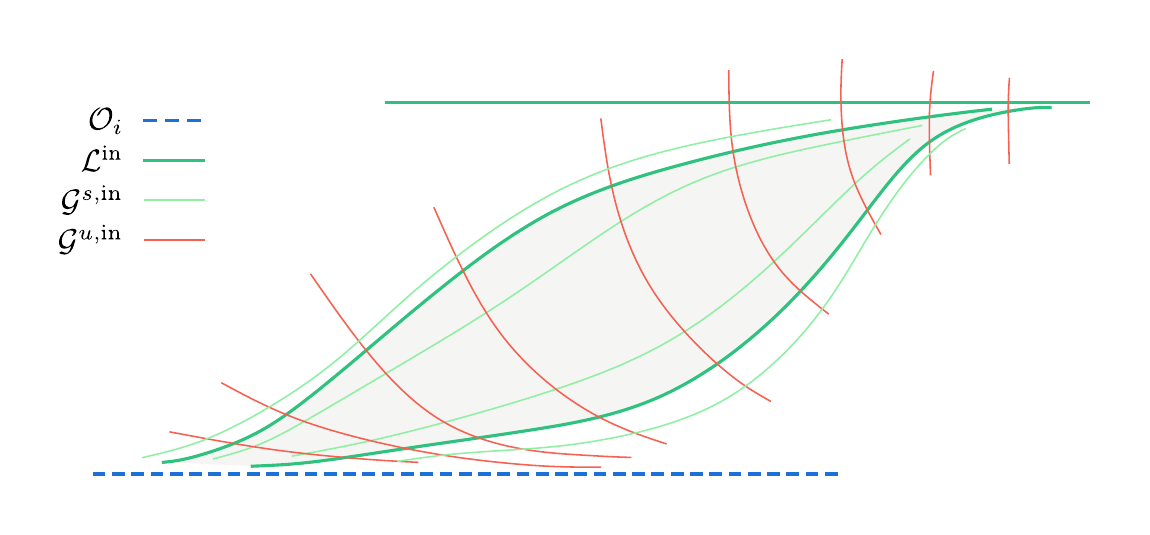}
    \vspace*{-1em}
    \caption{A strip which accumulates on one end on a periodic orbit $\cO_i \in \cO_*$ and on the other end on a compact leaf of $\cL^\iin$}
    \label{fig: entrance strip with foliations}
\end{figure}

The rest of the chapter is devoted to the proof of Lemma~\ref{lem: spread expansion strip}.
Let $B$ be a \cc{} of $\Pin \ssm \cL^\iin$.
Recall that $B$ is a strip bordered by two non-compact leaves $l_1$ and $l_2$ of the lamination $\cL^\iin$, which are asymptotic to each other at both ends (Figure~\ref{fig: entrance strip with foliations}).
Each of the two ends of the strip $B$ accumulates on a compact leaf $\cQ_i$ of the boundary lamination $\cL$ of $(P,X)$ (possibly the same in the non orientable case).
There are two possible cases for the compact leaf $\cQ_i$: either it is a compact leaf of the lamination $\cL^\iin$, or it is a periodic orbit of $(P,X)$.

\begin{rmk} \label{rmk: figure choice strip}
To fix the ideas, we make the choice to represent in our figures a strip $B$ with one end accumulating on a periodic orbit and the other on a compact leaf of $\cL^\iin$.
This strip is contained in a \cc{} $A$ of the complementary $\Pin \ssm \Gamma$, where $\Gamma$ is the union of compact leaves of $\cL$, this component $A$ is an annulus, and the foliation $\cG^{s, \iin}$ on $A$ is a suspension foliation.
Of course, one should keep in mind that the following situations are also possible:
\begin{itemize}
    \item $A$ is an annulus and the boundary compact leaves are two leaves of $\cL^\iin$, or two periodic orbits;
    \item $A$ is a Moebius strip bordered by a unique compact leaf of $\cL^\iin$, or a periodic orbit;
    \item the foliation $\cG^{s, \iin}$ forms a Reeb component on $A$.
\end{itemize}
\end{rmk}

\subsection{Action of the holonomy}
Before constructing the \diff{} $\psi \colon B \to B$, we will need some general results on the holonomy maps of foliations on \Pin{} and \Pout{}.

\subsubsection*{Bound for the distortion of the holonomy of the stable foliation in a strip on the entrance boundary}
Let $\sigma$ and $\sigma'$ be leaves of the foliation $\cG^{u, \iin} \cap B$.
These are open arcs contained in a leaf of $\cG^{u, \iin}$, and which cross the strip $B$ from (accessible) boundary to (accessible) boundary, i.e., the two ends of $\sigma$ and $\sigma'$ lie on the two leaves $l_1$ and $l_2$ of $\cL^\iin$ which form the accessible boundary of $B$.
Denote the holonomy map of the foliation $\cG^{s, \iin}$ between the leaves $\sigma$ and $\sigma'$ by $H^{s, \iin}_{\sigma', \sigma} \colon \sigma \to \sigma'$.
These maps are well defined because each leaf of $\cG^{s, \iin}$ in $B$ intersects one and only one time each leaf $\sigma$ of $\cG^{u, \iin} \cap B$.
Since the leaves are $\cC^1$, the holonomy maps are \diffs{}.
We will need to control the derivatives of these holonomies. 
The derivatives of the elements of the family $\{ H^{s, \iin}_{\sigma',\sigma} \}_{\sigma, \sigma'}$ are not uniformly bounded.
Indeed, the width of the strip $B$ tends to 0 at the ends, and the map $H^{s, \iin}_{\sigma', \sigma}$ is very contracting when $\sigma'$ tends to an end of $B$ which accumulates on a compact leaf of $\cL^\iin$.
However, the following lemma states that one can uniformly bound the distortion of the holonomy of the foliation $\cG^{s, \iin}$ between two leaves of $\cG^{u, \iin} \cap B$.

\begin{figure}[htb]
    \centering
    \vspace*{-1em}
    \hspace*{-1em}
    \includegraphics[height=0.28\textheight]{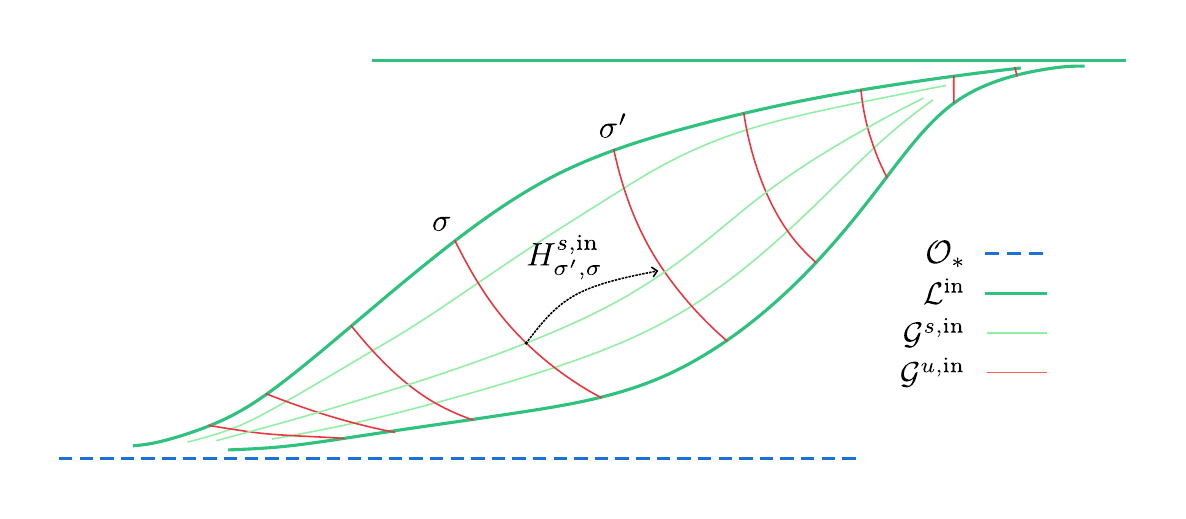}
    \vspace*{-1em}
    \caption{Holonomy of $\cG^{s, \iin}$ on a strip}
    \label{fig: holonomies in and out}
\end{figure}

\begin{lem} \label{lem: bounded distortion of the stable holonomy on entrance strip}
There exists a constant $\vartheta>0$ such that for any strip $B \subset \Pin \ssm \cL^\iin$, for any $\sigma, \sigma' \in \cG^{u, \iin} \cap B$, for any $p, q \in \sigma$
$$ \vartheta^\inv < \frac{\vert (H^{s, \iin}_{\sigma',\sigma})'(p) \vert}{\vert (H^{s, \iin}_{\sigma',\sigma})'(q) \vert} < \vartheta.$$
\end{lem}

We begin with the following fact, which is true on both \Pin{} and \Pout{}.
In order not to repeat it, we state it with the notation $^*$ which means $\iin$ or $\out$.

\begin{claim} \label{claim: constant holonomy neighborhood orbit}
Let $\cO_i$ be an orbit in $\cO_*$ and $\cV_i$ be a normalized neighborhood of $\cO_i$.
Then
\begin{enumerate}
    \item There exists a constant $\cst>0$ such that for all leaves $l_1$ and $l_2$ of $\cG^{s, *}$, if $\sigma$ and $\sigma'$ are two arcs of leaves of $\cG^{u, *} \cap \cV_i$ whose ends are on $l_1$ and $l_2$, then
    $\cst^\inv \leq {l(\sigma)}/{l(\sigma')} \leq \cst $.
    \label{claim: hol cst; it: length cst}
    
    \item The derivatives of the holonomy maps $\{ H^{s, *}_{\sigma',\sigma} \}_{\sigma, \sigma'}$ are uniformly bounded for $\sigma, \sigma' \in \cG^{u, *} \cap \cV_i$.
    \label{claim: hol cst; it: bounded derivative}
\end{enumerate}
\end{claim}

\begin{figure}[htb]
    \centering
    \vspace*{-1em}
    \includegraphics[height=0.3\textheight]{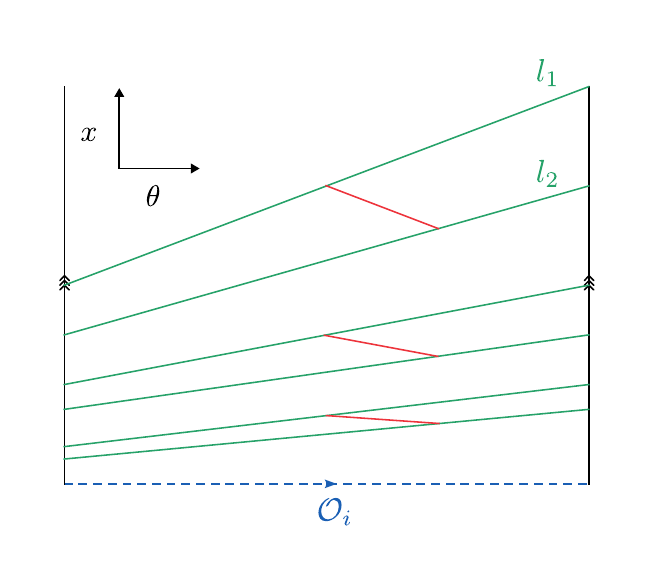}
    \vspace*{-1em}
    \caption{Arcs of leaves of $\cG^{u, *}$ between two leaves $l_1$ and $l_2$ of $\cG^{s, *}$ in a normalized neighborhood of $\cO_i$}
    \label{fig: constant holonomy}
\end{figure}

We refer to Figure~\ref{fig: constant holonomy} which illustrates Claim~\ref{claim: constant holonomy neighborhood orbit}, Item~\ref{claim: hol cst; it: length cst}.

\begin{proof}
Let us show this fact on \Pin{}, the proof is identical on \Pout{}.
Let $(\cV_i, \xi_i = (x,y,\theta))$ be a \ncs{} of $\cO_i$.
Let $\rho_i = (x, \theta) \in \R \times \R/\,\Z$ be the \ncs{} induced by $\xi_i$ on the boundary of $P$ (Remark~\ref{rmk: equation affine foliation}).
Let two leaves $l_1$ and $l_2$ of the foliation $\cG^{s, \iin}$ in $\cV$, and $p =(x, \theta) \in l_1$.
The lifts of the leaves in the lift $(x, \tilde \theta) \in \R \times \R$ of the coordinates $(x, \theta)$ have equation $x = c_i 2^{\tilde \theta}$.
It follows that the distance $\dist(p, l_2)$ decreases linearly in $x$.
Moreover the angle between the foliations $\cG^{s, \iin}$ and $\cG^{u, \iin}$ also decreases linearly in $x = \dist(p, \cO_i)$ according to Lemma~\ref{lem: angle foliation and vector field on boundary}.
We conclude that the ratio of the lengths ${l(\sigma')}/{l(\sigma)}$ of two leaf arcs $\sigma$ and $\sigma'$ of $\cG^{u, \iin}$ joining $l_1$ and $l_2$ in $\cV_i$ is uniformly bounded away from $0$ and $+\infty$ (Figure~\ref{fig: distance and angle}).

\begin{figure}[htb]
    \centering
    \vspace*{-1em}
    \includegraphics[height=0.28\textheight]{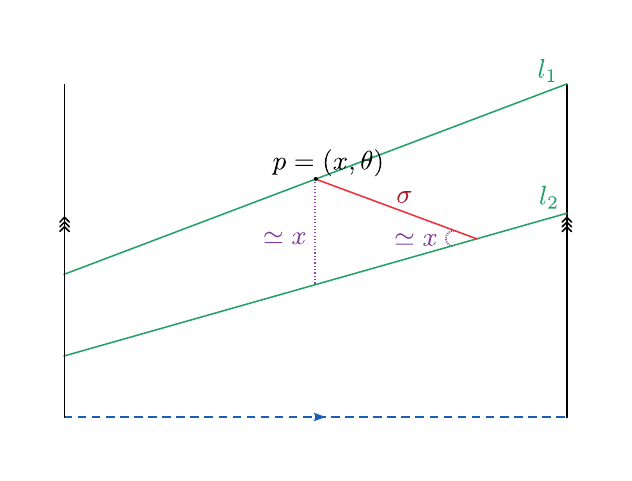}
    \vspace*{-1em}
    \caption{The distance between $l_1$ and $l_2$ and the angle between the foliations $\cG^{u, \iin}$ and $\cG^{s, \iin}$ decreases at the same rate}.
    \label{fig: distance and angle}
\end{figure}

\vspace*{-1em}

The second item is a consequence of the first.
Indeed, if $\sigma$ and $\sigma'$ are two arcs of leaves of $\cG^{u, \iin}$ in $\cV_i$, we can choose $\sigma_0$ an arc of leaf contained in $\sigma$ of arbitrarily small length around a point $p$, and $\sigma_0' = H^{s, \iin}_{\sigma', \sigma}$ its image contained in $\sigma'$.
Then the ratio of the lengths ${l(\sigma_0')}/{l(\sigma_0)}$ tends to the derivative of $H^{s, \iin}_{\sigma', \sigma}$ at point $p$ when the length of $\sigma_0$ tends to 0.
According to Item~\ref{claim: hol cst; it: length cst}, this ratio is uniformly bounded.
We deduce Item~\ref{claim: hol cst; it: bounded derivative}.
\end{proof}

\begin{figure}[tb]
    \centering
    \vspace*{-1em}
    \hspace*{-1.5em}
    \includegraphics[height=0.275\textheight]{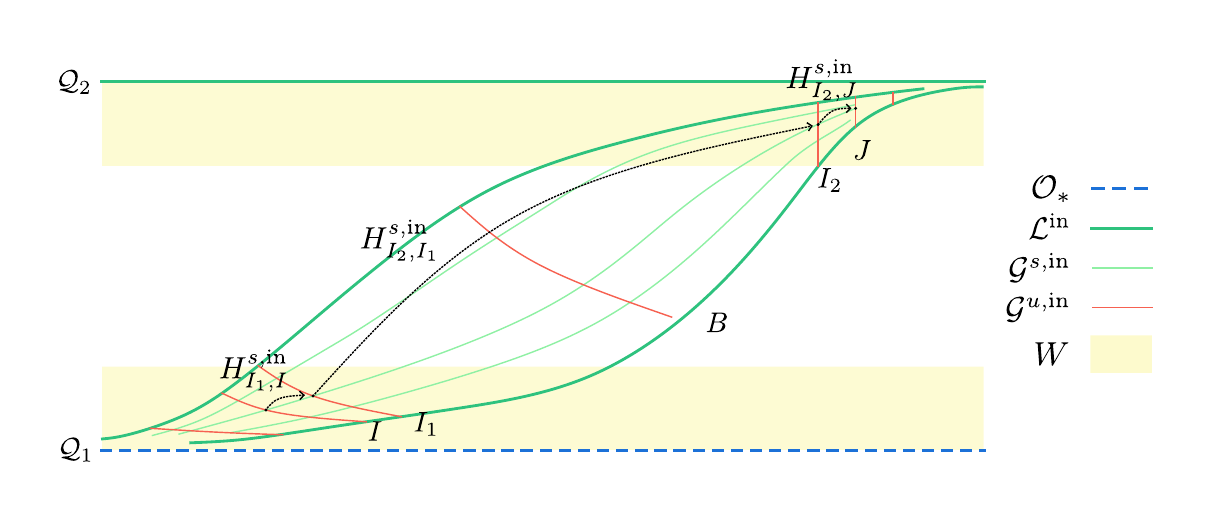}
    \vspace*{-2em}
    \caption{Decomposition of $H^{s, \iin}_{J,I}$ in $B$}
    \label{fig: decompo holonomy stable strip}
\end{figure}

\begin{proof}[Proof of Lemma~\ref{lem: bounded distortion of the stable holonomy on entrance strip}]
Let $\cQ$ be a compact leaf of $\cL = \cL^\iin \cup \cO_* \cup \cL^\out$ on which the strip $B$ accumulates.
Then $\cQ$ is a compact boundary leaf (Definition~\ref{def: boundary of a lamination}).
Let $W$ be a neighborhood of $\cQ$ on which the holonomy maps of the foliation $\cG^{s, \iin}$ are conjugate to affine maps (Proposition~\ref{prop: boundary foliation induced by paif}, Item~\ref{prop: boundary foliation, it: Gb}).
 If $\cQ \in \cO_*$ is a periodic orbit, $W$ is contained in a \ncs{} $(\cV_i, \xi_i)$ of $\cO_i = \cQ$.
    We use the fact that the holonomy maps $H^{s, \iin}_{\sigma',\sigma}$ have a uniformly bounded derivative in this neighborhood.
 If $\cQ$ is a compact leaf of $\cL^\iin$, the foliations $\cG^{u, \iin}$ and $\cG^{s, \iin}$ are uniformly transverse on $W$.
    By definition of $W$ the holonomy maps of a compact leaf of $\cL^\iin$ are conjugate to affine maps, hence the distortion of the holonomy maps $H^{s, \iin}_{\sigma',\sigma}$ is uniformly bounded for $\sigma, \sigma'$ in $W$.
    
Finally, it is enough to see that we can compose a holonomy map in the following way (Figure~\ref{fig: decompo holonomy stable strip}).
There exists a constant $l>0$, such that if $I, J$ are two leaves of $\cG^{u, \iin} \cap B$, then we can write
$$ H^{s, \iin}_{I, J} = H^{s, \iin}_{I, I_1} \circ H^{s, \iin}_{I_1, I_2} \circ H^{s, \iin}_{I_2, J} $$
with $H^{s, \iin}_{I, I_1}$ and $H^{s, \iin}_{I_2, J}$ holonomy maps in the neighborhood $W$,
and such that the leaves of $\cG^{s, \iin}$ have a length uniformly bounded by $l$ between $I_1$ and $I_2$ and are uniformly transverse to the leaves of $\cG^{u, \iin}$.
The derivatives of $H^{s, \iin}_{I_1, I_2}$ are uniformly bounded, and so is the composition.
Each map of the composition has a uniformly bounded distortion, and the bound is independent of the strip $B$.
The same is true for the composition.
\end{proof}

\pagebreak[4]

Let us quote the following lemma which is completely analogous to Lemma \ref{lem: bounded distortion of the stable holonomy on entrance strip}, where $H^{u, \out}_{\sigma', \sigma}$ denotes the holonomy map of the foliation $\cG^{u, \out}$ between two leaves $\sigma, \sigma'$ of $\cG^{s, \out} \cap B$.

\begin{lem} \label{lem sym: bounded distortion of the stable holonomy on entrance strip}
There exists a constant $\vartheta>0$ such that for any strip $B\!\subset\! \Pout \ssm \cL^\out$, for any $\sigma, \sigma' \in \cG^{s, \out} \cap B$, for any $p, q \in \sigma$
$$ \vartheta^\inv < \frac{\vert (H^{u, \out}_{\sigma',\sigma})'(p) \vert}{\vert (H^{u, \out}_{\sigma',\sigma})'(q) \vert} < \vartheta.$$
\end{lem}

\subsubsection*{Bound for the derivative of the holonomy of the stable foliation on a strip on the exit boundary}
\begin{figure}[htb]
    \centering
    \vspace*{-1em}
    \includegraphics[height=0.27\textheight]{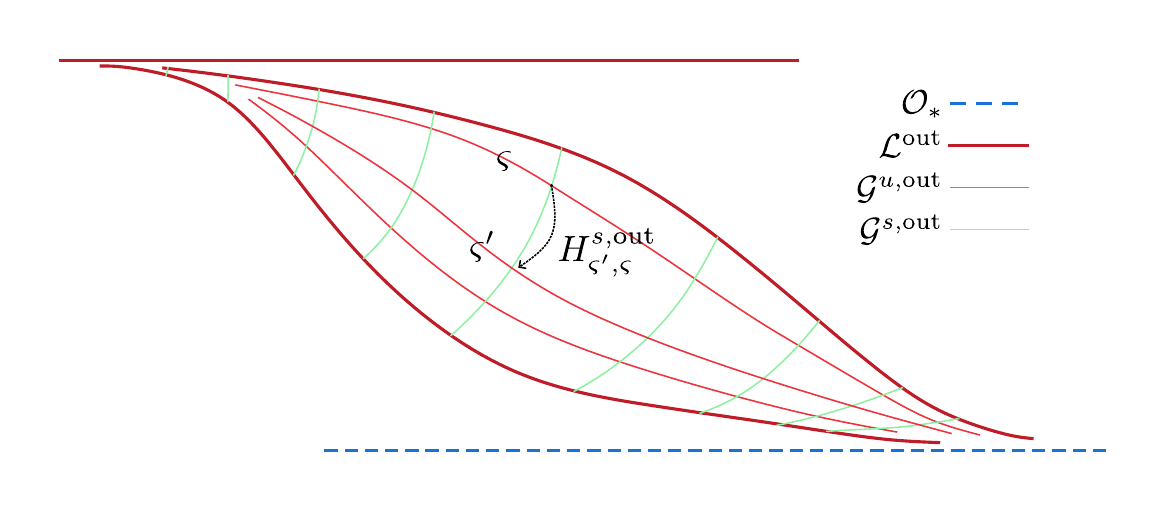}
    \vspace*{-1em}
    \caption{Holonomy map of the foliation $\cG^{s, \out}$ on a strip}
    \label{fig: holonomy out}
\end{figure}

.
Let $\varsigma$ and $\varsigma'$ be leaves of the foliation $\cG^{u, \out}$ on \Pout{} in a strip $B \subset \Pout \ssm \cL^\out$.
We define the holonomy map of the foliation $\cG^{s, \out} \cap B$ between the leaves $\varsigma$ and $\varsigma'$ (Figure~\ref{fig: holonomy out}) by $H^{s, \out}_{\varsigma', \varsigma} \colon \varsigma \to \varsigma'$.
These maps are well defined because each arc of leaf of $\cG^{s, \out} \cap B$ intersects each leaf $\varsigma$ of $\cG^{u, \iin} \cap B$.
These maps are \diffs{}.
Let us show that the derivatives of the holonomy maps $\big\{ H^{s, \out}_{\varsigma', \varsigma}\big \}$ for ${\varsigma', \varsigma}$ leaves of \Pout{} in a same strip $B \subset \Pout \ssm \cL^\out$ are uniformly bounded.

\begin{lem} \label{lem: bounded derivative of the stable holonomy on exit strip}
There exists a constant $\alpha\!>0\!$ such that for any strip $B\!\subset\! \Pout \ssm \cL^\out$, for any $\varsigma, \varsigma' \in \cG^{u, \out} \cap B$, for any $p \in \varsigma$, we have
$ \alpha^\inv < \vert (H^{s, \out}_{\varsigma', \varsigma})'(p) \vert < \alpha $.
\end{lem}

\begin{proof}
The family $\big\{ H^{s, \out}_{\varsigma', \varsigma} \big\}_{\varsigma',\varsigma}$ is a relatively compact family of \diffs{}.
Indeed, let $I$ be a leaf of $\cG^{s, \out} \cap B$.
Each leaf $\varsigma$ of $\cG^{s, \out}$ intersects $I$ once.
We can therefore parameterize the family $\big\{ H^{s, \out}_{\varsigma', \varsigma} \big\}_{\varsigma',\varsigma}$ by the product $I \times I$.
We conclude by noticing that $I$ is a relatively compact arc.
It is then sufficient to show that the derivative of a \diff{} $H^{s, \out}_{\varsigma', \varsigma}$ is bounded on $\varsigma$.
On the complementary $\Pout \ssm \cV_*$ of the normalized neighborhoods of the periodic orbits $\cO_*$, the angle of the foliation $\cG^{u, \iin}$ with $\varsigma$ is bounded, and the distance between $\varsigma$ and $\varsigma'$ is increased.
On a normalized neighborhood $\cV_i$ of an orbit $\cO_i$, the derivative of $H^{s, \out}_{\varsigma', \varsigma}$ is bounded according to Claim~\ref{claim: constant holonomy neighborhood orbit}.
We deduce that the derivative of $H^{s, \out}_{\varsigma', \varsigma}$ is bounded on $\varsigma$.
Finally, let us note that we can uniformly bound the derivatives of $\big\{H^{s, \out}_{\varsigma', \varsigma} \big\}_{\varsigma',\varsigma}$ independently of the strip $B$.
Indeed, there exists a relatively compact leaf arc $I$ of $\cG^{s, \out}$ such that each leaf $\varsigma$ of $\cG^{u, \out}$ intersects $I$ only once.
It follows that the family $\big\{ H^{s, \out}_{\varsigma', \varsigma} \big\}_{\varsigma',\varsigma}$ for $\varsigma, \varsigma' \in \cG^{u, \out}$ is a relatively compact family of \diffs{}, and the previous proof works.
\end{proof}

Let us quote the following lemma which is completely analogous to Lemma \ref{lem: bounded derivative of the stable holonomy on exit strip}, where $H^{u, \iin}_{\varsigma', \varsigma}$ denotes the holonomy of the foliation $\cG^{u, \iin}$ between two leaves of $\cG^{s, \iin} \cap B$.

\begin{lem} \label{lem sym: bounded derivative of the stable holonomy on exit strip}
There exists a constant $\alpha>0$ such that for any strip $B \subset \Pin \ssm \cL^\iin$, for any $\varsigma, \varsigma' \in \cG^{s, \iin} \cap B$, for any $p \in \varsigma$, we have
$ \alpha^\inv < \vert (H^{u, \iin}_{\varsigma', \varsigma})'(p) \vert < \alpha $.
\end{lem}

\subsection{Proof of Lemma~\ref{lem: spread expansion strip}}
 \label{sec: spread; subsec: proof lemma strip}

Let $B$ be a strip on $\Pin{}$. Let $\lambda>1$, $\epsilon>0$ and $\delta>0$ be fixed.
Let $l_1, l_2$ be the leaves of $\cL^\iin$ forming the accessible boundary of $B$.
The proof will be done in three steps.

\begin{itemize}[leftmargin=*]
    \item In Step 1, we fix a leaf $\sigma$ of $\cG^{u, \iin} \cap B$ that crosses the strip $B$ and we construct a \diff{} $h: \sigma \to \sigma$
    equal to the identity at the ends of the leaf,
    such that the composition $\foutin \circ h: \sigma \to \foutin(\sigma)$ has a derivative bounded from below by a given factor proportional to $\lambda$.
    This is Lemma~\ref{lem: existence h}.
    
    \item In Step 2, we push the diffeomorphism $h$ along the foliation $\cG^{s, \iin}$ by the holonomy maps $H^{s, \iin}_{\sigma', \sigma}$ to obtain a diffeomorphism $\hat\psi \colon B \to B$, equal to the identity on a neighborhood of $l_1$ and $l_2$.
    This diffeomorphism naturally preserves the foliation $\cG^{u, \iin}$ (leaf-to-leaf) and the foliation $\cG^{s, \iin}$ (globally) on $B$.
    Moreover, the expansion of the composition $\foutin \circ \hat\psi$ in the direction $\cG^{u, \iin}$ on $B$ is essentially equal to the expansion of the composition of $\foutin \circ h$.
    Indeed, this is a consequence of the uniform bounds on the distortion and the derivative of the holonomy maps of the foliation $\cG^s$ between two leaves of $\cG^u$ from Lemmas
   ~\ref{lem: bounded distortion of the stable holonomy on entrance strip} and~\ref{lem: bounded derivative of the stable holonomy on exit strip}.
    We chose $h$ in order to compensate the contraction of these holonomy maps.
    Thus the composition $\foutin \circ \hat\psi$ expands by a factor $\lambda$ the direction tangent to $\cG^{u, \iin}$ on the strip $B$.
    The \diff{} thus constructed on $B$ satisfies all items of Lemma~\ref{lem: spread expansion strip},
    except that its support is not a compact of the strip $B$ (Item~\ref{lem: spread strip; it: support}).
    
    \item Step 3 consists in \say{slowing down} the \diff{} $\hat\psi$ in order to obtain a \diff{} $\psi$ equal to the identity at the ends of $B$.
    This is where the parameters $\epsilon$ and $\delta$ come in.
    \begin{itemize}
        \item In order not to destroy what has been done previously we slow down by a leaf-to-leaf barycentric isotopy along $\cG^{u, \iin}$.
        \item The speed of the \say{slowing down} controls the perturbation of the foliation $\cG^{s, \iin}$, and will therefore be chosen according to the parameter $\epsilon$.
        \item If one end of $B$ is a periodic orbit $\cO_i$, we lose control on the expansion of the composition with $\foutin$ at the location of the slowdown.
        We make the slowing down in a neighborhood of the periodic orbit of size $\delta$.
        Otherwise, the ends of the strip $B$ accumulates on a compact leaf $\cQ_i$ of $\cL^\iin$, and we make the slowdown in a neighborhood of $\cQ_i$ small enough so that the \say{natural} expansion of $\foutin$ obtained by Proposition~\ref{prop: crossing map} compensates the contracting effects of the slowdown.
        This preserves the expansion of the composition by a factor $\lambda$ on this end.
    \end{itemize}
\end{itemize}

\subsubsection*{Step 1: On a leaf}
Let $\sigma$ be a leaf of $\cG^{u, \iin} \cap B$, of \emph{quasi-maximal} length among all the leaves of $\cG^{u, \iin} \cap B$, more precisely we ask that for any leaf $\sigma'$ of $\cG^{u, \iin} \cap B$, we have $l(\sigma)\geq \frac{1}{2} l(\sigma')$.
This is possible because the leaves of $\cG^{u, \iin} \cap B$ are of uniformly bounded length.
We fix $\sigma$ for the rest of the section.
It is an arc whose ends are on the non-compact leaves $l_1$ and $l_2$ of $\cL^\iin$ which border $B$.
Let $\vartheta>0$ and $\alpha>0$ be the constants which satisfy Lemmas~\ref{lem: bounded distortion of the stable holonomy on entrance strip} and~\ref{lem: bounded derivative of the stable holonomy on exit strip}.

\begin{lem} \label{lem: existence h}
There exists a \diff{} $h \colon \sigma \to \sigma$, equal to the identity on a neighborhood of the ends of $\sigma$ such that
$$ \forall p \in \sigma,
\qq \forall v \in T_p\sigma,
\qq  \Vert (\foutin \circ h)_* v \Vert > 2 \vartheta \alpha \lambda \Vert v \Vert.$$
\end{lem}

\begin{proof}
As $\partial \sigma \in \cL^\iin$, we know from Proposition~\ref{prop: crossing map} that the derivative of $\res{\foutin}{\sigma}$ tends to infinity at the ends of $\sigma$.
The unit-speed parametrization identifies $\sigma$ to $]0,l[$ for some $l>0$, and $\foutin(\sigma)$ to $\R$ by isometries.
We come to the following result, whose proof is left to the reader.

\begin{claim}
If $f \colon \, ]0,l[ \, \to \R$ is a \diff{} such that
$ \vert f'(x) \vert \underset{x \to 0, \, l}{\longrightarrow} + \infty $ then for all $A>1$, there exists a \diff{} $f_A \colon \, ]0,l[ \, \to \R$ such that
\begin{itemize}
    \item $f_A = f$ on a neighborhood of $0$ and $l$,
    \item $\vert f_A'(x) \vert > A$ for all $x$ in $]0,l[$.
\end{itemize}
\end{claim}
We apply this result with $A = 2 \vartheta \alpha \lambda >1$, and we put
$ h:= \foutin^\inv \circ (\foutin)_A $.
\end{proof}

\subsubsection*{Step 2: On the strip}
We extend the \diff{} $h$ to the strip $B$ by pushing through the holonomy of the foliation $\cG^{s, \iin}$ between the leaves of $\cG^{u, \iin} \cap B$.
More precisely, we define for any $p \in B$
\begin{equation} \label{eq: hat psi}
    \hat\psi (p) = H^{s, \iin}_{\sigma_p, \sigma} \circ h \circ H^{s, \iin}_{\sigma, \sigma_p} (p),
\end{equation}
where $\sigma_p$ denotes the leaf of $\cG^{u, \iin} \cap B$ passing through $p$.
This expression defines a \diff{} on $B$.
It satisfies the following properties:

\begin{claim}
\label{claim: psi hat}
\mb{}
\begin{enumerate}
    \item \label{claim: psi hat; it: egal id} $\hat\psi$ is equal to the identity on a neighborhood of the leaves $l_1$ and $l_2$ which border $B$;
    \item \label{claim: psi hat; it: preserve foliation u} $\hat\psi$ preserves the foliation $\cG^{u, \iin} \cap B$ leaf-to-leaf;
    \item \label{claim: psi hat; it: preserve foliation s} $\hat\psi$ preserve the foliation $\cG^{s, \iin} \cap B$ globally;
    \item \label{claim: psi hat; it: expansion composition} The composition $\foutin \circ \hat\psi$ expands the norm of the vectors in the direction $\cG^{u, \iin}$ by a factor $\lambda$:
    $ \forall p \in B, \quad \forall v \in T_p \cG^{u, \iin}, \quad \Vert (\foutin \circ \hat\psi)_* v \Vert \geq \lambda \Vert v \Vert. $
\end{enumerate}
\end{claim}

\begin{proof}$\,$
\begin{enumerate}
    \item The neighborhood in question is the saturated set of the complementary of the support of $h$ by the foliation $\cG^{s, \iin}$.
    It is a neighborhood of the leaves $l_1$ and $l_2$ by Lemma~\ref{lem: existence h}.
    \item Obvious.
    \item Obvious.
    \item Let $\varsigma = \foutin(\sigma)$.
The map $\foutin$ commutes with the holonomy maps of $\cG^{s, \iin}$ and $\cG^{s, \out}$ because it maps the pair $(\cG^{s, \iin}, \cG^{u, \iin})$ on the pair \linebreak[4]$(\cG^{s, \out}, \cG^{u, \out})$ on each strip.
We deduce that
$$(\foutin \circ \hat\psi) (p) = H^{s, \out}_{\varsigma_p,\varsigma} \circ \l( \foutin \circ h \r) \circ H^{s,in}_{\sigma, \sigma_p} (p)$$
where $\varsigma_p = \foutin(\sigma_p)$ denotes the leaf of $\cG^{u, \out}$ passing through $\foutin(p)$.
We refer to Figure~\ref{fig: spreading composition passage}.

\begin{figure}[htb]
    \centering
    \vspace*{-1em}
    \includegraphics[height=0.43\textheight]{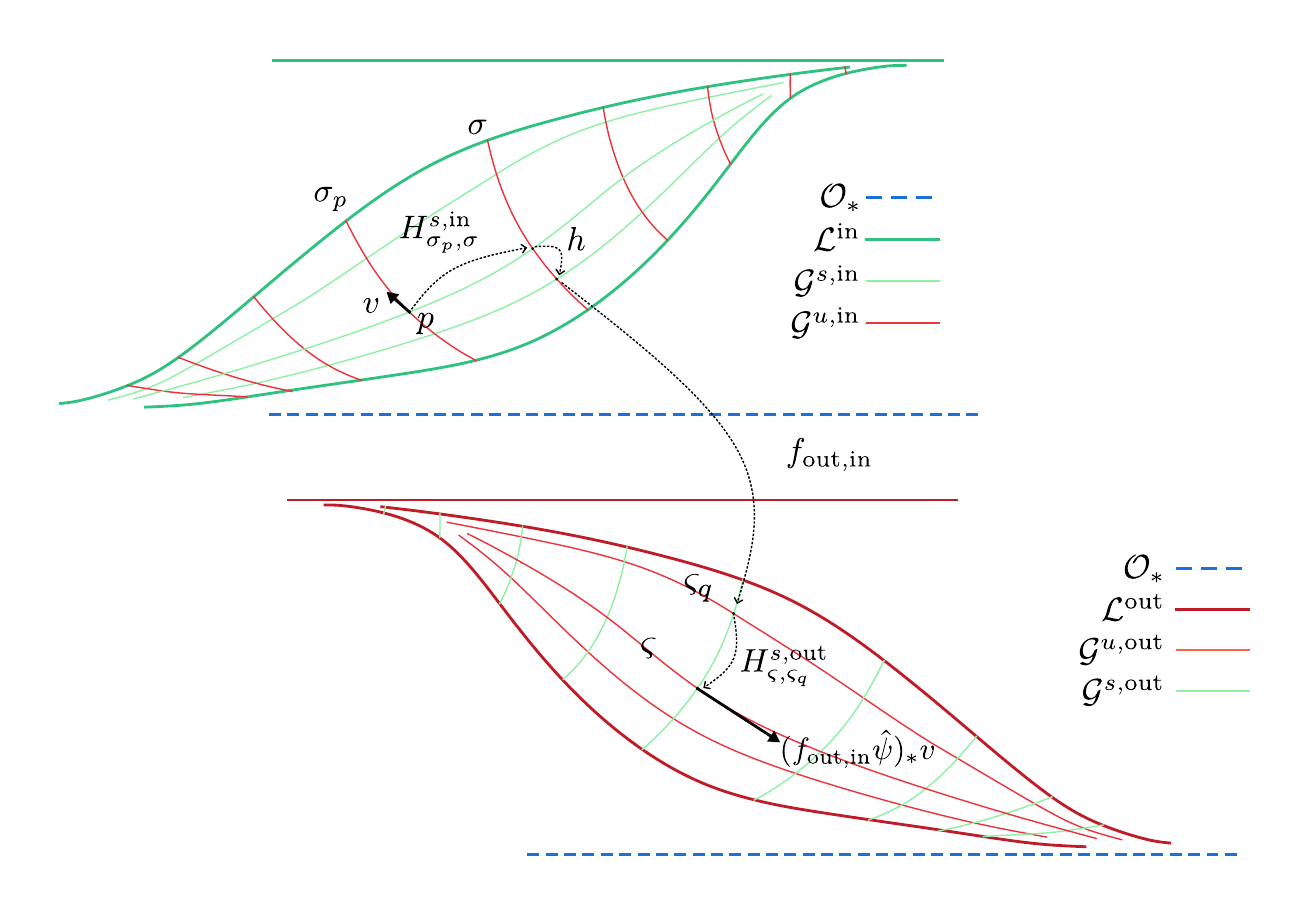}
    \vspace*{-1em}
    \caption{Action of $\foutin \circ \hat \psi$}
    \label{fig: spreading composition passage}
\end{figure}

Let $\vartheta$ be the constant given by Lemma~\ref{lem: bounded distortion of the stable holonomy on entrance strip}. 
Then, if $l(\sigma)$ denotes the length of the leaf $\sigma$, we have
$$
    | (H^{s,in}_{\sigma, \sigma_p})' | \  > \  \vartheta^\inv \frac{l(\sigma)}{l(\sigma_p)} \ \geq \  \frac{1}{2} \vartheta^\inv
$$
because $\sigma$ is of quasi-maximal length among the leaves of $\cG^{u, \iin} \cap B$ (see Step 1).
Moreover, by Lemma~\ref{lem: bounded derivative of the stable holonomy on exit strip},
we have $\vert (H^{s, \out}_{\varsigma,\varsigma_p})' \vert > \alpha^\inv$.
So for all $p \in B$, for all $v \in T_p \sigma_p$, if $q = H_{\sigma,\sigma_p}(p)$,
$$ \Vert (\foutin \circ \hat\psi)_* v \Vert > \frac{1}{2} \vartheta^\inv \alpha^\inv \vert (\foutin \circ h)'(q) \vert \Vert v \Vert.$$
We conclude with Lemma~\ref{lem: existence h}.\qedhere
\end{enumerate}
\end{proof}

\subsubsection*{Step 3: Slowing down at the ends of the strip}

The \diff{} $\hat\psi \colon B \to B$ (Claim~\ref{claim: psi hat}) satisfies all items of Lemma~\ref{lem: spread expansion strip} on the strip $B$, except Item~\ref{lem: spread strip; it: support}: it is not compactly supported, because it is not equal to the identity at the ends of $B$.
To solve this problem, we will \say{slow down} $\hat\psi$ on a neighborhood of the ends of $B$.
We can choose the slowing down area in order to keep the expansion of the composition $\foutin \circ \hat\psi$ on $\cG^{u, \iin}$ at the end of $B$ if it accumulates on a compact leaf of $\cL^\iin$.
If the end accumulates on a periodic orbit $\cO_i \in \cO_*$, one loses the expansion of $\foutin \circ \hat\psi$ on the slowing down zone.
We will then choose this zone in a $\delta$-neighborhood of $\cO_i$.
In fact, we will operate the slowing down on the $\hat\psi^\inv$ diffeomorphism in order to control the perturbation of the foliation $\cG^{s, \iin}$ by $\hat\psi^\inv$
(Item~\ref{lem: spread strip; it: foliation perturbed} of Lemma~\ref{lem: spread expansion strip}).
Recall that $l_1$ and $l_2$ are the non-compact leaves of $\cL^\iin$ which form the accessible boundary of $B$.
Let $s \colon \R \to l_1$ be a parameterization by length of the leaf $l_1$,
where $s(0)= l_1 \cap \sigma$.
Denote $\sigma_t$ the leaf of $\cG^{u, \iin} \cap B$ such that $\sigma_t \cap l_1 = s_1(t)$.
We have $\sigma = \sigma_0$.
As $\hat\psi$ preserves the foliation $\cG^{u, \iin}$ leaf-to-leaf on $B$, it induces the existence of a continuous one-parameter family of \diffs{} $\hat\psi_t = \res{\hat\psi}{\sigma_t} \colon \sigma_t \to \sigma_t$.

Let $W$ be a neighborhood of the ends of $B$ on which the holonomy maps of the foliation $\cG^{s, \iin}$ are conjugate to affine maps (Remark~\ref{rmk: nbh of end of strip with affine holonomy}).
For any leaf $\sigma_t$ contained in $W$, we define $\Psi = \{ \Psi_{t, \tau} \colon \sigma_t \to \sigma_t \tau \in [0,1], t \in \R \}$ a family of barycentric isotopies between $\hat\psi_t^\inv$ and the identity on $\sigma_t$ on each leaf $\sigma_t$: 
\begin{equation}\label{eq: Psi_t,tau}
    \forall p \in \sigma_t, \qq \forall \tau \in [0,1], \qq \Psi_{t, \tau}(p) = \tau \hat\psi_t^\inv(p) + (1-\tau) p.
\end{equation}
This formula makes sense because the leaf $\sigma_t$ contained in $W$ has an affine structure, and defines a family $\Psi$ of \diffs{}.

\begin{figure}[htb]
    \centering
    \vspace*{-1em}
    \includegraphics[height=0.20\textheight]{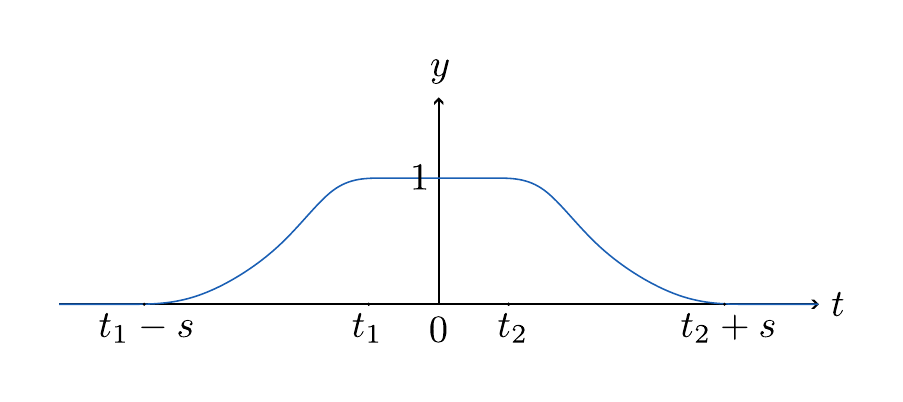}
    \vspace*{-1em}
    \caption{Graph of the function $\tau \colon \R \to [0,1]$}
    \label{fig: tau graph}
\end{figure}

Let $t_1<0$, $t_2>0$, such that $\sigma_t \in W \cap B$ for all $t \in \ ]- \infty, t_1[ \  \cup \ ]t_2, +\infty[$.
Let $s>0$.
Let \mb{$\tau = \tau(t_1, t_2, s) \colon \R \to [0,1]$} be a plateau function of class $\cC^1$ satisfying the following properties (Figure~\ref{fig: tau graph}):
\begin{enumerate}
    \item $\tau = 1$ on $[t_1, t_2]$ and $\tau = 0$ on $ \R \ssm [t_1-s, t_2 +s]$;
    \label{it: support of tau}
    \item $\tau$ is increasing on $\R_-$ and decreasing on $\R_+$, and the norm of its derivative is upper bounded by $2s^\inv$. \label{it: bounded derivative of tau}
\end{enumerate}

\begin{figure}[htb]
    \centering
    \vspace*{-1em}
    \includegraphics[height=0.27\textheight]{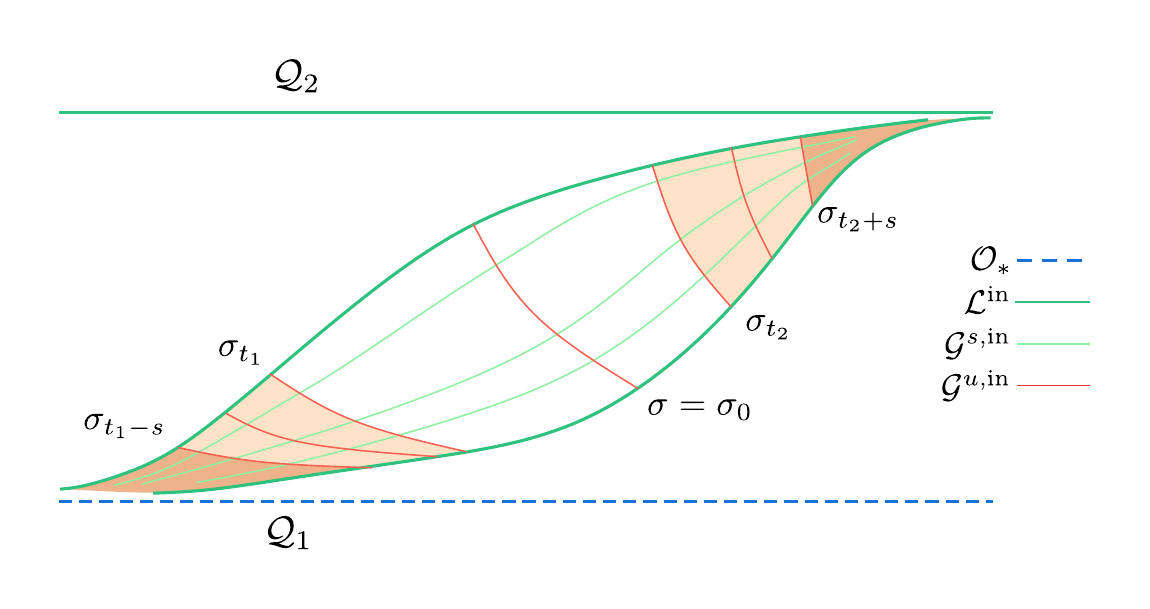}
    \vspace*{-1em}
    \caption{Area of the start of the isotopy (light) and the end of the isotopy (dark) on the strip $B$}
    \label{fig: slowing down area}
\end{figure}

The parameters $t_i$ determine the area of the beginning of the isotopy in the strip, and the parameters $s$ determines the speed of the isotopy (Figure~\ref{fig: slowing down area}).

For $p \in B$, if $\sigma_t$ is the leaf of $\cG^{u, \iin} \cap B$ passing through $p$, we define
\begin{equation} \label{eq: psi}
    \psi(p) = \l( \Psi_{t, \tau(t)} (p) \r)^\inv.
\end{equation}
This expression defines a diffeomorphism $\psi \colon B \to B$.

\begin{proof}[Proof of Lemma~\ref{lem: spread expansion strip}]
Let us show each item for the \diff{} $\psi \colon B \to B$ defined by Equation~\eqref{eq: psi}.
\begin{enumerate}[leftmargin=*]
    \item By definition of $\tau$ (Item~\ref{it: support of tau} above), the support of $\psi$ is a compact of the strip $B$ equal to
    \[\operatorname{Supp} (\psi) =  \operatorname{Supp}(\hat\psi) \ssm 
    \l( \l( \underset{-\infty}{\overset{t_1-s}{\bigcup}} \sigma_t \r)
    \cup \l( \underset{t_2 +s}{\overset{+\infty}{\bigcup}} \sigma_t \r) \r).\]
    \item It is clear that this diffeomorphism preserves the foliation $\cG^{u, \iin}$ leaf-to-leaf.
    \item Let $(\partial_u, \partial_s)$ be a pair of two unit vector fields tangent to $\cG^{u, \iin}$ and $\cG^{s, \iin}$ respectively.
    Let $p \in B$, let $t$ be a parameter such that $\sigma_t$ is the leaf of $\cG^{u, \iin} \cap B$ passing through $p$, let $q = \psi^\inv(p)$.
    Then the image of $\partial_s$ by the differential of $\psi^\inv$ is 
    $d_p \psi^\inv (\partial_s) = \tau'(t) \, \partial_u + (H^{u, \iin}_{\varsigma_p, \varsigma_q})'(p) \, \partial_s. $
    According to Lemma~\ref{lem sym: bounded derivative of the stable holonomy on exit strip}, the derivatives of these holonomy maps are uniformly bounded.
    It suffices to choose a function $\tau$ such that the derivative is close enough to 0, in other words a parameter $s$ large enough (Item~\ref{it: bounded derivative of tau} of the definition of $\tau$ above), so that
    $\vert \tau'(t) \vert \leq \epsilon \alpha^\inv$
    where $\alpha$ is the constant given by Lemma~\ref{lem sym: bounded derivative of the stable holonomy on exit strip} which bounds from below the derivative of $H^{u, \iin}_{\varsigma_p, \varsigma_q}$.
    Thus the image $\psi^\inv_* \partial_s$ is in a $(\epsilon,\cG^{s, \iin}/\, \cG^{u, \iin})$-cone field (Definition~\ref{def: cone slope}). This proves Item~\ref{lem: spread strip; it: foliation perturbed}.
    
    \item According to Equation~\eqref{eq: Psi_t,tau} of the barycentric isotopy in each leaf of $\cG^{u, \iin}$, the derivative of $\psi^\inv$ in the direction of $\cG^{u, \iin}$ is bounded by the upper and lower bounds of the derivatives of the family of \diffs{} $\{ \hat\psi_t \}_t$, in other words by the upper and lower bounds of the derivative of $\hat\psi$ in the direction $\cG^{u, \iin}$.
    Recall that $\hat\psi$ (Equation~\eqref{eq: hat psi}) is the conjugate by $H^{s, \iin}_{\sigma, \sigma'}$ of the \diff{} $h \colon \sigma \to \sigma$ (Lemma~\ref{lem: existence h}).
    Since the distortion of the holonomies $H^{s, \iin}_{\sigma', \sigma}$ is uniformly bounded by a constant $\vartheta>0$, (Lemma~\ref{lem: bounded distortion of the stable holonomy on entrance strip}), we deduce that the derivative of $\hat\psi$ in the direction of $\cG^{u, \iin}$ is uniformly controlled by the derivative of $h$,
    $$ \forall v \in T \cG^{u, \iin}, \qq \vartheta^\inv \vert h' \vert \leq \Vert \hat\psi_* v \Vert \leq \vartheta \vert h' \vert.$$
    Since $h$ is a diffeomorphism with compact support, we conclude that the differential of $\psi$ in the direction $\cG^{u,\iin}$ is uniformly bounded by a constant which only depends on $\lambda$.
    We thus obtain Item~\ref{lem: spread strip; it: derivative bounded}.
    \item Let us show that we can choose $t_1, t_2$ in order to obtain Item~\ref{lem: spread strip; it: composition expansion}.
    Let $c>0$ be a constant which bounds from below the derivative of $\psi$ in the direction $\cG^{u, \iin}$.
    This constant exists according to the previous Item~\ref{lem: spread strip; it: derivative bounded} and is independent of $t_i$.
    Denote $B_1$ the \nbh{} of the end of $B$ corresponding to $t\leq t_1$, and $B_2$ the \nbh{} of the end of $B$ corresponding to $t \geq t_2$.
    \begin{enumerate}
    \item If the end $B_i$ accumulates on a periodic orbit $\cQ_i \in \cO_*$, then choose $t_i$ such that $\sigma_{t_i}$ is contained in a $\delta$-neighborhood $\cV_\delta \subset \cV_*$ of $\cQ_i$, where we recall that $\cV_*$ is the union of the normalized neighborhoods of the periodic orbits of $\cO_*$.
    \item If the end $B_i$ accumulates on a compact leaf $\cQ_i$ of $\cL^\iin$ then choose $t_i$ such that $\sigma_{t_i}$ is contained in the neighborhood $\cW^\iin_{\lambda c^\inv}$ given by Proposition~\ref{prop: crossing map}.\footnote{We can have $\cQ_1 = \cQ_2$ in the non-orientable case.}
\end{enumerate}
Such a choice is always possible because $\sigma_{t} \to \cQ_i$ when $t \to \pm \infty$, and because the neighborhoods $\cW^\iin_\lambda$ of Proposition~\ref{prop: crossing map} are neighborhoods of compact leaves of $\cL^\iin$.
Then we have

\begin{itemize}
    \item On $\underset{[t_1, t_2]}{\bigcup} \sigma_t$, we have $\psi = \hat\psi$ and the expansion of the composition $\foutin \circ \psi$ in the direction $\cG^{u, \iin}$ is satisfied according to Claim~\ref{claim: psi hat}, Item~\ref{claim: psi hat; it: expansion composition}.
    \item By construction, we have
    $$B \ssm \l( \l( \underset{[t_1, t_2]}{\bigcup} \sigma_t \r) \cup \cV_\delta \r)  = \l(\underset{\R \ssm [t_1, t_2]}{\bigcup} \sigma_t \r) \ssm \cV_\delta \ \subset \ \cW_{\lambda c}.$$
    According to the previous Item~\ref{lem: spread strip; it: foliation preserved}, the foliation $\cG^{u, \iin}$ is preserved by $\psi$. 
    Hence for any $v \in T \cG^{u, \iin}$ at a point of this set, we have by Proposition~\ref{prop: crossing map} and the previously shown Item~\ref{lem: spread strip; it: derivative bounded} of Lemma~\ref{lem: spread expansion strip}:
    $$\Vert (\foutin \circ \psi)_*v \Vert > \lambda c^\inv \Vert \psi_*v \Vert > \lambda \Vert v \Vert.$$
\end{itemize}
    \item By construction, the \diff{} $\psi$ coincides with $\hat \psi$ on the complementary of the normalized $\cV_*$ neighborhoods of the tangent $\cO_*$ periodic orbits in $B$, hence it is independent of $\delta$ in this region.
    Item~\ref{lem: spread strip; it: independant delta} is thus verified.\qedhere
\end{enumerate}
\end{proof}

\section{Choice of the parameters and cone criterion}
\label{sec: parameters and cones}
In this section and for the rest of the proof, $\varphi$ denotes a normalized gluing map of $(P,X)$ (Definitions~\ref{def: normalized gluing map}). 
We refer to the beginning of Section~\ref{sec: crossing map} for a reminder of the notations.
By definition of a normalized block (Definition~\ref{def: normalized block}, Item~\ref{def: normalized, it: affine section}), the block $(P,X)$ is provided with an affine section $\Sigma$,\footnote{Such a section is defined in an extension of $(P,X)$, but we consider in this chapter its restriction to $P$} in other words a transverse local section of $\Lambda_s$ such that the first return map of the flow of $X$ is affine in the neighborhood of $\Sigma \cap \Lambda_s$.

In the previous section we built a family of \diffs{} \mb{$\psi^\iin_\led \colon \pP \to \pP$} with support in $\Pin$ and \mb{$\psi^\out_\led \colon \pP \to \pP$} with support in $\Pout$ indexed by parameters $\lambda>1$, $\epsilon>0$, $\delta>0$, which satisfy the hypotheses of Proposition~\ref{prop: spread expansion} and Proposition~\ref{prop sym: spread expansion}.
Recall that:
\begin{itemize}
    \item $\lambda>1$ is a lower bound for the expansion factor of the composition $\foutin \circ \psi^\iin$ in the direction $\cG^{u, \iin}$ outside a neighborhood of the orbits~$\cO_*$,
    \item $\epsilon>0$ measures the perturbation of the foliation $\cG^{s, \iin}$ under the action of $(\psi^\iin)^\inv$,
    \item $\delta>0$ is the size of the neighborhood of the orbits $\cO_*$ on which we no longer control the dilation of the composition $\foutin \circ \psi^\iin$.
\end{itemize}
Symmetrically,
\begin{itemize}
    \item $\lambda>1$ is a lower bound for the expansion factor of the composition $(\psi^\out \circ \foutin)^\inv$ in the direction $\cG^{s, \out}$ outside a neighborhood of the orbits $\cO_*$,
    \item $\epsilon>0$ measures the perturbation of the foliation $\cG^{s, \iin}$ under the action of $\psi^\out$,
    \item  $\delta>0$ is the size of the neighborhood of the $\cO_*$ orbits on which we no longer control the dilation of the composition $(\psi^\out \circ \foutin)^\inv$ on $\cG^{s, \out}$.
\end{itemize}

\subsection{Statement of the main proposition and summary of the proof}
\label{sec: parameters; subsec: statement proposition}

\subsubsection*{Gluing map}
For any parameter \led{}, define the composition 
$$\psi_\led := \psi^\iin_\led \, \varphi \, \psi^\out_\led \, \colon \pP \to \pP,$$
where $\varphi$ is a normalized gluing map of $(P,X)$.
For any \led{}, $\psi^\iin$ and $\psi^\out$ are equal to the identity on a neighborhood of $\cO_*$ in $\pP$, so the same holds for $\psi$ and $\psi$ is a dynamical gluing map of $(P,X)$ (Definition~\ref{def: dynamical gluing map} and Proposition~\ref{prop: normalized gluing map is dynamic}).
Moreover, $\psi^\iin$ and $\psi^\out$ are isotopic to the identity via \diffs{} equal to the identity on a neighborhood of the boundary lamination $\cL$, so $\psi_\led$ and $\varphi$ are isotopic via strongly transverse gluing maps.

\begin{claim} \label{claim: modify triple strongly isotopic}
For all $\led$, if $\psi = \psi_\led$, then
$(P,X,\psi)$ is a triple strongly isotopic to $(P,X,\varphi)$, and the quotient space $P_\psi := P/\,\psi$ is a closed manifold of dimension three with a vector field of class $\cC^1$ induced by $X$ on $P_\psi$ which is denoted by~$X_\psi$.
\end{claim}

\begin{proof}
Indeed, $\psi = \psi_\led$ is a gluing map of $(P,X)$ isotopic to $\varphi$.
Moreover, there is an isotopy $h_t \colon \pP \to \pP$ which joins $h_0 = \Id$ and $h_1 = \psi^\out \circ (\psi^\iin)^\inv$, which preserves the boundary lamination $\cL$, and such that $h_1$ maps leaves of $\psi_* \cL$ on leaves of $\varphi_* \cL$.
The triples are thus strongly isotopic in the sense of Definition~\ref{def: strongly isotopic triples}.
The last affirmation relies on the fact that $\psi$ is a dynamical gluing map and Claim~\ref{claim: dynamical gluing map induce vector field}.
\end{proof}

Let $\pi_\psi$ be the projection $P \to P_\psi$.
In order not to complicate the notations, denote $\Pin_\psi = \pi_\psi(\Pin)$, $\Sigma_\psi = \pi_\psi(\Sigma)$, $\pi_\psi(\cO_*) = (\cO_*)_\psi$, etc.
More generally the index $_\psi$ indicates that one considers an element projected in the quotient manifold~$P_\psi$.

\subsubsection*{Local section in $P_\psi$}
The projection
$ \pi_\psi(\Sigma \cup \Pin \cup \Pout) = \pi_\psi(\Sigma \cup \Pin) = \pi_\psi(\Sigma \cup \Pout) $
is a non-compact, immersed surface in $P_\psi$, transverse to the orbits of $X_\psi$.
Denote
$S_{0, \psi} := \pi_\psi (\Pin \cup \Sigma)$.
Let $f_{0, \psi} \colon S_{0, \psi} \to S_{0, \psi}$ the first return map of the flow of $X_\psi$ on $S_{0, \psi}$.
Notice that we have $S_{0, \psi} = \Pin_\psi \cup \Sigma_\psi$.
The section $\Pin_\psi$ is a union of closed surfaces transverse to $X_\psi$ and closed surfaces quasi-transverse to $X_\psi$ minus the periodic orbits $(\cO_*)_\psi$.
It follows that the section $S_{0, \psi}$ is not uniformly transverse to the vector field $X_\psi$.
Note that we have the equality $\Pin_\psi = \Pout_\psi$, but we will make the arbitrary choice to consider $\Pin_\psi$ in the following.

\begin{figure}[htb]
    \centering
    \vspace*{-2em}
    \hspace*{-1.5em}
    \includegraphics[height=0.23\textheight]{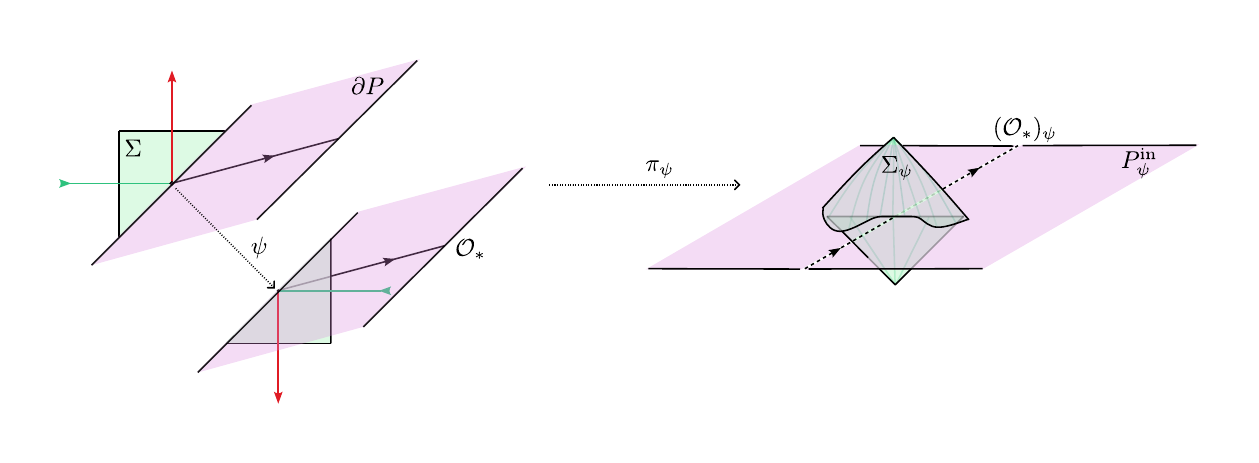}
    \vspace*{-2em}
    \caption{The section $S_{0, \psi}$ in the neighborhood of a periodic orbit of $(\cO_*)_\psi$}
    \label{fig: intermediate section in quotient}
\end{figure}

\subsubsection*{The cone fields condition}

\begin{defi}[Cone fields condition]
\label{def: cone fields condition}
Let $f \colon S \to S$ be a \diff{}.
We say that $f$ satisfies the \emph{cone fields condition} if there exists a pair of cone fields $(C^u, C^s)$ on $S$ and constants $\cst>0$, $\lambda>1$ such that for all $p \in S$:
\begin{itemize}[--]
    \item $f_* (C^u(p)) \subset \intr C^u(f(p))$, and for all $v \in C^u(p)$, $n \geq 0$,
    \vspace*{-0.2em}
    \[\Vert f^n_* v \Vert \geq \cst \lambda^n \Vert v \Vert.\]
    \item $f^\inv_* (C^s(f(p))) \subset \intr C^s(p)$,
    and for all $v \in C^s(f(p))$, $n \geq 0$,
    \vspace*{-0.2em}
    \[\Vert f^{-n}_* v \Vert \geq \cst \lambda^n \Vert v \Vert.\]
\end{itemize}
Analogously, we say that $(C^u, C^s)$ is a \emph{pair of $f$-invariant cone fields}.
We say that $C^u$ is an \emph{unstable cone field} and $C^s$ is a \emph{stable cone field} for $f$.
\end{defi}

\subsubsection*{Goal}

The goal of this section is to choose parameters \led{} such that if $\psi = \psi_\led$, the first return map $f_{0, \psi}$ has hyperbolic properties in the \nbh{} of its maximal invariant set.
This set captures all the $X_\psi$-orbit in $P_\psi$ which do not accumulate on the attractors $\pi_\psi(\gA)$ in the past or on the repellers $\pi_\psi(\gR)$ in the future (the projection is trivial in the \nbh{} of the attractors and repellers).
We will deal with those orbit separately in the next section.

We will show that an iterate (which depends on the starting point) of $f_{0, \psi}$ satisfies the cone fields condition in the \nbh{} of certain points.
However, the non-uniform transversality of the section $f_{0, \psi}$ with the vector field $X_\psi$ prevents the first return (or an iterate) from being uniformly hyperbolic.
In the following Section~\ref{sec: proof gluing thm}, we will show that for such a choice of parameters, we can truncate the section $S_{0, \psi}$ by removing a well-chosen neighborhood of the orbits $\pi_\psi (\cO_*)$ to make it compact and uniformly transverse to the \vf{} $X_\psi$ while catching the same set of orbits in uniformly bounded time, and such that the maximal invariant set is a hyperbolic set.
We have to specify some technical aspects before stating the result.

\subsubsection*{Self-intersections}

Notice that $S_{0, \psi}$ has double or triple self-intersection points (Figure~\ref{fig: intermediate section in quotient}).
Indeed the intersection $\Sigma \cap \Pin$ is a double point locus in $P$ so likewise for the image in $P_\psi$. Moreover, a point of $\pi_\psi(\Sigma)$ can be in the image of two disjoint disks of $\Sigma$ in the neighborhood of two periodic orbits $\cO_i$ and $\cO_j$ of $\cO_*$.
    The projections by $\pi_\psi$ of these disks intersect at some points and are not necessarily equal everywhere.
Then the surface $S_{0, \psi}$ has two or three tangent planes at its self-intersection points.
\begin{itemize}
    \item If $p \in \Pin_\psi \subset S_{0, \psi}$, we denote $T^\iin_p \subset T_p S_{0, \psi}$ the tangent plane induced by $T_{\tilde p} P^\iin$, where $\tilde p \in \Pin$ is the lift of $p$ in $\Pin$.
    \item If $q \in \Sigma_\psi \subset S_{0, \psi}$, we denote $T^\Sigma_q \subset T_p S_{0, \psi}$ the tangent plane(s) induced by $T_{\tilde q} \Sigma$, where $\tilde q \in \Sigma$ is a lift of $q$ in $\Sigma$.
\end{itemize}
The differential of $f_{0, \psi}$ in $p$ is well defined for a choice of tangent plane on $p$ and $f_{0, \psi}(p)$.
More precisely, for $p \in S_{0, \psi}$, once a tangent plane $T_p^*$ to $S_{0, \psi}$ has been chosen at point $p$ and $T^*_{f(p)}$, the differential $d_p f_{0, \psi} \colon T_p^* \to T_{f(p)}^*$ is well defined.

\begin{rmk}
One might want to adjust the section $\Sigma$ to avoid triple points of $S_{0, \psi}$.
However, the choice of the gluing maps $\psi = \psi_\led$ determined by the parameters $\led$ will depend on the data of $\Sigma$, and one cannot proceed in this way.
\end{rmk}

\subsubsection*{Metric on the section}
The main proposition of this section involves expansion results on non-invariant sets.
We must first choose a metric on $S_{0, \psi}$.
The metric $g$ on $P$ (see beginning of Section~\ref{sec: crossing map}) induces a metric $g^\iin$ on $\Pin$ and a metric $g^\Sigma$ on $\Sigma$.
For technical reasons, we will consider the metric $\hat g^\iin := \psi^\iin_*g^\iin$ which is the push-forward by $\psi^\iin$ of the metric $g^\iin$.
We then consider the quotient of $\hat g^\iin$ and $g^\Sigma$ onto $\pi_\psi(\Pin \cup \Sigma)$, denoted $g_\psi$.
These metrics do not coincide along self-intersections and the image $g_\psi$ does not define a continuous metric.
However, one can always lift everything in $P$ by the local homeomorphism $\pi_\psi$ and \say{pretend} that there are no intersections and all the pieces coming from $\Sigma$ and $\Pin$ are disjoint in $P_\psi$.
We don't write it but all these objects depend on the choice of the parameters $\led{}$.

\subsubsection*{Cone field}
A cone field $C$ on $S_{0, \psi}$ is the data, for each $p \in S_{0, \psi}$ and for each tangent plane $T^*$ to $S_{0, \psi}$ at the point $p$, of a cone $C(p, T^*) \subset T^* $.
\begin{itemize}
    \item When we specify that a point $p$ in $S_{0, \psi}$ belongs to $\Pin_\psi$, we consider the tangent space $T^\iin_p \subset T_p \Sigma_\psi$ induced by $T_{\tilde p}  \Pin$ in $p$, and the metric induced by $\hat g^\iin$, in other words we lift by the local homeomorphism $\pi_\psi$ in the neighborhood of the unique lift $\tilde p \in \Pin$.
    \item When we specify that a point $p \in S_{0, \psi}$ belongs to $\Sigma_\psi$, we consider the tangent space $T^\Sigma_p \subset T_p \Sigma_\psi$ induced by $T_{\tilde p} \Sigma$ in $p$, and the metric induced by $g^\Sigma$,
    in other words we lift by the local homeomorphism $\pi_\psi$ in the neighborhood of an (arbitrary) lift $\tilde p \in \Sigma$.
\end{itemize}

\subsubsection*{Statement of the main proposition}
The following proposition states that one can choose parameters \led{} and cones on $\Sigma$ and on $\Pin$ such that a (non-uniform) iterate of the first return map $f_{0, \psi} \colon S_{0, \psi} \to S_{0, \psi}$ of the flow of $X_\psi$ satisfies the cone fields condition in the neighborhood of some points.
Denote $\cW^s_\psi = \pi_\psi(\cW^s)$ the projection of the stable lamination of $\Lambda$ for $X$ in $P_\psi$, and similarly $\cW^u_\psi = \pi_\psi(\cW^u)$ for the projection of the unstable of $\Lambda$ in $P_\psi$.

\begin{prop} \label{prop: parameters and cones}
There are parameters $\lambda_0>1$, $\epsilon_0>0$, $\delta_0>0$, with the following properties.
We denote
$\psi^\iin = \psi^\iin_{\lambda_0, \epsilon_0, \delta_0}$,
$\psi^\out = \psi^\out_{\lambda_0, \epsilon_0, \delta_0}$,
$\psi = \psi^\iin \varphi \psi^\out$.
Denote $\cV_0 := \psi^\iin( \cV_{\delta_0})$.
There exists a pair $(C^u_\psi, C^s_\psi)$ of cone fields on $S_{0, \psi}$ such that 
\begin{enumerate}
    \item  \label{prop: parameters; it: directions}
    $C^u_\psi$ contains the direction tangent to $\cW^u_\psi \cap S_{0, \psi}$ in the interior and does not contain the direction tangent to $\cW^s_\psi \cap S_{0, \psi}$; 
    
 \noindent $C^s_\psi$ contains the direction tangent to $\cW^s_\psi \cap S_{0, \psi}$ in the interior does not contains the direction tangent to $\cW^u_\psi \cap S_{0, \psi}$.
\end{enumerate}
There exists integers $N_0$, $N_1$ and $N_\Sigma$ which satisfy the following properties. Let $p \in S_{0, \psi}$, and $n \geq 1$ such that $q := f_{0, \psi}^n(p)$ is well defined.

\begin{enumerate}[resume]
    \item \label{prop: parameters; it: v0}
    Suppose that $p \in\! (\Pin \cap \cV_0)_\psi$.~Then $f^k_{0, \psi}(p)\! \in\! \Sigma_\psi$ for $k= -N_0, \dots, -1, 1, \dots, N_0$.

    \item \label{prop: parameters; it: Pin minus v0 to sigma}
    Suppose $p \in (\Pin \ssm \cV_0)_\psi, q \in \Sigma_\psi,n \geq N_1$, and $f^k_{0,\psi}(p) \in \Sigma_\psi$ for $k=1, \dots, n-1$. Then
    \begin{itemize}[--]
        \item $(f_{0, \psi}^n)_* C^u_\psi(p, T^\iin_p) \subset \intr C^u_\psi(q, T^\Sigma_q)$,
        and for all $v \in C^u_\psi(p, T_p^\iin)$, $$\Vert (f_{0, \psi}^n)_* v \Vert \geq 2 \Vert v \Vert,$$
        \item $(f_{0, \psi}^n)^\inv_* C^s_\psi(q, T_q^\Sigma) \subset \intr C^s_\psi(p, T_p^\iin)$,
        and for all $v \in C^s_\psi(q, T_q^\Sigma)$, $$\Vert (f_{0, \psi}^n)^\inv_* v \Vert \geq 2 \Vert v \Vert.$$
    \end{itemize}

    \item \label{prop: parameters; it: Pin to Pin}
    Suppose $p \in \Pin_\psi$, $q \in \Pin_\psi$ and $f^k_{0, \psi}(p) \in \Sigma_\psi$ for $k=1, \dots, n-1$.
    Then
    \begin{itemize}[--]
        \item $(f_{0, \psi}^n)_* C^u_\psi(p, T_p^\iin) \subset \intr C^u_\psi(q, T_q^\iin)$,
        \item $(f_{0, \psi}^n)^\inv_* C^s_\psi(q, T_q^\iin) \subset \intr C^s_\psi(p, T_p^\iin)$,
    \end{itemize}
    Moreover, if both $p,q \in (\Pin \ssm \cV_0)_\psi$, then
    \begin{itemize}[--]
        \item for all $v \in C^u_\psi(p, T_p^\iin)$, $\Vert (f_{0, \psi}^n)_* v \Vert \geq 2 \Vert v \Vert$,
        \item for all $v \in C^s_\psi(q, T_q^\iin)$, $\Vert (f_{0, \psi}^n)^\inv_* v \Vert \geq 2 \Vert v \Vert$.
    \end{itemize}
    
    \item \label{prop: parameters; it: Pin minus v0 to sigma through Pin inter v0}
    Suppose that $p \in (\Pin \ssm \cV_0)_\psi$, $q \in \Sigma_\psi$, $n \geq N_0 +1$,
    there exists $1 \leq m \leq n-N_0$ such that $f^m_{0, \psi}(p) \in (\Pin \cap \cV_0)_\psi$ and $f^k_{0, \psi}(p) \in \Sigma_\psi$ for $k =1, \dots, m-1$, \linebreak[4]$m+1, \dots, n$.
    Then
    \begin{itemize}[--]
        \item $(f_{0, \psi}^{n})_* C^u_\psi(p, T_p^\iin) \subset \intr C^u_\psi(q, T_q^\Sigma)$,
        and for all $v \in C^u_\psi(p, T_p^\iin)$,\\[-1em] \[\Vert (f_{0, \psi}^n)_* v \Vert \geq 2 \Vert v \Vert,\]
        \item $(f_{0, \psi}^{n})^\inv_* C^s_\psi(q, T_q^\Sigma) \subset \intr C^s_\psi(p, T_p^\iin)$,
        and for all $v \in C^s_\psi(q, T_q^\Sigma)$, \\[-1em]\[\Vert (f_{0, \psi}^n)^\inv_* v \Vert \geq 2 \Vert v \Vert.\]
    \end{itemize}

    \item \label{prop: parameters; it: sigma to sigma}
    Suppose $p \in \Sigma_\psi$, $q \in \Sigma_\psi$, and $f^k_{0, \psi}(p) \in \Sigma_\psi$ for $k=1, \dots, n-1$.
    Then
    \begin{itemize}[--]
    \item $(f_{0, \psi}^n)_* C^u_\psi(p, T_p^\Sigma) \subset \intr C^u_\psi(q, T_q^\Sigma)$,
    \item $(f_{0, \psi}^n)_*^\inv C^s_\psi(q, T_q^\Sigma) \subset \intr C^s_\psi(p, T_p^\Sigma)$.
    \end{itemize}
    \nn Moreover if $n \geq N_\Sigma$, then
    \begin{itemize}[--]
        \item for all $v \in C^u_\psi(p, T_p^\Sigma)$, $\Vert (f_{0, \psi}^n)_* v \Vert \geq 2 \Vert v \Vert$,
        \item for all $v \in C^s_\psi(q, T_q^\Sigma)$, $\Vert (f_{0, \psi}^n)^\inv_* v \Vert \geq 2 \Vert v \Vert$.
    \end{itemize}
    
    \item \label{prop: parameters; it: sigma to Pin}
    Suppose $p \in \Sigma_\psi$, $q \in \Pin_\psi$ and $f^k_{0, \psi}(p) \in \Sigma_\psi$ for $k=1, \dots, n-1$.
    Then
    \begin{itemize}[--]
        \item $(f_{0, \psi}^n)_* C^u_\psi(p, T_p^\Sigma) \subset \intr C^u_\psi(q, T_q^\iin)$,
        \item $(f_{0, \psi}^n)_*^\inv C^s_\psi(q, T_q^\iin) \subset \intr C^s_\psi(p, T_p^\Sigma)$.
    \end{itemize}

    \item \label{prop: parameters; it: sigma to sigma through Pin}
    Suppose that $p \in \Sigma_\psi$, $q \in \Sigma_\psi$, $n \geq N_0+1$, there exists $1 \leq m \leq n-N_0$ such that $f_{0, \psi}^m(p) \in \Pin_\psi$, and $f^k(p) \in \Sigma_\psi$ for $k=1, \dots, m-1, m+1, \dots, n-1$. 
    Then
    \begin{itemize}[--]
        \item $(f_{0, \psi}^n)_* C^u_\psi(p, T_p^\Sigma) \subset \intr C^u_\psi(q, T_q^\Sigma)$,
        and for all $v \in C^u_\psi(p, T_p^\Sigma)$, \\[-1em]$$\Vert (f_{0, \psi}^n)_* v \Vert \geq 2 \Vert v \Vert,$$
        \item $(f_{0, \psi}^n)^\inv_* C^s_\psi(q, T_q^\Sigma) \subset \intr C^s_\psi(p, T_p^\Sigma)$, and
         for all $v \in C^s_\psi(q, T_q^\Sigma)$, \\[-1em]$$\Vert (f_{0, \psi}^n)^\inv_* v \Vert \geq 2 \Vert v \Vert.$$
    \end{itemize}
    
    \item \label{prop: parameters; it: sigma to Pin minus v0}
    Suppose $p \in \Sigma_\psi$, $q \in (\Pin \ssm \cV_0)_\psi$, and there exist $m \geq N_1$, and integers \linebreak[4]$1 \leq k_1 < \dots <k_m = n$ such that
    and $f^{k_1}(p), \dots, f^{k_m}(p) \in (\Pin \ssm \cV_0)_\psi$
    and $f^k(p) \in \Sigma_\psi$ for $k \neq k_1, \dots, k_m$, $1 \leq k \leq n-1$.
    Then
    \begin{itemize}[--]
        \item $(f_{0, \psi}^n)_* C^u_\psi(p, T_p^\Sigma) \subset \intr C^u_\psi(q, T_q^\iin)$, and
    for all $v \in C^u_\psi(p, T_p^\Sigma)$, $$\Vert (f_{0, \psi}^n)_* v \Vert \geq 2 \Vert v \Vert.$$
        \item $(f_{0, \psi}^n)^\inv_* C^s_\psi(q, T_q^\iin) \subset \intr C^s_\psi(p, T_p^\Sigma)$, and
    for all $v \in C^s_\psi(q, T_q^\iin)$, $$\Vert (f_{0, \psi}^n)^\inv_* v \Vert \geq 2 \Vert v \Vert.$$
    \end{itemize}
\end{enumerate}
\end{prop}

We refer to Figure~\ref{fig: diagram orbit} and the associated sub-figures for a diagram of the orbit corresponding to each of the items in the proposition.

\begin{figure}[hp]
    \centering
    \begin{subfigure}{1\textwidth}
         \centering
         \includegraphics[width=0.95\textwidth]{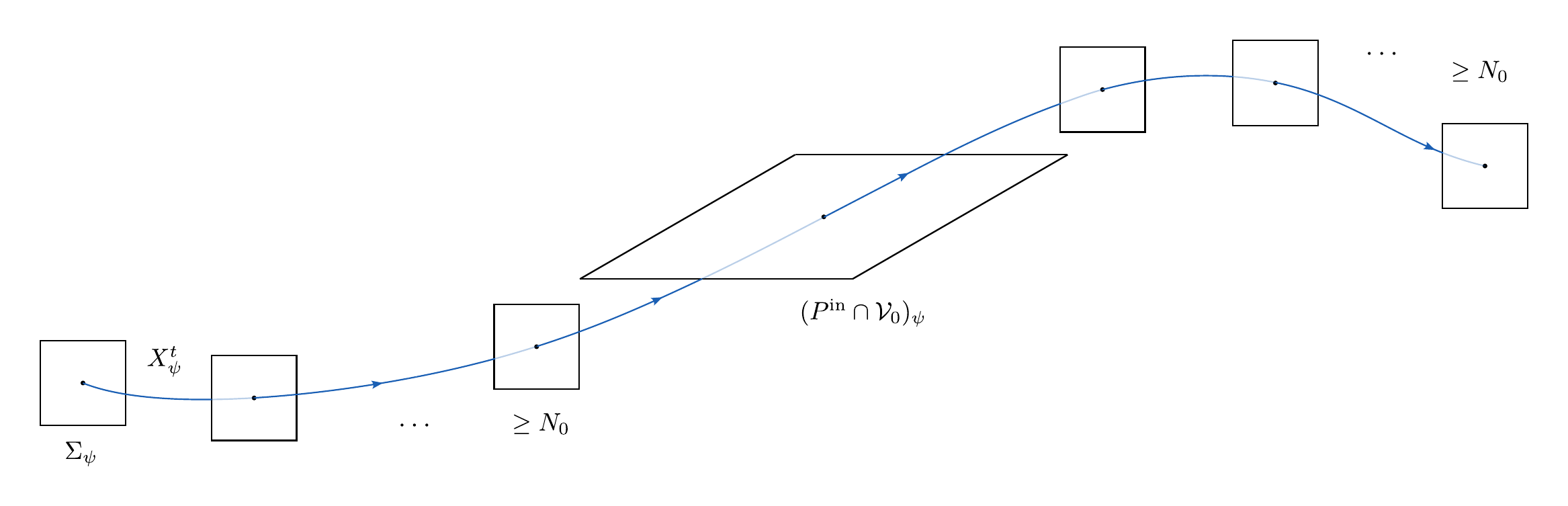}
         \vspace*{-1.5em}
         \caption{Item~\ref{prop: parameters; it: v0}}
         \label{fig: diagram orbite; subfig: v0}
    \end{subfigure}
    \vspace*{2em}
    \begin{subfigure}{1\textwidth}
         \centering
         \includegraphics[width=0.95\textwidth]{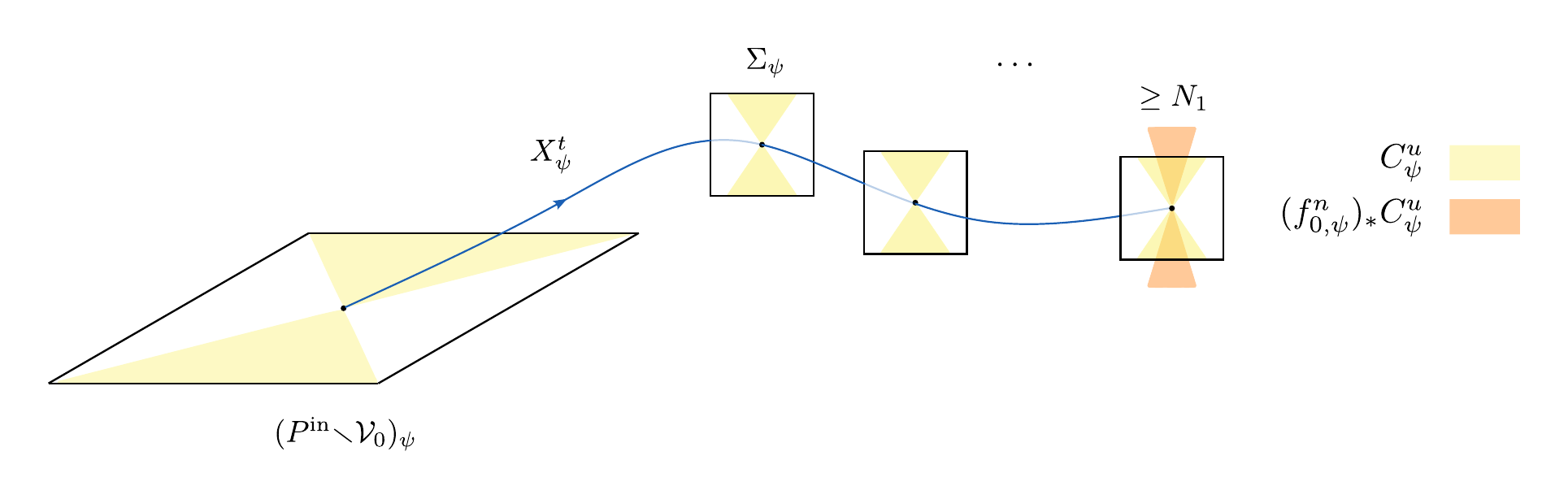}
         \vspace*{-1em}
         \caption{Item~\ref{prop: parameters; it: Pin minus v0 to sigma}}
         \label{fig: diagram orbite; subfig: Pin minus v0 to sigma}
    \end{subfigure}
    \vspace*{2em}
    \begin{subfigure}{1\textwidth}
         \centering
         \includegraphics[width=0.95\textwidth]{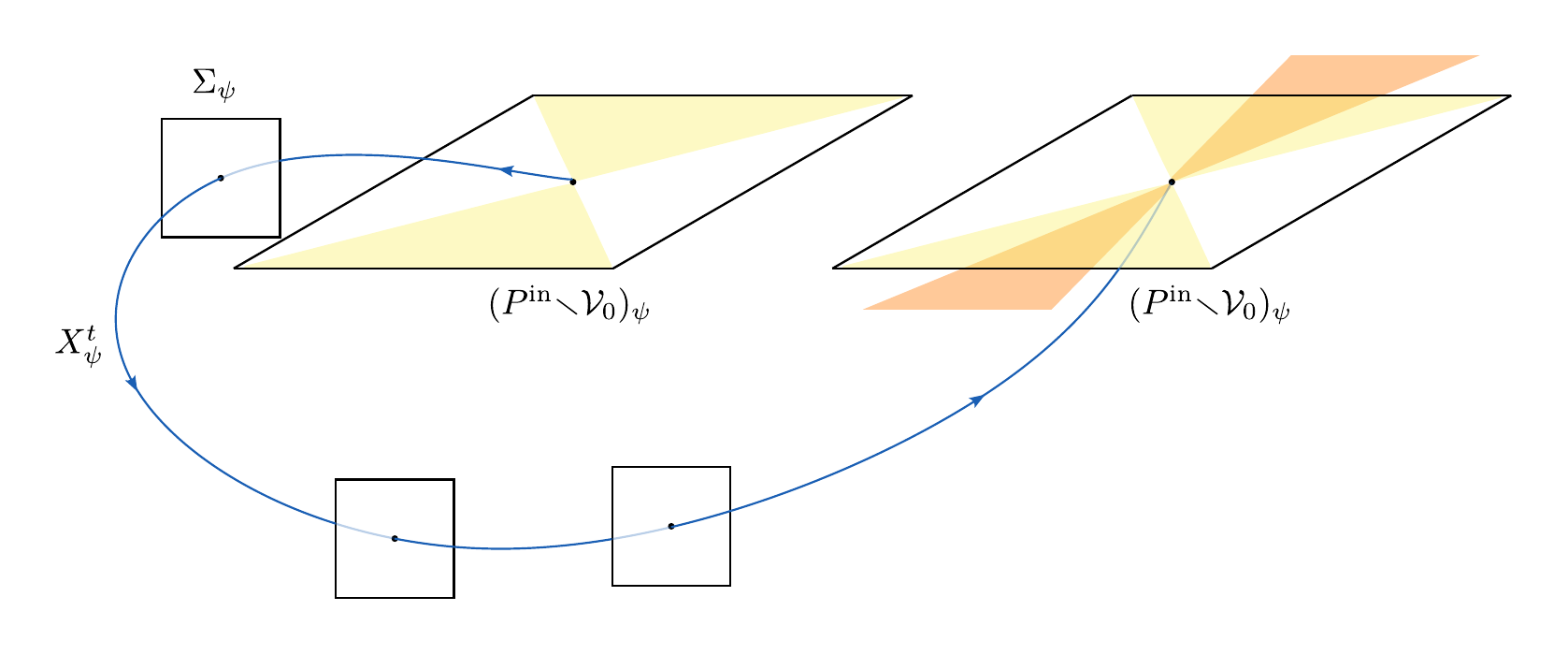}
         \vspace*{-1em}
         \caption{Item~\ref{prop: parameters; it: Pin to Pin}}
         \label{fig: diagram orbite; subfig: Pin minus v0 to Pin minus v0}
   \end{subfigure}
     \vspace*{2em}
      \begin{subfigure}{1\textwidth}
        \centering
        \includegraphics[width=0.95\textwidth]{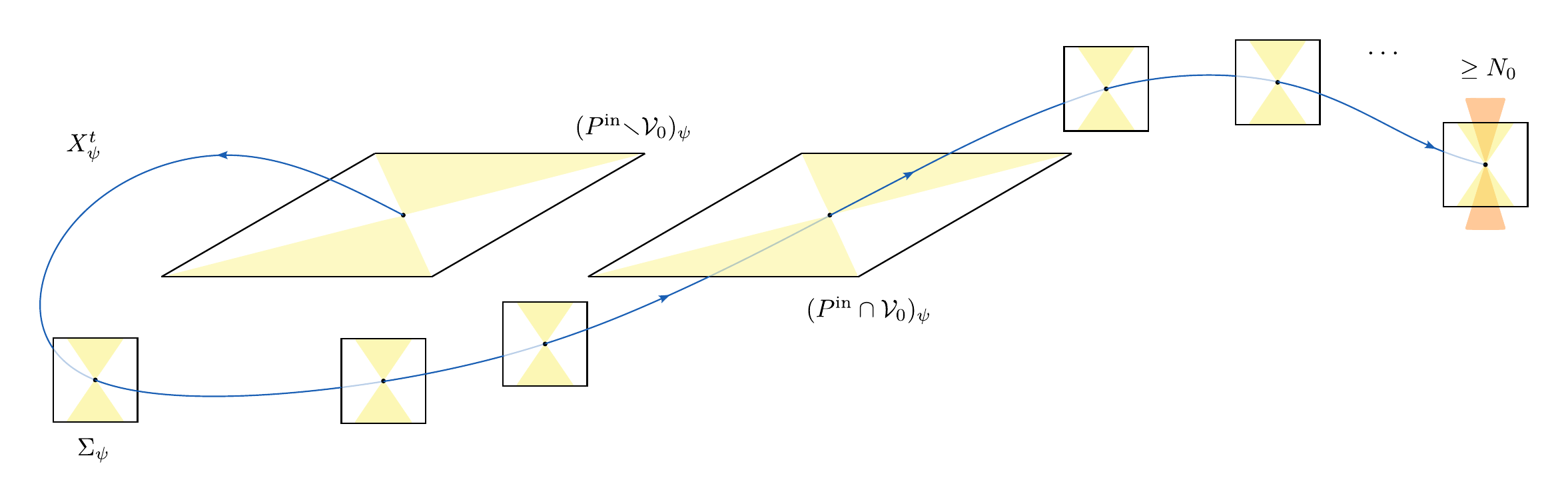}
         \vspace*{-1em}
        \caption{Item~\ref{prop: parameters; it: Pin minus v0 to sigma through Pin inter v0}}
        \label{fig: diagram orbite; subfig: Pin minus v0 to sigma through Pin inter v0}
    \end{subfigure}
    \vspace*{-3em} 
    \caption{Action of $f^n_{0,\psi}$ for each orbit of Proposition~\ref{prop: parameters and cones}}
    \label{fig: diagram orbit}
\end{figure}

\begin{figure}[bp]     
\ContinuedFloat
    \centering
    \begin{subfigure}{1\textwidth}
        \centering
        \includegraphics[width=0.95\textwidth]{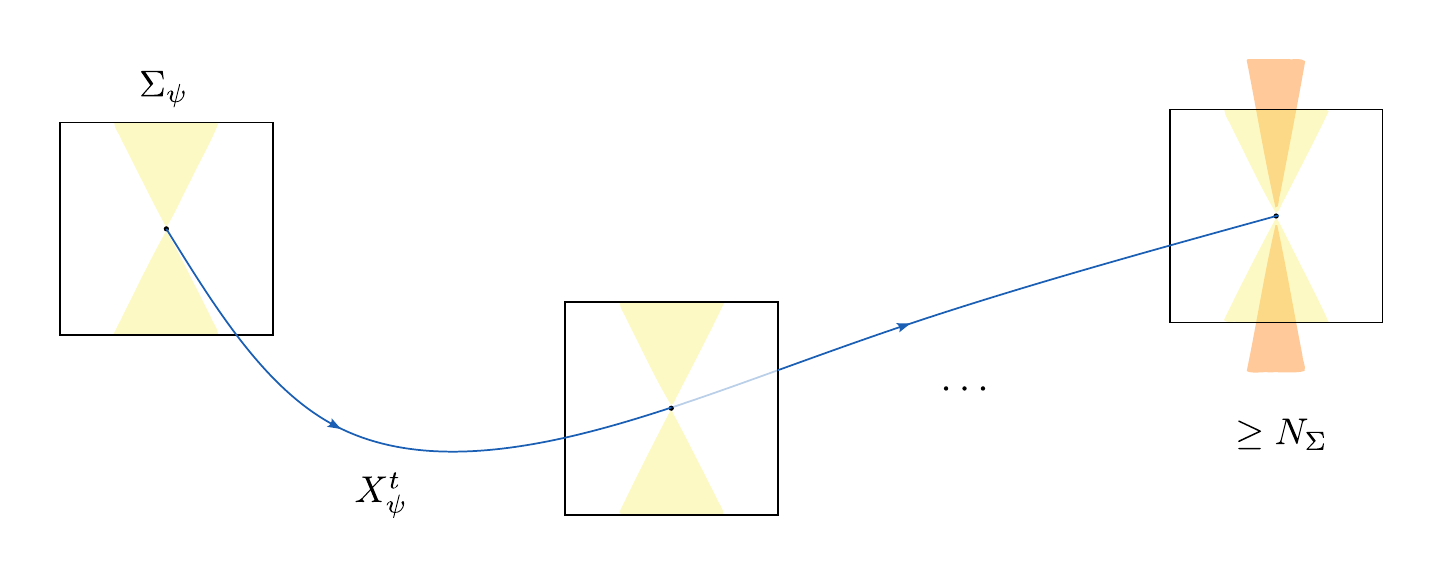}
         \vspace*{-1em}
        \caption{Item~\ref{prop: parameters; it: sigma to sigma}}
        \label{fig: diagram orbite; subfig: sigma to sigma}
    \end{subfigure} 
    \vspace*{2em}
    \begin{subfigure}{1\textwidth}
        \centering
        \includegraphics[width=0.7\textwidth]{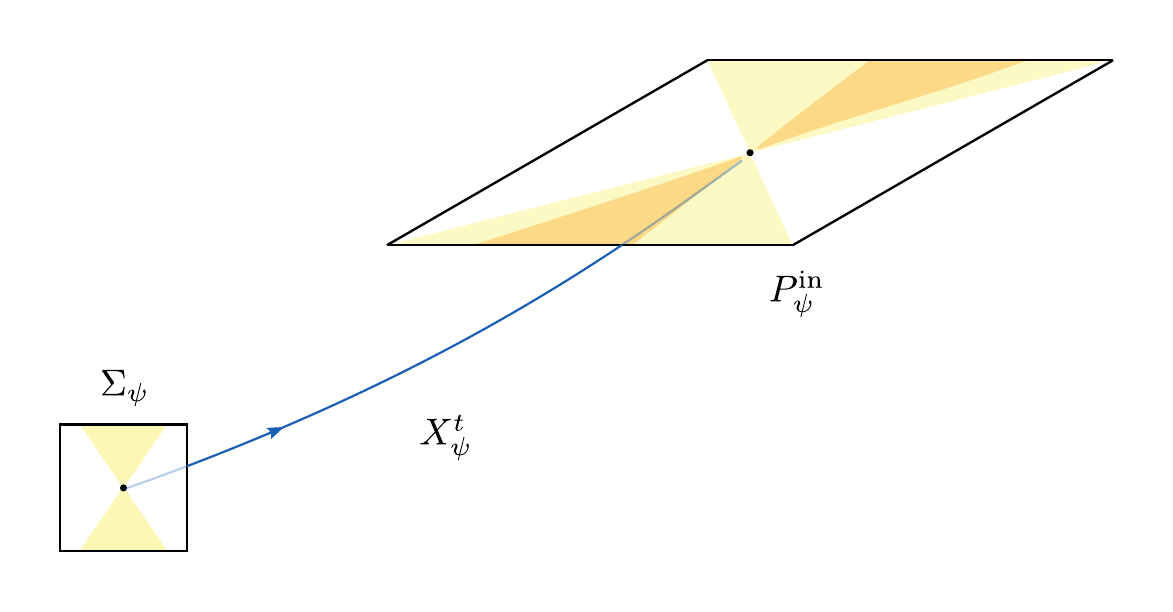}
         \vspace*{-1.5em}
        \caption{Item~\ref{prop: parameters; it: sigma to Pin}}
        \label{fig: diagram orbite; subfig: sigma to Pin}
    \end{subfigure}
    \vspace*{2em}
    \begin{subfigure}{1\textwidth}
         \centering
         \vspace*{-2em}
         \includegraphics[width=0.95\textwidth]{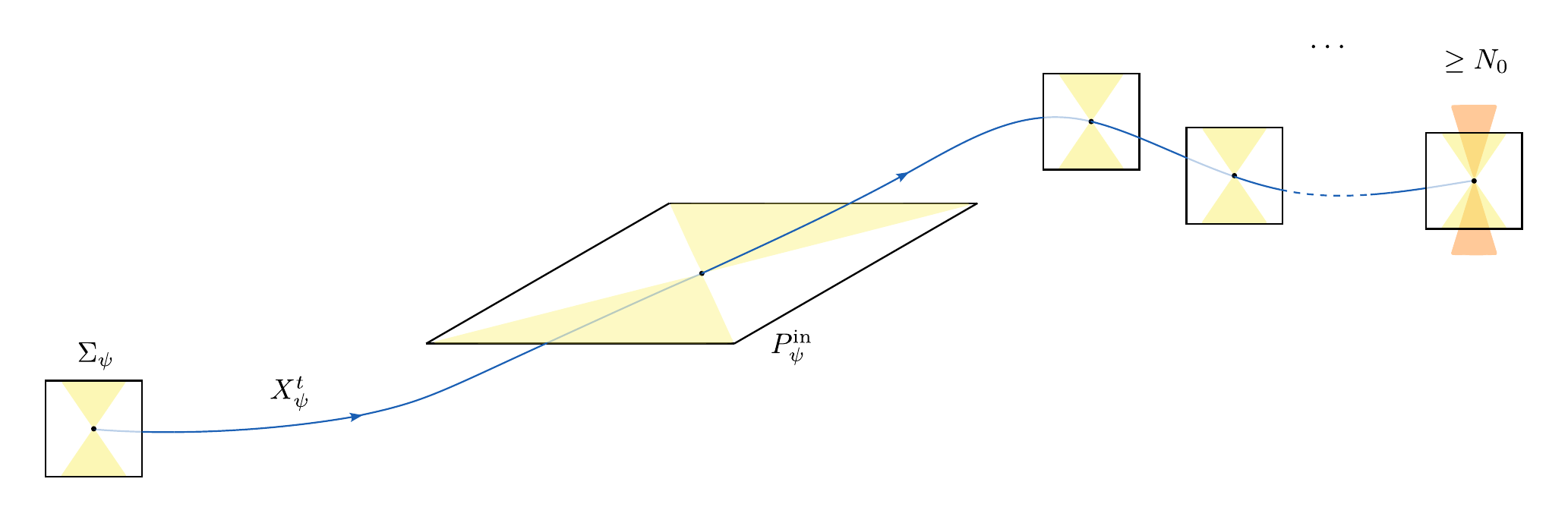}
          \vspace*{-1.5em}
         \caption{Item~\ref{prop: parameters; it: sigma to sigma through Pin}}
         \label{fig: diagram orbite; subfig: sigma to sigma through Pin}
    \end{subfigure}
    \vspace*{2em}
    \begin{subfigure}{1\textwidth}
         \centering
         \includegraphics[width=0.95\textwidth]{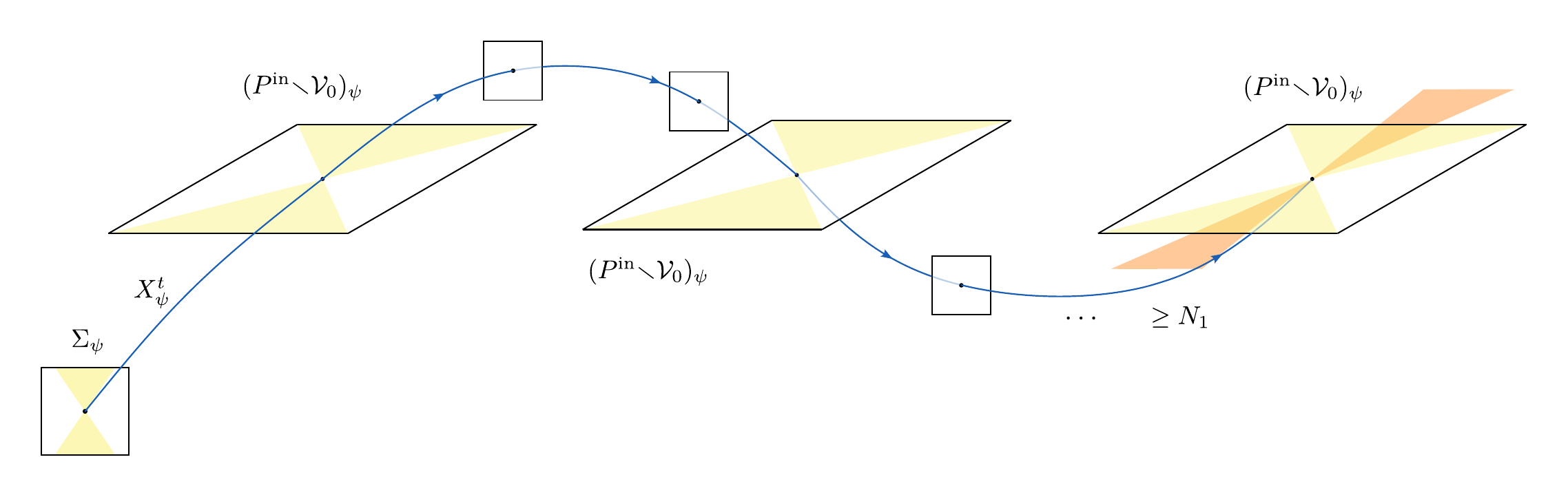}
          \vspace*{-2em}
         \caption{Item~\ref{prop: parameters; it: sigma to Pin minus v0}}
         \label{fig: diagram orbite; subfig: sigma to Pin minus v0}
    \end{subfigure}
     \vspace*{-2em} 
    \caption{Action of $f^n_{0,\psi}$ for each orbit of Proposition~\ref{prop: parameters and cones}}
\end{figure}


\subsubsection*{Lifting the return map in $P$}
For the proof, we will study lifts in $P$ of the iterates $f_{0, \psi}^n$.
If a lift $\tilde f_{0, \psi}^n$ maps a cone field $C$ on $\Sigma \cup \Pin$ on its interior in the neighborhood of a point $\tilde p$, and expands the norm of the vectors of $C$ by a factor $2$ for the metric $g$, then we can project into $P_\psi$ and we have the same result for the action $f_{0, \psi}^n$ on the cone field $C_\psi$ for the metric $g_\psi$ in $P_\psi$, in the neighborhood of $p= \pi_\psi(\tilde p)$. 
The lifts $\tilde f_{0, \psi}^n$ are compositions of the crossing and return maps of the $X$-flow and the gluing map $\psi$.

\begin{nota}\label{nota: return and crossing map in P}
Define
\begin{itemize}
    \item $f_\Sigma \colon \Sigma \to \Sigma$, the first return map of the flow of $X$ on $\Sigma$;
    \item $f_{\Sigma, \iin} \colon \Pin \to \Sigma$, the crossing map of the flow of $X$ from $\Pin$ to $\Sigma$;
    \item $f_{\out, \Sigma} \colon \Sigma \to \Pout$, the crossing map of the flow of $X$ from $\Sigma$ to $\Pout$;
    \item $\foutin \colon \Pin \to \Pout$ ,the crossing map of the flow of $X$ from $\Pin$ to $\Pout$.
\end{itemize}
\end{nota}

We have the following result.

\begin{claim} \label{claim: lift in P of return map}
\mb{}
\begin{itemize}[--]
    \item If $p \in \Sigma_\psi$ and $f_{0, \psi}(p) \in \Sigma_\psi$, then a lift of $f_{0, \psi}$ is $f_\Sigma \colon \Sigma \to \Sigma$.
    \item If $p \in \Sigma_\psi$ and $f_{0, \psi}(p) \in \Pin_\psi$, then a lift of $f_{0, \psi}$ is $\psi \, f_{\out, \Sigma}: \Sigma \to \Pin$.
    \item If $p \in \Pin_\psi$ and $f_{0, \psi}(p) \in \Sigma_\psi$, then a lift of $f_{0, \psi}$ is $f_{\Sigma, \iin}: \Pin \to \Sigma$.
    \item If $p \in \Pin_\psi$ and $f_{0, \psi}(p) \in \Pin_\psi$, then a lift of $f_{0, \psi}$ is $\psi \, \foutin: \Pin \to \Pin$.
\end{itemize}
\end{claim}

\subsubsection*{Section summary}
We will look for cones on $\Pin$ and on $\Sigma$ such that the lifts of the iterates $f^n_{0, \psi}$ satisfy the cone fields condition in the \nbh{} of some points.
The cones of Proposition~\ref{prop: parameters and cones} will be the projections in $P_\psi$ of the cones on $\Pin$ and on $\Sigma$.
The difficulty of the proof is to fix the parameters and the cones in the right order according to their relative dependence, in order to have the different compatibility conditions.
The steps are done in the following order.
Denote $\hat \cV_\delta = \psi^\iin(\cV_\delta)$.
\begin{itemize}[leftmargin=*]

    \item In Subsection~\ref{sec: parameters; subsec: lambda epsilon cone Pin},
    we show the existence of the parameters $\epsilon=\epsilon_0$, $\lambda = \lambda_0$ and the cone fields $(C^u_\iin, C^s_\iin)$ on $\Pin$ which will allow us to obtain the expansion property of Item~\ref{prop: parameters; it: Pin to Pin} of Proposition~\ref{prop: parameters and cones}, and this uniformly in $\delta$.
    Notice that we will deal separately with the invariance property of $\Pin_\psi$ and the expansion property on $(\Pin \ssm \hat \cV_\delta)_\psi$ of Item~\ref{prop: parameters; it: Pin to Pin}.
    The last one requires to choose correctly the cones on $\Pin$ and the parameters $\lambda$ and $\epsilon$ (uniformly in $\delta$), which is the work of this subsection.
    The first one requires to choose correctly $\delta$, which will be the work of Subsection~\ref{sec: parameters; subsec: delta}.
    
    We consider here the points $p \in (\Pin \ssm \hat \cV_\delta)_\psi$ whose future orbit by $X_\psi$ intersects again $(\Pin \ssm \hat \cV_\delta)_\psi$ in $q = f^n_{0,\psi}(p)$ (Figure~\ref{fig: diagram orbite; subfig: Pin minus v0 to Pin minus v0}).
    We therefore consider a well chosen restriction of the return map of the flow of $X_\psi$ on the section $(\Pin)_\psi$.
    A lift in $P$ of $f^n_{0, \psi}$ in the neighborhood of $p$ is $\psi \foutin \colon \Pin \to \Pin$, but for practical reasons we will study the composition
    \[\mb{$ f_{\iin, \psi} = \varphi \psi^\out \foutin \psi^\iin \colon \Pin \to \Pin $}\]
    which is its conjugate by $\psi^\iin$.
    We show that $f_{\iin, \psi}$ satisfies the cone fields condition for a good choice of cones $(C^u_\iin, C^s_\iin)$ on $\Pin$ and of parameters $\lambda$ and $\epsilon$ on the complementary of $\cV_\delta$, and this for all $\delta$.
    The result is summarized in Proposition~\ref{prop: lambda epsilon cones Pin}.
    The key point is that one can choose $\lambda$ large enough to have arbitrarily strong expansion and contraction of two transverse directions by $f_{\iin, \psi}$, uniformly in the parameters $\epsilon$ and $\delta$.
    
    The cone fields $(C^u_\iin, C^s_\iin)$ on $\Pin$ and the parameters $\epsilon=\epsilon_0$, $\lambda = \lambda_0$ will be fixed for the rest of the proof.
    
    \item In Subsection~\ref{sec: parameters; subsec: cone Sigma and N_Sigma}, we show the existence of cone fields $(C^u_\Sigma, C^s_\Sigma)$ on $\Sigma$ which will allow us to obtain Item~\ref{prop: parameters; it: sigma to Pin}, and which determine an integer $N_\Sigma$ which satisfies Item~\ref{prop: parameters; it: sigma to sigma} of Proposition~\ref{prop: parameters and cones}, uniformly in $\delta$.
    \begin{itemize}[leftmargin=*]
        \item Item~\ref{prop: parameters; it: sigma to Pin} concerns points $p \in \Sigma_\psi$ whose future orbit by $X_\psi$ intersects a first time $P^\iin_\psi$ in $q = f^n_{0,\psi}(p)$ (Figure~\ref{fig: diagram orbite; subfig: sigma to Pin}).
        A lift in $P$ of $f^n_{0, \psi}$ in the neighborhood of $p$ is the composition $\psi f_{\out, \Sigma} f_\Sigma^i \colon \Sigma \to \Pin$.
        We will study
        \mb{$\varphi \psi^\out f_{\out, \Sigma} f_\Sigma^i$}
        which is its conjugate by $\psi^\iin$ ($\psi^\iin$ is the identity outside $\Pin$).
        We show the existence of two cone fields $(C^u_\Sigma, C^s_\Sigma)$ on $\Sigma$ such that a compatibility criterion between the cones on $\Sigma$ and the cones on $\Pin$ is satisfied for the composition $\varphi \psi^\out f_{\out, \Sigma} f_\Sigma^i$, and this uniformly in $\delta$.
        This is Proposition~\ref{prop: cones Sigma}.
        The key point is that the projection $\pi \colon T\Sigma \to T\Pout$ parallel to $\R.X$ between $\Sigma$ and $\Pout$ (or $\Pin$) has a uniformly bounded effect on the slope of the cones (Lemma~\ref{lem: slope cone from Sigma to boundary}).
        \item Item~\ref{prop: parameters; it: sigma to sigma} concerns the points $p$ in $\Sigma_\psi$ whose future orbit by $X_\psi$ intersects $N_\Sigma$ consecutive times $\Sigma_\psi$ (Figure~\ref{fig: diagram orbite; subfig: sigma to sigma}).
        A lift in $P$ of $f^n_{0, \psi}$ in the neighborhood of $p$ is $f^n_\Sigma \colon \Sigma \to \Sigma$.
        Once the cones are fixed on $\Sigma$ containing the suitable tangent directions, the hyperbolicity of $f_\Sigma$ ensures that an iterate greater than a uniform integer $N_\Sigma$ satisfies the cone fields condition.
    \end{itemize}
    The cones $(C^u_\Sigma, C^s_\Sigma)$ on $\Sigma$ and the integer $N_\Sigma$ will be fixed for the rest of the proof.
    The pairs $(C^u_\Sigma, C^s_\Sigma)$ and $(C^u_\iin, C^s_\iin)$ are built so as to check Item~\ref{prop: parameters; it: directions} of Proposition~\ref{prop: parameters and cones}.
    
    \item In Subsection~\ref{sec: parameters; subsec: N0}, we show the existence of an integer $N_0$ which will allow us to obtain Items~\ref{prop: parameters; it: sigma to sigma through Pin} and~\ref{prop: parameters; it: Pin minus v0 to sigma through Pin inter v0} of Proposition~\ref{prop: parameters and cones}, and this uniformly in $\delta$.
    
    \begin{itemize}[leftmargin=*]
        \item Item~\ref{prop: parameters; it: sigma to sigma through Pin} concerns points $p \in \Sigma_\psi$ whose future orbit by $X_\psi$ intersects $\Pin_\psi$ once and then $N_0$ consecutive times $\Sigma_\psi$ (Figure~\ref{fig: diagram orbite; subfig: sigma to sigma through Pin}).
        A lift in $P$ of $f^n_{0, \psi}$ in the neighborhood of $p$ is the composition map 
        \mb{$f_\Sigma^j , f_{\Sigma, \iin} \, \psi \, f_{\out, \Sigma} \, f_\Sigma^i \colon \Sigma \to \Sigma$.}
        We show that the pair $(C^u_\Sigma, C^s_\Sigma)$ on $\Sigma$ fixed in the previous section form a pair of strictly invariant cone fields under this map, for an integer $j$ greater than a certain uniform integer $N_0$, uniformly in $\delta$.
        This property is summarized in Proposition~\ref{prop: N0 for sigma to sigma through Pin}. 
        The key point is that the contraction of the composition map
        $ f_{\Sigma, \iin} \, \psi \, f_{\out, \Sigma} $ 
        is bounded below uniformly in $\delta$ (Lemma~\ref{lem: contraction crossing boundary}).
        \item Item~\ref{prop: parameters; it: Pin minus v0 to sigma through Pin inter v0} is shown uniformly in $\delta$.
        It concerns points $p \in (\Pin \ssm \hat \cV_\delta)_\psi$ whose \linebreak[4]future orbit by $X_\psi$ intersects once $(\Pin \cap \hat \cV_\delta )_\psi$ and then $N_0$ consecutive times $\Sigma_\psi$ (Figure~\ref{fig: diagram orbite; subfig: Pin minus v0 to sigma through Pin inter v0}).
        A lift in $P$ of $f^n_{0, \psi}$ in the neighborhood of $p$ is the composition map
        $f_\Sigma^j \, f_{\Sigma, \iin} \, \psi \, f_{\out, \Sigma} \, f_\Sigma^i \, f_{\Sigma, \iin} \colon \Pin \ssm \cV_\delta \to \Sigma.$
        We will study the composi-\linebreak[4]tion
        $f_\Sigma^j \, f_{\Sigma, \iin} \, \psi \, f_{\out, \Sigma} \, f_\Sigma^i \, f_{\Sigma, \iin} \, \psi^\iin$
        which is its conjugate by $\psi^\iin$.
        We show that it satisfies the cone fields condition for the pair $(C^u_\iin \cup C^u_\Sigma, C^s_\iin \cup C^s_\Sigma)$ for an integer $j$ greater than a certain uniform integer $N_0$, and this uniformly in $\delta$.
        This is Proposition~\ref{prop: N0 for Pin minus V0 to Sigma through Pin inter VO}.
        The key result is that the contraction along an orbit of the flow of $X$ from $\Pin \ssm \cV_\delta$ to $\Sigma \cap \cV_\delta$ is bounded below, uniformly in $\delta$ (Lemma~\ref{lem: contraction from Pin minus V_delta to V_delta uniform in delta}).
    \end{itemize}
    The integer $N_0$ is now fixed for the rest of the proof.

    \item In Subsection~\ref{sec: parameters; subsec: delta}, we show the existence of a parameter $\delta_0$ which will allow us to obtain Item~\ref{prop: parameters; it: v0} and Item~\ref{prop: parameters; it: Pin to Pin} of Proposition~\ref{prop: parameters and cones}.

    \begin{itemize}
        \item Item~\ref{prop: parameters; it: v0} concerns the $X_\psi$ orbits crossing $(\cV_0)_\psi= (\psi^\iin (\cV_{\delta_0}))_\psi$ (Figure~\ref{fig: diagram orbite; subfig: v0}).
        We will study the composition maps
        $f_\Sigma^{n} \, f_{\Sigma, \iin} \, \psi^\iin$ and $f_\Sigma^{-n} \, f_{\out, \Sigma}^\inv \, (\varphi \psi^\out)^\inv$ from $\Pin$ to $\Sigma$.
        We show that for $\delta$ small enough, the maps restricted to $\Pin \cap \cV_\delta$ are well defined for arbitrarily big $n$ (Claim~\ref{claim: delta(n)}).        
        A key argument is that a point $p$ close enough to $\cO_*$ in $P$ has an orbit which intersects the section $\Sigma$ an arbitrarily large number of times in the past and in the future, as long as it is well defined.
        We will use a key lemma which states that the action of $\psi^\iin$ and $\psi^\out$ on the distance to the set $\cO_*$ is bounded uniformly in $\delta$ (Lemma~\ref{lem: action diffs distance orbit}).
        It then suffice to choose $\delta$ small enough depending on the integer $N_0$ fixed at the last section.
    
        \item Item~\ref{prop: parameters; it: Pin to Pin} concerns points $p \in \Pin$ whose future orbits by $X_\psi$ intersects $\Pin_\psi$ a first time at point $q = f^n_{0, \psi}(p)$.
        We know from Subsection~\ref{sec: parameters; subsec: lambda epsilon cone Pin} that if $p$ and $q$ are in $(\Pin \cap \ssm \cV_\delta)_\psi$ then the cones are invariant and expanded by $f^n_{0, \psi}$ (uniformly in $\delta$).
        It remains to check the cone invariance if $p$ or $q$ are in $(\Pin \cap \hat \cV_\delta)_\psi$.
        We will study the map $\varphi \, \psi^\out \, f_{\out, \Sigma} \, f_\Sigma^{n} \, f_{\Sigma, \iin} \, \psi^\iin \colon \Pin \to \Pin$.
        We show that for a uniformly large $n$, the cones $(C^u_\iin, C^s_\iin)$ are invariant (Lemma~\ref{lem: m0 for invariance from pin to sigma}), and this uniformly in $\delta$.
        We then use our last point which states that the map restricted to $\Pin \cap \cV_\delta$ at the start or at the destination are well defined for a given $n$ if $\delta$ is small enough.
    \end{itemize}
    The parameter $\delta_0$ is now fixed for the rest of the proof.

    \item In Subsection~\ref{sec: parameters; subsec: N1}, we show the existence of an integer $N_1$ which will allow us to have Items~\ref{prop: parameters; it: Pin minus v0 to sigma} and~\ref{prop: parameters; it: sigma to Pin minus v0} of Proposition~\ref{prop: parameters and cones}.
    
    At this point, the parameters $(\led{}) =( \lambda_0,\epsilon_0,\delta_0)$ are now fixed as well as the cone fields $(C^u_\Sigma, C^s_\Sigma)$ and $(C^u_\iin, C^s_\iin)$ on $\Sigma$ and on $\Pin$.
    Denote $\cV_0= \psi^\iin (\cV_0)$.
    \begin{itemize}[leftmargin=*]
        \item Item~\ref{prop: parameters; it: Pin minus v0 to sigma} concerns the points $p \in (\Pin \ssm \cV_0)_\psi$ whose future orbit by $X_\psi$ intersects $N_1$ consecutive times the section $\Sigma_\psi$ (Figure~\ref{fig: diagram orbite; subfig: Pin minus v0 to sigma}).
        A lift in $P$ of $f^n_{0, \psi}$ in the neighborhood of $p$ is the composition map $f^j_\Sigma , f_{\Sigma, \iin} \colon \Pin \ssm \cV_\delta \to \Sigma$.
        We will study the composition map
        \mb{$f^j_\Sigma \, f_{\Sigma, \iin} \, \psi^\iin$}
        which is its conjugate by $\psi^\iin$.
        We show that it satisfies the cone fields condition for the pair $(C^u_\iin \cup C^u_\Sigma, C^s_\iin \cup C^s_\Sigma)$ for an integer $j$ greater than some uniform integer $N_1$.
        The key points are the uniform transversality between the vector field $X$ and the surface $\Pin \ssm \cV_0$, the hyperbolicity of $f_\Sigma$, and the fact that the cone fields $(C^u_\iin, C^s_\iin)$ and $(C^u_\iin \cup C^s_\iin)$ contain \say{compatible} directions.
        \item Item~\ref{prop: parameters; it: sigma to Pin minus v0} concerns the points $p \in \Sigma_\psi$ whose future orbit by $X_\psi$ intersects $(\Pin \ssm \cV_0)_\psi$ a number $N_1$ of times (Figure~\ref{fig: diagram orbite; subfig: sigma to Pin minus v0}).
        A lift in $P$ of $f^n_{0, \psi}$ in the neighborhood of $p$ is the composition maps $(\foutin \psi)^j f_{\out, \Sigma} f_\Sigma^i \colon \Sigma \to \Pin{}$.
        We will study the conjugate by $\psi^\iin$.
        We show that it satisfies the cone fields condition for the pair $(C^u_\iin \cup C^u_\Sigma, C^s_\iin \cup C^s_\Sigma)$, for an integer $j$ greater than some uniform integer $N_1$.
        The key points are the hyperbolicity of $f_{\iin, \psi}$, and the bounded effect of the composition map $\varphi \psi^\out f_{\out, \Sigma} f_\Sigma^i$ on the slope of the cones and on the norm of the vectors.
    \end{itemize}
    The integer $N_1$ will be fixed for the rest of the proof.
    
    \item In Subsection~\ref{sec: parameters; subsec: proof}, we show each item of Proposition~\ref{prop: parameters and cones} using the results of the previous subsections.
    Each previous result is true for the parameter choice $\lambda=\lambda_0$, $\epsilon = \epsilon_0$ and $\delta= \delta_0$.
    We have to restore the conjugation defect by $\psi^\iin$ and then project into $P_\psi$.
    In other words, we push by $\psi^\iin$ the cones and the metric on $\Pin$ and project on $P_\psi$ these cones and this metric, as well as the cones and the metric on $\Sigma$.
    All the results follow afterwards.
\end{itemize}

\subsection{Parameters \texorpdfstring{$\lambda_0$}{¼ª\_0}, \texorpdfstring{$\epsilon_0$}{¼µ\_0} and cone fields on \texorpdfstring{$\Pin$}{P\^{in}}}
\label{sec: parameters; subsec: lambda epsilon cone Pin}

In this subsection
we show the existence of the parameters $\lambda = \lambda_0$, $\epsilon = \epsilon_0$ and the cone fields on $\Pin$ which will allow us to obtain the expansion property of Item~\ref{prop: parameters; it: Pin to Pin} of Proposition~\ref{prop: parameters and cones}, and this uniformly in $\delta$.
This item concerns the points $p \in (\Pin \ssm \hat \cV_\delta)_\psi$ whose future orbit by $X_\psi$ intersects again $(\Pin \ssm \hat \cV_\delta)_\psi$.
Define the following composition map:
\begin{equation} \label{eq: f_in,psi}
    f_{\iin, \psi} := \varphi \psi^\out \foutin \psi^\iin \colon \Pin \to \Pin.
\end{equation}
Let us show the following

\begin{prop} \label{prop: lambda epsilon cones Pin}
There exist parameters $\lambda_0 >1$, $\epsilon_0>0$,
a pair of continuous (closed) cone fields $(C^u_\iin,C^s_\iin)$ on $\Pin$,
there exist $K>0$, $k>0$, $k'>0$
such that for all $\delta >0$,
if $\psi^\iin = \psi_{\lambda_0, \epsilon_0, \delta}^\iin$ and $\psi^\out = \psi_{\lambda_0, \epsilon_0, \delta}^\out$ are the diffeomorphisms given by Propositions~\ref{prop: spread expansion} and~\ref{prop sym: spread expansion}, then
\begin{enumerate}
\item \label{prop: cones Pin; it: directions} \emph{(directions)}
the cone field $C^u_\iin$:
\begin{itemize}[--]
    \item is contained in a $(K, \cG^{u, \iin} /\, (\psi^\iin)^\inv_* \cG^{s, \iin})$-cone field,
    \item contains a $(k, \cG^{u, \iin} /\, (\psi^\iin)^\inv_* \cG^{s, \iin})$-cone field,
    \item contains a $(k', \ (\varphi \psi^\out)_* \cG^{u, \out} /\,\varphi_* \cG^{s, \out})$-cone field.
\end{itemize}
and the cone field $C^s_\iin$:
\begin{itemize}[--]
    \item is contained in a $(K,\ (\psi^\iin)^\inv_* \cG^{s, \iin}/\, \cG^{u, \iin})$-cone field.
    \item contains a $(k,\ (\psi^\iin)^\inv_* \cG^{s, \iin}/\, \cG^{u, \iin})$-cone field.
    \item contains a  $(k',\ \varphi_* \cG^{s, \out}/\,(\varphi \psi^\out)_* \cG^{u, \out})$-cone field.
\end{itemize}

\item \label{prop: cones Pin; it: invariance} \emph{(invariance and expansion)}
For any $p \in (\Pin \ssm \cV_\delta) \cap f_{\iin, \psi}^\inv (\Pin \ssm \cV_\delta)$,
\begin{itemize}[--]
    \item $\l(f_{\iin, \psi}\r)_*C^u_\iin(p) \subset \intr  C^u_\iin (f_{\iin, \psi}(p))$, and for all $v \in C^u_\iin(p)$, $$\Vert (f_{\iin, \psi})_* v \Vert \geq 2 \Vert v \Vert;$$
    \item$\l(f_{\iin, \psi}^\inv\r)_*C^s_\iin( f_{\iin, \psi} (p) ) \subset \intr  C^s_\iin (p)$, and for all $v \in C^s_\iin(f_{\iin, \psi} (p))$, $$\Vert (f_{\iin, \psi})^\inv_* v \Vert \geq 2 \Vert v \Vert.$$
\end{itemize}
\end{enumerate}
\end{prop}
See Figure~\ref{fig: orbit f_in,psi and cone}.

\begin{figure}[htb]
    \centering
    \vspace*{-1em}
    \includegraphics[width=\textwidth]{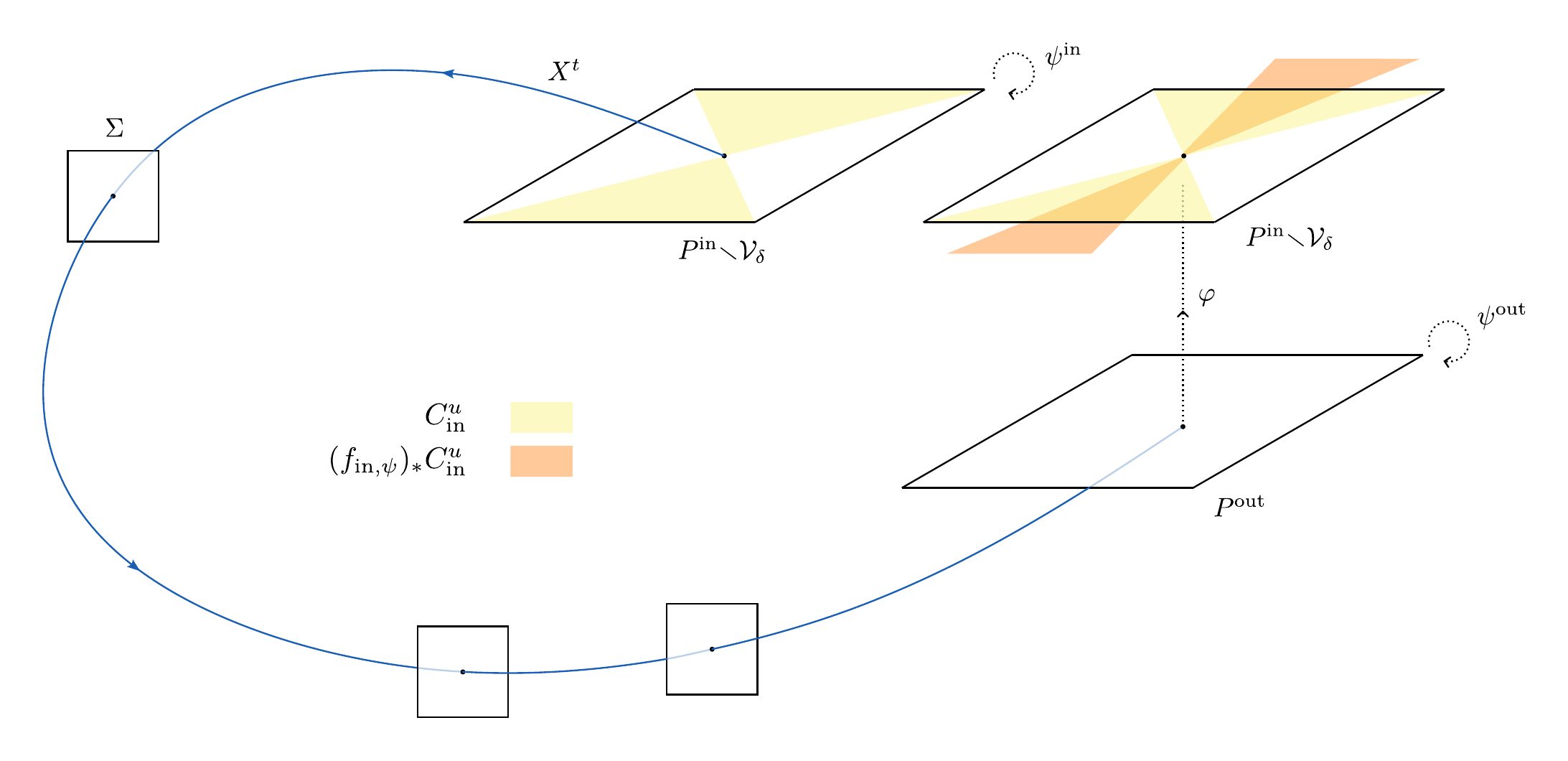}
     \vspace*{-1em}
    \caption{Decomposition of the map $f_{\iin, \psi}$ and action on the cone field $C^u_\iin$ on $(\Pin \ssm \cV_\delta)$}
    \label{fig: orbit f_in,psi and cone}
\end{figure}

We will only show the existence of the unstable cone field $C^u_\iin$. 
The proof to find a stable cone field is symmetric up to reversing the flow.
We do not require that these cones are disjoint, or complementary, so the two proofs are independent.
To show the existence of an unstable cone field for $f_{\iin, \psi}$ it is enough to find two transverse directions (not necessarily invariant), one uniformly expanded and the other uniformly contracted by $f_{\iin, \psi}$.
A cone field $C$ which contains the expanded direction and its image by $f_{\iin, \psi}$ in its interior, and whose adhesion does not contain the contracted direction is a good candidate.
A strong enough expansion and contraction factor ensures that such a cone field will be an unstable cone field, in other words strictly invariant and expanded by $f_{\iin, \psi}$.

The following lemma gives us the two transverse directions which are respectively expanded and contracted by a factor proportional to $\lambda$ (respectively $\lambda^\inv)$, on the complementary of a $\delta$-neighborhood of $\cO_*$ and its reciprocal image, and this uniformly in \led.

\begin{lem} \label{lem: directions expanded and contracted by f_in,psi}
There exists a constant $\cst{}>0$ such that for any \led{},
\begin{enumerate}
    \item for all $p \in \Pin \ssm \cV_\delta$, for all $v \in T_p\cG^{u, \iin}$,
    $$\Vert (f_{\iin, \psi})_* v \Vert > \cst \lambda \Vert v \Vert ;$$
    \label{prop: direction f_in,psi; it: expansion}
    \vspace*{-0.5em}
    \item for all $p \in f_{\iin, \psi}^\inv(\Pin \ssm \cV_\delta)$, for all $v \in T_p (\psi^\iin)^\inv_* \cG^{s, \iin}$,
    $$\Vert (f_{\iin, \psi})_* v \Vert < (\cst \lambda)^\inv \Vert v \Vert .$$
    \label{prop: direction f_in,psi; it: contraction}
\end{enumerate}
\end{lem}

\vspace*{-1em}
We deduce the following corollary.

\begin{coro} \label{coro: image of a k cone by f_in,psi}
There exists a constant $\cst>0$, such that for any \led, for any 
\mb{$p \in (\Pin \ssm \cV_\delta) \cap f_{\iin, \psi}^\inv (\Pin \ssm \cV_\delta)$}, 
if $C \subset T_p\Pin$ is a $(K,\cG^{u, \iin}/\, (\psi^\iin)^\inv_* \cG^{s, \iin})$-cone, and \mb{$q = f_{\iin, \psi}(p)$}, then
\begin{enumerate}
    \item $(f_{\iin, \psi})_* (C) \subset \intr C'$ where $C'$ is a $\l(\frac{K}{(\cst \, \lambda)^2},\   \varphi\psi^\out_* \cG^{u, \out}/\varphi_* \cG^{s, \out}\r)$-cone in\linebreak[4]$T_q\Pin$,
    \item $\Vert (f_{\iin, \psi})_* v \Vert \geq \cst \lambda K^\inv \Vert v \Vert$ for all $v \in C$.
\end{enumerate}
\end{coro}

\begin{proof}[Proof of Corollary~\ref{coro: image of a k cone by f_in,psi}]
It is enough to check that the mapped directions are indeed those announced.
Recall that $\psi^\iin$ preserves $\cG^{u, \iin}$ (Proposition \ref{prop: spread expansion}, Item~\ref{prop: spread; it: foliation preserved}), $\foutin \colon \Pin \to \Pout$ maps the pair $(\cG^{u, \iin}, \cG^{s, \iin})$ to the pair $(\cG^{u, \out}, \cG^{s, \out})$, and $\psi^\out$ preserves $\cG^{s, \out}$.
Thus we have:
\begin{align*}
    (f_{\iin, \psi})_* \cG^{u, \iin} &= (\varphi \psi^\out \foutin \psi^\iin)_* \cG^{u,\iin} \\
    &= (\varphi \psi^\out \foutin)_* \cG^{u,\iin}\\
    &= (\varphi\psi^\out)_* \cG^{u,\out}
\end{align*}
and 
\begin{align*}
    (f_{\iin, \psi})_* (\psi^\iin)^\inv_* \cG^{s, \iin} &= (\varphi \psi^\out \foutin \psi^\iin)_* \, (\psi^\iin)^\inv_* \cG^{s, \iin} \\
    &= (\varphi \psi^\out \foutin)_* \cG^{s,\iin}\\
    &= \varphi_* \cG^{s,\out}
\end{align*} 
The rest follows directly from Lemma~\ref{lem: directions expanded and contracted by f_in,psi} when standing on the intersection of $\Pin \ssm \cV_\delta$ and its reciprocal image by $f_{\iin, \psi}$.
\end{proof}

To show Lemma~\ref{lem: directions expanded and contracted by f_in,psi}, we need the following lemma, which states that the derivative of $\psi_\led^\out$ in the direction tangent to $\cG^{u, \out}$ is uniformly (in \led{}) bounded, and symmetrically, the derivative of $\psi_\led^\iin$ in the direction tangent to $\cG^{s, \iin}$ is uniformly (in \led{}) bounded.
\begin{lem}
\label{lem: bounds for derivative of diffs in perturbed direction}
There exists a constant $\cst>0$, such that for any \led{}, if $\psi^\out = \psi^\out_\led$ and $\psi^\iin = \psi^\iin_\led$, then
\begin{itemize}[--]
    \item for all $p \in \Pout$ and $v \in T_p \cG^{u, \out}$, we have 
    $\cst^\inv \Vert v \Vert \leq \Vert \psi^\out_* v \Vert \leq \cst \Vert v \Vert$;
    \item for all $p \in \Pin$ and $v \in T_p \cG^{s, \iin}$, we have 
    $\cst^\inv \Vert v \Vert \leq \Vert \psi^\iin_* v \Vert \leq \cst \Vert v \Vert$.
\end{itemize}
\end{lem}

\begin{proof}
Let us show the lemma for $\psi^\out$.
Let $v^u$ be a tangent vector to $\cG^{u, \out}$ and 
\mb{$w := \psi^\out_*(v^u) = w^s + w^u \in T\cG^{s, \out} \oplus T \cG^{u, \out}$}.
From Proposition~\ref{prop sym: spread expansion}, Item~\ref{prop sym: spread; it: foliation perturbed}, we know that $w$ is in a $(\epsilon,\cG^{u, \out}/, \cG^{s, \out})$-cone, which means that
$\Vert w^s \Vert \leq \epsilon \Vert w^u \Vert$.
From Item~\ref{prop sym: spread; it: foliation preserved}, Proposition~\ref{prop sym: spread expansion}, we know that $\psi^\out$ preserves the foliation $\cG^{s, \out}$ leaf-to-leaf, and it follows that $w^u = (H^{s, \out}_{\varsigma_{q}, \varsigma_p})_* (v^u)$
where $H^{s, \out}_{\varsigma_{q}, \varsigma_p}$ is the holonomy map of $\cG^{s,\out}$ between the leaf $\varsigma_p \in \cG^{u, \out}$ passing through $p$ and the leaf $\varsigma_{q} \in \cG^{u, \out}$ passing through $q = \psi^\out(p)$.
As the support of $\psi^\out$ is contained in strips $B_i \subset \Pout \ssm \cL^\out$ (Item~\ref{prop sym: spread; it: support}, Proposition~\ref{prop sym: spread expansion}), we have $p, q \in B_i$ and the maps $H^{s, \out}_{\varsigma_{q}, \varsigma_p}$ are holonomy maps in the strips $B_i$, we can apply Lemma~\ref{lem: bounded derivative of the stable holonomy on exit strip} which states that the ratio $\frac{\Vert w^u \Vert}{\Vert v^u \Vert}$ is uniformly bounded.
It follows that the ratio $\frac{\Vert w \Vert}{\Vert v^u \Vert}$ is bounded uniformly in the parameters \led{} for $\epsilon$ small enough.
\end{proof}

\begin{proof}[Proof of Lemma~\ref{lem: directions expanded and contracted by f_in,psi}]
Recall that $f_{\iin, \psi} = \varphi \psi^\out \foutin \psi^\iin$.
From Proposition~\ref{prop: spread expansion}, Item~\ref{prop: spread; it: composition expansion} and~\ref{prop: spread; it: foliation preserved}, the composition $\foutin \psi^\iin$ expands by a factor $\lambda$ the vectors tangent to $\cG^{u, \iin}$ on $\Pin \ssm \cV_\delta$ and maps the direction tangent to $\cG^{u, \iin}$ to the direction tangent to $\cG^{u, \out}$.
According to the previous Lemma~\ref{lem: bounds for derivative of diffs in perturbed direction}, the derivative of $\psi^\out$ in the direction $\cG^{u, \out}$ is bounded by a constant uniform in \led{}.
Moreover, $\varphi$ is a \diff{} with compact support, its differential is therefore uniformly bounded.
We deduce that there exists a constant $\cst>0$ uniform in \led{} which satisfies Item~\ref{prop: direction f_in,psi; it: expansion} of Lemma~\ref{lem: directions expanded and contracted by f_in,psi}.

Let us show Item~\ref{prop: direction f_in,psi; it: contraction}.
Let $p \in f_{\iin, \psi}^\inv (P^\iin \ssm \cV_\delta)$ and $v \in T_p (\psi^\iin)^\inv \cG^{s, \iin}$.
Then $w = (f_{\iin, \psi})_* v$ is a vector tangent to $T_q \varphi_* \cG^{s, \out}$ at the point $q = f_{\iin, \psi}(p) \in \Pin \ssm \cV_\delta$.
This item is equivalent to showing that $(f_{\iin, \psi}^\inv)_*$ expands the vector $w$ by a factor uniformly proportional to $\lambda$.
We have $f_{\iin, \psi}^\inv = (\psi^\iin)^\inv \psi^\out (\psi^\out \foutin)^\inv \varphi^\inv$. 
We study the effect of each map on the vector $w \in T_q \varphi_* \cG^{s, \out}$.
\begin{itemize}
    \item $\varphi^\inv$ is a compactly supported \diff{}, so its differential is uniformly bounded and maps the direction tangent to $\varphi_* \cG^{s, \out}$ to the direction tangent to $\cG^{s, \out}$.
    As $\varphi$ is a normalized gluing map, it maps $\Pin \ssm \cV_\delta$ on $\Pout \ssm \cV_\delta$ (Definition~\ref{def: normalized gluing map}, Item~\ref{def: normalized gluing; it: trivial}).
    \item According to Proposition~\ref{prop sym: spread expansion}, Item~\ref{prop sym: spread; it: composition expansion} and~\ref{prop sym: spread; it: foliation preserved}, the composition $\!(\psi^\out \foutin)^\inv\!$ expands by a factor $\lambda$ the vectors tangent to $\cG^{s, \out}$ on $\Pout \ssm \cV_\delta$ and maps the direction tangent to $\cG^{s, \out}$ to the direction tangent to $\cG^{s, \iin}$.
    \item By the previous Lemma~\ref{lem: bounds for derivative of diffs in perturbed direction}, the derivative of $(\psi^\iin)^\inv$ in the direction $\cG^{s, \iin}$ is bounded by a constant uniform in \led{},
\end{itemize}
We deduce that $(f_{\iin, \psi}^\inv)_*$ expands the vector $w$ by a factor uniformly proportional to $\lambda$, so there exists a constant $\cst>0$ uniform in \led{} which satisfies Item~\ref{prop: direction f_in,psi; it: contraction} of Lemma~\ref{lem: directions expanded and contracted by f_in,psi}.
\end{proof}

\begin{figure}[h]
    \centering
    \vspace*{-1em}
    \includegraphics[height=0.32\textheight]{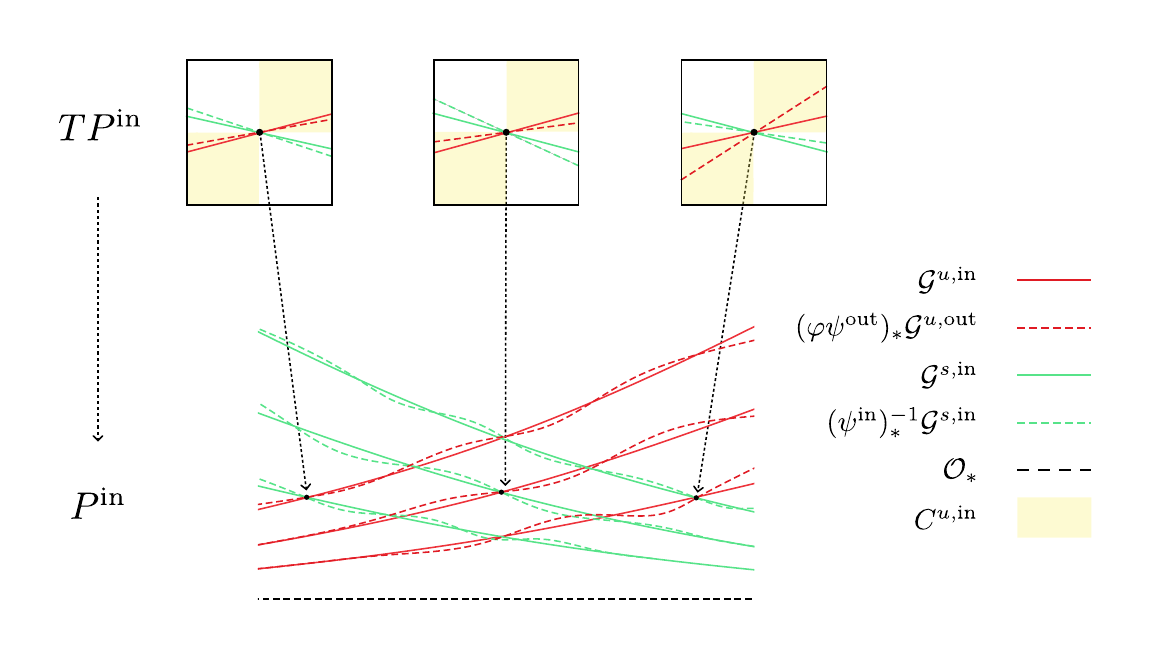}
    \vspace*{-1em}
    \caption{Cone field $C^{u, \iin}$ on $\Pin$ in the neighborhood of~$\cO_*$}
    \label{fig: C^u_in}
\end{figure}

Let us show the following lemma, which allows us to separate the suitable directions by a cone field (Figure~\ref{fig: C^u_in}).

\begin{lem} \label{lem: C^u_in and epsilon}
There exists a closed cone field $C^u_\iin$ on $\Pin$ and a parameter $\epsilon_0$ such that for any $\lambda>1$, $\delta>0$, if $\psi^{\iin} = \psi_{\lambda, \epsilon_0, \delta}^{\iin}$ and $\psi^{\out} = \psi_{\lambda, \epsilon_0, \delta}^{\out}$,
then
\begin{itemize}[--]
    \item $\intr C^u_\iin$ contains the directions tangent to ${\psi^\iin}^\inv \cG^{u, \iin}$ and $\varphi \psi^\out \cG^{u, \out}$;
    \item $C^u_\iin$ does not contain the direction tangent to $(\psi^\iin)^\inv_* \cG^{s, \iin}$;
    \item $C^u_\iin$ is the bisector $(\cG^{u, \iin}/\, \cG^{s, \iin})$-cone field on $\Pin \cap \cV_*$.
\end{itemize}
\end{lem}

Let us start by showing the following 
\begin{lem} \label{lem: C^u_in}
There exists a closed cone field $C^u_\iin$ on $\Pin$ such that
\begin{itemize}[--]
    \item $\intr (C^u_\iin)$ contains the directions tangent to $\cG^{u, \iin}$ and $\varphi_* \cG^{u, \out}$;
    \item $C^u_\iin$ does not contain the direction tangent to $\cG^{s, \iin}$;
    \item $C^u_\iin$ is the bisector $(\cG^{u, \iin}/ \, \cG^{s, \iin})$-cone field on $\Pin \cap \cV_*$.
\end{itemize}
\end{lem}

\begin{proof}[Proof of Lemma~\ref{lem: C^u_in}]
Recall that the gluing map $\varphi$ is normalized for the normalized block $(P,X)$ equipped with the \paif{} $(\cG^s, \cG^u)$, hence $\varphi_*(\cG^{u,\out})$ is transverse to $\cG^{s, \iin}$ on $\Pin$ (Item~\ref{def: normalized gluing; it: foliation transverse}, Definition~\ref{def: normalized gluing map}).
Moreover $\varphi$ is the reflection in the \ncss{} $(\cV_i, \xi_i)$, which implies that we have the relation $\varphi_* \cG^{u, \out} = \cG^{u, \iin}$ and $\varphi_* \cG^{s, \out} = \cG^{s, \iin}$ on $\Pin \cap \cV_*$.
According to this remark and by transversality of $\cG^{u, \iin}$ and $\cG^{s, \iin}$ on \Pin{}, there exists an open (thin) cone field around the direction tangent to $\cG^{s, \iin}$, which does not contain either the direction tangent to $\cG^{u, \iin}$ or the direction tangent to $\varphi_* \cG^{u, \out}$.
Moreover, we can choose this cone field as a bisector of $(\cG^{s, \iin}/\, \cG^{u, \iin})$ on $\Pin \cap \cV_*$.
The cone field $C^u_\iin$ is the complementary cone field.
\end{proof}

\begin{coro}\label{coro: bounded slope for C^u_in}
There exists $\hat K>0$, $\hat k>0$, $\hat k'>0$, such that the field $C^u_\iin$
\begin{itemize}[--]
    \item is contained in a $(\hat K,\ \cG^{u, \iin}/\, \cG^{s, \iin})$-cone field,
    \item contains a $(\hat k,\ \cG^{u, \iin}/\, \cG^{s, \iin})$-cone field,
    \item contains a $(\hat k',\ \varphi_* \cG^{u, \out}/\, \varphi_* \cG^{s, \out})$-cone field.
\end{itemize}
\end{coro}
\begin{proof}
$\Pin$ is not compact, but on $\cV_* \cap \Pin$, the cone field $C^u_\iin$ coincides with the $(1,\cG^{u, \iin}/, \cG^{s, \iin}) = (1,\ \varphi_* \cG^{u, \out}/, \varphi_* \cG^{s, \out})$-cone field.
The complementary of these neighborhoods is compact in \Pin{}, which allows us to conclude.
\end{proof}

We can now prove Lemma~\ref{lem: C^u_in and epsilon}.

\begin{proof}[Proof of Lemma~\ref{lem: C^u_in and epsilon}]
According to Propositions~\ref{prop: spread expansion} and~\ref{prop sym: spread expansion}, we have $$({\psi^\iin})_*^\inv\cG^{u, \iin} = \cG^{u, \iin},$$ and the direction tangent to $({\psi^\iin})_*^\inv\cG^{s, \iin}$ is in a $(\epsilon,\cG^{s, \iin}/\,\cG^{u, \iin})$-cone field.
Similarly $(\varphi \psi^\out)_* \cG^{s, \out} = \varphi_* \cG^{s, \out}$, and the direction tangent to $(\psi^\out)_* \cG^{u, \out}$ is in a $(\varphi_* \cG^{u, \out}, \varphi_* \cG^{s, \out})$-cone field of slope uniformly proportional to $\epsilon$ (the uniform constant comes from the action of $\varphi$ which is a compactly supported \diff{}).
Let $C^u_\iin$ be the cone field given by Lemma~\ref{lem: C^u_in}.
According to Corollary~\ref{coro: bounded slope for C^u_in},
it suffices to choose $\epsilon_0$ small before $\hat k$ and $\hat k'$.
\end{proof}

Let $\epsilon=\epsilon_0$ satisfying Lemma~\ref{lem: C^u_in and epsilon}. We have the following corollary.

\begin{coro}\label{coro: bounds for slope of perturbed direction in C^u_in}
There exists $K>0$, $k>0$, $k'>0$, such that for any $\lambda>1, \delta>0$, if $\psi^{\iin} = \psi_{\lambda, \epsilon_0, \delta}^{\iin}$ and $\psi^{\out} = \psi_{\lambda, \epsilon_0, \delta}^{\out}$, then
the cone field $C^u_\iin$
\begin{itemize}[--]
    \item is contained in a $(K,\ \cG^{u, \iin}/\, (\psi^\iin)^\inv_* \cG^{s, \iin})$-cone field,
    \item contains a $(k,\ \cG^{u, \iin}/\, (\psi^\iin)^\inv_* \cG^{s, \iin})$-cone field,
    \item contains a $(k',\ (\varphi^\iin)_* \cG^{u, \out}/\, \varphi_* \cG^{s, \out})$-cone field.
\end{itemize}
\end{coro}

\begin{proof}
This is a direct consequence of the first two items of Lemma~\ref{lem: C^u_in and epsilon}, which states that the directions tangent to $\cG^{u, \iin} = (\psi^\iin)^\inv_* \cG^{u, \iin}$ and to $(\varphi \psi^\out)_* \cG^{u, \out}$ are contained in the interior of $C^u_\iin$ and the direction tangent to $(\psi^\iin)^\inv_* \cG^{s, \iin}$ is everywhere disjoint from $C^u_\iin$.
\end{proof}

\begin{figure}[htb]
    \centering
    \vspace*{-3em}
    \includegraphics[width=\textwidth]{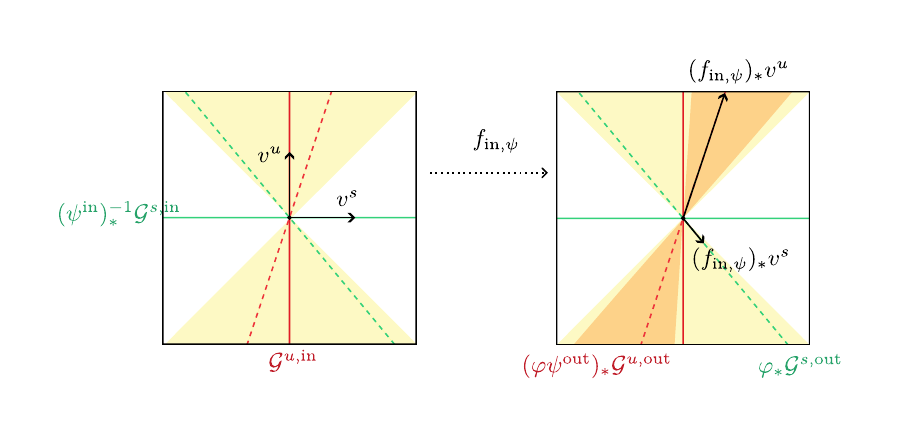}
    \vspace*{-3em}
    \caption{Directions in the cone field $C^u_\iin$ and $(f_{\iin, \psi})_* C^u_\iin$}
    \label{fig: directions in C^u_in and image}
\end{figure}

The cone field $C^u_\iin$ contains the suitable directions for $\epsilon= \epsilon_0$, and this is true for any choice of $\lambda$ and $\delta$.
It remains to show that we can choose $\lambda$ large enough so that the cone field $C^u_\iin$ is $f_{\iin, \psi}$-invariant and its vectors are expanded by $f_{\iin, \psi}$ on the well-chosen domain, and this uniformly in $\delta$.
We refer to Figure~\ref{fig: directions in C^u_in and image} for an image of the action of $f_{\iin, \psi}$ on the cone field $C^u_\iin$.

\begin{proof}[Proof of Proposition~\ref{prop: lambda epsilon cones Pin}]
We fix $C^u_\iin$ given by Lemma~\ref{lem: C^u_in and epsilon}.
Item~\ref{prop: cones Pin; it: directions} of Proposition~\ref{prop: lambda epsilon cones Pin} is satisfied by Corollary~\ref{coro: bounds for slope of perturbed direction in C^u_in} and let $ K>0$ and $ k'>0$ be the constants of the corollary.
Let $\cst>0$ be a constant which satisfies Corollary~\ref{coro: image of a k cone by f_in,psi}.
Then for $\lambda$ large enough compared to $\cst^\inv K k^\inv$, the field $C^u_\iin$ is strictly invariant by $f_{\iin, \psi}$ and its vectors are uniformly expanded on $(\Pin \ssm \cV_\delta) \cap f_{\iin, \psi}^\inv(\Pin \ssm \cV_\delta)$, and this uniformly in $\delta$.
The proof for the existence of a stable cone field $C^s_\iin$ is symmetric.
\end{proof}

\textbf{We fix the parameters $\lambda_0$, $\epsilon_0$ and a pair of cone fields $(C^u_\iin, C^s_\iin)$ satisfying Proposition~\ref{prop: lambda epsilon cones Pin}}.

\subsection{Cone fields on \texorpdfstring{$\Sigma$}{¼£} and integer \texorpdfstring{$N_\Sigma$}{N\_¼£}}

\label{sec: parameters; subsec: cone Sigma and N_Sigma}

In this section we show the existence of cone fields on $\Sigma$ which will allow us to obtain Item~\ref{prop: parameters; it: sigma to Pin}, and which determine an integer $N_\Sigma$ which satisfies Item~\ref{prop: parameters; it: sigma to sigma} of Proposition~\ref{prop: parameters and cones}.
Recall that Item~\ref{prop: parameters; it: sigma to Pin} concerns the points $p \in \Sigma_\psi$ whose future orbit by $X_\psi$ intersects a first time $P^\iin_\psi$ in $q = f^n_{0,\psi}(p)$,
and Item~\ref{prop: parameters; it: sigma to sigma} concerns the points $p \in \Sigma_\psi$ whose future orbit by $X_\psi$ intersects $n$ consecutive times $\Sigma_\psi$.
Recall that the parameters $\lambda = \lambda_0$ and $\epsilon= \epsilon_0$ were fixed in the previous Subsection~\ref{sec: parameters; subsec: lambda epsilon cone Pin}, as well as the cone fields $(C^u_\iin, C^s_\iin)$ on $\Pin$. The results will be uniform in $\delta$.
We study the maps
$$f^i_\Sigma \colon \Sigma \to \Sigma, \qq \mb{ and } \qq \varphi \psi^\out f_{\out, \Sigma} f_\Sigma^i \colon \Sigma \to \Pin.$$ 
Recall that $\Sigma$ is an affine section of $(P,X)$ (Definition~\ref{def: normalized block}, Item~\ref{def: normalized, it: affine section}).
Up to restrict $\Sigma$, there exists a coordinate system on $\Sigma$ in which the map $f_\Sigma \colon \Sigma \to \Sigma$ is an hyperbolic diagonal affine map, wherever it is well defined.
Let $T\Sigma = E^\uu \oplus E^\ss$ be the hyperbolic splitting of $f$ on $\Sigma$, and
$(\zeta^s, \zeta^u)$ the pair of stable and unstable foliation of $f_\Sigma$ on $\Sigma$.
The metric $g_\Sigma$ induced by the ambient metric $g$ on $\Sigma$ is adapted to the hyperbolic splitting of $f_\Sigma$.

\begin{nota}
Let $\mu>1$ be a uniform factor of expansion of $f_\Sigma$, in other words such that
\begin{itemize}[--]
    \item for all $v \in E^\uu$, $n \geq 0$, $\Vert (f_\Sigma^n)_* v \Vert \geq \mu^n \Vert v \Vert$,
    \item for all $v \in E^\ss$, $n \geq 0$, $\Vert (f_\Sigma^{-n})_* v \Vert \geq \mu^n \Vert v \Vert$.
\end{itemize}
\end{nota}

It follows that any closed cone field $C^u$ of constant slope whose interior $\intr C^u$ contains the direction $E^\uu$ and is disjoint from the direction $E^\ss$ is strictly invariant under the action of $f_\Sigma$ (at the first iteration) and its vectors are uniformly expanded by a (uniformly) large enough iteration $f^n_\Sigma$, which depends on the slope of the cone.
In other words, we have:

\begin{claim} \label{claim: cones return on sigma}
For any $K$, if $C^u$ is a $(K,\ \zeta^u/\,\zeta^s)$-cone field on $\Sigma$, then
\begin{itemize}[--]
    \item $(f_\Sigma)_* C^u$ is contained in a $(K \mu^{-2},\ \zeta^u/\, \zeta^s)$-cone field;
    \item there exists a constant $\cst = \cst(K) >0$ such that $\forall v \in C^u$, $\forall n \geq 0$, $$\Vert (f^n_\Sigma)_* v \Vert \geq \cst \mu^n \Vert v \Vert.$$
\end{itemize}
Similarly, if $C^s$ is a $(K,\ \zeta^s/\, \zeta^u)$-cone field on $\Sigma$, then
\begin{itemize}[--]
    \item $(f_\Sigma)^\inv_* C^s$ is contained in a $(K \mu^{-2},\ \zeta^s/\, \zeta^u)$-cone field;
    \item there exists a constant $\cst = \cst(K) >0$ such that $\forall v\in C^s$, $\forall n \geq N$, $$\Vert (f^{-n}_\Sigma)_* v \Vert \geq \cst \mu^n \Vert v \Vert.$$
\end{itemize}
\end{claim}

\begin{figure}[htb]
    \centering
    \captionsetup{width=.85\linewidth}
    \vspace*{-1em}
    \includegraphics[height=0.3\textheight]{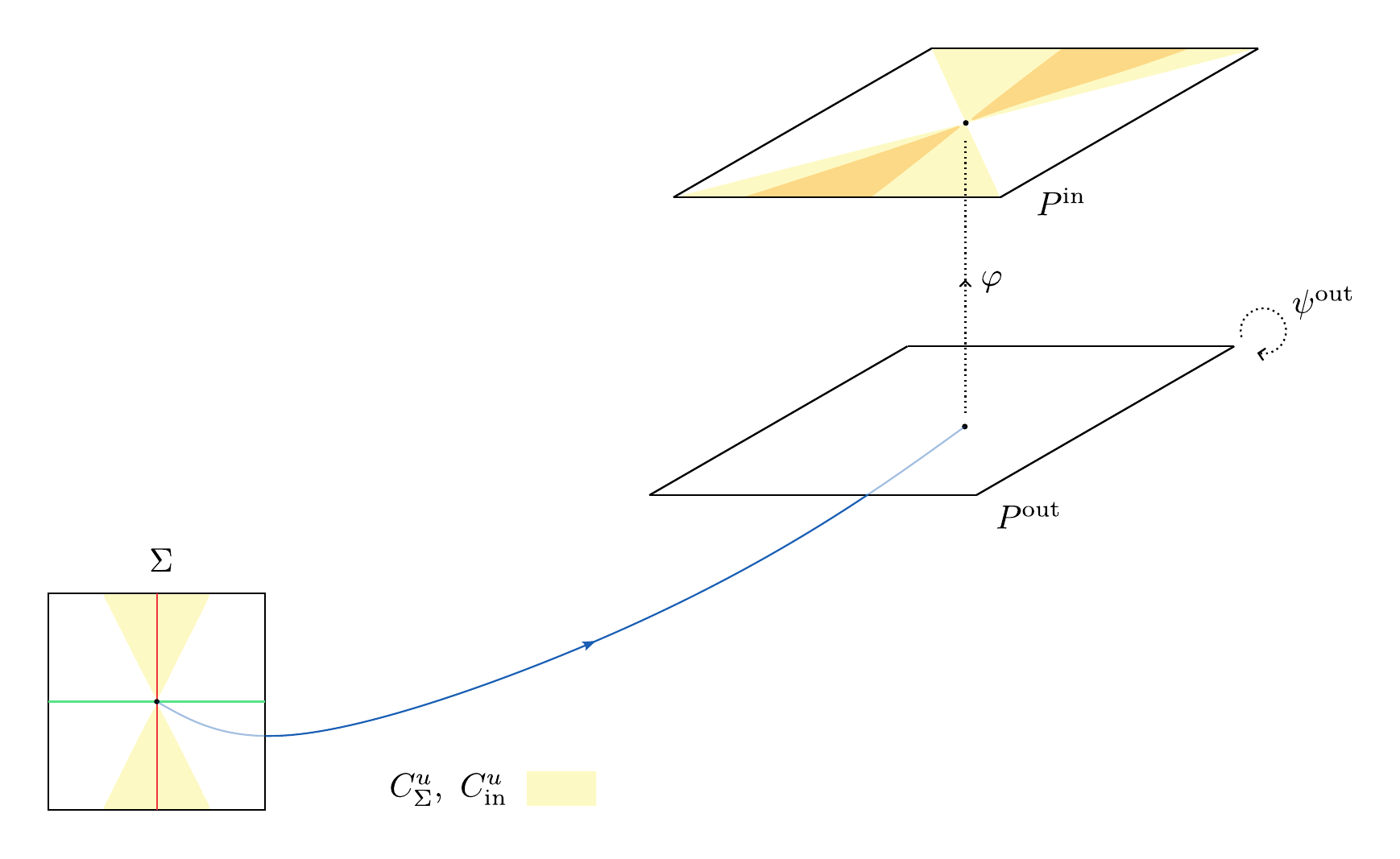}
    \vspace*{-1em}
    \caption{Trajectory of a cone $C^u_\Sigma$ under the map $\varphi, \psi^\out \, f_{\out, \Sigma} \, f_\Sigma^n$ on $\Pin$}
    \label{fig: cone Sigma to Pin compatible}
\end{figure}

We will choose the slope $K$ of a cone field $C^u_\Sigma$ of $(\zeta^u/\,\zeta^s)$-cone field such that we have a compatibility of cones on $\Sigma$ and cones on $\Pin$ for the composition map $\varphi, \psi^\out \, f_{\out, \Sigma} \, f_\Sigma^n : \Sigma \to \Pin$, for any $n \geq 0$, when it is well defined (Figure~\ref{fig: cone Sigma to Pin compatible}).
It will then suffice to choose an iterate of $f_\Sigma$ large enough for the vectors of $C^u_\Sigma$ to be uniformly expanded.
Denote $\psi^\iin_\delta = \psi^\iin_{\lambda_0, \epsilon_0, \delta}$ and $\psi^\out_\delta = \psi^\out_{\lambda_0, \epsilon_0, \delta}$.
This is summed up in the following main proposition.

\begin{prop} \label{prop: cones Sigma}
There exists a pair of cone fields $(C^u_\Sigma, C^s_\Sigma)$ on $\Sigma$
\begin{enumerate}
    \item \emph{(direction)} \label{prop: cones Sigma; it: directions} 
    There exist constants $K^u>0$ and $K^s>0$ such that
    
    $C^u_\Sigma$ is a $(K^u,\ \zeta^u/\,\zeta^s)$-cone field and $C^s_\Sigma$ is a $(K^s,\ \zeta^s/\,\zeta^u)$-cone field.

     \item \emph{(compatibility on $\Pin$)} \label{prop: cones Sigma; it: compatibilty pin} 
     For all $\delta>0$, if $\psi^\out = \psi^\out_\delta$, and for all $p \in \Sigma$ such that the point $q = \varphi \, \psi^\out \, f_{\out, \Sigma} \, f_\Sigma^n \, (p) \in \Pin$ is well defined for some integer $n \geq 0$, then
     \begin{itemize}[--]
     \item $(\varphi \, \psi^\out \, f_{\out, \Sigma} \, f_\Sigma^n \,)_* C^u_\Sigma(p) \subset \intr C^u_\iin(q)$, and
     
     \item $(\varphi \, \psi^\out \, f_{\out, \Sigma} \, f_\Sigma^n \,)^\inv_* C^s_\iin(q) \subset \intr C^s_\Sigma(p) $.
    \end{itemize}

     \item \emph{(invariance and expansion on $\Sigma$)} \label{prop: cones Sigma; it: invariance sigma}  
     There exists an integer $N_\Sigma \geq 0$ such that
      \begin{itemize}[--]
     \item $(f_\Sigma)_* C^u_\Sigma \subset \intr C^u_\Sigma, \ \mb {and} \ \forall v \in C^u_\Sigma$, $\forall n \geq N_\Sigma$, $\Vert (f_\Sigma^n)_* v \Vert \geq 2 \Vert v \Vert$,
     \item $(f_\Sigma)^\inv_* C^s_\Sigma \subset \intr C^s_\Sigma, \ \mb {and} \ \forall v \in C^s_\Sigma$, $\forall n \geq N_\Sigma$, $\Vert (f_\Sigma^n)^\inv_* v \Vert \geq 2 \Vert v \Vert$.
      \end{itemize}
\end{enumerate}
\end{prop}

\begin{rmk}
There is no requirement that the fields $C^s_\Sigma$ and $C^u_\Sigma$ be complementary to each other.
For the lemma to be true, the slope $K^s$ of $C^s_\Sigma$ will be large in front of 1, and the two cone fields can have a large intersection.
The important point of this lemma is that the slope $K^s$ and $K^u$ is independent of~$\delta$.
\end{rmk}

Recall that the map $f_{\out, \Sigma}: \Sigma \to \Pout$ is the flow of $X$ from the section $\Sigma$ to the boundary $\Pout$.
In order to show this lemma we need to show that we keep a control on the slope of the cones when we apply $f_{\out, \Sigma} \colon \Sigma \to \Pout$.
This is the following lemma, which is symmetric for $f_{\out, \Sigma}$ and $f_{\Sigma, \iin}$.

\begin{lem}
\label{lem: slope cone from Sigma to boundary}
There exists a constant $\cst>0$ such that for all $K>0$,
\begin{itemize}[--]
    \item if $C^u$ is a $(K,\ \zeta^u/\,\zeta^s)$-cone field on $\Sigma$, then $(f_{\out, \Sigma})_* C^u$ is contained in a $(\cst K,\ \cG^{u, \out}/\,\cG^{s, \out} )$-cone field on $\Pout$;
    \item if $C^s$ is a $(K,\ \cG^{s, \out}/\,\cG^{u, \out} )$-cone field on $\Pout$, then $(f_{\out, \Sigma})^\inv_* C^s$ is contained in a $(\cst K,\ \zeta^s /\, \zeta^u)$-cone field on $\Sigma$.
\end{itemize}
Likewise,
\begin{itemize}[--]
    \item if $C^u$ is a $(K,\ \cG^{u, \iin}/\,\cG^{s, \iin} )$-cone field on $\Pin$, then $(f_{\Sigma, \iin})_* C^u$ is contained in a $(\cst K,\ \zeta^u/\,\zeta^s)$-cone field on $\Sigma$;
    \item if $C^s$ is a $(K,\ \zeta^s/\,\zeta^u)$-cone field on $\Sigma$, then $(f_{\Sigma, \iin})^\inv_* C^s$ is contained in a $(\cst K,\ \cG^{s, \iin}/\,\cG^{u, \iin} )$-cone field on $\Pout$.
\end{itemize}
\end{lem}

We will just show the result for $f_{\out, \Sigma}$.
This result is not trivial because $\Pout$ is not uniformly transverse to the vector field, so the projection $T P \to T\Pout$ parallel to $\R.X$ has an unbounded effect on the norm of vectors.
The key point is that the surface $\Pout$ is invariant by a symmetry which exchanges the strong stable direction $E^\ss$ and strong unstable direction $E^\uu$ in a neighborhood of $\cO_*$.
As a result, the effect of the projection on the norm of the vectors of $E^\uu$ and $E^\ss$ is the same, which ensures that the slope of a $(E^\uu/\,E^\ss)$-cone is preserved.

\begin{proof}[Proof of Lemma~\ref{lem: slope cone from Sigma to boundary}]
Let us show the first item. The proof of the others are analogous.
This lemma concerns the segments of orbits which intersect $\Sigma$ then $\Pout$. 
They are endowed with a hyperbolic splitting $E^\ss \oplus \R.X \oplus E^\uu$ for the flow of $X$ (Remark~\ref{rmk: extend hyperbolic splitting}). 

First observe that the lemma is true for all orbits intersecting $\Pout$ outside $\cV_*$.
Indeed the sections $\Sigma$ and $\Pout \ssm \cV_*$ are uniformly transverse to the vector field $X$, the flow of $X$ is hyperbolic along these orbits, and the cones in question are uniformly far from the center-stable bundle.
It remains to check this for the orbits intersecting $\Pout$ in a \ncs{} $(\cV_i, \xi_i)$.

Let $C^u$ be a $(K,\ \zeta^u/\,\zeta^s)$-cone on $\Sigma$, and $v$ a vector tangent to $\Sigma$ at a point $p$.
Then $$v = v^\uu + v^\ss + 0 \in E^\uu \oplus E^\ss \oplus \R. X$$
with $\Vert v^\ss \Vert \leq K \Vert v^\uu \Vert$.
The transport by the (linear) flow $X^t$ during a time $t>0$ contributes to expand exponentially the component $v^\uu$ and contract exponentially the component $v^\ss$, and has no effect on the component $v^X=0$.
So the action of $X^t$ contributes to close a $(E^\uu/\,E^\ss)$-cone field exponentially fast, in other words to contract the slope of the cones exponentially fast.
    
    It is therefore sufficient to look at the effect of the projection on $T\Pout$ parallel to $X$ of a $(K,\ E^\uu/\,E^\ss)$-cone field transverse to the vector field $X$ (contained in the plane $E^\uu \oplus E^\ss$) in a \ncs{} $(\cV_i, \xi_i = (x, y, \theta))$.
    Let us note this projection $\pi \colon E^\uu \oplus E^\ss \to T \Pout$.
    To say that $\pi$ maps a $(K,\ E^\uu/\,E^\ss)$-cone onto a $(K,\ \cG^{u, \out}/\,\cG^{s, \out})$-cone is equivalent to saying that
    $\pi$ maps the decomposition $E^\uu \oplus E^\ss$ onto $T\cG^{u, \out} \oplus T\cG^{s, \out}$ and that
    $||| \res{\pi}{E^\uu} ||| = ||| \res{\pi}{E^\ss} |||$.
    The first assertion is true because by definition we have $T\cG^{u, \out} = (E^\uu \oplus \R.X) \cap T\Pout$, and $T\cG^{s, \out} = (E^\ss \oplus \R.X) \cap T\Pout$.
    The second one is shown with the following argument.
    By definition of a normalized block (Definition~\ref{def: normalized block}, Item~\ref{def: normalized, it: straight boundary}), the sets $\Sigma$, $\Pout$ and $\R.X$ are invariant in $\cV_i$ by the symmetry $s \colon (x, y, \theta) \mapsto (y, x, \theta)$, so
    we have the commutation relation $\pi \circ s_* = s_* \circ \pi$.
    We have $\Vect (\partial_x) = E^\ss$ and $\Vect (\partial_y) = E^\uu$,
    and $s_* \partial_x = \partial_y$, so $s_*$ permutes the directions $E^\uu$ and $E^\ss$.
    Finally, $s$ is an isometry for the metric $g$ which coincides with $d x^2 + d y^2 + d \theta^2$ on $\cV_i$.
    We conclude by writing that $\res{\pi}{E^\uu} = \res{\pi \circ s_*}{E^\ss} = s_* \circ \res{\pi}{E^\ss}$.
\end{proof}

To measure the action of $\psi^\out_\delta$, we need the following lemmas. 
The first lemma states that the action of $\psi^\out_\delta$ on the distance to the orbit is bounded uniformly in $\delta$.

\begin{lem}\label{lem: action diffs distance orbit}
There exists a constant $\cst>1$ such that for all $\delta>0$,
\begin{itemize}[--]
    \item $\forall \, p \in \Pout{}$, 
$\ \cst^\inv \, \dist (p, \cO_*) \leq \dist (\psi^\out_\delta(p), \cO_*) \leq \cst\, \dist (p, \cO_*)$;
\item $\forall \, p \in \Pin{}$,
$\ \cst^\inv \, \dist (p, \cO_*) \leq \dist (\psi^\iin_\delta(p), \cO_*) \leq \cst\, \dist (p, \cO_*)$.
\end{itemize}
\end{lem}

\begin{proof}
We show the lemma for $\psi^\out_\delta$.
The \diff{} $\psi^\out_\delta$ is independent of the parameter $\delta$ on the complementary of the union $\cV_*$ (Proposition~\ref{prop sym: spread expansion}, Item~\ref{prop sym: spread; it: independant delta}).
The complementary of $\cV_*$ in $\Pout$ is compact, so the property is true for any point $p \in \Pout \ssm \cV_*$.
It is sufficient to check this property in a \ncs{} $(\cV_i, \xi_i = (x, y, \theta))$.
Let $(x, \theta)$ be the coordinates induced by $\xi_i$ on $\Pout$.
A point $p=(x, \theta)$ is in the closure of a strip $B$ which is a \cc{} of $\Pout \ssm \cL^\out$.
According to Items~\ref{prop sym: spread; it: support} and~\ref{prop sym: spread; it: foliation preserved} of Proposition~\ref{prop sym: spread expansion}, the diffeomorphism $\psi^\out$ moves the point $p$ along the leaf arc of $\cG^{s, \out}$ inside $B$ and this for any parameter $\led$.
Let $\sigma^s$ be the leaf of $\cG^{s, \iin} \cap B$ passing through $p = (x,\theta)$.
Then $\sigma^s$ intersects one of the boundary leave of strip $B$ at a point $(x_1,\theta_1)$ with $x_1 > x/\sqrt{2}$ and the other boundary leave of $B$ at a point $(x_2, \theta_2)$ with $x_2 < \sqrt{2} x$ (Remark~\ref{rmk: equation affine foliation})
We conclude that the variation of the coordinate $x$ of a point $p \in B$ under the action of $\psi^\out$ is uniformly bounded.
The same is true for the distance to the orbit $\cO_i$.
\end{proof}

The second lemma states that the differentials of the maps $\psi^\out_\delta$ and $\psi^\iin_\delta$ are bounded uniformly in $\delta$.

\begin{lem}
\label{lem: differential of diffs bounded uniform in delta}
There exists $\cst>0$, such that for all $\delta>0$ and for all $v \in T\Pout{}$, 
\begin{itemize}[--]
    \item for all $v \in T\Pout$, $\cst^\inv \Vert v \Vert \leq \Vert (\psi^\out_\delta)_* v \Vert \leq \cst\,\Vert v \Vert$,
    \item for all $v \in T\Pin$, $\cst^\inv \Vert v \Vert \leq \Vert (\psi^\iin_\delta)_* v \Vert \leq \cst\,\Vert v \Vert$.
\end{itemize}
\end{lem}

\begin{proof}
We show the lemma for $\psi^\out$.
Let $(\partial_u, \partial_s)$ be a basis of unit vectors tangent to $\cG^{u, \out}$ and $\cG^{s, \out}$ respectively on \Pout{}.
The following fact is a direct consequence of Item~\ref{prop sym: spread; it: derivative bounded} of Proposition~\ref{prop sym: spread expansion} and of Lemma~\ref{lem: bounds for derivative of diffs in perturbed direction}.

\begin{claim} \label{claim: psi^out bound basis u s}
The norm of the matrix of $(\psi^\out_\delta)_*$ in the basis $(\partial_u, \partial_s)$ is uniformly (in $\delta$) bounded.
\end{claim}

This fact is not enough to bound the action of the differential because the angle between the vectors $(\partial_u, \partial_s)$ tends to 0 as we approach a periodic orbit of the boundary.
If $B$ is a field of basis on $\Pout$, and $f \colon \Pout \to \Pout$ is a \diff{}, we denote $M_B(d_pf)$ the matrix of the differential from the basis $B(p)$ on $T_p \Pout$ to the basis $B(f(p))$ on $T_{f(p)} \Pout$.
Let $M^{B'}_{B, p}$ be the change of basis matrix from basis $B$ to basis $B'$ on $T_p \Pout$.
Let $B := (\partial_u, \partial_s)$ and $B'$ be a field of orthonormal basis on $T\Pout{}$.
We want to bound uniformly in $\delta$ the operator
$\Pout{} \ni p  \mapsto M_{B'} \l( d_p \psi^\out_\delta \r) \in GL_2(\R) $.
We have the following change of basis relation in $GL_2(\R)$, where $q := \psi^\out_\delta(p)$:
    \begin{align*}
        {M_{B'}\l( d_p \psi^\out_\delta \r)}
        = {M^{B'}_{B, q}} \times {M_{B}\l( d_p \psi^\out_\delta \r)} \times {M^{B}_{B', p}}.
    \end{align*}
    As the operator
    $ p \mapsto {M_{B} (d_p\psi^\out_\delta)}$
    is bounded on $\Pout{}$ uniformly in $\delta$ (Claim \ref{claim: psi^out bound basis u s}), it is sufficient to show that
    $p \mapsto {M^{B'}_{B, q}} \times {M^{B}_{B', p}}$
    is bounded on $\Pout{}$.
    We can assume without loss of generality that the vector $\partial_u$ is the first vector of the orthonormal basis $B'$.
    Therefore the product $ {M^{B'}_{B, q}} \times {M^{B}_{B', p}} $ is a $GL_2(\R)$ matrix of determinant
    $\frac{\sin \theta_{B}(p)}{\sin \theta_{B}(q)}$,
    where $\sin \theta_{B}$ denotes the sine of the angle of the vectors of the basis $B$.
    We want to bound this ratio uniformly.
    We compare the angle
    $\measuredangle(\cG^{u, \out}, \cG^{s, \out})$ at the point $p \in \Pout$ and the point $q= \psi^\out_\delta(p)$:
    \begin{align*}
        \frac{\measuredangle(\cG^{u, \out}, \cG^{s, \out})_p}{\measuredangle(\cG^{u, \out}, \cG^{s, \out})_q} & 
        = \frac{\measuredangle(\cG^{u, \out}, \cG^{s, \out})_p}{\dist(p, \cO_*)} \times \frac{\dist(q, \cO_*)}{\measuredangle(\cG^{u, \out}, \cG^{s, \out})_q} \times  \frac{\dist(p, \cO_*)}{\dist(\psi^\out_\delta(p), \cO_*)}\,.
    \end{align*}
    According to Lemma~\ref{lem: angle foliation and vector field on boundary} and Lemma~\ref{lem: action diffs distance orbit}, each of these ratios is bounded uniformly in $\delta$.
    As explained above, Lemma~\ref{lem: differential of diffs bounded uniform in delta} follows from this.
\end{proof}

We deduce the following lemma.

\begin{lem} \label{lem: action diffs slope C^u,s_in}
There exists a constant $\hat K$, such that for any $\delta$,
\begin{itemize}[--]
    \item $(\psi^\iin_\delta)_* C^u_\iin$ is contained in a $(\hat K,\cG^{u, \iin}/, \cG^{s, \iin})$-cone field;
    \item $(\varphi \psi^\out_\delta)^\inv_* C^s_\iin$ is contained in a $(\hat K,\ \cG^{s, \out}/\, \cG^{u, \out})$-cone field.
\end{itemize}
\end{lem}

\begin{proof}
Let us show the first item.
The differential of $\psi^\iin_\delta$ is bounded uniformly in $\delta$ (Lemma~\ref{lem: differential of diffs bounded uniform in delta}) and
$C^u_\iin$ is contained in a $(K,\cG^{u, \iin}/\, (\psi^\iin)_*^\inv(\cG^{s, \iin}))$-cone field (Proposition~\ref{prop: lambda epsilon cones Pin}, Item~\ref{prop: cones Pin; it: directions}).
As $\psi^\iin_\delta$ preserves the foliation $\cG^{u, \iin}$ (Proposition~\ref{prop: spread expansion}, Item~\ref{prop: spread; it: foliation preserved}), we conclude that $\psi^\iin_\delta$ maps $C^u_\iin$ inside a $(\cG^{u, \iin}/\,\cG^{s, \iin})$-cone field of uniformly increased slope in $\delta$.

The second item is similar.
The field $C^s_\iin$ is contained in a $$(K,\ \varphi_* \cG^{s,\out}/\, (\varphi \psi^\out)_* \cG^{u, \out})\text{-cone field}$$
(Proposition~\ref{prop: lambda epsilon cones Pin}, Item~\ref{prop: cones Pin; it: directions}).
The \diff{} $\varphi^\inv$ is supported on a compact set so it maps a $(K,\ \varphi_* \cG^{s,\out}/\, (\varphi \psi^\out)_* \cG^{u, \out})$-cone field of $\Pin$ inside a $(K',\ \cG^{u, \out}/\, \psi^\out_* \cG^{u, \out})$-cone field where $K'$ is a uniform constant.
According to Lemma~\ref{lem: differential of diffs bounded uniform in delta},
the differential of $(\psi^\out_\delta)^\inv$ is bounded on $\Pout$ uniformly in $\delta$.
Since $(\psi^\out_\delta)^\inv$ preserves the foliation $\cG^{s, \out}$ (Proposition~\ref{prop sym: spread expansion}, Item~\ref{prop sym: spread; it: foliation preserved}), we deduce that it maps a $(K',\ \cG^{u, \out}/\, \psi^\out_* \cG^{u, \out})$-cone field inside a $(\cG^{s, \out}/\, \cG^{u, \out})$-cone field of slope bounded below, uniformly in $\delta$.
\end{proof}

We are now able to show Proposition~\ref{prop: cones Sigma}.

\begin{proof}[Proof of Proposition~\ref{prop: cones Sigma}]
Let us show the existence of the cone field $C^u_\Sigma$ satisfying Item~\ref{prop: cones Sigma; it: compatibilty pin}.
The effects of the differential of $\varphi \, \psi^\out \, f_{\out, \Sigma} \, f_\Sigma^n$ on the slope of a $(\zeta^u/\,\zeta^s)$-cone of $\Sigma$ are the following:
\begin{itemize}
    \item $f_\Sigma^n$ maps a $(K,\ \zeta^u/\,\zeta^s)$-cone field on a $(K,\ \ \zeta^u/\, \zeta^s)$-cone field for $n \geq 0$ (Claim~\ref{claim: cones return on sigma}),
    
    \item $f_{\out, \Sigma}$ maps a $(K,\ \zeta^u/\, \zeta^s)$-cone field on $\Sigma$ inside a $(\cG^{u, \out}/\, \cG^{s, \out})$-cone field, of slope $K_1 \leq \cst K$ on $\Pout$ (Lemma~\ref{lem: slope cone from Sigma to boundary}),
    
    \item the differential of $\varphi \psi^\out$ is bounded on $\Pout$ uniformly in $\delta$ (Lemma~\ref{lem: differential of diffs bounded uniform in delta}), so it maps a $(K_1, \cG^{u, \iin}, \cG^{s, \iin})$-cone field inside a $$(K_2,\ (\varphi \psi^\out)_* \cG^{u, \out}/\, \varphi_* \cG^{s, \out})\text{-cone field}$$ on \Pin{}, $K_2 \leq \cst K_1$ and the constant is uniform in $\delta$.
\end{itemize}

To summarize, a $(K,\ \zeta^u/\zeta^s)$-cone field on $\Sigma$ is mapped by $\varphi \, \psi^\out_\delta \, f_{\out, \Sigma} \, f_\Sigma^n$ inside a $(K',\ (\varphi \psi^\out)_* \cG^{u, \out}/\, \varphi_* \cG^{s, \out})$-cone field on $\Pout$, with $K' \leq \cst K$, where the constant is uniform in $\delta >0$.
Finally, according to Item~\ref{prop: cones Pin; it: directions} of Proposition~\ref{prop: lambda epsilon cones Pin}, it suffices to choose $C^u_\Sigma$ a $(K^u, \zeta^u/\,\zeta^s)$-cone field of small enough slope $K^u$ (in front of the constant $k'$ given by the proposition).

The existence of the $C^s_\iin$ field is shown symmetrically, and requires to choose a $(K^s,\ \zeta^s/\, \zeta^u)$-cone field of fairly large $K^s$ opening.
Indeed:
\begin{itemize}
    \item $(\varphi \psi^\out)^\inv$ maps $C^s_\iin$ inside a $(\hat K,\cG^{s, \out}/\, \cG^{u, \out})$-cone field (Lemma~\ref{lem: action diffs slope C^u,s_in}),
    \item $f_{\out, \Sigma}^\inv$ maps a $(\hat K,\cG^{s, \out}/\,\cG^{u, \out})$-cone field on $\Pout$ inside a $(K_1,\ \zeta^u/\, \zeta^s)$-cone field on $\Sigma$, $K_1 \leq \cst \hat K$ (Lemma~\ref{lem: slope cone from Sigma to boundary}),
    \item $f_\Sigma^{-n}$ maps a $(K_1,\zeta^s/\,\zeta^u)$-cone field on a $(K_1,\zeta^s/\,\zeta^u)$-cone field for $n \geq 0$ (Claim~\ref{claim: cones return on sigma})
\end{itemize}

To summarize, the field $C^s_\iin$ is mapped by $(\varphi, \psi^\out_\delta \, f_{\out, \Sigma} f_\Sigma^n)^\inv$ inside a \linebreak[4]$(K',\ \zeta^s/\,\zeta^u)$-cone field on $\Sigma$, where $K'$ is a constant uniform in $\delta >0$. So it suffices to choose $K^s = K'$ and $C^s_\Sigma$ is the $(K^s,\zeta^s/\,\zeta^u)$-cone field.
The cone fields $C^u_\Sigma$ and $C^s_\Sigma$ have constant slope $K^u$ and $K^s$, i.e., satisfies Item~\ref{prop: cones Sigma; it: directions}.
Item~\ref{prop: cones Sigma; it: invariance sigma} follows directly from Claim~\ref{claim: cones return on sigma}.
\end{proof}

\textbf{We fix a pair of cone fields $(C^u_\Sigma, C^s_\Sigma)$ on $\Sigma$ and integer $N_\Sigma \geq 0$ which satisfies Proposition~\ref{prop: cones Sigma}}.

\begin{figure}[htb]
    \centering
     \vspace*{-1em}
    \includegraphics[height=0.34\textheight]{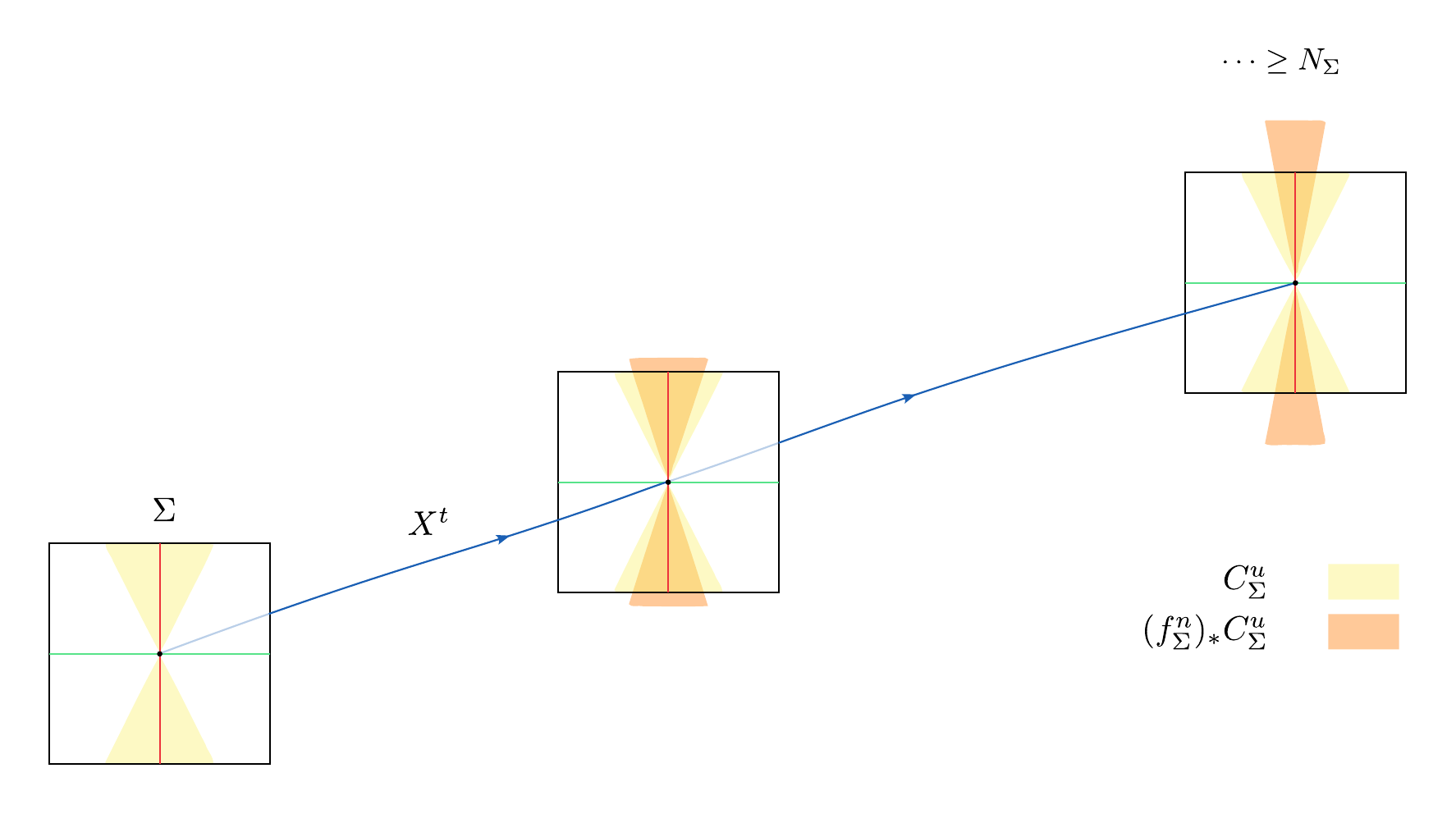}
     \vspace*{-1em}
    \caption{Trajectory of a cone $C^u_\Sigma$ under $f_\Sigma^n$}
    \label{fig: cone sigma to sigma}
\end{figure}

\subsection{Integer \texorpdfstring{$N_0$}{N\_0}}
\label{sec: parameters; subsec: N0}

In this subsection, we show the existence of an integer $N_0$ which will allow us to obtain Items~\ref{prop: parameters; it: sigma to sigma through Pin} and~\ref{prop: parameters; it: Pin minus v0 to sigma through Pin inter v0} of Proposition~\ref{prop: parameters and cones} and this uniformly in $\delta$.
We recall that Item~\ref{prop: parameters; it: sigma to sigma through Pin} concerns the points $p \in \Sigma_\psi$ whose future orbit by $X_\psi$ intersects once $\Pin_\psi$ and then $N_0$ consecutive times $\Sigma_\psi$,
and Item~\ref{prop: parameters; it: Pin minus v0 to sigma through Pin inter v0} concerns points $p \in (\Pin \ssm \hat \cV_\delta)_\psi$ whose future orbit by $X_\psi$ intersects a single time $(\Pin \cap \hat \cV_\delta )_\psi$ then $N_0$ consecutive times $\Sigma_\psi$.

\subsubsection*{From $\Sigma$ to $\Sigma$ through $\Pin$}

Let $p \in \Sigma$ and $i, j \geq 0$ such that the following map is well defined in the neighborhood of $p$:
$$f_\Sigma^j \, f_{\Sigma, \iin} \, \psi_\delta \, f_{\out, \Sigma} \, f_\Sigma^i(p) \colon \Sigma \to \Sigma.$$
We look for an integer $N_0$, such that if $j \geq N_0$, the composition map above satisfies the cone fields condition on $\Sigma$ for the pair $(C^u_\Sigma, C^s_\Sigma)$.
We want to prove the following main proposition.

\begin{prop} \label{prop: N0 for sigma to sigma through Pin}
There exists $N_0$ such that for any $\delta>0$, if
$f_\Sigma^j \, f_{\Sigma, \iin} \, \psi \, f_{\out, \Sigma} \, f_\Sigma^i$ is well defined in $p \in \Sigma$ for some $i \geq 0$ and $j \geq N_0$, then, denoting $q$ the image of~$p$,
\begin{itemize}[--]
    \item $\l(f_\Sigma^j \, f_{\Sigma, \iin} \, \psi \, f_{\out, \Sigma} \, f_\Sigma^i\r)_*$ maps $C^u_\Sigma(p)$ inside $C^u_\Sigma(q)$ and expands the norm of the vectors of $C^u_\Sigma(p)$ by a factor 2,
    \item $\l( f_\Sigma^j \, f_{\Sigma, \iin} \, \psi \, f_{\out, \Sigma} \, f_\Sigma^i \r)_*^\inv$ maps $C^s_\Sigma(q)$ inside $C^s_\Sigma(p)$ and expands the norm of the vectors of $C^s_\Sigma(q)$ by a factor 2.
\end{itemize}
\end{prop}

\begin{figure}[htb]
    \centering
    \vspace*{-2em}
    \captionsetup{width=.85\linewidth}
    \includegraphics[width=\textwidth]{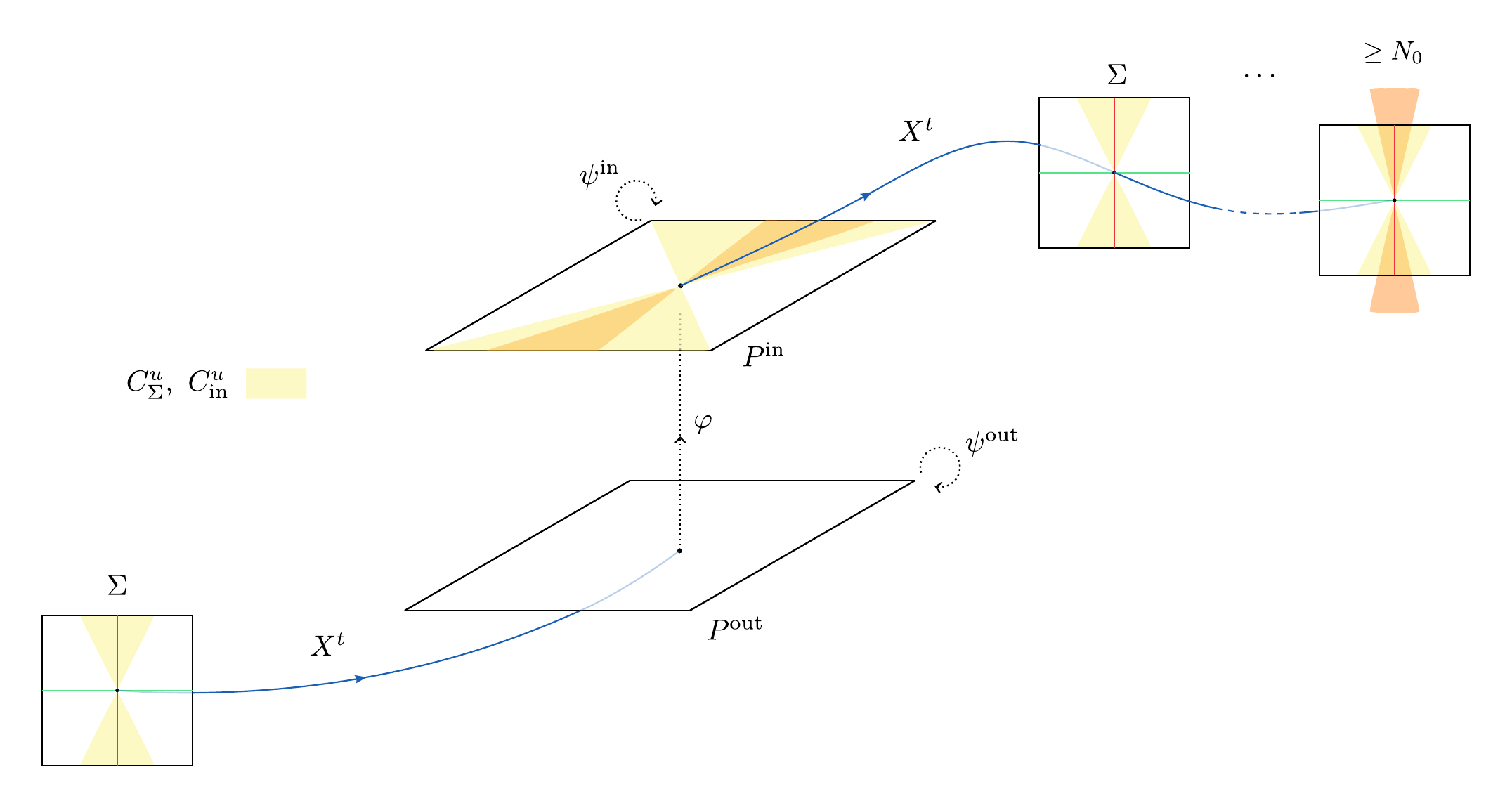}
     \vspace*{-1em}
    \caption{Trajectory of a cone $C^u_\Sigma$ under the map $f_\Sigma^{N_0} \, f_{\Sigma, \iin} \, f_{\out, \Sigma} \, f_\Sigma^i$}
    \label{fig: N0 for sigma to sigma through Pin}
\end{figure}

We know that the iterates of $f_\Sigma$ contribute to close $C^u_\Sigma$ and expand its vectors exponentially (Claim~\ref{claim: cones return on sigma}).
We will first show a result (Lemma~\ref{lem: contraction crossing boundary}) which gives a lower bound for the contraction of the composition map $$f_{\Sigma, \iin} \, \psi_\delta \, f_{\out, \Sigma} \colon \Sigma \to \Sigma$$ on the vectors of $C^u_\Sigma$ and on the slope of a $(\zeta^u/\,\zeta^s)$-cone, and this uniformly in $\delta$.
Let us start by estimating the effect of the crossing maps $f_{\out, \Sigma} \colon \Sigma \to \Pout$ and $f_{\Sigma, \iin} \colon \Pin \to \Sigma$.

\begin{lem} \label{lem: contraction sigma to boundary is distance to orbit}
There exists a constant $\cst >0$ such that:
\begin{itemize}[--]
    \item If $p \in \Sigma$ such that $q = f_{\out, \Sigma}(p) \in \Pout$ is well defined, then for any $v \in T_p \Sigma$
$$ \cst^\inv \dist(q, \cO_*)^\inv \leq \frac{\Vert (f_{\out, \Sigma})_* v\Vert}{\Vert v \Vert} \leq \cst \dist(q, \cO_*)^\inv.$$
    \item If $p \in \Pin$ such that $q = f_{\Sigma, \iin}(p) \in \Sigma$ is well defined, then for any $v \in T_p \Pin$
    $$ \cst^\inv \dist(p, \cO_*)
\leq \frac{\Vert (f_{\Sigma, \iin})_* v\Vert}{\Vert v \Vert} \leq 
\cst \dist(p, \cO_*).$$
\end{itemize}
\end{lem}

\begin{proof}
Let $p \in \Sigma$ and $v \in T_p \Sigma$. Let $\tau = \tau(p)$ be the time such that
$X^\tau(p) = f_{\out, \Sigma}(p) = q \in \Pout$.
The time $\tau$ is uniformly bounded.
Let $\pi \colon TP \to T\Pout$ be the projection on $T\Pout$ parallel to $\R.X$.
We have $\l( f_{\out, \Sigma} \r)_* v = \pi(X^\tau_*v)$.
Since the time $\tau$ is bounded, the norm $\Vert X^\tau_*v \Vert$ and the angle $\measuredangle(X^\tau_*v, X)_q$ are proportional to the norm $\Vert v \Vert$ and the angle $\measuredangle(X, v)_p$, respectively, and the proportionality constant is uniform in $p$.
From this observation and from the general Claim~\ref{claim: projection bounded by angle}, there exists a uniform constant $\cst >0$ such that
$$  
\cst^\inv \l\vert \frac{\sin \measuredangle (v, X)_p}{\sin \measuredangle (f_{\out, \Sigma})_*v , X)_q} \r\vert
\leq \frac{\Vert \l( f_{\out, \Sigma} \r)_* v \Vert}{\Vert v \Vert} \leq
\cst \l\vert \frac{\sin \measuredangle (v, X)_p}{\sin \measuredangle (f_{\out, \Sigma})_*v , X)_q} \r\vert.
$$
On the other hand, by uniform transversality we know that the angle between a vector $v \in T\Sigma$ and $X$ is uniformly bounded,
and that the angle between a vector $v \in T_q\Pout$ and $X$ is bounded below by a factor uniformly proportional to the distance $\dist(\cdot, \cO_*)$ (Lemma~\ref{lem: angle foliation and vector field on boundary}).
We deduce the existence of a constant such that
\[
\frac{\cst^\inv }{\dist(q, \cO_*)} 
\leq \frac{\Vert \l( f_{\out, \Sigma} \r)_* v \Vert}{\Vert v \Vert} \leq
\frac{\cst }{\dist(q, \cO_*)}.
\]
The same arguments applied to $f_{\Sigma, \iin}^\inv$ allow us to prove the second item.
\end{proof}

The following lemma bounds uniformly in $\delta$ the differential of the composition map $f_{\Sigma, \iin} \, \psi_\delta \, f_{\out, \Sigma}$ (when well defined) and its action on the slope of the cone field $C^u_\Sigma$ and $C^s_\Sigma$.

\begin{lem} \label{lem: contraction crossing boundary}
\mb{}
\begin{enumerate}
    \item The differentials of
$f_{\Sigma, \iin} \, \psi_\delta \, f_{\out, \Sigma} \colon \Sigma \to \Sigma$
are bounded uniformly in $\delta$;
    \item There exists a constant $K>0$ uniform in $\delta$ such that
\begin{itemize}[--]
    \item $(f_{\Sigma, \iin} \, \psi_\delta \, f_{\out, \Sigma})_* C^u_\Sigma$
    is inside a $(K, \zeta^u/\,\zeta^s)$-cone field,
    \item $(f_{\Sigma, \iin} \, \psi_\delta \, f_{\out, \Sigma})^\inv_* C^s_\Sigma$
    is inside a $(K, \zeta^s/\,\zeta^u)$-cone field.
\end{itemize}
\end{enumerate}
\end{lem}

\begin{proof}
Let $\psi = \psi_\delta$.
Let $p \in \Sigma$ such that $q = f_{\Sigma, \iin} \, \psi \, f_{\out, \Sigma} \in \Sigma$ is well defined.
Let us write
$$ 
\frac{\Vert (f_{\Sigma, \iin} \, \psi \, f_{\out, \Sigma})_* v \Vert}{\Vert v \Vert } = 
\underset{(a)}{\underbrace{\frac{\Vert (f_{\Sigma, \iin} \, \psi \, f_{\out, \Sigma})_* v \Vert}{\Vert ( \psi \, f_{\out, \Sigma})_* v \Vert }}} \ 
\underset{(b)}{\underbrace{\frac{\Vert (\psi \, f_{\out, \Sigma})_* v \Vert }{\Vert (f_{\out, \Sigma})_* v \Vert }}} \ 
\underset{(c)}{\underbrace{\frac{\Vert (f_{\out, \Sigma})_* v \Vert }{\Vert v \Vert }}}
$$ 
Let us note $p_1 = f_{\out, \Sigma}(p) \in \Pout$ and $p_2 = \psi(p_1) \in \Pin$.
There exists a constant $\cst>0$ such that
\begin{itemize}
    \item from Lemma~\ref{lem: contraction crossing boundary}: $ \cst^\inv \dist(p_2, \cO_*)
\leq (a) \leq \cst \dist(p_2, \cO_*) $;
    \item from Lemma~\ref{lem: differential of diffs bounded uniform in delta}: $\cst^\inv \leq (b) \leq \cst $;
    \item from Lemma~\ref{lem: contraction sigma to boundary is distance to orbit}: $ \cst^\inv \dist(p_1, \cO_*)^\inv
\leq (c) \leq
\cst \dist(p_1, \cO_*)^\inv $.
\end{itemize}

Up to change the constant, we have
$$ \cst^\inv \frac{\dist(p_2, \cO_*)}{\dist(p_1, \cO_*)} \leq
\frac{\Vert (f_{\Sigma, \iin} \, \psi \, f_{\out, \Sigma})_* v \Vert}{\Vert v \Vert } 
\leq \cst \frac{\dist(p_2, \cO_*)}{\dist(p_1, \cO_*)}
$$
Since $p_2 = \psi(p_1)$, according to Lemma~\ref{lem: action diffs distance orbit}, up to change the constant, we have a constant $\cst$ uniform in $\delta$ such that
$$\cst^\inv
\leq \frac{\vert (f_{\Sigma, \iin} \, \psi \, f_{\out, \Sigma})_* v \vert }{\vert v \vert } \leq
\cst. $$

Let us show the second item.
We prove only the statement concerning the image of the cone field $C^u_\Sigma$, the one for $C^s_\Sigma$ being symmetric.
\begin{itemize}
    \item We know from Proposition~\ref{prop: cones Sigma} that the composition map $$\varphi \psi^\out f_{\out, \Sigma} \colon \Sigma \to \Pin$$ maps $C^u_\Sigma$ inside $C^u_\iin$,  for any $\delta>0$.
    \item The differential of $\psi^\iin_\delta$ is bounded uniformly in $\delta$ (Lemma~\ref{lem: differential of diffs bounded uniform in delta}) and
    $C^u_\iin$ is contained in a $(K,\cG^{u, \iin}/,(\psi^\iin)_*^\inv(\cG^{s, \iin}))$-cone field (Proposition~\ref{prop: lambda epsilon cones Pin}, Item~\ref{prop: cones Pin; it: directions}).
    Since $\psi^\iin_\delta$ preserves $\cG^{u, \iin}$, we have that $\psi^\iin_\delta$ maps $C^u_\iin$ inside a $(\cG^{u, \iin}/\,\cG^{s, \iin})$-cone field with a slope bounded uniformly in $\delta$.
    \item Using Lemma~\ref{lem: slope cone from Sigma to boundary}, we deduce that the crossing map $f_{\Sigma, \iin}$ maps the cone field $\psi^\iin_* C^u_\iin$ inside a $(\zeta^u/\,\zeta^s)$-cone field on $\Sigma$ with a slope bounded uniformly in $\delta$.
\qedhere\end{itemize}
\end{proof}

We are now ready to prove Proposition~\ref{prop: N0 for sigma to sigma through Pin}.

\begin{proof}[Proof of Proposition~\ref{prop: N0 for sigma to sigma through Pin}]
We decompose $f_\Sigma^j \, f_{\Sigma, \iin} \, \psi \, f_{\out, \Sigma} \, f_\Sigma^i$ in three pieces:
\begin{itemize}
    \item $f_\Sigma^i$ leaves the cones $C^u_\Sigma$ invariant and does not contract the vectors of $C^u_\Sigma$ for $i \geq 0$ (Proposition~\ref{prop: cones Sigma}, Item~\ref{prop: cones Sigma; it: invariance sigma});
    \item $ f_{\Sigma, \iin} \, \psi \, f_{\out, \Sigma}$ maps $C^u_\Sigma$ to a $(\zeta^u/\,\zeta^s)$-cone on $C^u_\Sigma$, of uniformly bounded slope, and contracts the vectors of $C^u_\Sigma$ by a factor uniformly bounded below, and the constants are uniform in $\delta$ (Lemma~\ref{lem: contraction crossing boundary});
    \item $f_\Sigma^j$ maps a $(\zeta^u/\,\zeta^s)$-cone field of bounded slope inside a $(\zeta^u/\,\zeta^s)$-cone field of arbitrarily small slope for $j$ large enough, and expands its vectors arbitrarily strongly (Claim~\ref{claim: cones return on sigma}).
\end{itemize}
We deduce the existence of an integer $N_0$ which satisfies the first item if $j \geq N_0$.
The proof of the second item is symmetric.
\end{proof}

\subsubsection*{From $\Pin \ssm \cV_\delta$ to $\Sigma$ through $\Pin \cap \cV_\delta$}

\begin{figure}[htb]
    \centering
    \vspace*{-5em}
    \includegraphics[width=\textwidth]{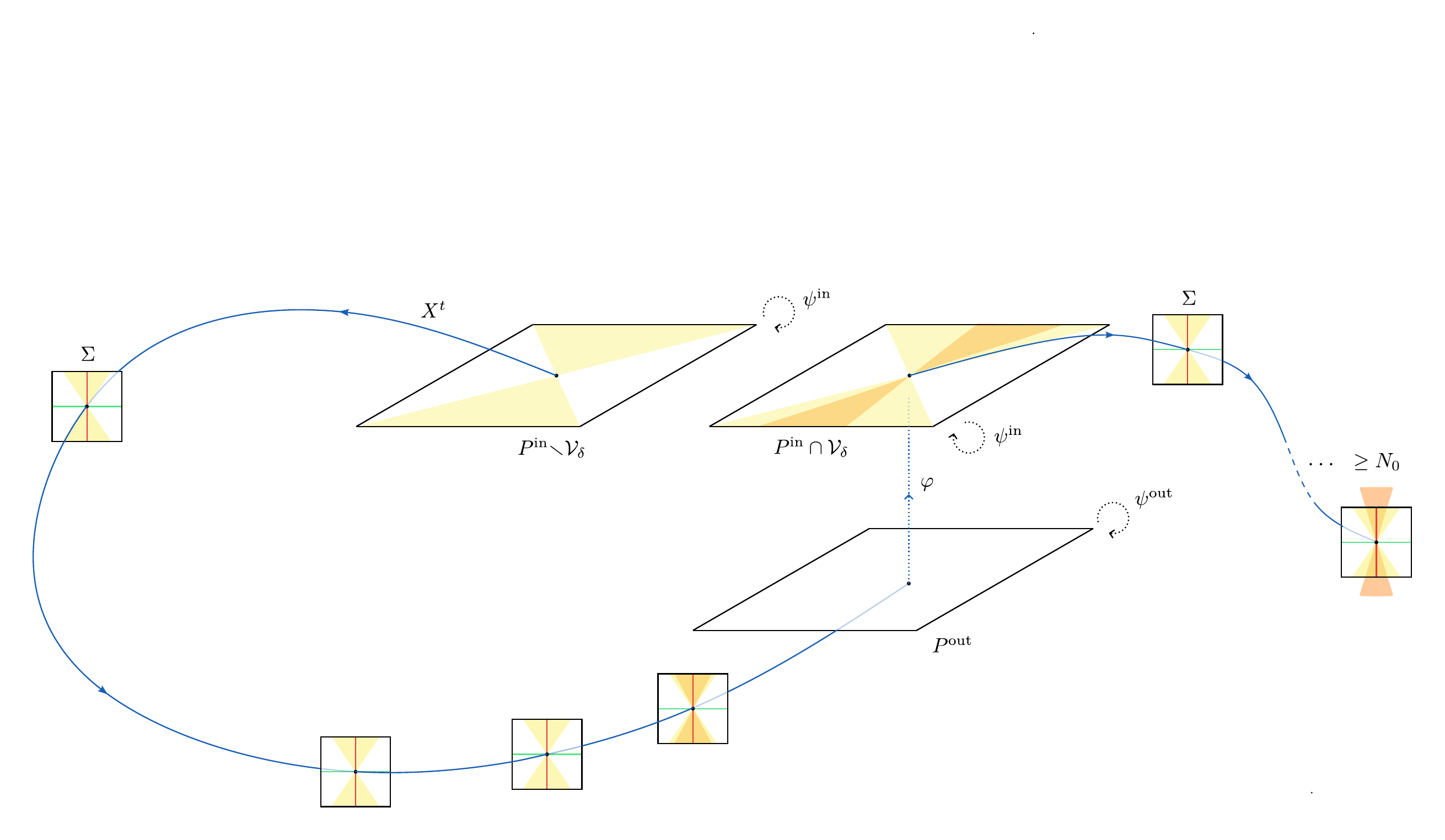}
     \vspace*{-1em}
    \caption{Trajectory of a cone $C^u_\iin$ under the map $f_\Sigma^j \, f_{\Sigma, \iin} \, \psi \, f_{\out, \Sigma} \, f_\Sigma^i \, f_{\Sigma, \iin} \, \psi^\iin \colon \Pin \ssm \cV_\delta \to \Sigma$ crossing $\Pin \cap \cV_\delta$}
    \label{fig: Pin minus V0 to Sigma through Pin inter V0}
\end{figure}

Now we want to show that the composition map
\begin{equation} \label{eq: map from Pin minus V_delta to Sigma through Pin}
    f_\Sigma^j \, f_{\Sigma, \iin} \, \psi \, f_{\out, \Sigma} \, f_\Sigma^i \, f_{\Sigma, \iin} \, \psi^\iin \colon \Pin \ssm \cV_\delta \to \Sigma
\end{equation}
satisfies the cone fields condition for the pair $(C^u_\iin \cup C^u_\Sigma, C^s_\iin \cup C^s_\Sigma)$
on the set of points $p \in \Pin \ssm \cV_\delta$ whose image by $\varphi \, \psi^\out  \, f_{\out, \Sigma} \, f_\Sigma^i \, f_{\Sigma, \iin} \, \psi^\iin$ belongs to $\Pin \cap \cV_\delta$.
We refer to Figure~\ref{fig: Pin minus V0 to Sigma through Pin inter V0} and the following proposition.

\begin{prop}
\label{prop: N0 for Pin minus V0 to Sigma through Pin inter VO}
There exists $N_0$ such that
for all $\delta>0$ small enough,
the following condition is satisfied.
Let $p \in \Pin \ssm \cV_\delta$ such that the following are well defined for $i\geq 0$ and $j \geq N_0$
\begin{itemize}[--]
    \item $p' = \varphi \, \psi^\out  \, f_{\out, \Sigma} \, f_\Sigma^i \, f_{\Sigma, \iin} \, \psi^\iin(p) \in \Pin \cap \cV_\delta $, and
    \item $q =  f_\Sigma^j \, f_{\Sigma, \iin} \, \psi \, f_{\out, \Sigma} \, f_\Sigma^i \, f_{\Sigma, \iin} \, \psi^\iin(p) \in \Sigma$.
\end{itemize}
Then
\begin{itemize}[--]
    \item $f_\Sigma^j \, f_{\Sigma, \iin} \, \psi \, f_{\out, \Sigma} \, f_\Sigma^i \, f_{\Sigma, \iin} \, \psi^\iin$ maps $C^u_\iin(p)$ inside $C^u_{\Sigma}(q)$ and expands the norm of the vectors of $C^u_\iin(p)$ by a factor $2$,
    \item $(f_\Sigma^j \, f_{\Sigma, \iin} \, \psi \, f_{\out, \Sigma} \, f_\Sigma^i \, f_{\Sigma, \iin} \, \psi^\iin)^\inv$ maps $C^s_\Sigma(q)$ inside $C^s_\iin(p)$ and expands the norm of the vectors of $C^s_\Sigma(q)$ by a factor $2$.
\end{itemize}
\end{prop}

The key result is that the contraction along an orbit of the flow of $X$ from $\Pin \ssm \cV_\delta$ to $\Sigma \cap \cV_\delta$ is bounded below, uniformly in $\delta$ (Lemma~\ref{lem: contraction from Pin minus V_delta to V_delta uniform in delta}).
More precisely, let us show that the composition $f_\Sigma^n, f_{\Sigma, \iin} \, \psi^\iin$ restricted to $\Pin \ssm \cV_\delta$ at the domain and $\Sigma \cap \cV_\delta$ at the destination has a contraction bounded below on the norm of the vectors of $C^u_\iin$, and the bound is uniform in $\delta$, moreover the map contributes to close a $(\zeta^u/\,\zeta^s)$-cone field.

\begin{lem} \label{lem: contraction from Pin minus V_delta to V_delta uniform in delta}
There exists a constant $\cst>0$ such that for any $\delta>0$ uniformly small enough, if $\psi^\iin = \psi^\iin_\delta$,
for any $p \in \Pin \ssm \cV_\delta$ such that there exists $n \geq 1$ with the image
$q = f_\Sigma^n \, f_{\Sigma, \iin} \, \psi^\iin(p) \in \Sigma \cap \cV_\delta$ well defined, then
\begin{itemize}[--]
\item $(f_\Sigma^n \, f_{\Sigma, \iin} \, \psi^\iin)_* C^u_\iin(p)\, \subset \, \intr C^u_\Sigma(q), \ \mb{and} \ \forall \, v \in C^u_\iin(p)$, $$\Vert (f_\Sigma^n \, f_{\Sigma, \iin} \, \psi^\iin)_* v \Vert \geq \cst \Vert v \Vert;$$
\item $(f_\Sigma^n \, f_{\Sigma, \iin} \, \psi^\iin)^\inv_* C^s_\Sigma(q) _, \subset \, \intr C^s_\iin(p), \ \mb{and} \ \forall v \in C^s_\Sigma(q)$, $$\ \Vert (f_\Sigma^n \, f_{\Sigma, \iin} \, \psi^\iin)^\inv_* v \Vert \geq \cst \Vert v \Vert.$$
\end{itemize}
\end{lem}

\begin{proof}
Let us show the first item, the proof of the second being similar.
In order not to introduce new constants, we can assume that each normalized neighborhood $\cV_i$ of $\cO_i$ contains a $1$-neighborhood of $\cO_i$.
We denote $\cV_1$ the union of $1$-neighborhoods of $\cO_i$ in $\cV_i$.
The proof relies on the following key points.
\begin{itemize}
    \item The differential of $f_{\Sigma, \iin}$ is bounded below by a factor uniformly proportional to $\dist(\cdot, \cO_*)$ from Lemma~\ref{lem: contraction sigma to boundary is distance to orbit}, so greater than $\delta$ on $\Pin \ssm \cV_\delta$.
    \item The orbit of any point $p \in \Pin$ exits a $1$-neighborhood of $\cO_*$ before exiting $P$.
    \item The number $n$ of crossings through $\Sigma$ of the orbit of a point crossing from $\cV_1(\cO_i)$ to $\cV_\delta(\cO_i)$ is of the order of $\log_2(\delta^\inv)$.
    The corresponding iteration $f_\Sigma^{n}$ of the return map on $\Sigma$ expands the norm of the vectors $v \in E^\uu$ by a factor greater than $\delta^\inv$ (up to a uniform multiplicative constant), which will compensate for the contraction of $f_{\Sigma, \iin}$.
\end{itemize}
Let us insist on the fact that we will have an exact compensation of a contracting effect of order $\delta$ and a expanding effect of order $\delta^\inv$.
The contraction of the composition map $f_\Sigma^n , f_{\Sigma, \iin} \, \psi^\iin \colon \Pin \ssm \cV_\delta \to \Sigma \cap \cV_\delta$ will thus be bounded below by a constant \emph{independent of $\delta$}.
We formally state the last two items.
Let $l$ be the minimum diameter of the disks of section $\Sigma$ and $0< \delta< l$.

\begin{claim} \label{claim: orbit exit V_1}
Let $p \in \Pin \ssm \cV_\delta$. Then there exists $T \geq 0$ such that for all $t \leq T$, $X^t(p) \in P \ssm \cV_\delta$ and $X^T(p) \in P \ssm \cV_1$.
\end{claim}

\begin{proof}
If $p \in \Pin \ssm \cV_1$, the fact is true for $T=0$.
Otherwise, $p \in \Pin \cap \cV_i$, where $\cV_i$ is the normalized neighborhood of the orbit $\cO_i \in \cO_*$ for some $i$, given coordinates $\xi_i = (x,y,\theta)$.
We have $p = (x,x,\theta)$, with $\vert x \vert > \delta$ and $$X^t(x,y,\theta) = (2^{-t} x, 2^t x, \theta +t).$$
The distance of $X^t(p)$ to the orbit $\cO_i$ is therefore greater than $2^t \delta$.
We deduce that the orbit of $p$ leaves $\cV_1$ in time $t=\log_2(\delta^\inv)$ without intersecting $\cV_\delta$.
\end{proof}

\begin{claim} \label{claim: from V_1 to V_delta}
Let $q \in \Sigma \cap \cV_\delta$ be such that there exists $T \geq 0$ with $X^{-T}(q) \in P \ssm \cV_1$.
Then the orbit segment $[X^{-T}(q), q]$ intersects the section $\Sigma$ a number of times $n_0 \geq \log_2( l \delta^\inv) -1$  in a normalized neighborhood $\cV_j$ of $\cO_j \in \cO_*$.
\end{claim}

\begin{proof}
Let $(\cV_j, \xi_j)$ be an \ncs{} of $\cO_j \in \cO_*$ such that $q \in \cV_j$.
Let $\cV_\delta \subset \cV_j$ be the $\delta$-neighborhood of $\cO_j$.
By definition of $l$, and by Definition~\ref{def: normalized block}, Item~\ref{def: normalized, it: affine section}, the section $\Sigma$ contains the set $R_l := \{\theta=0, \, \vert x \leq l, \, \vert y \vert \leq l\}$.
Let $q = (x, y, \theta)$ with $x^2 + y^2 \leq \delta^2$ and we have $X^{-t} (p) = (2^{t} x, 2^{-t} y, \theta -t)$ for $t \geq 0$.
Then the negative orbit $X^{-t}(p)$ intersects the set $R_l \subset \Sigma$ a number greater than $n = \log_2 (l \delta^\inv) -1$ times before exiting $\cV_1$.
\end{proof}

Now let $p \in \Pin \ssm \cV_\delta$ such that $q = f^n_\Sigma f_{\Sigma, \iin} \psi^\iin (p)$ is well defined in $\Sigma \cap \cV_\delta$.
According to Claim~\ref{claim: orbit exit V_1}
the orbit of $p$ exits a $1$-neighborhood at a point $p' = X^T(p)$.
According to Claim~\ref{claim: from V_1 to V_delta},
the orbit segment from $p'$ to $q$ intersects $n_0 \geq \log_2(l \delta^\inv) -1$ times the $\Sigma$ section in a normalized $\cV_j$ neighborhood of a $\cO_j \in \cO_*$ orbit.
In the neighborhood of $p$ we have the following equality
$$ f^n_\Sigma f_{\Sigma, \iin} \psi^\iin = f^{n_0}_\Sigma f^{n-n_0}_\Sigma f_{\Sigma, \iin} \psi^\iin$$
with $n-n_0 \geq 0$ and the map $f^{n_0}_\Sigma$ is the iterate $n_0$ of the restriction of $f_\Sigma$ to $\Sigma \cap \cV_j$.
We can study the effect on the norm of the vectors of each map of the composition $f^n_\Sigma f_{\Sigma, \iin} \psi^\iin$ in the neighborhood of a point $p \in \Pin \ssm \cV_\delta$.

\begin{enumerate}[leftmargin=*]

    \item \label{it: composition; psi^in} 
    The \diff{} $\psi^\iin = \psi^\iin_\delta$ maps $C^u_\iin$ inside a $(\hat K,\cG^{u, \iin}/\, \cG^{s, \iin})$-cone field, where $\hat K$ is a constant uniform in $\delta$ (Lemma~\ref{lem: action diffs slope C^u,s_in}).
    Its differential is bounded (below) by a constant $\cst$ uniform in $\delta$ (Lemma~\ref{lem: differential of diffs bounded uniform in delta}).
    
    \item \label{it: composition; f_sigma,in} The crossing map $f_{\Sigma, \iin}$ maps a $(\hat K,\cG^{u, \iin}/\cG^{s, \iin})$-cone field inside a $(\zeta^u /\,\zeta^s)$-cone field on $\Sigma$ of slope $K_1 \leq \cst \hat K$ (Lemma~\ref{lem: slope cone from Sigma to boundary}).
    Its differential at a point $\psi^\iin(p)$ is bounded below by $\cst \dist(\psi^\iin(p),\cO_*)$
    (Lemma~\ref{lem: contraction sigma to boundary is distance to orbit}).
    According to~\ref{lem: action diffs distance orbit}, this factor is bounded below by $\cst \dist(p,\cO_*) > \cst \delta$.
    
    \item \label{it: composition; f_sigma^n-n0}
    The iteration $f_\Sigma^{n-n_0} \colon \Sigma \to \Sigma$ maps a $(K_1,\zeta^u /\,\zeta^s)$-cone field to a $(K_1, \zeta^u /\,\zeta^s)$-cone field (Claim~\ref{claim: cones return on sigma}) and contracts the vectors of a $(K_1,\zeta^u /\,\zeta^s)$-cone by a factor bounded below by a $\cst>0$ which depends only on $K_1$.
    
    \item \label{it: composition; f_sigma^n0} Finally, we can use Claim~\ref{claim: cones return on sigma} with an expansion factor $\mu=2$ and an iterate $n_0 \geq \log_2(l \delta^\inv) -1$ to show that $f_\Sigma^{n_0}$ maps a $(K_1,\zeta^u/\, \zeta^s)$-cone inside a $(\zeta^u/\,\zeta^s)$-cone field of slope $K_2 \leq \cst \, \delta^2 \, K_1$, where the constant depends only on $l$.
    It expands the vectors of a $(K_1,\ \zeta^u/\, \zeta^s)$-cone by a factor greater than $\cst \delta$, where the constant depends only on $K_1$ and $l$.
\end{enumerate}

We deduce that one can choose $\delta$ small enough so that the image \linebreak[4]$(f^n_\Sigma f_{\Sigma, \iin} \psi^\iin)_* C^u_\iin$ is contained inside $C^u_\Sigma$.
Moreover, up to change the constant $\cst$ which is uniform in $\delta$, we have in the neighborhood of $p$:
\begin{align*}
    \l\Vert \l( f_\Sigma^n f_{\Sigma, \iin} \psi^\iin \r)_* v \r\Vert &= \l\Vert \l( f_\Sigma^{n_0} f_\Sigma^{n-n_0} f_{\Sigma, \iin} \psi^\iin \r)_* v \r\Vert\\  
    & \geq \cst\, \delta^\inv \l\Vert \l(f_{\Sigma, \iin} \psi^\iin \r)_* v \r\Vert & \mb{(Item~\ref{it: composition; f_sigma^n-n0} and~\ref{it: composition; f_sigma^n0})} \\  
    & \geq \cst\, \delta^\inv  \delta \Vert (\psi^\iin)_*v \Vert
    & \mb{(Item~\ref{it: composition; f_sigma,in})} \\
    & \geq \cst\, \Vert v \Vert.
    & \mb{(Item~\ref{it: composition; psi^in})}
\end{align*} 
This completes the proof of Claim~\ref{claim: from V_1 to V_delta}.
\end{proof}

\begin{rmk} \label{rmk: contraction from Pin minus V_delta to V_r.delta uniform in delta}
It is easy to check that Lemma~\ref{lem: contraction from Pin minus V_delta to V_delta uniform in delta} is still true (up to change the constant) for any orbit trajectory starting from $p \in \Pin \ssm \cV_\delta$ and arriving at
$q = f_\Sigma^n \, f_{\Sigma, \iin} \, \psi^\iin(p) \in \Sigma \cap \cV_{r.\delta}$, where $r>0$ is a uniform constant previously fixed.
\end{rmk}

We are now ready to prove Proposition~\ref{prop: N0 for Pin minus V0 to Sigma through Pin inter VO}.

\begin{proof}[Proof of Proposition~\ref{prop: N0 for Pin minus V0 to Sigma through Pin inter VO}]
Let us show the first item concerning unstable cones, the proof of the second is similar.
The fact that the point $$p' = \varphi \, \psi^\out  \, f_{\out, \Sigma} \, f_\Sigma^i \, f_{\Sigma, \iin} \, \psi^\iin(p)$$ belongs to $\Pin \cap \cV_\delta$ (by assumption of the corollary) implies that the point $f_\Sigma^i \, f_{\Sigma, \iin} \, \psi^\iin (p)$ belongs to $\Sigma \cap \cV_{r.\delta}$, in other words the positive orbit of $p$ intersects the $\Sigma$ section one last time in the neighborhood $\cV_{r.\delta}$ before exiting through $\Pout$, where $r>0$ is a (uniform) constant which satisfies Lemma~\ref{lem: action diffs distance orbit}.

We decompose the map 
\mb{$f_\Sigma^j \, f_{\Sigma, \iin} \, \psi \, f_{\out, \Sigma} \, f_\Sigma^i \, f_{\Sigma, \iin} \, \psi^\iin$}
in the neighborhood of $p$ in three pieces:
\begin{itemize}
    \item The map $f_\Sigma^i \, f_{\Sigma, \iin} \, \psi^\iin \colon \Pin \ssm \cV_\delta \to \Sigma \cap \cV_{r .\delta}$ is addressed in Lemma~\ref{lem: contraction from Pin minus V_delta to V_delta uniform in delta} and Remark~\ref{rmk: contraction from Pin minus V_delta to V_r.delta uniform in delta}.
    For $\delta$ small enough, it maps the cone $C^u_\iin(p)$ inside $C^u_\Sigma$ and contracts the vectors of $C^u_\iin$ by a factor bounded below by a constant $c_0>0$ uniform in $\delta$.
    \item The map $f_{\Sigma, \iin} \psi f_{\out, \Sigma} \colon \Sigma \to \Sigma$ is addressed in Lemma~\ref{lem: contraction crossing boundary}. Its differential is bounded by a constant $c_1>0$ and it maps the cone field $C^u_\Sigma$ inside a $(K_1,\ \zeta^u/\, \zeta^s)$-cone field, with $c_1$ and $K_1$ uniform in $\delta$.
    \item From Claim~\ref{claim: cones return on sigma}, there exists an integer $N_0$, such that if $j \geq N_0$, then $f_\Sigma^j$ maps a $(K_1,\ \zeta^u/\,\zeta^s)$-cone field inside $C^u_\Sigma$ and expands the norm of the vectors of a $(K_1,\ \zeta^u/\,\zeta^s)$-cone by a factor greater than $2 (c_1 c_0)^\inv$, which completes the proof.
\qedhere\end{itemize}

\end{proof}

\textbf{We fix the integer $N_0$ which satisfies Propositions~\ref{prop: N0 for sigma to sigma through Pin} and~\ref{prop: N0 for Pin minus V0 to Sigma through Pin inter VO}}.

\subsection{Parameter \texorpdfstring{$\delta_0$}{¼¥\_0}}
\label{sec: parameters; subsec: delta}

In this Subsection, we show the existence of a parameter $\delta_0$ which will allow us to obtain Item~\ref{prop: parameters; it: v0} and Item~\ref{prop: parameters; it: Pin to Pin} of Proposition~\ref{prop: parameters and cones}.

\begin{figure}[htb]
    \centering
    \captionsetup{width=.83\linewidth}
    \includegraphics[width=\textwidth]{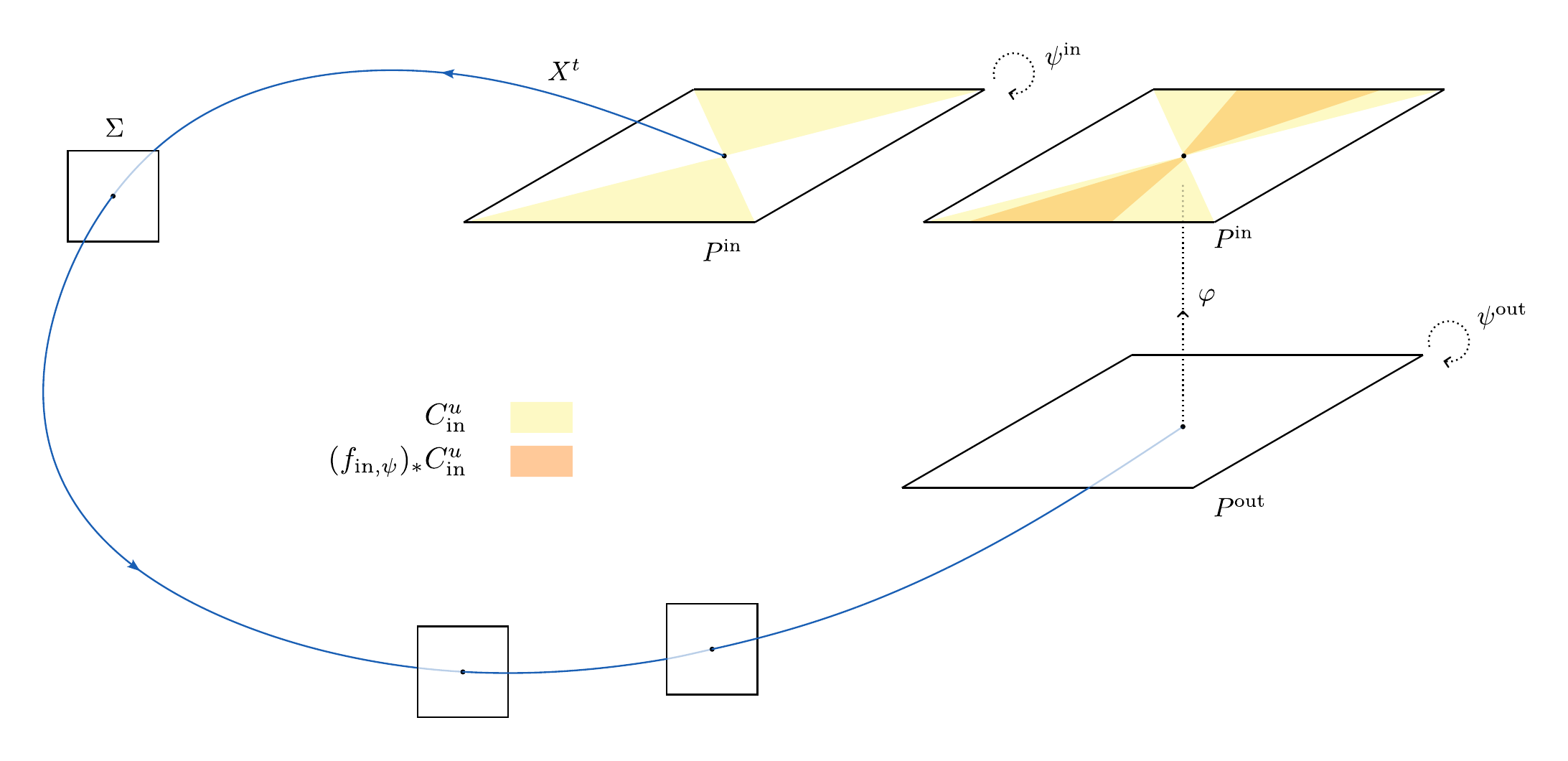}
    \vspace*{-1em}
    \caption{Trajectory of a cone $C^u_\iin$ under the map $f_{\iin, \psi}$ on~$\Pin$}
    \label{fig: cone pin to pin compatible}
\end{figure}

    \begin{itemize}
        \item Item~\ref{prop: parameters; it: v0} concerns the $X_\psi$ orbits crossing $(\cV_0)_\psi= (\psi^\iin (\cV_{\delta_0}))_\psi$ (Figure~\ref{fig: diagram orbite; subfig: v0}).
        We will study the composition maps
        $f_\Sigma^{n} \, f_{\Sigma, \iin} \, \psi^\iin$ and $f_\Sigma^{-n} \, f_{\out, \Sigma}^\inv \, (\varphi \psi^\out)^\inv$ from $\Pin \cap \cV_\delta$ to $\Sigma$.
    
        \item Item~\ref{prop: parameters; it: Pin to Pin} concerns points $p \in \Pin_\psi$ whose future orbits by $X_\psi$ intersects $\Pin_\psi$ a first time at point $q = f^n_{0, \psi}(p)$ (Figure~\ref{fig: cone pin to pin compatible}).
        We know from Subsection~\ref{sec: parameters; subsec: lambda epsilon cone Pin} that if $p$ and $q$ are in $(\Pin \cap \ssm \cV_\delta)_\psi$ then the cones are invariant and expanded by $f^n_{0, \psi}$ (uniformly in $\delta$).
        It remains to check the cone invariance if $p$ or $q$ are in $(\Pin \cap \hat \cV_\delta)_\psi$.
        We will study the map $$f_{\iin, \psi} = \varphi \, \psi^\out \, f_{\out, \Sigma} \, f_\Sigma^{n} \, f_{\Sigma, \iin} \, \psi^\iin \colon \Pin \to \Pin.$$
    \end{itemize}

The conditions are summarized in the following proposition.

\pagebreak[4]
\begin{prop} \label{prop: delta_0}
There exists $\delta_0$ such that:
\begin{enumerate}
    \item \label{prop: delta_0; it: crossing Sigma N0}
    The maps $f_\Sigma^{N_0} \, f_{\Sigma, \iin} \, \psi^\iin$ and $f_\Sigma^{-N_0} \, f_{\out, \Sigma}^\inv \, (\varphi \psi^\out)^\inv$ are well defined from \linebreak[4]$\Pin \cap \cV_{\delta_0}$ to $\Sigma$.
    \item \label{prop: delta_0; it: Pin invariance}
    For all $p \in \Pin$ such that $f_{\iin, \psi}(p)=:q$ is well defined in $\Pin$, the map $f_{\iin, \psi}$ maps $C^u_\iin(p)$ inside $C^u_\iin(q)$, and the map $f_{\iin, \psi}^\inv$ maps $C^s_\iin(q)$ inside $C^s_\iin(p)$.
\end{enumerate}
\end{prop}

\begin{proof}
Let us start with the following result.
\begin{claim} \label{claim: delta(n)}
    For every integer $n>0$, there exist $\delta(n)>0$ such that for all $\delta \leq \delta(N)$ the maps $f_\Sigma^n \, f_{\Sigma, \iin} \, \psi^\iin$ and $f_\Sigma^{-n}  = \, f_{\out, \Sigma}^\inv \, (\varphi \psi^\out)^\inv$ are well defined from $ \Pin \cap \cV_\delta$ to $\Sigma$. 
\end{claim}
\begin{proof}
    A point $p \in \cap \cV_\delta$ is in a $\delta$-neighborhood of a periodic orbit $\cO_i \in \cO_*$.
For arbitrarily small $\delta$, the orbit by the flow of $X$ of a point $p \in \Pin \cap \cV_\delta$ spends an arbitrarily long time in the neighborhood of $\cO_* \subset \Lambda_s$ so intersects an arbitrarily large number of times the section $\Sigma$ before intersecting $\Pout$.
It suffices to recall that the diffeomorphisms $\psi^\iin$ and $\varphi \psi^\out$ have a bounded effect on the distance $\dist(\cdot, \cO_*)$, uniformly in $\delta$ (Lemma~\ref{lem: action diffs distance orbit}).
\end{proof}

We now state a lemma about the projective action of $f_\Sigma^n f_{\Sigma, \iin} \psi^\iin : \Pin \to \Sigma$.

\begin{lem} \label{lem: m0 for invariance from pin to sigma}
    There exist an integer $M_0 >0$ such that for all $\delta>0$, if $p \in \Pin$ and $q = f_\Sigma^n f_{\Sigma, \iin} \psi^\iin \in \Sigma$, with $n \geq M_0$, then
    $f_\Sigma^n f_{\Sigma, \iin} \psi^\iin$ maps $C^u_\iin(p)$ inside $C^u_\Sigma(q)$ and the inverse $(f_\Sigma^n f_{\Sigma, \iin} \psi^\iin)^\inv$ maps $C^s_\iin(q)$ inside $C^s_\iin(p)$.
\end{lem}

\begin{proof}
We write the proof for the unstable cone field.
We study each part of the composition in the following way:
\begin{itemize}
    \item The \diff{} $\psi^\iin = \psi^\iin_\delta$ maps $C^u_\iin$ inside a $(\hat K,\cG^{u, \iin}/\, \cG^{s, \iin})$-cone field, where $\hat K$ is a constant uniform in $\delta$ (Lemma~\ref{lem: action diffs slope C^u,s_in}). 
    
    \item The crossing map $f_{\Sigma, \iin}$ maps a $(\hat K,\cG^{u, \iin}/\cG^{s, \iin})$-cone field inside a $(\zeta^u /\,\zeta^s)$-cone field on $\Sigma$ of slope $K_1 \leq \cst \hat K$ (Lemma~\ref{lem: slope cone from Sigma to boundary}).

    \item $f_\Sigma^n$ maps a $(\zeta^u/\,\zeta^s)$-cone field of bounded slope inside a $(\zeta^u/\,\zeta^s)$-cone field of arbitrarily small slope for $n$ large enough (uniformly in $\delta$) (Claim~\ref{claim: cones return on sigma}).

\end{itemize}
We deduce the existence of an integer $M_0$ such that the lemma is satisfied.
\end{proof}

Now we can prove Proposition~\ref{prop: delta_0}.
From Claim~\ref{claim: delta(n)} we can choose $\delta_0 < \min (\delta(N_0), \delta(M_0))$ where $M_0$ is given by Lemma~\ref{lem: m0 for invariance from pin to sigma} and $N_0$ is the integer of the Proposition.
Item~\ref{prop: delta_0; it: crossing Sigma N0} is then a direct consequence of Claim~\ref{claim: delta(n)} with $N=N_0$.

Now let us prove Item~\ref{prop: delta_0; it: Pin invariance}.
If $p \in \Pin \ssm \cV_{\delta_0}$ and $q \in \Pin \ssm \cV_{\delta_0}$, it follows from Proposition~\ref{prop: lambda epsilon cones Pin}, Item~\ref{prop: cones Pin; it: invariance}.
If $p \in \Pin \cap \cV_{\delta_0}$ or $q \in \Pin \cap \cV_{\delta_0}$.
then we write
$f_{\iin, \psi}(p) = \varphi \psi^\out f_{\out, \Sigma}  f_\Sigma^n  f_{\Sigma, \iin}  \psi^\iin (p) $, and we know that $n \geq M_0$ from Claim~\ref{claim: delta(n)}.
Define $p' = f_\Sigma^n  f_{\Sigma, \iin}  \psi^\iin \in \Sigma$.
We decompose in the following way
\begin{itemize}
    \item $f_\Sigma^n  f_{\Sigma, \iin}  \psi^\iin $ maps $C^u_\iin(p)$ inside $C^u_\Sigma(p')$ according to Lemma~\ref{lem: m0 for invariance from pin to sigma}.
    \item $ \varphi \psi^\out f_{\out, \Sigma} $ maps $C^u_\Sigma(p')$ to $C^u_\iin(q)$ according to Proposition~\ref{prop: cones Sigma}, Item~\ref{prop: cones Sigma; it: compatibilty pin}.
\end{itemize}
This completes the proof.
\end{proof}

\textbf{We fix the parameter $\delta_0$ which satisfies Proposition~\ref{prop: delta_0}.}


\subsection{Integer \texorpdfstring{$N_1$}{N\_1}}

\label{sec: parameters; subsec: N1}

The parameters $\lambda = \lambda_0$, $\epsilon = \epsilon_0$ and $\delta = \delta_0$ are fixed, as well as the cones $(C^u_\iin, C^s_\iin)$ on $\Pin$ and $(C^u_\Sigma, C^s_\Sigma)$ on $\Sigma$.
From now on we call $\psi^\iin = \psi^\iin_{\lambda_0, \epsilon_0, \delta_0}$, $\psi^\out = \psi^\out_{\lambda_0, \epsilon_0, \delta_0}$ the \diffs{} of $\pP$ given by Proposition~\ref{prop: spread expansion} and~\ref{prop sym: spread expansion}.
Let $\psi := \psi^\iin \varphi \psi^\out$ be the composition map.
Let $\cV_0 := \psi^\iin(\cV_{\delta_0})$.
In this section, we show the existence of an integer $N_1$ which will allow us to have Items~\ref{prop: parameters; it: Pin minus v0 to sigma} and~\ref{prop: parameters; it: sigma to Pin minus v0} of Proposition~\ref{prop: parameters and cones}.
Recall that Item~\ref{prop: parameters; it: Pin minus v0 to sigma} concerns the points $p \in (\Pin \ssm \cV_0)_\psi$ whose future orbit by $X_\psi$ intersects $\Sigma_\psi$ a number $N_1$ of consecutive times,
and Item~\ref{prop: parameters; it: sigma to Pin minus v0} concerns the points $p \in \Sigma_\psi$ whose future orbit by $X_\psi$ intersects $(\Pin \ssm \cV_0)_\psi$ a number $N_1$ of times.

\begin{figure}[htb]
    \centering
    \vspace*{-1em}
    \captionsetup{width=.85\linewidth}
    \includegraphics[width=0.9\textwidth]{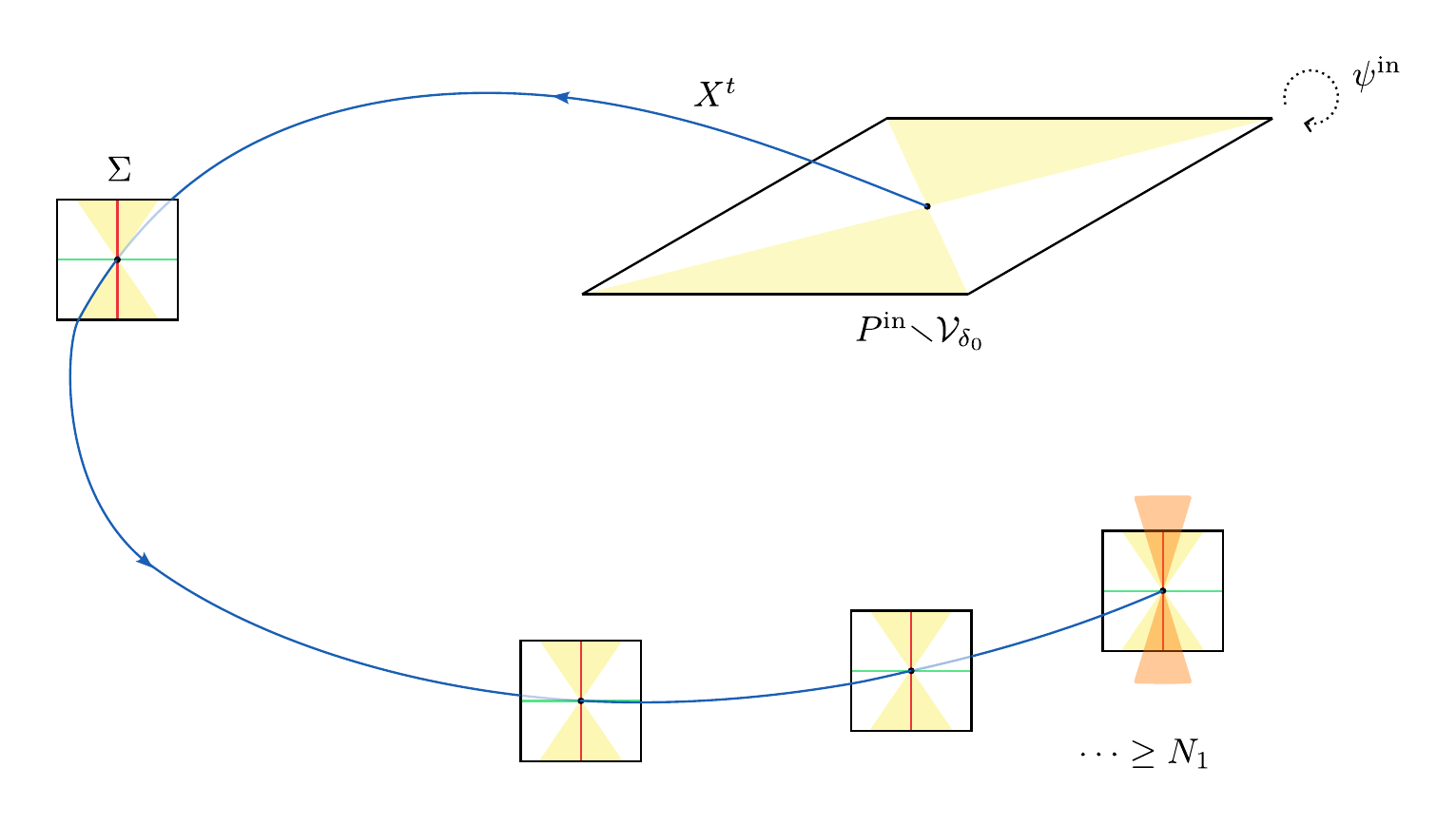}
    \vspace*{-1em}
    \caption{Cone $C^u_\iin$ mapped by $f^j_\Sigma \, f_{\Sigma, \iin} \si^\iin$ restricted to $\Pin \ssm \cV_{\delta_0}$ }
    \label{fig: Pin minus v0 to sigma}
\end{figure}

Let us start by showing that the map $f^j_\Sigma \, f_{\Sigma, \iin} \, \psi^\iin \colon \Pin \to \Sigma$
restricted to $\Pin \ssm \cV_{\delta_0}$ satisfies the cone fields condition for the pair $(C^u_\iin \cup C^u_\Sigma, C^s_\iin \cup C^s_\Sigma)$ for a sufficiently large iterate $j$.
We refer to Figure~\ref{fig: Pin minus v0 to sigma}.

\begin{lem} \label{lem: N1 for Pin minus v0 to Sigma}
There exists $N_1$ such that for any $p \in \Pin \ssm \cV_{\delta_0}$ if the image $q = f^n_\Sigma \, f_{\Sigma, \iin} \, \psi^\iin(p) \in \Sigma$ is well defined for a $n \geq N_1$,
then
\begin{itemize}[--]
    \item $(f^n_\Sigma \, f_{\Sigma, \iin}  \psi^\iin)_* C^u_\iin(p) \subset \intr C^u_\Sigma (q)$ and 
    $\Vert (f^n_\Sigma \, f_{\Sigma, \iin}  \psi^\iin)_* v \Vert \geq 2 \Vert v \Vert$ for all $v \in C^u_\iin(p)$
    \item $(f^n_\Sigma \, f_{\Sigma, \iin}  \psi^\iin)^\inv_* C^s_\Sigma(q) \subset \intr C^s_\Sigma (p)$ and
    $\Vert (f^n_\Sigma \, f_{\Sigma, \iin}  \psi^\iin)^\inv_* v \Vert \geq 2 \Vert v \Vert$ for all $v \in C^s_\Sigma(p)$
\end{itemize}
\end{lem}

\vspace*{-1em}
\begin{proof}
We show the existence of an integer $N_1$ such that the first item is true. 
The second item is shown in a similar way.
Let us study each map of the composition.
\begin{itemize}
    \item The cone field $\psi^\iin(C^u_\iin)$ is inside a $(\hat K,\cG^{u,\iin}/, \cG^{s, \iin})$-cone field (Lemma~\ref{lem: action diffs slope C^u,s_in}) and $\psi^\iin$ has a differential uniformly bounded by a constant $c_0>0$.
    
    \item The crossing map $f_{\Sigma, \iin} \colon \Pin \to \Sigma$ maps a $(\hat K,\cG^{u, \iin}/\, \cG^{s, \iin})$-cone field on $\Pin$ inside a $(\zeta^u /\,\zeta^s)$-cone field on $\Sigma$ of uniform slope $K_1$ (Lemma~\ref{lem: slope cone from Sigma to boundary}).
    Moreover, the differential of $f_{\Sigma, \iin}$ restricted to $\Pin \ssm \cV_{\delta_0}$ is uniformly bounded by a constant $c_1>0$ because it is a crossing map of the flow $X^t$ between a compact surface $\Pin \ssm \cV_{\delta_0}$ transverse to the \vf{} $X$ and a compact surface $\Sigma$  transverse to the \vf{} $X$ and the time of the crossing is uniformly bounded.
    
    \item By hyperbolicity of $f_\Sigma$ (Claim~\ref{claim: cones return on sigma}), there exists an integer $N_1$ such that if $n \geq N_1$, then $f_\Sigma^n$ maps a $(K_1,\zeta^u/\,\zeta^s)$-cone field inside $C^u_\Sigma$ and expands the norm of the vectors of a $(K_1,\zeta^u/\,\zeta^s)$-cone by a factor greater than $2 (c_1 c_0)^\inv$.
\end{itemize}
Putting together we obtain the existence of an integer $N_1$ such that the composition $f^j_\Sigma \, f_{\Sigma, \iin} \, \psi^\iin \colon \Pin \to \Sigma$ with $n \geq N_1$ maps the cone $C^u_\iin(p)$ inside $C^u_\Sigma(q)$ and expands the vectors of $C^u_\iin$ by a factor $2$.
\end{proof}

We recall the definition in Subsection~\ref{sec: parameters; subsec: lambda epsilon cone Pin} of the following composition
$f_{\iin, \psi} = \varphi \psi^\out \foutin \psi^\iin \colon \Pin \to \Pin$.
Denote $\hat f_{\iin, \psi} \colon \Pin \ssm \cV_{\delta_0} \to \Pin \ssm \cV_{\delta_0}$
the restriction of $f_{\iin, \psi}$ to the set $\Pin \ssm \cV_{\delta_0}$ at the domain and destination. 
{\sloppy Let us show that the composition map
$$ \hat f_{\iin, \psi}^j \varphi \psi^\out f_{\out, \Sigma} f_\Sigma^i \colon \Sigma \to \Pin \ssm \cV_{\delta_0} $$ 
satisfies the cone fields condition for the pair $(C^u_\Sigma \cup C^u_\iin, C^s_\Sigma \cup C^s_\iin)$ if $j$ is large enough.}

\begin{figure}[htb]
    \centering
    \vspace*{-1em}
    \captionsetup{width=.95\linewidth}
    \includegraphics[width=\textwidth]{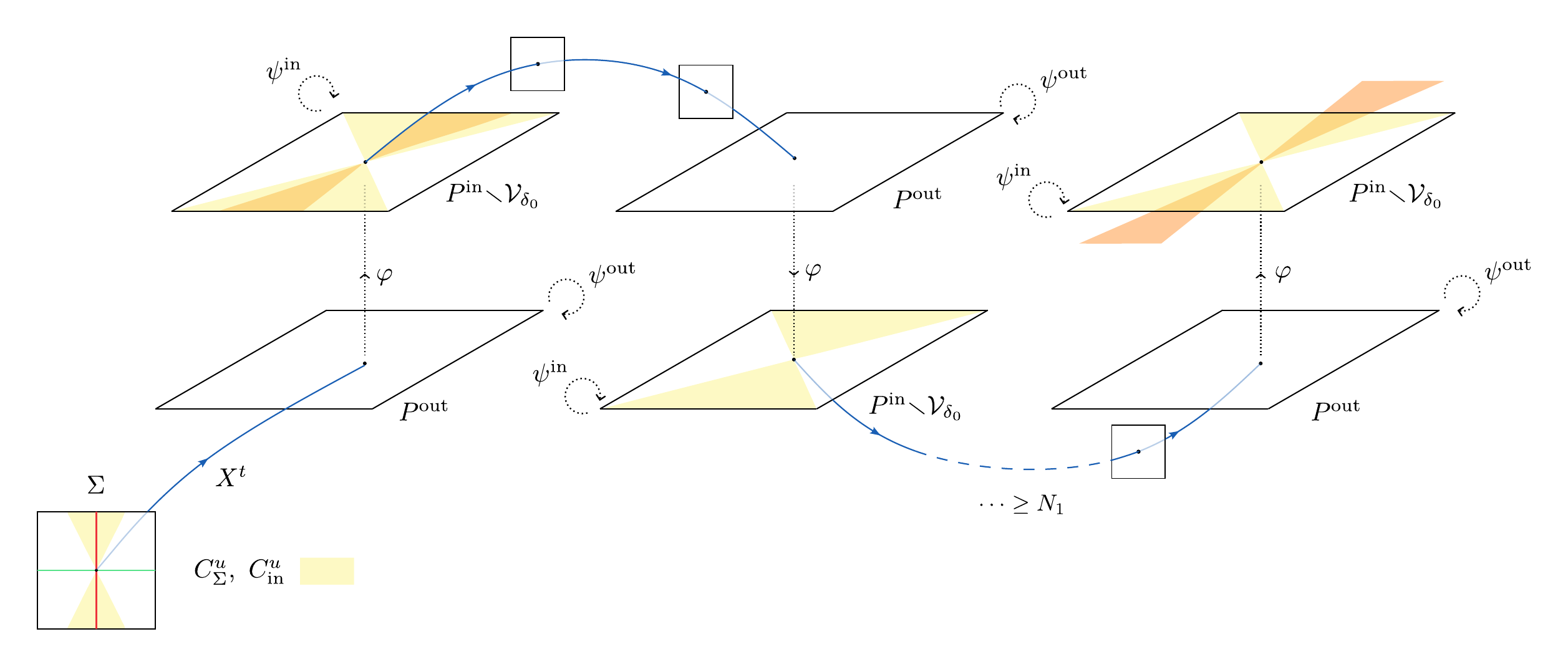}
    \vspace*{-2em}
    \caption{Trajectory of a cone $C^u_\Sigma$ under the map $\hat f_{\iin, \psi}^{j} \varphi \psi^\out f_{\out, \Sigma} f_\Sigma^i$}
    \label{fig: Sigma to Pin minus V0}
\end{figure}

\begin{lem} \label{lem: N1 for Sigma to Pin minus V0}
There exists an integer $N_1 \geq 0$ such that for all
$p \in \Sigma$, if the image $\hat f_{\iin, \psi}^j \varphi \psi^\out f_{\out, \Sigma} f_\Sigma^i(p) \in \Pin \ssm \cV_{\delta_0}$ is well defined for $j \geq N_1$ and $i \geq 0$, in other words $f_{\iin, \psi}^k \varphi \psi^\out f_{\out, \Sigma} f_\Sigma^i(p) \in \Pin \ssm \cV_\delta$ for all $k=0, \dots, j$,
then
\begin{itemize}[--]
    \item $\hat f_{\iin, \psi}^j \, \varphi \psi^\out \, f_{\out, \Sigma} \, f_\Sigma^i$ maps $C^u_\Sigma(p)$ inside $C^u_\iin(q)$ and expands the norm of the vectors of $C^u_\Sigma(p)$ by a factor $2$;
    \item $(\hat f_{\iin, \psi}^j \, \varphi \psi^\out \, f_{\out, \Sigma} \, f_\Sigma^i)^\inv$ maps $C^s_\iin(p)$ inside $C^s_\Sigma(q)$ and expands the norm of the vectors of $C^s_\iin(p)$ by a factor $2$.
\end{itemize}
\end{lem}

Let us insist on the fact that each crossing of the orbit of such a point $p$ on the surface $\Pin$ takes place on the complementary of $\cV_{\delta_0}$ by assumption.

\begin{proof}
The composition is decomposed into three pieces.
\begin{itemize}
    \item $f_\Sigma^i$ maps $C^u_\Sigma$ on $C^u_\Sigma$ and the contraction of the norm of the vectors is greater than a uniform constant $c_0>0$ for $i \geq 0$ (Proposition~\ref{prop: cones Sigma}, Item~\ref{prop: cones Sigma; it: invariance sigma}).
    \item We know that $\l( \varphi \psi^\out \, f_{\out, \Sigma} \r)_* C^u_\Sigma \subset \intr C^u_\iin$ (Proposition~\ref{prop: cones Sigma}, Item~\ref{prop: cones Sigma; it: compatibilty pin}).
    Moreover, by hypotheses the codomain of $f_{\out, \Sigma}$ for our restriction is contained in $\Pout \ssm (\varphi \psi^\out)^\inv (\cV_{\delta_0})$, so uniformly far from the set $\cO_*$.
    Using Lemma \ref{lem: contraction sigma to boundary is distance to orbit}, we deduce that the differential of $f_{\out, \Sigma}$ is bounded, uniformly in $p$.
    The diffeomorphisms $\varphi$, $\psi^\out$ are compactly supported so the differential of $\varphi \, \psi^\out \, f_{\out, \Sigma}$ is bounded by a constant $c_1>0$ (uniform in $p$).
   \item The map $\hat f_{\iin, \psi} \colon \ssm \cV_{\delta_0} \to \Pin \ssm \cV_{\delta_0}$, which coincides with the restriction of $f_{\iin, \psi}$ to the set $(\Pin \ssm \cV_{\delta}) \cap f_{\iin, \psi}^\inv( \Pin \ssm \cV_\delta)$ maps the field $C^u_\iin$ inside $C^u_\iin$ and expands the norm of the vectors of $C^u_\iin$ by a factor $2$ (Proposition~\ref{prop: lambda epsilon cones Pin}).
   It is enough to choose $N_1 \geq -\log_2(c_0 c_1)$ so that the expansion of $\hat f_{\iin, \psi}$ compensates the contraction of the previous maps.
   This integer is uniform in~$p$.
\end{itemize}
The second item is shown in a similar way.
\end{proof}

\textbf{We fix an integer $N_1$ which satisfies Lemmas~\ref{lem: N1 for Pin minus v0 to Sigma} and~\ref{lem: N1 for Sigma to Pin minus V0}.}

\subsection{Proof of Proposition~\ref{prop: parameters and cones}}
\label{sec: parameters; subsec: proof}

We are now able to prove Proposition~\ref{prop: parameters and cones}.
We have fixed in the previous sections the parameters $\lambda= \lambda_0$, $\epsilon = \epsilon_0$ (Subsection~\ref{sec: parameters; subsec: lambda epsilon cone Pin}), $\delta=\delta_0$ (Subsection~\ref{sec: parameters; subsec: delta}),
a pair of cone fields $(C^u_\iin, C^s_\iin)$ on $\Pin$ (Subsection~\ref{sec: parameters; subsec: lambda epsilon cone Pin}) and a pair of cone fields $(C^u_\Sigma, C^s_\Sigma)$ on $\Sigma$ (Subsection~\ref{sec: parameters; subsec: cone Sigma and N_Sigma}),
and we defined integers $N_\Sigma$ (Subsection~\ref{sec: parameters; subsec: cone Sigma and N_Sigma}), $N_0$ (Subsection~\ref{sec: parameters; subsec: N0}), $N_1$ (Subsection~\ref{sec: parameters; subsec: N1}).
These choices have been made in the \say{right order}, which means that all the results from Subsections~\ref{sec: parameters; subsec: lambda epsilon cone Pin} to~\ref{sec: parameters; subsec: N1} are true for the parameters $\lambda_0, \epsilon_0, \delta_0$.
Denote $\psi^\iin = \psi^\iin_{\lambda_0, \epsilon_0, \delta_0}$ and $\psi^\out = \psi^\out_{\lambda_0, \epsilon_0, \delta_0}$.
We use the notations introduced in Section~\ref{sec: parameters and cones}.
Let $\psi = \psi^\iin \varphi \psi^\out$ be the modified gluing map, $P_\psi = P/\,\psi$ the closed manifold and $X_\psi$ the vector field induced by $X$ on $P_\psi$. The projection is $\pi_\psi \colon P \to P_\psi$. Define the section $S_{0, \psi} = \pi_\psi( \Pin \cup \Sigma )$, $f_{0, \psi}$ the first-return map of the flow of $X_\psi$ on $S_{0, \psi}$ and $\cV_0 = \psi^\iin (\cV_{\delta_0})$.

Define the following pair of cone fields on $S_{0, \psi}$:
$$(C^u_\psi, C^s_\psi) = \l( \pi_\psi( \psi^\iin_* (C^u_\iin) \sqcup C^u_\Sigma),\  \pi_\psi(\psi^\iin_* (C^s_\iin) \sqcup C^s_\Sigma) \r) $$
and the metric
$g_\psi = \pi_\psi(\psi^\iin_* g^\iin) \sqcup \pi_\psi(g^\Sigma)$ on $P_\psi$,
where $g^\Sigma$ and $g^\iin$ are restrictions of the metric fixed on $P$.
Recall that the metric $g$ is the metric fixed at the beginning of Section~\ref{sec: crossing map}, which coincides with $dx^2 + dy^2 + d\theta^2$ in each $(\cV_i, \xi_i = (x,y,\theta))$, and which is adapted to the hyperbolic splitting of the flow.
Let us note $\Vert \cdot \Vert_\psi, \Vert \cdot \Vert_\Sigma, \Vert \cdot \Vert_\iin$ the associated norms.
Recall that for practical reasons, we have studied lifts in $P$ \emph{conjugated by $\psi^\iin$} of the map $f^n_{0, \psi}$.
This explains that we compose the projection $\pi_\psi$ by $\psi^\iin$ on $\Pin$ to fit.
Recall that:
\begin{itemize}
    \item $S_{0, \psi}$ has points of self-intersections $p$, at these points it has two or three tangent planes $T^* = T^*_p \subset T_p S_{0, \psi}$.
    The differential of $f_{0, \psi}$ in $p$ is well defined for a choice of tangent plane on $p$ and on $f_{0, \psi}(p)$.
    A cone field $C$ on $S_{0, \psi}$ is a collection of cones $C(p, T^i)$ on each of the tangent planes $T^i \subset T_p S_{0, \psi}$.
    \item If $p \in \Pin_\psi \subset S_{0, \psi}$, we denote $T^\iin_p \subset T_p S_{0, \psi}$ the tangent plane induced by $T_{\tilde p} P^\iin$, where $\tilde p \in \Pin$ is the lift of $p$ in $\Pin$. We denote $(C^u_\psi(p, T^\iin_p), C^s_\psi(p, T^\iin_p))$ the pair $(C^u_\psi, C^s_\psi)$ on $T^\iin_p$ in $p$.
    
  \item  If $q \in \Sigma_\psi \subset S_{0, \psi}$, we denote $T^\Sigma_\psi \subset T_p S_{0, \psi}$ the (or one of the two) tangent plane(s) induced by $T_{\tilde q} P^\Sigma$, where $\tilde q \in \Sigma$ is a lift of $q$ in $\Sigma$, and we denote $(C^u_\psi(q, T^\Sigma_q), C^s_\psi(q, T^\Sigma_q)$ the pair $(C^u_\psi, C^s_\psi)$ on $T^\Sigma_q$ in $q$.
\end{itemize}

\subsubsection*{Lifts in \texorpdfstring{$P$}{P}}
We summarize the lifts in $P$ of the points, cones, vectors and norms that we will use.

\begin{claim} \label{claim: lifts in P of points}
\mb{}
\begin{enumerate}
    \item \emph{(maps)} \label{claim: lifts; it: maps}
    \begin{itemize}[--]
    \item If $p \in \Sigma_\psi$ and $f_{0, \psi}(p) \in \Sigma_\psi$, then a lift of $f_{0, \psi}$ is $f_\Sigma \colon \Sigma \to \Sigma$.
    \item If $p \in \Sigma_\psi$ and $f_{0, \psi}(p) \in \Pin_\psi$, then a lift of $f_{0, \psi}$ is $\psi \, f_{\out, \Sigma}: \Sigma \to \Pin$.
    \item If $p \in \Pin_\psi$ and $f_{0, \psi}(p) \in \Sigma_\psi$, then a lift of $f_{0, \psi}$ is $f_{\Sigma, \iin}: \Pin \to \Sigma$.
    \item If $p \in \Pin_\psi$ and $f_{0, \psi}(p) \in \Pin_\psi$, then a lift of $f_{0, \psi}$ is $\psi \, \foutin: \Pin \to \Pin$
\end{itemize}
    \item \emph{(points)} \label{claim: lifts; it: points}
    \begin{itemize}[--]
        \item If $p \in \Pin_\psi$, there exists a unique lift $\tilde p$ of $p$ in $\Pin$, then denote $\check p := (\psi^\iin)^\inv (\tilde p)$, in other words we have $\pi_\psi (\psi^\iin (\check p)) =p$.
        Moreover if $p \in (\Pin \ssm \cV_0)_\psi$ then $\check p \in \Pin \ssm \cV_{\delta_0}$ by definition of $\cV_0 = \psi^\iin (\cV_{\delta_0})$.
        \item If $q \in \Sigma_\psi$, there exists $\tilde q \in \Sigma$ a lift of $q$ in $P$.
    \end{itemize}
    \item \emph{(cones)} \label{claim: lifts; it: cones}
    \begin{itemize}[--]
        \item If $p \in \Pin_\psi$, a lift of the cone $C^u_\psi(p, T^\iin_p)$ is the cone $\psi^\iin_* C^u_\iin ( \check p )$.
        \item If $q \in \Sigma$, a lift of the cone $C^u_\psi(q, T^\Sigma_q)$ is the cone $C^u_\Sigma (\tilde q)$.
    \end{itemize}
    \item \emph{(vectors)} \label{claim: lifts; it: vectors}
    \begin{itemize}[--]
        \item If $v \in T \Pin_\psi$, there exists a unique lift of $v$ in $TP$ of the form $\psi^\iin_* \check v$ with $\check v \in T \Pin$.
        Moreover, if $v \in C^u_\psi(p, T^\iin_p)$, then $\check v \in C^u_\iin(\check p)$.
        \item If $w \in T \Sigma_\psi$, there exists $\tilde w \in T \Sigma$ a lift of $w$ in $P$.
        Moreover, if $w \in C^u_\psi(q, T^\Sigma_q)$, then $\tilde w \in C^u_\Sigma(\tilde q)$.
    \end{itemize}
    \item \emph{(norm)} \label{claim: lifts; it: norm}
    \begin{itemize}[--]
        \item If $v \in T \Pin_\psi$, then $\Vert v \Vert_\psi = \Vert \check v \Vert_\iin = \Vert \check v \Vert$ where $\check v$ is given by the previous item and $\Vert \cdot \Vert$ is the norm from the Riemannian metric on $P$.
        \item If $w \in T \Sigma_\psi$, then $\Vert w \Vert_\psi = \Vert \tilde w \Vert_\Sigma = \Vert \tilde w \Vert$ where $\tilde w$ is given by the previous item and $\Vert \cdot \Vert$ is the norm from the Riemannian metric on~$P$.
    \end{itemize}
\end{enumerate}
\end{claim}
All the above items are direct consequences of the definitions of the cones and the metric on $P_\psi$.
Let us prove Proposition~\ref{prop: parameters and cones}.

\begin{proof}[Proof of Proposition~\ref{prop: parameters and cones}]
Let $p \in S_{0, \psi}$, $n \geq 1$ and $q=f_{0, \psi}^n(p)$.
Let us show each item of the proposition.
We only show the properties for the unstable cones $C^u_\psi$.
The proofs for the stable cones $C^s_\psi$ are symmetric.

\begin{enumerate}[leftmargin=*]
    \item[\ref{prop: parameters; it: directions}.] 
    Suppose $p \in \Pin_\psi$ and let $\tilde p := \psi^\iin \check p$ be a lift of $p$ in $P$ with $\check p \in \Pin$ (Claim~\ref{claim: lifts in P of points}, Item~\ref{claim: lifts; it: points}). Then
    \begin{align*}
        \intr C^u_\psi(p, T_p^\iin) & = (\pi_\psi)_* \psi^\iin_*( \intr C^u_\iin(\check p))\\
        & \supset (\pi_\psi)_* \psi^\iin_* (\varphi \psi^\out)_* \cG^{u, \out}
        & \mb{(Proposition~\ref{prop: lambda epsilon cones Pin}, Item~\ref{prop: cones Pin; it: directions})} \\
        & =  (\pi_\psi)_* \psi_* \cG^{u, \out} = (\pi_\psi)_* \cG^{u, \out} \\
        & \supset (\pi_\psi)_* \cL^\out = \cW^u_\psi \cap \Pin_\psi,
    \end{align*}
    \begin{align*}
    \intr C^u_\psi(p, T_p^\iin) \cap \cW^s_\psi & =  (\pi_\psi)_* \psi^\iin_*( \intr C^u_\iin(\check p)) \cap (\pi_\psi)_*\cW^s \\
    & \subset (\pi_\psi)_* \psi^\iin_*( \intr C^u_\iin(\check p)) \cap (\pi_\psi)_* \cG^{s, \iin} \\
    & = (\pi_\psi \circ \psi^\iin)_* \l( \intr C^u_\iin(\check p) \cap (\psi^\iin)^\inv_*\cG^{s, \iin} \r) \\
    & = (\pi_\psi \circ \psi^\iin)_* \{O\} = \{O\}. \qquad\quad \mb{(Proposition~\ref{prop: lambda epsilon cones Pin}, Item~\ref{prop: cones Pin; it: directions})}
\end{align*}
    Suppose now $p \in \Sigma_\psi$ and let $\tilde p :=\check p$ be a lift of $p$ in $P$ with $\check p \in \Sigma$  (Claim~\ref{claim: lifts in P of points}, Item~\ref{claim: lifts; it: points}). Then
    \begin{align*}
    \intr C^u_\psi(p, T_p^\Sigma) & = (\pi_\psi)_* (C^u_\Sigma(\check p))\\
    & \supset (\pi_\psi)_* \zeta^u
    & \mb{(Proposition~\ref{prop: cones Sigma}, Item~\ref{prop: cones Sigma; it: directions})} \\
    & \supset (\pi_\psi)_* (\cW^u \cap \Sigma) = \cW^u_\psi \cap \Sigma_\psi,
\end{align*}
\begin{align*}
    \intr C^u_\psi(p, T_p^\Sigma) \cap \cW^s_\psi & =  (\pi_\psi)_* (\pi_\psi)_* (C^u_\Sigma(\check p)) \cap (\pi_\psi)_*\cW^s \\
    & \subset (\pi_\psi)_* (C^u_\Sigma(\check p)) \cap (\pi_\psi)_* \zeta^s \\
    & = (\pi_\psi)_* \l( C^u_\Sigma(\check p) \cap \zeta^s \r)\\
    & = (\pi_\psi)_* \{O\} = \{O\}. \qquad\qquad \mb{(Proposition~\ref{prop: cones Sigma}, Item~\ref{prop: cones Sigma; it: directions})}
\end{align*}

    \item[\ref{prop: parameters; it: v0}.]
    Suppose that $p \in (\Pin \cap \cV_0)_\psi$.
    Let $\tilde p := \psi^\iin \check p$ be a lift of $p$ in $P$ with $\check p \in \Pin \cap \cV_{\delta_0}$ (Claim~\ref{claim: lifts in P of points}, Item~\ref{claim: lifts; it: points}).
    According to Proposition~\ref{prop: delta_0}, the maps
    $f_\Sigma^{N_0} , f_{\Sigma, \iin} \psi^\iin \colon \Pin \to \Sigma$ and 
    $f_\Sigma^{-N_0} f_{\out, \Sigma}^\inv \, (\varphi \psi^\out)^\inv \colon \Pin \to \Sigma$ are well defined in $\check p$.
    From Claim~\ref{claim: lifts in P of points}, Item~\ref{claim: lifts; it: maps}, the lifts of the iterates $f^{k}_{0, \psi}(p)$ in $P$ are the iterates:
    \begin{itemize}[--]
        \item $ f_\Sigma^{k-1} \, f_{\Sigma, \iin} (\tilde p) = f_\Sigma^{k-1} \, f_{\Sigma, \iin} \psi^\iin (\check p) $ for $k = 1, \dots, N_0$, and
        \item $f_\Sigma^{k-1} \, f_{\out, \Sigma}^\inv \, \psi^\inv (\tilde p) =f_\Sigma^{k-1} \, f_{\out, \Sigma}^\inv \, (\varphi \psi^\out)^\inv (\check p) $ for $k=-N_0, \dots, -1$.
    \end{itemize}
    We deduce that the iterates $f^k_{0, \psi}(p)$ belong to $\Sigma_\psi$ for $k= -N_0, \dots, -1, 1, \dots, N_0$.
 
    \item[\ref{prop: parameters; it: Pin minus v0 to sigma}.]
    Suppose $p \in (\Pin \ssm \cV_0)_\psi$, $q \in \Sigma_\psi$, $n \geq N_1$, and $f^k_{0,\psi}(p) \in \Sigma_\psi$ for $k=1, \dots, n-1$.
    
    Let $\tilde p := \psi^\iin \check p$ be a lift of $p$ in $P$ with $\check p \in \Pin \cap \cV_{\delta_0}$ (Claim~\ref{claim: lifts in P of points}, Item~\ref{claim: lifts; it: points}).
    According to Claim~\ref{claim: lifts in P of points}, Item~\ref{claim: lifts; it: maps}, $f^{n-1}_\Sigma, f_{\Sigma, \iin}$ is a lift of $f^n_{0, \psi}$ in the neighborhood of $p$.
    Let $\tilde q := f^{n-1}_\Sigma \, f_{\Sigma, \iin} \psi^\iin (\check p) \in \Sigma $.
    This is a lift of $q = f^n_{0, \psi}(p)$ in~$P$.
    
    According to Item~\ref{claim: lifts; it: cones}, Claim~\ref{claim: lifts in P of points}, a lift of $C^u_\psi(p, T^\iin_p)$ is $\psi^\iin_* C^u_\iin (\check p)$ and a lift of $C^u_\psi(q, T^\Sigma_q)$ is $C^u_\Sigma(\tilde q)$
    We use Lemma~\ref{lem: N1 for Pin minus v0 to Sigma}.\pagebreak[4]
    We have
\begin{align*}
    (f_{0, \psi}^n)_* C^u_\psi(p, T^\iin_p)
    & = (\pi_\psi)_* (f^{n-1}_\Sigma, f_{\Sigma, \iin} \, \psi^\iin)_* C^u_\iin (\check p) \\
    & \subset (\pi_\psi)_* \, \intr C^u_\Sigma(\tilde q)
    & \mb{(Lemma~\ref{lem: N1 for Pin minus v0 to Sigma})} \\
    & = \intr C^u_\psi(q, T^\Sigma_q).
\end{align*}

Let $v \in C^u_\psi(p, T^\iin_p)$.
Let $\psi^\iin_* \check v$ be a lift of $v$ in $T \Pin$, with $v \in C^u_\iin(\check p)$ (Claim~\ref{claim: lifts in P of points} Item~\ref{claim: lifts; it: vectors}).
Then $(f^{n-1}_\Sigma f_{\Sigma, \iin} \psi^\iin)_* \check v$ is a lift of $(f_{0, \psi}^n)_* v$ in $T \Sigma$, and we have
\begin{align*}    
    \dfrac{\Vert (f_{0, \psi}^n)_* v \Vert_\psi}{\Vert v \Vert_\psi} & = \dfrac{\Vert (f^{n-1}_\Sigma \, f_{\Sigma, \iin} \,\psi^\iin)_* \check v \Vert}{\Vert \check v \Vert}
    & \mb{(Claim~\ref{claim: lifts in P of points}, Item~\ref{claim: lifts; it: norm})}\\
    & \geq  2. &  \mb{(Lemma~\ref{lem: N1 for Pin minus v0 to Sigma})} 
\end{align*}
    
    \item[\ref{prop: parameters; it: Pin to Pin}.]
    Suppose $p \in \Pin_\psi$, $q \in \Pin_\psi$ and $f^k_{0, \psi}(p) \in \Sigma_\psi$ for $k=1, \dots, n-1$. 
    
    Let $\tilde p := \psi^\iin \check p$ be a lift of $p$ in $P$ with $\check p \in \Pin$ (Claim~\ref{claim: lifts in P of points}, Item~\ref{claim: lifts; it: points}).
    According to Claim~\ref{claim: lifts in P of points}, Item~\ref{claim: lifts; it: maps}, $\psi \foutin$ is a lift of $f^n_{0, \psi}$ in the neighborhood of $p$.
    Let $\tilde q := \psi \foutin (\psi^\iin (\check p)) = \psi^\iin \, f_{\iin, \psi} (\check p) \in \Pin $ (Equation~\eqref{eq: f_in,psi}).
    This is a lift of $q = f^n_{0, \psi}(p)$ in $P$.
    Let $\check q := (\psi^\iin)^\inv \tilde q = f_{\iin, \psi} (\check p)$.
    From Item~\ref{claim: lifts; it: cones}, Claim~\ref{claim: lifts in P of points}, $\psi^\iin_* C^u_\iin (\check p)$ is a lift of $C^u_\psi(p, T^\iin_p)$ and $\psi^\iin_* C^u_\iin(\check q)$ is a lift of $C^u_\psi(q, T^\iin_q)$.
    We use Proposition~\ref{prop: delta_0}, Item~\ref{prop: delta_0; it: Pin invariance}.
    We have
\begin{align*}    
    (f_{0, \psi}^n)_* C^u_\psi(p, T^\iin_p) & = (\pi_\psi)_* (\psi \, \foutin \, \psi^\iin)_* C^u_\iin (\check p) \\
    & =  (\pi_\psi \psi^\iin)_* (f_{\iin, \psi})_* C^u_\iin (\check p) \\ 
    & \subset (\pi_\psi \psi^\iin)_* (\intr C^u_\iin(\check q)) 
    & \mb{(Proposition~\ref{prop: delta_0}, Item~\ref{prop: delta_0; it: Pin invariance})} \\
    & = (\pi_\psi)_* (\intr \psi^\iin_* C^u_\iin(\check q))) \\
    & = \intr C^u_\psi( q, T^\iin_q).
\end{align*}
Now suppose $p$ and $q$ are in $(\Pin \ssm \cV_0)_\psi$, and $q \in (\Pin \ssm \cV_0)_\psi$.
Then $\check p$ and $\check q$ are in $\Pin \cap \cV_{\delta_0}$ (Claim~\ref{claim: lifts in P of points}, Item~\ref{claim: lifts; it: points}).
Let $v \in C^u_\psi(p, T^\iin_p)$.
Let $\psi^\iin_* \check v$ be a lift of $v$ in $T \Pin$, with $v \in C^u_\iin(\check p)$ (Claim~\ref{claim: lifts in P of points} Item~\ref{claim: lifts; it: vectors}).
Then $\psi^\iin_* (f_{\iin, \psi})_* \check v$ is a lift of $(f_{0, \psi}^n)_* v$ in $T \Pin$, and we have
\begin{align*}    
    \dfrac{\Vert (f_{0, \psi}^n)_* v \Vert_\psi}{\Vert v \Vert_\psi} & = \dfrac{\Vert (f_{\iin, \psi})_* \check v \Vert}{\Vert \check v \Vert}
    & \mb{(Claim~\ref{claim: lifts in P of points}, Item~\ref{claim: lifts; it: norm})}\\
    & \geq  2. & \mb{(Proposition~\ref{prop: lambda epsilon cones Pin})}
\end{align*}

    \item[\ref{prop: parameters; it: Pin minus v0 to sigma through Pin inter v0}.]
    Suppose $p \in (\Pin \ssm \cV_0)_\psi$, $q \in \Sigma_\psi$, $n \geq N_0 +1$,
    there exists $1 \leq m \leq n-N_0$ such that $f^m_{0, \psi}(p)\in (\Pin \cap \cV_0)_\psi$ and $f^k_{0, \psi}(p) \in \Sigma_\psi$ for $k=1, \dots, m-1, m+1, \dots$,$n-1$.
    
    Let $\tilde p := \psi^\iin \check p$ be a lift of $p$ in $P$ with $\check p \in \Pin \cap \cV_{\delta_0}$ (Claim~\ref{claim: lifts in P of points}, Item~\ref{claim: lifts; it: points}).
    From Claim~\ref{claim: lifts; it: points}, Item~\ref{claim: lifts; it: maps}, $f_\Sigma^{n-m} \, f_{\Sigma, \iin} \, \psi \, f_{\out, \Sigma} \, f_\Sigma^{m-1} \, f_{\Sigma, \iin}$ is a lift of $f^n_{0, \psi}$ in the neighborhood of $p$.
    Let $\tilde q := f_\Sigma^{n-m} \, f_{\Sigma, \iin} \, \psi \, f_{\out, \Sigma} \, f_\Sigma^{m-1} \, f_{\Sigma, \iin} \psi^\iin (\check p) \in \Sigma $.
    This is a lift of $q = f^n_{0, \psi}(p)$ in $P$.
    
    According to Item~\ref{claim: lifts; it: cones}, Claim~\ref{claim: lifts in P of points}, a lift of $C^u_\psi(p, T^\iin_p)$ is $\psi^\iin_* C^u_\iin (\check p)$ and a lift of $C^u_\psi(q, T^\Sigma_q)$ is $C^u_\Sigma(\tilde q)$.
    We have $n-m \geq N_0$.
    By Proposition~\ref{prop: N0 for Pin minus V0 to Sigma through Pin inter VO}
\begin{align*}    
    (f_{0, \psi}^n)_* C^u_\psi(p, T^\iin_p)
    & = (\pi_\psi)_* (f_\Sigma^{n-m} \, f_{\Sigma, \iin} \, \psi \, f_{\out, \Sigma} \, f_\Sigma^{m-1} \, f_{\Sigma, \iin} \psi^\iin)_* C^u_\iin (\check p) \\
    & \subset (\pi_\psi)_* \, \intr C^u_\Sigma(\tilde q) \\
    & = \intr C^u_\psi(q, T^\Sigma_q).
\end{align*}

Let $v \in C^u_\psi(p, T^\iin_p)$.
Let $\psi^\iin_* \check v$ be a lift of $v$ in $T \Pin$, with $v \in C^u_\iin(\check p)$ (Claim~\ref{claim: lifts in P of points} Item~\ref{claim: lifts; it: vectors}).
Then $f_\Sigma^{n-m} \, f_{\Sigma, \iin} \, \psi \, f_{\out, \Sigma} \, f_\Sigma^{m-1} \, f_{\Sigma, \iin} \psi^\iin)_* \check v$ is a lift of $(f_{0, \psi}^n)_* v$ in $T \Sigma$, and we have
\begin{align*}   
    \dfrac{\Vert (f_{0, \psi}^n)_* v \Vert_\psi}{\Vert v \Vert_\psi} & = \dfrac{\Vert ( f_\Sigma^{n-m} \, f_{\Sigma, \iin} \, \psi \, f_{\out, \Sigma} \, f_\Sigma^{m-1} \, f_{\Sigma, \iin} \psi^\iin)_* \check v \Vert}{\Vert \check v \Vert}
    & \mb{(Item~\ref{claim: lifts; it: norm})}\\
    & \geq  2. & \mb{(Proposition~\ref{prop: N0 for Pin minus V0 to Sigma through Pin inter VO})}
\end{align*}

    \item[\ref{prop: parameters; it: sigma to sigma}.]
    Suppose $p \in \Sigma_\psi$, $q \in \Sigma_\psi$, and $f^k_{0, \psi}(p) \in \Sigma$ for $k=1, \dots, n-1$.
    Let $\tilde p \in \Sigma$ be a lift of $p$ in $P$.
    According to Claim~\ref{claim: lifts in P of points}, Item~\ref{claim: lifts; it: maps}, $f_\Sigma^n$ is a lift of $f^n_{0, \psi}$ in the neighborhood of $p$.
    Let $\tilde q := f_\Sigma^n (\tilde p)$.
    This is a lift of $q = f^n_{0, \psi}(p)$ in~$P$.
    
    According to Item~\ref{claim: lifts; it: cones}, Claim~\ref{claim: lifts in P of points}, a lift of $C^u_\psi(p, T^\Sigma_p)$ is $C^u_\Sigma (\tilde p)$ and a lift of $C^u_\psi(q, T^\Sigma_q)$ is $C^u_\Sigma(\tilde q)$.
    We use Proposition~\ref{prop: cones Sigma}.
    We have
\begin{align*}    
    (f_{0, \psi}^n)_* C^u_\psi(p, T^\Sigma_p) & = (\pi_\psi)_* (f_\Sigma^n)_* C^u_\Sigma (\tilde p) \\
    & \subset (\pi_\psi)_* (\intr C^u_\Sigma(\tilde q)) 
    & \mb{(Proposition~\ref{prop: cones Sigma})} \\
    & = \intr C^u_\psi( q, T^\Sigma_q).
\end{align*}
Moreover, if $n \geq N_\Sigma$:
Let $v \in C^u_\psi(p, T^\Sigma_p)$ and $\tilde v \in C^u_\Sigma(\tilde p)$ be a lift of $v$.
We have
\begin{align*}    
    \dfrac{\Vert (f_{0, \psi}^n)_* v \Vert_\psi}{\Vert v \Vert_\psi} & = \dfrac{\Vert (f_\Sigma^n)_* \check v \Vert}{\Vert \check v \Vert}
    & \mb{(Claim~\ref{claim: lifts in P of points}, Item~\ref{claim: lifts; it: norm})}\\
    & \geq  2. & \mb{(Proposition~\ref{prop: cones Sigma})} 
\end{align*}

    \item[\ref{prop: parameters; it: sigma to Pin}.]
    Suppose $p \in \Sigma_\psi$, $q \in \Pin_\psi$ and $f^k_{0, \psi}(p) \in \Sigma_\psi$ for $k=1, \dots, n-1$.
    Let $\tilde p \in \Sigma$ be a lift of $p$ in $P$.
    From Claim~\ref{claim: lifts in P of points}, Item~\ref{claim: lifts; it: maps}, $ \psi \, f_{\out, \Sigma} \, f_\Sigma^n$ is a lift of $f^n_{0, \psi}$ in the neighborhood of $p$.
    Let $\tilde q := \psi \, f_{\out, \Sigma} \, f_\Sigma^n \, (\tilde p) \in \Pin $.
    This is a lift of $q = f^n_{0, \psi}(p)$ in $P$.
    Denote $\check q := (\psi^\iin)^\inv \tilde q = \varphi \, \psi^\out \, f_{\out, \Sigma} \, f_\Sigma^n (\tilde p)$.
    
    According to Item~\ref{claim: lifts; it: cones}, Claim~\ref{claim: lifts in P of points}, a lift of $C^u_\psi(p, T^\Sigma_p)$ is $C^u_\Sigma (\tilde p)$ and a lift of $C^u_\psi(q, T^\iin_q)$ is $\psi^\iin_* C^u_\iin(\check q)$.
    We use Proposition~\ref{prop: cones Sigma},
\begin{align*}    
    (f_{0, \psi}^n)_* C^u_\psi(p, T^\Sigma_p) & = (\pi_\psi)_* (\psi \, f_{\out, \Sigma} \, f_\Sigma^n)_* C^u_\iin (\tilde p) \\
    & =  (\pi_\psi \psi^\iin)_* ( \varphi \, \psi^\out \, f_{\out, \Sigma} \, f_\Sigma^n )_* C^u_\Sigma (\tilde p) \\ 
    & \subset (\pi_\psi \psi^\iin)_* (\intr C^u_\iin(\check q)) 
    & \mb{(Proposition~\ref{prop: cones Sigma})} \\
    & = (\pi_\psi)_* \intr \psi^\iin_* C^u_\iin(\check q)) \\
    & = \intr C^u_\psi(q, T^\iin_q).
\end{align*}

    \item[\ref{prop: parameters; it: sigma to sigma through Pin}.]
    Suppose that $p \in \Sigma_\psi$, $q \in \Sigma_\psi$, $n \geq N_0+1$, and there exists $1 \leq m \leq n-N_0$ such that $f_{0, \psi}^m(p) \in \Pin_\psi$, and $f^k(p) \in \Sigma_\psi$ for $k=1, \dots, m-1, m+1, \dots, n-1$.

Let $\tilde p \in \Sigma$ be a lift of $p$ in $P$.
    From Claim~\ref{claim: lifts in P of points}, Item~\ref{claim: lifts; it: maps}, $f_\Sigma^{n-m} \, f_{\Sigma, \iin} \, \psi \, f_{\out, \Sigma} \, f_\Sigma^{m-1}$ is a lift of $f^n_{0, \psi}$ in the neighborhood of $p$.
    We set $\tilde q := f_\Sigma^{n-m} \, f_{\Sigma, \iin} \, \psi \, f_{\out, \Sigma} \, f_\Sigma^{m-1} (\tilde p)$.
    This is a lift of $q = f^n_{0, \psi}(p)$ in $P$.
    
    According to Item~\ref{claim: lifts; it: cones}, Claim~\ref{claim: lifts in P of points}, a lift of $C^u_\psi(p, T^\Sigma_p)$ is $C^u_\Sigma (\tilde p)$ and a lift of $C^u_\psi(q, T^\Sigma_q)$ is $C^u_\Sigma(\tilde q)$.
    Since $n-m \geq N_0$, we can use Proposition~\ref{prop: N0 for sigma to sigma through Pin}.
    We have
\begin{align*}       
    (f_{0, \psi}^n)_* C^u_\psi(p, T^\Sigma_p) & = (\pi_\psi)_* (f_\Sigma^{n-m} \, f_{\Sigma, \iin} \, \psi \, f_{\out, \Sigma} \, f_\Sigma^{m-1})_* C^u_\Sigma (\tilde p) \\
    & \subset (\pi_\psi)_* (\intr C^u_\Sigma(\tilde q)) 
    & \mb{(Proposition~\ref{prop: N0 for sigma to sigma through Pin})} \\
    & = \intr C^u_\psi( q, T^\Sigma_q).
\end{align*}

Let $v \in C^u_\psi(p, T^\Sigma_p)$ and $\tilde v \in C^u_\Sigma(\tilde p)$ be a lift of $v$.
We have
\begin{align*}       
    \dfrac{\Vert (f_{0, \psi}^n)_* v \Vert_\psi}{\Vert v \Vert_\psi} & = \dfrac{\Vert (f_\Sigma^{n-m} \, f_{\Sigma, \iin} \, \psi \, f_{\out, \Sigma} \, f_\Sigma^{m-1})_* \check v \Vert}{\Vert \check v \Vert}
    & \mb{(Claim~\ref{claim: lifts in P of points}, Item~\ref{claim: lifts; it: norm})}\\
    & \geq  2. & \mb{(Proposition~\ref{prop: N0 for sigma to sigma through Pin})}
\end{align*}

    \item[\ref{prop: parameters; it: sigma to Pin minus v0}.] Suppose
    $p \in \Sigma_\psi$, $q \in (\Pin \ssm \cV_0)_\psi$, and there exist $m \geq N_1$, and integers $1 \leq k_1 < \dots <k_m = n$ such that
    and $f_{0, \psi}^{k_1}(p), \dots, f_{0, \psi}^{k_m}(p) \in (\Pin \ssm \cV_0)_\psi$
    and $f_{0, \psi}^k(p) \in \Sigma$ for $k \neq k_1, \dots, k_m$, $1 \leq k \leq n-1$.
    Let $\tilde p \in \Sigma$ be a lift of $p$ in $P$.
   According to Claim~\ref{claim: lifts in P of points}, Item~\ref{claim: lifts; it: maps}, $(\foutin \psi)^{m-1} f_{\out, \Sigma} f_\Sigma^{k_1-1} $ is a lift of $f^n_{0, \psi}$ in the \nbh{} of $p$, and we write in the following decomposition:
     $$(\foutin \psi)^{m-1} f_{\out, \Sigma} f_\Sigma^{k_1 -1} = \psi^\iin \hat f_{\iin, \psi}^{m-1} \varphi \psi^\out f_{\out, \Sigma} f_\Sigma^{k_1 -1}  $$
    with $\hat f_{\iin, \psi} \colon \Pin \ssm \cV_{\delta_0} \to \Pin \ssm \cV_{\delta}$ is the restriction of $f_{\iin, \psi}$ to the corresponding sets.
    
    Let $\tilde q := \psi^\iin \hat f_{\iin, \psi}^{m-1} \varphi \psi^\out f_{\out, \Sigma} f_\Sigma^{k_1 -1} \in \Pin$.
    This is a lift of $q = f^n_{0, \psi}(p)$ in $P$.
    We set $\check q := (\psi^\iin)^\inv \tilde q = (\hat f_{\iin, \psi})^{m-1} \varphi \psi^\out f_{\out, \Sigma} f_\Sigma^{k_1 -1} (\tilde p)$.
    
    According to Item~\ref{claim: lifts; it: cones}, Claim~\ref{claim: lifts in P of points}, a lift of $C^u_\psi(p, T^\Sigma_p)$ is $C^u_\Sigma (\tilde p)$ and a lift of $C^u_\psi(q, T^\iin_q)$ is $\psi^\iin_* C^u_\iin(\check q)$.
    Using Lemma~\ref{lem: N1 for Sigma to Pin minus V0}, we have
     \begin{align*}       
    (f_{0, \psi}^n)_* C^u_\psi(p, T^\Sigma_p) & = (\pi_\psi)_* (\psi^\iin \hat f_{\iin, \psi}^{m-1} \varphi \psi^\out f_{\out, \Sigma} f_\Sigma^{k_1 -1})_* C^u_\iin (\tilde p) \\
    & =  (\pi_\psi \psi^\iin)_* (\hat f_{\iin, \psi}^{m-1} \varphi \psi^\out f_{\out, \Sigma} f_\Sigma^{k_1 -1})_* C^u_\Sigma (\tilde p) \\ 
    & \subset (\pi_\psi \psi^\iin)_* (\intr C^u_\iin(\check q)) 
    & \mb{(Lemma~\ref{lem: N1 for Sigma to Pin minus V0})} \\
    & = (\pi_\psi)_* \intr \psi^\iin_* C^u_\iin(\check q)) \\
    & = \intr C^u_\psi(q, T^\iin_q).
\end{align*}

Let $v \in C^u_\psi(p, T^\Sigma_p)$, and $\tilde v \in C^u_\Sigma(\tilde p)$ be a lift of $v$ in $T \Sigma$.
Then $$(\psi^\iin)_* (\hat f_{\iin, \psi}^{m-1} \varphi \psi^\out f_{\out, \Sigma} f_\Sigma^{k_1 -1})_* \tilde v$$ is a lift of $(f_{0, \psi}^n)_* v$ in $T \Pin$, and we have
\begin{align*}       
    \dfrac{\Vert (f_{0, \psi}^n)_* v \Vert_\psi}{\Vert v \Vert_\psi} & = \dfrac{\Vert (\hat f_{\iin, \psi}^{m-1} \varphi \psi^\out f_{\out, \Sigma} f_\Sigma^{k_1 -1})_* \tilde v \Vert}{\Vert \tilde v \Vert}
    & \mb{(Claim~\ref{claim: lifts in P of points}, Item~\ref{claim: lifts; it: norm})}\\
    & \geq  2. & \mb{(Lemma~\ref{lem: N1 for Sigma to Pin minus V0})}
\end{align*}
\end{enumerate}
This completes the proof of Proposition \ref{prop: parameters and cones}.
\end{proof}

\subsection{Additional result}
\label{sec: parameters; subsec: additional result}
We can improve Proposition~\ref{prop: parameters and cones} by adding (for the same parameters) the existence of another pair of cone fields on $S_{0, \psi}$, denoted $(\hat C^u_\psi, \hat C^s_\psi)$, and which will have a complementary role to the one of the pair $(C^u_\psi, C^s_\psi)$ in the next section. Let us explain.

In the next section, we will show (Proposition~\ref{prop: future iterate for f_psi on S_psi}) that we can remove a neighborhood of the boundary of $\Pin_\psi$ in $S_{0, \psi}$ to obtain a compact section $S_\psi$ with boundary, transverse to the vector field $X$, and catching the same orbit in uniformly bounded time.
The first return map $f_\psi$ which is a restriction of $f_{0, \psi}$ then satisfies the following property.
For \emph{any $p \in S_\psi$}, there exists an integer $n^\star = n^\star(p) \geq 0$, such that $f^n_\psi$ satisfies the cone fields condition in the neighborhood of $p$ for the pair $(C^u_\psi, C^s_\psi)$ (restricted on $S_\psi$).
Since the integer $n^\star = n^\star(p)$ depends on the point $p$, this result is not symmetric by inversion of the time direction, in other words it does not follow that there exists an integer $\hat n^\star = \hat n^\star (p) \leq 0$ such that $f^{\hat n^\star}_\psi$ satisfies the cone fields condition in the neighborhood of $p$ for the pair $(C^s_\psi, C^u_\psi)$.
Now this is necessary to prove the existence of a hyperbolic splitting on $S_\psi$ using the invariant cone fields criterion, which will then give the Anosov structure for $X_\psi$ on $P_\psi$ (Proposition~\ref{prop: hyperbolic splitting on S_psi}).
In fact, this is not the case for the pair $(C^s_\psi, C^u_\psi)$ given by Proposition~\ref{prop: parameters and cones}. Let us explain.

We chose the pair $(C^s_\psi, C^u_\psi)$ so as to satisfy Item~\ref{prop: parameters; it: sigma to Pin} of Proposition~\ref{prop: parameters and cones}, which gives a compatibility of cones from $\Sigma_\psi$ to $\Pin_\psi$ at the first crossing, and allows us to keep the control of the image of a cone during the trajectory of an orbit crossing $\Pin_\psi$ several times.
On the other hand, there is no compatibility of cones from $\Pin_\psi$ to $\Sigma_\psi$ on the first crossing.
To get around this problem, we use a pair $(\hat C^u_\psi, \hat C^s_\psi)$ which satisfies this property, i.e., such that we have strict invariance of the cones by $f_\psi$ at the first crossing from $\Pin_\psi$ to $\Sigma_\psi$.
This pair will also continue to satisfy all the items of Proposition~\ref{prop: parameters and cones}, except that of Item~\ref{prop: parameters; it: sigma to Pin}.

Let us state the result. Let the parameters $\lambda_0>1$, $\epsilon_0>0$, $\delta_0>0$ and the integers $N_0$, $N_1$ and $N_\Sigma$ be given by Proposition~\ref{prop: parameters and cones}.
Let
$\psi^\iin = \psi^\iin_{\lambda_0, \epsilon_0, \delta_0}$,
$\psi^\out = \psi^\out_{\lambda_0, \epsilon_0, \delta_0}$,
$\psi = \psi^\iin \varphi \psi^\out$,
the \diffs{} of $\pP$ associated with this choice of parameters.

\begin{prop} \label{prop sym: parameters and cones}
There exists a pair $(\hat C^u_\psi, \hat C^s_\psi)$ of cone fields on $S_{0, \psi}$ which satisfy the following properties.
Let $p \in S_{0, \psi}$, let $n \geq 1$ and let $q = f_{0, \psi}^n(p)$.
\begin{enumerate}
    \item All the items of Proposition~\ref{prop: parameters and cones} except Item~\ref{prop: parameters; it: sigma to Pin} are satisfied.
    
    \item Suppose $p \in \Pin_\psi$, $q \in \Sigma_\psi$ and $f^k_{0, \psi}(p) \in \Sigma_\psi$ for $k=1, \dots, n-1$.
    Then
    \begin{itemize}[--]
        \item $(f_{0, \psi}^n)_* \hat C^u_\psi(p, T_p^\iin) \subset \intr \hat C^u_\psi(q, T_q^\Sigma)$, and
        \item $(f_{0, \psi}^n)_*^\inv \hat C^s_\psi(q, T_q^\Sigma) \subset \intr \hat C^s_\psi(p, T_p^\iin)$.
    \qedhere
    \end{itemize}
\end{enumerate}
\end{prop}

The proof is the same, except for Section~\ref{sec: parameters; subsec: cone Sigma and N_Sigma}, Proposition~\ref{prop: cones Sigma} (this is the part which allows us to obtain Item~\ref{prop: parameters; it: sigma to Pin} of Proposition~\ref{prop: parameters and cones} in the conclusion).
The analogous of Proposition~\ref{prop: cones Sigma} is

\begin{prop} \label{prop sym: cones Sigma}
There exists a pair of cone fields $(\hat C^u_\Sigma, \hat C^s_\Sigma)$ on $\Sigma$
\begin{enumerate}
    \item \emph{(direction)} \label{prop sym: cones Sigma; it: directions} 
    There exist constant $\hat K^u>0$ and $\hat K^s>0$ such that\\[0.4em]
    \hspace*{2em}$\hat C^u_\Sigma$ is a $(\hat K^u,\ \zeta^u/\,\zeta^s)$-cone field, and $\hat C^s_\Sigma$ is a $(\hat K^s,\ \zeta^s/\,\zeta^u)$-cone field.

     \item \emph{(compatibility on $\Pin$)} \label{prop sym: cones Sigma; it: compatibilty pin} 
     For all $\delta>0$, if $\psi^\out = \psi^\out_\delta$, and for all $p \in \Sigma$ such that the point $q = \varphi \, \psi^\out \, f_{\out, \Sigma} \, f_\Sigma^n \, (p) \in \Pin$ is well defined for some integer $n \geq 0$, then
      \begin{itemize}[--]
     
     \item $(\varphi \, \psi^\out \, f_{\out, \Sigma} \, f_\Sigma^n \,)_* \hat C^u_\Sigma(p) \subset \intr C^u_\iin(q)$, and
     
     \item $(\varphi \, \psi^\out \, f_{\out, \Sigma} \, f_\Sigma^n \,)^\inv_* C^s_\iin(q) \subset \intr \hat C^s_\Sigma(p)$.
      \end{itemize}

     \item \emph{(invariance and expansion on $\Sigma$)} \label{prop sym: cones Sigma; it: invariance sigma}  
     There exists an integer $N_\Sigma \geq 0$ such that
      \begin{itemize}[--]
     \item $(f_\Sigma)_* \hat C^u_\Sigma \subset \intr \hat C^u_\Sigma, \ \mb {and} \ \forall v \in \hat C^u_\Sigma$, $\forall n \geq N_\Sigma$, $\Vert (f_\Sigma^n)_* v \Vert \geq 2 \Vert v \Vert$,
     
     \item $(f_\Sigma)^\inv_* \hat C^s_\Sigma \subset \intr \hat C^s_\Sigma, \ \mb {and} \ \forall v \in \hat C^s_\Sigma$, $\forall n \geq N_\Sigma$, $\Vert (f_\Sigma^n)^\inv_* v \Vert \geq 2 \Vert v \Vert$.
     \end{itemize}
\end{enumerate}
\end{prop}

The proof is based on the same arguments.
We have to choose a slope $\hat K^u$ of $\hat C^u_\Sigma$ big enough and a slope $\hat K^s$ of $\hat C^s_\Sigma$ small enough to satisfy the compatibility of Item~\ref{prop sym: cones Sigma; it: compatibilty pin}. Essentially, the slope $\hat K^u$ is of the order of the slope $K^s$ and conversely the slope $\hat K^s$ is of the order of the slope $K^u$, hence the integer $N_\Sigma$ is indeed the same in both proposition.
The cone fields $(\hat C^u_\psi, \hat C^s_\psi)$ will be the quotient by $\pi_\psi$ of cone fields on $\Pin$ and on $\Sigma$ in $P$.
The cone fields on $\Pin$ are the same as for the pair $(C^u_\psi, C^s_\psi)$, namely $(\psi^\iin_* C^u_\iin, \psi^\iin_* C^s_\iin)$ where $C^u_\iin$ and $C^s_\iin$ are the cone fields given in Proposition~\ref{prop: lambda epsilon cones Pin}.
We replace the cones $(C^u_\Sigma, C^s_\Sigma)$ by the cones $(\hat C^u_\Sigma, \hat C^s_\Sigma)$ on $\Sigma$.
By reproducing point by point the previous proof and replacing the cone fields $(C^u_\Sigma, C^s_\Sigma)$ by the cone fields $(\hat C^u_\Sigma, \hat C^s_\Sigma)$ we obtain the result announced in Proposition~\ref{prop sym: parameters and cones}.

\section{Modification of the gluing map and proof of Theorem~\ref{thmintro: gluing theorem}}
\label{sec: proof gluing thm}

In Section~\ref{sec: parameters and cones}, we showed that there exist parameters $\lambda_0, \epsilon_0, \delta_0$ which satisfy Proposition~\ref{prop: parameters and cones}.
We refer to the notations introduced in Section~\ref{sec: parameters and cones}.
Let $\psi = \psi^\iin \varphi \psi^\out$ be the modified gluing map for such parameters and $X_\psi$ be the vector field induced by $X$ on $P_\psi = P/\,\psi$. 
Recall that $S_{0, \psi} = \pi_\psi( \Pin \cup \Sigma )$ and $f_{0, \psi}$ is the first-return map of the flow of $X_\psi$ on $S_{0, \psi}$.
According to Proposition~\ref{prop: parameters and cones}, there exists a pair of cone fields
$(C^u_\psi, C^s_\psi)$ on $S_{0, \psi}$ for which a (non-uniform) iterate of $f_{0, \psi}$ satisfies the cone fields condition in the neighborhood of some points, \emph{outside a small neighborhood $(\cV_0)_\psi$ of the periodic orbits $(\cO_*)_\psi$ in $S_{0, \psi}$.}
Indeed, the section is not uniformly transverse and the expansion of the vector of the cones is missing along some orbits when approaching the periodic orbits $(\cO_*)_\psi$ contained in $\Pin_\psi$ (Item~\ref{prop: parameters; it: sigma to Pin} and Item~\ref{prop: parameters; it: Pin to Pin}).
The goal of this section is therefore to remove a well-chosen neighborhood of the orbits $(\cO_*)_\psi$ in $S_{0, \psi}$ in order to have the hyperbolicity of the return map everywhere while catching the same set of orbit.

\subsection{Gluing procedure and invariant cone fields on a global section}
\label{sec: proof; subsec: invariant cone fields global section}

\subsubsection*{Compact local section}
We remove from $S_{0, \psi}$ the projection of the neighborhood $\cV_0 := \psi^\iin (\cV_{\delta_0})$ of the orbits $\cO_*$. We obtain the surface
\begin{equation} \label{eq: S_psi}
    S_\psi := \pi_\psi \l( (\Pin \ssm \cV_0 ) \cup \Sigma \r) = (\Pin \ssm \cV_0)_\psi \cup \Sigma_\psi \subset S_{0, \psi}.
\end{equation}
It is an immersed compact section catching the same set of $X_\psi$-orbit.
Indeed, according to Item~\ref{prop: parameters; it: v0} of Proposition~\ref{prop: parameters and cones}, an orbit of $X_\psi$ which intersects $(\Pin \cap \cV_0)_\psi$ intersects the section $\Sigma_\psi \subset S_\psi$ in the past and in the future in uniformly bounded time.

\begin{figure}[htb]
    \centering
    \vspace*{-2em}
    \includegraphics[width=\textwidth]{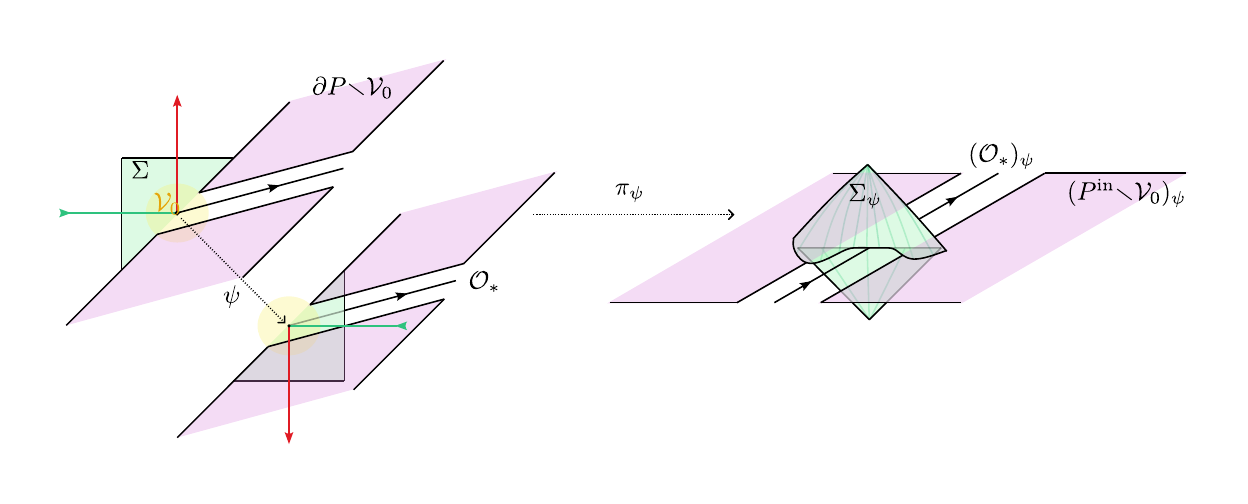}
    \vspace*{-3em}
    \caption{Section $S_\psi$}
    \label{fig: S_psi}
\end{figure}
We recall (see Subsection~\ref{sec: parameters; subsec: statement proposition}) that as for the section $S_{0, \psi}$, the surface $S_\psi$ admits points $p$ of self-intersection on which it can have two or three tangent planes $T^* = T^*_p \subset T_p S_\psi$.
Denote $T^\Sigma_p$ the tangent plane(s) to $\Sigma_\psi$ at $p$ and $T^\iin_q$ the tangent plane to $\Pin_\psi$ at $q$.
The surface $S_\psi$ has a metric $g_\psi$, which is a collection of metrics $g_\psi (p, T^*)$ on each of the tangent planes $T^* \subset T_p S_\psi$.
A cone field $C$ on $S_\psi$ is the data for each $p \in S_\psi$ and for each tangent plane $T^*$ to $S_\psi$ at the point $p$ of a cone $C(p, T^*) \subset T^*$.

\subsubsection*{Statement and proof of the proposition}
If $f_\psi$ is the first return map of the flow of $X_\psi$ on $S_\psi$,
then we show that it has a hyperbolic behavior in the \emph{future} of the orbit of the points $p \in \Sigma_\psi$.
We will deal separately with the maximal \emph{positive} invariant set in this section and the maximal \emph{negative} invariant set in the next subsection. More precisely

\begin{defi}[Positive, Negative, Maximal invariant set]
\label{def: positive, negative, maximal invariant set f_psi S_psi}
    Define 
    \begin{itemize}
        \item $K_\psi^+ := \bigcap_{n\geq 0} f_\psi^n(S_\psi)$ the \emph{maximal positive invariant set for $f_\psi$ in $S_\psi$},
        \item $K_\psi^- := \bigcap_{n\geq 0} f_\psi^n(S_\psi)$ the \emph{maximal negative invariant set for $f_\psi$ in $S_\psi$},
        \item $K_\psi := \bigcap_{n\in \Z} f_\psi^n(S_\psi)$ the \emph{maximal invariant set for $f_\psi$ in $S_\psi$}.
    \end{itemize}
\end{defi}

The following proposition gives the existence of a \emph{positive} integer $n^\star(p)\geq 0$ for any point $p \in K_\psi^+$, uniformly bounded from above, such that the iterate $f^{n^\star}_\psi$ satisfies the cone fields condition in the neighborhood of $p$ for the pair $(C^s_\psi, C^u_\psi)$ of Proposition~\ref{prop: parameters and cones}.

\begin{prop} \label{prop: future iterate for f_psi on S_psi}
There exists an integer $\bN \geq 0$ such that,
for all $p \in K_\psi^+$, if $T_p$ is (one of) the tangent plane(s) to $S_\psi$ at $p$,
there exists an integer $n^\star = n^\star(p, T_p)$ with $0 \leq n^\star \leq \bN$,
and a tangent plane $T_q = T_q(p, T_p)$ to $S_\psi$ in
$q = f_\psi^{n^\star}(p) \in S_\psi$
such that
\begin{enumerate}
    \item $(f_\psi^{n^\star})_* C^u_\psi(p, T_p) \subset \intr C^u_\psi(q, T_q)$,
    and for all $v \in C^u_\psi(p, T_p)$, $$\Vert (f_\psi^{n^\star})_* v \Vert \geq 2 \Vert v \Vert,$$
    \item $(f_\psi^{n^\star})^\inv_* C^s_\psi(q, T_q) \subset \intr C^s_\psi(p, T_p)$,
    and for all $v \in C^s_\psi(q, T_q)$, $$\Vert (f_\psi^{n^\star})^\inv_* v \Vert \geq 2 \Vert v \Vert.$$
\end{enumerate}
\end{prop}

\begin{proof}

Note that $f_\psi$ is a restriction of $f_{0, \psi}$.
We can use exclusively the results of Proposition~\ref{prop: parameters and cones} on the return map $f_{0, \psi}$ of the flow of $X_\psi$ on the section $S_{0, \psi}$.
If $p \in K_\psi^+$, then there exist $n$ and $m$ such that $f^n_\psi(p) = q = f^m_{0, \psi}(p)$ and $n \leq m$.
Indeed $f_\psi$ does not \say{see} the crossings on $(\Pin \cap \cV_0)_\psi$.
We will make dichotomy according to the starting point $p$ on the set $K_\psi^+$, and the trajectory of its orbit.
We represent the dichotomy by the two diagrams of Figure~\ref{fig: graph orbit}.

\begin{figure}[htb]
    \centering
    \begin{subfigure}{0.92\textwidth}
    \vspace*{-1em}
    \includegraphics[width=\textwidth]{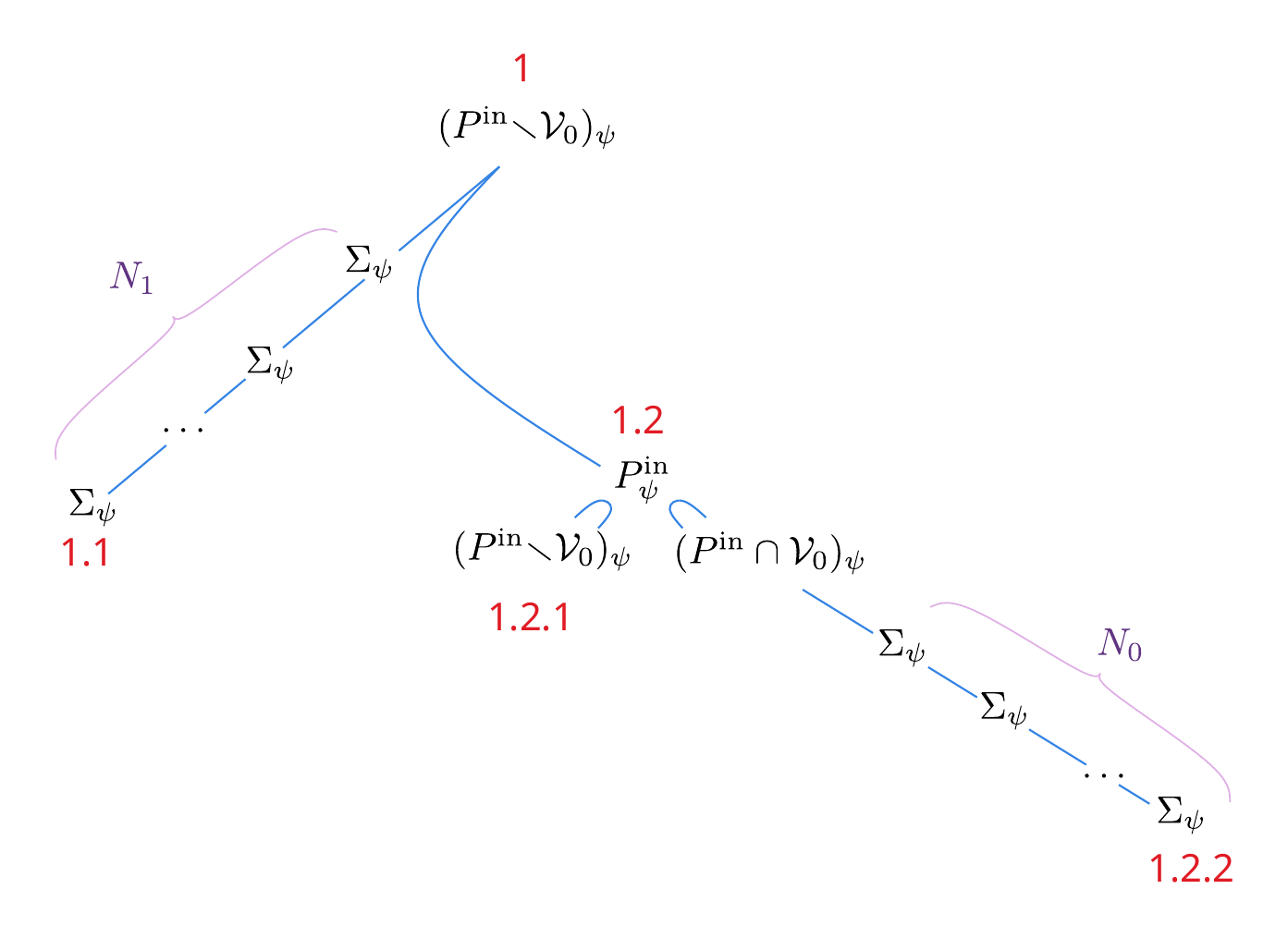}
    \vspace*{-1em}
    \caption{Dichotomy diagram of orbits in Case~\ref{it: start Pin}}
    \label{fig: graph orbit start Pin}
    \end{subfigure}
    \caption{Dichotomy diagram of $f_\psi$-orbits on $S_\psi$}
    \label{fig: graph orbit}
\end{figure}

\begin{enumerate} [label=\arabic*., leftmargin=*]

    \item \textbf{Suppose $p \in (\Pin \ssm \cV_0)_\psi$ and $T_p = T_p^\iin$.}
    \label{it: start Pin}
    Two cases are possible depending on whether the future orbit of $p$ by the flow of $X_\psi$ crosses the section $\Sigma_\psi$ enough times without crossing $\Pin_\psi$ to have strict invariance and expansion of the cones, or whether it crosses $\Pin_\psi$ before.
    This return number is the integer $N_1$ which satisfies Item~\ref{prop: parameters; it: Pin minus v0 to sigma} of Proposition~\ref{prop: parameters and cones}.
 We have the following dichotomy:
    
    \begin{enumerate} [label=\theenumi\arabic*.,  leftmargin=*]
        \item \textbf{The orbit of $p$ by the flow of $X_\psi$ crosses  $\Sigma_\psi$ $N_1$ consecutive times without crossing $\Pin_\psi$.}
        In other words $f_{0, \psi}^k(p) \in \Sigma_\psi$ for $k=1, \dots, N_1$.
        Then according to Item~\ref{prop: parameters; it: Pin minus v0 to sigma} Proposition~\ref{prop: parameters and cones},
        Proposition~\ref{prop: future iterate for f_psi on S_psi} is satisfied for $n^\star(p)=N_1$ and $T_q = T_q^\Sigma$.

        \item \textbf{The orbit of $p$ by the flow of $X_\psi$ crosses $\Pin_\psi$ before crossing $\Sigma_\psi$ $N_1$ times.}
        In other words there exists $k_0 \leq N_1$ minimal such that $f_{0, \psi}^{k_0}(p) \in \Pin_\psi$.
        Then two cases are possible:
         \begin{enumerate} [label=\theenumii\arabic*., wide=0pt, itemsep=1pt]
            \item ${\boldsymbol{f_{0, \psi}^{k_0}(p) \in (\Pin \ssm \cV_0)_\psi}}$.
            According to Item~\ref{prop: parameters; it: Pin to Pin} of Proposition~\ref{prop: parameters and cones}, Pro\-position~\ref{prop: future iterate for f_psi on S_psi} is satisfied for $n^\star(p) = k_0 \leq N_1$ and $T_q = T_q^\iin$.
            
            \item $\boldsymbol{f_{0, \psi}^{k_0}(p) \in (\Pin \cap \cV_0)_\psi}$.
            Then according to Item~\ref{prop: parameters; it: v0} of Proposition~\ref{prop: parameters and cones}, the future orbit of $f_{0, \psi}^{k_0}(p)$ by the flow of $X_\psi$ crosses $N_0$ consecutive times the section $\Sigma_\psi$ without crossing $\Pin_\psi$.
            In other words $f^K_\psi\l(f^{k_0}_{0, \psi}(p)\r) \in \Sigma_\psi$ for $k=1, \dots, N_0$.
            Then according to Item~\ref{prop: parameters; it: Pin minus v0 to sigma through Pin inter v0} of Proposition~\ref{prop: parameters and cones}, Proposition~\ref{prop: future iterate for f_psi on S_psi} is satisfied for $n^\star(p) = N_0 + k_0-1 \leq N_0 + N_1 -1$ and $T_q = T_q^\Sigma$.
            Note that the endpoint
            $q$ on $\Sigma_\psi$ corresponds to the iteration $N_0 + k_0$ of the orbit of $p$ by $f_{0, \psi}$, and the iteration $N_0 + k_0 -1$ of the orbit of $p$ by $f_\psi$ because in the second case we do not count the crossing of the orbit through $(\Pin \cap \cV_0)_\psi$.
         \end{enumerate}
    \end{enumerate}
    
     \begin{figure}[htb]
    \ContinuedFloat
    \centering
    \begin{subfigure}{0.93\textwidth}
     \vspace*{-1em}
    \includegraphics[width=\textwidth]{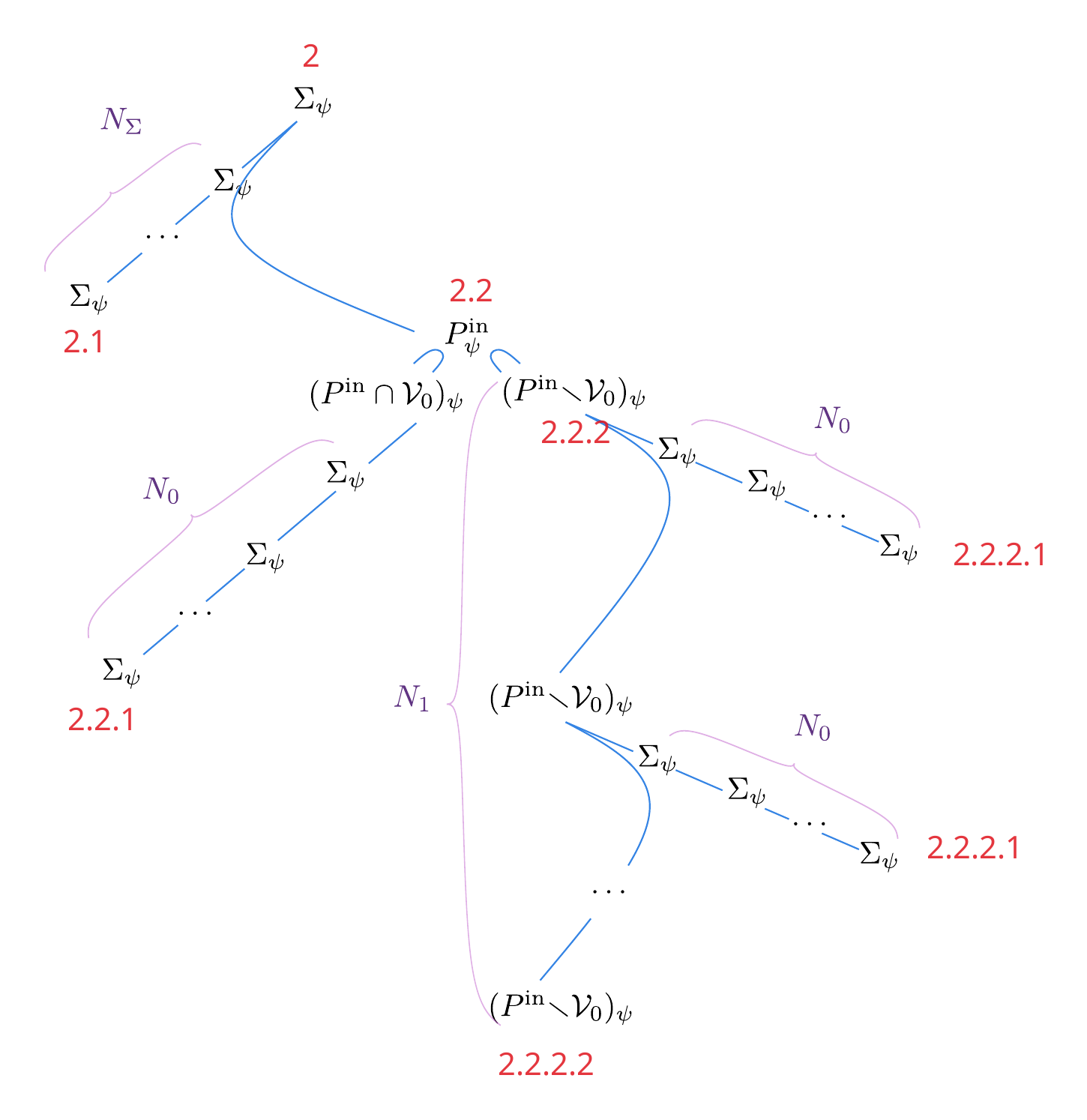}
     \vspace*{-1em}
    \caption{Dichotomy diagram of orbits in Case~\ref{it: start Sigma}}
    \label{fig: graph orbit start Sigma}
    \end{subfigure}
    \caption{Dichotomy diagram of $f_\psi$-orbits on $S_\psi$}
\end{figure}

    \vspace*{0.2em}
    \item \textbf{Suppose $p \in \Sigma_\psi$, and $T_p = T_p^\Sigma$.}
    \label{it: start Sigma}
   We have a dichotomy depending on whether the orbit of $p$ by the flow of $X_\psi$ crosses $\Sigma_\psi$ enough times without crossing $\Pin_\psi$ to have strict invariance and expansion of the cones, or if it crosses $\Pin_\psi$ before.
This return number is the integer $N_\Sigma$ which satisfies Item~\ref{prop: parameters; it: sigma to sigma} of Proposition~\ref{prop: parameters and cones}.
    
    \begin{enumerate} [label=\theenumi\arabic*., leftmargin=*]
        \item \textbf{The orbit of $p$ by the flow of $X_\psi$ crosses $N_\Sigma$ consecutive times without crossing $\Pin$.}   
        In other words $f^K_\psi(p) \in \Sigma_\psi$ for $k=0, \dots, N_\Sigma$.
        Then according to Item~\ref{prop: parameters; it: sigma to sigma} of Proposition~\ref{prop: parameters and cones}, Proposition~\ref{prop: future iterate for f_psi on S_psi} is satisfied for $n^\star(p) = N_\Sigma$ and $T_q = T_q^\Sigma$.
        
        \item \textbf{The orbit of $p$ by the flow of $X_\psi$ crosses $\Pin_\psi$ before}.
        In other words there exists $k_1 \leq N_\Sigma$ minimal such that $f^{k_1}_{0, \psi}(p) \in \Pin_\psi$.
        Then two cases are possible, depending on whether this intersection point is in the neighborhood $(\cV_0)_\psi$ or not.
        
        \begin{enumerate} [label=\theenumii\arabic*., leftmargin=*, wide=0pt]
        
            \item \label{case: sigma to Pin inter v0}
            $\boldsymbol{f^{k_1}_{0, \psi}(p) \in (\Pin \cap \cV_0)_\psi}$.    
            Then according to Item~\ref{prop: parameters; it: v0} of Proposition~\ref{prop: parameters and cones}, the future orbit of $f^{k_1}(p)$ by the flow of $X_\psi$ crosses $N_0$ consecutive times the section $\Sigma_\psi$ without crossing $\Pin_\psi$.
            In other words $f^K_\psi \l( f^{k_1}_{0, \psi}(p) \r) \in \Sigma_\psi$ for $k=1, \dots, N_0$.
            Then according to Item~\ref{prop: parameters; it: sigma to sigma through Pin} of Proposition~\ref{prop: parameters and cones}, Proposition~\ref{prop: future iterate for f_psi on S_psi} is satisfied for $n^\star(p) = N_0 + {k_1}-1 \leq N_0 + N_\Sigma -1$ and $T_q = T_q^\Sigma$.
            Note that the endpoint $q$ on $\Sigma_\psi$ corresponds to the iteration $N_0 + k_1$ of $p$ by $f_{0, \psi}$, and the iteration $N_0 + k_1 -1$ of $p$ by $f_\psi$, because in the second case we do not count the crossing of the orbit by $(\Pin \cap \cV_0)_\psi$.
           
            \item \label{case: sigma to Pin minus v0}
          $\boldsymbol{f^{k_1}_{0, \psi}(p) \in (\Pin \ssm \cV_0)_\psi}$.
            Let $q_1 = f^{k_1}_{0, \psi}(p) = f^{k_1}_\psi(p)$.
            Either the orbit of $q_1$ crosses the $\Sigma$ section enough consecutive times without crossing $\Pin_\psi$.
            This return number is the integer $N_0$ which satisfies Item~\ref{prop: parameters; it: sigma to Pin minus v0} of Proposition~\ref{prop: parameters and cones}, and we will have strict invariance and expansion of the cones.
            Either the orbit of $q_1$ crosses $\Pin_\psi$ at a point $q_2$ with less than $N_0$ crossing through $\Sigma_\psi$ between each.
            According to Item~\ref{prop: parameters; it: v0} of Proposition~\ref{prop: parameters and cones}, the point of intersection $q_2$ is in $(\Pin \ssm \cV_0)_\psi$.
            In this case, we can reiterate this dichotomy by replacing $q_2$ by $q_1$.
            This recurrence is finite because if the orbit of $q_1$ crosses $(\Pin \ssm \cV_0)_\psi$ enough times without crossing $(\Pin \cap \cV_0)_\psi$, then we have the strict invariance and expansion of cones. This return number is the integer $N_1$ which satisfies Proposition~\ref{prop: parameters and cones}, Item~\ref{prop: parameters; it: sigma to Pin minus v0}.
            Formally, the dichotomy is as follows.
            
            \begin{enumerate} [label=\theenumiii\arabic*.,leftmargin=16pt, labelwidth=-36pt]
                \item \textbf{The orbit of $p$ crosses $\Pin_\psi$ $N_1$ times, with strictly less than $N_0$ crossings through $\Sigma_\psi$ between each}.
                Then each crossing through $\Pin_\psi$ occurs in $(\Pin \ssm \cV_0)_\psi$.
                There exists $k_1< \dots < k_{N_1}$ such that
                \begin{itemize}[--]
                    \item $k_i \leq i . N_0$,
                    \item $f_{0, \psi}^{k_i}(p) \in (\Pin \ssm \cV_0)_\psi$,
                    \item $f^k_{0, \psi}(p) \in \Sigma_\psi$ for $k \leq k_{N_1}$ and $k \neq k_1, \dots, k_{N_1}$.
                \end{itemize}
                According to Item~\ref{prop: parameters; it: sigma to Pin minus v0} of Proposition~\ref{prop: parameters and cones}, Proposition~\ref{prop: future iterate for f_psi on S_psi} is satisfied for $n^\star(p)= k_{N_1} \leq N_1 \times (N_0-1)$ and $T_q = T_q^\iin$.
                \item \label{case: sigma to sigma through Pin minus v0}
                \textbf{The orbit of $p$ crosses $\Pin_\psi$ strictly less than $N_1$ times,  with strictly less than $N_0$ crossings through $\Sigma_\psi$ between each,
                then its orbit crosses the section $\Sigma_\psi$ $N_0$ consecutive times without crossing $\Pin$.}
                Then each crossing through $\Pin_\psi$ occurs in $(\Pin \ssm \cV_0)_\psi$.
                There exists $m \leq N_1 -1$ and $k_1< \dots <k_m$ such that
                \begin{itemize}[--]
                    \item $k_i \leq i . N_0 \in (\Pin \ssm \cV_0)_\psi$,
                    \item $f^k_{0, \psi}(p) \in \Sigma_\psi$ for $k\leq k_m$ and $k \neq k_1 \dots, k_m$,
                    \item $f_{0, \psi}^k \l( f_{0, \psi}^{k_m}(p) \r) \in \Sigma_\psi$ for $k=1, \dots, N_0$,
                \end{itemize}
                
                We set $q_i = f^{k_i}(p)$ for $i=1, \dots, m$.
                Note that according to Item~\ref{prop: parameters; it: sigma to Pin}, $f^{k_1}_\psi$ maps $C^u_\psi(p)$ inside $C^u_\psi(q_1)$ and its inverse maps $C^s_\psi(q_1)$ inside $C^s_\psi(p)$.
                Then, each crossing of the orbit from $(\Pin \ssm \cV_0)_\psi$ to $(\Pin \ssm \cV_0)_\psi$ contributes to map the cone $C^u_\psi(q_i)$ inside the cone $C^u_\psi(q_{i+1})$ and expand the norm of the vectors by a factor of 2 according to Item~\ref{prop: parameters; it: Pin to Pin}.
                So it is enough to show that an orbit of a point of $\Sigma_\psi$ which crosses $\Pin_\psi$ and then crosses $N_0$ times the section $\Sigma_\psi$ without crossing $\Pin_\psi$ leaves the cones of $C^u_\psi$ strictly invariant and expanded, and similarly for the action of the inverse on the cones $C^s_\psi$.
                This is Item~\ref{prop: parameters; it: sigma to sigma through Pin} of Proposition~\ref{prop: parameters and cones}.
                We deduce that Proposition~\ref{prop: future iterate for f_psi on S_psi} is satisfied for $n^\star(p)= N_0 + k_m \leq N_0 + (N_1 -1) \times N_0$, and $T_q = T_q^\Sigma$.\qedhere          
            \end{enumerate}     
        \end{enumerate}
    \end{enumerate}
\end{enumerate}
\end{proof}

\subsubsection*{Additional result}

As announced in Subsection~\ref{sec: parameters; subsec: additional result}, we have an additional complementary result on the the maximal negative invariant set $K_\psi^-$ for $f_\psi$ in $S_\psi$ (Definition~\ref{def: positive, negative, maximal invariant set f_psi S_psi}).
It gives the existence of an integer $\hat n^\star(p) \leq 0$ for any point $p \in K_\psi^-$, uniformly bounded from below, such that the iterate $f^{\hat n^\star}_\psi$ satisfies the cone fields condition in the neighborhood of $p$ for the pair of cone fields $(\hat{C}^s_\psi, \hat{C}^u_\psi)$ given by Proposition~\ref{prop: parameters and cones}.
This means that the behavior of the map $f^{\hat n^\star}$ is hyperbolic in the \emph{past} orbits.

\begin{prop} \label{prop: past iterate for f_psi on S_psi}
There exists an integer $\bN \geq 0$ such that,
for all $p \in K_\psi^-$, if $T_p$ is (one of) the tangent plane(s) to $S_\psi$ in $p$,
there exists an integer $\hat n^\star = \hat n^\star(p, T_p)$ with $-\bN \leq \hat n^\star \leq 0$,
and a tangent plane $T_q = T_q(p, T_p)$ to $S_\psi$ in
$q = f_\psi^{\hat n^\star}(p) \in S_\psi$
such that:
\begin{enumerate}
    \item $(f_\psi^{\hat n^\star})_* \hat C^s_\psi(p, T_p) \subset \intr \hat C^s_\psi(q, T_q)$,
    and for all $v \in \hat C^s_\psi(p, T_p)$, 
    $$\Vert (f_\psi^{\hat n^\star})_* v \Vert \geq 2 \Vert v \Vert,$$
    \item $(f_\psi^{\hat n^\star})^\inv_* \hat C^u_\psi(q, T_q) \subset \intr \hat C^u_\psi(p, T_p)$,
    and for all $v \in \hat C^u_\psi(q, T_q)$, 
    $$\Vert (f_\psi^{\hat n^\star})^\inv_* v \Vert \geq 2 \Vert v \Vert.$$
\end{enumerate}
\end{prop}

The proof is identical using Proposition~\ref{prop sym: parameters and cones}.

\subsection{Hyperbolic splitting on the local section}
\label{sec: proof; subsec: hyperbolic splitting global section}

Propositions~\ref{prop: future iterate for f_psi on S_psi} and~\ref{prop: past iterate for f_psi on S_psi} will allow us to construct a hyperbolic splitting of the tangent bundle of $S_\psi$ on the maximal invariant set $K_\psi$ for the return map $f_\psi \colon S_\psi \to S_\psi$.
Let us show the following proposition which gives the existence of a hyperbolic splitting for the return map $f_\psi$ on $K_\psi$.
More precisely, we construct the stable direction on all of $K_\psi^+$, and the unstable direction on all of $K_\psi^-$.
We will also need that the fibers of the invariant direction on two different tangent planes generates the same sum with $\R. X_\psi$.
A line bundle $F$ on $S_\psi$ is the data, for any $p\in S_\psi$ and any tangent plane $T^*=T^*_p$ to $S_\psi$ at $p$, of a 1-dimensional vector space $F(p, T^*) \subset T^*$.

\begin{prop}
\label{prop: hyperbolic splitting on S_psi}
There exists a $f_\psi$-invariant direction $F^s \subset TS_\psi$ on $K_\psi^+$, and a transverse $f_\psi^\inv$-invariant direction $F^u \subset TS_\psi$ on $K_\psi^-$, and constants $\lambda >1$ and $C>0$ such that: 
\begin{enumerate}
    \item  \label{prop: hyp splitting S_psi; it: cones}
    $F^s \subset \intr C^s_\psi$ and $F^u \subset \intr C^u_\psi$.
    \item \label{prop: hyp splitting S_psi; it: contraction}
    For all $p \in K_\psi^+, \ \forall n \geq 0, \ \forall v \in F^s(p),\ \Vert (f^n_\psi)_* v \Vert \leq C \lambda^{-n} \Vert v \Vert$.\\
    For all $p \in K_\psi^-, \ \forall n \leq 0, \ \forall v \in F^u(p),\ \Vert (f^n_\psi)_* v \Vert \leq C \lambda^n \Vert v \Vert$.
    \item \label{prop: hyp splitting S_psi; it: generated planes}
    If $T$ and $T'$ are two planes tangent to $S_\psi$ at $p \in K_\psi^+$, then $$F^s(p,T) \oplus \R. X_\psi = F^s(p, T') \oplus \R. X_\psi .$$
    If $T$ and $T'$ are two planes tangent to $S_\psi$ at $p \in K_\psi^-$, then $$F^u(p,T) \oplus \R. X_\psi = F^u(p, T') \oplus \R. X_\psi .$$
\end{enumerate}
\end{prop}

\begin{proof}

In the proof we will just show the existence of the stable bundle $F^s$ on $K_\psi^-$, the proof for $F^u$ being similar.

\subsubsection*{Step 1: Infinite intersection of stable cones}
Let $p$ be a point of $K_\psi^+$ and $T$ (one of) the tangent plane(s) to $S_\psi$ at $p$.
Then there exists a sequence of particular points on the positive orbit of $p$ by $f_\psi$ defined as follows.
Let $(p_k, T_k)$ and $n_k$, defined by recurrence:
\begin{equation} \label{eq: p_k, T_k, n_k}
    \left\{ \begin{array}{l}
    (p_0, T_0) = (p, T), \\
    (p_{k+1}, T_{k+1})= \l(f_\psi^{n^\star(p_k, T_k)}(p_k), \ T_{p_{k+1}}(p_k, T_k)\r)
\end{array}\right. ,\ 
\left\{ \begin{array}{l}
    n_0 =  n^\star(p,T),\\
    n_{k+1} = n^\star( p_k, T_k) 
\end{array}\right.,
\end{equation}
where $n^\star(p, T)$ is the integer and $T_q (p, T)$ is the tangent plane satisfying Proposition~\ref{prop: future iterate for f_psi on S_psi} for the given $(p, T)$.
Recall that by Proposition~\ref{prop: future iterate for f_psi on S_psi} the sequence $n_k$ is bounded by an integer $\bN$.
Moreover, $f_\psi^{n_k}$ maps the cone $C^u_\psi(p_k, T_k)$ inside $C^u_\psi(p_{k+1}, T_{k+1})$ and expands the norm of the vectors of $C^u_\psi(p_k, T_k)$ by a factor 2, and $\l(f_\psi^{n_k}\r)^\inv$ maps the cone $C^s_\psi(p_{k+1}, T_{k+1})$ inside $C^s_\psi(p_k, T_k)$ and expands the norm of the vectors of $C^s_\psi(p_{k+1}, T_{k+1})$ by a factor 2.
Define $L_k$ to be the differential of $f_\psi^{n_k}$ from the plane $T_k$ to the plane $T_{k+1}$, in other words
\begin{equation} \label{eq: L_k}
    L_k := d_{p_k} f_\psi^{n_k} \colon T_k \to T_{k+1}.
\end{equation}
It is a family of invertible linear maps. Define
\begin{equation}
\label{eq: F^s}
    F^s(p,T) := \underset{k=0}{\overset{+\infty}{\bigcap}} \l( L_0^\inv \circ \dots \circ L_k^\inv \r) C^s_\psi(p_k, T_k) \subset C^s_\psi(p, T).
\end{equation}
By strict invariance of the stable cones $C^s_\psi$ under the iterates $f_\psi^{-n^\star}$, this infinite intersection is a cone contained in $C^s_\psi(p, T)$.
Repeating this operation for each point $p$ of $K_\psi$ and tangent plane $T$ of $S_\psi$ in $p$, we obtain a cone field $F^s = \{ F^s(p, T) \}$ on $S_\psi$ on $K_\psi^+$.
Then we have for all $k \geq 0$:
\begin{equation} \label{eq: invariance and contraction F^s}
F^s(p_{k+1}, T_{k+1}) = L_k F^s(p_k, T_k) \qq \mb{and} \qq \forall v \in F^s_k, \ \Vert L_k v \Vert \leq \frac{1}{2} \Vert v \Vert.
\end{equation}
The first equality follows from Equation~\eqref{eq: F^s} defining $F^s$ by an infinite intersection.
The following inequality is true because $F^s(p_k, T_k) \subset L_k^\inv C_\psi^s(p_{k+1}, T_{k+1})$ by~\eqref{eq: F^s}, and $L_k^\inv$ expands the norm of the vectors of $C_\psi^s(p_{k+1}, T_{k+1})$ by a factor~2.

\subsubsection*{Step 2: Line bundle}
\begin{claim} \label{claim: F^s line bundle}
$F^s$ is a line bundle.
\end{claim}

\begin{proof}
Let $(p, T)$ be fixed, $(p_k, T_k, n_k)$ the sequence associated by Equation \eqref{eq: p_k, T_k, n_k}.
Let $F^s_k = F^s(p_k, T_k)$ be the family of cones defined by Equation~\eqref{eq: F^s}, and \linebreak[4]$L_k \colon T_k \to T_{k+1}$ the family of invertible linear maps defined by Equation~\eqref{eq: L_k}.
According to Equation~\eqref{eq: invariance and contraction F^s}, $\{F^s_k\}_k$ is a collection of cones (exactly) invariant by $\{L_k\}_k$, and whose vectors are uniformly contracted by $\{L_k\}_k$.
Suppose by contradiction that $F^s_0$ is not reduced to a line.
Then the same is true for $F^s_k$ and we denote $F^1_k$ and $F^2_k$ the two disjoint lines which constitute the boundary of $F^s_k$, and such that $L_k F^i_k = F^i_{k+1}$.
Therefore, any vector $v \in T$ decomposes uniquely as $v = v_1 + v_2 \in F^1 \oplus F^2$, and we have by iteration of the inequality in Equation~\eqref{eq: invariance and contraction F^s}
\begin{equation} \label{eq: L_k contracts everything}
    \Vert L_n \circ L_{n-1} \circ \dots \circ L_0 (v) \Vert \leq 2^{-n} ( \Vert v_1 \Vert + \Vert v_2 \Vert ) \underset{n \to +\infty}{\longrightarrow} 0.
\end{equation}
In other words, the norm of any vector tends to 0 under the action of $\{L_k\}_k$.
Now there is an invariant direction expanded by $\{ L_k\}_k$ : it is enough to iterate any direction in the unstable field $C^u_\psi(p, T)$.
More precisely, according to Proposition~\ref{prop: future iterate for f_psi on S_psi}, and by definition of $(p_k, T_k, n_k)$ to Equation~\eqref{eq: p_k, T_k, n_k} and $L_k$ to Equation~\eqref{eq: L_k}, we have $L_k \, C^u_\psi(p_k, T_k) \subset \intr C^u_\psi(p_{k+1}, T_{k+1})$ and for all $v_k \in C^u_\psi(p_k, T_k)$, $\Vert L_k v_k \Vert \geq 2 \Vert v \Vert$.
We conclude by iterating the inequality for $v \in C^u_\psi(p, T) \subset T$:
$$ \Vert L_n \circ L_{n-1} \circ \dots \circ L_0 (v) \Vert \geq 2^n \Vert v \Vert \underset{n \to +\infty}{\longrightarrow} +\infty .$$
This is a contradiction with Equation~\eqref{eq: L_k contracts everything}.
\end{proof}

\subsubsection*{Step 3: Exponential contraction}
Let us show that the line bundle $F^s$ is contracted by a uniform iteration of $f_\psi$.

\begin{claim} \label{claim: uniform iterate contracts F^s}
There exists an integer $\bN_0$, such that for any $p \in K_\psi$ and $T \subset T_p S_\psi$ a tangent plane in $p$,
$$ \forall v \in F^s(p, T), \qq \Vert (f^{\bN_0}_\psi)_* v \Vert \leq \frac{1}{2} \Vert v \Vert.$$ 
\end{claim}

\begin{proof}
Let $(p, T)$ be fixed, $(p_k, T_k, n_k)$ the sequence associated by Equation \eqref{eq: p_k, T_k, n_k}.
Let $F^s_k = F^s(p_k, T_k)$ be the family lines defined by Equation~\eqref{eq: F^s}, and $L_k \colon T_k \to T_{k+1}$ be the family of invertible linear maps defined by Equation~\eqref{eq: L_k}.
Let $\bN$ be the uniform integer which bounds the sequence $\{ n_k \}_k$.
There exists a constant $\cst >0$ which bounds the differential of $f_\psi^i$ for $0 \leq i \leq \bN$.
Let be an integer $k_0$ satisfying
\begin{equation}\label{eq: k_0}
    k_0 > \log_2 \cst + 1
\end{equation}
and let $\bN_0$ be an integer greater than $k_0 . \bN$.
Then we can decompose $f_\psi^{\bN_0}$ in the neighborhood of $p$ in the form
$f_\psi^{\bN_0} = f_\psi^i \circ f_\psi^{n_k} \circ \dots \circ f_\psi^{n_0}$,
with $k \geq k_0$ and $0 \leq i \leq \bN$.
Then by iterating the inequality of Equation~\eqref{eq: invariance and contraction F^s} (with $L_k = d_{p_k} f^{n_k}_\psi$) for a $v \in F^s$, we have
\begin{align*}
    \Vert (f_\psi^{\bN_0})_* v \Vert &= \Vert (f_\psi^i \circ f_\psi^{n_k} \circ \dots \circ f_\psi^{n_0})_* v \Vert \\
     &\leq \cst \Vert (f_\psi^{n_k} \circ \dots \circ f_\psi^{n_0})_* v \Vert\\
    & \leq \cst 2^{-k} \Vert v \Vert & \mb{Equation~\eqref{eq: invariance and contraction F^s}} \\
    & \leq \frac{1}{2} \Vert v \Vert & \mb{Equation~\eqref{eq: k_0}}
\end{align*}
This completes the proof of the claim.
\end{proof}

\subsubsection*{Step 4: Invariance}
Let us show that $F^s$ is an $f_\psi$-invariant bundle.

\begin{claim} \label{claim: F^s invariant}
For any $p \in K_\psi^+$ and $T \subset T_p S_\psi$ a tangent plane at $p$, for any $j \geq 0$, if $q = f_\psi^j(p)$ and $T' \subset T_q S_\psi$ is a tangent plane at $q$, then $ F^s(q, T')= (f_\psi^j)_* F^s(p, T) $.
\end{claim}

\begin{proof}
Suppose by contradiction that $F^s(q, T') \neq (f_\psi^j)_* F^s(p, T)$.
Let $F^1 = F^s(q, T')$ and $F^2 =(f_\psi^j)_* F^s(p, T) $.
These are two lines in $T'$.
We can assume without loss of generality that $0 \leq j \leq \bN_0$.
According to Claim~\ref{claim: uniform iterate contracts F^s}, we have for $ v_1 \in F^1$ 
and $n \geq 0$, \,
$\Vert (f^{n \bN_0})_* v_1 \Vert \leq 2^{-n} \Vert v_1 \Vert$,
and for $v_2 \in F^2$ and $n \geq 0$,
$$\Vert (f_\psi^{n \bN_0})_* v_2 \Vert = \Vert f^j_\psi (f_\psi^{\bN_0})^n_* (f^{-j}_\psi)_*v_2  \Vert
\leq \cst 2^{-n} \Vert v_1 \Vert .$$
Any vector $v \in T'$ decomposes uniquely as $v = v_1 + v_2 \in F^1 \oplus F^2$.
We deduce that the norm of any vector $v \in T'$ under the action of $\{ (f^{n \bN_0}_\psi)_* \}$ tends to 0.

Now there exists an invariant direction expanded by the family $\{ (f_\psi^{n \bN_0})_* \}_n$, it is enough to iterate any direction in the unstable cone field $C^u_\psi(q, T')$.
More precisely, let $(q_k, T'_k, n_k)$ be the sequence defined by Equation~\eqref{eq: p_k, T_k, n_k} from $(q, T')$.
According to Proposition~\ref{prop: future iterate for f_psi on S_psi} we have for all $k \geq 0$,\, $$(f_\psi^{n_k})_* \, C^u_\psi(q_k, T'_k) \subset \intr C^u_\psi(q_{k+1}, T'_{k+1}),$$ and for all $v \in C^u_\psi(q_k, T'_k)$,\, $\Vert (f_\psi^{n_k})_* v \Vert \geq 2 \Vert v \Vert$.
For all $n \geq 0$, there exists a unique decomposition of $f_\psi^{n \bN_0}$ in the neighborhood of $q'$ of the form
$f_\psi^{n \bN_0} = f_\psi^i \circ f_\psi^{n_{k_n}} \circ \dots \circ f_\psi^{n_0} $ with $0 \leq i < \bN_0$.
Moreover, $k_n \underset{n \to \infty}{\longrightarrow} + \infty$.
There is also a uniform constant $\cst >0$ which bounds the differential of $f_\psi^i$ for $0 \leq i \leq \bN_0$.

We conclude that for $v \in C^u_\psi(q, T') \subset T'$:
\begin{align*}
    \Vert (f_\psi^{n \bN_0})_* v \Vert & = \Vert (f_\psi^i \circ f_\psi^{n_{k_n}} \circ \dots \circ f_\psi^{n_0})_* v \Vert \\
    & \geq \cst \Vert (f_\psi^{n_{k_n}} \circ \dots \circ f_\psi^{n_0})_* v \Vert \\
    & \geq \cst^\inv 2^{k_n} \Vert v \Vert  \underset{n \to +\infty}{\longrightarrow} +\infty.
\end{align*}
This is a contradiction.
\end{proof}

\subsubsection*{Step 5: Generated tangent planes.}
We have shown that $F^s$ is a $f_\psi$-invariant line bundle of $S_\psi$ on the subset $K_\psi^+$, and exponentially contracted by $f_\psi$, i.e., $F^s$ satisfies Item~\ref{prop: hyp splitting S_psi; it: contraction} of Proposition~\ref{prop: hyperbolic splitting on S_psi}.
Let us show that $F^s$ satisfies Item~\ref{prop: hyp splitting S_psi; it: generated planes} of Proposition~\ref{prop: hyperbolic splitting on S_psi}.

\begin{claim}
Let $T$ and $T'$ be two tangent planes at $p \in K_\psi^+$ in $T_p S_\psi$, then $$F^s(p, T) \oplus \R.X_\psi = F^s(p, T') \oplus \R.X_\psi.$$
\end{claim}

\begin{proof}
Let $q = f_\psi(p)$ and $T_q$ be a tangent plane to $S_\psi$ at $q$.
Let $L \colon T \to T_q$ and $L' \colon T' \to T_q$ be the two restrictions of the differential $d_p f_\psi$.
Then by Claim~\ref{claim: F^s invariant},
$L (F^s(p, T)) = F^s(q, T_q) = L (F^s(p, T'))$.
Now $L$ and $L'$ differ at a projection parallel to $\R.X_\psi$, in other words if $\pi \colon T \to T'$ is the projection on $T'$ parallel to $\R.X_\psi$, then we have $L = L' \circ \pi$.
Since $L'$ is invertible, we deduce that $\pi F^s(p, T) = F^s(p, T')$.
The equality of sums with $\R. X$ is therefore true.
\end{proof}

It follows that the line bundle $F^s$ on $S_\psi$ satisfies Proposition~\ref{prop: hyperbolic splitting on S_psi}.
The proof for the line bundle $F^u$ is symmetric. 
We have to use Proposition~\ref{prop: past iterate for f_psi on S_psi}, construct a sequence in the past orbit of a point $p \in K_\psi^-$, the iteration is given by the integer $\hat{n}^\star$, and the action of the inverse differential on the pair $(\hat C^s_\psi, \hat C^u_\psi)$.
\end{proof}

\subsection{Proof of Theorem~\ref{thmintro: gluing theorem}}
\label{sec: proof; subsec: proof}
We want to show that the vector field $X_\psi$ on $P_\psi$ is Anosov.
We will use the following general criterion, see for example \cite[Theorem~6.2.24]{fisherHyperbolicFlows2019}.

\begin{lem}[Ma\~n\'e Criterion]
    A $\cC^1$ flow on a compact manifold is an Anosov flow if and only if the chain-recurrent set is hyperbolic, the stable and unstable manifolds intersect transversely at one (hence every) point of each orbit, and their dimension is constant.
\end{lem}

To use the Ma\~n\'e Criterion, we will need to prove that the chain-recurrent set of $X_\psi$ in $P_\psi$ is hyperbolic. 
To do so, we use the following general lemma.

\begin{lem} \label{lem: expanded plane field criteria hyperbolic set}
Let $\Lambda$ be a compact invariant set for a $\cC^1$ flow $Y$ on a compact 3-manifold $M$.
Suppose that the differential of the flow of $Y$ preserves two plane fields $E^\cs$ and $E^\cu$ on $P$ on $\Lambda$, tangent to $Y$, and transverse to each other, such that $Y^t$ contracts exponentially the area of $E^\cs$ and expands exponentially the area of $E^\cu$ for $t \geq 0$.
Then $\Lambda$ is a hyperbolic set of index (1,1).
\end{lem}

\begin{proof}
Let us show the existence of the strong unstable direction $E^\uu \subset E^\cu$.
Let $V$ be a vector field tangent to $E^\cu$ and transverse to the vector field $Y$ on $E^\cu$, such that the pair $(V,Y)$ forms a vector field basis of $E^\cu$.
Let $L^t := \res{Y^t_*}{E^\cu}$.
The matrix of $L^t$ in the basis $(V,Y)$ is
$$ M_{(V,Y)}L^t = \begin{pmatrix}
    a_t & 0 \\
    b_t & 1
\end{pmatrix},$$
where $a_t$ and $b_t$ are continuous functions of $\cM$ and bounded at $t$ fixed.
Let us show that $L^t$ satisfies the cone fields condition.
By hypothesis, there exist constants $C>0$ and $\lambda >1$ such that
$\operatorname{Jac} (L^t) = \vert a_t \vert \geq C \lambda^t  $.
Up to change the metric, we can suppose that
$\inf \vert a_1 \vert > 1$.
We suppose that $b_1$ is not uniformly zero, otherwise the lemma is true : the strong unstable direction $E^\uu$ is the direction generated by $V$.
We define
\begin{equation} \label{eq: slope unstable cone}
    K := 2 \, \frac{\sup \vert b_1 \vert}{\inf \vert a_1 \vert -1} >0.
\end{equation}
Let $C^u$ be the $K$-cone field on $E^\cu$ in the $(V,Y)$ basis.

\begin{claim}
For all $n \geq 1$, we have
$L^n C^u \subset \intr C^u$,
and $L^n$ exponentially expands the norm of the vectors of $C^u$.
\end{claim}

\begin{proof}
Let $v = v^u V + v^Y Y \in C^u$, then we have $\vert v^Y \vert \leq K \vert v^u \vert$.
Then $L^1 v = a_1 v^u V + (b_1 v^u + v^Y) Y$ and
$$ \frac{\vert b_1 v^u + v^Y \vert}{\vert a_1 v^u \vert} 
\leq \frac{ \vert b_1 \vert K^\inv +1 }{\vert a_1 \vert} K  < K.$$
So $L^1 v$ is inside $C^u$.
We deduce the property for $L^n$ by writing $L^n = L^1 \circ \dots \circ L^1$.

Let $n$ be a positive integer and $L^n v = v^u_n V + v^Y_n Y = a_n v^u V + (b_n v^u + v^Y) Y$.
We have $$\l\lbrace
\begin{array}{rcl l}
    \vert v^u_n \vert & \geq & (\inf \vert a_1 \vert)^n \vert v^u \vert> 2 \vert v^u \vert , & \forall n \geq \log_{\inf \vert a_1 \vert}(2) \\
    \vert v^u_n \vert &\geq & (\inf \vert a_1 \vert)^n \frac{1}{K} \vert v^Y \vert > 2 \vert v^Y \vert, & \forall n \geq \log_{\inf \vert a_1 \vert}(2K).
\end{array} \r.$$

Therefore $\max (\vert v_1^n \vert, \vert v_Y^n \vert) > 2 \max(\vert v_1 \vert, \vert v_Y \vert)$ for $n$ (uniformly) large enough, which proves that $L^n$ expands the norm of the vectors of $C^u$ by a factor of $2$ for $n$ greater than some integer $n_0$.
Since the operator $L^i$ is uniformly bounded for $0 \leq i \leq n_0$,
we deduce that $L^n$ expands exponentially the norm of the vectors of $C^u$ by iterating $L^{k n_0+i} = L^i \circ L^{n_0} \circ \dots \circ L^{n_0}$.
\end{proof}

Define
$E^\uu (p) := \underset{n \geq 0}{\bigcap} 
L^n C^u (Y^{-n}(p)) \subset C^u(p)$.

\begin{claim}
$E^\uu$ is an $L^t$-invariant line bundle, exponentially expanded by $L^t$ for $t \geq 0$.
\end{claim}

\begin{proof}
Indeed, by construction $E^\uu$ is a cone field in the plane field $E^\cu$ on $K$, invariant by $L^{-n}$ and exponentially contracted by $L^{-n}$ for $n \geq 0$.
Now we know that $L^{-n}$ acts isometrically on the direction generated by $Y$, so there exists an $L^{-n}$-invariant direction which is not contracted.
Since $\{ L^{-n} \}_n$ is a family of invertible linear maps in dimension two,
we deduce that $E^\uu$ cannot contain two distinct invariant line fields (these are the same arguments as the proof of Proposition~\ref{prop: future iterate for f_psi on S_psi}).
It is therefore a line bundle exponentially contracted by $L^{-n}$.

Let us assume by contradiction that $E_1 := E^\uu(p)$ and $E_2 := L^t E^\uu(Y^{-t}(p))$ are distinct for a certain $t \in \R$ (we can assume $\vert t \vert \leq 1$),
then the families $\{ L^{-n} E_1 \}_n$ and $\{ L^{-n} E_2 \}_n$ are two line bundles invariant by $ L^{-n} $ and exponentially contracted for $n \geq 0$.
Indeed, for the second family, it suffices to write that if $v \in E^2$ then
$ L^{-n} v =  L^{-t} (L^{-n} v')$ with $v' \in E^\uu(Y^{-t}(p))$ 
and use the fact that $L^t$ is uniformly bounded for $\vert t \vert \leq 1$ and $L^{-n}$ contracts exponentially $E^\uu$.
By the same argument as in the previous paragraph, we have a contradiction because the $\R.Y$-direction is an $L^n$-invariant direction which is not contracted.
So $E^\uu$ is a line bundle $L^t$-invariant for $t \in \R$, and exponentially expanded for $t\geq 0$.
Indeed, it suffices to write for $v \in E^\uu$,
$L^t v = L^r \circ L^{E(t)} v = L^r \circ (L^{-E(t)})^\inv v$ ($E(t)$ is the integer part of $t$) and use that
$L^r$ is uniformly bounded for $r \in [0,1]$ and $L^{-n}$ contracts exponentially $E^\uu$ for $n \geq 0$.
\end{proof}

We conclude that $E^\uu$ is the strong unstable bundle of $\Lambda$ for $Y$.
We show in the same way the existence of the strong stable bundle $E^\ss$ of $\Lambda$ for $Y$.
We deduce that $\Lambda$ is a hyperbolic set for $Y$.
\end{proof}

Now is the final step of the proof of Theorem~\ref{thm: gluing theorem}.

\begin{proof}[Proof of Theorem~\ref{thm: gluing theorem}]
Let $(P_0, X_0, \varphi_0)$ be a filled \bb{}
with a \sqt{} gluing map.
Let $(P,X,\varphi)$ be the normalized strongly isotopic triple given by Proposition~\ref{prop: normalization of triple}.
Let $\lambda_0, \epsilon_0, \delta_0$, be the parameters which satisfy Proposition~\ref{prop: parameters and cones}, and
$\psi = \psi_{\lambda_0, \epsilon_0, \delta_0} \colon \pP \to \pP$
the modified gluing map.
According to Claim~\ref{claim: modify triple strongly isotopic}, the triple
$(P, X, \psi)$ is strongly isotopic to $(P,X, \varphi)$ thus to $(P_0, X_0, \varphi_0)$, and $P_\psi = P/\psi$ is a compact manifold of dimension 3 equipped with a vector field $X_\psi$ of class $\cC^1$ induced by $X$ on $P_\psi$.

Let us show that $X_\psi$ is an Anosov vector field.
Recall that $\gA$ and $\gR$ denotes the union of, respectively, attracting and repelling basic pieces of $\Lambda$.
Denote $\gA_\psi$ and $\gR_\psi$ the projection in $P_\psi$. Note that the projection is trivial in the \nbh{} of $\gA$ and $\gR$, because uniformly away from the boundary, hence $\gA_\psi$ and $\gR_\psi$ are union of hyperbolic attracting and repelling basic pieces respectively for the flow of $X_\psi$.
A filtration adapted to the splitting $\Lambda = \gR \cup \Lambda_s \cup \gA$ in $P$ for the flow of $X$ gives a filtration for the flow of $X_\psi$ in $P_\psi$ (see Figure~\ref{fig: filtration in P}).

\begin{figure}[htb]
    \centering
    \vspace*{-1em}
    \includegraphics[height=0.58\textheight]{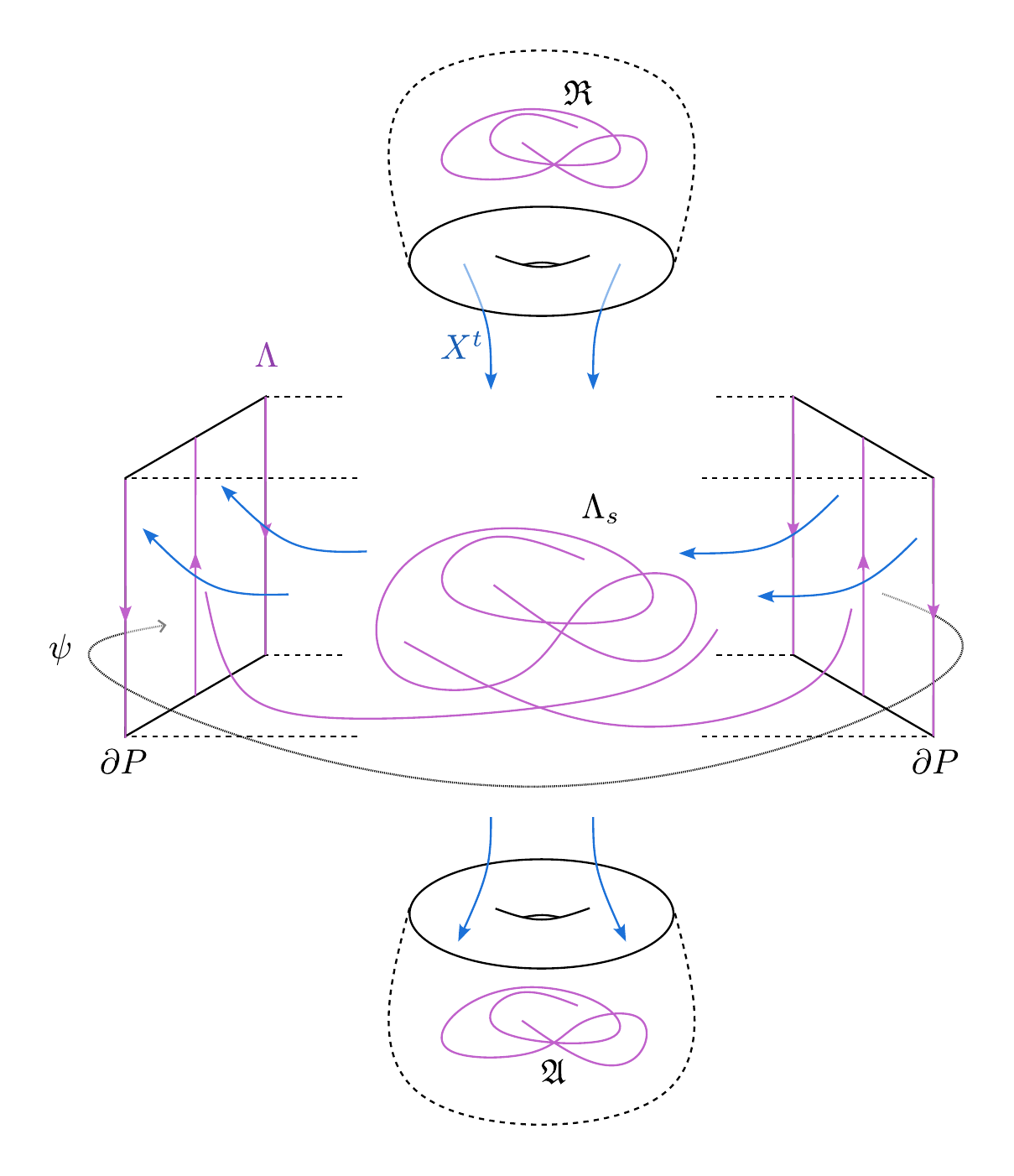}
    \vspace*{-1em}
    \caption{Filtration in $P$ adapted for the decomposition $\Lambda = \gR \cup \Lambda_s \cup \gA$ for the flow $X^t$}
    \label{fig: filtration in P}
\end{figure}

\begin{nota} \label{nota: Lambda'_psi}
    Define $P'_\psi$ the complementary of the filtrating \nbhs{} of $\gA_\psi$ and $\gR_\psi$ in $P_\psi$ for $X_\psi$, and $\Lambda'_\psi$ the maximal invariant set of $P'_\psi$ for the flow of $X_\psi$.
\end{nota}

Let us show the following result.

\begin{claim}
    The set $\gA_\psi \, \cup  \, \Lambda'_\psi \, \cup \, \gR_\psi$ is hyperbolic for $X_\psi$.
\end{claim}

\begin{proof}
It is clear that $\gA_\psi$ and $\gR_\psi$ are hyperbolic sets for $X_\psi$ because they inherits the hyperbolic structure of $\gA$ and $\gR$ for $X$ in $P$ (the projection restricted to small \nbhs{} of $\gA$ and $\gR$ is trivial, because the sets are contained in the interior of $P$).
Let us show that $\Lambda'_\psi$ is a hyperbolic set.
Let $S_\psi$ the section
given by Equation~\eqref{eq: S_psi} and $f_\psi$ the first return map of $X_\psi$ on $S_\psi$.
First remark that $S_\psi$ intersect every orbit of $\Lambda'_\psi$ in uniformly bounded time.
Indeed the filtration adapted to $\gA_\psi \cup \Lambda'_\psi \cup \gR_\psi$ is the projection of a filtration adapted to 
$\gR \cup \Lambda_s \cup \gA$.
Hence an orbit $\Lambda'_\psi$ is either an orbit of $(\Lambda_s)_\psi$, or crosses $\Pin_\psi$, and never escapes to $\gA_\psi$ in the future of $\gR_\psi$ in the past.
It follows that every segment of this orbit will cross $\Sigma_\psi \cup \Pin_\psi$, hence $\Sigma_\psi \cup (\Pin \ssm \cV_0)_\psi = S_\psi$, in uniformly bounded time, .
Hence it suffices to use the fact that $K_\psi$ is a hyperbolic set for $f_\psi$.

Let $\res{T S_\psi}{K_\psi} = F^u \oplus F^s$ be the hyperbolic splitting of $S_\psi$ for $f_\psi$ given by Proposition~\ref{prop: hyperbolic splitting on S_psi} above $K_\psi = \Lambda'_\psi \cap S_\psi$.
We define at any point $p \in K_\psi$ the tangent plane in $T_p P_\psi$ by
$E^\cu(p) := F^u(p) \oplus \R.X_\psi(p)$.
Proposition~\ref{prop: hyperbolic splitting on S_psi} ensures that this definition does not depend on the choice of a tangent plane $T^* \subset T_p S_\psi$ in $p$.
Let $p \in \Lambda'_\psi$.
Then the orbit of $p$ intersects $S_\psi$ in uniformly bounded time, and let $\tau(p) \geq 0$ be the first positive time such that $q = X^{\tau(p)}_\psi(p) \in S_\psi$.
Define for all $p\in P_\psi$,
$E^\cu(p) := (X_\psi^{-\tau(p)})_* E^\cu(q)$.
This definition gives a plane field $E^\cu$ at any point of $\Lambda'_\psi$.
Let us show that 
$E^\cu$ is a $X_\psi^t$-invariant plane field, and whose area is exponentially expanded by the differential of the flow of $X_\psi^t$.

Indeed, the invariance of $E^\cu$ by the differential of the flow of $X_\psi$ is ensured by the invariance of $F^u$ by the first return map $f_\psi$ on $S_\psi$ according to Proposition~\ref{prop: hyperbolic splitting on S_psi}.
Let $V$ be a unit vector field tangent to $F^u$ at any point of $S_\psi$.
By transversality of $S_\psi$ and $X_\psi$, the pair $(V,X_\psi)$ forms a basis of vector fields on the plane field $E^\cu$ at any point of $S_\psi$.
The matrix of the differential of $X^t_\psi$ restricted to $E^\cu$ in the basis $(V,X)$ is of the form
$$ M_{(V,X_\psi)} \res{(X^t_\psi)_*}{E^\cu} = \begin{pmatrix}
    a_t & 0 \\
    b_t & 1
\end{pmatrix},$$ 
where $a_t$ and $b_t$ are continuous functions on $S_\psi$.
By compactness, these functions are uniformly bounded at $t$.
By Proposition~\ref{prop: hyperbolic splitting on S_psi},
there exist constants $C>0$ and $\lambda>0$ such that
for all $p \in S_\psi$, if $\tau$ is the time of first return on the section $S_\psi$, then
$\vert a_{\tau^n(p)} \vert > C \lambda^n$.
The flow returns in uniformly bounded time on $S_\psi$, so the differential of $X^t_\psi$ is uniformly bounded between two crossings on the section.
We deduce that
for all $p \in P_\psi$ and $t \geq 0$, up to change the constant $C>0$, we have
$\vert a_{t} \vert > C \lambda^t$.
So the differential of $X^t_\psi$ expands exponentially the area of $E^\cu$ for all $t \geq 0$.
We similarly show the existence of a plane field $E^\cs$ on $\Lambda'_\psi$, invariant by the differential of the flow of $X_\psi$, transverse to $E^\cu$ and whose area is exponentially contracted by $X^t_\psi$ for $t \geq 0$.
It follows that $\Lambda'_\psi$ is a hyperbolic set (of index $(1,1)$) for the flow of $X_\psi$ in $P_\psi$, from Lemma~\ref{lem: expanded plane field criteria hyperbolic set}.
\end{proof}

The existence of a filtration adapted to $\gA_\psi \, \cup \, \gR_\psi \, \cup  \, \Lambda'_\psi$ implies that the chain-recurrent set of $X_\psi$ in $P_\psi$ is contained in
$\gA_\psi \, \cup \, \gR_\psi \, \cup  \, \Lambda'_\psi$, hence hyperbolic.
It is clear that the dimension of the stable and unstable manifolds of these sets are of constant dimension 2.
To use Ma\~n\'e Criterion, it remains to show

\begin{claim}
    The stable and unstable manifolds of $\gA_\psi \, \cup \, \gR_\psi \, \cup  \, \Lambda'_\psi$ intersects transversally.
\end{claim}

\begin{proof}
We consider orbits $O$ which are not contained in $\Lambda_\psi$, where the intersection is naturally transverse (via lift in $P$).
Every such orbit intersects at most once $\Pin_\psi$.
We study the following disjoint cases.

\begin{enumerate}[leftmargin=*]
    \item Let $O$ be an $X_\psi$-orbit in $\cW^s_{X_\psi}(\gA_\psi) \cap \cW^u_{X_\psi}(\gR_\psi)$.
There is a first intersection point $p \in \Pin_\psi$, and a last intersection point $q = X_\psi^\tau(p) = f_{0, \psi}^n(p) \in \Pin_\psi$ in $O$.
The segment of $X_\psi$-orbit in the past of $p$ lift to a segment of $X$-orbit, so the direction tangent to $\cW^u_{X_\psi}$ on $p$ coincide with the direction tangent to $\cW^u_\psi$.
The trace on $\Pin_\psi$ is contained in the interior of $C^u_\psi$ according to Proposition~\ref{prop: parameters and cones}, Item~\ref{prop: parameters; it: directions}.
Similarly, the segment of $X_\psi$-orbit in the future of $q$ lift to a segment of $X$-orbit, so the direction tangent to $\cW^s_{X_\psi}$ on $q$ coincide with the direction tangent to $\cW^s_\psi$.
The trace on $\Pin_\psi$ is disjoint from the closure of $C^u_\psi$ according to Proposition~\ref{prop: parameters and cones}, Item~\ref{prop: parameters; it: directions}. 
Moreover, we know from Proposition~\ref{prop: parameters and cones}, Item~\ref{prop: parameters; it: Pin to Pin} that $(f_{0, \psi}^n)_*(C^u_\psi)$ is contained in the interior $C^u_\psi$.
As $f_{0, \psi}^n$ maps the trace of the direction tangent to $\cW^u_{X_\psi}$ on $p$ to the trace of the direction tangent to $\cW^u_{X_\psi}$ on $q$ to, it follows that $\cW^s_{X_\psi}$ is transverse to $\cW^u_{X_\psi}$ in $q$, hence everywhere in~$O$.

\item Let $O$ an $X_\psi$-orbit in $\cW^s_{X_\psi}(\Lambda'_\psi) \cap \cW^u_{X_\psi}(\gR_\psi)$.
There is a first intersection point $p \in S_{0, \psi}$.
Remark that this point necessarily belong to $S_\psi = S_{0, \psi} \ssm (\Pin \cap \cV_0)$.
Indeed, if it belongs to $(\Pin \cap \cV_0)_\psi$, then the orbit of $p$ has crossed $\Sigma_\psi$ before by Item~\ref{prop: parameters; it: v0} of Proposition~\ref{prop: parameters and cones}, hence it is not the first intersection point.
Moreover $p$ belongs to $K_\psi^+$.
The direction tangent to $\cW^u_{X_\psi}$ on $p$ coincide with the direction tangent to $\cW^u_\psi$.
The trace on $S_\psi$ of the direction tangent to $\cW^s_{X_\psi}$ on coincide with $F^s$ and is contained in the interior of $C^s_\psi$ from Proposition~\ref{prop: hyperbolic splitting on S_psi}, Item~\ref{prop: hyp splitting S_psi; it: cones}.
We know in addition that the closure of $C^s_\psi$ is disjoint from $\cW^u_\psi$ from Proposition~\ref{prop: parameters and cones}, Item~\ref{prop: parameters; it: directions}.
It follows that $\cW^s_{X_\psi}$ is transverse to $\cW^u_{X_\psi}$ in $p$, hence everywhere in $O$.

\item Let $O$ an $X_\psi$-orbit in $\cW^s_{X_\psi}(\gA_\psi) \cap \cW^u_{X_\psi}(\Lambda'_\psi)$.
There is a last intersection point $q \in S_{0, \psi}$.
Remark that this point necessarily belong to $S_\psi$, as for the previous item.
Moreover $q$ belongs to $K_\psi^-$.
The direction tangent to $\cW^s_{X_\psi}$ on $p$ coincide with the direction tangent to $\cW^s_\psi$.
The trace on $S_\psi$ of the direction tangent to $\cW^u_{X_\psi}(O)$ on $p$ coincide $F^u$ and is contained in the interior of $C^u_\psi$ from Proposition~\ref{prop: hyperbolic splitting on S_psi}, Item~\ref{prop: hyp splitting S_psi; it: cones}.
We know in addition that the closure of $C^u_\psi$ is disjoint from $\cW^s_\psi$ from Proposition~\ref{prop: parameters and cones}, Item~\ref{prop: parameters; it: directions}.
This shows that $\cW^s_{X_\psi}$ is transverse to $\cW^u_{X_\psi}$ in $p$, hence everywhere in $O$.

\end{enumerate}

They are the only cases to study, because $\cW^u_{X_\psi}(\gA_\psi) = \gA_\psi$ and $\cW^s_{X_\psi}(\gR_\psi) = \gR_\psi$ as they are respectively a union of basic attractors and repellers.
\end{proof}

This completes the proof of Theorem~\ref{thm: gluing theorem}.
\end{proof}

\subsubsection*{Generalization to partial gluing maps}

The paragraph is devoted to state an easy generalization of Theorem~\ref{thmintro: gluing theorem} for \say{partial gluing maps}, a similar construction but where we only want to glue a \emph{part} of the boundary components.
More precisely 

\begin{defi}[Partial gluing map] \label{def: partial gluing map}
A \emph{partial gluing map} of a block $(P,X)$ is an involution $\varphi$ of an union $\ppP$ of connected components of $\partial P$, such that $\varphi$ satisfies the condition of a gluing map on $\ppP$, i.e., it pairs \ccs{} of $\ppP$, maps $\Pout \cap \ppP$ on $\Pin \cap \ppP$ and $\cO_* \cap \ppP$ on $\cO_* \cap \ppP$ respecting the flow orientation.
We say equivalently that $(P,X,\varphi: \ppP \to \ppP)$ is a partial triple.
\end{defi}

We generalize the Gluing Theorem~\ref{thm: gluing theorem} with the following statement.
The proof is essentially the same, except that we deal with the hyperbolicity of the orbit of the new maximal invariant set of $(P/\varphi, X_\varphi)$ which is not the entire set of orbits.

\begin{thm} \label{thm: partial gluing theorem}
Let $(P,X, \varphi: \ppP \to \ppP)$ be a filled building block, equipped with a \emph{partial} \sqt{} gluing map.
Then there exists a partial triple $(P', X', \varphi')$ strongly isotopic to $(P,X,\varphi)$ such that $P'/ \varphi'$ is a smooth manifold with boundary equipped with a vector field $X'_{\varphi'}$ induced by $X'$, and the pair $(P'/\varphi', X'_{\varphi'})$ is a building block, whose boundary is $\pP \ssm \ppP$.
\end{thm}

\subsection{Criterion of transitivity}
\label{sec: proof; subsec: transitivity criterion}
Recall that a \bb{} $(P,X)$ is said to be \emph{transitive} if the maximal invariant set $\Lambda$ is a transitive set for the flow of $X$.
An Anosov vector field $Y$ on a closed manifold $\cM$ is said to be \emph{transitive} if $\cM$ is a transitive set for the flow of $Y$.
In this subsection, we give a criterion analogous to \cite[Proposition~1.6]{beguinBuildingAnosovFlows2017}, which allows us to determine the transitivity of an Anosov flow, or a \bb{}, obtained by gluing the boundary of a \bb{} under the hypothesis of Theorem~\ref{thm: gluing theorem}.
Let $(P,X)$ be a \bb{} and $\Lambda$ be the maximal invariant set of $X$ in $P$, and let $\Lambda_1, \dots, \Lambda_n$ be the decomposition of the non-wandering set of $\Lambda$ into \emph{basic pieces} (\cite{smaleDifferentiableDynamicalSystems1967}).

\begin{defi}[Graph of a triple] \label{def: graph g(P,X,varphi)}
Let $G = G(P,X,\varphi)$ be the directed graph having as vertex the basic pieces $i = \Lambda_i$ and a oriented edge $(i,j)$ if:
\begin{itemize}[label=--]
    \item $\cW^u(\Lambda_i) \cap \cW^s(\Lambda_j) \neq \varnothing$, or
    \item $\varphi(\cW^u(\Lambda_i)) \cap \cW^s(\Lambda_j) \neq \varnothing$, or
    \item $\cW^u(\Lambda_i) \cap \varphi(\cW^s(\Lambda_j)) \neq \varnothing$.
\end{itemize}
\end{defi}

We have the following criterion.

\begin{prop} \label{prop: transitivity criterion}
\mb{}
\begin{enumerate}
    \item If $(P_0, X_0, \varphi_0)$ and $(P_1, X_1, \varphi_1)$ are two strongly isotopic (partial) triples, then the graphs $G(P_0, X_0, \varphi_0)$ and $G(P_1, X_1, \varphi_1)$ are isomorphic.
    \item Let $(P,X,\varphi)$ be a triple such that $\varphi$ is a \sqt{} gluing map and $X_\varphi$ is an Anosov vector field on $P/\varphi$.
    If $G(P,X,\varphi)$ is \emph{strongly connected}\footnote{each pair of vertex can be connected by a oriented path of edges}, then $X_\varphi$ is transitive.
    \item Let $(P,X,\varphi)$ be a partial triple such that $\varphi$ is a \sqt{} gluing map and $(P/\varphi, X_\varphi)$ is a \bb{}.
    If $G(P,X,\varphi)$ is strongly connected, then $X_\varphi$ is transitive.
\end{enumerate}
\end{prop}

Notice that it is not assumed in the criterion of Proposition~\ref{prop: transitivity criterion} that the \bb{} is saddle or filled.

\begin{proof}
Up to make an orbit equivalence of one of the two blocks, $(P_0, X_0)$ and $(P_1, X_1)$ have a common \minc{} (Proposition~\ref{prop: isotopy vs orbit eq}, Item~\ref{prop: isotopy vs orbit eq; it: isotopy imples same minc}).
It follows that the maximal invariant set $\Lambda_{X_1}$ of $X_1$ and $\Lambda_{X_2}$ of $X_2$ is the same.
Moreover, there exists a \homeo{} $h: \pP_1 \ssm \cO_* \to \pP_0 \ssm \cO_*$ which maps the pair $(\cL^\iin_{X_1}, (\varphi_1)_*(\cL^\out_{X_1}))$ to the pair $(\cL^\iin_{X_0}, (\varphi_0)_*(\cL^\out_{X_0}))$ and the pair $(\cL^\out_{X_1}, (\varphi_1)_*(\cL^\iin_{X_1}))$ to the pair $(\cL^\out_{X_0}, (\varphi_0)_*(\cL^\iin_{X_0}))$ (Definition~\ref{def: strongly isotopic triples}, Item~\ref{def: strongly isotopic triples; it: lamination}).
This shows that the graphs $G(P_1, X_1, \varphi_1)$ and $G(P_0, X_0, \varphi_0)$ are isomorphic and the first item is true.

Let us show the second item.
Let $P_\varphi:= P/\varphi$ be the quotient manifold, $\pi_\varphi \colon P \to P_\varphi$ be the projection, and $\Lambda_\varphi := \pi_\varphi(\Lambda)$.
It is a hyperbolic invariant compact for the flow of $X_\varphi$ in $P_\varphi$.
Note $\cL = \cL^\out \cup \cO_* \cup \cL^\iin$ the boundary lamination on $\pP$.
Let $\Lambda_1, \dots, \Lambda_n$ be the basic pieces of the Smale decomposition of the non wandering set of $\Lambda$.
These are the maximal transitive hyperbolic compact sets of $\Lambda$.
Let $\Lambda_{\varphi,i} = \pi_\varphi(\Lambda_i)$ be the projection of the basic pieces of $\Lambda$ onto $\Lambda_\varphi$.
These are transitive hyperbolic sets for the flow of $X_\varphi$.
The stable manifold of $\Lambda_{\varphi, i}$ contains the projection of
$\cW^s(\Lambda_i) \cup \varphi (\cW^s(\Lambda_i))$
and its unstable manifold contains the projection of $\cW^u(\Lambda_i)\cup \varphi (\cW^u(\Lambda_i))$.
The strong connectedness of the graph $G(P,X,\varphi)$ is equivalent to the fact that the $\Lambda_{\varphi,i}$ are all connected by a cycle: for all $i, j$, there exists a sequence $i_0 = i, \dots, i_n=j$ such that
$\cW^s(\Lambda_{\varphi, i_k}) \cap \cW^u(\Lambda_{\varphi, i_{k+1}}) \neq \emptyset$ ($i_{n+1} = i_0$).
We deduce that they form a single basic piece (see \cite{smaleDifferentiableDynamicalSystems1967}), in other words $\Lambda_\varphi$ is transitive.
Finally, let $O$ be an orbit of the flow of $X_\varphi$ which is not contained in $\Lambda_\varphi$.
Then
\begin{claim}
The stable manifold $\cW^s(O)$ intersects $\cL^\out_\varphi$ and the unstable manifold $\cW^u(O)$ intersects $\cL^\iin_\varphi$.
\end{claim}
Indeed, if $O$ is not contained in $\Lambda_\varphi$, then it intersects transversally $P^\out_\varphi$.
The surface $P^\out_\varphi$ is foliated by a pair $(f^s, f^u)$ of transverse foliations which are the intersections of the stable $\cF^s$ and unstable $\cF^u$ laminations of the Anosov vector field $X_\varphi$ on $P_\varphi$.
Moreover, $f^u$ contains the lamination $\cL^\out_\varphi$ and $f^s$ contains $\cL^\iin_\varphi$.
The stable manifold $\cW^s(O)$ intersects $P^\out_\varphi$ along a leaf $l \in f^s$.
The leaf $l$ intersects the closure of a \cc{} $R$ of the complementary $\Pout_\varphi \ssm (\cL^\iin_\varphi \cup \cL^\out_\varphi)$, i.e., the projection of $\Pout \ssm (\varphi_* \cL^\iin \cap \cL^\out)$.
By strong transversality of the gluing map $\varphi$, the closure of $R$ is a compact rectangle bounded by two horizontal segments of leaves of $\cL^\iin_\varphi$ and two vertical segments of leaves of $\cL^\out_\varphi$.
It follows that the segment of leaf of $l$ crosses $R$ horizontally, and therefore intersects $\cL^\out_\varphi$.
The proof for the unstable manifold of $O$ is symmetric.
Therefore, each orbit $O$ of $X_\varphi$ has its stable manifold intersecting $\cW^u(\Lambda_\varphi)$ and its unstable manifold intersecting $\cW^s(\Lambda_\varphi)$.
We conclude that the manifold $P_\varphi$ is a unique basic piece, and therefore the Anosov flow $X_\varphi$ on $P_\varphi$ is transitive.

An analogous argument works for the last Item, by considering the orbits in the new maximal invariant set of the \bb{}.
\end{proof}

\section{Building block with prescribed boundary lamination}
\label{sec: prescribed boundary lamination}

In \cite[Theorem~1.10]{beguinBuildingAnosovFlows2017}, it is shown that any Morse--Smale foliation (in other words, a \qms{} foliation without marked leaves) on the torus can be realized (up to topological equivalence) as a boundary foliation of an attracting transitive orientable \bby{} block $(P,X)$, in other words such that the maximal invariant set $\Lambda$ is a transitive attractor for the flow of $X$.
The building blocks considered in this section will be connected and orientable.

\begin{defi}[Topological equivalence of laminations] \label{def: topo eq lamination}
Two laminations $\cL_1$ and $\cL_2$ on oriented surfaces $S_1$ and $S_2$ are said to be topologically equivalent if there exists a homeomorphism $h: S_1 \to S_2$ which preserves the orientation and such that $h_* \cL_1 = \cL_2$.
\end{defi}

In this section we prove the following analogous main result.

\begin{prop} \label{prop: transitive block with prescribed lam}
Let $\cF_1$ and $\cF_2$ be two \qms{} foliations on an oriented torus $S_1$ and an oriented torus $S_2$ respectively, such that $\cF_1$ and $\cF_2$ have the same nonzero number of marked leaves.
Then there exists $(P,X)$ a transitive filled \bb{}, with $P$ oriented and connected, such that $\partial P$ is the union of two torus $T_1$ and $T_2$ \qt{} to $X$, and the boundary lamination $\cL_X$ restricted to $T_i$ extends into a foliation topologically equivalent to $\cF_i$.
\end{prop}

\begin{rmk} \label{rmk: precision orientation boundary}
Recall that an orientation of $P$ induces a canonical orientation on the boundary $\pP$.
The proposition states that there exists a homeomorphism $h_i : S_i \to T_i$ which maps $\cF_i$ to a foliation containing the boundary lamination and preserves the orientation for $i=1, 2$.
\end{rmk}

Recall that the boundary lamination $\cL_X$ of a filled block is, by definition, filling, so it extends to a foliation which is unique up to topological equivalence (the \ccs{} of the complementary of the lamination are strips, Definition~\ref{def: strip and filling lam}).
We will also show the main Theorem~\ref{thmintro: realize qms bifoliation}, which states that any pair of \qt{} foliations on the torus can be realized as the trace of the stable and unstable foliation of a transitive Anosov flow on a \qt{} embedded torus.

\subsubsection*{Section summary}
The section is organized as follows.

\begin{itemize}[leftmargin=*]

    \item In a preliminary Subsection~\ref{sec: prescribed boundary lam; subsec: gluing without cycle}, we give a simple result which allow us to build new \bbs{} from partial gluing of \bbs{}, without the condition of being \emph{filled} and equipped with a \emph{\sqt{}} gluing map. It is the particular case where the gluing pattern has \say{no cycle}.

    \item In Subsection~\ref{sec: prescribed boundary lam; subsec: combinatorial type}, we study \qms{} prefoliation (Definition~\ref{def: qms lam}).
    We associate to every \qms{} prefoliation a combinatorial type (Definition~\ref{def: qms combi type} and~\ref{def: qms combi type})
    It is a complete invariant of the topological equivalence class of a \qms{} foliation (Proposition~\ref{prop: qms foliation for given combi type}).
    
    \item In Subsection~\ref{sec: prescribed boundary lam; subsec: non-transitive block}, we show a preliminary result which states that any \qms{} foliation can be realize (up to topological equivalence) as the boundary foliation of a filled \emph{non-transitive} \bb{} (Proposition~\ref{prop: non-transitive block with prescribed lam}).
    We first construct a \bb{} $(P,X)$ whose boundary consist of the union of two quasi-transverse tori $T_1$ and $T_2$, containing periodic orbits whose number and orientation is compatible with the given combinatorial type and which constitute the boundary lamination on $T_i$, and a collection of tori $T^\iin$ and $T^\out$ transverse to the \vf{} $X$ (Lemma~\ref{lem: empty block}).
    We then glue attracting and repelling \bbs{} along the boundary tori $T^\iin$ and $T^\out$ with prescribed boundary foliation (Lemma~\ref{lem: BBY attractors with prescribed MS lam}).
    This gluing operation must be done in a way to induce foliations on the \qt{} tori whose combinatorial type is the one prescribed.
    We use properties of building blocks obtained by partial gluing with no cycle (Subsection~\ref{sec: prescribed boundary lam; subsec: gluing without cycle}).
    
    \item In Subsection~\ref{sec: prescribed boundary lam; subsec: transitivity by blow up} we show a general result which allows, from a non-transitive building block $(P,X)$ satisfying a condition on its Smale's graph, to create a filled transitive building block $(P',X')$ and such that the boundary laminations of $(P,X)$ and $(P',X')$ extends to topologically equivalent foliations.
    We use the \emph{Blow-up -- Excise -- Glue} surgeries coming from \cite[Section~8]{beguinBuildingAnosovFlows2017}.
    The idea is to create a collection of transverse boundaries by bifurcation \emph{Derived from Anosov} on periodic orbits, then excising small tubular neighborhoods whose boundaries are transverse to the flow, and this in each of the non-trivial basic pieces of $(P,X)$.
    We then glue these boundaries back together so as to create cycles in the graph $G$ associated to the triple (Definition~\ref{def: graph g(P,X,varphi)}).
    We then use the transitivity criterion (Proposition~\ref{prop: transitivity criterion}).
    We show that Proposition~\ref{prop: transitive block with prescribed lam} follows from this general result and the result of the previous subsection.
    
    \item In Subsection~\ref{sec: prescribed boundary lam; subsec: bifoliation}, we show Theorem~\ref{thmintro: realize qms bifoliation} which allows us to realize any pair of \qt{} foliations on the torus as the trace of the stable and unstable foliation of a transitive Anosov flow on a \qt{} embedded torus.
    It is a corollary of Proposition~\ref{prop: transitive block with prescribed lam} and of the Gluing Theorem~\ref{thm: gluing theorem}.
\end{itemize}

\subsection{Partial gluing without cycle}
\label{sec: prescribed boundary lam; subsec: gluing without cycle}

There is an easy way to construct blocks from partial gluing (Definition~\ref{def: partial gluing map}) of \bbs{} without the condition of being \emph{filled}, and equipped with a \emph{\sqt{}} gluing map.
This is the particular case where there is \emph{no cycle} in the gluing pattern, and the \cc{} of the boundary glued together are transverse to the vector field.
It is \cite[Proposition~1.1]{beguinBuildingAnosovFlows2017}, that we reformulate here.
The argument is then elementary, and relies essentially on the $\lambda$-lemma.

\begin{prop} \label{prop: block obtained by partial gluing without cycle}
Let $(P,X)$ be a \bb{} and $\varphi: \ppP \to \ppP$ a partial gluing map of $(P,X)$, such that:
\begin{enumerate}
    \item $\ppP$ is a union of \ccs{} of $\partial P$ transverse to $X$;
    \label{prop: partial without cycle; it: transverse boundary}
    \item each orbit of the flow of $X$ in $P$ intersects at most once $\ppP$ (we say that the partial gluing map $\varphi$ has \emph{no cycle});
    \label{prop: partial without cycle; it: without cycle}
    \item $\varphi_* \cL$ is transverse to $\cL$ on $\ppP$.
    \label{prop: partial without cycle; it: lam transverse}
\end{enumerate}
Then up to modify the partial triple $(P,X,\varphi)$ by strong isotopy, the vector field $X$ induces a \vf{} $X_\varphi$ on the manifold $P_\varphi:= P/\varphi$, such that the couple $(P_\varphi, X_\varphi)$ is a \bb{}.
\end{prop}

Before showing the proposition, let us state a number of properties for such \bb{} in the following lemma.
We express the boundary, maximal invariant set, and lamination of a block $(N,Y)$ obtain by a partial gluing of a block $(P,X)$ without cycle.

\begin{lem} \label{lem: block obtained by partial gluing without cycle}
Let $(P,X,\varphi:\ppP \to \ppP)$ be a partial triple without cycle.
Suppose the manifold $N := P/ \varphi$ is provided with a vector field $Y$ of class $\cC^1$ induced by $X$. Then
\begin{enumerate}
    \item \label{lem: block by partial gluing; it: boundary}
    \emph{(boundary)} $\partial N = \pP \ssm \ppP$,
    \item \label{lem: block by partial gluing; it: max inv}
    \emph{(maximal invariant set)} $\Lambda_Y = \pi_\varphi \l( \Lambda_X \cup \underset{t \in \R}{\bigcup} X^t (\cL_X \cap \varphi_* \cL_X) \r)$.
\end{enumerate}
Moreover, if $(N,Y)$ is a \bb{} and if $\cW^s_Y$, $\cW^u_Y$ denotes the stable and unstable manifolds of $\Lambda_Y$ and, $\cL_Y$ the boundary lamination of $(N,Y)$, then
\begin{enumerate}[resume]
    \item \label{lem: block by partial gluing; it: invariant manifolds}
    \emph{(invariant manifolds)}
    \begin{itemize}[--]
        \item $\cW^s_Y = \pi_\varphi \l( \cW^s_X \, \cup \, \underset{t \leq 0}{\bigcup} X^t \l( \varphi _* \cL^\iin_X \r) \r)$,
        \item $\cW^u_Y = \pi_\varphi \l( \cW^u_X \, \cup \, \underset{t \geq 0}{\bigcup} X^t \l( \varphi_* \cL^\out_X \r) \r)$;
    \end{itemize}
    \item \label{lem: block by partial gluing; it: lamination}
    \emph{(boundary lamination)} $\cL_Y = \cL_X \cup \foutin{}_* (\varphi_* \cL^\out_X) \cup \foutin^\inv{}_* (\varphi_*\cL^\iin_X)$.
\end{enumerate}
\end{lem}
For simplicity, we omit the domain where the previous expressions make sense. In other words
\begin{itemize}[--]
    \item $\varphi (*)$ means $\varphi (* \cap \ppP)$,
    \item $\foutin (*)$ means $\foutin (* \cap \Pin \ssm \cL^\iin))$,
    \item $\foutin^\inv (*)$ means $\foutin (* \cap (\Pout \ssm \cL^\out))$, etc.
\end{itemize}

\begin{proof} $\, $
\begin{enumerate}
    \item It is obvious.
    \item It is sufficient to see that if $\tilde p \in \cL_X \cap \Pin = \cL^\iin_X$, then the positive orbit by the flow $X^t$ is defined for all $t\geq 0$, and if $\tilde p \in (\varphi_* \cL_X) \cap \Pin$, then $\tilde q = \varphi (\tilde p) \in \cL_X \cap \Pout = \cL^\out_X$ and the negative orbit of $\tilde q$ by the flow $X^t$ is defined for all $t \leq 0$.
    We deduce the inclusion $\supset$ by projection in $N = \pi_\varphi (P)$.
    Suppose now that there is an orbit $\gamma$ of $\Lambda_Y$ which is not in $\pi_\varphi(\Lambda_X)$.
    Such an orbit crosses transversely the surface $\pi_\varphi(\ppP)$ in a point $p$.
    Let $\tilde p$ be the lift of $p$ in $\Pin \cap \ppP$.
    By the no cycle assumption, the positive $X$-orbit of $\tilde p$ never intersects $\ppP$ again. 
    If it intersects $\partial P$, by projection in $N$ we conclude the same for $\gamma$, which contradicts that $\gamma \in \Lambda_Y$.
    Hence $\tilde p \in \cL^\iin_X$, and the positive $X$-orbit of $\tilde p$ is the positive $Y$-orbit of $p$ (up to trivial lift in $N$). 
    Analogously, if $\tilde q = \varphi (\tilde p)$ is the lift of $p$ in $\Pout \cap \ppP$, the negative $X$-orbit of $\tilde q$ cannot intersects $\partial P$, hence belongs to $\cL^\out_X$, and the negative $X$ orbit of $\tilde q$ is the negative $Y$ orbit of $p$ (up to trivial lift in $N$). 
    
    \item It is sufficient to use the expression of the maximal invariant set from the previous item.
    \item It is sufficient to use the expression for invariant manifolds from the previous item.\qedhere
\end{enumerate}
\end{proof}

\begin{proof}[Proof of Proposition~\ref{prop: block obtained by partial gluing without cycle}]
Item~\ref{prop: partial without cycle; it: transverse boundary} ensures that $\varphi$ is a dynamical gluing map and $P_\varphi$ is a compact manifold of dimension 3 with a vector field $X_\varphi$ induced by $X$.
Let $(N,Y) := (P_\varphi, X_\varphi)$ and denote $\Lambda_Y$ the maximal invariant set.
Let $p \in \Lambda_Y$.
Then either $p \in \pi_\varphi(\Lambda_X)$, or the orbit of $p$ by the flow of $Y$ transversely intersects the set $\pi_\varphi (\ppP) =: (\ppP)_\varphi$ (Lemma~\ref{lem: block obtained by partial gluing without cycle}).
By Item~\ref{prop: partial without cycle; it: without cycle} of Proposition~\ref{prop: block obtained by partial gluing without cycle}, this intersection is unique.
In the first case, we have a decomposition of the tangent space given by the projection of the hyperbolic splitting of $\Lambda_X$ for $X$.
In the second case, $p \in \Lambda_Y \cap (\ppP)_\varphi$, and let $\tilde p \in \Pin \cap \ppP$ be a lift of $p$ in $P$.
Then $\tilde p \in \cL^\iin_X \cap \varphi_* \cL^\out_X$ (Lemma~\ref{lem: block obtained by partial gluing without cycle}, Item~\ref{lem: block by partial gluing; it: max inv}) and we define the following tangent planes on $\tilde p$:
\begin{equation}\label{eq: E^cu E^cs on partial_1P}
   E^\cs(\tilde p) = T \cL_X^s \oplus \R. X, \qq \mb{ and } \qq E^\cu(\tilde p) = T \varphi_* \cL_X^u \oplus \R. X.
\end{equation}
We obtain by projection in $P_\varphi$ the data of two tangent planes $E^\cu_Y$, $E^\cs_Y$ on the set $\Lambda_Y \cap (\ppP)_\varphi$ tangent to the vector field $Y$.
By Item~\ref{prop: partial without cycle; it: lam transverse}, $\cL^\iin$ is transverse to $\varphi_* \cL^\out$ on $\ppP$ so the planes $E^\cu$ and $E^\cs$ are transverse to each other.
We push these tangent planes by the flow of $X_\varphi$ on $\Lambda_Y \ssm (\Lambda_X)_\varphi$ in the following way.
If $q \in \Lambda_Y \ssm (\Lambda_X)_\varphi$, there exists a unique $p = Y^t(q) \in (\ppP)_\varphi$ and we define
\begin{equation} \label{eq: E^cu E^cs general}
    E_Y^\cu(p) =Y^{-t}_* E_Y^\cu(Y^t(q)), \qq \mb{ and }  \qq E_Y^\cs(p) =Y^{-t}_* E_Y^\cs(Y^t(q)).
\end{equation}
This defines an $Y^t$-invariant collection of tangent planes of $TN$ on $\Lambda_Y \ssm (\Lambda_X)_\varphi$.
Let us show that the area of $E^\cu_Y$ is exponentially expanded by $Y^t$ for $t \geq 0$.
Let $p \in (\ppP)_\varphi \cap \Lambda_Y$, and $v \in E^\cu_Y(p)$ be a vector tangent to $E^\cu$ and transverse to $X$ (uniformly in $p$).
Let $\tilde p \in \Pin \cap \ppP$ and $\tilde v \in E^\cu(\tilde p)$ be lifts in $P$.
Note that $\tilde p \in \cL^\iin$.
By definition of $E^\cu(\tilde p)$ in Equation~\eqref{eq: E^cu E^cs on partial_1P}, $\tilde v$ is tangent to $\varphi_* \cL^\out$ and transverse to $\cL^\iin$.
The positive $X$-orbit of $\tilde p$ accumulates on the maximal invariant set $\Lambda$, hence enters in uniformly bounded time a \nbh{} $\cU$ of $\Lambda$ where the flow is hyperbolic (we extend the hyperbolic splitting by continuity on $\cU$) and the vector $\tilde v$ is (uniformly) transverse to the weak stable direction $E^\ss \oplus \R.X$ in $\tilde p$.
According to the $\lambda$-lemma, $X^t$ exponentially expands the norm of the vector $\tilde v$ for $t \geq 0$.
Projecting into $P_\varphi$ it follows that $Y_t$ exponentially expands the norm of $v$ for $t \geq 0$.
Let us now look at the effect on the norm of the vector pushed backward by the flow of $Y$.
The vector $\varphi_* (v)$ is tangent to $\cL^\out$, so to $\cW^u_X$ in $\varphi (\tilde p) \in \Pout$.
Note that $\varphi (\tilde p) \in \cL^\out$.
The negative orbit of $\varphi (\tilde p)$ accumulates on the maximal invariant set $\Lambda$, hence enters in uniformly bounded time the \nbh{} $\cU$, and the vector $\tilde v$ is transverse to the weak stable direction $E^\ss \oplus \R.X$ in $\varphi (\tilde p)$.
According to the $\lambda$-lemma, $X^t$ exponentially contracts the norm of the vector $\tilde v$ for $t \leq 0$.
Projecting into $P_\varphi$ it follows that $Y^t$ exponentially contracts the norm of $v$ for $t \leq 0$.
Since $Y^t$ acts isometrically on the direction $\R. Y \subset E^\cu$, we deduce that the area of the plane field $E^\cu_Y$ is exponentially expanded by $Y^t$ for $t \geq 0$.
We show in the same way that the area of $E^\cs_Y$ is exponentially contracted in the future, and we conclude with Lemma~\ref{lem: expanded plane field criteria hyperbolic set}.
\end{proof}

\subsection{Combinatorial type of \qms{} prefoliation}
\label{sec: prescribed boundary lam; subsec: combinatorial type}

We refer to Subsection~\ref{sec: preli; subsec: boundary lamination} for general facts about \qms{} laminations and boundary lamination of \bb{}.
Let $\cL$ be a \qms{} lamination on a closed surface $S$, $\Gamma$ the compact leaves and $\Gamma_*$ the marked compact leaves with $\cL$, and $S \ssm \Gamma_* = S^\iin \cup S^\out$ the (in,out)-splitting of $\cL$ (Definition~\ref{def: qms lam}).

\begin{defi}[Dynamical orientation] \label{def: dynamical orientation compact leaf}
There exists an orientation of the elements of $\Gamma$ such that
the holonomy of each compact leaf is expanding on $S^\iin$ and contracting on $S^\out$.
We say that this is a \emph{dynamical orientation} of $\Gamma$.
We say that the lamination $\cL$ equipped with a dynamical orientation of the compact leaves $\Gamma$ is a \emph{dynamically oriented lamination}.
\end{defi}

\begin{rmk} \label{rmk: contracting and dilating dyn orientation}
In the case where the set $\Gamma_*$ is empty,
we have either $S = S^\iin$, or $S= S^\out$.
In the first case, the dynamical orientation is such that the holonomy of all compact leaves is expanding.
In the second case, the dynamical orientation is such that the holonomy of all compact leaves is contracting.
\end{rmk}

\begin{rmk} \label{rmk: canonical dynamical orientation block}
Such a lamination can a priori contain isolated compact leaves, and there are then several dynamical orientations.
In the following two cases, there is a canonical dynamical orientation.
\begin{itemize}
    \item If $\cL$ is a filling lamination (Definition~\ref{def: strip and filling lam}) or a foliation, then the dynamical orientation of the compact leaves $\Gamma_*$ is unique, because each compact leaf $\gamma$ is accumulated on both sides by non-compact leaves of $\cL$.
    
    \item If $\cL = \cL_X$ is the boundary lamination of a \bb{} $(P,X)$, then it has a canonical dynamical orientation.
    Recall that $\cL = \cL^\iin \cup \cO_* \cup \cL^\out$, with $\cL^\iin = \cW^s \cap \Pin$, and $\cL^\out = \cW^u \cap \Pout$.
    If $\gamma$ is a compact leaf of $\cL$, then we have the following situations:

    \begin{itemize}
    \item  $\gamma = \cO_i \in \cO_*$ is a periodic orbit of $X$, and it is oriented by the flow.
    
    \item $\gamma \in \cL^\iin$, and it is a circle transverse to the orbits of $X$ contained in a stable separatrix of a periodic orbit $O$.
    There exists a unique orientation of $\gamma$ such that $\gamma$ is freely homotopic to the oriented orbit $O$ in $\cW^s(O)$.
    
    \item $\gamma \in \cL^\out$, and it is a circle transverse to the orbits of $X$ contained in an unstable separatrix of a periodic orbit $O$.
    There exists a unique orientation of $\gamma$ such that $\gamma$ is freely homotopic to the oriented orbit $O$ in $\cW^u(O)$.
    \end{itemize}
    
    We refer to the proof of Lemma~\ref{lem: free separatrix}.
    This orientation is a dynamical orientation of $\gamma$ as a compact leaf of a \qms{} lamination in the sense of Definition~\ref{def: dynamical orientation compact leaf} (see the proof of Proposition~\ref{prop: block boundary lam are qms}).
\end{itemize}
\end{rmk}

\begin{lem} \label{lem: extend qms prefoliation by qms foliation}
Let $\cL$ be a dynamically oriented \qms{} prefoliation on a torus $T$.
Let $\Gamma_\cL$ be the compact leaves, $\Gamma_{\cL,*}$ the marked compact leaves, and $T^\iin_\cL \cup T^\out_\cL$ the (in,out)-splitting of $T$ for $\cL$.
There exists a dynamically oriented \qms{} foliation $\cF$ on $T$, which contains $\cL$ as a sublamination and such that, if we denote $\Gamma_\cF$ the compact leaves, $\Gamma_{\cF,*}$ the marked compact leaves, and $T^\iin_\cF \cup T^\out_\cF$ the (in,out)-splitting of $T$ for $\cF$, then
\begin{enumerate}
    \item $\Gamma_\cF = \Gamma_\cL$ and $\Gamma_{\cF,*} = \Gamma_{\cL, *}$ (as a set of oriented closed curves)\\[-0.8em]
    \item $T^\iin_\cL = T^\iin_\cF$ and $T^\out_\cL = T^\out_\cF$
\end{enumerate}
\end{lem}

\begin{proof}
According to Proposition~\ref{prop: complementary of a prefoliation}, a \cc{} of $T \ssm \cL$ is either homeomorphic to an annulus bounded by compact leaves (possibly the same one), or it is a strip (Definition~\ref{def: strip and filling lam}), in other words homeomorphic to $\R^2$ and bordered by two noncompact leaves asymptotic to each other at both end (Figure~\ref{fig: prefoliation}).
Let $C$ be a \cc{} of $T \ssm \cL$.
If $C$ is a strip, then we extend the lamination $\cL$ on $C$ by a trivial foliation $\cF$ on $C$ by (non-compact) leaves going from one end of the strip to the other.
If $C$ is an annulus, let $\gamma_1$ and $\gamma_2$ be the compact leaves of $\cL$ which border $C$, equipped with the dynamical orientation.
    The component $C$ is either contained in $T^\iin_\cL$ or in $T^\out_\cL$.
    Extend $\cL$ by a foliation $\cF$ on $C$, without compact leaf, such that each leaf of $\cF$ accumulates at both ends on a compact leaf $\gamma_i$, and such that the holonomy of the oriented compact leaf $\gamma_i$ of the boundary of $C$ for this foliation is expanding if $C$ is contained in $T^\iin_\cL$, contracting if $C$ is contained in $T^\out_\cL$.
The foliation $\cF$ thus constructed satisfies Lemma~\ref{lem: extend qms prefoliation by qms foliation}.
\end{proof}

Using the Poincar\'e--Hopf theorem, we have the following corollary.

\begin{coro}
\label{coro: homotopic compact leaf in qms prefoliation}
Let $\cL$ be a \qms{} prefoliation on the torus $T$.
The compact leaves of $\cL$ are non-contractible in $T$, and they are pairwise freely homotopic in $T$ as \emph{non-oriented} closed curves.
\end{coro}

In particular, there is a cyclic order of the elements of $\Gamma_\cL$ on $T$.
If $T$ is oriented, then giving a dynamically oriented compact \say{first} leaf determines a cyclic order.

\begin{defi}[Geometric enumeration] \label{def: geometric enumeration}
Let $\cL$ be a dynamically oriented prefoliation on the oriented torus $T$.
We will say that $\{ \gamma_0, \dots, \gamma_{n-1} \}$ is a \emph{geometric enumeration} of the compact leaves of $\cL$ if the leaf $\gamma_{i+1}$ is the next leaf of $\gamma_i$ for the order determined by the dynamical orientation of $\gamma_0$ and the orientation of~$T$.
\end{defi}

There are as many geometric enumerations as there are compact leaves.

\begin{defi}[Abstract combinatorial type] \label{def: abstract combi type}
An abstract combinatorial type $\sigma$ is a map
$$ \sigma \colon \Z/n\Z \; \longrightarrow \; \{ \ra{}, \la{} \} \times \{\ua{}, \da{} \} \times \{ \ra{}, \la{} \} $$
satisfying the following two conditions:
\begin{enumerate}[i)]
    \item $\sigma(0)= (*,\ua, *)$,
    \label{def: abstract combi type; it: value 0}
    \item $ \sigma(i) = (*,*,\ra{}) \iff \sigma(i+1) = (\la{},*,*) $ for all $i \in \Z/n\Z$.
    \label{def: abstract combi type; it: g and d}
\end{enumerate}
\end{defi}

\begin{nota}\label{nota: alphabet H and O}
Define the following alphabets by $\scA_H := \{ \ra{}, \la{} \}$ and $\scA_O :=\linebreak[4] \{ \ua{}, \da{} \}$.
If $\sigma \colon \Z/n\Z \to \scA_H \times \scA_O \times \scA_H$ is an abstract combinatorial type, we denote its three components by $\sigma = (\sigma_g, \sigma_o, \sigma_d)$, \footnote{For left, orientation, right (see Remark~\ref{rmk: combi type and holonomy}).}
and we denote $-\sigma := (\sigma_g, -\sigma_o, \sigma_d)$ where the sign $-$ means that we reverse the direction of the arrow.
\end{nota}

\begin{defi} [Combinatorial type of a \qms{} prefoliation] \label{def: qms combi type}
Let $\cL$ be a dynamically oriented \qms{} prefoliation on the oriented torus $T$.
Let $\Gamma= \{ \gamma_0, \dots, \gamma_{n-1} \}$ be a geometric enumeration of the compact leaves and $T \ssm \Gamma_* = T^\iin \cup T^\out$ be the (in,out)-splitting of $T$ for $\cL$.
Let the map
$$ \sigma = \sigma_\cL \colon \Z/n\Z \; \longrightarrow \;  \{ \ra{}, \la{} \} \times \{\ua{}, \da{} \} \times \{ \ra{}, \la{} \}  $$
defined by: 
\begin{enumerate}[i)]
    \item $ \sigma(i) = (*, \ua, *)$ if and only if the oriented leaf $\gamma_i$ is freely homotopic to the oriented leaf $\gamma_0$;
    \label{def: qms type; it: leaf orientation}
    \item  \label{def: qms type; it: holonomy left} $ \sigma(i) = (\ra, *, *)$ if and only if the left (local) side of $\gamma_i$ is contained in $T^\out$;\footnote{The dynamical orientation of $\gamma_0$ and the orientation of $T$ determines a left and a right side for each compact leaf.}
    \item  \label{def: qms type; it: holonomy right} $ \sigma(i) = (*, *, \la)$ if and only if the right (local) side of $\gamma_i$ is contained in $T^\out$.
  
\end{enumerate}
We say that $\sigma$ is the \emph{combinatorial type of $\cL$ on $T$}.
\end{defi}

\begin{example}
For example, a combinatorial type of the foliation of Figure~\ref{fig: qms lamination} is
$$ \la \ua \ra; \qq \la \da \ra; \qq \la \ua \la; \qq \ra \ua \ra $$
where we give the list of the values of $\sigma$ ordered from $i=0$ to $i=3$ following the figure from left to right.
\end{example}

\begin{rmk} \label{rmk: combi type and holonomy}
If the leaf $\gamma_i$ is accumulated on the left by the lamination $\cL$, then Item~\ref{def: qms type; it: holonomy left} is equivalent to saying that the holonomy of $\cL$ to the left of $\gamma_i$ is contracting for the dynamical orientation of $\gamma_i$.
Similarly, if the leaf $\gamma_i$ is accumulated on the right by the lamination $\cL$, then Item~\ref{def: qms type; it: holonomy right} is equivalent to saying that the holonomy of $\cL$ to the right of $\gamma_i$ is contracting for the dynamical orientation of $\gamma_i$.
This justifies the use of arrow symbols, and the notation $\scA_H$ for Holonomy and $\scA_O$ for Orientation.
\end{rmk}

\begin{claim} \label{claim: dictionary}
The map $\sigma_\cL$ is an abstract combinatorial type.
We have the following relations, via the identification $\Gamma \simeq \Z/n\Z$:
\begin{enumerate}
    \item $\Gamma_* = \{ \sigma_g = \sigma_d \}$,
    \item $\Gamma \cap T^\iin = \{\sigma=  (\la{}, *, \ra{}) \}$,
    \item $\Gamma \cap T^\out = \{\sigma= (\ra{}, *, \la{}) \}$.
\end{enumerate}
\end{claim}

\begin{proof}
It is clear that $\sigma_\cL(0) = (*, \ua, *)$ by definition.
Let $A$ be a connected component of $T \ssm \Gamma_\cL$.
It is an annulus bounded by two compact leaves $\gamma_i$ and $\gamma_{i+1}$, such that $\gamma_{i+1}$ is to the right of $\gamma_i$ for the orientation of $A$ induced by the orientation of $T$.
The annulus $A$ is either contained in $T^\iin$ or in $T^\out$.
It follows that $\sigma(i) = (*,*,\ra)$ if and only if $\sigma(i+1) = (\la,*,*)$.
So $\sigma_\cL$ is an abstract combinatorial type.
The following items follow directly from the definition.
\end{proof}

For an abstract combinatorial type $\sigma$, we will note $\Gamma_{*, \sigma} := \{ \sigma_g = \sigma_d \}$, and we will call it by use of language the set of \emph{marked leaves} of $\sigma$.
It follows from Definition~\ref{def: abstract combi type} that this set is of even cardinal.
Recall that a \qms{} lamination $\cL$ whose set of marked leaves is empty is said to be a Morse--Smale lamination (Remark~\ref{rmk: QMS lam and MS lam}).
Similarly, an abstract combinatorial type $\sigma$ whose set of marked leaves is empty will be a Morse--Smale combinatorial type.
If it is the combinatorial type of a Morse--Smale lamination $\cL$ on the torus $T$, then there are two possible cases:
\begin{itemize}
    \item either $\sigma(i) = (\ra, *, \la)$ for all $i$, this case corresponds to $T= T^\iin$,
    \item either $\sigma(i) = (\la, *, \ra)$ for all $i$, this case corresponds to $T= T^\out$.
\end{itemize}

\subsubsection*{The combinatorial type determines the foliation}
For what follows, we fix an orientation on the torus $\T^2$.

\begin{prop}
\label{prop: qms foliation for given combi type}
Any abstract combinatorial type $\sigma$ is the combinatorial type of a dynamically oriented foliation $\cF$ on the oriented torus $\T^2$, and $\cF$ is unique up to topological equivalence.
\end{prop}

\begin{proof}
Let $\sigma = (\sigma_g, \sigma_o, \sigma_d) \colon \Z/n\Z \to \scA_H \times \scA_O \times \scA_H$ be an abstract combinatorial type.
For any $i \in \Z/n\Z$, let $A_i = [0,1] \times \R/\Z$ be an annulus, and let $\gamma^i_0 = \{0\} \times \R/\Z \subset$ and $\gamma^i_1 = \{ 1 \} \times \R/\Z$ the boundaries of $A_i$.
We orient $\gamma^i_0$ by the increasing orientation of $\R/\Z$ if $\sigma_o(i)=\ua$, and by the decreasing orientation otherwise.
Similarly, we orient $\gamma^i_1$ by the increasing orientation of $\R/\Z$ if $\sigma_o(i+1)=\ua$, and by the decreasing orientation otherwise.
We have the following polychotomy:
\begin{enumerate}
\item $\sigma_o(i) = \sigma_o(i+1)$, and $\sigma_d(i) = \, \ra$
\item $\sigma_o(i) = \sigma_o(i+1)$, and $\sigma_d(i) = \, \la$
\item $\sigma_o(i) = -\sigma_o(i+1)$, and $\sigma_d(i) = \, \ra$
\item $\sigma_o(i) = -\sigma_o(i+1)$, and $\sigma_d(i) = \, \la$.
\end{enumerate}
Let $\cF_i$ be a foliation on $A_i$, with no compact leaves other than $\gamma_0^i$ and $\gamma_1^i$, and such that each noncompact half-leaf of $\cF_i$ accumulates at both ends on $\gamma^i_0$ or $\gamma^i_1$ in a expanding way if $\sigma_d(i)=\ra$, and contracting otherwise.
We obtain a cyclic sequence of closed annuli $A_i$ indexed by $i \in \Z/n\Z$, equipped with a foliation $\cF_i$.
We identify the oriented boundary $\gamma^i_1 = \{1\} \times \R/\Z$ of $A_i$ with the oriented boundary $\gamma^{i+1}_0 = \{0\} \times \R/\Z$ of $A_{i+1}$ via a \diff{} of the circle $\varphi_i$, equal to the identity in coordinate $\theta \in \R/\Z$ of $A_i$ and $A_{i+1}$.
We obtain a cyclic sequence of gluing \diffs{} $\varphi_i$ indexed by $i \in \Z/n\Z$.
The quotient $T = \bigcup_i A_i /\varphi$, where $\varphi$ denotes the product of $\varphi_i$, is a torus oriented by the canonical orientation on each $A_i$, and the foliations $\cF_i$ come together in a foliation $\cF$ on the torus $T \simeq \T^2$ of combinatorial type $\sigma$.

For uniqueness, we refer to the proof of \cite[Proposition~7.7]{beguinBuildingAnosovFlows2017} which deals with the case of Morse--Smale foliation for a given abstract Morse--Smale combinatorial type, i.e., when the set of marked leaves is empty.
Reducing to uniqueness on the closure of each connected component of the complement of compact leaves, it is clear that the proof is completely similar for \qms{} and for Morse--Smale foliations.
\end{proof}

\begin{rmk} \label{rmk: same type implies same filling lam}
If $\cL$ and $\cL'$ are two filling laminations on the oriented torus, it is equivalent to say that $\cL$ and $\cL'$ have the same combinatorial type to say that $\cL$ and $\cL'$ complete each other in foliations $\cF$ and $\cF'$ which are topologically equivalent (in the sense of Definition~\ref{def: topo eq lamination}).
This remark follows from Lemma~\ref{lem: extend qms prefoliation by qms foliation} and from Proposition~\ref{prop: qms foliation for given combi type}.
\end{rmk}

\begin{defi}[Restriction] \label{def: restriction type}
Let $\sigma$ be a combinatorial type on $\Z/n\Z$ and let $I = \{k_1, \dots, k_m \} \subset \Z/n\Z$ a subset of cardinal $m$, where the $k_i$ are cyclically ordered.
The restriction $\res{\sigma}{I}$ is the map
$\res{\sigma}{I} \colon \Z/m\Z \to \scA_H \times \scA_O \times \scA_H$ defined by
$\res{\sigma}{I}(i) = \sigma(k_i)$.
\end{defi}

A restriction of a combinatorial type is not a priori a combinatorial type.
This definition allows us to compare restrictions to subsets of the same cardinal of any two combinatorial types $\sigma$ and $\sigma'$.

\subsection{Non-transitive building block with prescribed boundary lamination}
\label{sec: prescribed boundary lam; subsec: non-transitive block}

In this section, we prove the following propositon.

\begin{figure}[h]
    \centering
    \vspace{-1em}
    \hspace*{-1.5em}\includegraphics[scale=0.65]{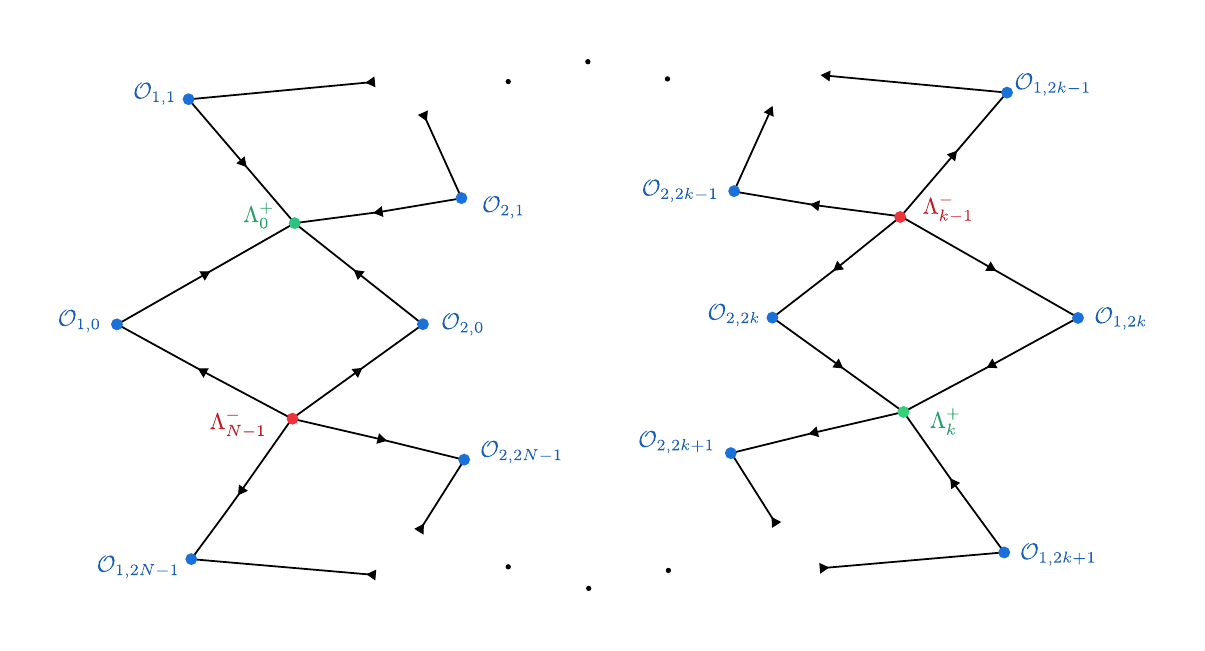}
    \hspace*{-1.5em}
    \vspace*{-1em}
    \caption{Smale's graph of $\Lambda_X$}
    \label{fig: smale graph non-transitive block}
\end{figure}

\begin{prop} \label{prop: non-transitive block with prescribed lam}
Let $\cF_1$ and $\cF_2$ be two foliations on the oriented torus $S_1$ and the oriented torus $S_2$ respectively, with the same non-zero number of marked leaves, i.e., $ \operatorname{Card} \Gamma_{*, \cF_1} = \operatorname{Card} \Gamma_{*, \cF_2} = 2N \neq 0$.
Then there exists $(P,X)$ a oriented connected \bb{} such that:
\begin{enumerate}
    \item \emph{(boundary)} \label{prop: non-transitive block; it: boundary}
    $\partial P$ is the union of two tori $T_1$ and $T_2$ quasi-transverse to~$X$.

    \item \emph{(boundary lamination)} \label{prop: non-transitive block; it: boundary lamination}
    The boundary lamination $\cL_{X}$ restricted to $T_i$ is a foliation topologically equivalent to $\cF_i$.
    
    \item \emph{(Smale's graph)} \label{prop: non-transitive block; it: smale graph}
    The basic pieces of the maximal invariant set $\Lambda_X$ are the $4N$ periodic orbits $\cO_{i,j}$ contained in $T_i \subset \pP$, $N$ attractors $\Lambda_k^+$, and $N$ repellers $\Lambda_k^-$, each containing infinitely many periodic orbits with negative multipliers.
    The Smale's graph of $\Lambda_X$ is given by Figure~\ref{fig: smale graph non-transitive block}.
\end{enumerate}
\end{prop}

For the rest of the proof, let $\Gamma_i$ be the set of compact leaves of $\cF_i$, $\Gamma_{*, i}$ the set of marked compact leaves.
We fix a geometric enumeration $\Gamma_i = \{ \gamma_{i, 0}, \dots, \gamma_{i,{n_i}} \}$ such that $\gamma_{i, 0} \in \Gamma_{i,*}$, in other words the first compact leaf is a marked leaf.
Let us denote $\sigma_i \colon \Z/n_i\Z \to \scA_H \times \scA_O \times \scA_H$ the corresponding combinatorial type of $\cF_i$.
Let us prove the following lemma.

\begin{figure}[htb]
    \centering
    \vspace*{-1.5em}
    \includegraphics[height=0.43\textheight]{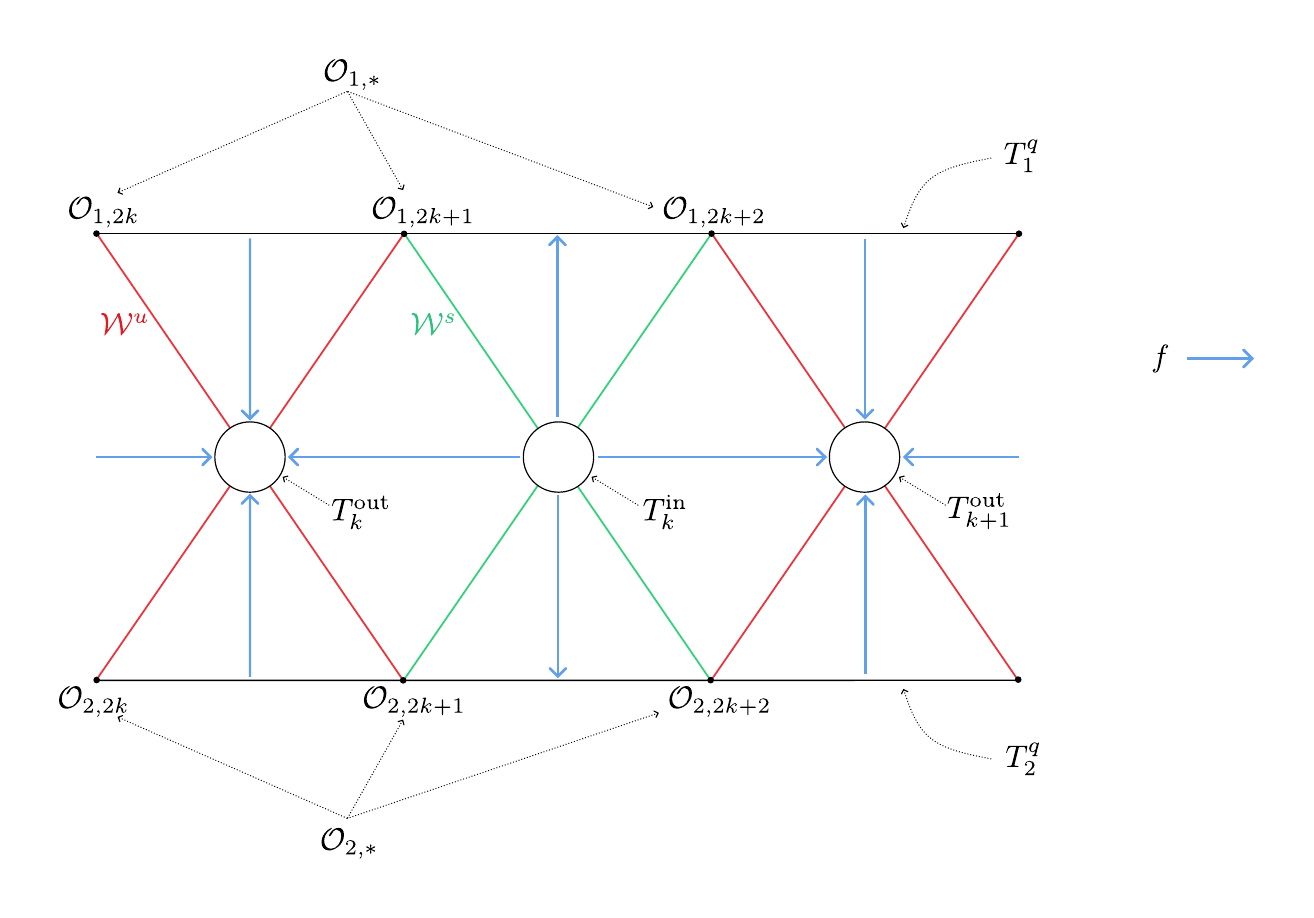}
    \vspace*{-1em}
    \caption{Block $(P, X)$ and crossing map $f$}
    \label{fig: empty block crossing map}
\end{figure}

\begin{lem} \label{lem: empty block}
There exists an oriented connected \bb{} $(P, X)$ which satisfies the following properties.

\begin{enumerate}
    \item \emph{(boundary)} \label{lem: empty block; it: boundary}
    The boundary $\partial P$ is the union of two tori $T_1^q$, $T_2^q$ quasi-trans\-verse to $X$ each containing a collection $\cO_{i,*} \subset \cO_*$ of $2N$ periodic orbits, $N$ tori $T^\iin_k \subset \Pin$ and $N$ tori $T^\out_k \subset \Pout$ transverse to $X$;
    
    \item \emph{(maximal invariant set)} \label{lem: empty block; it: max inv}
    The maximal invariant set $\Lambda$ is the union of the $4N$ isolated saddle periodic orbits $\cO_*$ of $X$ contained in the boundary~$\pP$;
    
    \item \emph{(boundary lamination on transverse tori)} \label{lem: empty block; it: lam t}
    The boundary lamination $\cL$ on $T^\iin_k$ and on $T^\out_k$ is the union of 4 compact leaves;
    
    \item \emph{(boundary lamination on \qt{} tori)} \label{lem: empty block; it: lam qt}
    The boundary lamination $\cL$ on $T_i^q$ is reduced to the periodic orbits $\cO_{i,*}$, and one of its combinatorial types coincides with the restriction of $\sigma_i$ on the set of marked leaves $\Gamma_{i,*}$;

    \item \emph{(crossing map)} \label{lem: empty block; it: crossing map}
    The crossing map $f: P^\iin \to \Pout$ pairs the \ccs{} of $\pP \ssm \cL$ according to the diagram given by Figure~\ref{fig: empty block crossing map}. 
\end{enumerate}
\end{lem}

\begin{rmk}
The manifold $P$ is a circle bundle on an annulus minus $2N$ disks, and the fiber is homotopic to a periodic orbit of $X$.
\end{rmk}

\begin{proof}
Let us note $\{p_{i,0}, \dots, p_{i,2N-1} \} \subset \Z/ 2N \Z$ the ordered set of marked leaves of $\sigma_i$, which we also denote (by use of notation) $\Gamma_{i,*}$.
To fix the ideas we suppose that $p_{i,0} = 0$, and that the left side of $\gamma_{i,0}$ is contained in $S_i^\out$ and the right side of $\gamma_{i,0}$ in $S_i^\iin$, which is equivalent to the following relation
\begin{equation} \label{eq: enumeration Fi}
    \sigma_i(p_{i,0}) = \sigma_i(0) = (\ra, \ua, \ra).
\end{equation}

Let $A = \R \times [-1,1]$ be the band of $\R^2$ with coordinates $(x,y)$.
Let $Y$ be a Morse--Smale gradient vector field on $A$ described by Figure~\ref{fig: morse smale field}, $2$-periodic along the variable $x$, i.e., $Y(x+2,y) = Y (x,y)$, and whose non-wandering set consists of

\begin{itemize}[--]
    \item hyperbolic fixed points of index $(1,1)$ on $(i, \pm 1)$ for $i \in \Z$,
    \item attracting fixed point on $(i + \frac{1}{2}, 0)$ for $i$ even integer,
    \item repelling fixed point on $(i + \frac{1}{2}, 0)$ for $i$ odd integer.
\end{itemize}

\begin{figure}[htb]
    \centering
    \vspace*{-1em}
    \includegraphics[height=0.33\textheight]{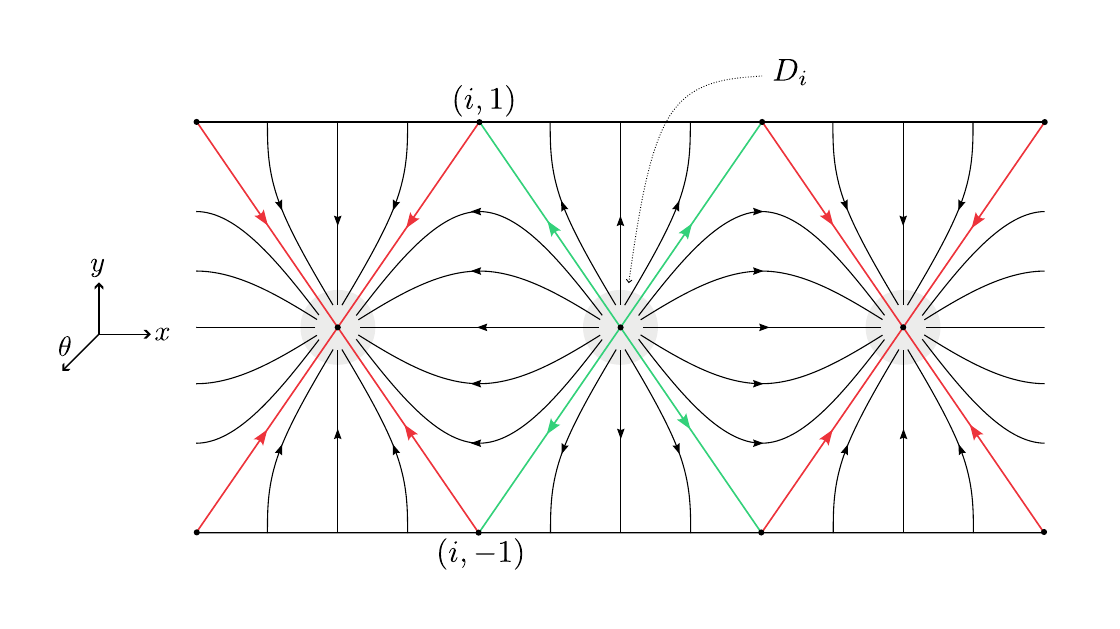}
    \vspace*{-1em}
    \caption{Morse--Smale vector field $Y$ on the band $\R \times [-1,1]$}
    \label{fig: morse smale field}
\end{figure}

We cut in $A$ disks $D_i = D( (i+\frac{1}{2}, 0), \epsilon)$ of center $(i+\frac{1}{2}, 0)$ and radius $\epsilon$, such that $Y$ is transverse to the boundary of $D_i$.
Let $\hat A := A \ssm (\bigcup_i D_i) \subset \R^2$, and let $\tilde P := \hat A \times \R/\Z$ be the product on the circle, given coordinates $(x,y,\theta)$, and oriented by the orientation of the canonical basis in these coordinates.
Let $\tilde X$ be the vector field on $\tilde P$ defined by
$$ \tilde X = Y + g(x,y) \partial_\theta,$$
where $g \colon \R^2 \to \R$ is a function of class $\cC^1$, which satisfies the following conditions:
\begin{enumerate}[i)]
    \item \label{it: support g}
    The support of $g$ is contained in disks $\D((i,\pm 1), \frac{1}{4})$ centered in $(i, \pm 1)$ of radius $\frac{1}{4}$, its value depends only on the radius $r$, and is monotonic in $r$.
    
    \item \label{it: periode abs g}
    The absolute value of $g$ is $1$-periodic along the variable $x$.
    
    \item \label{it: sign g}
    $g$ is constant equal to $\pm 1$ on the disk $\D((i,\pm 1), \frac{1}{8})$ centered at $(i, \pm 1)$ and of radius $\frac{1}{8}$:
    \begin{itemize}[--]
        \item $g=1$ on $\D((i,1), \frac{1}{8})$ if and only if $\sigma_1(p_{1,i})= (*,\da,*)$,
        \item $g=1$ on $\D((i,-1), \frac{1}{4})$ if and only if $\sigma_2(p_{2,i})= (*,\ua,*)$.
    \end{itemize}
\end{enumerate}

See Figure~\ref{fig: support g}.

\begin{figure}[htb]
    \centering
    \vspace*{-2em}
    \includegraphics[height=0.36\textheight]{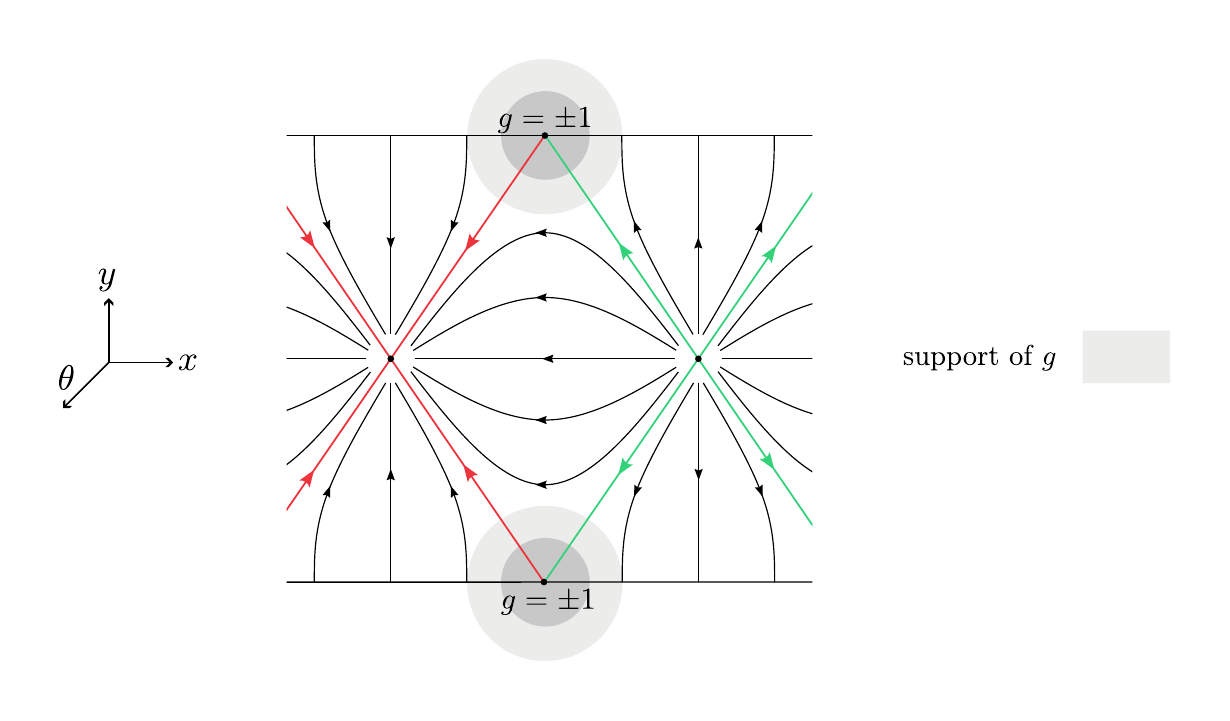}
    \vspace*{-2em}
    \caption{Support of $g$}
    \label{fig: support g}
\end{figure}

The choice on the sign of $g$ (Item~\ref{it: sign g}) will allow us to obtain the orientations of the periodic orbits compatible with the combinatorial type, which will be detailed at the end of the proof.
The reason for the mirror effect is that the canonical orientation of the boundary $\{y=1\}$ and $\{y=-1\}$ induced by the orientation of $\tilde P$ is \say{reversed}: the first one coincides with the orientation given by the basis field $(-\partial_x, \partial_\theta)$ and the other one by the basis field $(\partial_x, \partial_\theta)$.

The $g$-map is $2N$-periodic along the variable $x$ and the vector field $\tilde X$ quotients into a field $X$ on the manifold
$$ P := \tilde P / (x,y,\theta) \sim (x + 2N, y, \theta).$$
This manifold is compact, connected, with boundary, oriented by the basis field $(\partial_x, \partial_y, \partial_\theta)$.
Let us show that $(P, X)$ satisfies the lemma.

\begin{enumerate}[leftmargin=*]

    \item The boundary $\partial P$ is the union of two tori $T_1^q = \{ y=1\}$ and $T_1^q = \{ y=-1\}$, of $N$ tori $T^\out_k = \partial D_{2k} \times \R/\Z$ and of $N$ tori $T^\iin_k = \partial D_{2k+1} \times \R/\Z$ with $k \in \{0, \dots, N-1\}$.
    The vector field $X$ is transverse to $T^\iin_i$ and $T^\out_i$, and points inward along $T^\iin_i$ and outward along $T^\out_i$.
    Each torus $T_i^q$ contains a finite collection $\cO_{i,*}$ of hyperbolic periodic orbits of index $(1,1)$ of $X$ induced by the saddle hyperbolic fixed points $\Z \times \pm 1\}$ of the vector field $Y$.
    The vector field $X$ is transverse to the boundary on the complementary $T_i^q \ssm \cO_{i,*}$, alternately inward and outward on two connected components of $T_i^q \ssm \cO_{i,*}$ adjacent along a periodic orbit of $\cO_{i,*}$, in other words the tori $T_i^q$ are \qt{} to $X$ in the sense of Definition~\ref{def: qt surface}.
    Item~\ref{lem: empty block; it: boundary} is true.
    
    \item All orbits of $X$ which are not in $\cO_* = \cO_{1,*} \cup \cO_{2,*}$ intersect the boundary of $P$ in the future or in the past.
    It follows that the maximal invariant set $\Lambda$ coincides with the set $\cO_*$, and is therefore provided with a hyperbolic structure of index $(1,1)$ for $X$.
    We conclude that $(P, X)$ is a \bb{} in the sense of Definition~\ref{def: building block}.
    Item~\ref{lem: empty block; it: max inv} is true.
    
    \item Let $\cL$ be the boundary lamination of $(P, X)$.
    Each orbit of a point of $T_i^q \ssm \cO_*$ intersects a torus $T^\iin_k$ or $T^\out_k$.
    It follows that the boundary lamination $\cL$ restricted to $T_i^q$ is reduced to the periodic orbits $\cO_{i,*}$, the boundary lamination $\cL$ on a torus $T^\iin_i$ is the union of 4 compact leaves which are the intersection of the free stable separatrices of two consecutive periodic orbits of $T_1^q$ and two consecutive periodic orbits of $T_2^q$.
    Similarly, the boundary lamination $\cL$ on a torus $T^\out_i$ is the union of 4 compact leaves which are the intersection of the free unstable separatrices of two consecutive periodic orbits of $T_1^q$ and two consecutive periodic orbits of $T_2^q$.
    Item~\ref{lem: empty block; it: lam t} is true.
    
    \item The lamination $\cL$ is a dynamically oriented prefoliation (Remark~\ref{rmk: canonical dynamical orientation block}) and the boundary $\partial P$ is provided with the orientation induced by the orientation of the canonical basis $(\partial_x, \partial_y, \partial_\theta)$ in coordinates $(x, y, \theta)$ on $P$.

    We enumerate the periodic orbits of $T_1^q$ by $\cO_{1,i} := \{(i, 1)\} \times \R/\Z$ in coordinates $(x,y,\theta)$.
    Let $\sigma(\res{\cL}{T_1^q})$ be the combinatorial type of $\cL$ on $T_1^q$ associated with the first compact leaf $\cO_{1,0}$.
    The orientation of $T_1^q$ coincides with the orientation defined by the base field $(-\partial_x, \partial_\theta)$.
    The choice made between $+1$ or $-1$ when defining the map $g$ at the point $(i,1)$ determines the orientation of the vector field $\tilde X$ in the direction $\partial_\theta$ or $-\partial_\theta$ along the orbit $\cO_{1,i}$.
    By construction $g(0,1) = -1$, so the orbit $\cO_{1,0}$ is oriented in the decreasing direction of $\theta$, and its left side for the orientation of $T_1^q$ coincides with the side of the lower $x$.
    The vector field $X$ points outward on this side of $\cO_{1,0}$ (Figure~\ref{fig: morse smale field}), so $\sigma(\res{\cL}{T_1^q})(0) = (\ra, \ua, \ra)$.
    According to Equation~\eqref{eq: enumeration Fi}, we check that we have $\sigma(\res{\cL}{T_1^q})(0) = \sigma_1(p_{1,0})= \sigma_1(0)$.
    Finally, the enumeration of the collection $\cO_{1,*}$ coincides with the geometric enumeration determined by the choice of the first leaf $\cO_{1,0}$.
    Consequently, the sign of $g$ has been chosen so that $\sigma(\res{\cL}{T_1^q})(i) = \sigma_1 (p_{1,i}) = \res{\sigma_i}{\Gamma_{i,*}}(i)$ for all $i \in \Z/2N\Z$.
    
    The argument is the same for the combinatorial type of $\cL$ on $T_2^q$, except that the orientations are in \say{mirror}.
    Let us explain.
    We enumerate the periodic orbits of $T_2^q$ by $\cO_{2,i} := \{(i, -1)\} \times \R/\Z$ in coordinates $(x,y,\theta)$.
    Let $\sigma(\res{\cL}{T_2^q})$ be the combinatorial type of $\cL$ on $T_2^q$ associated with the given first compact leaf $\cO_{2,0}$.
    The orientation of $T_2^q$ coincides with the orientation defined by the base field $(\partial_x, \partial_\theta)$.
    The choice made between $+1$ or $-1$ when defining the map $g$ at the point $(i,-1)$ determines the orientation of the vector field $\tilde X$ in the direction $\partial_\theta$ or $-\partial_\theta$ along the orbit $\cO_{2,i}$.
    By construction $g(0,-1) = 1$, so the orbit $\cO_{2,0}$ is oriented in the increasing direction of $\theta$ and its left side for the orientation of $T_2^q$ therefore coincides with the side of the lower $x$.
    The vector field $X$ points outward on this side of $\cO_{2,0}$ (Figure~\ref{fig: morse smale field}), so $\sigma(\res{\cL}{T_2^q})(i) = (\ra, \ua, \ra)$.
    According to Equation~\eqref{eq: enumeration Fi}, we check that we have $\sigma(\res{\cL}{T_2^q})(0) = \sigma_2(p_{2,0})= \sigma_2(0)$.
    Finally, the enumeration chosen on the collection $\cO_{2,*}$ coincides with the geometric enumeration determined by the choice of the first leaf $\cO_{2,0}$.
    Consequently, the sign of $g$ has been chosen so that $\sigma(\res{\cL}{T_2^q})(i) = \sigma_2 (p_{2,i}) = \res{\sigma_i}{\Gamma_{i,*}}(i)$ for all $i \in \Z/2N\Z$.
    Item~\ref{lem: empty block; it: lam qt} is true.
    
    \item By construction.
\end{enumerate}

Finally, let us explain what to do if we change the choice made in Equation~\eqref{eq: enumeration Fi} for the value of the combinatorial type of the first marked leaf.
There are two possible cases (up to reversing $\sigma_1$ and $\sigma_2$):
\begin{itemize}
    \item $\sigma_1(0) = \sigma_2(0) = (\la, \ua, \la)$; then it is sufficient to construct the \vf{} $Y$ on the strip $A$ by exchanging the attracting fixed points and the repelling fixed points, which exchange the components of the entrance boundary and the components of the exit boundary on the two \qt{} tori~$T_{1,2}^q$.
    \item $\sigma_1(0) = (\la, \ua, \la)$ and $\sigma_2(0) = (\ra, \ua, \ra)$; then we change the definition of the map $g$ so that we have $g(i, 1) = 1 \iff \sigma_{1}(n-i) = (*,\ua,*)$ at Item~\ref{it: sign g}, which reverse the left and right sides of the first orbit on $T_1^q$ as well as the geometric enumeration, but not on $T_2^q$.\qedhere
\end{itemize}
\end{proof}

We say that a \bb{} $(P,X)$ is \emph{attracting} if the maximal invariant set $(P,X)$ is an attractor for the flow of $X$, and \emph{repelling} if the maximal invariant set is a repeller.
The rest of the proof consists of gluing together attracting transitive \bbs{} $(U_k^+, X_k^+)$ and repelling transitive \bbs{} $(U_k^-, X_k^-)$ along the boundary components $T^\out_k$ and $T^\iin_k$ of the block $(P, X)$ which are transverse to the \vf{} $X$, so as to obtain a new \bb{} whose boundary lamination on $T_1$ and $T_2$ has a prescribed combinatorial type.
We will use a result of \cite{beguinBuildingAnosovFlows2017} to obtain attracting and repelling \bby{} blocks realizing a given Morse--Smale boundary lamination on the torus.
We recall that a \bby{} block $(U,Y)$ is a \bb{} such that the vector field $Y$ is transverse to the boundary $\partial U$, or equivalently such that the collection of periodic orbits of $Y$ contained in $\partial U$ is empty (Remark~\ref{rmk: block and bby block}).
A Morse--Smale lamination is a \qms{} lamination such that the set of marked leaves is empty (Remark~\ref{rmk: QMS lam and MS lam}).

\begin{lem}[{\cite[Theorem~1.10]{beguinBuildingAnosovFlows2017}}]
\label{lem: BBY attractors with prescribed MS lam}
For any Morse--Smale foliation $\cF$ on the torus $\T^2$, there exists a connected transitive attracting orientable \bby{} block $(U,Y)$, and a \homeo{} \mb{$h: \partial U \to \T^2$} such that $h_* \cL_Y = \cF$, where $\cL_Y$ denotes the boundary lamination of $(U,Y)$.
Moreover, we can choose $(U,Y)$ such that the maximal invariant set $\Lambda_Y$ contains infinitely many periodic orbits with negative multipliers.\footnote{See \cite[Lemma~9.8]{beguinBuildingAnosovFlows2017} in the proof of \cite[Theorem~1.10]{beguinBuildingAnosovFlows2017}.}
\end{lem}

\begin{rmk}
\label{rmk: BBY repellers with prescribed MS lam}
We have a completely analogous result where $(U,Y)$ is a repelling block. We just need to reverse the vector field of the block given by Lemma~\ref{lem: BBY attractors with prescribed MS lam}.
\end{rmk}

\begin{proof}[Proof of Proposition~\ref{prop: non-transitive block with prescribed lam}]
Let $(P, X)$ be the oriented \bb{} given by \linebreak[4]Lemma~\ref{lem: empty block}.
Let $(P, X)$ be a geometric enumeration of the periodic orbits $\cO_{i,*} = \{ \cO_{i,0}, \dots, \cO_{i,2N-1}\}$ of $X$ contained in $T_i^q$, such that the combinatorial type of the lamination $\cL \cap T_i^q = \cO_{i,*}$ is equal to the restriction of $\sigma_i$ to the marked leaves $\Gamma_{i,*}$ of $\cF_i$ (Item~\ref{lem: empty block; it: lam qt}, Lemma~\ref{lem: empty block}).
Recall that $\Gamma_{i,*} = \{ p_{i,0}, \dots, p_{i, 2N-1}\}$ is the ordered collection of marked leaves with $\sigma_i$, with $p_{i,0} = 0$.

\subsubsection*{Step 1: Gluing attractors.}
Let $0 \leq k \leq N-1$ be an integer, and let $T^\out = T^\out_k \subset P^\out$ be a torus of $\partial P$ along which $Y$ is transverse outward.
According to Item~\ref{lem: empty block; it: crossing map} of Lemma~\ref{lem: empty block}, and Figure~\ref{fig: empty block crossing map}, there exists a \cc{} $A_1$ of $T_1^q \ssm \cL$ and a \cc{} $A_2$ of $T_2^q \ssm \cL$ such that $f (A_i) \subset T^\out$ for $i=1,2$.
The component $A_1 \subset P^\iin$ is an annulus bounded by the periodic orbits $\cO_{1,j}$ and $\cO_{1, j+1}$ of $\cO_{1,*}$, and the annulus $f(A_1)$ is bounded by the trace of the unstable manifolds of $\cO_{1,j}$ and $\cO_{1,j+1}$.
Let $B_1$ be the connected component of $S^\iin_1$ bounded by the marked leaves $p_{1,j}$ and $p_{1,j+1}$ of $\cF_1$.
Similarly, the component $A_2 \subset P^\iin$ is a annulus bounded by the periodic orbits $\cO_{2,l}$ and $\cO_{2, l+1}$ of $\cO_{2,*}$, and the oriented annulus $f(A_2)$ is bounded by the trace of the unstable manifolds of $\cO_{2,l}$ and $\cO_{2,l+1}$.
Let $B_2$ be the annulus of $S^\iin_2$ bounded by the marked leaves $p_{2,l}$ and $p_{2,l+1}$ of $\cF_2$.

\begin{claim} \label{claim: MS foliation F_k on T^out_k}
There exists a Morse--Smale foliation $\cF = \cF_k$ on $T^\out_k$ such that:
\begin{enumerate}
    \item $\cF$ is transverse to $\cL$,
    \item the compact leaves of $\cF$ are parallel to the compact leaves of $\cL$,
    \item \label{claim: F_k on T^out_k; it: homeo} 
    there exists homeomorphisms $h_i : B_i \to f(A_i)$ for $i=1,2$ which maps the foliation $\cF_i$ to the foliation $\cF$ and preserves orientation
\end{enumerate}
\end{claim}

\begin{figure}[htb]
    \centering
    \vspace*{-2em}
    \includegraphics[width=\textwidth]{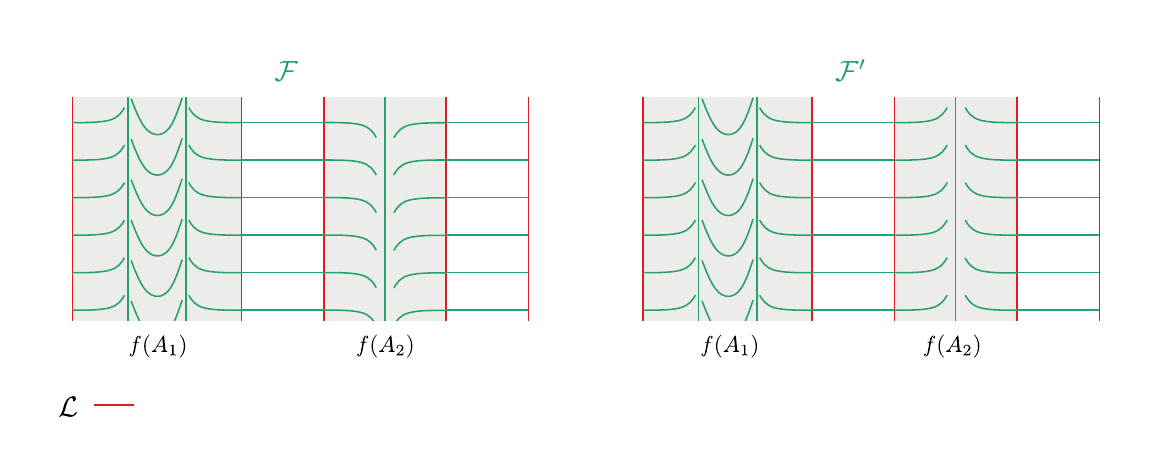}
    \vspace{-4em}
    \caption{Foliations $\cF$ and $\cF'$ satisfying Claim~\ref{claim: MS foliation F_k on T^out_k}}
    \label{fig: foliation concatenation type combi}
\end{figure}

\begin{proof}
The lamination $\cL$ on $T^\out = T^\out_k$ is the union of four compact leaves (Item~\ref{lem: empty block; it: lam t}, Lemma~\ref{lem: empty block}).
If $B_1$ and $B_2$ do not contain any compact leaf of $\cF_1$ and $\cF_2$, we choose for $\cF$ a foliation containing a single compact leaf $\gamma$ inside one of the \ccs{} of $T^\out \ssm \cL$ different from $f(A_1)$ and $f(A_2)$, and such that each noncompact half-leaf of $\cF$ accumulates on the compact leaf $\gamma$, and transversely intersects each compact leaf of $\cL$.
It satisfies the lemma.
Otherwise, let $n_1 \geq 1$ be the number of compact leaf of $\cF_1$ in $B_1$.
We construct a foliation $\cF$ on the annulus $f(A_1)$, transverse to the boundary of $f(A_1)$, which contains $n_1$ compact leaves parallel to the boundary of the annulus, which we list from left to right.
We arbitrarily choose an orientation on the first compact leaf.
We orient the following compact leaves in a way compatible with the combinatorial type of the foliation $\cF_1$ restricted to $B_1$.
The other leaves of $\cF$ are non-compact leaves such that each half-leaf $\cF$ accumulates on the compact leaves in a expanding way for this orientation.
We do the same on the annulus $f(A_2)$.
We complete the foliation $\cF$-foliation on $T^\out$ by the trivial foliation on the annuli $T^\out \ssm \bigcup_i f(A_i)$ and transverse to the boundary.
This is a Morse--Smale foliation transverse to $\cL$.
Analogously to the uniqueness result of Proposition~\ref{prop: qms foliation for given combi type}, there exists a \homeo{} $h_i : B_i \to f(A_i)$ which maps the restriction of $\cF$ to the annulus $f(A_i)$ to the restriction of $\cF_i$ on $B_i$, for $i=1,2$, and preserves orientation.
We refer to Figure~\ref{fig: foliation concatenation type combi}.
\end{proof}

We make an arbitrary choice of the orientation of the first compact leaf of $\cF$ in $f(A_i)$ for $i=1,2$.
Up to change this orientation (see Figure~\ref{fig: foliation concatenation type combi} where we reverse the first compact leaf of $f(A_2)$), we deduce :
\begin{rmk} \label{rmk: F_k' reverse sign}
Let $\cF$ be a foliation which satisfies Claim~\ref{claim: MS foliation F_k on T^out_k} above.
Then there exists a foliation $\cF'$ which satisfies Claim~\ref{claim: MS foliation F_k on T^out_k}, which coincides with $\cF$ along their respective compact leaves, and such that
\begin{enumerate}
    \item \label{rmk: F_k'; it: - sigma 2}
    either the dynamical orientations of the compact leaves of $\cF$ and $\cF'$ coincide on $f(A_1)$ and are opposite on $f(A_2)$;
    \item \label{rmk: F_k'; it: - sigma 1}
    either the dynamical orientations of the compact leaves of $\cF$ and $\cF'$ are opposite on $f(A_1)$ and coincide on $f(A_2)$;
    \item \label{rmk: F_k'; it: - sigma 1 and 2}
    either the dynamical orientations of the compact leaves of $\cF$ and $\cF'$ are opposite on $f(A_1)$ and on $f(A_2)$;
\end{enumerate}
\end{rmk}

By Lemma~\ref{lem: BBY attractors with prescribed MS lam}, there exists an orientable block $(U_k^+, X_k^+)$ such that the maximal invariant set $\Lambda_k^+$ is a transitive attractor which admits infinitely many periodic orbits with negative multipliers, and there exists a \diff{} \linebreak[4]$\varphi_k^+ \colon \partial U_k^+ \to T_k^\out$ such that $(\varphi_k^+)_*(\cL_k^+) =: \cF_k$ satisfies Claim~\ref{claim: MS foliation F_k on T^out_k}.
Let $$P' := \l( P \cup \l( \bigcup_{k=0}^{N-1} U_k^+ \r) \r) \big/ \varphi^+$$
where $\varphi^+$ is the product of \diffs{} $\varphi^+_k \colon \partial U_k^+ \to T^\out_k$.
Then:
\begin{claim} \label{claim: empty block attractors}
The vector fields $X$ and $X_k^+$ induce on the manifold $P'$ a vector field $X'$ such that $(P', X')$ is a connected oriented \bb{}.
It satisfies the following
\begin{enumerate}
    \item \emph{(boundary)} \label{claim: empty block attractors; it: boundary}
    the boundary $\partial P'$ is the union of two tori $T_1^q$, $T_2^q$ quasi-trans\-verse to $X'$ each containing a collection $\cO_{i,*} \subset \cO_*$ of $2N$ periodic orbits and $N$ tori $T^\iin_k \subset P'{}^\iin$ transverse to $X$;
    
    \item \emph{(Smale's graph)} \label{claim: empty block attractors; it: smale's graph}
    The basic pieces of the maximal invariant set $\Lambda'$ are $4N$ saddle periodic orbits $\cO_{i,j}$ contained in $T_i^q \subset \pP'$ and $N$ attractors $\Lambda_k^+$, each containing infinitely many periodic orbits with negative multipliers.
    The Smale's graph is given in Figure~\ref{fig: smale graph with attractors}

    \item \emph{(lamination)} \label{claim: empty block attractors; it: lam}
    the boundary lamination is $\cL' = \cL \cup \l( \bigcup_k (f^\inv)_* \l( \cF_k \ssm \cL \r) \r)$
    
    \item \emph{(combinatorial type)} \label{claim: empty block attractors; it: combi type}
    There exists a combinatorial type of $\cL'$ on $T_i^q$ such that
    $\sigma(\res{\cL'}{T_i^q})$ coincides with the restriction of $\sigma_i$ to the set $\Gamma_{i,*} \cup \{ \sigma_i =(\la{},*,\ra{}) \} $
    
    \item \emph{(crossing map)} \label{claim: empty block attractors; it: crossing map}
    The crossing map $f'$ of the block $(P', X')$ is the restriction of the crossing map $f$ of $(P, X)$.
\end{enumerate}
\end{claim}

\begin{figure}[htb]
    \centering
    \vspace*{-1em}
    \includegraphics[scale=0.66]{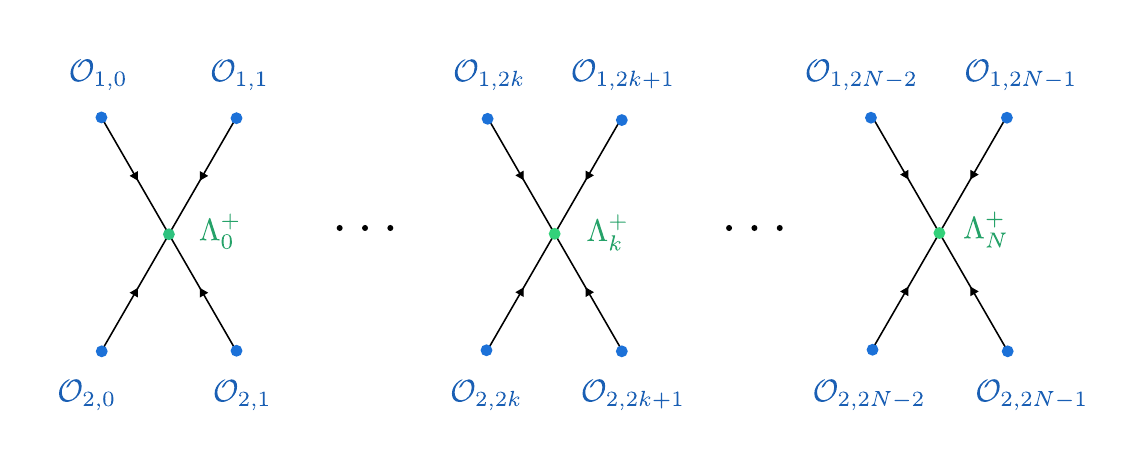}
    \vspace*{-1em}
    \caption{Smale's graph of the block $(P',X')$}
    \label{fig: smale graph with attractors}
\end{figure}

\begin{proof}
Consider the \bb{} given by the union $$(M,Y):= \textstyle{(P \bigcup_k U_k^+, X \bigcup_k X_k^+),}$$ equipped with the orientation of $P$ and the orientation on $U_k^+$ such that $\varphi_k^+$ reverses the orientation.
The \diff{} $\varphi^+$ is then a partial gluing of $(M,Y)$ which reverses the orientation, associated to a subset $\partial_1 M := \bigcup_k (\partial U^+_k \cup T^\out_k) $ of $\partial M$ which satisfies the hypotheses of Proposition~\ref{prop: block obtained by partial gluing without cycle}.
Indeed, $\partial_1 M$ is transverse to $Y$, each orbit of the flow of $Y$ intersects $\partial_1 M$ at at most one point, and $\varphi$ maps the boundary lamination of $(M,Y)$ to a transverse lamination by construction.
We deduce that the vector field $Y$ induces on the (compact, boundary oriented) manifold $P' = M/\varphi$ a vector field $X'$, such that $(P', X')$ is a \bb{}.
Items~\ref{claim: empty block attractors; it: boundary},~\ref{claim: empty block attractors; it: smale's graph},~\ref{claim: empty block attractors; it: lam},~\ref{claim: empty block attractors; it: crossing map} follow from Lemma~\ref{lem: block obtained by partial gluing without cycle}.
The maximal invariant set $\Lambda'$ coincide with the union of $\cO_*$, the attractors $\Lambda_k^+$ and the unstable separatrices (for the flow of $X'$) of the four periodic orbits contained in $T_1^q \cup T_2^q$ which intersect the torus $T^\out_k$.
The basic pieces of the maximal invariant set $\Lambda'$ are the periodic orbit $\cO_*$ and the attractors $\Lambda_k^+$ and the Smale's graph follow naturally from the gluing graph $G(M,Y,\varphi)$.

Let us show Item~\ref{claim: empty block attractors; it: combi type}.
We fix the geometric enumeration of compact leaves of $\cL'$ on $T_i^q$ such that the first leaf is the periodic orbit $\cO_{i,0} \in \cO_{i,*}$.
The laminations $\cL'$ and $\cL$ on $T_i^q$ coincide along the oriented marked leaves, in other words the periodic orbits $\cO_{i,*}$, and the (in,out)-splitting of $T_i^q \ssm \cO_{i,*}$ coincide for $\cL'$ and $\cL$.
It follows (Item~\ref{lem: empty block; it: lam qt}, Lemma~\ref{lem: empty block}) that $\sigma(\res{\cL'}{\cO_{i,*}})= \sigma(\res{\cL}{T_i^q}) = \res{\sigma_i}{\Gamma_{i,*}}$.
Let $0 \leq k \leq N-1$ be a fixed integer.
For $i=1, 2$, there exists a unique $A_i = A_{i,k}$ of $T_i^q \ssm \cO_{i,*}$ along which $X$ is transverse inward, and such that $f(A_i)\in T^\out_k$ (Figure~\ref{fig: empty block crossing map}).
On the interior of $A_i$, the lamination $\cL'$ coincides with the image $(f^\inv)_* \l( \cF_k \cap f(A_i) \r) $ (Item~\ref{claim: empty block attractors; it: lam}, Claim~\ref{claim: empty block attractors}) where $\cF_k$ is the lamination given by Claim~\ref{claim: MS foliation F_k on T^out_k}.
We refer to Figure~\ref{fig: gluing attractors}.
It follows from Claim~\ref{claim: MS foliation F_k on T^out_k}, Item~\ref{claim: F_k on T^out_k; it: homeo} that the combinatorial type of $\cL'$ restricted to the compact leaves on the interior of $A_i$ coincides with the restriction of $\pm \sigma_i$ on an interval of type $I_i = I_{i,k} = [[ p_{i,j} +1,\, p_{i,j+1} -1 ]]$ contained in the domain $\{ \sigma_i = (\la, *, \ra) \}$ where $j = j_k$ is an integer that depends on $k$, and the integers $j_k$ are consecutive in $\Z/2N\Z$.

\begin{figure}[htb]
    \centering
    \vspace*{-1em}
    \includegraphics[width=\textwidth]{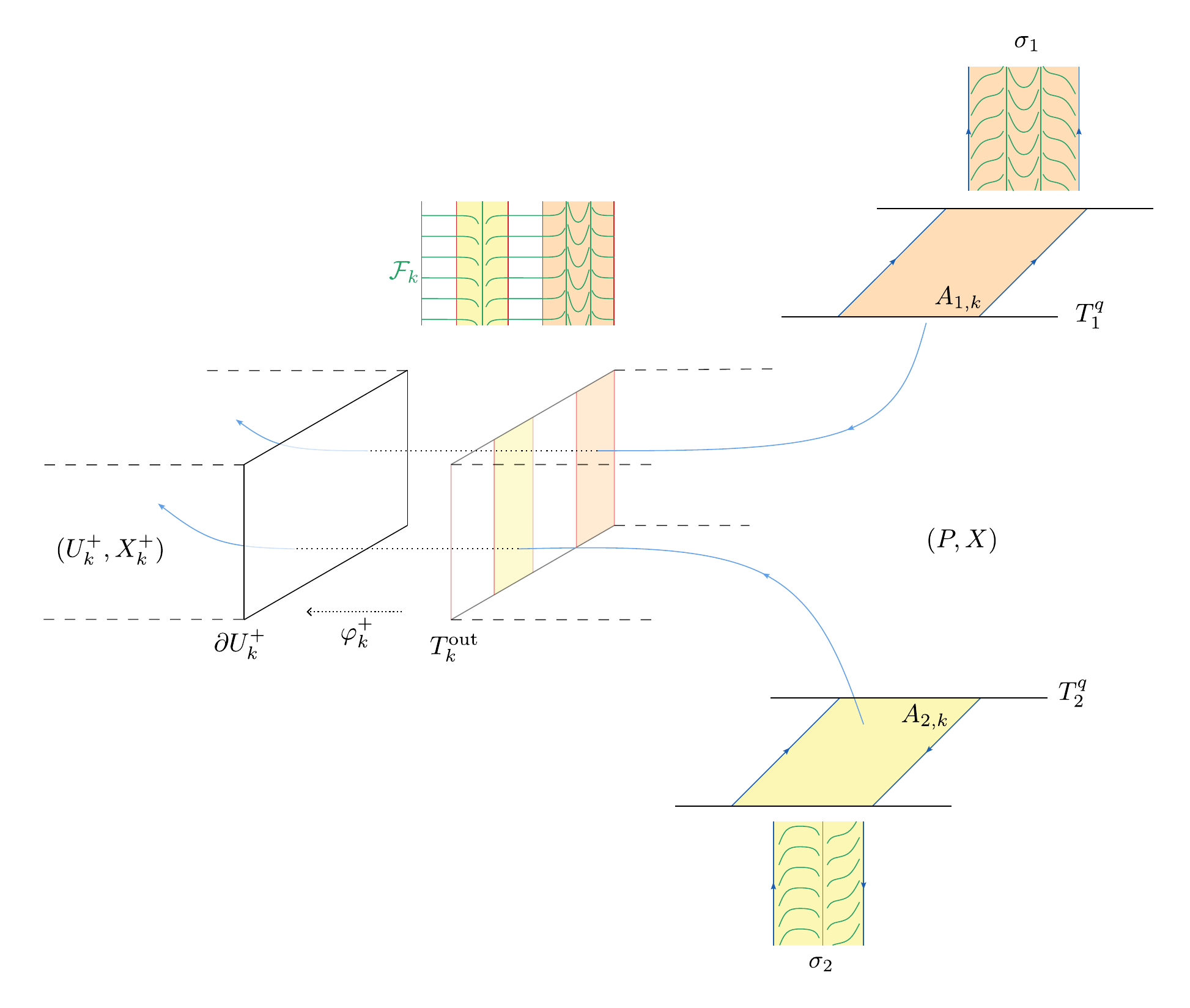}
    \vspace*{-1em}
    \caption{Gluing an attractor $(U_k^+, X_k^+)$ and induced foliation on the tori $T^i_q$ using the foliation of Figure~\ref{fig: foliation concatenation type combi}}
    \label{fig: gluing attractors}
\end{figure}

There are four possible cases, depending on the sign of $\sigma_1$ and the sign of~$\sigma_2$.
\begin{enumerate} \setcounter{enumi}{-1}
    \item \emph{($+\sigma_1$ and $+\sigma_2$):} \label{it: + 1 and + 2}
    In this case, we do nothing.
    \item \emph{($+\sigma_1$ and $-\sigma_2$):} \label{it: + 1 and - 2}
    In this case, we replace the foliation $\cF_k$ by the foliation $\cF_k'$ foliation given by Remark~\ref{rmk: F_k' reverse sign}, Item~\ref{rmk: F_k'; it: - sigma 2}.
    This has the effect of changing the sign of $\sigma_2$ but not of $\sigma_1$.
    \item \emph{($-\sigma_1$ and $+\sigma_2$):} \label{it: - 1 and + 2}
    In this case, we replace the foliation $\cF_k$ by the foliation $\cF_k'$ given by Remark~\ref{rmk: F_k' reverse sign}, Item~\ref{rmk: F_k'; it: - sigma 1}.
    This has the effect of changing the sign of $\sigma_1$ but not of $\sigma_2$.
    \item \emph{($-\sigma_1$ and $-\sigma_2$):} \label{it: - 1 and - 2}
    In this case, we replace the foliation $\cF_k$ by the foliation $\cF_k'$ given by Remark~\ref{rmk: F_k' reverse sign}, Item~\ref{rmk: F_k'; it: - sigma 1 and 2}.
    This has the effect of changing the sign of $\sigma_1$ and of $\sigma_2$.
\end{enumerate}
Repeating this argument for all $k= 0, \dots, N-1$, we show that the combinatorial type of $\cL'$ on each component $A_{i,k}$ of $T_i^q \ssm \cO_{i,*}$ coincides with the restriction of $\sigma_i$ to a corresponding interval $I_{i,k}$.
The union $\bigcup_k I_{i,k}$ is equal to the set $\sigma_i = (\la, *, \ra)$.
Putting this result together with the combinatorial type on the marked leaves, we deduce that the combinatorial type of $\cL'$ on $T_i^q$ coincide with the restriction of $\sigma_i$ to $\Gamma_{i,*} \cup \{ \sigma_i = (\la, *, \ra) \}$.
\end{proof}

\subsubsection*{Step 2: Gluing repellers.}
This second step consists of gluing transitive repelling blocks $(U_k^-, X_k^-)$ along the entrance tori $T^\iin_k$ on $(P', X')$, in order to induce the desired foliation on each connected component of $T_1^q \cap P'{}^\out$ and $T_2^q \cap P'{}^\out$.
The proof is essentially the same by replacing $(P, X)$ by $(P', X')$.
The only difference is that the lamination $\cL'$ along the torus $T^\iin_k$ is no longer the same as the lamination $\cL$. 
Recall that we have
$\cL' = \cL \bigcup_k( f^\inv)_* \l( \cF_k \ssm \cL \r)$ 
(Figure~\ref{fig: lamination T^in_k after gluing attractors}).

\begin{figure}[htb]
    \centering
    \vspace*{-1em}
    \includegraphics[width=\textwidth]{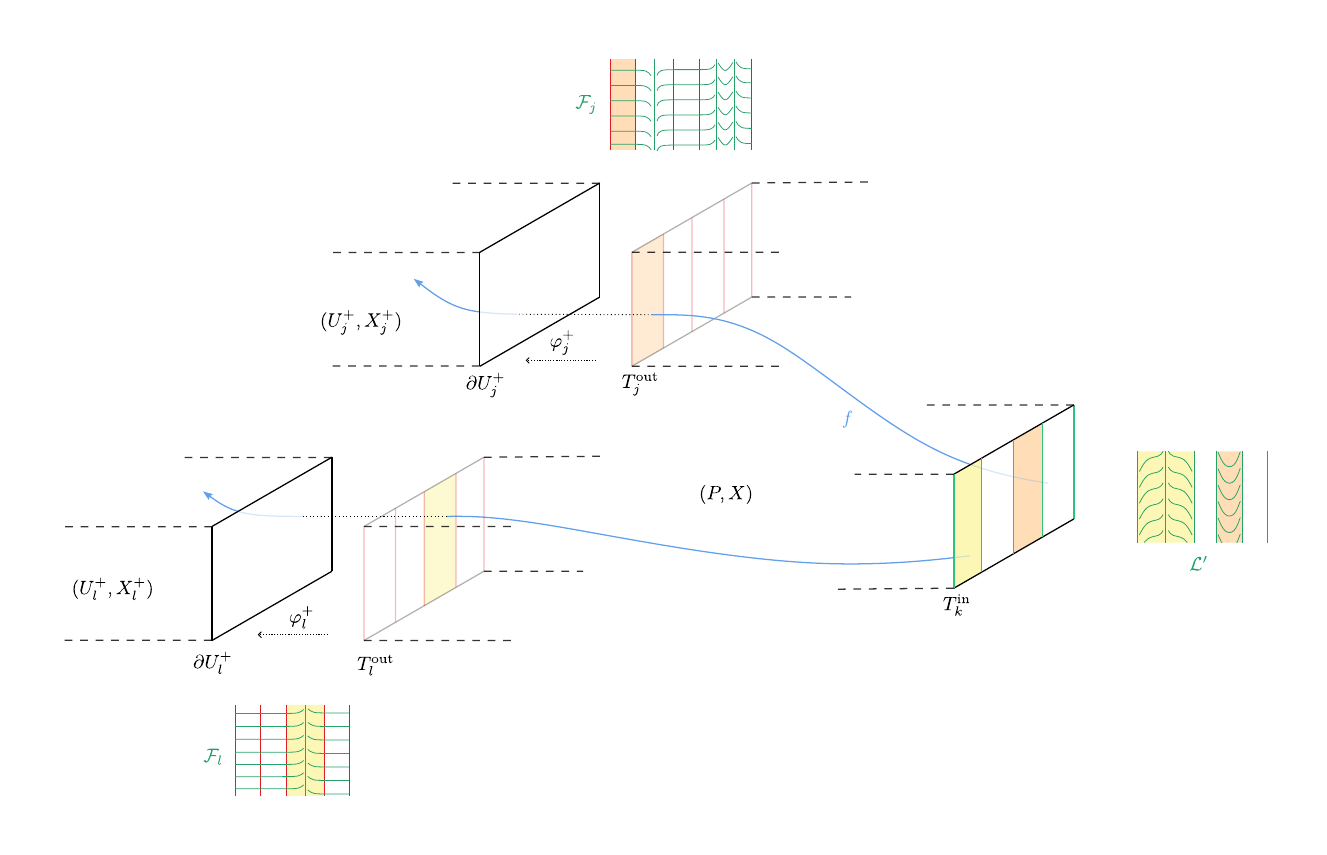}
    \vspace*{-1em}
    \caption{Lamination $\cL'$ on a torus $T^\iin_k$ induced by the gluing procedure of attractors $(U_l^+, X_l^+)$ and $(U_j^+, X_j^+)$}
    \label{fig: lamination T^in_k after gluing attractors}
\end{figure}

Let $0 \leq k \leq N-1$ be an integer, and let $T^\iin = T^\iin_k \subset P'{}^\iin$ a torus of $\partial P'$ along which $X'$ is transverse inward.
According to Item~\ref{claim: empty block attractors; it: crossing map} of Claim~\ref{claim: empty block attractors}, and Figure~\ref{fig: empty block crossing map}, there exists a \cc{} $A_1$ of $T_1^q \ssm \cL'$ and a \cc{} $A_2$ of $T_2^q \ssm \cL'$ such that $f'^\inv (A_i) \subset T^\iin $ for $i=1,2$.
The component $A_1 \subset P'{}^\out$ is an annulus bounded by the periodic orbits $\cO_{1,j}$ and $\cO_{1, j+1}$ of $\cO_{1,*}$.
Let $B_1$ be the connected component of $S^\out_1$ bounded by the marked leaves $p_{1,j}$ and $p_{1,j+1}$ of $\cF_1$.
Similarly, the component $A_2 \subset P'{}^\out$ is an annulus bounded by the periodic orbits $\cO_{2,l}$ and $\cO_{2, l+1}$ of $\cO_{2,*}$.
Let $B_2$ be the connected component of $S^\out_2$ bounded by the marked leaves $p_{2,l}$ and $p_{2,l+1}$ of $\cF_2$.

\begin{claim}
There exists a Morse--Smale foliation $\cG = \cG_k$ on $T_k^\iin$ such that:
\begin{enumerate}
    \item $\cG$ is transverse to $\cL'$;
    \item the compact leaves of $\cG$ are parallel to the compact leaves of $\cL'$;
    \item there exists homeomorphisms $g_i : B_i \to (f')^\inv(A_i)$ for $i=1,2$ which and maps the foliation $\cF_i$ on the foliation $\cF$ and preserves orientation.
\end{enumerate}
\end{claim}

\begin{proof}
The sets $f'^\inv(A_i) \subset T^\iin$ are annuli bounded by two compact leaves of the lamination $\cL'$ and disjoint from $\cL'$ on their interior.
The only difference with the proof of Claim~\ref{claim: MS foliation F_k on T^out_k} is that these are the only connected components of $T^\iin \ssm \cL'$.
Indeed, the lamination $\cL'$ is a foliation on $T^\iin \ssm \bigcup_i f'^\inv(A_i)$, because all points of $T^\iin \ssm \bigcup_i f'^\inv(A_i)$ have a positive orbit which intersects a torus $T^\out_k$, so accumulates on an attractor $\Lambda_k^+$ in the future (Figure~\ref{fig: empty block crossing map}).
We construct the foliation $\cG$ on the annuli $f'^\inv(A_i)$ in a similar way to the proof of Claim~\ref{claim: MS foliation F_k on T^out_k}: we place the compact leaves in the annuli $f'^\inv(A_i)$ parallel to the boundary and enumerated in the right way, orient them in a way compatible with the combinatorial type of $\cF_i$ on $B_i$ and complete with non-compact leaves that accumulate on the compact leaves oriented for the expanding orientation.
We complete the foliation thus constructed on $T^\iin$ by a transverse foliation on $\cL'$.
We can choose for example the foliation generated by the orthogonal vector field $Z^\perp$, for a certain metric, to a vector field $Z$ generating $\cL'$ on $T^\iin \ssm \bigcup_i f'^\inv(A_i)$.
By density of Morse--Smale vector fields, up to make a small perturbation on $T^\iin \ssm \bigcup_i f'^\inv(A_i)$, we obtain a Morse--Smale foliation $\cG$ on $T^\iin$, transverse to $\cL$, and which satisfies the fact.
\end{proof}

By Lemma~\ref{lem: BBY attractors with prescribed MS lam} and Remark~\ref{rmk: BBY repellers with prescribed MS lam}, there exists an orientable connected block $(U_k^-, X_k^-)$ such that the maximal invariant set $\Lambda_k^-$ is a transitive repeller which admits infinitely many periodic orbits with negative multipliers, and there exists a \diff{} $\varphi_k^- \colon \partial U_k^- \to T_k^\iin$ such that $(\varphi_k^-)_*(\cL_k^-) =: \cG_k$ satisfies Claim~\ref{claim: MS foliation F_k on T^out_k}.
Finally, we define
$$P'' := \l( P' \cup \l( \bigcup_k U_k^- \r) \r) \big/ \varphi^-$$
where $\varphi^-$ is the product of $\varphi^-_k \colon \partial U_k^- \to T^\iin_k$.
\begin{claim} \label{claim: empty block with attractors and repellers}
The \vfs{} $X'$ and $X_k^-$ induce on the manifold $P''$ a vector field $X''$ such that $(P'', X'')$ is a \bb{} which satisfies Proposition~\ref{prop: non-transitive block with prescribed lam}.
\end{claim}

\begin{proof}
The proof is the same as the proof of Claim~\ref{claim: empty block attractors}.
\end{proof}
This completes the proof of Proposition~\ref{prop: non-transitive block with prescribed lam}.
\end{proof}

\subsection{Transitivity by the \emph{Blow-up -- Excise -- Glue} surgery}
\label{sec: prescribed boundary lam; subsec: transitivity by blow up}

Proposition~\ref{prop: non-transitive block with prescribed lam} gives us the existence of a building block $(P,X)$ which satisfies Proposition~\ref{prop: transitive block with prescribed lam} for a given pair of \qms{} foliations, except transitivity.
We show the following general

\begin{prop} \label{prop: transitive block by blow up}
Let $(P,X)$ be a filled oriented connected \bb{}, we denote $\Lambda$ the maximal invariant set.
Assume that
for any basic piece $\Lambda_i$ of $\Lambda$, there exists a basic piece $\Lambda_k$ and $\Lambda_j$ containing infinitely many periodic orbits with negative multipliers, and such that there exists as oriented (possibly trivial) path in the Smale's graph connecting $\Lambda_j$ to $\Lambda_i$ and $\Lambda_i$ to $\Lambda_k$.

Then there exists a \emph{transitive} filled oriented connected block $(P', X')$ such that the boundary lamination $\cL_{X'}$ and the boundary lamination $\cL_X$ complete into topologically equivalent foliations.
\end{prop}

Let us first show that Proposition~\ref{prop: transitive block with prescribed lam} is an immediate corollary.

\begin{proof}[Proof of Proposition~\ref{prop: transitive block with prescribed lam}]
Let $\cF_1$ and $\cF_2$ be two \qms{} foliations on the oriented torus $S_1$ and $S_2$, and let $(P,X)$ be the block given by Proposition~\ref{prop: non-transitive block with prescribed lam}.
It suffices to check that $(P,X)$ satisfies the hypotheses of Proposition~\ref{prop: transitive block by blow up}.
The basic pieces of $\Lambda$ are attracting and repelling basic pieces $\Lambda_k^\pm$ and saddle periodic orbits $\cO_{i,j}$.
Each basic piece $\Lambda_k^\pm$ contains infinitely many periodic orbits with negative multipliers (Item~\ref{prop: non-transitive block; it: smale graph}, Proposition~\ref{prop: non-transitive block with prescribed lam}).
According to the Smale's graph on Figure~\ref{fig: smale graph non-transitive block} each orbit $\cO_{i,j}$ is connected by an oriented edge to an attracting basic piece $\Lambda_k^+$ and by a reversed oriented edge to a repelling basic piece $\Lambda_k^+$.
The hypotheses are thus satisfied, there exists a block $(P', X')$ associated by Proposition~\ref{prop: transitive block by blow up} which satisfies Proposition~\ref{prop: transitive block with prescribed lam}.
\end{proof}

Let us now show Proposition~\ref{prop: transitive block by blow up}.
The idea is to perform surgeries on the block $(P,X)$ to obtain a transitive block while keeping the combinatorial type of boundary lamination.
We use the so-called \emph{Blow-up -- Excise -- Glue} surgery of \cite[Section~8]{beguinBuildingAnosovFlows2017}.
We will need the following lemma.

\begin{lem} \label{lem: attracting DA block}
Let $(P,X)$ be a connected oriented filled \bb{}.
Let $O$ be a periodic orbit in the interior of $P$, without free separatrices and with negative multiplier.
Denote $\Lambda$ the maximal invariant set and $\cL$ the boundary lamination of $(P,X)$.
There exists a connected oriented block $(\hat P, \hat X)$ such that, 
\begin{enumerate}
    \item \emph{(manifold)} \label{lem: DA block; it: block}
    $\hat P = P \ssm \cV$ where $\cV$ is a tubular neighborhood of the orbit $O$.
    
    \item \emph{(boundary)} \label{lem: DA block; it: boundary}
    $\partial \hat P = \partial P \cup T$, where $T= \partial \cV$ and $\hat X$ is transverse outward on~$T$.
    
    \item \emph{(semi-conjugation)} \label{lem: DA block; it: conjugation}
    There exists a continuous surjective map $\pi \colon \hat P \ssm \hat T \to P$ which maps the oriented orbits of the flow of $\hat X$ onto the oriented orbits of the flow of $X$, where $\hat T$ is the flow-saturated set of $T$.
    There exists a periodic orbit $O' \subset \intr \hat P$ of $\hat X$ such that the restriction
    $\pi \colon \hat P \ssm (\hat T \cup \cW^s_{\hat X}(O')) \to \hat P \ssm \cW^s_X(O)$ is a homeomorphism.
    
    \item \emph{(lamination)} \label{lem: DA block; it: boundary lamination}
    The boundary lamination $\hat \cL$ is filling, $\hat \cL$ and $\cL$ completes on $\pP$ as topologically equivalent foliations on $\pP$, and $\hat \cL$ has a single compact leaf on $T$. 
    
    \item \emph{(Smale's graph)} \label{lem: DA block; it: smale graph}
    The Smale's graphs of $(\hat P, \hat X)$ and $(P,X)$ are isomorphic.\footnote{Here, we speak of the Smale's graph of the maximal invariant set of the block.}
\end{enumerate}

We say that $(\hat P, \hat X)$ is a \bb{} obtained by an \emph{attracting bifurcation Derived from Anosov (DA)} on the orbit $O$.
\end{lem}

\begin{proof}[Proof idea]
This block is obtained by an attracting DA bifurcation on the periodic orbit $O$.
We refer to \cite[Subsection~2.2.2]{ghristKnotsLinks3Dimensional1997} for a detailed description of this bifurcation, and to Figure~\ref{fig: attracting DA}.
Let us also mention \cite[Section~3.3]{barthelmeAnomalousAnosovFlows2021}, where the authors give an explicit construction of a DA bifurcation on a codimension 2 geodesic flow, but which can be adapted in our case.

\begin{figure}[h]
    \centering
    \vspace*{-2em}
    \hspace*{-2.5em}
    \includegraphics[height=0.25\textheight]{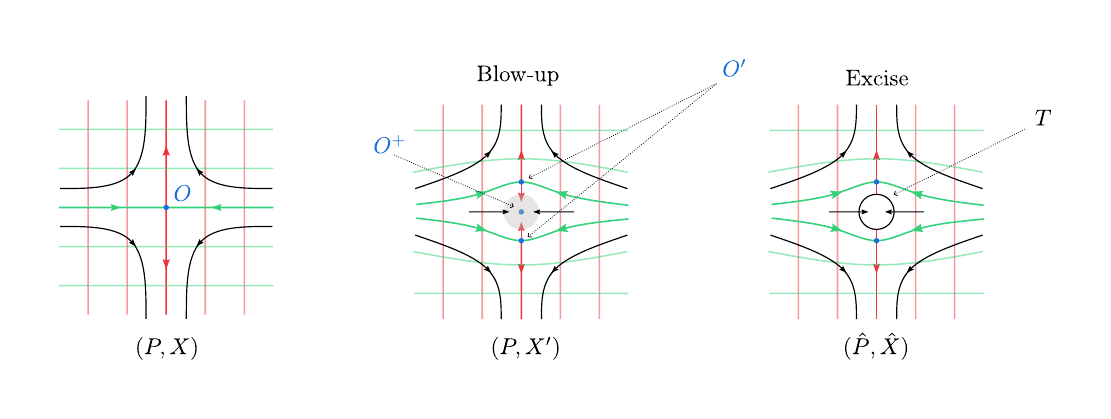}
    \vspace*{-2em}
    \caption{Attracting DA bifurcation on a periodic orbit $O$ with negative multipliers}
    \label{fig: attracting DA}
\end{figure}

This bifurcation consists in \say{opening} the stable manifold of $O$, which creates an attracting periodic orbit $O^+$ instead of $O$, and transforms $O$ into a saddle periodic orbit $O'$ with positive multipliers (corresponding to the boundary of the local unstable manifold of $O$ which is a Moebius strip).
More precisely, there exists a small perturbation $X'$ of $X$ on $P$ such that $X'$ coincides with $X$ outside a linearizing neighborhood $\cV_0$ of the orbit $O$, and such that the vector field $X'$ on $\cV_0$ admits an attracting periodic orbit $O^+$ with negative multipliers, and a periodic orbit $O$ with positive multipliers.
We remove a small tubular neighborhood $\cV = \cV_\delta$ of size $\delta$ from the orbit $O^+$, in the basin of attraction of $O^+$ and whose boundary is transverse to the \vf{} $X'$.
Denote $\hat P = P \ssm \cV$ and $\hat X$ the restriction of $X''$ on $\hat P$.
Since everything takes place far from the boundary components of $(P,X)$, we can use \cite[Propositions~8.1~and~8.3]{beguinBuildingAnosovFlows2017}, which states that if $\delta$ is small enough, the pair $(\hat P,\hat X)$ is a \bb{} which satisfies Items~\ref{lem: DA block; it: block},~\ref{lem: DA block; it: boundary}, and~\ref{lem: DA block; it: boundary lamination}.

Moreover (see \cite[Subsection~2.2.2]{ghristKnotsLinks3Dimensional1997}) there exists a continuous surjective map $\pi \colon P \to P$, which maps the oriented orbits of the flow of $X'$ onto the oriented orbits of the flow of $X$, in other words it is a semi-conjugation between the flow of $X'$ and a reparametrization of the flow of $X$.
Let $B$ be the basin of attraction of the orbit $O^+$ for the flow of $X'$, and $\hat B = B \cup \cW^s(O')$.
Then $\pi$ maps $\hat B$ onto $\cW^s(O)$, and the pre-image of a point $x \in \cW^s(O)$ is an arc through $B$ whose ends are on the manifold $\cW^s(O)$.
The restriction $\pi \colon P \ssm \hat B \to P \ssm \cW^s(O)$ is a homeomorphism, hence an orbit equivalence.
Since the set of orbits of $X'$ which intersects the boundary of $\cV$ coincide with the basin of attraction of $O^+$, Item~\ref{lem: DA block; it: conjugation} follows.

Let us show the last item, by summarize the proof of \cite[Propositions~8.1]{beguinBuildingAnosovFlows2017}.
Let $\Lambda_i$ be the basic piece of the maximal invariant set $\Lambda$ of $(P,X)$ containing the orbit $O$.
The piece $\Lambda_i$ is nontrivial because the orbit $O$ is assumed to have no free separatrices, and there exists $\gamma \in \Lambda_i$ a periodic orbit dense in $\Lambda_i$ and disjoint from $\cW^s(O)$.
Since $\pi \colon P \ssm \hat B \to P \ssm \cW^s(O)$ is a conjugation, the pre-image of $\Lambda_i \ssm \cW^s(O)$ contains a dense orbit, hence is a transitive set for the flow of $X'$.
Let $\hat \Lambda_i$ be the basic piece of $\hat \Lambda$ containing this set.
Since the unstable separatrices of $O$ are not free, the stable manifold $\cW^s(O)$ is accumulated by leaves of $\cW^s(\Lambda_i)$, and it follows by orbit equivalence of $X$ on $P \ssm \cW^s(O)$ and $X'$ on $P \ssm \hat B = P \ssm (B \cup \cW^s(O'))$ that $\cW^s(O')$ is accumulated by $\cW^s(\hat \Lambda_i)$.
It follows that $O' \subset \hat \Lambda_i$, and the map $\pi$ induces an isomorphism of Smale's graphs.
\end{proof}

\begin{rmk} \label{rmk: repelling DA block}
There exists an analogue of Lemma~\ref{lem: attracting DA block} where $(\hat P, \hat X)$ is such that $\partial \hat P = \partial P \cup T$, and $\hat X$ is transverse inward on $T$.
We say that $(\hat P, \hat X)$ is a  \bb{} obtained by \emph{repelling DA bifurcation} on the orbit $O$.
In this case, in Item~\ref{lem: DA block; it: conjugation}, we have an orbit equivalence between $X$ on $P \ssm \cW_X^u(O)$ and $\hat X$ on $\hat P \ssm \cW_{\hat X}^u(O')$.
\end{rmk}

\begin{figure}[htb]
    \centering
    \vspace*{-2em}
    \includegraphics[height=0.27\textheight]{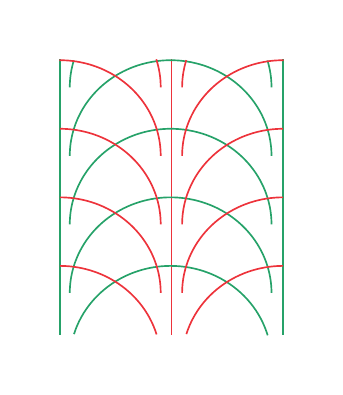}
   \vspace*{-1em}
    \caption{Pair of the strongly transverse Morse--Smale laminations $(\hat \cL, (\varphi_i)_* \hat \cL)$ on a torus $T^\iin$}
    \label{fig: transverse reeb zip}
\end{figure}

\begin{proof}[Proof of Proposition~\ref{prop: transitive block by blow up}]
Let $(P,X)$ be a block satisfying the hypothesis of Proposition~\ref{prop: transitive block by blow up}.
Denote $\Lambda_1, \dots, \Lambda_n$ the basic pieces of $\Lambda$ containing infinitely many periodic orbits with negative multipliers.
For any $\Lambda_i$, choose $2(n-1)$ periodic orbits $O_{i,j}^+$ and $O_{i,j}^-$ in $\Lambda_i$ for $j=1, \dots, n$ and $j \neq i$, without free separatrices, with negative multipliers.
Such orbits exist according to Proposition~\ref{prop: non-transitive block with prescribed lam}, Item~\ref{prop: non-transitive block; it: smale graph}, and because the periodic orbits with a free separatrix are in finite number (Lemma~\ref{lem: free separatrix}).
Let $(\hat P, \hat X)$ be the oriented connected block obtained by a \linebreak[4]repelling DA bifurcation on the orbits $O^+_{i,j}$ and an attracting DA bifurcation on the orbits $O^-_{i,j}$ by iteration of Lemma~\ref{lem: attracting DA block} and Remark~\ref{rmk: repelling DA block} for every $i,j$.
It satisfies the following properties, deduced from the previous lemma.
The boundary $\partial \hat P$ is the union $\partial P \bigcup_{i,j} (T^\iin_{i,j} \cup T^\out_{i,j})$, with $T^\iin_{i,j}$ the torus created by repelling bifurcation on $O_{i,j}^+$ where $\hat X$ is transverse inward, and $T^\out_{i,j}$ the torus created by attracting bifurcation on $O_ {i,j}^-$ where $\hat X$ is transverse outward (Item~\ref{lem: DA block; it: boundary}).
The lamination $\hat \cL$ is filling, completes on $\pP$ into a foliation topologically equivalent to $\cL$, and contains a unique compact leaf on each torus $T_{i,j}^\iin$ and $T_{i,j}^\out$ (Item~\ref{lem: DA block; it: boundary lamination}).
By composing each semi-conjugation obtained by iterating the construction (Item~\ref{lem: DA block; it: conjugation}), we obtain a semi-conjugation $h$ between the vector field $\hat X$ on the complementary of the flow-saturated set $\hat T$ of the tori $T^\out_{i,j}$ and $T^\iin_{i,j}$, and a renormalization of the vector field $X$, which is a homeomorphism on the complementary of the union of the stable manifolds $\cW^s(O^+_{i,j})$ and the unstable manifolds $\cW^u(O^-_{i,j})$.
The Smale's graph of $(\hat P, \hat X)$ is isomorphic to the Smale's graph of $(P,X)$ (Item~\ref{lem: DA block; it: smale graph}).
Let $\hat \Lambda_1, \dots, \hat \Lambda_n$ be the basic pieces of $\hat \Lambda$ associated to $\Lambda_1, \dots, \Lambda_n$.
\begin{claim} \label{claim: lamination and basic piece}
    for all $(i,j)$ we have
    $\hat \cL \cap T^\out_{i,j} = \cW^u(\hat \Lambda_i) \cap T^\out_{i,j}$,
    and
    $\hat \cL_Y \cap T^\iin_{i,j} = \cW^s(\hat \Lambda_i) \cap T^\iin_{i,j}$.
\end{claim}

\begin{proof}
It follows from Item~\ref{lem: DA block; it: conjugation} of Lemma~\ref{lem: attracting DA block} that the pre-image of $\Lambda \cap W^s(O_{i,j}^-) = \Lambda_i \cap W^s(O_{i,j}^-)$ by $h$ is the set of orbits of $\hat \Lambda$ whose unstable manifold intersects the torus $T^\out_{i,j}$, and is contained in $\hat \Lambda_i$.
Similarly, the pre-image of $\Lambda \cap W^u(O_{i,j}^+) = \Lambda_i \cap W^u(O_{i,j}^+)$ by $h$ is the set of orbits of $\hat \Lambda$ whose unstable manifold intersects the torus $T^\iin_{i,j}$, and is contained in $\hat \Lambda_i$.
\end{proof}

Let us go to the gluing operation. 
We want a \sqt{} partial gluing map $\varphi$ of $(\hat P,\hat X)$ so that the graph $G=G(\hat P,\hat X,\varphi)$ contains a cycle between $\Lambda_i$ and $\Lambda_j$ for any pair $(i,j)$.
Let $\varphi_{i,j} \colon T^\out_{i,j} \to T^\iin_{j,i}$ be a strongly (quasi)-transverse \diff{} which maps the single compact leaf of the lamination $\res{\hat \cL}{T^\out_i}$ to a leaf parallel (disjoint) to the single compact leaf of the lamination $\res{\hat \cL}{T^\iin_j}$ (Figure~\ref{fig: transverse reeb zip}).
By symmetry of the lamination, we can choose a diffeomorphism which reverses the orientation.
Let $\varphi$ be the product of \diffs{} $\varphi_{i,j}$ (made involutive).
Then $\varphi$ is a \sqt{} partial gluing map of the filled \bb{} $(\hat P, \hat X)$ which reverses the orientation.
Moreover, there is now a cycle between $\hat \Lambda_i$ and $\hat \Lambda_j$ in the graph $G = G(\hat P,\hat X,\varphi)$ for every pair $(i,j)$.
Indeed by the gluing operation and Claim~\ref{claim: lamination and basic piece}, the intersection $\varphi(\cW^u(\Lambda_i)) \cap \cW^s(\Lambda_j) \supset \varphi_{i,j} (\res{\hat \cL}{T^\out_i}) \cap \res{\hat \cL}{T^\iin_j}$ is non-empty.
Similarly, the intersection $\varphi(\cW^u(\Lambda_j)) \cap \cW^s(\Lambda_i) \supset \varphi_{j,i} (\res{\hat \cL}{T^\out_j}) \cap \res{\hat \cL}{T^\iin_i}$ is non-empty.
By assumption of Proposition~\ref{prop: transitive block by blow up}), and isomorphism of Smale's graphs, every basic piece of $\hat \Lambda$ is connected by an oriented path and a reverses oriented path to a basic piece with infinitely many negative multipliers $\hat \Lambda_i$, and every such pair $\hat \Lambda_i, \hat \Lambda_j$ is connected by a cycle in $G(\hat P, \hat X,\varphi)$.
It follows that any two basic pieces are connected by a path of oriented edges in $G(\hat P, \hat X,\varphi)$, in other words the graph is strongly connected.

According to Theorem~\ref{thm: partial gluing theorem} (generalization of Theorem~\ref{thm: gluing theorem} for partial gluing maps), up to make a strong isotopy of the triple $(\hat P, \hat X, \varphi)$,
the quotient $P' := \hat P /\varphi$ is a compact oriented manifold with boundary equipped with a vector field $Z$ induced by $\hat X$, such that the pair $(P', X')$ is a connected oriented \bb{}.
The boundary of $P'$ is $\partial P' = \partial P$, and the boundary lamination $\cL'$ contains the filling lamination $\hat \cL$ (Item~\ref{lem: block by partial gluing; it: lamination}, Lemma~\ref{lem: block obtained by partial gluing without cycle}), so it is a filling lamination of the same combinatorial type.
It follows that these two laminations complete each other in topologically equivalent foliations on $\partial P'$ (Lemma~\ref{lem: extend qms prefoliation by qms foliation} and Proposition~\ref{prop: qms foliation for given combi type}), so the same occurs for $\cL'$ and $\cL$.
By strong connectedness of the graph $G(\hat P, \hat X,\varphi)$, the block $(P',X')$ is transitive (Proposition~\ref{prop: transitivity criterion}).
We conclude that $(P',X')$ satisfies Proposition~\ref{prop: transitive block by blow up}.
\end{proof}

\subsection{Realizing quasi-transverse bifoliation in an Anosov flow}
\label{sec: prescribed boundary lam; subsec: bifoliation}

\begin{figure}[htb]
    \centering
    \vspace*{-1em}
    \includegraphics[scale=0.55]{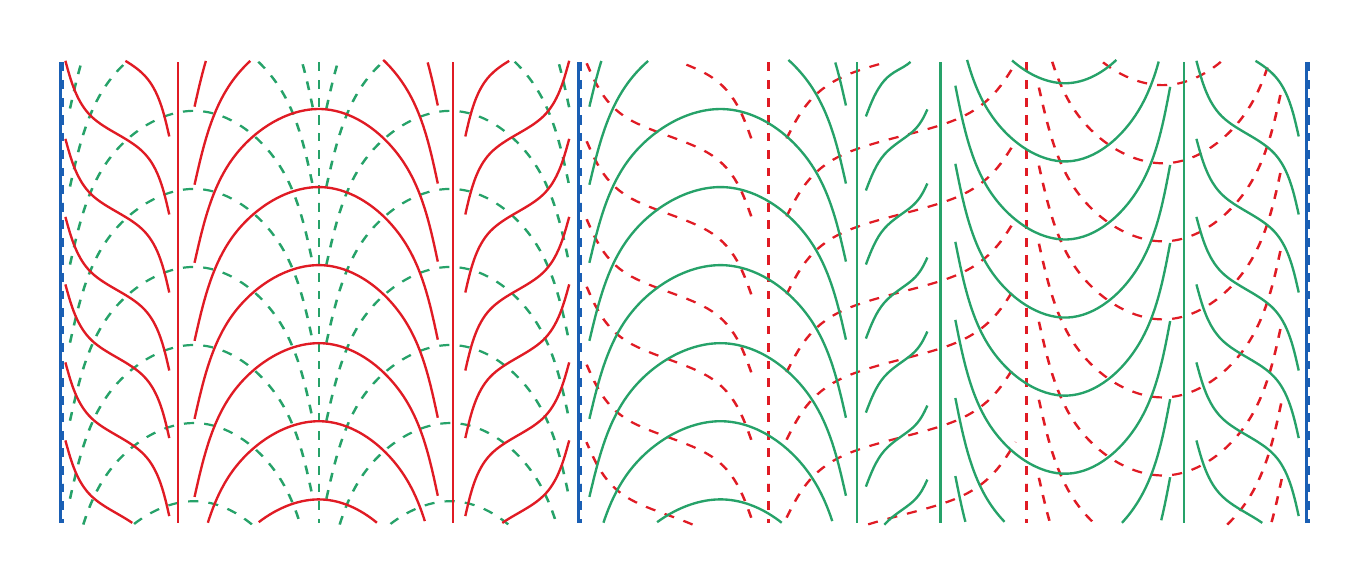}
    \vspace*{-2.5em}
    \caption{A \qt{} bifoliation on the torus}
    \label{fig: qt bifoliation}
\end{figure}

In this subsection we show that we can realize any pair of \qt{} foliations as the trace of the stable and unstable foliations of a transitive Anosov flow on an embedded \qt{} torus (Theorem~\ref{thmintro: realize qms bifoliation})
Recall that a pair of laminations $(\cF_1, \cF_2)$ on torus $T$ is said to be \emph{quasi-transverse} if it is a pair of two \qms{} foliations which coincide along their oriented marked leaves, are transverse on the complementary of the marked leaves, and such that $T^\iin_{\cF_1} = T^\out_{\cF_2}$ (Definition~\ref{def: sqt laminations} and Figure~\ref{fig: qt bifoliation}).

\begin{lem} \label{lem: qt bifoliation and compact leaves}
Let $(\cF_1, \cF_2)$ be a pair of \qt{} foliations on torus $T$, such that the set $\Gamma_*$ of marked (common) compact leaves is \emph{non-empty}.
Let $\Gamma = \Gamma(\cF_1 \cup \cF_2) = \gamma_0, \dots, \gamma_{n-1} $ be the collection of compact leaves of $\cF_1 \cup \cF_2$.
Then the compact leaves $\Gamma$ of $\cF_1 \cup \cF_2$ are non-contractible in $T$ and pairwise freely homotopic as non-oriented closed curves in $T$.
The compact leaves of $\Gamma \ssm \Gamma_*$ are pairwise disjoint.
\end{lem}

\begin{proof}
This is a consequence of the fact that the foliations $\cF_i$ have a non-empty set of marked leaves, and coincide along these leaves.
It is then sufficient to apply the corollary~\ref{coro: homotopic compact leaf in qms prefoliation}.
By transversality of the laminations on the complementary of the marked leaves, the compact unmarked leaves are pairwise disjoint.
\end{proof}

\vspace*{-0.1em}

There is thus a cyclic order on the elements of $\Gamma$.
We provide the elements of $\Gamma$ with the dynamical orientation induced by the dynamical orientation of the respective compact leaves of $\cF_1$ and $\cF_2$. 
It coincides along the common marked leaves.
Let be a geometric enumeration (Definition~\ref{def: geometric enumeration}) $\Gamma= \{\gamma_0, \dots, \gamma_{n-1} \}$ of the compact leaves of $\cF_1 \cup \cF_2$. 
Denote $\Gamma_* \subset \Gamma$ the collection of common marked leaves.
The following definition allows us to characterize the pattern of the quasi-transverse intersection of the two foliations.

\begin{defi}[Combinatorial type of a bifoliation] \label{def: combi type of bifoliation}
Define the map
$$ \sigma = \sigma(\cF_1, \cF_2) \colon \Z/n\Z \; \longrightarrow \; \{1,2\} \times \l( \scA_H \times \scA_O \times \scA_H \r) $$
by:
\begin{enumerate}[i)]
    \item $ \sigma(k) = (1, * )$ if and only if $\gamma_k$ belongs to the foliation $\cF_1$, and $ \sigma(k) = (2, * )$ if and only if $\gamma_k$ belongs to the foliation $\cF_2 \ssm \Gamma_*$;
    \item $ \sigma(k) = (*, (* \ua *) )$ if and only if the leaf $\gamma_i$ is freely homotopic to $\gamma_0$ as oriented paths;
    \item $ \sigma(k) = (i, (\ra * *))$ if and only if the holonomy from $\cF_i$ on the left\footnote{the orientation of the first leaf and of the torus $T$ determines a left and a right for each compact leaf} of $\gamma_k$ is contracting;
    \item $ \sigma(k) = (i, (* * \la))$ if and only if the holonomy of $\cF_i$ on the right of $\gamma_k$ is contracting;
\end{enumerate}
We say that $\sigma$ is a \emph{combinatorial type of the bifoliation $(\cF_1, \cF_2)$ on $T$}.
\end{defi}
This definition is not symmetrical because we make the arbitrary choice to privilege $\cF_1$, and $\sigma(\cF_1, \cF_2)$ is not equal to $\sigma(\cF_2, \cF_1)$.
Indeed, we consider the marked leaves as compact leaves of $\cF_1$.
If $\sigma_1$ is a combinatorial type \qms{} of $\cF_1$, then there exists an injection $\phi_1 \colon \Z/n_1 \Z \to \Z/n\Z$ such that
$\sigma(\cF_1, \cF_2) \circ \phi_1 = (1, \pm \sigma_1)$.
This is not true for $\cF_2$.

\begin{example}
We give the combinatorial type of the bifoliation of Figure~\ref{fig: qt bifoliation}, where the geometrical enumeration on the compact leaves is the one from left to right, the green refers to a \say{in} component, the red refers to a \say{out} component, $\cF_1$ is the full line foliation and $\cF_2$ the dashed line foliation
\begin{align*}
    \sigma(\cF_1, \cF_2) = 
& \{ (1, (\la, \ua, \la)); \ 
(1, (\la, \da, \ra); \  
(2, (\ra, \ua, \la);\ 
(1, (\la, \da, \ra);\ 
(1, (\ra, \ua, \ra);\ \\ 
&(2, (\ra, \da, \la);\ 
(1, (\la, \ua, \ra);\ 
(1, (\la, \da, \ra);\ 
(2, (\ra, \ua, \la);\ 
(1, (\la, \da, \ra) \}.
\end{align*}

\end{example}

\begin{mainthm}[Theorem~\ref{thmintro: realize qms bifoliation}] \label{thm: realize qms bifoliation}
Let $\sigma$ be a combinatorial type of a \qt{} bifoliation.
There exists a transitive Anosov vector field $Z$ on an oriented $3$-manifold $\cM$ and an incompressible torus $T$ embedded in $\cM$, \qt{} to $Z$, such that the trace of the stable and unstable foliation $\cF^s$ and $\cF^u$ on $T$ induces a bifoliation $(\cF_1, \cF_2)$ on $T$ of combinatorial type $\sigma$.
\end{mainthm}

\begin{proof}
Let $(\cG_1, \cG_2)$ be a bifoliation of combinatorial type $\sigma$ on the oriented torus $T$.
Denote by $\check T$ the torus $T$ with the opposite orientation.
Let $(P,X)$ be the transitive filled block given by Proposition~\ref{prop: transitive block with prescribed lam} associated to the foliations $\cG_1$ on the oriented torus $T$ and $\cG_2$ on the oriented torus $\check T$.
Denote $\cL$ the boundary lamination of $(P,X)$.
The boundary $\pP$ is the union of two tori $T_1$ and $T_2$ such that the lamination $\cL$ on $T_1$ is filling and extends into a topologically equivalent foliation $\cG_1$ on $T$ and the lamination $\cL$ on $T_2$ is filling and extends into a topologically equivalent foliation $\cG_2$ on $\check T$.
By construction, there exists an involution $\varphi \colon \partial P \to \partial P$ which pairs $T_1$ and $T_2$, which reverses the orientation and such that $(\varphi_* \cL_X, \cL_X)$ is a pair of \sqt{} filling laminations, and completes on $T_1$ into a pair of \sqt{} foliations topologically equivalent to the pair $(\cG_1, \cG_2)$.
This diffeomorphism induces a \sqt{} gluing map of the block $(P,X)$.
By the Gluing Theorem~\ref{thm: gluing theorem},
there exists a strongly isotopic triple $(P',X', \varphi')$ such that the quotient $\cM := P'/\varphi'$ is an oriented closed 3-manifold with a vector field $Z$ induced by $X'$ which is Anosov.
As $(P,X)$ is transitive, so is $(P',X')$ and the Anosov flow $Z$ is transitive (Proposition~\ref{prop: transitivity criterion}).
The projection of $\partial P'$ in $\cM$ is an oriented torus $T$ \qt{} to the Anosov \vf{} $Z$ so it is incompressible (\cite{barbotMisePositionOptimale1995, brunellaSeparatingBasicSets1993}).
The stable and unstable foliations $(\cF^s, \cF^u)$ of the Anosov flow induce on $T$ a pair of foliations $(\cF_1, \cF_2)$, which contains the projection of the pair of filling laminations $(\varphi'_* \cL_{X'}, \cL_{X'})$.
The pair $(\cF_1, \cF_2)$ is not a priori topologically equivalent to the pair $(\cG_1, \cG_2)$ (this is lost by strong isotopy, see Remark~\ref{rmk: no usual sqt isotopy for equivalent triples}) but it follows from Proposition~\ref{prop: properties of strongly isotopic triples}, Item~\ref{prop: strongly isotopic trip, it: lamination pattern} that the combinatorial type of the bifoliation which extends the pair $(\varphi'_* \cL_{X'}, \cL_{X'})$ and of the bifoliation which extends the pair $(\varphi_* \cL_{X}, \cL_{X})$ are equal.
It follows that the combinatorial type of the bifoliation of $(\cF_1, \cF_2)$ is equal to $\sigma$.
\end{proof}

\section{Embedding a filled block in an Anosov flow}
\label{sec: embed block in anosov flow}
The goal of this section is to show Theorem~\ref{thmintro: embed block in anosov}.
We recall that a block is said to be \emph{filled} if the boundary lamination is filling.
We refer to Section~\ref{sec: preliminaries} for the definitions.
In this section, we consider orientable blocks.

\begin{mainthm}[Theorem~\ref{thmintro: embed block in anosov}] \label{thm: embed block in anosov}
For any (transitive) orientable, filled block $(P,X)$, there exists a (transitive) Anosov vector field $Z$ on a closed orientable 3-manifold $\cM$, such that $(P,X)$ is embedded in $(\cM, Z)$.
More precisely, there exists a finite collection of incompressible tori $\cT$ embedded in $\cM$, \qt{} to $Z$, such that the closure of one \cc{} of $\cM \ssm \cT$ is a compact submanifold diffeomorphic to $P$ and such that the restriction of $Z$ on $P$ is orbit equivalent to~$X$.
\end{mainthm}

Let us summarize the proof.
First, we notice that for any \qms{} foliation $\cF$ on $\T^2$, there exists a \qms{} foliation $\cG$ \qt{} to $\cF$ on $\T^2$ (Lemma~\ref{lem: foliation qt to a given qms foliation}).
For a given block $(P,X)$, we associate to each (quasi-transverse) boundary component $T_k$ of $\pP$ a block $(N_k, Y_k)$ given by Proposition~\ref{prop: transitive block with prescribed lam}, whose filling boundary lamination $\cL_{Y_k}$ extends to a foliation topologically equivalent to a foliation \qt{} to the boundary lamination $\cL_X$ of $(P,X)$ on $T_k$.
The block $(P,X)$ and $(N_k, Y_k)$ are glued together along the boundary $T_k$ by a \sqt{} (partial) gluing map which reverses the orientation, and we proceed in this way for each boundary component $T_k$.
Up to make a strong isotopy of the triples, the Gluing Theorem~\ref{thm: gluing theorem} tells us that the vector fields of the different blocks induce an Anosov vector field on the quotient (orientable) manifold, and transitive according to the criterion of Proposition~\ref{prop: transitivity criterion}.

\subsection{Construction of a quasi-transverse foliation}

\label{dec: embed block; subsec: qt foliation}
A first step is to show that for any \qms{} foliation $\cF$ on $\T^2$, there exists a \qms{} foliation $\cG$ \qt{} to $\cF$ on $\T^2$.

\begin{lem}
\label{lem: foliation qt to a given qms foliation}
Let $\cF$ be a \qms{} foliation on the torus $T$.
Then there exists a \qms{} foliation $\cG$ on $T$ such that $\cG$ is \qt{} to~$\cF$.
\end{lem}

\begin{proof}
Let $\Gamma_\cF$ denote the compact leaves of $\cF$, $\Gamma_*$ the marked compact leaves of $\cF$, and $T^\iin \cup T^\out$ the (in,out)-splitting of $T$ for $\cF$.
Let $A$ be the closure of a $T \setminus \Gamma_\cF$.
Suppose that this \cc{} belongs to $T^\iin$.
Let $\gamma \in \Gamma_\cF$ be a connected component of $\partial A$.
It is a compact leaf of $\Gamma_\cF$ oriented by the dynamical orientation (Definition~\ref{def: dynamical orientation compact leaf}, Remark~\ref{rmk: canonical dynamical orientation block}), such that the holonomy of the oriented leaf $\gamma$ of the foliation $\cF$ is expanding.
Let $N = N(\gamma)$ be a small collar neighborhood of $\gamma$ in $\partial A$ and $(x, \theta) : N \to [0,1] \times \R/\Z$ coordinates on $A$ such that $\gamma = \{ x = 0 \}$ is oriented in the positive direction of $\theta$.
Two foliations whose holonomy is contracting are topologically equivalent, so we can assume that the direction of $\cF$ is given by $dx + x d\theta = 0$.
We distinguish two cases:
\begin{enumerate}
    \item If $\gamma \in \Gamma_*$ is a marked leaf of $\cF$, we define $\cG$ on $N$ by the foliation of direction \mb{$dx - x d\theta = 0$.}
    \item If $\gamma \in \Gamma_\cF \ssm \Gamma_*$ is an unmarked compact leaf of $\cF$, we define $\cG$ on $N$ by the direction \mb{$d\theta - x dx = 0$}.
\end{enumerate}
We do a similar construction on a collar neighborhood $N(\gamma')$ of the second boundary component $\gamma'$ of $\partial A$.
We obtain a foliation $\cG$ on a collar neighborhood $V = N(\gamma) \cup N(\gamma')$ of $\partial A$ in $A$.
Choose the metric $\Vert \cdot \Vert^2 = dx^2 + d\theta^2$ on $N(\gamma)$ and $N(\gamma')$.
Then by construction, $\cG$ is orthogonal for this metric to $\cF$ on the boundary of $V$, and we extend $\cG$ on $A$ by the foliation defined by the vector field $X'$ orthogonal to a vector field $X$ generating the foliation $\cF$ in the metric $\Vert \cdot \Vert$.
By density of Morse--Smale vector fields, there exists a Morse--Smale \vf{} $X''$ close to $X'$ in $\cC^1$ topology on $A$, which coincides with $X'$ on a neighborhood of $\partial A$.
This vector field defines a foliation $\cG'$ on $A$, which contains a finite number of compact leaves, all of which have an orientation such that the holonomy is expanding (this is the dynamical orientation), and such that $\cG'$ is transverse to $\cF$ on the interior of $A$ (transversality is an open property).
We thus construct a foliation $\cG'_k$ on each of the \ccs{} $A_k$ of $T \ssm \Gamma_{\cF}$, all of which reconnect into a foliation $\cG'$ on $T$, which is a \qms{} foliation \qt{} to $\cF$ on $T$.
\end{proof}

\subsection{Proof of Theorem~\ref{thmintro: embed block in anosov}}
\label{dec: embed block; subsec: proof}

\begin{figure}[h]
    \centering
    \vspace*{-2em}
    \includegraphics[width=\textwidth]{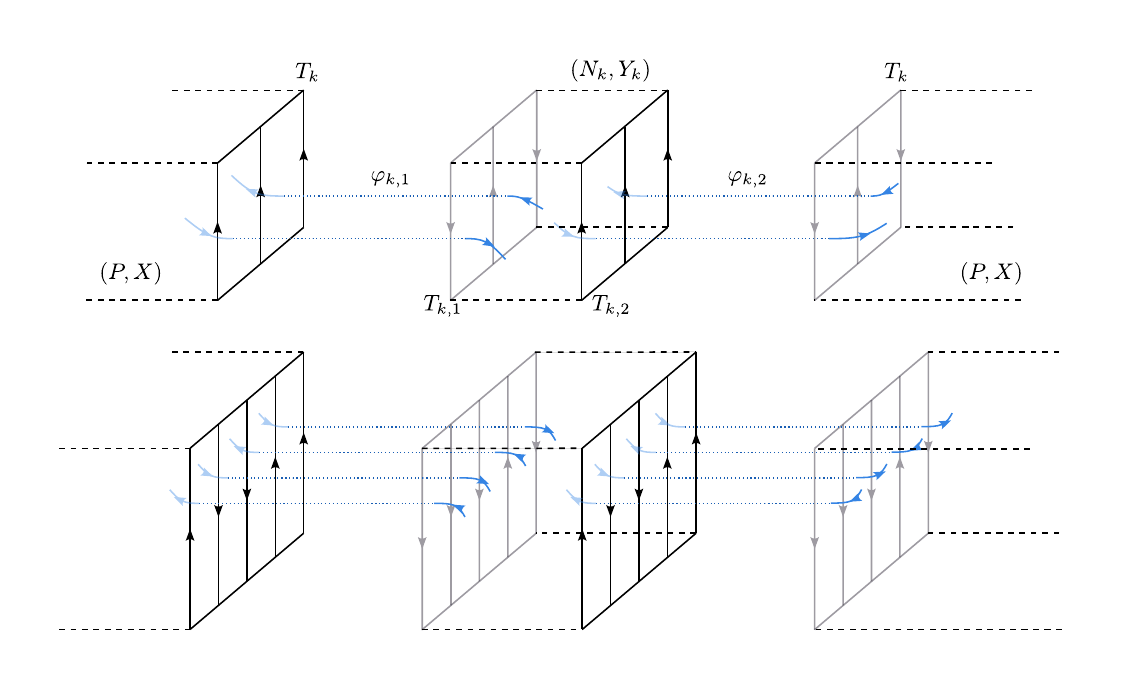}
    \vspace*{-1em}
    \centering \caption{Gluing two copies of $(P,X)$ and $(N_k, Y_k)$}
    \label{fig: embed block}
\end{figure}

\begin{proof}[Proof of Theorem~\ref{thm: embed block in anosov}]

Let $(P,X)$ be a transitive filled orientable \bb{}.
Let $\cL$ be the boundary lamination of $(P,X)$.
We denote $T_1, \dots, T_n$ the \ccs{} of $(P,X)$, and assume to begin with that none of them is transverse to the \vf{} $X$.
For each $k=1, \dots, n$, let $\cF_k$ be a \qms{} foliation on $T_k$ containing $\cL$ as sublamination (Lemma~\ref{lem: extend qms prefoliation by qms foliation}).
Let $\cG_k$ be a \qms{} foliation on $T_k$ \qt{} to $\cF_k$ given by Lemma~\ref{lem: foliation qt to a given qms foliation}.
Let $(N_k, Y_k)$ be the \bb{} given by Proposition~\ref{prop: transitive block with prescribed lam} associated to two copies of the lamination $\cG_k$ on the oriented torus $\check T_k$ (which denotes the torus $T_k$ provided with the reversed orientation).
The boundary $\partial N_k$ is formed by two connected components $T_{k,1}$ and $T_{k,2}$ and the boundary lamination of $(N_k, Y_k)$ on $T_{k,i}$ extends to a foliation topologically equivalent to the foliation $\cG_k$ on $\check T_k$.
We deduce that there exists a \diff{} $\varphi_{k,i} : T_{k,i} \to T_k$ which maps the boundary lamination of $(N_k, Y_k)$ to a lamination \sqt{} to the boundary lamination of $(P,X)$, and which reverses the orientation.
The product $\varphi_k$ of the $\varphi_{k,i}$ (made involutive) is a \sqt{} {partial} gluing map of the block formed from the union of two copies of $(P,X)$ and the block $(N_k, Y_k)$, and which reverses the orientation.
The product $\varphi$ of $\varphi_k$ (made involutive) is a (non-partial) \sqt{} gluing map of the filled block $(M,Y)$ formed from two copies of $(P,X)$ and the union $\bigcup_k (N_k, Y_k)$ (Figure~\ref{fig: embed block}), and which reverses the orientation.
According to the Theorem~\ref{thm: gluing theorem}, up to modify $(M,Y, \varphi)$ by strong isotopy, the quotient $\cM := M/\varphi$ is a closed, orientable manifold, with a vector field $Z := Y_\varphi$ induced by $Y$ which is Anosov.
The transitivity is deduced from the strong connectedness of the associated graph $G = G(M,Y,\varphi)$ (Proposition~\ref{prop: transitivity criterion}).
Indeed, the maximal invariant set of the block $(M,Y)$ is the union of the maximal invariant sets of the two copies of $(P,X)$ and the maximal invariant sets of the blocks $(N_k, Y_k)$.
Each of these sets is transitive, thus a basic piece.
The gluing map $\varphi_k$ pairs a boundary component $T_k$ of each of the two copies of $(P,X)$ with one of the two boundary components of $(N_k, Y_k)$.
By assumption, each boundary component of $(P,X)$ and $(N_k, Y_k)$ contains a non-zero number of periodic orbits, so the stable and unstable manifolds of the maximal invariant set of each of the blocks intersects each of the boundary components of the block.
It follows that the image of the entrance lamination $\cL^\iin_k$ of $(N_k, Y_k)$ intersects the exit boundary lamination $\cL^\out_X$ of each of the two copies of $(P,X)$ via the gluing map $\varphi$, and similarly the image of the exit lamination $\cL^\out_k$ of $(N_k, Y_k)$ intersects the entrance boundary lamination $\cL^\iin_X$ of each of the two copies of $(P,X)$ via the gluing map $\varphi$ (all boundary laminations are filling).
We deduce that the maximal invariant set of the two copies of $(P,X)$ and $(N_k, Y_k)$ are connected by a cycle in the graph $G(M,Y,\varphi)$.
Since this result is true for all $k$, the graph $G$ is strongly connected.
The block $P$ is naturally embedded in $\cM$ in a submanifold whose boundary is the union of tori \qts{} to the Anosov vector field $Z$.
Finally, we know from \cite{barbotMisePositionOptimale1995} and \cite{brunellaSeparatingBasicSets1993} that any torus embedded in a \qt{} position in an Anosov flow is incompressible.

The proof generalizes to any $(P,X)$ by using the results of \cite{beguinBuildingAnosovFlows2017} to pairs the connected components that are transverse to the vector field.
\end{proof}

\section{Realize an abstract geometric type in a filled block}
\label{sec: geometric type}

In this section we are interested in the dynamics in the neighborhood of the maximal invariant set of a filled building block.
All the manifolds we consider will be orientable.
We show Theorem~\ref{thmintro: geometric type in block} which gives a necessary and sufficient criterion for an orientable filled block to realize an \emph{abstract geometric type} as the geometric type of a Markov partition of its maximal invariant set.
This will be done using a correspondence between \bbs{} and \bby{} blocks given by operations of handles attaching and removal (Subsection~\ref{sec: geometric type; subsec: bb and bby}), and results from \cite{beguinFlotsSmaleDimension2002} and \cite{beguinConstructionFlotsSmale1999} which give the existence and uniqueness of a canonical \bby{} block (the \emph{model}) which realizes an abstract geometric type.
We then describe a simple algorithmic method which allows us to check these criteria from the geometric type (Subsection~\ref{sec: geometric type; subsec: example model}).
Let us start with some precise definitions.

\begin{defi}[Model block] \label{def: model block}
An orientable building block $(U,X)$ is said to be a \emph{model block} if, denoting $\Lambda$ the maximal invariant set and $\cL$ the boundary lamination of $(U, X)$, we have
\begin{enumerate}
    \item \label{def: model block; it: transverse boundary}
    $X \pitchfork \partial U$,
    \item \label{def: model block; it: connected components}
    any connected component of $U$ intersects $\Lambda$,
    \item \label{def: model block; it: saddle}
    $\Lambda$ does not contains neither attractors or repellers (we say that $\Lambda$ is \emph{saddle}),
    \item \label{def: model block; it: lam}
    any simple closed curve $c$ embedded in $\partial U$ and disjoint from $\cL$ is the boundary of a disk $D \subset \partial U$ disjoint from $\cL$.
\end{enumerate}
\end{defi}
In particular, a model block is a \bby{} block.

\begin{rmk} \label{rmk: model of the germ and model block}
This definition is compatible with the definition \cite[Definition~0.2]{beguinFlotsSmaleDimension2002} of the model of the germ of a saddle saturated set $\Lambda$ of a Smale \vf{} $X$ on a compact orientable manifold $\cM$ of dimension 3.
A compact invariant set $\Lambda$ for a vector field $X$ is said to be \emph{saddle saturated} if it is a hyperbolic set of index $(1,1)$, without attractor or repeller, and contains the intersection of its stable and unstable manifolds, in other words $\Lambda = \cW^s(\Lambda) \cap \cW^u(\Lambda)$.
The \emph{germ} $[X, \Lambda]$ is the equivalence class of the couple $(X, \Lambda)$ for the relation $(X, \Lambda) \sim (X', \Lambda')$ if there exists a neighborhood of $V$ of $\Lambda$ and $V'$ of $\Lambda'$ such that the vector fields $X$ on $V$ and $X'$ on $\Lambda'$ are orbit equivalent.
A \emph{model of germ of a saddle saturated set} $[X, \Lambda]$ is then a pair $(U,Y)$ where $U$ is a compact orientable manifold with boundary of dimension 3, provided with a vector field $Y$ transverse to the boundary, such that if $\Lambda_Y$ is the maximal invariant set of the flow of $Y$ in $U$, then $[Y,\Lambda_Y] = [X, \Lambda]$, and Items~\ref{def: model block; it: connected components} and~\ref{def: model block; it: lam} of Definition~\ref{def: model block} are satisfied.
We can indeed speak of boundary lamination of $(U,Y)$ because $(U,Y)$ is then a \bb{} in the sense of Definition~\ref{def: building block}.
A model block $(U,X)$ according to Definition~\ref{def: model block} is then a model of the germ $[X, \Lambda]$ of its maximal invariant set.
In \cite[Theorem~0.3]{beguinFlotsSmaleDimension2002} the authors show that the model of the germ of a saddle saturated set is unique up to orbit equivalence, which justifies the term of \emph{model}.
\end{rmk}

We follow the definition of \cite[Definition~2.3 and Definition~2.4]{beguinFlotsSmaleDimension2002} for the geometric type.

\begin{defi}[Abstract geometric type]
An abstract geometric type $\scT$ is the data of
\begin{itemize}[--]
    \item an integer $n \in \N^*$,
    \item for any $1 \leq i \leq n$, the data of two integers $h_i$ and $v_i$ such that $\sum_i h_i = \sum_i v_i$,
    \item a map $\phi$ from 
    $\{ i \in \{1, \dots, n\}, j \in \{1, \dots, h_i \} \}$ to $\{ k \in \{1, \dots, n\}, l \in \{1, \dots, v_k \} \} \times \{\pm 1\}$
    which induces a bijection by forgetting the signs.
\end{itemize}
\end{defi}

Let $\cR =(\Sigma, R)$ be a Markov partition in the sense of Definition~\ref{def: markov partition}.
We consider partitions embedded in orientable manifolds, with a choice of orientation on $\Sigma$.
Then a choice of an orientation of the verticals of $\cR$ induces an orientation of the horizontals of $\cR$.

\begin{defi}[Geometric type of a Markov partition]
Let $\cR =(\Sigma, R)$ be a Markov partition with an orientation of verticals and horizontals, and let $f : R \to \Sigma$ be the return map of the Markov partition.
The geometric type $\scT$ of $\cR$ is an abstract geometric type defined by
\begin{itemize}[--]
    \item $n$ is the number of rectangle of $R$;
    \item $h_i$ is the number of connected components of the intersection of $R_i \cap f^\inv (R)$, these are horizontal subrectangles $H_i^1,\dots H_i^{h_i}$ of $R_i$, numbered in the order induced by the orientation of the verticals of $R_i$;
    \item $v_k$ is the number of connected components of the intersection $R_k \cap f(R)$, they are vertical subrectangles $V_k^1,\dots V_k^{v_k}$ of $R_k$, numbered in the order induced by the orientation of the horizontals of $R_k$;
    \item $\phi(i,j) = ((k,l), \epsilon)$ if $f(H_i^j) = V_k^l$, and $\epsilon = +$ if the orientation of the verticals of $f(H_i^j)$ coincides with the orientation of the verticals of $V_k^l$, and $\epsilon = -$, otherwise.
\end{itemize}
\end{defi}

\begin{figure}[htb]
    \centering
    \vspace*{-1em}
    \includegraphics[height=0.3\textheight]{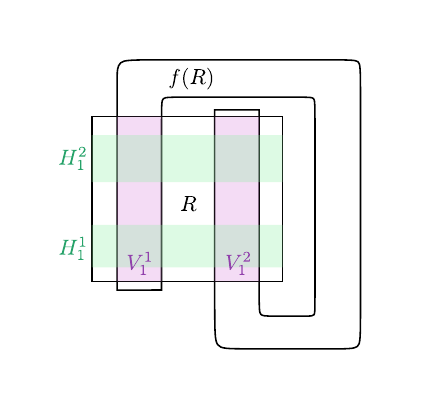}
    \vspace*{-1em}
    \centering \caption{Partition of the fake horseshoe}
    \label{fig: partition fake horseshoe}
\end{figure}

\begin{example} \label{ex: type of horseshoe}
The geometric type of the fake horseshoe (Figure~\ref{fig: partition fake horseshoe}) is
$$ \scT = \{ n=1,\ h_1=2,\ v_1=2,\ \ \phi(1,1) = (1,1,+), \, \phi(1,2) =(1,2,+) \}.$$
\end{example}

\begin{thm}[{\cite[Theorem~0.1]{beguinConstructionFlotsSmale1999} and \cite[Theorem~0.4]{beguinFlotsSmaleDimension2002}}]
\label{thm: model geometric type}
For any abstract geometric type $\scT$, there exists a model block $(U, X)$, whose maximal invariant set admits a Markov partition of geometric type $\scT$, and this one is unique up to orbit equivalence.
We say that $(U, X)$ is \emph{the model of the geometric type $\scT$}.
\end{thm}

The aim of this section is to show the following main proposition.

\begin{mainthm}[Theorem~\ref{thmintro: realize qms bifoliation}] \label{thm: geometric type in block}
Let $\scT$ be an abstract geometric type.
There exists an orientable filled block $(P,X)$, admitting a Markov partition of geometric type $\scT$ if and only if the model $(U, Y)$ of $\scT$ satisfies the following conditions:
\begin{enumerate}
    \item \label{prop: geometric type; it: boundary components}
    $\partial U$ is a union of tori and spheres, each sphere contains exactly two open disks $D_i$ and $D_j$ disjoint from the boundary lamination $\cL_Y$, bounded by two distinct compact leaves of $\cL_Y$;
    \item \label{prop: geometric type; it: prefoliation}
    $\cL_Y$ is a pre-foliation on the complementary $\partial U \ssm \bigcup_i D_i$.
\end{enumerate}
Such a block $(P,X)$ is then unique up to isotopy.
Moreover, the lamination $\cL_X$ is topologically equivalent on $\overline \Pin$ to the lamination $\cL_Y$ on $U^\iin \ssm \bigcup_i D_i$, and the lamination $\cL_X$ is topologically equivalent on $\overline \Pout$ to the lamination $\cL_Y$ on $U^\out \ssm \bigcup_i D_i$.
\end{mainthm}

The proof consists in going from a building block to a \bby{} block (and vice versa) by handles attaching and removal.
We analyze the boundary lamination of these corresponding blocks obtained by such surgeries.

\begin{itemize}
    \item In Subsection~\ref{sec: geometric type; subsec: bb and imb}, we build a correspondence between \bbs{} and \say{intermediate} blocks (Definition~\ref{def: intermediate block}), which will be \say{standard} \minc{}s of \bbs{}, whose boundary decomposes into an entrance boundary, an exit boundary and a tangent boundary.
    These will be the pieces of our surgeries.
    We will be able to glue handles along the tangent boundary to build blocks whose boundary is transverse to the vector field.
    
    \item In Subsection~\ref{sec: geometric type; subsec: bb and bby}, we describe two mutually reciprocal surgical operations which allow us to go from \bbs{} to \bby{} blocks and vice-versa.
    More precisely, we show that every \bbs{} is associated with a single \bby{} block obtained by gluing spheres together along the boundary (Lemma~\ref{lem: bb to bby}).
    The reciprocal operation consists of, given a \bby{} block and a finite collection of connected components of the complementary of the lamination which are bordered by compact leaves (satisfying some assumptions), to construct a \bb{} by a finite number of removal of handles which has the effect of \say{killing} the components of the collection to \say{replace} them by periodic orbits contained in the boundary (Lemma~\ref{lem: bby to bb}).
    
    \item Finally, we use these two reciprocal operations, to show the necessary and sufficient criterion of Theorem~\ref{thm: geometric type in block}.
    
\end{itemize}

\subsection{Intermediate blocks}

\label{sec: geometric type; subsec: bb and imb}

This subsection is a preliminary technical section in order to perform surgery between building blocks and \bby{} blocks in the next subsection.

\begin{defi}[Intermediate block] \label{def: intermediate block}
Let $\check P$ be a3-manifold with boundary, and $\check X$ a vector field of class $\cC^1$ on $\check P$.
We say that $(\check P, \check X)$ is an \emph{\imb{}} if
\begin{enumerate}
    \item \label{def: imb; it: edges}
    $\check P$ is a smooth manifold outside a finite number of simple closed curves $c_1^\iin, \dots, c_n^\iin$ and $c_1^\out, \dots, c_n^\out$ contained in $\partial \check P$.
    \item \label{def: imb; it: boundary}
    $\partial \check P$ splits into the union $\partial \check P = \check P^\iin \cup \check P^\tan \cup \check P^\out$, where
    \begin{itemize}[--]
        \item $\check P^\iin$ is a surface with boundary, bordered by $c_1^\iin, \dots, c_n^\iin$, the \vf{} $\check X$ is transverse to $\check P^\iin$ and points inward;
        \item $\check P^\out$ is a surface with boundary, bordered by $c_1^\out, \dots, c_n^\out$, the \vf{} $\check X$ is transverse to $\check P^\out$ and points outward;
        \item $\check P^\tan$ is a union of annuli $A_1, \dots, A_n$ tangent to the \vf{} $\check X$, such that each $A_i$ is bounded by the curves $c^\iin_i$ and $c^\out_i$, and the orbits of $\check X$ on $A_i$ go from $c^\iin_i$ to $c^\out_i$.
    \end{itemize}
    \item \label{def: imb; it: hyp max inv}
    The maximal invariant set $\check \Lambda := \bigcap_{t \in \R} \check X^t(P)$ of $\check X$ in $\check P$ is a hyperbolic set of index $(1,1)$ for $\check X$.    
    \item \label{def: imb; it: annuli in boundary}
    If we denote $\check \cW^s$ and $\check \cW^u$ the stable and unstable manifolds of $\check \Lambda$ in $\check P$, then for any $i$ there exist two compact leaves $\gamma^s_i$ of $\check \cL^\iin := \cW^s \cap \check P^\iin$ and $\gamma^u_i$ of $\check \cL^\out:= \cW^u \cap \check P^\out$, such that $A_i$ is contained in an annulus $\cA_i \subset \partial \check P$, bounded by $\gamma^s_i$ and $\gamma^u_i$ and disjoint from $\check \cW^s \cup \check \cW^u$.
\end{enumerate}
\end{defi}

\begin{figure}[htb]
    \centering
    \vspace*{-1em}
    \includegraphics[height=0.28\textheight]{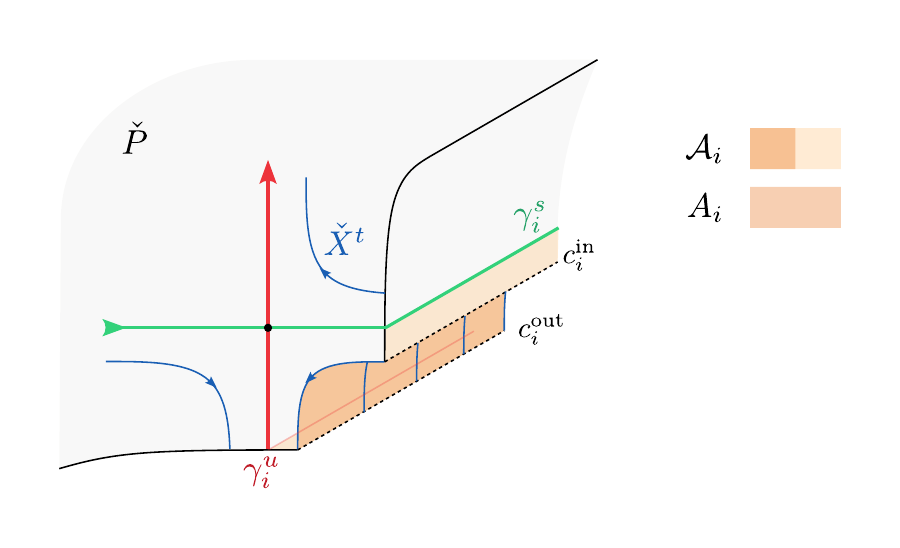}
    \vspace*{-1em}
    \caption{Boundary of an intermediate block $(\check P,\check X)$ in the neighborhood of an annulus $A_i$ tangent to $\check X$}
    \label{fig: intermediate block boundary}
\end{figure}

We refer to Figure~\ref{fig: intermediate block boundary}.

\begin{defi}[Entrance, Exit, Tangent boundary] \label{def: in out tangent boundary of imb}
Let $(\check P,\check X)$ be an \imb{} and $\partial \check P = \check P^\iin \cup \check P^\tan \cup \check P^\out$ be the decomposition of $\partial \check P$ given by Item~\ref{def: imb; it: boundary}, Definition~\ref{def: intermediate block}.
We say that
\begin{itemize}[--]
    \item $\check P^\iin$ is the \emph{entrance boundary} of $(\check P,\check X)$;
    \item $\check P^\out$ is the \emph{exit boundary} of $(\check P,\check X)$;
    \item $\check P^\tan$ is the \emph{tangent boundary} of $(\check P,\check X)$.
\end{itemize}
\end{defi}

\begin{claim} \label{claim: boundary lam on imb}
With these same notations, we have
\begin{enumerate}
    \item The stable lamination $\check \cW^s$ induces a lamination $\check \cL^\iin := \check \cW^s \cap \partial \check P = \check \cW^s \cap \check P^\iin$ of dimension 1 on $\check P^\iin$, called \emph{entrance lamination} of $(\check P,\check X)$.
    \item The unstable lamination $\check \cW^u$ induces a lamination $\check \cL^\out := \check \cW^u \cap \partial \check P = \check \cW^u \cap \check P^\out$ of dimension 1 on $(\check P,\check X)$, called \emph{exit lamination} of $(\check P,\check X)$.
    \item The union $\check \cL = \check \cL^\iin \cup \check \cL^\out$ is a lamination of dimension 1 on $\partial \check P$, called \emph{boundary lamination} of $(\check P,\check X)$.
\end{enumerate}
\end{claim}

\begin{proof}
Each orbit of a point of a connected component $A_i$ of $\check P^\tan$ is contained in the annulus $A_i$ and crosses $A_i$ from one boundary of the annulus to the other in finite time (Item~\ref{def: imb; it: boundary}, Definition~\ref{def: intermediate block}).
We deduce that the stable lamination $\check \cW^s$ of the maximal invariant set $\check \Lambda$ is disjoint from the tangent boundary $\check P^\tan$.
It is disjoint from $\check P^\out$, because all orbits intersecting $\check P^\out$ exit $\check P$ by definition, so cannot accumulate on $\check \Lambda$ in the future.
We deduce that $\check \cW^s \cap \partial \check P = \check \cW^s \cap \check P^\iin$.
Since the vector field $\check X$ is transverse to $\check P^\iin$, and $\check \cW^s$ is tangent to the vector field $\check X$, we deduce that the leaves of $\check \cW^s$ transversely intersect the surface $\check P^\iin$, and the set $\check \cL^\iin := \check \cW^s \cap \check P^\iin$ is a lamination of dimension 1 on $\check P^\iin$.
The other statements of the fact are proved in a similar way.
\end{proof}

Let $A_* = \{A_1, \dots, A_n \}$ be the collection of the \ccs{} of the tangent boundary of $(\check P,\check X)$,
and $\cA_* = \{ \cA_1, \dots, \cA_n \}$ the collection of the \ccs{} $\partial \check P \ssm \check \cL$ such that $A_i \subset \cA_i$, where $\check \cL$ is the boundary lamination of $(\check P, \check X)$.
With these notations, the curves $\gamma^s_i$ and $\gamma^u_i$ are compact leaves of $\check \cL^\iin$ and $\check \cL^\out$ respectively, and $\cA_i$ is the connected component of $\partial \check P \ssm \check \cL$ which contains the annulus $A_i$.

Let $(P,X)$ be a \bb{}. 
Let $\Lambda$ be its maximal invariant set, $\cO_*$ the collection of periodic orbits of $X$ contained in $\partial P$, $\Pin$ the entrance boundary, and $\Pout$ the exit boundary of $(P,X)$.
Recall that the closure of $\Pin$, denoted $\overline \Pin$, coincides with $\Pin \cup \cO_*$.
Similarly, $\overline \Pout = \Pout \cup \cO_*$.
The following lemma states that we can canonically embed $(P,X)$ into an intermediate block $(\check P,\check X)$, which is then a \say{standard} \minc{} of $(P,X)$.

\begin{lem} \label{lem: bb to imb}
There exists an \imb{} $(\check P,\check X)$ such that, if $\check \Lambda$ denotes the maximal invariant set, $(\check \cW^s, \check \cW^u)$ the stable and unstable manifolds of $\check \Lambda$, and $\cA_*$ the collection of \ccs{} of $\partial \check P \ssm \check \cL$ containing the tangent boundary (Definition~\ref{def: in out tangent boundary of imb}), we have
\begin{enumerate}
    \item \emph{(\minc{})} \label{lem: bb to imb, it: minc}
    $(\check P,\check X)$ is a \minc{} of $(P,X)$ (Definition~\ref{def: minc});
    \item \emph{(lamination)} \label{lem: bb to imb; it: lamination}
    $\check P^\iin \ssm \cA_*$ is isotopic to $ \overline \Pin$ along the lamination $\check \cW^s$ and \linebreak[4]$\check P^\out \ssm \cA_*$ is isotopic to $ \overline \Pout$ along the lamination $\check \cW^u$;
    \item \emph{(uniqueness)} \label{lem: bb to imb; it: uniqueness}
    such an \imb{} is unique up to topological equivalence.
\end{enumerate}
\end{lem}

\begin{proof}
We repeat the proof of Lemma~\ref{lem: from bb to bby} which gives a method of constructing an \imb{} $(\check P, \check X)$ associated to a \bb{} $(P,X)$.
Let $(\tilde P, \tilde X)$ be a \minc{} (Definition~\ref{def: minc}) of $(P,X)$.
Let $\cO_* = \{ \cO_1, \dots, \cO_n\}$ be the collection of periodic orbits of $X$ contained in $\pP$.
For each $\cO_i \in \cO_*$, up to local orbit equivalence of the flow of $\tilde X$, we can consider $\cV_i$ a linearizing tubular neighborhood of $\cO_i$ for the flow of $\tilde X$ in $\tilde P$, provided with a \lcs{} $(x,y,\theta) \in \R^2 \times \R/\Z$.
The boundary $\partial P$ crosses the opposite quadrants $\{x > 0, y > 0\}$ and $\{x < 0, y < 0\}$ (Claim~\ref{claim: boundary quadrant and multipliers}).
Let $S$ be a topological surface, smooth outside a finite number of simple closed curves $c_1^\iin, \dots, c_n^\iin$ and $c_1^\out, \dots, c_n^\out$, which coincides with $\pP$ outside of the $\cV_i$ neighborhoods, and decomposes into the union $S = S^\iin \cup S^t \cup S^\out$ (Figure~\ref{fig: bb to imb}), where
\begin{itemize}
    \item $S^\iin$ is a surface with boundary transverse to the \vf{} $X$ and bordered by $c_1^\iin, \dots, c_n^\iin$, which coincides with $\Pin$ outside the union of neighborhoods $\cV_i$;
    \item $S^\out$ is a surface with boundary transverse to the \vf{} $X$ and bordered by $c_1^\out, \dots, c_n^\out$, which coincides with $\Pout$ outside the union of neighborhoods $\cV_i$;
    \item $S^t$ is a union of annuli $A_i$ tangent to the vector field $X$, each contained in the neighborhood $\cV_i$ of $\cO_i$ in the quadrant $\{ x > 0, y < 0\}$, and bounded by the curves $c^\iin_i$ and $c^\out_i$.
\end{itemize}

\begin{figure}[htb]
    \centering
    \vspace*{-2em}
    \captionsetup{width=.84\linewidth}
    \includegraphics[height=0.28\textheight]{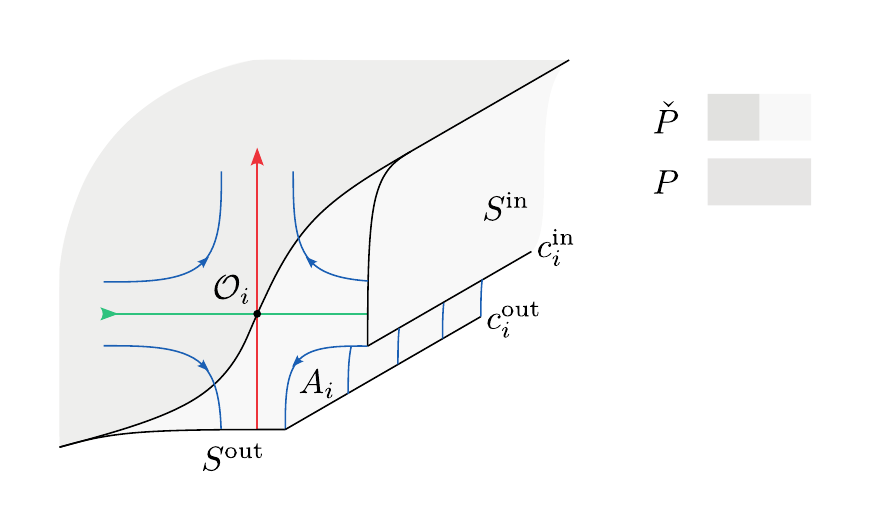}
    \vspace*{-1em}
    \caption{The manifold $\check P$ in a linearizing neighborhood $\cV_i$ of $\cO_i$}
    \label{fig: bb to imb}
\end{figure}

The surface $S$ cuts into $\tilde P$ a compact submanifold $\check P$ with boundary, which contains $P$, and coincides with $P$ outside the union of the neighborhoods $\cV_i$.
Let $\check X$ be the restriction of $\tilde X$ to $\check P$.
The pair $(\check P, \check X)$ is by construction a \minc{} of $(P,X)$.
The maximal invariant set $\check \Lambda$ of $(\check P, \check X)$ coincides with the maximal invariant set $\Lambda$ of $(P,X)$ because an orbit of $\check X$ is either an orbit of the flow of $X$ contained in $P$, or intersects the boundary of $\check P$.
Moreover, $\partial \check P = S$, so $(\check P, \check X)$ satisfies Items~\ref{def: imb; it: edges},~\ref{def: imb; it: boundary} and~\ref{def: imb; it: hyp max inv} of Definition~\ref{def: intermediate block} of an \imb.
Let $\check \cW^s$ and $\check \cW^u$ be the stable and unstable manifolds of the set $\check \Lambda = \Lambda$ for the flow of $\check X$ in $\check P$.
Then $\check \cW^s$ is the union of the flow-saturated set of the laminations $\cW^s$ under the flow of $\check X$ and the local stable manifold of periodic orbits $\cO_i$ in the linearizing neighborhoods $\cV_i$.
Similarly $\check \cW^u$ is the union of the flow-saturated set of the laminations $\cW^u$ under the flow of $\check X$ and the local unstable manifolds of the periodic orbits $\cO_i$ in the linearizing neighborhoods $\cV_i$.
For any $i$, the local stable and unstable manifold of orbit $\cO_i$ in neighborhood $\cV_i$ intersect $\partial \check P$ along two compact leaves $\gamma^s_i$ and $\gamma^u_i$, which borders a connected component $\cA_i$ of $\partial \check P \ssm \check \cL$.
This is an annulus which contains the annulus $A_i \subset S$ tangent to the vector field $\check X$ (Figure~\ref{fig: bb to imb}).
It follows that $(\check P, \check X)$ satisfies Item~\ref{def: imb; it: annuli in boundary} of Definition~\ref{def: intermediate block}, so it is an \imb{}.

Let us denote $\cA_*$ the collection of annuli $\cA_i$.
The surfaces $\check P^\iin \ssm \cA_*$ and $\overline \Pin$ coincide outside the union of neighborhoods $\cV_i$, and in each $\cV_i$ the surface $\check P^\iin \ssm \cA_*$ is an annulus whose interior is transverse to the lamination $\check \cW^s$.
There is therefore an isotopy between $\check P^\iin \ssm \cA_*$ and $\overline \Pin$, supported in the linearized \nbhs{} $\cV_i$, which preserves the lamination $\check \cW^s$ leaf-to-leaf.
The argument is symmetric for $\check P^\out \ssm \cA_*$ and $\overline \Pout$, and $(\check P, \check X)$ satisfies Item~\ref{lem: bb to imb; it: lamination} of Lemma~\ref{lem: bb to imb}.
It remains to prove Item~\ref{lem: bb to imb; it: uniqueness}, in other words the uniqueness of such an \imb{} up to orbit equivalence.
Let $(Q,Y)$ be another \imb{} satisfying Lemma~\ref{lem: bb to imb}.
Up to reduce the initial \minc{} used to construct $(\check P, \check X)$ and make an orbit equivalence, we can make the hypothesis that $\check P$ is contained in $Q$, and that they are contained in a common \minc{} of $(P,X)$ according to Item~\ref{lem: bb to imb, it: minc} of Lemma~\ref{lem: bb to imb}.
It is sufficient to show the existence of a local orbit equivalence on a linearizing neighborhood $\cV_i$ of a periodic orbit $\cO_i$ of $X$ contained in $\pP$, because the \imbs{} $(Q,Y)$ and $(\check P, \check X)$ coincide each with $(P,X)$ outside of the union of $\cV_i$, so coincide between them.
By transversality of $\Pin$ and $Q^\iin$ with the vector field $Y$, up to an isotopy of $\partial Q$ along the (positive) orbits of the flow of $Y$ (which does not change the orbit equivalence class), we can suppose that $Q^\iin$ contains $\check P^\iin$, and up to an isotopy of $\partial \check P$ along the positive orbits of the flow of $Y$, we can suppose that $Q^\out$ contains $\check P^\out$ (Figure~\ref{fig: uniqueness imb}). 

\begin{figure}[htb]
    \centering
    \vspace*{-1.5em}
    \includegraphics[width=\textwidth]{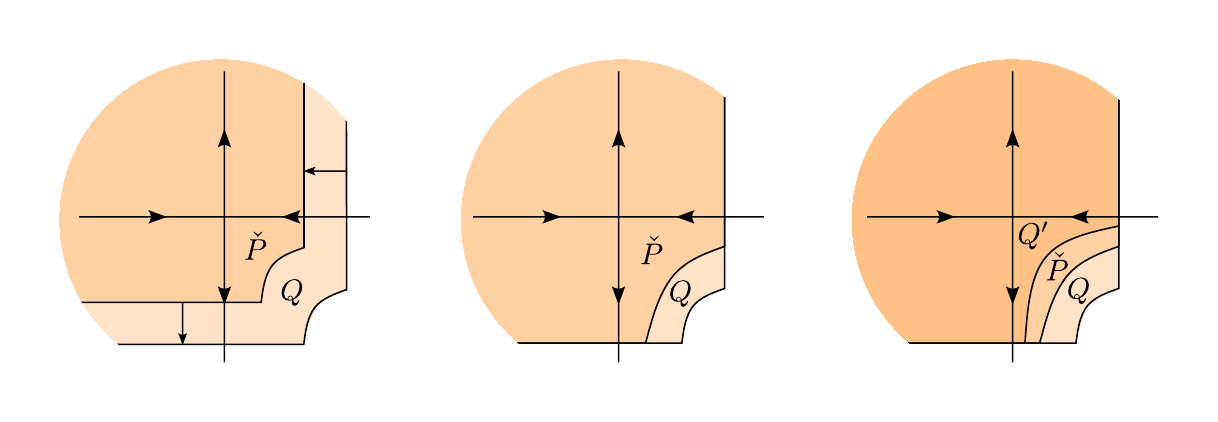}
    \vspace*{-1.5em}
    \caption{Orbit equivalence of \imbs{} in a\\linearizing neighborhood of a periodic orbit $\cO_i$}
    \label{fig: uniqueness imb}
\end{figure}

We introduce a third \imb{} $(Q', Y')$ which satisfies Lemma~\ref{lem: bb to imb}, which is strictly contained in $\check P$, and such that $(Q')^\iin \subset \check P^\iin$ and $(Q')^\out \subset \check P^\out$ (Figure~\ref{fig: uniqueness imb}).
Then the vector field $\check X$ restricted to $\check P \ssm Q'$ and the vector field $Y$ restricted to $Q \ssm Q'$ are both orbit equivalent to the trivial field $\partial /\partial_t$ on $I \times J \times \R/\Z$ given coordinates $(s,t,\theta)$.
Since the \imbs{} $(\check P, \check X)$ and $(Q,Y)$ coincide on $Q'$, we deduce that they are orbit equivalent.
\end{proof}

Let $(\check P, \check X)$ be an \imb{} and $\cA_*$ be the collection of \ccs{} of $(\check P, \check X)$ containing the tangent boundary of $(\check P, \check X)$ (Definition~\ref{def: in out tangent boundary of imb}), where $\check \cL$ is the boundary lamination of $(\check P, \check X)$.
Recall that a stable separatrix is a connected component of $\check \cW^s(O) \ssm O$, where $O$ is an orbit of $\check X$.

\begin{lem} \label{lem: tangent annuli imb}
For any $\cA_i$, if $\gamma^s_i$ and $\gamma^u_i$ are the two compact leaves of $\check \cL$ which border $\cA_i$ (Item~\ref{def: imb; it: boundary}, Definition~\ref{def: intermediate block}), then there exists a periodic orbit $\cO_i$ of $\check X$ such that
$\gamma^s_i$ is the intersection of a free stable separatrix of $\cO_i$ with $\check P^\iin$ and $\gamma^u_i$ is the intersection of a free unstable separatrix of $\cO_i$ with $\check P^\out$.
Up to orbit equivalence of $(\check P, \check X)$, we can suppose that $\cA_i$ is contained in a linearizing neighborhood of~$\cO_i$.
\end{lem}

\begin{proof}
We can adapt Lemma~\ref{lem: free separatrix} for intermediate blocks, and $\gamma^s_i$ is a compact leaf of the lamination $\check \cL^\iin$, thus the intersection of a free stable separatrix $\cW^s_+ (\cO_i)$ of a periodic orbit $\cO_i$ of $\check X$.
If $p \in \check P^\iin$ is a point on $\check P^\iin$ arbitrarily close to $\gamma^s_i$, then the orbit of $p$ by the flow of $\check X$ exits $\check P$ at a point $q \in \check P^\out$ arbitrarily close to the unstable manifold $\cW^u(\cO_i)$ of $\cO_i$.
Now we know that the flow of $\check X$ on the annulus $A_i \subset \cA_i$ crosses $A_i$ from one boundary $c^\iin_i$ to the other boundary $c^\out_i$ (Item~\ref{def: imb; it: boundary}, Definition~\ref{def: intermediate block}).
Since $c^\iin_i$ is the boundary of $\check P^\iin \cap \cA_i$, and $c^\out_i$ is the boundary of $\check P^\out \cap \cA_i$, we deduce by continuity that the orbit of a point of $\check P^\iin \cap \cA_i$ exits $P$ by the component $\check P^\out \cap \cA_i$.
It follows that $\gamma^u_i \in \cW^u(\cO_i)$, and $\gamma^u_i$ is the intersection of a free unstable separatrix $\cW^u_+ (\cO_i)$ of the same periodic orbit $\cO_i$.

Let us show the second statement of Lemma~\ref{lem: tangent annuli imb}. 
For each $i$, consider $c_i$ a simple closed curve in (the interior of) $\check P^\iin \cap \cA_i$ parallel to $\gamma^s_i$ and close enough to $\gamma^s_i$ so that the orbit of $c_i$ by the flow of $\check X$ intersects a linearizing neighborhood $\cV_i$ of $\cO_i$.
Let $C_i := \bigcup_{t \in \R} \check X^t(c_i)$ the orbit of $c_i$ by the flow of $\check X$.
Each point of $c_i$ exits $\check P$ in uniformly bounded time, because $c_i$ is uniformly far from $\check \cL^\iin$.
It follows that $C_i$ is homeomorphic to the cylinder $c_i \times [0,1]$, and the vector field $\check X$ is orbit equivalent to the vertical field $\partial_t$ on $C_i$.
The union of cylinders $C_i$ cuts into $\check P$ a submanifold with boundary $Q$, such that if $Y$ is the restriction of $\check X$ on $Q$, the pair $(Q,Y)$ is an \imb{}, embedded in $(\check P, \check X)$, such that the maximal invariant set of $(Q,Y)$ coincides with the maximal invariant set of $X$.
Let us show that the intermediate block $(Q,Y)$ is orbit equivalent to $(\check P, \check X)$.
Let $(Q', Y')$ be a third \imb{} contained in $(Q,Y)$ obtained by cutting $(P,X)$ along the orbit of a simple closed curve $c'_i$ contained in the interior of the annulus bounded by $\gamma^s_i$ and $c_i$ and parallel to the boundary of this annulus.
Then $Q \ssm Q'$ is homeomorphic to the product $I \times J \times \R/\Z$, and the vector field $Y$ on $Q \ssm Q'$ is orbit equivalent to the vector field $\partial_t$ on $I \times J \times \R/\Z$ given coordinates $(s,t,\theta)$.
Similarly, $\check P \ssm Q'$ is homeomorphic to the product $I \times J \times \R/\Z$, and the vector field $Y$ on $Q \ssm Q'$ is orbit equivalent to the vector field $\partial_t$ on $I \times J \times \R/\Z$ with coordinates $(s,t,\theta)$ (Figure~\ref{fig: uniqueness imb})
It follows that the vector fields $Y$ on $Q \ssm Q'$ and $\check X$ on $\check P \ssm Q'$ are orbit equivalent.
Since they coincide outside $\check P \ssm Q$, we deduce that $(\check P, \check X)$ is orbit equivalent to $(Q,Y)$.
The tangent boundary $Q^\tan$ of $\partial Q$ coincides with the union of the cylinders $\bigcup_i C_i$ defined previously, and we denote $\partial Q = Q^\iin \cup Q^\tan \cup Q^\out$ the decomposition of the boundary given by Item~\ref{def: imb; it: boundary}, Definition~\ref{def: intermediate block}.
By construction, the boundary components $C_i$ intersect a linearizing neighborhood $\cV_i$ of the periodic orbit $\cO_i$.
Let $\cC_i$ be the connected component of $\partial Q \ssm \cL_Y$ containing $C_i$.
By transversality of $Q^\iin$ with the vector field $Y$, up to push $Q^\iin$ along the (positive) orbits of $Y$ and $Q^\out$ along the (negative) orbits of $Y$, we can assume that $\cC_i$ is contained in $\cV_i$.
This operation does not change the orbit equivalence class of $(Q,Y)$.
\end{proof}

\begin{defi}[(h)-hypothesis]
\label{def: (h) hypothesis for imb}
We will say that the $(\check P, \check X)$ satisfies the \emph{$(h)$ hypothesis} if the periodic orbits $\cO_i$ associated to the annuli $\cA_i$ by Lemma~\ref{lem: tangent annuli imb} are pairwise distinct, and have positive multipliers.
\end{defi}

\begin{rmk} \label{rmk: bb to imb (h) hyp}
In the case of a $(\check P, \check X)$ associated to a $(P,X)$ by Lemma~\ref{lem: bb to imb}, the periodic orbits $\cO_i$ are the periodic orbits of $X$ contained in $\partial P$.
Indeed, $(\check P, \check X)$ is a \minc{} of $(P,X)$ (Item~\ref{lem: bb to imb, it: minc}, Lemma~\ref{lem: bb to imb}) and the annuli $\cA_i$ are contained in $\partial \check P \ssm \partial P$, thus in the linearizing neighborhoods $\cV_i$ of the periodic orbits $\cO_*$.
Consequently, Claim~\ref{claim: boundary quadrant and multipliers} implies that each $(\check P, \check X)$ satisfies the $(h)$-hypothesis.
There are some \imbs{} which do not satisfy the $(h)$-hypothesis : see for that Lemma~\ref{lem: elementary block} and Figure~\ref{fig: elementary block}.
\end{rmk}

The following lemma is the counterpart of Lemma~\ref{lem: bb to imb}.
Let $(\check P,\check X)$ be an \imb{} satisfying the $(h)$-hypothesis.
Let $\check \Lambda$ be its maximal invariant set, $(\check \cW^s, \check \cW^u)$ be the pair of stable and unstable manifolds of $\check \Lambda$, $\check P^\iin$ be the entrance boundary, $\check P^\out$ be the exit boundary, and
$\cA_*$ the collection of \ccs{} of $\partial \check P \ssm \check \cL$ containing the tangent boundary of $(\check P, \check X)$ (Definition~\ref{def: in out tangent boundary of imb}).

\begin{lem} \label{lem: imb to bb}
There exists a \bb{} $(P,X)$ such that, if $\Lambda$ denotes the maximal invariant set of $(P,X)$, $\cO_*$ the collection of periodic orbits of $X$ contained in $\pP$, $\Pin$ the entrance boundary and $\Pout$ the exit boundary of $(P,X)$, then:
\begin{enumerate}
    \item \emph{(embedding)} \label{lem: imb to bb; it: embedding}
    there exists an embedding $h: P \to \check P$ which maps the oriented orbits of $X$ on the oriented orbits of $\res{\check X}{\check P}$ and $\Lambda$ on $\check \Lambda$;
    \item \emph{(lamination)} \label{lem: imb to bb; it: lamination}
    $\check P^\iin \ssm \cA_*$ is isotopic to $\overline \Pin$ along the lamination $\check \cW^s$ and \linebreak[4]$\check P^\out \ssm \cA_*$ is isotopic to $ \overline \Pout$ along the lamination $\check \cW^u$;
    \item \emph{(uniqueness)} \label{lem: imb to bb; it: uniqueness}
    such a \bb{} $(P, X)$ is unique up to isotopy among the \bbs{}.
    \end{enumerate}
\end{lem}
We refer to Definition~\ref{def: isotopic blocks} for the isotopy of \bbs{}.

\begin{proof}
In order to make an orbit equivalence of $(\check P, \check X)$, we assume that each annulus $\cA_i$ is contained in the linearizing neighborhood $\cV_i$ of a periodic orbit $\cO_i$ (Lemma~\ref{lem: tangent annuli imb}).
By assumption, the multipliers of $\cO_i$ are positive, and the orbits $\cO_i$ are pairwise disjoint, so we can assume that the linearizing \nbhs{} $\cV_i$ are pairwise disjoint.
The entrance boundary $\check P^\iin$ and the exit boundary $\check P^\out$ of $(\check P, \check X)$ lie in two opposite quadrants of $\cO_i$ on the neighborhood $\cV_i$ (see the proof of Claim~\ref{claim: boundary quadrant and multipliers}).
So there exists a smooth surface $S$ in $\check P$, which coincides with $\partial \check P$ outside the union of $\cV_i$, and such that in each $\cV_i$, $S$ is an annulus which contains the orbit $\cO_i$, which coincides with $\partial \check P$ in the \nbh{} of the boundary of $\cV_i$ and crosses $\cV_i$ in two opposite quadrants of $\cO_i$, and which is transverse to the vector field $\check X$ on the complementary of $\cO_i$.
Then $S$ cuts into $\check P$ a submanifold $P$ such that if $X$ is the restriction of $\check X$ to $P$, then $(P,X)$ is a \bb{}, embedded in $(\check P, \check X)$ and such that the maximal invariant set $\Lambda$ of $(P,X)$ coincides with the maximal invariant set $\check \Lambda$ of $(\check P, \check X)$.
We show Item~\ref{lem: imb to bb; it: lamination} in an analogous way to the proof of Lemma~\ref{lem: bb to imb}, Item~\ref{lem: bb to imb; it: lamination}.
Finally, two \bbs{} $(P,X)$ and $(P',X')$ which satisfy Lemma~\ref{lem: imb to bb} have a common \minc{} (up to orbit equivalence of $(\check P, \check X)$, according to Lemma~\ref{lem: tangent annuli imb}).
They are therefore isotopic among the \bbs{} (Proposition~\ref{prop: isotopy vs orbit eq}).
\end{proof}

\subsection{Building block and \bby{} block}

\label{sec: geometric type; subsec: bb and bby}

Let $(P,X)$ be a \bb{}.
Let $\cO_*$ be the set of periodic orbits of $X$ contained in $\pP$, $\Lambda_X$ the maximal invariant set, $P^\iin$ the entrance boundary and $P^\out$ the exit boundary.

\begin{lem} \label{lem: bb to bby}
There exists a \bby{} block $(U,Y)$ such that, if $\Lambda_Y$ denotes the maximal invariant set of $(U,Y)$, 
$(\cW_Y^s, \cW_Y^u)$ the pair of stable and unstable manifolds of $\Lambda_Y$, 
$U^\iin$ the entrance boundary, $U^\out$ the exit boundary, and $\cL_Y = \cL^\iin_Y \cup \cL^\out_Y$ the boundary lamination, then:
\begin{enumerate}
    \item \emph{(embedding)} \label{lem: bb to bby; it: embedding}
    there exists an embedding $h : P \to U$ which maps the oriented orbits of $X$ to the oriented orbits of $Y$ in $h(P)$ and $\Lambda_X$ to $\Lambda_Y$;
    \item \emph{(lamination)}\label{lem: bb to bby; it: lamination}
    there exists a finite collection $D_* = \{D_1, \dots, D_n\}$ of open disks contained in $\partial U$, disjoint from $\cL_Y$, bounded by pairwise distinct compact leaves of $\cL_Y$, such that $ U^\iin \ssm D_*$ is isotopic to $\overline \Pin$ along the lamination $\cW^s_Y$ and $U^\out \ssm D_*$ is isotopic to $ \overline \Pout$ along the lamination $\cW^u_Y$;
    \item \emph{(uniqueness)}  \label{lem: bb to bby; it: uniqueness}
    such a block $(U,Y)$ is unique up to orbit equivalence
\end{enumerate}
\end{lem}

\begin{proof}
We refer to the proof of Lemma~\ref{lem: from bb to bby} for Items~\ref{lem: bb to bby; it: embedding} and~\ref{lem: bb to bby; it: lamination}.
Let us show the last Item~\ref{lem: bb to bby; it: uniqueness}.
Let $(U,Y)$ be a \bby{} block which satisfies Lemma~\ref{lem: bb to bby}.
Let $\bar D_1, \dots, \bar D_n$ be a collection of compact disks such that $\bar D_i \subset D_i$.
Since each point $p$ is uniformly far from the boundary lamination $\cL_Y$, the orbit of $p$ exits $U$ in uniformly bounded time.
We deduce that the orbit of $\bar D_i$ by the flow of $Y$ into $U$, denoted $\bar C_i$, is homeomorphic to a compact cylinder $\D^2 \times I$, on which the vector field $Y$ is orbit equivalent to the vertical vector field $\partial/\partial t$, where $t$ is the coordinate on $I$.
The complementary $\check P := U \ssm \bigcup_i \bar C_i$ equipped with the restriction $\check X$ of the vector field $Y$ is an \imb{} which satisfies Lemma~\ref{lem: bb to imb} associated to the block $(P,X)$.
The orbit equivalence class of the block $(\check P,\check X)$ obtained by such a surgery does not depend on the choice of the compact disk $\bar D_i$ contained in $D_i$ since they are all isotopic.
We conclude with Item~\ref{lem: bb to imb; it: uniqueness} of Lemma~\ref{lem: bb to imb}.
\end{proof}

We now describe the reciprocal operation.

\begin{defi}[Simple collection] \label{def: simple collection}
Let $(U,Y)$ be a \bb{} and $D_* = \{D_1, \dots, D_n\}$ be a collection of connected components of the complementary of the boundary lamination $\partial U \ssm \cL_Y$ bounded by compact leaves of $\cL_Y$.
We say that $D_*$ is a \emph{simple collection} if two distinct compact leaves of the lamination $\cL_Y$ contained in the boundary of $D_*$ are on distinct leaves of $\cW_Y^s$ or $\cW_Y^u$, and the periodic orbits contained in these leaves have positive multipliers.
\end{defi}

Recall that any compact leaf $\gamma$ of the boundary lamination $\cL_Y$ is the trace of a free separatrix of a periodic orbit $O$ on $\partial U$ (Lemma~\ref{lem: free separatrix}).
This justifies that we can speak of (the only) periodic orbit contained in the leaf of $\cW^s_Y$ (or $\cW^u_Y$) which contains $\gamma$.

\begin{rmk} \label{rmk: collection bb to bby is simple}
Let $D_*$ be the collection given by Lemma~\ref{lem: bb to bby}, Item~\ref{lem: bb to bby; it: lamination}.
Then it is a simple collection.
The isotopy of Item~\ref{lem: bb to bby; it: lamination} induces a \diff{} $H^s: U^\iin \ssm D_* \to \overline P^\iin$ which preserves the stable lamination $\cW^s_Y$ leaf to leaf.
If $\gamma^s$ is a compact leaf of $\cL^\iin_Y$, boundary of a component $D_i$, then $\gamma^s$ is mapped by $H^s$ on a periodic orbit of $\cO \in \cO_*$ along the stable leaf $\cW_Y^s(\cO)$.
The multipliers of $\cO$ are positive according to Claim~\ref{claim: boundary quadrant and multipliers}.
By injectivity of $H^s$ it is the only boundary component of $D_*$ mapped on $\cO$.
It follows that the leaves bordering $D_*$ are on disjoint leaves of $\cW_Y^s$.
Similarly, the compact leaves of $\cL^\out_Y$ which border a component of $D_*$ are on disjoint leaves of $\cW_Y^u$.
\end{rmk}

The following lemma is a reciprocal of Lemma~\ref{lem: bb to bby}.
Let $(U,Y)$ be a \bby{} block, and let $\Lambda_Y$ be the maximal invariant set, $\cL_Y$ be the boundary lamination of $(U,Y)$, $U^\iin$ be the entrance boundary and $U^\out$ be the exit boundary.
Let $D_* = \{D_1, \dots D_n \}$ be a simple collection (Definition~\ref{def: simple collection}) of connected components of $\partial P \ssm \cL$ bounded by compact leaves of $\cL_Y$ and invariant by the crossing map $f_Y: U^\iin \to U^\out$ of the flow of $Y$ from the entrance boundary to the exit boundary.

\begin{lem}
\label{lem: bby to bb}
There exists a \bb{} $(P,X)$ such that, if $\cO_*$ denotes the set of periodic orbits of $X$ contained in $(P,X)$, $\cL_X$ the boundary lamination of $(P,X)$ and $\Lambda_X$ the maximal invariant set of $(P,X)$, we have:
\begin{enumerate}
    \item \emph{(embedding)} \label{lem: bby to bb; it: embedding}
    there exists an embedding $h : P \to U$ which maps the oriented orbits of $X$ to the oriented orbits of $Y$ in $h(P)$ and $\Lambda_X$ to $\Lambda_Y$;
    \item \emph{(lamination)}\label{lem: bby to bb; it: lamination}
    $U^\iin \ssm D_*$ is isotopic to $ \overline \Pin$ along the lamination $\cW^s_Y$, and $U^\out \ssm D_*$ is isotopic to $ \overline \Pout$ along the lamination $\cW^u_Y$;
    \item \emph{(uniqueness)} \label{lem: bby to bb; it: uniqueness}
    such a building block $(P,X)$ is unique up to isotopy among the \bbs{}.
\end{enumerate}
\end{lem}

\begin{rmk} \label{rmk: bb and bby reciprocal}
Let us justify why Lemma~\ref{lem: bby to bb} does give a reciprocal of Lemma \ref{lem: bb to bby}.
If $(U,Y)$ is the \bby{} block associated to a building block $(P,X)$ by Lemma~\ref{lem: bb to bby}, and $D_*$ is the collection of disks given by Item~\ref{lem: bb to bby; it: lamination} of this lemma, then it is a simple collection, and we can apply Lemma~\ref{lem: bby to bb} to $(U,Y)$ and to the collection $D_*$.
We obtain a \bb{} $(P', X')$ isotopic to $(P,X)$ by Item~\ref{lem: bby to bb; it: uniqueness}.
\end{rmk}

\begin{proof}[Proof of Lemma~\ref{lem: bby to bb}]
By $f_Y$-invariance of $D_*$, we can rename the components of $D_*$ as follows $D_* = \{ D^\iin_1, \dots, D^\iin_n, D^\out_1, \dots, D^\out_n \}$ such that $D^\iin_i \subset U^\iin$, and $D^\out_i = f_Y (D^\iin_i) \subset U^\out$.
Let $\gamma_i^s$ be a compact leaf that borders $D^\iin_i$.
It is the intersection of a free stable separatrix of a periodic orbit $O_i$.
By $f_Y$-invariance of $D_*$, there exists a compact leaf $\gamma_i^u$ in the boundary of $D^\out_i$ which is on the unstable manifold of the periodic orbit $O_i$ (it suffices to notice that the orbit of a point in $D^\iin_i$ arbitrarily close to $\gamma_i^s$ exits through $U^\out$ on $D^\out_i$ at a point arbitrarily close to $\cW^u(O_i)$).
If we push $U^\iin$ and $U^\out$ along the orbits of the flow of $Y$ we can assume that a neighborhood of $\gamma^s$ in $\partial U$ (contained in $U^\iin$) and a neighborhood of $\gamma^u$ in $\partial U$ (contained in $U^\out$) are contained in a linearizing neighborhood $\cV_i$ of $O_i$, and any orbit that enters through $D^\iin_i$ close to $\gamma_i^s$ exits through $D^\out_i$ close to $\gamma_i^u$ and is entirely contained in $\cV_i$.
We refer to Figure~\ref{fig: bby to imb}.

\begin{figure}[htb]
    \centering
    \vspace*{-1em}
    \captionsetup{width=.85\linewidth}
    \includegraphics[height=0.27\textheight]{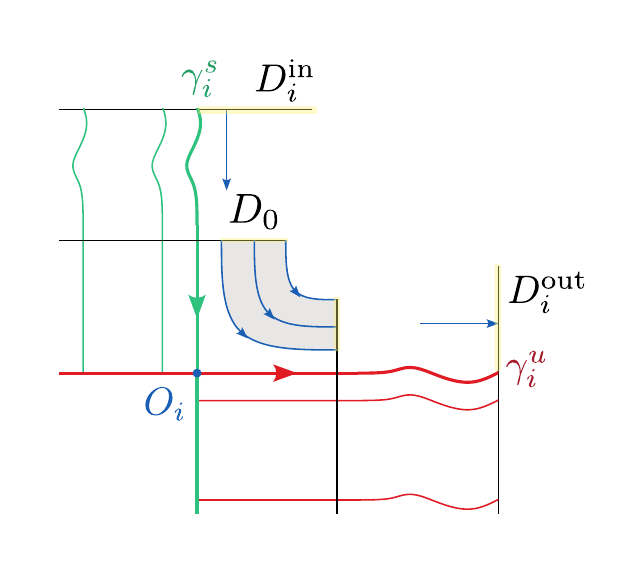}
    \vspace*{-1em}
    \caption{Isotopy along the flow and handles removal in $(U,Y)$}
    \label{fig: bby to imb}
\end{figure}

By repeating this same argument, we obtain a collection $O_1, \dots, O_n$ of periodic orbits and linearizing neighborhoods $\cV_1, \dots, \cV_n$.
Since the collection is simple, the periodic orbits and the linearizing neighborhoods are pairwise distinct.
Let $D_0 \subset U^\iin$ be a compact set contained in the union $\bigcup_i D^\iin_i$, whose interior is isotopic to $\bigcup_i D^\iin_i$ and whose boundary is a union of simple closed curves contained in $\bigcup_i \cV_i$.
By construction, the orbit by the flow of $Y$ of each of the connected components of $\partial D_0$ is a cylinder contained in a single $\cV_i$.
The orbit $\cD_0 := \bigcup_{t \in \R} Y^t(D_0)$ of $D_0$ by the flow of $Y$ in $U$ is homeomorphic to $D_0 \times [0,1]$ and the vector field $Y$ is orbit equivalent to the trivial vector field $\partial /\partial_t$, $t \in [0,1]$.
Let $\check P := U \ssm \cD_0$ and $\check X$ be the restriction of $Y$ to $\check P$.
Then $(\check P, \check X)$ is an \imb{} embedded in $(U,Y)$ and whose maximal invariant set coincides with $\Lambda_Y$.
Indeed by construction, the tangent boundary $\check P^\tan$ of $(\check P, \check X)$ is a union of annuli $A_i$ homeomorphic to $\partial D_0 \times [0,1]$, each contained in a single linearizing neighborhood $\cV_i$.
The connected component $\cA_i$ of $\check P \ssm \check \cL$ containing $A_i$ is an annulus bounded by the compact leaves $\gamma^s_i$ and $\gamma^u_i$ arising from the free stable and local free unstable separatrices of $O_i$.
Moreover, $(\check P, \check X)$ satisfies the $(h)$-hypothesis (Definition~\ref{def: (h) hypothesis for imb}).
We use Lemma~\ref{lem: imb to bb} which gives the existence of a \bb{} $(P,X)$ embedded in $(\check P, \check X)$, such that $(\check P, \check X)$ is a $(P,X)$ \minc{}.
Let us show that $(P,X)$ satisfies Lemma~\ref{lem: bby to bb}.

\begin{enumerate}[leftmargin=*]
    \item By construction, the \imb{} $(\check P, \check X)$ is contained in $(U,Y)$.
    Recall that $\check P = U \ssm \cD_0$, and each orbit of $Y$ in $\cD_0$ enters and exits $\cD_0$, hence $U$, in finite time.
    The maximal invariant set $\Lambda_Y$ of $(U,Y)$ is thus contained in $\check P$, and coincides with the maximal invariant set of $(\check P, \check X)$.
    Item~\ref{lem: bby to bb; it: embedding} follows from Item~\ref{lem: imb to bb; it: embedding} of Lemma~\ref{lem: imb to bb}.
    
    \item We have the equality $\partial \check P \ssm \cA_* = \partial U \ssm D_*$, and the lamination $\check \cL$ on $\partial \check P \ssm \cA_*$ coincides with the lamination $\cL_Y$ on $\partial U \ssm D_*$, because $\check \cL$ is the intersection of the manifolds $\cW^s_Y$ and $\cW^u_Y$ with $\partial \check P$.
    Item~\ref{lem: bby to bb; it: lamination} follows from Item~\ref{lem: imb to bb; it: lamination} of Lemma~\ref{lem: imb to bb}.
    
    \item Let $(P', X')$ be another \bb{} which satisfies Lemma~\ref{lem: bby to bb}.
    Then, up to make an orbit equivalence, according to Item~\ref{lem: bby to bb; it: embedding} of Lemma~\ref{lem: bby to bb}, we can assume that $P$ and $P'$ are contained in $U$, and $X$ and $X'$ coincide with the restrictions of $Y$ to $P$ and $P'$, and the vector fields coincide on a neighborhood of $\Lambda_X = \Lambda_{X'} = \Lambda_Y$.
    In particular, $X$ and $X'$ coincide on linearizing neighborhoods $\cV_i$ of $\cO_i \in \cO_*$ periodic orbits contained in the boundary of $P$.
    According to Item~\ref{lem: bby to bb; it: lamination}, $P^\iin$ is isotopic to $U^\iin \ssm D_*$ along the $\cW^s_Y$ lamination.
    On the complementary of the $\cV_i$ neighborhoods, the surfaces $P^\iin$ and $U^\iin$ are compact and transverse to the \vf{} $Y$, so we can choose an isotopy which preserves the orbits of the $Y$ flow on the complementary of the \nbhs{} $\cV_i$.
    With the same argument for $P'$, we can assume (up to orbit equivalence) that the boundary of $P$ and the boundary of $P'$ coincide outside the neighborhoods $\cV_i$.
    It follows that the blocks $(P,X)$ and $(P', X')$ have a common \minc{}.
    According to Proposition~\ref{prop: isotopy vs orbit eq} they are isotopic.\qedhere
\end{enumerate}
\end{proof}

\begin{rmk} \label{rmk: generalize bby to bb}
We can generalize the previous lemma to a \bb{} $(U,Y)$ which is not a \bby{} block, in other words which admits a non-empty collection $\cO_*$ of periodic orbits contained in the boundary.
We obtain the same result if we suppose moreover that the closure of the elements of $D_*$ are disjoint from $\cO_*$.
The proof is identical because the surgeries are far from~$\cO_*$.
\end{rmk}

\subsection{Proof of Theorem~\ref{thmintro: geometric type in block}}
\label{sec: geometric type; subsec: proof}

We show the realization criterion of an abstract geometric type in a filled block.

\subsubsection*{Step 1: Necessary condition.}
Suppose that $(P,X)$ is a filled orientable \bb{} which admits a Markov partition of geometric type $\scT$.
Let $(U,Y)$ be the \bby{} block associated by Lemma~\ref{lem: bb to bby}.
Let $D_* = \{D_1, \dots, D_n \}$ be the collection of connected components of $\partial U \ssm \cL_Y$ given by Item~\ref{lem: bb to bby; it: lamination} of Lemma \ref{lem: bb to bby}.
Let us show that $(U,Y)$ is a model block.
Let $c$ be a simple closed curve embedded in $\partial U$ disjoint from $\cL_Y$.
If $c$ is contained in a disk $D_i \in D_*$, then $c$ borders a disk disjoint from the lamination $\cL_Y$.
Otherwise $c$ is contained in $U^\iin \ssm D_*$ or $U^\out \ssm D_*$.
Suppose that $c$ is contained in $U^\iin \ssm D_*$.
It follows from Item~\ref{lem: bb to bby; it: lamination} that the lamination $\cL_Y$ restricted to $U^\iin \ssm D$ is topologically equivalent to the lamination $\cL^\iin_X \cup \cO_*$ on $\overline{\Pin}$.
Since the lamination $\cL_X$ is filling on $\partial P$, it follows that any connected component of the complementary of the lamination $\cL^\iin_X \cup \cO_*$ on $\overline{\Pin}$ is a strip.
By topological equivalence, the same is true for the lamination $\cL_Y$ restricted to $U^\iin \ssm D$, and the curve $c$ is contained in a strip.
It follows that $c$ borders a disk disjoint from the lamination.
Since the boundary lamination of $(P,X)$ is filling, it intersects every connected component of $\pP$.
The same is true for the boundary lamination $\cL_Y$ of $(U,Y)$ from Item~\ref{lem: bb to bby; it: lamination}, and we deduce that every \cc{} of $U$ intersects the maximal invariant set $\Lambda_Y$.
Therefore, $(U,Y)$ is a model block.
According to Item~\ref{lem: bb to bby; it: embedding}, the block $(P,X)$ is embedded in $(U,Y)$, it follows that the block $(U,Y)$ admits a Markov partition of geometric type $\scT$.
Let us show that it satisfies the two conditions of Theorem~\ref{thm: geometric type in block}.
The lamination $\cL^\iin_Y$ restricted to $U^\iin \ssm D_*$ is topologically equivalent to the lamination $\cL^\iin_X \cup \cO_*$ on $\overline{\Pin}$, hence is a prefoliation, which shows Item~\ref{prop: geometric type; it: prefoliation}.
It follows that $U^\iin \ssm D_*$ is a union of tori and closed annuli bounded by compact leaves of $\cL^\iin_Y$.
Let $A$ be such an annulus, and let $\gamma_i, \gamma_j$ be the two (distinct) compact leaves of $\cL^\iin_Y$ which border $A$.
There exist two disks $D_i$ and $D_j$ of $D_*$ such that $\gamma_i = \partial D_i$ and $\gamma_j = \partial D_j$.
We conclude that the annulus $A$ is contained in a $\partial U$ which is a sphere $S_{i,j}$, obtained by joining $A$ and the two disks $D_i$ and $D_j$ along their boundary.
Item~\ref{prop: geometric type; it: boundary components} is satisfied.

\subsubsection*{Step 2: Sufficient condition}

Let $\scT$ be a geometric type whose model $(U,Y)$ satisfies Theorem~\ref{thm: geometric type in block}.

\begin{claim} \label{claim: filling lam on model block}
The lamination $\cL_Y$ on $\partial U \ssm D_*$ is filling.
\end{claim}

\begin{proof}
According to Item~\ref{prop: geometric type; it: prefoliation}, the \ccs{} of $\partial U \ssm D_*$ are annuli bordered by compact leaves of $\cL_Y$ or tori, and according to Proposition~\ref{prop: complementary of a prefoliation}, the connected components of the complementary of the lamination $\cL_Y$ on $\partial U \ssm D_*$ are annuli bordered by compact leaves, or tori.
Now let us notice that if $A$ is a \cc{} of $\partial U \ssm \cL_Y$ which is an annulus, then there exists a simple closed curve $c$ embedded in $A$ disjoint from the lamination $\cL_Y$ and which does not border any disk disjoint from the lamination $\cL_Y$, which contradicts with the fact that $(U,Y)$ is a model block (Item~\ref{def: model block; it: lam}, Definition~\ref{def: model block}).
\end{proof}

\begin{claim}
The collection $D_*$ given by Item~\ref{prop: geometric type; it: prefoliation} is a simple and $f_Y$-invariant collection (Definition~\ref{def: simple collection})
\end{claim}

\begin{proof}
The collection $D_*$ is the collection of all connected components of the complementary of $\partial U \ssm \cL_Y$ which are not strips, so it is $f_Y$-invariant.
Suppose there are two compact leaves $\gamma^s_i = \partial D_i$ and $\gamma^s_j = \partial D_j$ in the boundary of $D_*$ from the same stable manifold $\cW^s(O)$ of a periodic orbit $O$.
Then $f_Y (D_i)$ and $f_Y (D_j)$ are two disks of $D_*$ bordered by two distinct compact leaves, both arising from the unstable manifold $\cW^u(O)$ of $O$.
It follows that all separatrices of $O$ are free, so $O$ is an isolated periodic orbit.
The compact leaves induced on the boundary by the separatrices are therefore isolated in the boundary lamination.
This is impossible because we have shown that the lamination $\cL_Y$ is filling on $\partial U \ssm D_*$.
\end{proof}

We can therefore apply Lemma~\ref{lem: bby to bb}, with the collection $D_* = \{D_1, \dots, D_n \}$ of disks given by Item~\ref{prop: geometric type; it: boundary components} of Theorem~\ref{thm: geometric type in block}.
It is a simple and $f_Y$-invariant collection.
We obtain a block $(P,X)$ embedded in $(U,Y)$ by an embedding which maps the maximal invariant set $\Lambda_X$ to $\Lambda_Y$ and the oriented orbits of $X$ to the oriented orbits of $Y$.
It follows that $(P,X)$ is an orientable block, and $\Lambda$ admits a Markov partition of geometric type $\scT$.
It follows from Item~\ref{lem: bby to bb; it: lamination} that the lamination $\cL^\iin_X \cup \cO_*$ on $\overline{\Pin}$ is topologically equivalent to the lamination $\cL^\iin_Y$ restricted to $U^\iin \ssm D_*$, hence is filling according to Claim~\ref{claim: filling lam on model block}.
It is the same for $\cL^\out_X \cup \cO_*$, so for the boundary lamination $\cL_X = \cL^\iin_X \cup \cO_* \cup \cL^\out_X$ on $\partial P$, so $(P,X)$ is a filled block.

\subsubsection*{Step 3: Uniqueness and topological equivalence of laminations.}
Let $(P,X)$ and $(P',X')$ be two \bbs{} satisfying Theorem~\ref{thm: geometric type in block}.
Let $(U,Y)$ and $(U', Y')$ be the \bby{} blocks associated by Lemma~\ref{lem: bb to bby} to $(P,X)$ and $(P',X')$ respectively.
Then as noticed in the first step, $(U,Y)$ and $(U', Y')$ are model blocks.
Moreover, they admit a Markov partition of geometric type $\scT$.
By Theorem~\ref{thm: model geometric type}, the model blocks $(U,Y)$ and $(U',Y')$ are orbit equivalent.
By Lemma~\ref{lem: bby to bb}, seen as a reciprocal of Lemma~\ref{lem: bb to bby} (Remark~\ref{rmk: bb and bby reciprocal}) we deduce that the blocks $(P,X)$ and $(P',X')$ are isotopic.
The topological equivalence of the laminations follows directly from Item~\ref{lem: bby to bb; it: lamination} of Lemma~\ref{lem: bby to bb} applied to $(U,Y)$ and to the collection $D_*$.

\subsection{Example of model of Markov partition}
\label{sec: geometric type; subsec: example model}

Theorem~\ref{thm: geometric type in block} is useful because the criteria can be checked by a algorithm consisting of gluing rectangles.
It allows us, given an abstract geometric type, to build the entrance boundary and the entrance lamination of the model.
We refer to the proof of \cite[Lemma~2.13]{beguinFlotsSmaleDimension2002}, and repeat the idea here.

Let $\scT = \{n, \{h_i\}, \{v_i\}, \phi\}$ be an abstract geometric type.
Consider a union $R$ of $n$ rectangles $R_1, \dots, R_n$, provided with an orientation of the verticals and horizontals of $R_i$, such that each $R_i$ contains $h_i$ horizontal subrectangles $H_i^1, \dots, H_i^{h_i}$ and $v_i$ vertical subrectangles $V_i^1, \dots, V_i^{v_i}$.
Denote $R$ the union of $R_i$, $H$ the union of $H_i^j$ and $V$ the union of $V_i^j$.
Let $\partial^u$ denote the vertical boundary of a rectangle or subrectangle of $R$, and let $\partial^u_g$ and $\partial^u_d$ denote the left and right vertical boundary of the subrectangle respectively relative to the orientations of the verticals and horizontals of $R$.
Let $(i,j)$ be a fixed pair with $1 \leq i \leq n$ and $1 \leq j \leq v_i$.
Let $(k,l)$ be the unique pair such that
$\phi(k,l) = ((i,j), \epsilon)$.
Let $h_{(i,j)}$ be a homeomorphism which maps $\partial^u_g V_i^j$ on $\partial^u_g H_k^l$ and $\partial^u_d V_i^j$ on $\partial^u_d H_k^l$ if $\epsilon =+$ and which maps $\partial^u_g V_i^j$ on $\partial^u_d H_k^l$ and $\partial^u_d V_i^j$ on $\partial^u_s H_k^l$ if $\epsilon = -1$.
The set $A := R \ssm V$ is a finite union of semi-open vertical subrectangles of $R$, whose closure are pairwise disjoint.
The closure $\bar A$ contains the set $\partial^u V$ and the set $\partial^u H$ in its boundary.
Consider the set $\Sigma^s:= \bar A / h$ where $h$ is the product of homeomorphisms $h_{(i,j)}$.
Then $h$ connects the rectangles of $\bar A$ along a finite number of segments of their boundaries.
The surface $\Sigma^s$ is a surface with corners and boundary.
It is provided with a foliation $\cF^s$, image by the gluing map $h$ of the horizontal foliation $R$ induced on $\bar A$.

A Markov partition such that the intersection of every orbit with the suspension \nbh{} of $R$ is connected is called \emph{essential}.
Such a partition always exists in a block without attractors or repellers (Lemma~\ref{lem: existence of markov partition for saddle block}).

\begin{claim}
Let $(U, Y)$ be a model block with an essential Markov partition $\cR$ of geometric type $\scT$.
Then $\Sigma^s$ is a compact surface, with boundary and corners, homeomorphic to the intersection of the flow-saturated set of $\cR$ under the flow of $Y$ on the entrance boundary $U^\iin$.
The lamination $\cL_Y$ on $U^\iin$ is topologically equivalent to a sublamination of the foliation $\cF^s$ on $\Sigma^s$.
The boundary $U^\iin$ is obtained by gluing disks to the boundary of $\Sigma^s$.
\end{claim}

\begin{proof}[Proof idea]
The idea of the proof is based on the following two observations, using the previous notations
\begin{itemize}
    \item The past orbit of a point of $\bar A = R \ssm V$ intersects the boundary $U^\iin$ in uniformly bounded time and at a single point.
    \item The orbits of the points of $\intr A$ are pairwise distinct,
    \item The negative orbit of a point of $\bar A$ reconnects with a negative orbit of a point of the vertical boundary of $A$ after time 1.
\end{itemize}
With these observations, we see that the gluing map $h$ on the set $\bar A \ssm A$ corresponds to the identification of points on the same orbit.
It then follows from the first two items that the flow-saturated set of $R$ under the negative flow of $Y$ intersects the boundary of $U$ on a surface homeomorphic to the quotient of $\bar A$ by $h$.
A surface obtained by gluing a finite number of segments, immersed in a smooth surface, is a compact surface with boundary and corner.
We refer to the proof \cite[Lemma~2.13]{beguinFlotsSmaleDimension2002} for a complete argument.
The stable lamination on $R$ is a sublamination of the horizontal foliation on $R$, which justifies the second statement.
It follows from the criterion~\ref{def: model block; it: lam} of Definition~\ref{def: model block} of a model block that the complementary $U^\iin \ssm \Sigma^s$ is a finite union of disks (otherwise we would have a circle disjoint from the lamination which borders no disk).
\end{proof}

\begin{example}[Fake horseshoe] \label{ex: model fake horseshoe}
Recall that the geometric type of the fake horseshoe is given in Figure~\ref{fig: partition fake horseshoe}, and the example~\ref{ex: type of horseshoe}.
The following figure shows the construction of the entrance boundary $U^\iin$ and the entrance lamination $\cL^\iin$ of the horseshoe model by gluing rectangles.
We see that $U^\iin$ is a sphere containing exactly two disks $D_1$ and $D_2$ of interiors disjoint from the boundary lamination $\cL$, and bordered by two distinct compact leaves of $\cL^\iin$.
The complementary of these disks in $U^\iin$ is an annulus bordered by two compact leaves, on which the lamination is a (filling) prefoliation, with no compact leaves other than the two boundary leaves, and each noncompact leaf accumulates on a compact leaf with no Reeb component.
In particular, the fake horseshoe model satisfies the criterion of Theorem~\ref{thm: geometric type in block}, so the geometric type of the fake horseshoe is realizable in a orientable filled \bb{} $(P,X)$.

\begin{figure}[htb]
    \centering
     \vspace*{-1em}
    \includegraphics[width=\textwidth]{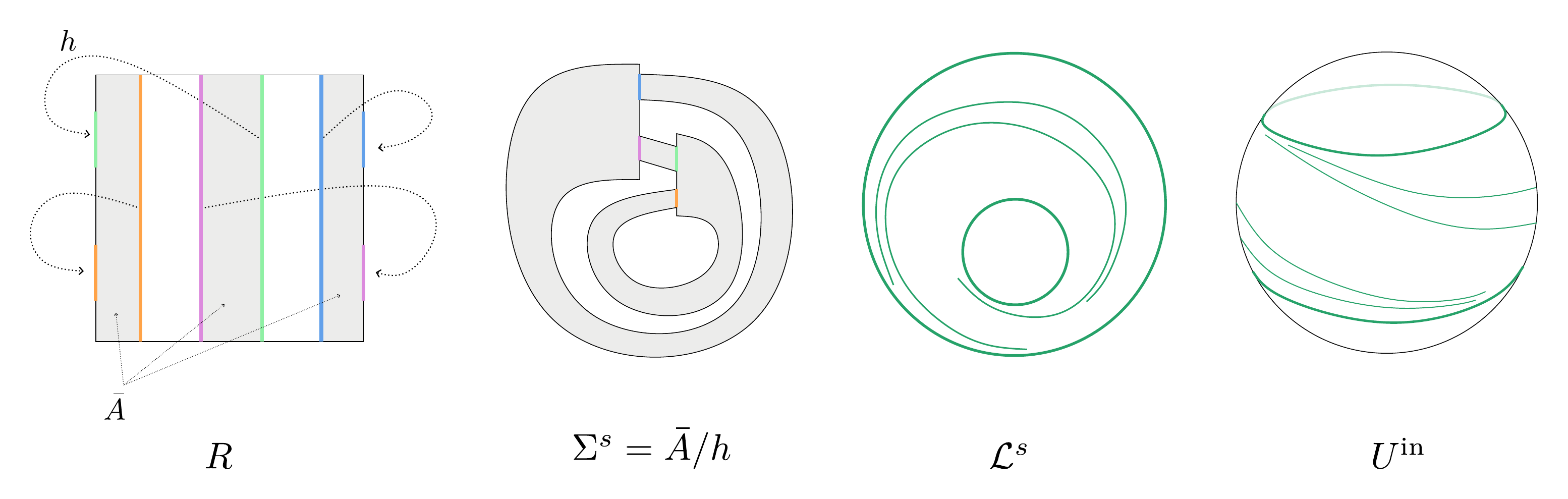}
     \vspace*{-1em}
    \caption{Constructing the surface $\Sigma^s$ and the boundary lamination on the entrance boundary $U^\iin$ of the model $(U,Y)$ of the fake horseshoe}
    \label{fig: entrance boundary fake horseshoe}
\end{figure}

The boundary lamination on the exit boundary $U^\out$ is symmetric to that on $U^\iin$.
According to Theorem~\ref{thm: geometric type in block}, we can also compute the boundary lamination of the associated block $(P,X)$.
The boundary of $(P,X)$ is a torus which contains exactly two periodic orbits $\cO_1$ and $\cO_2$ of the vector field $X$ (corresponding to the two compact leaves of $U^\iin$).
The lamination $\cL_X$ on $\partial P$ is a filling lamination, with no compact leaves other than the periodic orbits $\cO_i$ and no Reeb component.
Therefore, there exists a diffeomorphism $\varphi: \pP \to \pP$ which maps the periodic orbit $\cO_1$ to the orbit $\cO_2$, the exit boundary $\Pout$ to the entrance boundary $\Pin$, and the lamination $\cL_X$ to a lamination \sqt{} to $\cL_X$.
This diffeomorphism induces a gluing map for the building block formed by two copies $(P_1, X_1)$ and $(P_2, X_2)$ of $(P,X)$.
Up to a strong isotopy of the triple $(Q,Y, \varphi) := (P_1 \cup P_2, X_1 \cup X_2, \varphi)$, the \vfs{} $X_1$ and $X_2$ induce an Anosov \vf{} $Z$ on a closed orientable manifold $\cM = P_1 \cup P_2/\varphi$ according to the Gluing Theorem~\ref{thm: gluing theorem}.
The maximal invariant set of the fake horseshoe is transitive, and the graph $G(Q,Y,\varphi)$ is strongly connected, therefore the Anosov \vf{} is transitive (Proposition~\ref{prop: transitivity criterion}).

This construction is already known.
Indeed, the block given by Theorem~\ref{thm: geometric type in block} for the fake horseshoe is orbit equivalent to a block of geodesic flow on the compactified \emph{modular orbifold}.
The modular orbifold is the quotient $S_\Mod := \H^2 / \mathrm{PSL}_2(\Z)$.
It is a hyperbolic orbifold, which is a sphere with two singular points of order 2 and 3 and a \emph{cusp}.    
It is a known result (\cite{milnorHyperbolicGeometryFirst1983}) that the unitary tangent bundle $T^1 S_\Mod \simeq \mathrm{PSL}_2(\R) / \mathrm{PSL}_2(\Z)$ is homeomorphic to the complementary of the trefoil knot $K$ in $\S^3$, and the geodesic flow is the flow transverse to the fibration of $K$.
In \cite{ghysKnotsDynamics2007}, the author shows that one can deform the surface $S_\Mod$ so as to obtain a compact hyperbolic orbifold $\overline S_\Mod$ with geodesic boundary.
The geodesic flow $X^t$ induced on the tangent unitary fibered $P := T^1 \overline S_\Mod$ is an Axiom A flow whose basic piece admits a Markov partition of geometric type the fake horseshoe (see for example \cite{bonattiLorenzAttractorsModular2021}).
It is a building block with a single boundary component, two incoherently oriented periodic orbits in the boundary and a filling boundary lamination with no compact leaves other than the periodic orbits and no Reeb components.
It follows from the uniqueness in Theorem~\ref{thm: geometric type in block} that the block given by this proposition is isotopic to $(P,X)$.
In \cite{clayGraphManifoldsThat2021}, the authors show, in the manner of Handel-Thurston, the existence of a family of Anosov flows obtained by gluing back finite covers of two copies of $(P,X)$.
\end{example}

\section[Periodic orbit complement as JSJ piece]{Periodic orbit complement as JSJ piece of a transitive Anosov flow}
\label{sec: orbit complement}

The Jaco-Shalen-Johansson (JSJ) theorem (\cite{jacoNewDecompositionTheorem1978}, \cite{johannsonHomotopyEquivalences3manifolds1979}) allows us to split a compact, irreducible, orientable manifold of dimension three into components \linebreak[4]which admit a Seifert fibration, or which are atoroidal.
This splitting is done along embedded incompressible tori, and there is a minimal collection of such tori which is unique up to isotopy.
We recall that an embedded surface in a manifold is said to be \emph{incompressible} if the embedding induces an injective morphism of the fundamental groups.
We recall that a Seifert fibration (\cite{hempel3manifolds2004}) is a circle fibration over a 2-orbifold (possibly with boundary).
A manifold which admits such a structure (there may be several of them) is a \emph{Seifert fibered space}.
There is a finite number of \emph{singular fibers}, which are the fibers above the singular points of the orbifold.
The others are called \emph{regular}.
A manifold is said to be \emph{atoroidal} if every incompressible embedding of the torus is homotopic to a boundary component.

\begin{thm}[Jaco-Shalen, Johansson] \label{thm: JSJ decomposition}
Let $\cM$ be a closed, orientable, irreducible manifold of dimension three. 
Then there exists a finite collection $\cT = \{T_1, \dots, T_n\}$ of incompressible pairwise disjoint embedded tori,
such that the closure of the connected components $P_i$ of $\cM \ssm \bigcup_i T_i$ are either Seifert fibers or atoroidal.
The decomposition is unique up to isotopy for $T_i$ and $P_i$ if the collection $\cT$ is minimal.
\end{thm}

Any orientable 3-manifold $\cM$ which carries an Anosov vector field $X$ is irreducible (since its universal covering is $\R^3$, \cite{calegariFoliationsGeometry3manifolds2007}, \cite{barbotHyperbolicityGlobalHyperbolicity2005}), and thus admits a JSJ decomposition along incompressible tori.
Let $P$ be a Seifert piece of the JSJ decomposition of $\cM$. We say that
\begin{itemize}
    \item $P$ is a \emph{periodic Seifert piece} if there exists a Seifert fibration on $P$ such that the regular fiber is homotopic to a power of a periodic orbit of $X$.
    \item $P$ is a \emph{free Seifert piece} otherwise.
\end{itemize}

\begin{defi}[Pseudo-Anosov flow] \label{def: pseudo anosov flow}
A pseudo-Anosov flow $X^t$ on a closed 3-manifold $\cM$ is a flow which is locally modeled on a semi-branched covering of an Anosov flow.
\end{defi}

\begin{figure}[htb]
    \centering
     \vspace*{-2em}
    \includegraphics[height=0.27\textheight]{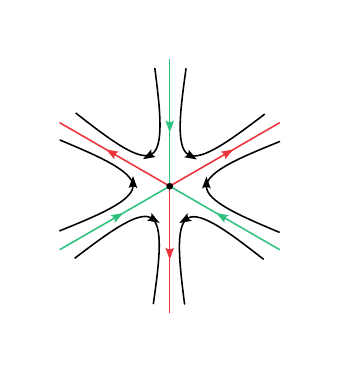}
     \vspace*{-2em}
    \caption{Transverse image of a $3$-prong singularity in a pseudo-Anosov flow}
    \label{fig: p prong}
\end{figure}

We follow the definition of \cite[Definition~6.41]{calegariFoliationsGeometry3manifolds2007}.
A pseudo-Anosov flow is thus a continuous flow, differentiable outside a finite collection of periodic orbits $\Gamma_* = \{ \gamma_1, \dots, \gamma_n \}$, called \emph{singular orbits}.
On the complementary of these orbits, the flow admits a hyperbolic splitting, in other words it preserves a splitting $T \cM = E^\ss \oplus \R.X \oplus E^\uu$ of the tangent bundle above $\cM \ssm \Gamma_*$, uniformly expands the vectors of $E^\uu$ and uniformly contracts the vectors of $E^\ss$ in the future.
Equivalently, we say that $X$ is a \emph{pseudo-Anosov \vf{}}.
The bundles $E^\ss \oplus \R X$ and $E^\uu \oplus \R X$ are uniquely integrable into $X^t$-invariant foliations $\cF^s$ and $\cF^u$ of dimension 2 called respectively the (weak) stable and (weak) unstable foliations, which admit $p_i$-prong singularities with $p_i \geq 3$ along the $\gamma_i$ orbits (Figure~\ref{fig: p prong}).
In the neighborhood of a singular orbit $\gamma_i$, the coordinates of the semi-branched covering are of class $\cC^1$ outside the singularity.

\begin{rmk}
There is another common definition of a \emph{topological pseudo-Anosov flow}, for example used in \cite[Definition~2.1]{barbotPseudoAnosovFlowsToroidal2013}, where one considers a flow only continuous on $\cM$.
We refer to \cite[Definition~5.8]{agolDynamicsVeeringTriangulations2022} to compare the two definitions.
The authors show in this paper (\cite[Theorem~5.10]{agolDynamicsVeeringTriangulations2022}) that a transitive topological pseudo-Anosov flow is orbit equivalent to a transitive pseudo-Anosov flow (in the smooth sense).
This is a generalization (more or less immediate) of M. Shannon's theorem \cite{shannonDehnSurgeriesSmooth2020} which states that a topologically transitive Anosov flow is orbit equivalent to a transitive Anosov flow.
\end{rmk}

The purpose of this chapter is to show Theorem~\ref{thmintro: orbit complement and JSJ piece}, which allows us to realize atoroidal JSJ pieces in a transitive Anosov flow as a complement of periodic orbits of transitive Anosov or pseudo-Anosov flow.

\begin{mainthm}[Theorem~\ref{thmintro: orbit complement and JSJ piece}] \label{thm: orbit complement and JSJ piece}
Let $\Gamma = \{ \gamma_1, \dots, \gamma_n \}$ be a finite collection of periodic orbits of a transitive pseudo-Anosov vector field $X$ on an orientable 3-manifold $\cM$.
Assume that all the singular orbits of $X$ are contained in $\Gamma$ and that the complementary $\cM \ssm \Gamma$ is atoroidal.
Then there exists an orientable 3-manifold $\cN$ carrying a transitive Anosov vector field $Y$ such that the JSJ decomposition of $\cN$ is made of two atoroidal pieces $P$ and $P'$, both homeomorphic to $\cM \ssm \Gamma$, and a periodic Seifert piece.
The restriction of $Y$ to $P$ and $P'$ is obtained from $X$ by a DA bifurcation on the orbits of $\Gamma$.
\end{mainthm}

\begin{rmk} \label{rmk: position of JSJ toruq}
The tori $T_i$ in the JSJ collection of $\cN$ are each isotopic to a \qt{} torus $T_i'$ in $\cN$, but the collection of $T_i'$ is not embedded: the tori intersect along periodic orbits.
This property will be checked during the proof, but we refer to \cite{barbotPseudoAnosovFlowsToroidal2013} for an analysis of the quasi-transverse position of JSJ tori in a dimension 3 Anosov flow.
An obstruction to embed the \say{modified} collection into quasi-transverse position occurs for tori adjacent to periodic Seifert pieces.
\end{rmk}

In fact, we will prove the following more general statement:
\begin{thm}\label{thm: orbit complement non-atoroidal}
Let $\Gamma$ be a finite collection of periodic orbits of a transitive pseudo-Anosov vector field $X$ on an orientable 3-manifold $\cM$, containing the set of singular periodic orbits of $X$.
Then there exists an orientable 3-manifold $\cN$ with a transitive Anosov \vf{} $Y$, and a collection $\cT$ of incompressible tori \linebreak[4]
embedded in $\cN$ and not homotopic,
which cut $\cN$ into three components $P_1, P_2, P_3$, where $P_1$ and $P_2$ are homeomorphic to $\cM \ssm \Gamma$ and $P_3$ admits a Seifert fibration with a fiber homotopic to a multiple of a periodic orbit.
The restriction of $Y$ to $P_1$ and $P_2$ is obtained from $X$ by a DA bifurcation on the orbits of $\Gamma$.
\end{thm}

Theorem~\ref{thm: orbit complement and JSJ piece} follows directly from Theorem~\ref{thm: orbit complement non-atoroidal} in the case $\cM \ssm \Gamma$ is atoroidal, since the collection $\cT$ is then the JSJ collection of tori, and $P_3$ is a periodic Seifert fibration.

\subsubsection*{Section summary}
The section is organized as follows.

\begin{itemize}[leftmargin=*]
    \item In Subsection~\ref{sec: orbit complement; subsec: double blow up} we give a general procedure for constructing a filled transitive block by \emph{Double blow up -- Excision} of a collection $\Gamma$ of periodic orbits of a pseudo-Anosov flow on a manifold $\cM$ containing the set of singular periodic orbits (Proposition~\ref{prop: double DA block})
    The resulting block is homeomorphic to the complementary $\cM \ssm \Gamma$ of the union of periodic orbits of the collection, and the boundary lamination is \emph{coherent elementary} on each boundary component.
    
    \item In Subsection~\ref{sec: orbit complement, subsec: add compact leaf}, we show a general result that allows us to \say{add} a compact (unmarked) leaf on an boundary lamination of a building block $(P,X)$ by gluing an \say{elementary} intermediate block, which is a tubular neighborhood of a periodic saddle orbit with negative multipliers.    
    Topologically, the block thus obtained is the union of $P$ and a \say{simple} Seifert fibered piece, adjoined by an incompressible torus.
    This result is summarized in Proposition~\ref{prop: add compact leaf}.
    
    \item In Subsection~\ref{sec: orbit complement; subsec: proof}, we show Theorem~\ref{thm: orbit complement non-atoroidal} using the previous results.
    We construct a block $(P,X)$ by \emph{Double Blow-up -- Excision} along the collection of orbits given by Proposition~\ref{prop: double DA block}, and add compact leaves on each connected component of the complement of the periodic orbits of the boundary with a well-chosen dynamical orientation using Proposition~\ref{prop: add compact leaf}.
    The resulting block $(P', X')$ is filled, transitive, and glues back to a copy of $(P,X)$ via a \sqt{} gluing map $\varphi$.
    We use the Gluing Theorem and the transitivity criterion, and the Anosov flow obtained will satisfy Theorem~\ref{thm: orbit complement non-atoroidal}.
    
    \item In Subsection~\ref{sec: orbit complement; subsec: knots}, we give as a corollary a criterion for the complementary of a \emph{hyperbolic knot} in $\S^3$ to be realized as an atoroidal JSJ piece of a transitive Anosov flow.
    We describe a construction by \emph{Hopf plumbing} which gives a family of knots satisfying this criterion.

\end{itemize}

\subsection{Building block by \emph{Double blow up -- Excision}}
\label{sec: orbit complement; subsec: double blow up}

In this preliminary section, we present the general result of creating \bb{} by \emph{Double blow-up -- Excision} of periodic orbits.
We refer to Subsection~\ref{sec: prescribed boundary lam; subsec: combinatorial type} for the definition of a \qms{} combinatorial type (Definitions~\ref{def: abstract combi type} and~\ref{def: qms combi type}).

\begin{defi}[Elementary, Coherent, Alternating \qms{} prefoliation] \label{def: elementary coherent alternating}
A \qms{} prefoliation on the torus is said to be
\begin{itemize}[--]
    \item \emph{elementary} if it admits no compact leaves other than marked leaves,
    \item \emph{coherent} if the oriented marked leaves are all freely homotopic,
    \item \emph{alternating} if any oriented marked leaf is freely homotopic to the successive marked leaf for a cyclic order.
\end{itemize}
\end{defi}

\begin{prop} \label{prop: double DA block}
Let $\Gamma = \{\gamma_1, \dots, \gamma_n\}$ be a collection of periodic orbits of a pseudo-Anosov vector field $X$ on an orientable closed manifold $\cM$ of dimension 3, which contains all (possible) singular orbits of $X$.
There exists a \bb{} $(P, Y)$ such that
\begin{enumerate}
    \item \emph{(manifold)}
    \label{prop: double DA block; it: block}
    $P = \cM \ssm \cV$, where $\cV$ is a union of disjoint tubular neighborhoods of the orbits $\gamma_i$.
    \item \emph{(boundary lamination)}
    \label{prop: double DA block; it: lamination}
    The boundary lamination $\cL$ is filling elementary coherent, each connected component $T_i$ of $\partial P$ contains $2 p_i$ periodic orbits, where $p_i$ is the number of stable separatrices of the orbit $\gamma_i$.
    \item \emph{(transitivity)}
    \label{prop: double DA block; it: transitivity}
    If $X$ is transitive, then $(P,Y)$ is a transitive block.
    \end{enumerate}
We say that $(P,Y)$ is a \bb{} obtained by \emph{Double blow-up -- Excision} on the collection $\Gamma$ of periodic orbits of the pseudo-Anosov \vf{} $X$.
\end{prop}
In particular, $P$ is homeomorphic to the complementary $\cM \ssm \Gamma$.

\begin{proof}[Idea of the proof]
We explain this construction for a single periodic orbit $\gamma$ of an Anosov flow, the generalization for a finite collection $\Gamma$ being similar.
Let us assume for the sake of argument that the multipliers of $\gamma$ are negative.
We make an attracting DA bifurcation on the periodic orbit $\gamma$.
We refer to \cite[Subsection~2.2.2]{ghristKnotsLinks3Dimensional1997} and \cite[Section~8]{beguinBuildingAnosovFlows2017} for a detailed description of this bifurcation, and to Figure~\ref{fig: attracting DA}.
Let us also mention \cite[Section~3.3]{barthelmeAnomalousAnosovFlows2021}, where the authors give an explicit construction of an DA bifurcation on a codimension 2 geodesic flow, but which can be adapted in our case.

This bifurcation consists in \say{opening} the stable manifold of $\gamma$, and creates an attracting periodic orbit $\gamma_+$ in place of $\gamma$ and a saddle periodic orbit $\gamma'$ with positive multipliers in place of the boundary of the local stable manifold of $\gamma$.
More precisely, there exists a small perturbation $X'$ of $X$ on $\cM$ such that $X'$ coincides with $X$ outside a small tubular neighborhood $\cV_0$ of the $\gamma$ orbit, and such that the vector field $X'$ on $\cV_0$ admits an attracting periodic orbit $\gamma_+$ with negative multipliers, and a saddle periodic orbit $\gamma$ with negative multipliers homotopic to the double of $\gamma_+$.
We then make a repelling Anosov derivative bifurcation on the periodic orbit $\gamma'$, which consists in \say{opening} the unstable manifold of $\gamma'$.
More precisely, there exists a small perturbation $X''$ of $X'$ on $\cM$ such that $X'$ coincides with $X$ outside a small tubular neighborhood $\cV'_0$ of the orbit $\gamma'$, and such that the vector field $X''$ on $\cV'_0$ admits a repelling periodic orbit $\gamma'_-$ of positive multipliers instead of $\gamma'$ and two saddle periodic orbits $\gamma'_1$ and $\gamma'_2$, instead of the boundary of the local unstable manifold of $\gamma'$, freely homotopic to $\gamma'_-$ and with positive multipliers, and the \vf{} has not been modified in the neighborhood of $\gamma^+$ (Figure~\ref{fig: double da}).
There exists a torus $T$ which is the union of two annuli $A_+$ and $A_-$ such that:
\begin{itemize}
    \item the boundary of $A_+$ is the union of the periodic orbits $\gamma'_1$ and $\gamma'_2$ and the interior of $A_+$ is contained in the basin of attraction $\cW^s(\gamma_+)$ of $\gamma_+$, and is transverse to the \vf{} $X''$;
    \item the boundary of $A_-$ is the union of the periodic orbits $\gamma'_1$ and $\gamma'_2$ and its interior is contained in the basin of repulsion $\cW^u(\gamma'_-)$ of $\gamma'_-$, and is transverse to the \vf{} $X''$.
\end{itemize}
The torus $T$ is drawn in Figure~\ref{fig: double da}.
It follows that $T$ is a torus quasi-transverse to the \vf{} $X''$ which contains the two saddle periodic orbits $\gamma'_i$.
It is the boundary of a solid torus $\cV$ contained in $\cV_0$ and which contains in its interior the attracting periodic orbit $\gamma_+$ and the repelling periodic orbit $\gamma'_-$.
Let $P:= \cM \ssm \cV$, and denote $Y$ the restriction of $X''$ to $P$.

\begin{figure}[htb]
    \centering
     \vspace*{-1em}
     \hspace*{-2em}
    \includegraphics[height=0.22\textheight]{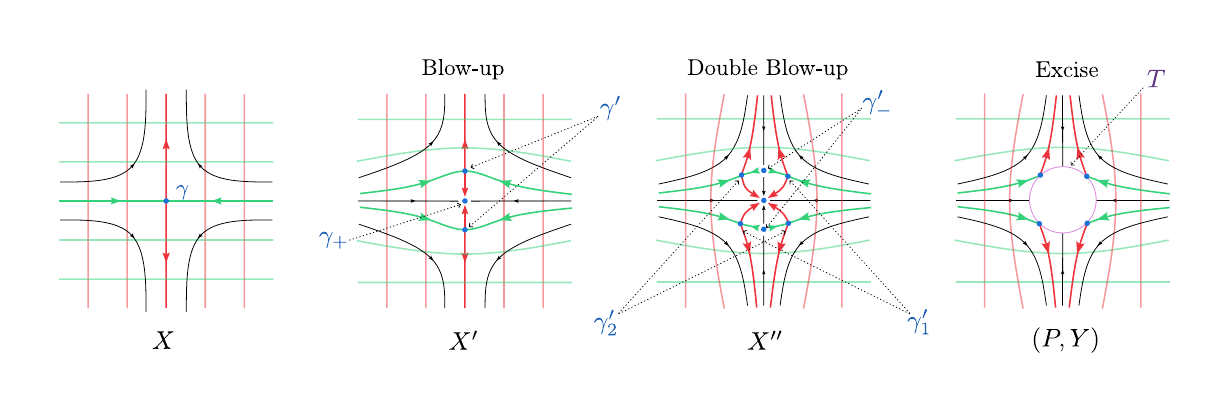}
     \vspace*{-2em}
    \centering \caption{Double blow-up and Excision on an orbit $\gamma$ with negative multipliers}
    \label{fig: double da}
\end{figure}

Then $(P,Y)$ is a building block satisfying Proposition~\ref{prop: double DA block}.
We refer to \cite[Proposition~3.8]{barthelmeAnomalousAnosovFlows2021} for a detailed proof of the hyperbolicity of the maximal invariant set which we can easily adapt to our case.
For the rest, we refer to the proof of Lemma~\ref{lem: attracting DA block}.

If $\gamma$ has positive multipliers, the first attracting DA bifurcation produces two saddle periodic orbits $\gamma_1$ and $\gamma_2$ and one attracting periodic orbit $\gamma_+$, all with positive multipliers.
We then perform a repelling DA bifurcation on $\gamma_1$ and on $\gamma_2$.
This creates two repelling periodic orbits $\gamma_{1,-}$ and $\gamma_{2,-}$ and four saddle periodic orbits $\gamma_{1,1}$ and $\gamma_{1,2}$ and $\gamma_{2,1}$ and $\gamma_{2,2}$.
There exists then a solid torus $\cV$ whose boundary is a torus $T$ quasi-transverse to the vector field $X''$, which contains the four saddle periodic orbits $\gamma_{i,j}$ and the interior of $\cV$ contains the attracting and repelling periodic orbits.
We show in the same way that the manifold $P= \cM \ssm \cV$ equipped with the \vf{} $Y$ equal to the restriction of $X''$ is a building block which satisfies the lemma.

In the case where $X$ is a pseudo-Anosov \vf{} and $\gamma$ is the unique singular orbit of $p$-prong type, the procedure is essentially the same and this operation will have the effect of \say{blowing-up} the $p$-prong singularity.
Let us explain briefly.
Suppose that $\gamma$ admits $p$-stable separatrices and $p$-unstable separatrices.
We perform an attracting bifurcation \emph{Derived from Pseudo-Anosov (DpA)} on the orbit $\gamma$, in other words a bifurcation analogous to the DA bifurcation, which consists in opening each of the stable separatrices of $\gamma$.
More precisely, there exists a tubular neighborhood $\cV$ of $\gamma$ on which the flow of $X$ is orbit equivalent to a semi-branched covering of a saddle hyperbolic periodic orbit.
This neighborhood is partitioned into $p$ regions $\cV_i$, $i=1, \dots, p$, each bordered by two unstable separatrices of $\gamma$, and on each of these regions there are coordinates $(x,y,t): \cV_i \to \R_+ \times \R \times \R/\Z$, of regularity $\cC^1$ outside $x=y=0$, in which the flow of $X$ is the suspension of $(x,y) \mapsto (\lambda^\inv x, \lambda y)$, with $\lambda >1$, restricted to the half space.
We can then return to the attracting DA bifurcation of the previous case on each region $\cV_i$ and to the explicit descriptions of \cite[Section~3.3]{barthelmeAnomalousAnosovFlows2021} or \cite[Section~8]{beguinBuildingAnosovFlows2017} in linearizing coordinates restricted to the half-space $x\geq0$.
We obtain a new \vf{} $X_i'$ on each $\cV_i$, of class $\cC^1$ outside $x=y=0$.
The blow-up \vfs{} $X_i'$ on the region $\cV_i$ reconnect in a $\cC^1$-way along the unstable separatrices by a reflection, and we obtain a field $X'$ of class $\cC^1$ outside $\gamma$.
This \vf{} $X'$ coincides with $X$ outside a small tubular neighborhood of $\gamma$, admits an attracting periodic orbit $\gamma^+$ in place of $\gamma$ and $p$ periodic $\gamma_i$ orbits in place of each boundary component of the local unstable separatrices of $\gamma$, and which are now saddle hyperbolic periodic orbits in the \nbh{} of which the flow is of class $\cC^1$.
The hyperbolicity (of index $(1,1)$) for the orbits of the flow of $X'$ on the complementary of the basin of attraction of $\gamma^+$ can be shown by following the proof of \cite[Proposition~3.8]{barthelmeAnomalousAnosovFlows2021} by cone criterion.
We can reduce to coordinates $(x,y,t)$ on each region $\cV_i$ and adapt the proof for half-spaces $x \geq 0$.
We then perform a repelling DA bifurcation on each of these $p$ saddle hyperbolic orbits in small and disjoint neighborhoods, far from $\gamma^+$.
We can then return to the classical construction, and we obtain a \vf{} $X''$ which coincides with $X'$ outside these small neighborhoods, of class $\cC^1$ outside $\gamma^+$, which admits $p$ repelling periodic orbits $\gamma_i^-$ in place of $\gamma_i$, and $2p$ saddle hyperbolic periodic orbits in place of the two boundaries of the local stable separatrices of the $p$ orbits $\gamma_i$.
The saddle hyperbolicity for the orbits in the complementary of the union of the repulsion basins of $\gamma_i^-$ and the attraction basin of $\gamma_i^+$ is shown by following the proof \cite[Proposition~3.8]{barthelmeAnomalousAnosovFlows2021}, since we are far from the singularity.
There exists a torus $T$ quasi-transverse to the vector field $X''$ described in Figure~\ref{fig: double da prong}, which contains each of the $2p$ saddle hyperbolic periodic orbits, and the complementary of these orbits in $T$ are annuli contained in the basins of attraction of $\gamma^+$ and basins of repulsion of $\gamma_i^-$.

\begin{figure}[h]
    \centering
    \vspace*{-2em}
    \captionsetup{width=.84\linewidth}
    \includegraphics[height=0.28\textheight]{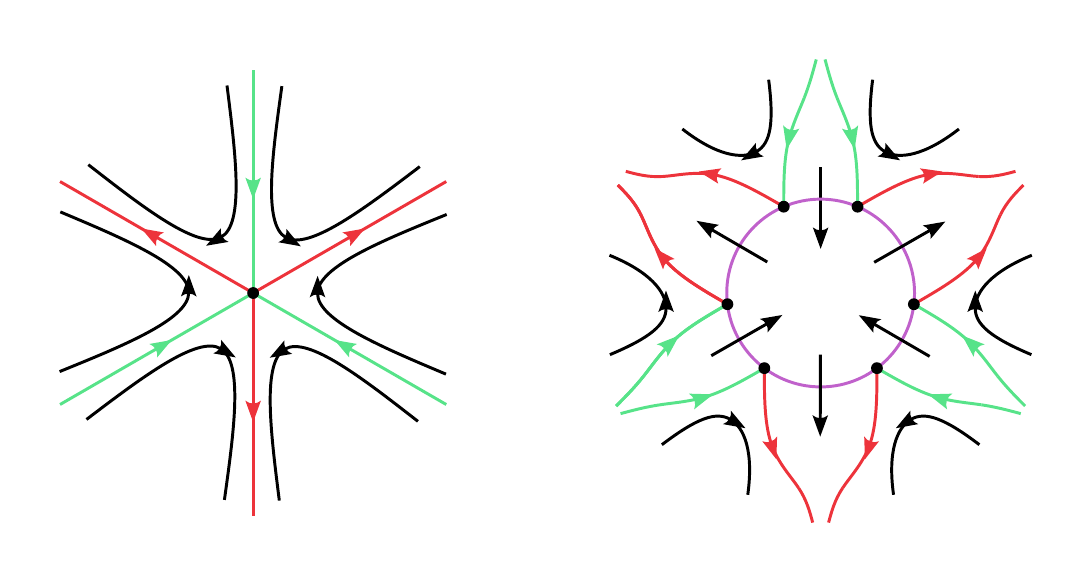}
     \vspace*{-1em}
    \caption{Double Blow-up and Excision on a $3$-prong singularity}
    \label{fig: double da prong}
\end{figure}

The torus $T$ is the boundary of a solid torus $\cV$ and contains in its interior the attracting periodic orbit $\gamma^+$ and the repelling periodic orbits $\gamma_i^-$.
Let $P:= \cM \ssm \cV$, and denote $Y$ the restriction of $X''$ to $P$.
We then show by analogous arguments that the pair $(P,Y)$ is a \bb{} which satisfies the proposition.
\end{proof}

\subsection{Addition of compact leaf}
\label{sec: orbit complement, subsec: add compact leaf}

A method is described for \say{adding} a compact leaf to the boundary lamination of a building block $(P,X)$. 
Topologically, this operation has the effect of gluing a solid torus along an annulus on a boundary component of $P$.
We will describe this operation using partial gluing maps of intermediate blocks (Definition~\ref{def: intermediate block}).
We use the results of Section~\ref{sec: geometric type} for the correspondence between building blocks and intermediate blocks.

Let $(P,X)$ be an orientable filled \bb{}, and $T$ be a boundary component of $\pP$ which contains a nonempty set of periodic orbits of $X$.
We denote $\cL_X$ the boundary lamination and $\Lambda_X$ the maximal invariant set.
Let $B$ be an annulus of $T$ bounded by two compact leaves $\gamma$ and $\gamma'$ of $\cL_X$ and containing no compact leaves of $\cL$ on its interior.
Equip $P$ with the orientation such that a geometric enumeration of the compact leaves of $\cL_X$ satisfy $\gamma_{n-1} = \gamma$ and $\gamma_0 = \gamma'$, i.e., $\gamma'$ is the first leaf, and $\gamma$ is to the left of $\gamma'$ for the dynamical orientation of $\gamma'$ and the orientation of $T$.
Let $\sigma$ be the associated combinatorial type of $\cL_X$.

\begin{prop} \label{prop: add compact leaf}
For any combinatorial type $\sigma'$ on $\Z/(n+1)\Z$ such that the restriction of $\sigma'$ to $\{ 0,\dots, n-1 \}$ coincides with $\sigma$, there exists $(Q, Y)$ a oriented filled \bb{}, such that
\begin{enumerate}
    \item \emph{(boundary)} \label{prop: add compact leaf; it: boundary}
    $\partial P \simeq \partial Q$.

    \item \emph{(Smale's graph)} \label{prop: add compact leaf; it: smale graph}
    The Smale's graph of $(Q,Y)$ is the Smale's graph of $(P,X)$ where we add a vertex $O$ and $n$ oriented edges $O \rightarrow \Lambda_i$ if $B \subset \Pin$, or $n$ oriented edges $\Lambda_i \rightarrow O$ if $B \subset \Pout$, where $O$ is a saddle periodic orbit, and $\Lambda_1, \dots, \Lambda_n$ are the basic pieces of $\Lambda_X$ whose invariant manifold intersect the annulus $B$.
    
    \item \emph{(combinatorial type)} \label{prop: add compact leaf; it: combi type}
    There exists $T'$ a \cc{} of $\partial Q$ such that:
    \begin{itemize}[--]
        \item the boundary lamination $\cL_Y$ on $T'$ has combinatorial type $\sigma'$, and the $(n+1)$-th compact leaf is the intersection of the invariant manifold of $O$ with $\partial Q$;
        \item $\cL_Y$ and $\cL_X$ extends as topologically equivalent foliations on $\partial Q \ssm T'$ and $\partial Q \ssm T$ respectively.
        
    \end{itemize}
    
    \item \emph{(topology)} \label{prop: add compact leaf; it: topology}
    there exists a torus $T''$ embedded in the interior of $Q$, isotopic to a torus \qt{} to $Y$, which cuts $Q$ into two pieces $P_0$ and $P_1$ such that $P_0 \simeq P$ and $P_1$ is a Seifert fibered space adjacent to $T'$ and $\partial Q$, whose regular fiber is freely homotopic to the periodic orbits of $Y$ contained in $T'$.
\end{enumerate}
\end{prop}

\begin{figure}[htb]
    \centering
    \vspace*{-1em}
    \includegraphics[height=0.28\textheight]{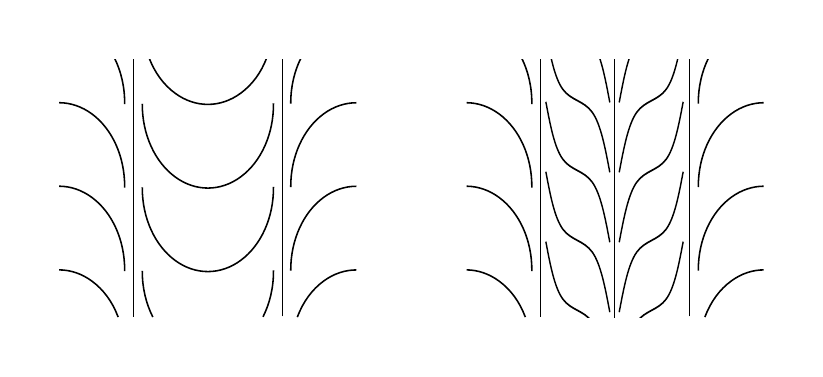}
    \vspace*{-1em}
    \caption{An addition of compact leaf}
    \label{fig: compact leaf addition}
\end{figure}

Item~\ref{prop: add compact leaf; it: combi type} means that the boundary lamination of $(P',X')$ is essentially the same as the boundary lamination of $(P,X)$, where a compact leaf has been added to the middle of the annulus $B$ bounded by $\gamma$ and $\gamma'$.

\begin{rmk}
If $\sigma'$ is a combinatorial type which coincides with a combinatorial type $\sigma$ on the complementary of an integer in $\Z/(n+1)\Z$, then this integer does not correspond to a marked leaf, by parity of the marked leaves.
This process thus allows us to add compact unmarked leaves to a boundary lamination.
\end{rmk}

A first step is to build an \say{elementary} intermediate block (Definition~\ref{def: intermediate block}), containing a single periodic orbit.

\begin{lemdef}
\label{lem: elementary block}
There exists an \imb{} $(V_e,X_e)$, such that:
\begin{enumerate}
     \item \emph{(boundary)} \label{lem: elem block; it: boundary}
    $\partial V_e = V^\iin_e \cup V^\tan_e \cup V_e^\out$, where the entrance boundary $V_e^\iin$ is a single annulus, the exit boundary $V^\out$ is a single annulus, and the tangent boundary $V_e^\tan$ is the union of two annuli tangent to the vector field $X_e$.
    \item \emph{(maximal invariant set)} \label{lem: elem block; it: max inv}
    the maximal invariant set is a single saddle periodic orbit $O$ with negative multipliers,
    \item \emph{(boundary lamination)} \label{lem: elem block; it: lamination}
    the boundary lamination $\cL_{X_e}$ is reduced to two compact leaves, one $\gamma^s$ on $V^\iin_e$ and one $\gamma^u$ on $V^\out_e$, corresponding to the intersections of the stable and unstable separatrix of $O$ respectively.
    \item \emph{(topology)} \label{it: block elem_ topology}
    $V_e$ is homeomorphic to a solid torus $\D^2 \times \S^1$, it is a Seifert fibered space whose regular fiber is freely homotopic to the trace of the stable manifold $\gamma^s$, and has a unique singular fiber homotopic to the periodic orbit $O$.
\end{enumerate}
Such a block is unique up to orbit equivalence. 
It is called the \emph{elementary intermediate block}.
\end{lemdef}

\begin{figure}[htb]
    \centering
    \vspace*{-1em}
    \includegraphics[height=0.35\textheight]{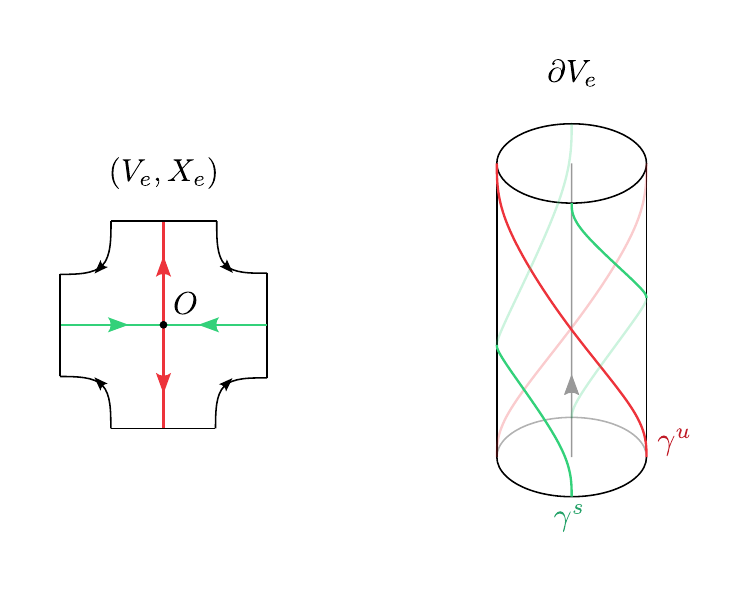}
    \vspace*{-2em}
    \caption{Elementary block $(V_e, X_e)$ and boundary lamination}
    \label{fig: elementary block}
\end{figure}

\begin{proof}
Consider the suspension of the matrix $A \colon (x,y) \mapsto (-\frac{1}{2}x, 2y)$ on $\R^2 \times \R/\Z$, and $O$ the periodic orbit induced by $(0,0)$.
Let $V_e$ be a tubular neighborhood of $O$ whose boundary consists of two tangent annuli and two transverse annuli, whose trace in the transverse plane $\R^2 \times \{0\}$ is described by Figure~\ref{fig: elementary block}.
The block $(V_e, X_e)$ where $X_e$ is the restriction of the suspension \vf{} satisfies the lemma.
If $(V,X)$ is another intermediate block which satisfies Lemma~\ref{lem: elementary block}, up to push the boundary of $V$ along the vector field we can assume that it is contained in a linearizing neighborhood of the unique periodic orbit $O$.
Saddle hyperbolic periodic orbits with negative multipliers are all orbit equivalent on a linearizing neighborhood.
We show with arguments similar to the proof of~\ref{lem: bb to imb} that neighborhoods described by Figure~\ref{fig: elementary block} are orbit equivalent.
\end{proof}

We will need the following lemma to describe partial gluing of intermediate blocks.
Let $(\check P, \check X)$ and $(\check Q, \check Y)$ be two \imbs{}.
Let $\partial \check P = \check P^\iin \cup \check P^\tan \cup \check P^\out$ denote the decomposition of the boundary into entrance boundary, tangent boundary and exit boundary, and $\check \cL_X$ the boundary lamination of $(\check P, \check X)$.
We denote $\check f_X : \check P^\iin \to \check P^\out$ the crossing map of $\check X$.
Similarly, $\partial \check Q = \check Q^\iin \cup \check Q^\tan \cup \check Q^\out$ the boundary decomposition, $\check \cL_Y$ the boundary lamination of $(\check Q, \check Y)$.
Let us denote $\check f_Y : \check Q^\iin \to \check Q^\out$ the crossing map of $\check Y$.

\begin{lem} \label{lem: gluing imb}
Let $C^\out$ be a union of \ccs{} of $\check P^\out$ and $C^\iin$ be a union of \ccs{} of $\check Q^\iin$, let $\varphi: C^\out \to C^\iin$ be a diffeomorphism which reverses the orientation and such that $\varphi_* \check \cL_X$ is transverse to $\check \cL_Y$ on $C^\iin$.
Then the vector fields $\check X$ and $\check Y$ induce a vector field $\check Z$ of class $\cC^1$ on the manifold
$\check N:= (\check P \cup \check Q)\, /\varphi$, and we have
\begin{enumerate}
    \item \emph{(intermediate block)} \label{lem: gluing imb; it: imb}
    $(\check N, \check Z)$ is an \imb{}.
    \item \emph{(boundary)} \label{lem: gluing imb; it: boundary}
    $\partial \check N = (\partial \check P \ssm C^\out) \cup_\varphi (\partial \check Q \ssm C^\iin)$, and
        \begin{itemize}[--]
    \item the entrance boundary is
    $\check N^\iin = \partial \check P^\iin \sqcup (\partial \check Q^\iin \ssm C^\iin)$
    \item the exit boundary is
    $\check N^\out = (\partial \check P^\out \ssm C^\out) \sqcup \partial \check Q^\out$
    \item the tangent boundary is 
    $\check N^\tan = \check P^\tan \cup_\varphi \check Q^\tan$.
    \end{itemize}
    \item \emph{(lamination)} \label{lem: gluing imb; it: lam}
    If $\check \cL_Z$ is the boundary lamination of $\check N$ then:
    $$\check \cL_Z = \check \cL_X \cup \check \cL_Y \cup (\check f_Y)_* (\varphi_* \cL_X \ssm \cL_Y) \cup (\check f_X^\inv)_* (\varphi^\inv_* \cL_Y \ssm \cL_X).\footnote{We omit by use of notation to write the domains where the maps are well defined.}$$
\end{enumerate}
\end{lem}

\begin{figure}[htb]
    \centering
    \vspace*{-2em}
    \includegraphics[height=0.27\textheight]{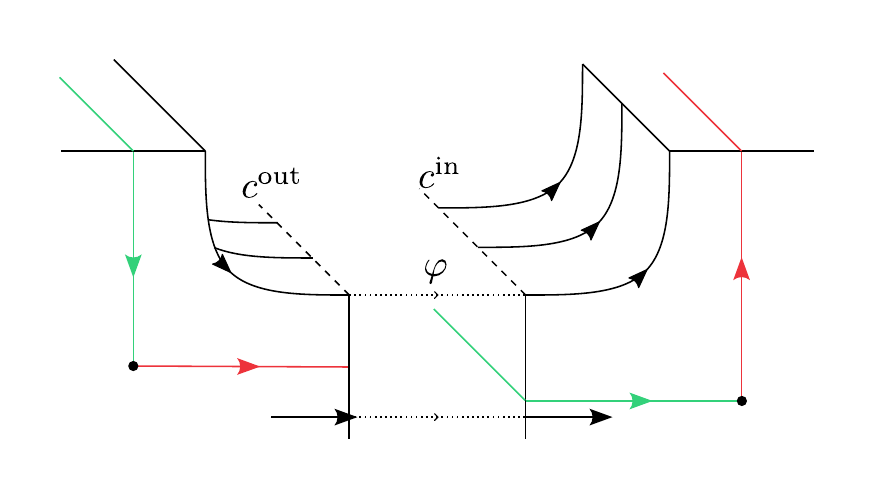}
    \vspace*{-1em}
    \caption{Gluing \imbs{}}
    \label{fig: gluing imb}
\end{figure}

\begin{proof}
For an \imb{} $(\check P, \check X)$, if $c^\iin$ is a curve which borders a connected component of $\check P^\iin$ and a connected component of $\check P^\tan$, then there exist local coordinates $(x, y, \theta)$ in the neighborhood of any point of $c^\iin$ such that the vector field $\check X$ is the vertical vector field $\partial_y$, the tangent boundary $\check Q^\tan$ is the vertical plane $\{x=0\}$ and the entrance boundary $\check P^\iin$ is the horizontal plane $\{y=0\}$.
The same is true in the \nbh{} of the points of a curve $c^\out$ which borders a connected component of $\check Q^\iin$ and $\check Q^\tan$ (Figure~\ref{fig: gluing imb}).
It follows that the quotient space
$\check N:= \check P \cup \check Q /\varphi$ is a manifold with boundary, smooth outside a finite number of simple closed curves contained in the boundary, and a vector field $\check Z$ of class $\cC^1$ induced by the vector fields $\check X$ and $\check Y$ and the pair $(\check N, \check Z)$ satisfies Item~\ref{def: imb; it: edges} and~\ref{def: imb; it: boundary} of Definition~\ref{def: intermediate block} of an intermediate block.
Item~\ref{def: imb; it: boundary} of Lemma~\ref{lem: gluing imb} is deduced from the construction.
Let us then remark that any orbit of the flow of $\check Z$ which intersects the projection of $C^\iin$ in $\check N$ never passes through it again, because each orbit of the flow of $\check X$ intersects the exit boundary at most once, each orbit of the flow of $\check Y$ intersects the entrance boundary at most once, and the blocks $(\check P, \check X)$ and $(\check Q, \check Y)$ are disjoint. 
It follows that we can see $\varphi$ as a \emph{partial gluing map without cycle} (see Subsection~\ref{sec: prescribed boundary lam; subsec: gluing without cycle}).
We can then adapt the proof of Proposition~\ref{prop: block obtained by partial gluing without cycle} for intermediate blocks: if $\varphi$ maps the boundary lamination $\check \cL$ of $(\check P, \check X)$ on $C^\iin$ on a lamination transverse to the boundary lamination $\check \cL_Y$ of $(\check Q, \check Y)$ on $C^\iin$, then the maximal invariant set $\check \Lambda_Y$ of the pair $(\check N, \check Z)$ is a hyperbolic set and the pair is a \imb{}.
Item~\ref{lem: gluing imb; it: imb} is verified.
Item~\ref{lem: gluing imb; it: lam} is deduced with considerations analogous to the proof of Lemma~\ref{lem: block obtained by partial gluing without cycle} concerning partial gluing without cycle of \bbs{}.
\end{proof}

\begin{proof}[Proof of Proposition~\ref{prop: add compact leaf}]
We recall that $(P,X)$ is a filled orientable block, $T$ a boundary component of $\pP$ which contains at least one periodic orbit of $X$, and we take the notations of Proposition~\ref{prop: add compact leaf}.
Recall that we have fixed an annulus $B$ of $T$ bounded by two compact leaves $\gamma$ and $\gamma'$ of the boundary lamination $\cL_X$ and containing no compact leaves on its interior.
We chose the first leaf $\gamma_0 = \gamma'$ and an orientation on $P$ such that $\gamma$ is to the right of $\gamma_0$ for its dynamical orientation.
The corresponding combinatorial type of lamination $\cL_X$ on $T$ is denoted $\sigma$ on $\Z/n\Z$.
The vector field $X$ is transverse to $B$ and say exiting to fix the ideas.
Let $(\check P, \check X)$ be the \imb{} associated to $(P,X)$ by Lemma~\ref{lem: bb to imb}.
Denote $\check P^\iin$ the entrance boundary, $\check P^\out$ the exit boundary, and $\check \cL_X$ the boundary lamination (Definition~\ref{def: intermediate block} and~\ref{def: in out tangent boundary of imb}).
Denote $\cA_*$ the collection of connected components of $\partial \check P \ssm \check \cL_X$ containing the tangent boundary $\check P^\tan$ (Definition~\ref{def: intermediate block}, Item~\ref{def: imb; it: annuli in boundary}).

By Item~\ref{lem: bb to imb; it: lamination} of Lemma~\ref{lem: bb to imb}, the closed annulus $\adh(B)$ is isotopic to a closed annulus $\adh(\check B) \subset \check P^\out$ via an isotopy which maps the restricted lamination $\cL_X$ to the lamination $\check \cL_X$.
Let $C^\out$ be the connected component of $\check P^\out$ which contains $\check B$.
Let $(V_e, X_e)$ be the elementary \imb{} given by Lemma~\ref{lem: elementary block}, and let $V^\iin_e$ be the entrance boundary of $V_e$.
It is a annulus and the boundary lamination $\cL_e$ on $V^\iin_e$ is reduced to a single compact leaf $\gamma^s \subset \intr(V^\iin_e)$ parallel to the boundary of the annulus.
Let $\varphi : C^\out \to V^\iin_e$ be a diffeomorphism such that the leaf $\varphi_*^\inv(\gamma^s)$ is inside $\check B \subset C^\out$, parallel (and disjoint) to the compact leaves of the lamination $\check \cL$ on $C^\out$ (Figure~\ref{fig: gluing to add compact leaf}).
Let $\check Q = (\check P \cup V_e)/\varphi$.
Then by Lemma~\ref{lem: gluing imb}, the vector fields $\check X$ and $X_e$ induce a vector field $\check Y$ of class $\cC^1$ on the manifold $\check Q$, such that $(\check Q, \check Y)$ is a \imb{} with orientation coinciding with the one chosen on $P$.

We want to compute the Smale's graph of $(\check Q, \check Y)$, from the graph $G(\check P \cup V_e, \check X \cup V_e, \varphi)$ given by Definition~\ref{def: graph g(P,X,varphi)} (the analogous for intermediate block).
Let $\Lambda_1, \dots$, $\Lambda_n$ be the basic pieces of $\Lambda_X$ such that $\cW_X^s(\Lambda_i) \cap B \neq \emptyset$.
By Item~\ref{lem: bb to imb, it: minc}, the germ of $X$ on $\Lambda_X$ is the same than the germ of $\check X$ on the maximal invariant set $\check \Lambda_X$ of $(\check P, \check X)$, hence the Smale's graph of $\Lambda_X$ is isomorphic to the Smale's graph of $\check \Lambda_X$.
Let $\check \Lambda_1, \dots, \check \Lambda_n$ be the basic pieces of $\check \Lambda_X$ corresponding to $\Lambda_1, \dots, \Lambda_n \subset \Lambda_X$, in other words which satisfy $\check \cW_X^s(\check \Lambda_i) \cap \check B \neq \emptyset$.
By definition of $\varphi$, the curve $\varphi^\inv_* \gamma^s$ intersects every leaf of $\check \cW^u_X$ intersecting $\check B$, and no one else.
It follows that $\cW^s(O) \cap \varphi_* (\check \cW^u_X(\check \Lambda_i)) \neq \emptyset$.
We deduce:
\begin{claim} \label{claim: smale graph after gluing imb elementary}
    The Smale's graph of $(\check Q, \check Y)$ is isomorphic to the Smale's graph of $(P,X)$ where we add a vertex $O$ corresponding to a saddle periodic orbit, and $n$ oriented edges $\Lambda_i \rightarrow O$, $i=1, \dots, n$.
\end{claim}

We want now to compute the boundary lamination of $(\check Q, \check Y)$.
The exit boundary $\check Q^\out$ is the disjoint union $(\check P^\out \ssm C^\out) \sqcup V^\out_e$ and the entrance boundary $\check Q^\iin$ is equal to $\check P^\iin$ (Item~\ref{lem: gluing imb; it: boundary} Lemma~\ref{lem: gluing imb}).
Let $\cA'_*$ be the collection of connected components of $\partial \check Q \ssm \check \cL_Y$ containing the tangent boundary $\check Q^\tan$.

\begin{claim} \label{claim: description gluing elementeray block}
\mb{}
\begin{enumerate}
    \item $\check \cL_Y$ is filling on $(\partial \check Q \ssm \cA'_*) \ssm V^\out_e$ and $\check \cL_Y$ and $\check \cL_X$ complete into topologically equivalent foliations on $(\partial \check Q \ssm \cA'_*) \ssm V^\out_e$ and $(\partial \check P \ssm \cA_*) \ssm C^\out$ respectively.
    \item $\check \cL_Y$ is filling on $V^\out_e \ssm \cA'_*$, and there exists a annulus $\check B' \subset V^\out_e$ bounded by two compact leaves of $\check \cL_Y$, such that $\check \cL_X$ on $(C^\out \ssm \cA_*) \ssm \check B$ and $\check \cL_Y$ on $(V^\out_e \ssm \cA'_*) \ssm \check B'$ are topologically equivalent, and $\check \cL_Y$ contains a unique compact leaf on $\check B'$, which is the unique intersection of the invariant manifolds of $O$ with $\partial \check Q$.
    \end{enumerate}
\end{claim}

\begin{figure}[htb]
    \centering
    \vspace*{-1em}
    \includegraphics[height=0.325\textheight]{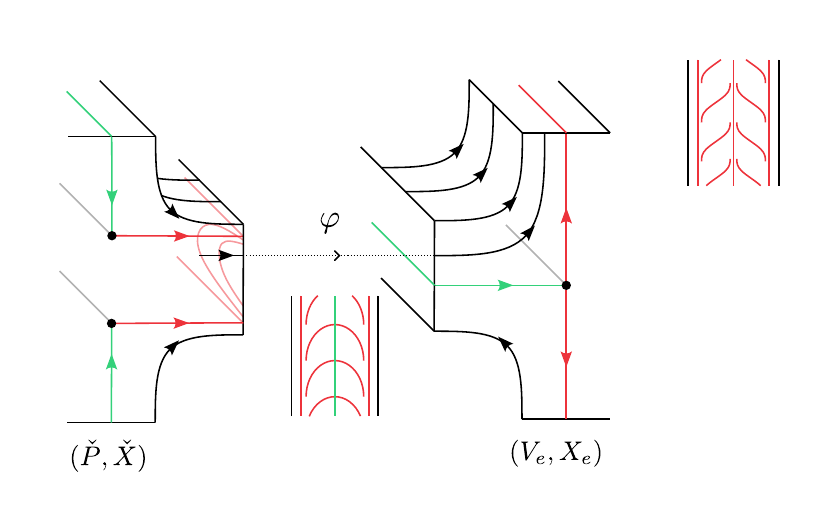}
    \vspace*{-1em}
    \caption{Gluing an elementary block $(V_e, X_e)$ along $V^\iin_e$ together with the intermediate block $(\check P, \check X)$ along $C^\out$, and induced lamination on $V^\out_e$}
    \label{fig: gluing to add compact leaf}
\end{figure}

\begin{proof}
It is enough to use Lemma~\ref{lem: gluing imb}, Item~\ref{lem: gluing imb; it: boundary} and~\ref{lem: gluing imb; it: lam} and to compute the lamination on each set.
Let $\check f_X: \check P^\iin \to \check P^\out$ denote the crossing map of the flow of $\check X$, and $f_e: V^\iin_e \to V^\out_e$ the crossing map of the flow of $X_e$.
Let $\gamma^u$ and $\gamma^s$ be the two leaves of the boundary lamination of $(V_e, X_e)$ on $V^\out_e$ and $V^\iin_e$ respectively.
These are the intersection of $\cW^u(O)$ and $\cW^s(0)$ with $\partial V_e$, where $O$ is a saddle periodic orbit (Item~\ref{lem: elem block; it: lamination} of Lemma~\ref{lem: elementary block}).
The lamination $\check \cL_Y$ on $\check Q^\out \ssm V^\out_e = \check P^\out$ is equal to the lamination $\check \cL_X$ on $\check P^\out$.
The lamination $\check \cL_Y$ on $\check Q^\iin \ssm V^\out_e = \check P^\iin$ is equal to the union
$\check \cL_X \cup (\check f_X)^\inv_* (\gamma^s \ssm \cL_X)$.
Since the lamination $\check \cL_X$ is filling on $\check P \ssm \cA_*$, then the lamination $\check \cL_Y$ on $\check Q^\iin \ssm \cA'_*$ is a lamination which contains a filling lamination, and the first item follows.
The lamination $\check \cL_Y$ on $C'^\out = V^\out_e$ coincides with the union \mb{$\gamma^u \cup (f_e)_* (\varphi_* \check \cL_X \ssm \gamma^s)$}.
As the lamination $\check \cL_X$ is filling, so is the image by $f_e \circ \varphi$, so the lamination $\check \cL_Y$ is filling on $V^\out_e$.
The lamination $\check \cL_Y$ on the complementary of $\gamma^u$ is topologically equivalent to the lamination $\check \cL$ on the complementary of $\varphi^\inv_* \gamma^s$.
Recall that we chose $\varphi$ so that $\varphi^\inv_* \gamma^s$ is contained in the interior of the annulus $\check B$ defined earlier, bounded by two compact leaves of $\check \cL_X$ and without compact leaves of $\check \cL_X$ on the interior (Figure~\ref{fig: gluing to add compact leaf}).
Denoting $\check B' = f_e(\check B)$, the last item follows.
\end{proof}

Finally, let $(Q,Y)$ be the \bb{} associated to the intermediate block $(\check Q, \check Y)$ by Lemma~\ref{lem: imb to bb}.
Note that the periodic orbits of $Y$ contained in $Q$ are the same as the periodic orbits of $X$ contained in $P$ (via the embedding of $(P,X)$ into $(Q,Y)$).
Indeed, it is enough to see on the associated intermediate blocks that the compact leaves which border the annuli $\cA'_*$ and the compact leaves which border the annuli $\cA_*$ in $\partial \check P$ come from the separatrices of the same periodic orbits by the gluing construction.
We recall that these separatrices correspond to a single periodic orbit in the intermediate block, and this orbit is a periodic orbit contained in the boundary of the associated building block (see Figure~\ref{fig: bb to imb} and Lemma~\ref{lem: tangent annuli imb}).
There is a bijective correspondence between the pairs of compact leaves bordering $\cA_*$ in the intermediate block, and the periodic orbits in the boundary of the building block (Remark~\ref{rmk: bb to imb (h) hyp}).

Item~\ref{prop: add compact leaf; it: boundary} of Proposition~\ref{prop: add compact leaf} is deduced from the fact that the surgery performed to go from $(P,X)$ to $(Q,Y)$ neither adds nor removes boundary components.
Topologically, a solid torus is glued together with $P$ along a annulus in the component $T$ of $\partial P$, resulting in a component $T'$ homeomorphic to $T$, and the other components are left intact.
Item~\ref{prop: add compact leaf; it: smale graph} directly follows from Claim~\ref{claim: smale graph after gluing imb elementary} and Item~\ref{lem: imb to bb; it: embedding}, Lemma~\ref{lem: imb to bb}.

Recall that Lemma~\ref{lem: imb to bb}, Item~\ref{lem: imb to bb; it: lamination} gives a \diff{} between $\check Q^\iin \ssm \cA'_*$ and $\overline Q^\iin$ which maps the lamination $\cL^\iin_X$ on the lamination $\check \cL^\iin_X$ and a \diff{} between $\check Q^\out \ssm \cA'_*$ and $\overline Q^\out$ which maps the lamination $\cL^\out_X$ on the lamination $\check \cL^\out_X$.
The laminations are topologically equivalent on these domains.
The annulus $\check B'$ is therefore isotopic to a annulus $B'$ in the exit boundary $Q^\out$, bounded by two compact leaves of the lamination $\cL_Y$ and containing a single compact leaf in its interior.
Let $T'$ be the connected component of $\pP$ which contains $B'$. 
Choose the first compact leaf of $\cL_Y$ on $T'$ to be the compact leaf which borders $B'$ and which corresponds to the leaf $\gamma'$ bordering $B$ (by the correspondence given by topological equivalences and isotopies to go from $B$ to $B'$).
Denote $\sigma_Y$ the associated combinatorial type.
Then according to Claim~\ref{claim: description gluing elementeray block} and the topological equivalence of laminations, $\sigma_Y$ is a combinatorial type on $\Z/(n+1)\Z$ which coincides with $\sigma$ on $\{ 0,\dots, n-1 \}$, and the $(n+1)$-th compact leaf is the unique intersection of $\cW^u(O)$ with $\partial Q$.
We don't know a priori if we have $\sigma_Y(n) = (*,\ua,*)$ or $\sigma_Y(n) = (*, \da, *)$, but up to replace the partial gluing map $\varphi : C^\out \to V^\iin_e$ in the proof by $\tau \circ \varphi$, where $\tau$ is a \diff{} of $V^\iin_e$ which fixes $\gamma^s$, reverses its orientation, and preserves the global orientation, then the value of $\sigma_Y(n)$ coincides with $\sigma'(n)$, where $\sigma'$ is the combinatorial type prescribed in Proposition~\ref{prop: add compact leaf}.
Item~\ref{prop: add compact leaf; it: combi type} of Proposition~\ref{prop: add compact leaf} follows.

Let us show Item~\ref{prop: add compact leaf; it: topology}.
By construction, $\check P$ and $V_e$ are embedded in $\check Q$, and the boundary $\partial \check P$ cuts $\check Q$ into a component homeomorphic to $\check P$, hence to $P$, and a component homeomorphic to $V_e$, hence to a Seifert fibered space.
Moreover, we know that the regular fiber of $V_e$ is homotopic to the leaf $\gamma^s$ of the boundary lamination of $V_e$, and $\gamma^s$ is connected to a curve (freely) homotopic to the compact leaves of the lamination $\check \cL_X$ on $\partial \check P$.
The block $P$ is embedded in $\check P$ and isotopic to $\check P$ (Lemma~\ref{lem: bb to imb}), and the compact leaves of the lamination of $\check \cL_X$ on $\partial \check P$ are homotopic to the periodic orbits of $\check \cL$ contained in the boundary of $T$ according to Item~\ref{lem: bb to imb; it: lamination}, Lemma~\ref{lem: bb to imb}.
Finally, $Q$ is isotopic to $\check Q$, it follows that we can embed $P$ into $Q$. Denote $T''$ the boundary of $P$ via this embedding.
Since $\check Q^\iin = \check P^\iin$ and the boundary laminations coincide along the compact leaves on the entrance boundary, it follows that the free homotopy class of the periodic orbits of $X$ contained in $T''$ coincides with the free homotopy class of the periodic orbits of $Y$ contained in $T' \subset \partial Q$.
By a small isotopy of $T''$ to make it disjoint from $T'$, it cuts in $Q$ a component homeomorphic to $P$ and a Seifert fibered space whose fiber is homotopic to the periodic orbits of $Y$ in $T'$. 
\end{proof}

\begin{rmk} \label{rmk: isotopic qt torus}
From the proof we can see that the torus $T''$ given by Item~\ref{prop: add compact leaf; it: topology} of Proposition~\ref{prop: add compact leaf} is isotopic to a torus quasi-transverse to $Y$ which contains the same periodic orbits as $T'$ (but is not isotopic to $T'$).
\end{rmk}

\subsection{Proof of Theorem~\ref{thmintro: orbit complement and JSJ piece}}
\label{sec: orbit complement; subsec: proof}

\begin{proof}[Proof of Theorem~\ref{thm: orbit complement non-atoroidal}]
We show the proposition for a collection $\Gamma$ reduced to a single orbit, the generalization for a collection being immediate.
Let $\gamma$ be a periodic orbit of a transitive Anosov vector field or the unique singularity of a transitive pseudo-Anosov flow on an orientable 3-manifold $\cM$.
Let $(P,Y)$ be the building block obtained by \emph{Double blow-up -- Excision} on the orbit $\gamma$ by Proposition~\ref{prop: double DA block}.
It is a connected transitive filled block with a single boundary component, and the boundary lamination is coherent and elementary, i.e., it admits no compact leaf other than the periodic orbits, and the oriented periodic orbits are freely homotopic in $\partial P$.
Let $2p$ be the number of periodic orbits contained in $\partial P$.
Then we can iterate Proposition~\ref{prop: add compact leaf} of compact leaf addition for every annulus of $\partial P$ bounded by two periodic orbits and without compact leaves inside, and we obtain a transitive filled saddle orientable connected block $(P', Y')$ such that if $\cL'$ denotes the boundary lamination, then $\cL'$ contains a collection $\cO_*'$ of $2p$ coherently oriented periodic orbits, and on each connected component of $\pP' \ssm \cO_*$ it admits a single compact leaf.
We further require that the dynamical orientation of each of these compact leaves is the reversed orientation of the periodic orbits.

\begin{claim}
There exists a \diff{} $\varphi : \pP \to \pP'$ which maps the oriented periodic orbits of $Y$ in $\pP$ to the oriented periodic orbits of $Y'$ in $\pP'$ and the exit boundary $\Pout$ to the entrance boundary $P'^\iin$, and such that $\varphi_* \cL$ is strongly quasi-transverse to $\cL'$ on $\pP'$.
\end{claim}

\begin{figure}[htb]
    \centering
    \vspace*{-2em}
    \includegraphics[height=0.28\textheight]{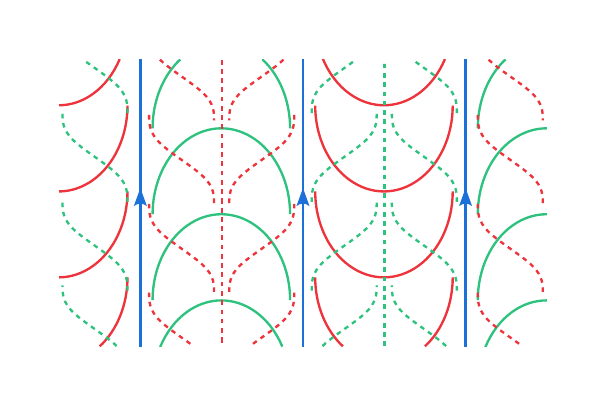}
   \vspace*{-1em}
    \caption{Strongly quasi-transverse gluing of the laminations $\cL$ (full line) and $\cL'$ (dashed line)}
    \label{fig: sqt gluing on elementary coherent}
\end{figure}
We refer to Figure~\ref{fig: sqt gluing on elementary coherent}.
The \diff{} $\varphi$ induces a strongly quasi-trans\-verse gluing map of the filled saddle block $(P \cup P', Y \cup Y')$.
Up to make a strong isotopy of the triple $(P \cup P', Y \cup Y', \varphi)$, we can apply the Gluing Theorem~\ref{thm: gluing theorem} and the vector fields $Y$ and $Y'$ induce on the closed orientable manifold $ \cN := (P \cup P') /\varphi$ an Anosov \vf{} $Z$.

Let us show that the graph $G = G(P \cup P', Y \cup Y', \varphi)$ is strongly connected.
By transitivity of $(P,Y)$, the set $\Lambda$ is a basic piece of $(P,Y)$, hence a vertex of $G$.
By iterating Item~\ref{prop: add compact leaf; it: smale graph} of Proposition~\ref{prop: add compact leaf}, the Smale's graph of $(P',Y')$ is formed by a vertex $\Lambda_0$ (corresponding to the maximal invariant set of a copy of $(P,Y)$) and $2p$ vertex $O_1, \dots, O_{2p}$, and an oriented edge $\Lambda_0 \rightarrow 0_i$ for $i$ even, and an oriented edge $O_i \rightarrow \Lambda_0$ for $i$ odd (up to change the numbering), where each $O_i$ is a saddle periodic orbit. 
We have the following claim.
\begin{claim} \label{claim: G(P cup P', Y cup Y', varphi)}
    The graph $G$ is the union of the Smale's graph of $(P,Y)$ and of $(P',Y')$ where we add the following edges:
\begin{itemize}
    \item $\Lambda \rightarrow \Lambda_0$ and $\Lambda_0 \rightarrow \Lambda$,
    \item $\Lambda \rightarrow O_i$ for $i$ odd,
    \item $O_i \rightarrow \Lambda$ for $i$ even.
\end{itemize}
\end{claim}
\begin{proof}
We refer to Figure~\ref{fig: graph add cpct leaf}.

\begin{figure}[htb]
    \centering
    \vspace*{-1em}
    \includegraphics[height=0.29\textheight]{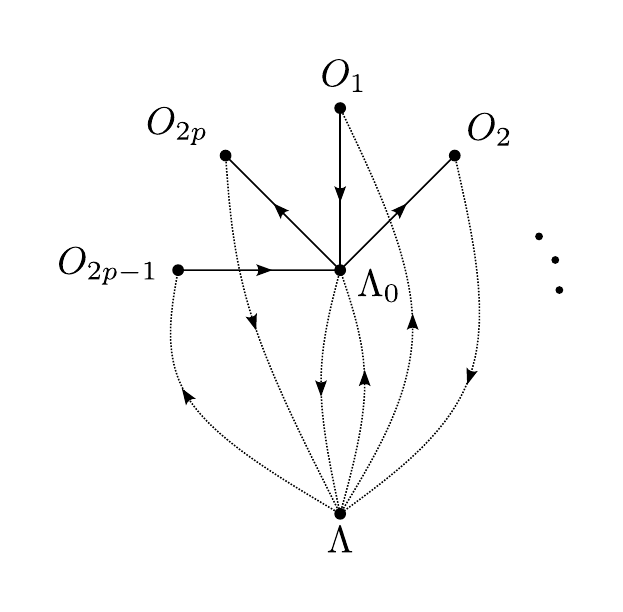}
    \vspace*{-1em}
    \caption{Graph $G(P\cup P', Y\cup Y', \varphi)$ with edges added by the gluing map in dashed line}
    \label{fig: graph add cpct leaf}
\end{figure}

If $O_i$ is an isolated saddle periodic orbit of $\Lambda'$ such that $O_i \rightarrow \Lambda_0$ (i.e., if $i$ odd), then its stable manifold intersect $\partial P'$ in a compact leaf $\gamma^s_i$.
According to Figure~\ref{fig: sqt gluing on elementary coherent}, and because the boundary laminations of both blocks are filling, we know that the image $\varphi (\gamma^s_i)$ intersect $\cW^u(\Lambda)$.
We deduce that there is an edge $\Lambda \rightarrow O_i$ for $i$ odd in the graph $G$.
The third item is deduced the same way.
Finally, the stable (unstable) manifold of $\Lambda_0$ intersects each annuli in the entrance (exit) boundary of $(P',X')$ along non-compact leaves.
Because the lamination is filling in both block, it follows that $\varphi (\cW^s(\Lambda_0))$ intersect $\cW^u(\Lambda)$ and $\varphi (\cW^u(\Lambda_0))$ intersect $\cW^s(\Lambda)$.
The first item follows.
\end{proof}

It is now easy to see from Claim~\ref{claim: G(P cup P', Y cup Y', varphi)} that the graph $G$ is strongly connected, and the Anosov flow $Z$ is transitive according to Proposition~\ref{prop: transitivity criterion}.

The manifolds $P$ and $P'$ are both embedded in $\cN$, by an embedding which realizes an orbit equivalence between the vector fields $Y$ and $Y'$ and the corresponding restriction of the vector field $Z$.
By iterating Proposition~\ref{prop: add compact leaf}, we can show that there exists a torus $T$ embedded in $P'$ which cuts in $P'$ a submanifold homeomorphic to $P$ and a submanifold which admits a Seifert fibration whose fiber is homotopic to the periodic orbits contained in $\partial P'$.
Indeed, in the second iteration of Proposition~\ref{prop: add compact leaf}, we glue a Seifert piece (the elementary block $(V_e, X_e)$) to a Seifert piece whose regular fiber is homotopic to the periodic orbits (this is the submanifold adjacent to the boundary given by Item~\ref{prop: add compact leaf; it: topology} of Proposition~\ref{prop: add compact leaf} applied once), and according to the construction we glue fiber to fiber.
The union forms a unique Seifert piece whose fiber is homotopic to the periodic orbits contained in $\partial P'$, and we continue by recurrence.
Moreover the torus $T$ is isotopic to a torus \qt{} in the vector field $Y'$ (Remark~\ref{rmk: isotopic qt torus}).

Let $T_1$ be the projection of $\partial P$ in $\cN$, and $T_2$ the projection of the torus $T \subset \intr P'$ in $\cN$.
The torus $T_1$ and $T_2$ are disjoint and both isotopic to a torus quasi-transverse to the Anosov field $Z$ and therefore incompressible (\cite{barbotMisePositionOptimale1995}, \cite{brunellaSeparatingBasicSets1993}).
They cut two submanifolds $P_1$ and $P_2$ both homeomorphic to $P$, and such that the vector field $Z$ on $P_1$ and on $P_2$ is orbit equivalent to $Y$ on $P$, and a third component $P_3$ which is a Seifert fibration, bounded by the torus $T_1$ and $T_2$.
The fiber of the fibration of $P_3$ is homotopic to a periodic orbit of the vector field, in other word $P_3$ is a periodic Seifert fibered space.
\end{proof}

\begin{rmk}
For a generalization to a collection of periodic orbits the proof is the same, but we have to justify that we can choose a gluing map which reverses the orientation.
We construct a building block $(P,X)$ by Double-blow up -- Excision on a collection of $n$ periodic orbits, this block has $n$ boundary components.
Let $(P', X')$ be the block associated by the (iterated) Proposition~\ref{prop: add compact leaf} where we have added a compact leaf between each periodic orbit of the lamination, with a dynamical orientation reversed to those of the periodic orbits.
We pair naturally the boundary components of $P$ and $P'$ and there is a \sqt{} gluing map $\varphi : \pP \to \pP'$ given by Figure~\ref{fig: sqt gluing on elementary coherent} which realizes this pairing.
Note that if $P$ has an orientation and $P'$ has the reversed orientation (the one which coincides with the reversed orientation of $P$ seen as a submanifold of $P'$), then the \diff{} $\varphi$ reverses the orientation (on each boundary component).
It follows that the glued manifold is orientable.
\end{rmk}

\subsection{Knot complement in \texorpdfstring{$S^3$}{S\^3} as atoroidal JSJ piece}
\label{sec: orbit complement; subsec: knots}
A knot is an embedding of the circle $\S^1$ to the 3-sphere $\S^3$.
A link is an embedding of $n$ copies of the circle $\S^1$ to $\S^3$.
We usually confuse a knot (link) with its \emph{isotopy class} in $\S^3$.
We refer to \cite{adamsKnotBookElementary2004} for a general reference.
Theorem~\ref{thm: orbit complement and JSJ piece} allows us to realize the complementary of some periodic orbits of transitive (pseudo-)Anosov flow as JSJ pieces of transitive Anosov flow.
One can then ask if there exist knots $K$ in $\S^3$ whose complementary $\S^3 \ssm K$ is homeomorphic to the complementary of a periodic orbit $\gamma$ of a transitive Anosov flow or of the unique singular orbit of a transitive pseudo-Anosov flow.
We are interested in this question for knots whose complementary in $\S^3$ is atoroidal.
Such knots are called \emph{hyperbolic}.

The Gordon-Luecke theorem states that the isotopy class of a knot $K$ is determined by its complementary in $\S^3$, i.e., if $K$ and $K'$ are two knots such that the complements $\S^3 \ssm K$ and $\S^3 \ssm K'$ are homeomorphic, then $K$ and $K'$ are isotopic.
Knot theory is a very rich theory for low dimensional topology.
Let us quote for example the Lickorish-Wallace theorem which states that any connected closed orientable manifold of dimension 3 is obtained by a \emph{Dehn surgery} on a link $L$ in $\S^3$ (\cite{adamsKnotBookElementary2004}).

\subsubsection*{Figure eight knot}
\begin{figure}[tb]
    \centering
    \vspace*{-2em}
    \includegraphics[height=0.265\textheight]{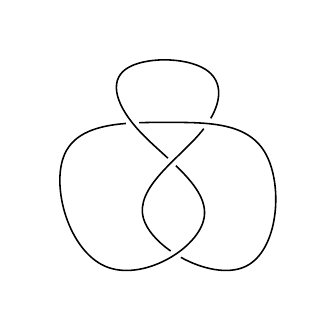}
    \vspace*{-1.5em}
    \caption{The figure eight knot}
    \label{fig: figure eight knot}
\end{figure}

It is a known fact (\cite{thurstonThreedimensionalGeometryTopology1997}) that the complementary of the figure eight knot $K_8$ (Figure~\ref{fig: figure eight knot}) is atoroidal, and homeomorphic to the complementary of the periodic orbit $\gamma$ of the matrix suspension flow
$$ A := \begin{pmatrix} 2 & 1 \\ 1 & 1 \end{pmatrix} : \T^2 \to \T^2 $$
induced by the fixed point $(0,0)$.
The suspension field $X_A$ is an Anosov vector field.
In \cite{franksAnomalousAnosovFlows1980}, the authors construct an Anosov flow by gluing two copies of a building block homeomorphic to the figure eight knot complement $\S^3 \ssm K_8$, obtained by blow-up and excision on the periodic orbit $\gamma$ of $X_A$.
This Anosov flow has the particularity of not being transitive.
Theorem~\ref{thm: orbit complement and JSJ piece} gives the existence of a transitive Anosov flow whose JSJ decomposition is formed by two copies of $\S^3 \ssm K_8$, adjoined by a periodic Seifert piece $P$.
Note that the periodic Seifert piece can be described completely: it is the gluing of two copies of the elementary intermediate block $(V_e, X_e)$ (Definition~\ref{lem: elementary block}).

\subsubsection*{Hyperbolic fibered knots}
Let $K$ be a knot in $\S^3$.
We say that a compact orientable surface $S$ embedded in $\S^3$ is a \emph{Seifert surface} for $K$ if $\partial S =K$.
We say that $K$ is a \emph{fibered knot} if the complementary $\S^3 \ssm K$ admits a fibration over the circle with fiber a Seifert surface $S$ of $K$.
The \emph{monodromy} of a fibered knot is the monodromy $h: S \to S$ of such a fibration, well-defined up to isotopy relative to $\partial S$.
A diffeomorphism $h: S \to S$ is said to be \emph{pseudo-Anosov} if there exists a pair of foliations with singularities $(\cF^s, \cF^u)$, transverse to each other, and provided with a transverse measures, which are uniformly contracted and expanded by $h$ respectively.
The laminations have a finite number of \emph{$p$-prong singularities} with $p \geq 3$ (Figure~\ref{fig: p prong}).
The suspension of a pseudo-Anosov diffeomorphism of a closed surface is a pseudo-Anosov flow in the sense of Definition~\ref{def: pseudo anosov flow}.
If the monodromy $h$ of a fibered knot $K$ admits a pseudo-Anosov representative, then the knot $K$ is a hyperbolic knot and the complementary of $K$ is homeomorphic to the complementary of a periodic orbit $\gamma$ of the pseudo-Anosov suspension flow of $\hat h: \hat S \to \hat S$, where $\hat h$ is the pseudo-Anosov diffeomorphism induced by $h$ on the surface $\hat S$ obtained by blowing down the boundary of $S$ (\cite{thurstonGeometryTopology3Manifolds1979}).
The boundary of $S$ induces a fixed point $\hat p$ of $\hat h$, and $\gamma$ is the suspension of $\hat p$.
If $h$ admits no singularities on the interior of $S$, then $\gamma$ is the unique singularity of the suspension pseudo-Anosov flow of $h$, and it follows as a corollary of Theorem~\ref{thm: orbit complement and JSJ piece}.

\begin{coro} \label{coro: hyperbolic knot no interior singularity}
If $K$ is a fibered knot whose monodromy is pseudo-Anosov and admits no singularity in the interior of the fiber, then the complementary $\S^3 \ssm K$ is an atoroidal JSJ piece of a transitive Anosov flow.
\end{coro}

We are looking for a family of knots which satisfies the hypotheses of Corollary~\ref{coro: hyperbolic knot no interior singularity}.

\subsubsection*{Plumbings of the Seifert surface of the figure eight knot}
I. Agol pointed out to me the following construction, allowing to exhibit a family of knots satisfying the hypotheses of Corollary~\ref{coro: hyperbolic knot no interior singularity}, which are the knots obtained by \emph{plumbings of the Seifert surface of the figure eight knot}.
Let us explain the ideas briefly.

\begin{figure}[htb]
    \centering
    \vspace*{-1em}
    \includegraphics[width=\textwidth]{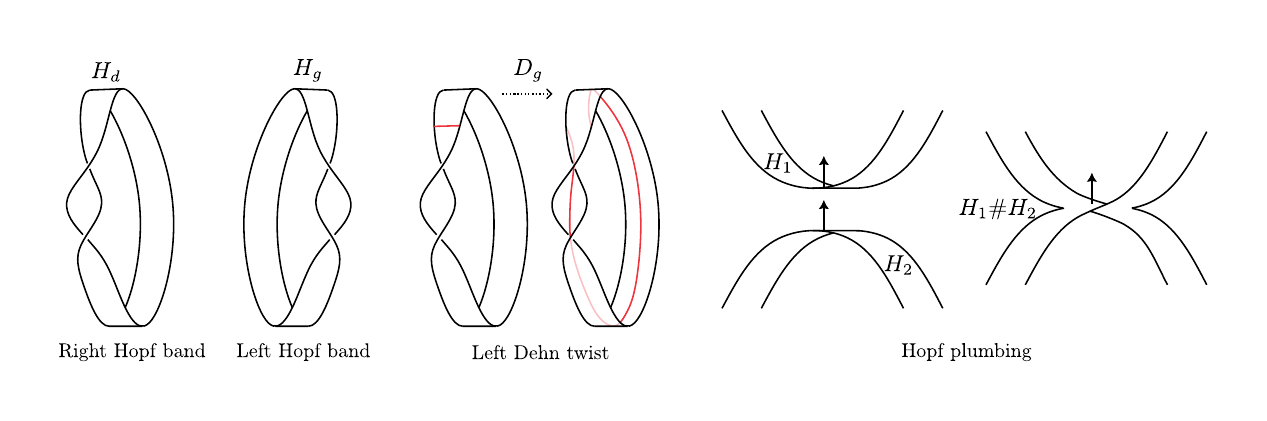}
    \vspace*{-2em}
    \caption{Dehn twist and Hopf plumbing}
    \label{fig: hopf band}
\end{figure}

In \cite{gabaiPseudoAnosovMapsSurgery1990}, the authors describe a construction by \emph{Hopf plumbing} allowing to compute the invariant foliations of the monodromy of a hyperbolic fibered \emph{2-bridge} knot.
Let $H_g$ be the left Hopf band, and $H_d$ the right Hopf band described in Figure~\ref{fig: hopf band}.
The left Hopf link $\partial H_g$ is a fibered knot with monodromy a right Dehn twist $D_d$ along the core, and the right Hopf link $\partial H_d$ is a fibered knot with monodromy a left Dehn twist $D_g$ (Figure~\ref{fig: hopf band}).
The \emph{Hopf plumbing} is the surface $S:= H_1 \# H_2$ obtained by gluing two unlinked Hopf bands in $\S^3$ together along a rectangular neighborhood of an arc which crosses the band from boundary to boundary (Figure~\ref{fig: hopf band}), with the orientation convention given by Figure~\ref{fig: hopf band}.
If $D_i$ is the monodromy of $H_i$, then $\partial (H_1 \# H_2)$ is a fibered knot of fiber $H_1 \# H_2$ and monodromy the product of $D_i$ (\cite{stallingsConstructionsFibredKnots1978}).

Consider $S:= H_g^1 \# H_d^1 \# \dots \# H_g^n \# H_d^n$ a $2n$ Hopf plumbing of alternating left band and right band, and denote $K:= \partial S$.

\begin{lem}
$K$ satisfies the assumptions of Corollary~\ref{coro: hyperbolic knot no interior singularity}, i.e., $K$ is a fibered knot whose monodromy is pseudo-Anosov and admits no singularity on the interior of the fiber.
\end{lem}

\begin{proof}[Proof idea]
The surface $S$ is a surface of genus $n$, and the cores of the Hopf bands are the curves drawn in Figure~\ref{fig: monodromy and covering}.
By construction, the boundary $K:= \partial S$ is a fibered knot of monodromy the product $h = D_d^1 D_g^1 \dots D_d^n D_g^n$ of $2n$ Dehn twists alternating right and left along the core of each of the Hopf strips.
The surface $S$ is a second order branching covering of the disk $D$ which is a quotient of $S$ by the involution around the horizontal axis of Figure~\ref{fig: monodromy and covering}, provided with $2n +1$ marked points corresponding to the branching points.
The cores of the Hopf bands are arcs connecting the marked points, and the \diff{} $h :S \to S$ induces a \diff{} $g : D \to D$, which is a product of \say{half}-twists of Dehn along the arcs, alternated on the left and on the right, denoted $\alpha_g^i$ and $\alpha_d^i$.

\begin{figure}[htb]
    \centering
    \vspace*{-1em}
    \includegraphics[width=0.9\textwidth]{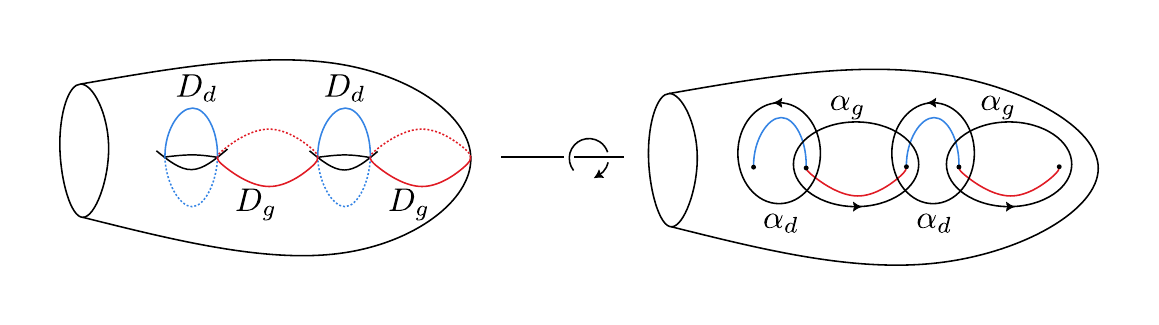}
    \vspace*{-1em}
    \caption{Monodromy on $S$ and covering $S \to D$}
    \label{fig: monodromy and covering}
\end{figure}

In \cite{gabaiPseudoAnosovMapsSurgery1990}, the authors use the method of \emph{invariant train tracks} (\cite{pennerCombinatoricsTrainTracks1992}) to construct a $g$-invariant {train track} obtained from the invariant \emph{subtracks} $\tau_g$ and $\tau_d$ of the half-twists $\alpha_g^i$ and $\alpha_d^i$.
They are described in Figure~\ref{fig: train track}.
The invariant {train track} $\tau$ is a recollection of the subtracks $\tau_g$ and $\tau_d$ which overlap on monogons (Figure~\ref{fig: train track}).

\begin{figure}[htb]
    \centering
    \vspace*{-1em}
    \captionsetup{width=.83\linewidth}
    \includegraphics[width=0.9\textwidth]{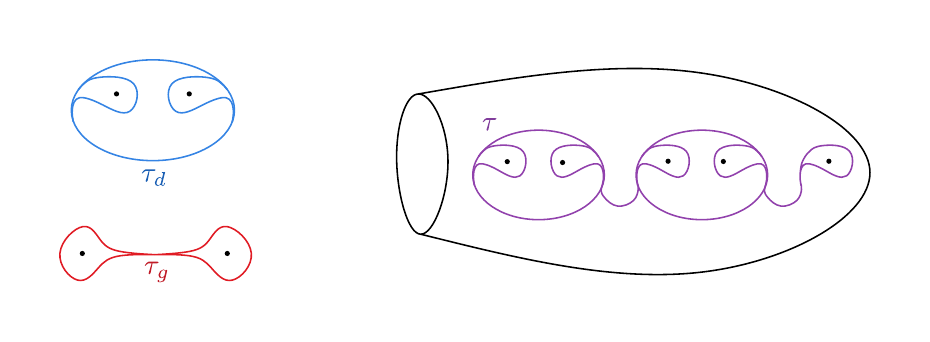}
    \vspace*{-1em}
        \caption{Train track $\tau$ as a recollection of subtracks $\tau_g$ and $\tau_d$}
    \label{fig: train track}
\end{figure}

The authors show (\cite[Proposition~4]{gabaiPseudoAnosovMapsSurgery1990}) that $g$ is isotopic relative to $\partial D$ to a pseudo-Anosov homeomorphism of $\hat g$ of $D$. Moreover if we denote $\lambda^s$ the stable invariant lamination of $\hat g$ on $D$, then a $n$-gone of $D \ssm \tau$ corresponds to a disk with $n$ cusps of $D \ssm \lambda^s$.
We then observe that, in the particular situation of an alternating left and right Hopf plumbing, these connected components are disks with one or two cusps, and it follows that the stable foliation with singularities $f^s$ of the pseudo-Anosov homeomorphism $\hat g$ that extends the $\lambda^s$ lamination admits only singularities of type $1$-prong.

Passing to the $2$ order covering on $S$, it follows that (up to isotopy relative to $\partial S$), the monodromy $f$ of the knot $K$ is pseudo-Anosov, and its foliations have no singularities on the interior $S$.
\end{proof}

Finally, let us notice that the Seifert surface of the figure eight knot is a plumbing of a right Hopf band and a left Hopf band.
We deduce:
\begin{coro} \label{coro: plumbing of figure eight knot}
Let $K = \partial (S_1 \# S_2 \# \dots \# S_n)$ be a plumbing of $n$ copies of the Seifert surface of the figure eight knot.
Then $K$ is a fibered hyperbolic knot and the complementary of $K$ is an atoroidal JSJ piece of a transitive Anosov flow.
\end{coro}

\section{Gluing pieces of skewed \texorpdfstring{$\R$}{R}-covered Anosov flows}
\label{sec: gluing skewed blocks}

In this section, we show a general result which allows us to glue pieces cut along incompressible tori in a skewed $\R$-covered Anosov flow on a toroidal 3-manifold.
An Anosov flow on a 3-manifold $\cM$ is said to be \emph{$\R$-covered} if the leaf space of the lifted stable foliations $\widetilde \cF^s$ on the universal cover $\widetilde \cM$ is separated (hence homeomorphic to $\R$).
It is said to be \emph{skewed $\R$-covered} if it is moreover not orbit equivalent to a suspension.
The proposition below allows us to cut building blocks out of a a skewed $\R$-covered Anosov flow along a collection of incompressible tori. 

\begin{prop}[{\cite[Theorem~A'~and~Theorem~E]{barbotMisePositionOptimale1995}}]
\label{prop: position tori skewed anosov}
Let $Z$ be a skewed $\R$-covered Anosov vector field on a closed orientable 3-manifold $\cM$, whose stable and unstable foliations are transversely orientable.
Let $\cT = \{ T_1, \dots, T_n \}$ be a finite collection of incompressible tori embedded in $\cM$, pairwise disjoint and pairwise non-homotopic.
Then there exists a collection $\cT' = \{ T'_1, \dots, T'_n \}$ of pairwise disjoint tori isotopic to $\cT$ and quasi-transverse to $Z$, and this collection is unique up to homotopy along the orbits of the flow.
As a consequence, if we set $P := \cM \ssm \cT'$, then $(P, \res{Z}{P})$ is a building block.
\end{prop}

\begin{proof}
According to \cite[Theorem~A']{barbotMisePositionOptimale1995} any incompressible torus $T_i$ embedded in an $\R$-covered Anosov flow in an orientable 3-manifold, and whose foliations are transversely oriented is isotopic to a quasi-transverse torus $T'_i$.
The torus $T'_i$ is unique up to homotopy along the orbits.
According to \cite[Theorem~E]{barbotMisePositionOptimale1995}, we can choose the tori $T'_1, \dots T'_n$ pairwise disjoint.
If we set $P := \cM \ssm \bigcup_i T'_i$, it is clear that the pair $(P, X) := (P, \res{Z}{P})$ is a \bb{}.
\end{proof}

\begin{defi}[Skewed $\R$-covered Anosov block]
In the setting of Proposition~\ref{prop: position tori skewed anosov}, we call a \emph{skewed $\R$-covered Anosov block} any union of \ccs{} of $(P,\res{Z}{P})$.
\end{defi}

We have the following result.

\begin{mainthm}[Theorem~\ref{thmintro: gluing skewed blocks}] \label{thm: gluing skewed blocks}
Let $(P,X)$ and $(P', X')$ be two skewed $\R$-covered Anosov blocks and $\varphi \colon \pP \to \pP'$ a gluing map.
There is a gluing map $\psi$ isotopic to $\varphi$ among gluing maps such that the \vfs{} $X$ and $X'$ induce an Anosov vector field $Z$ on the closed 3-manifold $\cN:= P \cup P' / \psi$.
\end{mainthm}

Recall that skewed $\R$-covered Anosov flow blocks are orientable, and we say that a \diff{} $\varphi \colon \pP \to \pP'$ reverses the orientation if there exists an orientation of $P$ and $P'$ such that $\varphi$ reverses the canonical orientation induced on the boundary $\pP$ and $\pP'$.

\begin{rmk}
Any \emph{piece of (finite cover of) geodesic flow} is a skewed $\R$-covered Anosov block.
This result then generalizes Handel-Thurston or Clay-Pinsky gluing surgeries, and generalizations by T. Barbot and S. Fenley and gets rid of the \emph{positivity} constraint on the isotopy class of the gluing map.
It allows us to use of much more general blocks than geodesic flow pieces.
For example, any Anosov flow obtained by Dehn-Goodman-Fried surgeries \emph{of coherent orientations} on a suspension or a geodesic flow is skewed $\R$-covered (\cite{fenleyAnosovFlows3manifolds1994}).
It also allows us to glue blocks $(P,X)$ and $(P',X')$ which are not cut in the same Anosov flow.
\end{rmk}

Note that the statement does not require any assumption on the boundary lamination and the gluing map.
We will see (Lemma~\ref{lem: cns sqt gluing skewed blocks}) that the existence of a \sqt{} gluing map $\varphi \colon \pP \to \pP'$ is equivalent to a simple criterion of matching boundary components containing the same number of periodic orbits.

\subsection{Boundary lamination}

The proof relies on the following key lemma.

\begin{lem} \label{lem: lamination on skewed blocks}
Let $(P,X)$ be a skewed $\R$-covered Anosov block.
Then $(P,X)$ is a filled block, and the boundary lamination $\cL$ of $(P,X)$ is alternating elementary (Definition~\ref{def: elementary coherent alternating}), i.e., the only compact leaves of $\cL$ are the periodic orbits of $X$, and two successive periodic orbits have an orientation given by the flow which is opposite.
\end{lem}

\begin{figure}[htb]
    \centering
    \vspace*{-2em}
    \includegraphics[height=0.27\textheight]{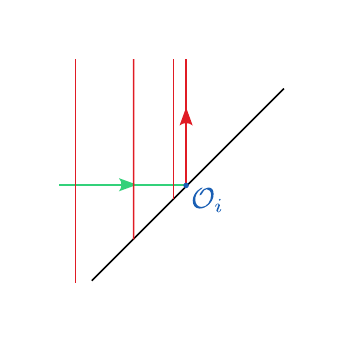}
    \vspace*{-2em}
    \caption{A periodic orbit $\cO_i$ which borders a component of $\pP \ssm \cL$ has a free separatrix}
    \label{fig: boundary orbit free separatrix}
\end{figure}

\begin{proof}
Let $Z$ be a skewed $\R$-covered Anosov vector field on a 3-manifold $\cM$ which contains the block $(P,X)$.
Let $T_i$ be a torus quasi-transverse at $Z$ embedded in $\cM$ corresponding to a \cc{} of $\pP$, and let $\cO_*$ be the collection of periodic orbits of $Z$ contained in $T_i$.
T. Barbot shows in \cite{barbotMisePositionOptimale1995} that $T_i$
is a gluing of elementary Birkhoff annuli, i.e., the traces of the stable and unstable foliations $\cF^s$ and $\cF^u$ of the Anosov \vf{} on the closure of each of the connected components of $T_i \ssm \cO_*$ are foliations of dimension 1 with no compact leaves other than the periodic orbits of the boundary and with no Reeb components.
In other words, the traces of the foliations $\cF^s$ and $\cF^u$ on $T_i$ are elementary alternating.
The periodic orbits bordering a connected component of $T_i \ssm \cO_*$ are then oppositely oriented by the flow of $X$, in other words they are freely homotopic to the inverse of the other.
The boundary lamination $\cL$ of the block $(P,X)$ is a sublamination of the \qms{} foliations induced by the pair $(\cF^s, \cF^u)$ on $T_i$, it follows that the lamination $\cL$ is elementary alternating, and so on any boundary component $T_i$.
Let us show that this is a filling lamination.
It is a prefoliation with no compact leaves other than the periodic orbits.
If it is not filling, there exists an annulus $A$ in, say, $\Pin$ bounded by two periodic orbits $\cO_i, \cO_j$, and disjoint from the lamination $\cL$ on its interior (Proposition~\ref{prop: complementary of a prefoliation}).
The periodic orbit $\cO_i$ thus admits an unstable separatrix $\cW^u_+(\cO_i)$, which is disjoint from the lamination $\cW^s$, and thus from $\Lambda$.
Indeed, an orbit of $\Lambda \cap (\cW^u_+(\cO_i))$ in a linearizing neighborhood of $\cO_i$ has its stable manifold intersecting the boundary of $P$ on the annulus $A$, which is a contradiction (Figure~\ref{fig: boundary orbit free separatrix}).
It follows that (the unique) unstable separatrix of $\cO_i$ in $P$ is free.
It intersects $\partial P$ along a compact leaf of lamination $\cL^\out$, which is not a periodic orbit.
This is a contradiction.
\end{proof}

\subsection{Criterion for the existence of a \sqt{} gluing map}
We show the following lemma, which gives a simple necessary and sufficient criterion for the existence of a \sqt{} gluing map.

\begin{lem} \label{lem: cns sqt gluing skewed blocks}
Let $(P,X)$ and $(P', X')$ be two oriented skewed $\R$-covered blocks.
Then the following are equivalent:
    \begin{enumerate}
        \item \label{lem: cns, it: sqt}
        There exists a strongly quasi-transverse gluing map $\varphi \colon \pP \to \pP'$.
        \item \label{lem: cns, it: gluing map}
        There is a gluing map $\varphi \colon \pP \to \pP'$.
        \item \label{lem: cns, it: pair}
        There is a pairing of boundary components of $P$ and $P'$ containing the same number of periodic orbits.
    \end{enumerate}

\end{lem}

Recall that $\varphi \colon \pP \to \pP'$ is a gluing map if $\varphi$ maps the oriented periodic orbits of $X$ in $\pP$ to the oriented periodic orbits of $X'$ in $\pP'$ and the exit boundary $\Pout$ to the entrance boundary $P'^\iin$.
The condition of Item~\ref{lem: cns, it: pair} means that $P$ and $P'$ have the same number of boundary components, denoted $T_1, \dots, T_n$ and $T_1', \dots, T_n'$, and for each $i= 1, \dots, n$ we require that there are as many periodic orbits of $X$ in $T_i$ as there are periodic orbits of $X'$ in $T_i'$.

\begin{figure}[htb]
    \centering
    \vspace*{-2em}
    \includegraphics[height=0.25\textheight]{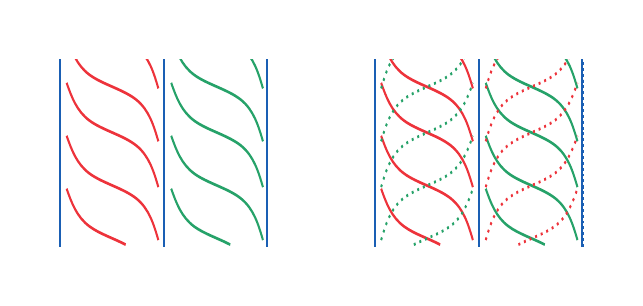}
    \vspace*{-2em}
    \caption{An alternating elementary lamination and a pair of strongly quasi-transverse alternating elementary laminations}
    \label{fig: alternating elementary lamination}
\end{figure}

\begin{proof}
Denote $P^\iin$ and $P^\out$ the entrance and exit boundary of $(P,X)$, $\cO_*$ the periodic orbits of $X$ contained in $\pP$, and $\cL$ the boundary lamination of $(P,X)$.
Similarly, denote $P'^\iin$ and $P'^\out$ the entrance and exit boundary of $(P',X')$, $\cO'_*$ the periodic orbits of $X'$ contained in $\pP'$, and $\cL'$ the boundary lamination of $(P',X')$.
 
It is clear that $\ref{lem: cns, it: sqt} \Rightarrow~\ref{lem: cns, it: gluing map} \Rightarrow~\ref{lem: cns, it: pair}$.
Conversely, suppose we have a pairing of boundary components $\{(T_i, T'_i)\}$ of $P$ and $P'$ respectively, containing the same number of periodic orbits.
According to Lemma~\ref{lem: lamination on skewed blocks}, the lamination $\cL$ on $T_i$ and $\cL'$ on $T_i'$ are elementary alternating, and according to the matching conditions, they have the same number of marked leaves.
Two such laminations can be glued together in a strongly quasi-transverse way (Figure~\ref{fig: alternating elementary lamination}), and we have for each $i=1, \dots, n$ the existence of a \diff{} $f_i \colon T_i \to T'_i$ which maps the oriented periodic orbits to the oriented periodic orbits, the exit boundary $P^\out$ to the entrance boundary $P'^\iin$, and the lamination $\cL$ to a strongly quasi-transverse lamination at $\cL'$.
The product $\varphi \colon \partial P \to \partial P'$ thus defines a \sqt{} gluing map of building blocks.
\end{proof}

\subsection{Proof of Theorem~\ref{thmintro: gluing skewed blocks}}
We will need the following general lemma.

\begin{lem} \label{lem: twist is sqt}
Let $(P,X)$ be a \bb{} and $\varphi$ be a \sqt{} gluing map of $(P,X)$.
Let $T$ be a boundary component which contains a nonzero number of periodic orbits of $X$.
Then for any Dehn twist $\tau \colon T \to T$ along the periodic orbits, there exists a strongly quasi-transverse gluing map isotopic to $\varphi \circ \tau$.
\end{lem}

\begin{proof}
We refer to the end of the proof of Proposition~\ref{prop: normalization of triple}, which shows that we can choose the Dehn twist $\tau$ in a well-chosen support so that $\varphi \circ \tau$ is still a gluing map \sqt{}.
Let us briefly recall the idea.
Denote $\cL$ the boundary lamination and $\cO_*$ the collection of periodic orbits contained in $T$.
It suffices to choose the support of $\tau$ in a compact annulus of $T \ssm \cO_*$ so that the composition is still a gluing map, in other words maps the oriented periodic orbits on the oriented periodic orbits and the exit boundary on the entrance boundary.
We then choose a small neighborhood $\cV$ of a periodic orbit $\cO \in \cO_*$, in which the holonomy of the lamination $\cL$ along a transversal of $\cO$ is contracting on $\Pout$ and expanding on $\Pin$, and the holonomy of the lamination $\varphi^\inv_* \cL$ is expanding on $\Pin$ and contracting on $\Pout$.
Following the sign of the Dehn twist, it suffices to choose the support of $\tau$ in a compact annulus contained either in $\Pin \cap \cV$ or in $\Pout \cap \cV$, so that $\tau_* (\varphi^\inv_* \cL)$ is (strongly) transverse to $\cL$.
\end{proof}

\begin{proof}[Proof of Theorem~\ref{thm: gluing skewed blocks}]
Let $(P,X)$ and $(P',X')$ be two skewed $\R$-covered Anosov\linebreak[4] blocks and $\varphi \colon \pP \to \pP'$ a \diff{} which maps the oriented periodic orbits of $X$ in $\pP$ to the oriented periodic orbits of $X'$ in $\pP'$ and the exit boundary $\Pout$ to the entrance boundary $P'^\iin$. 
It follows from Lemma~\ref{lem: cns sqt gluing skewed blocks} that there exists a \sqt{} gluing map $ \psi \colon \pP \to \pP'$.
The gluing map $\psi$ defines the same pairing of boundary components as $\varphi$.
Let $T_i$ and $T'_i$ be two boundary components of $\pP$ and $\pP'$ matched by $\varphi$ and $\psi$.
Since $\varphi(\cO_*) = \psi(\cO_*) = \cO_*$, it follows that the isotopy class of $\varphi$ and $\psi$ differ by a Dehn twist $\tau \colon T_i \to T_i$ along the periodic orbits $\cO_*$, in other words such that $\psi \circ \tau$ and $\varphi$ are isotopic.
According to Lemma~\ref{lem: twist is sqt}, $\psi \circ \tau$ is isotopic to a \sqt{} gluing map $\psi' \colon \pP \to \pP'$, which is then isotopic to $\varphi$.

We consider the triple $(P \cup P', X \cup X', \psi')$.
The block $(P \cup P', X \cup X')$ is a filled block according to Lemma~\ref{lem: lamination on skewed blocks}, and $\psi'$ is a \sqt{} gluing map.
Up to a strong triple isotopy, we can apply the Gluing Theorem~\ref{thm: gluing theorem} and the \vfs{} $X$ and $X'$ induce on the closed manifold $\cM := P \cup P' /\psi'$ a vector field $Z$ which is Anosov.
It satisfies Theorem~\ref{thm: gluing skewed blocks} by construction.
It is easy to check with the criterion of Proposition~\ref{prop: transitivity criterion} that if the blocks $(P,X)$ and $(P', X')$ are both transitive, the \vf{} $Z$ thus obtained is transitive.
\end{proof}

\begin{figure}[h]
    \centering
    \vspace*{-2em}
    \includegraphics[width=\textwidth]{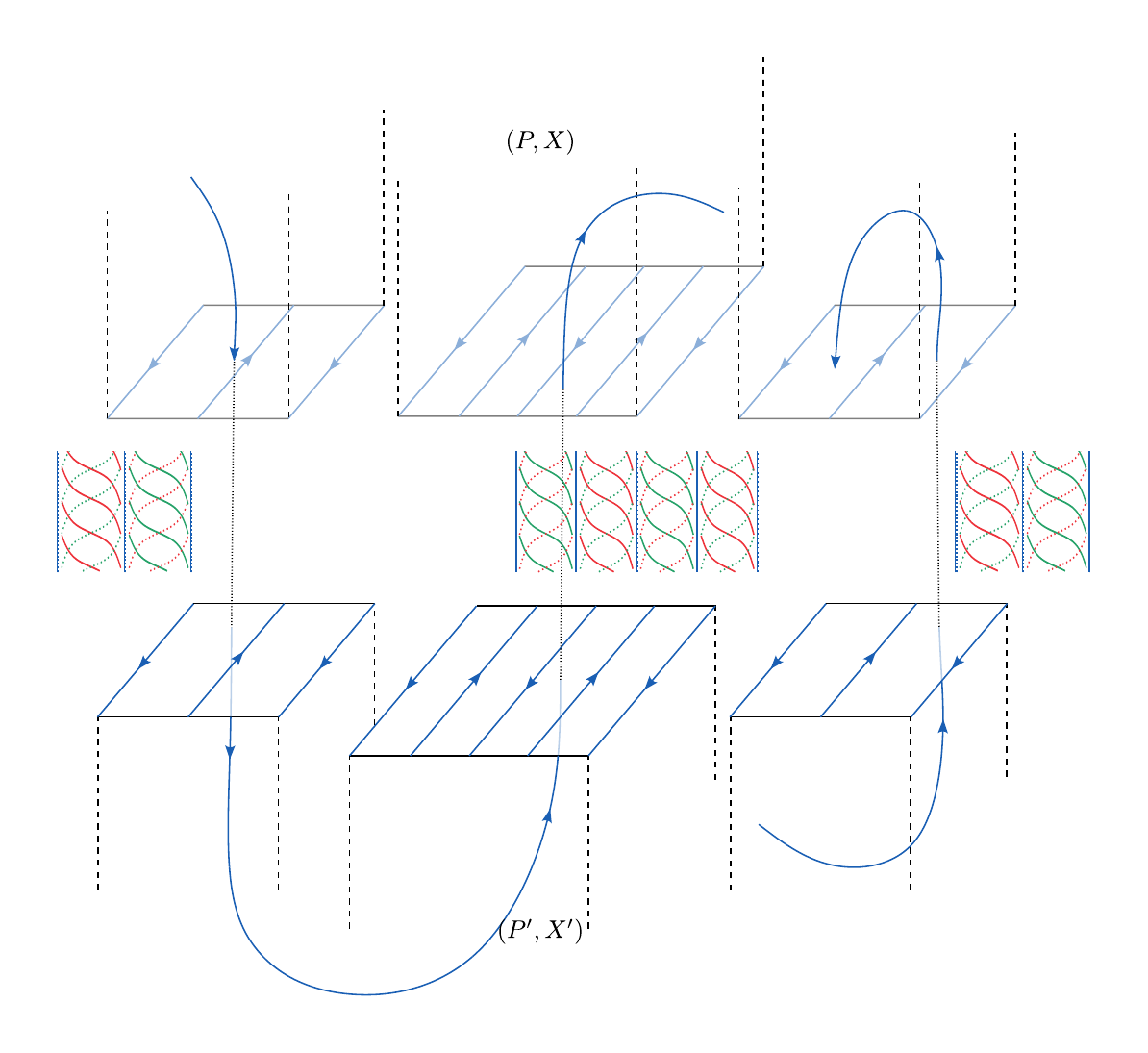}
     \vspace*{-1em}
    \caption{Gluing operation of two skewed $\R$-covered Anosov blocks}
    \label{fig: gluing skewed blocks}
\end{figure}

\subsection*{Acknowledgments}
We address a special thank to Fran\c{c}ois B\'eguin for his proofreading work, his technical and mental support. We also thank Christian Bonatti for his intuition in preliminary stages, Thomas Barthelm\'e for the sharing of many ideas on this work, and Boris Hasselblatt for the support.

\newcommand{\doititle}{}
\def\arXiv#1{\@ifundefined{href}{{\mdseries\ttfamily arxiv:#1}}{\href{http://arxiv.org/pdf/#1}{{\mdseries\ttfamily arXiv:#1}}}}
\let\arxiv\arXiv

\bibliographystyle{plain}

\end{document}